\documentclass[twoside,12pt,a4paper]{book}
\pdfoutput=1

\newcommand{\titleString}{Eisenstein series and automorphic representations}

\usepackage[centering, textwidth=16cm, top=1in, includehead, includefoot]{geometry}

\usepackage{amsmath,amssymb,amsthm,graphics,color,cite}
\usepackage[usenames,dvipsnames]{xcolor}
\usepackage{graphicx}
\usepackage[english]{babel}
\usepackage{fancyhdr}
\usepackage{mdframed}
\usepackage{mathtools}
\usepackage{wrapfig}
\usepackage[textfont=it]{caption}
\usepackage{enumitem}
\usepackage{booktabs}
\usepackage{xparse}

\usepackage{epigraph}
\usepackage{tocloft}
\ifpdf
\usepackage[protrusion=true,expansion=true]{microtype}
\fi

\usepackage{ifthen}
\usepackage{imakeidx}
\indexsetup{othercode=\small}
\makeindex[program=makeindex,columns=2,intoc=true,options=-s indexstyle]

\usepackage[bookmarksnumbered]{hyperref}
\usepackage{currfile}

\usepackage[vcentermath]{youngtab}
\usepackage{tensor}

\usepackage{chngcntr}
\counterwithout{footnote}{chapter}

\newcommand\nonumberfootnote[1]{%
  \begingroup
  \renewcommand\thefootnote{}\footnote{#1}%
  \addtocounter{footnote}{-1}%
  \endgroup
}

\usepackage{tikz}
\newcommand*\circled[1]{\tikz[baseline=(char.base)]{
            \node[shape=circle,draw,inner sep=2pt] (char) {#1};}}

\definecolor{darkred}{rgb}{0.5,0,0.5}
\definecolor{darkgreen}{rgb}{0,0.5,0}

\hypersetup{
    pdftitle={ \titleString },         %
    pdfauthor={ Philipp Fleig, Henrik P. A. Gustafsson, Axel Kleinschmidt, Daniel Persson },      %
    linktocpage=true,           % color toc pages, not text
    colorlinks=true,            % false: boxed links; true: colored links
    linkcolor=blue!75!black,    % color of internal links
    citecolor=green!50!black,   % color of links to bibliography
    filecolor=magenta,          % color of file links
    urlcolor=blue!75!black      % color of external links
}

\numberwithin{equation}{chapter}

\newcommand{\emphindex}[2][]{\emph{#2}%
        \ifthenelse{ \equal{#1} {} } %
        {\index{#2}} %
        {\index{#1}}}

\newcommand{\nn}{\nonumber}

\newcommand{\nats}{\mathbb{N}}
\newcommand{\ints}{\mathbb{Z}}
\newcommand{\reals}{\mathbb{R}}
\newcommand{\cx}{\mathbb{C}}
\newcommand{\rats}{\mathbb{Q}}

\newcommand{\adeles}{\mathbb{A}}
\newcommand{\ads}{\mathbb{A}}
\newcommand{\cmplx}{\mathbb{C}}

\newcommand{\bs}{\backslash}
\newcommand{\height}{\mathrm{ht}}
\newcommand{\DEVIII}[8]{
\ensuremath{\left[\begin{smallmatrix*}[c]
       &      & {#2} & & & & & \\
  {#1} & {#3} & {#4} & {#5} & {#6} & {#7} & {#8}
\end{smallmatrix*}\!\!\right]}
}
\newcommand{\DEVII}[8]{
\ensuremath{\left[\begin{smallmatrix*}[c]
       &      & {#2} & & & & & \\
  {#1} & {#3} & {#4} & {#5} & {#6} & {#7} 
\end{smallmatrix*}\!\!\right]}
}

\newcommand{\DDIV}[4]{
\ensuremath{\left[\begin{smallmatrix*}[c]
      & {#4} & \\
 {#1} & {#2} & {#3}
\end{smallmatrix*}\right]}
}
\newcommand{\beq}{\begin{equation}}
\newcommand{\eeq}{\end{equation}}
\newcommand{\beqa}{\begin{eqnarray}}
\newcommand{\eqa}{\end{eqnarray}}

\DeclareMathOperator{\WF}{WF}
\newcommand{\cO}{{\mathcal{O}}}

\DeclareMathOperator{\Ind}{Ind}
\newcommand{\ba}{{\bar{a}}}

\newcommand{\br}{{\bar{r}}}
\newcommand{\normord}[1]{:\mathrel{#1}:}

\DeclareMathAlphabet{\mathbbold}{U}{bbold}{m}{n}  
\newcommand{\id}{\mathbbold{1}}

\newcommand{\iso}{\ensuremath{\mathrel{\cong}}}

\newcommand{\prodp}{\sideset{}{'}\prod}

\NewDocumentCommand{\intl}{g}{
    \IfNoValueTF{#1}{\int\limits}{\int\limits_{\mathclap{#1}}}
}
\NewDocumentCommand{\suml}{g}{
    \IfNoValueTF{#1}{\sum\limits}{\sum\limits_{\mathclap{#1}}}
}
\NewDocumentCommand{\prodl}{g}{
    \IfNoValueTF{#1}{\prod\limits}{\prod\limits_{\mathclap{#1}}}
}

\newcommand{\lb}{\left[}
\newcommand{\rb}{\right]}

\newcommand{\abs}[1]{\left|#1\right|}

\newcommand{\wlong}{w_{\textrm{long}}}
\newcommand{\Weyl}{\mathcal{W}}
\newcommand{\eps}{\epsilon}

\newcommand{\supp}{\mathrm{supp}}
\newcommand{\stab}{\mathrm{stab}}
\newcommand{\UHP}{\mathbb{H}}

\newcommand{\wt}{\mathrm{wt}}
\newcommand{\BZL}{\mathrm{BZL}}
\newcommand{\field}{\mathbb{F}}
\newcommand{\fieldext}{\mathbb{E}}
\newcommand{\ells}{\ell_{\mathrm{s}}}
\newcommand{\gs}{g_{\mathrm{s}}}
\newcommand{\idm}{\id}%{1\!\!1}
\newcommand{\compl}{\star}
\newcommand{\gammaE}{\gamma_{\mathrm{E}}}

\newcommand{\st}{\ifnum\currentgrouptype=16 \mathrel{}\middle|\mathrel{}\else\mid\fi}
\newcommand{\gra}[2]{{\scriptscriptstyle (#1 , #2 )}}

\newcommand{\ordd}[1]{{\scriptscriptstyle (#1)}}

\newcommand{\lint}{\int\limits}

\newcommand{\mf}[1]{{\mathfrak{#1}}}
\newcommand{\lie}[1]{{\mathfrak{#1}}}
\newcommand{\emp}[1]{{\em #1}}

\makeatletter
\let\Re\@undefined
\let\Im\@undefined
\makeatother
\DeclareMathOperator{\Re}{Re}
\DeclareMathOperator{\Im}{Im}

\DeclareMathOperator{\GKdim}{GKdim}
\DeclareMathOperator{\Vol}{Vol}
\DeclareMathOperator{\sgn}{sgn}

\makeatletter
\def\th@plain{%
  \thm@notefont{}% same as heading font
  \itshape % body font
}
\def\th@definition{%
  \thm@notefont{}% same as heading font
  \normalfont % body font
}
\makeatother

\newtheorem{theorem}{Theorem}[chapter]
\newtheorem{lemma}[theorem]{Lemma}
\newtheorem{proposition}[theorem]{Proposition}
\newtheorem{corollary}[theorem]{Corollary}

\newtheorem{conjecture}[theorem]{Conjecture}

\theoremstyle{definition}
\newtheorem{definition}[theorem]{Definition}
\newtheorem{remark}[theorem]{Remark}

\mdfdefinestyle{examplestyle}{topline=false,bottomline=false,rightline=false,linecolor=gray!50,linewidth=2pt, font=\footnotesize, frametitleaboveskip=10pt, innerbottommargin=10pt, usetwoside=false, backgroundcolor=gray!7, % skipabove=12pt, skipbelow=12pt,
    leftmargin=\dimexpr-\mdflength{linewidth}-\mdflength{innerleftmargin},
    rightmargin=\dimexpr-\mdflength{innerrightmargin}}

\mdtheorem[style=examplestyle]{example}[theorem]{Example}[chapter]

\addto\captionsenglish{}

\graphicspath{ {gfx/} }
\makeatletter
\def\input@path{ {gfx/} }
\makeatother

\begin{document}

\frontmatter
\setcounter{page}{0}
\thispagestyle{empty}

\begin{center} {\bf
\Huge Eisenstein series and \\[0.5em] automorphic representations 
}

\vspace{8mm}

Philipp Fleig${}^{1}$, Henrik P. A. Gustafsson${}^2$, Axel Kleinschmidt${}^{3,4}$, Daniel Persson${}^{2}$

\footnotesize
\vspace{8mm}

${}^1${\it Institut des Hautes \'Etudes Scientifiques, IHES\\
Le Bois-Marie, 35, Route de Chartres, 91440 Bures-sur-Yvette, France\\[3mm]}
${}^2${\it Department of Physics, Chalmers University of Technology\\
 412 96 Gothenburg, Sweden\\[3mm]} 
${}^3${\it Max Planck Institute for Gravitational Physics, Albert Einstein Institute\\
Am M\"uhlenberg 1, 14476 Potsdam, Germany\\[3mm]}
${}^4${\it  International Solvay Institutes\\ Campus Plaine C.P. 231, Boulevard du
Triomphe, 1050 Bruxelles, Belgium}

\vspace{15mm}

\hrule

\vspace{5mm}

\parbox{130mm}{

\noindent \footnotesize 

We provide an introduction to the theory of Eisenstein series and automorphic forms on real simple Lie groups $G$, emphasising the role of representation theory. It is useful to take a slightly wider view and define all objects over the (rational) adeles $\mathbb{A}$, thereby also paving the way for connections to number theory, representation theory and the Langlands program. Most of the results we present are already scattered throughout the mathematics literature but our exposition collects them together and is driven by examples. Many interesting aspects of these functions are hidden in their Fourier coefficients with respect to unipotent subgroups and a large part of our focus is to explain and derive general theorems on these Fourier expansions. Specifically, we give complete proofs of the Langlands constant term formula for Eisenstein series on adelic groups $G(\mathbb{A})$ as well as the Casselman--Shalika formula for the $p$-adic spherical Whittaker function associated to unramified automorphic representations of $G(\mathbb{Q}_p)$. In addition, we explain how the classical theory of Hecke operators fits into the modern theory of automorphic representations of adelic groups, thereby providing a connection with some key elements in the Langlands program, such as the Langlands dual group ${}^LG$ and automorphic $L$-functions. Somewhat surprisingly, all these results have natural interpretations as encoding physical effects in string theory. We therefore also introduce some basic concepts of string theory, aimed toward mathematicians, emphasising the role of automorphic forms. In particular, we provide a detailed treatment of supersymmetry constraints on string amplitudes which enforce differential equations of the same type that are satisfied by automorphic forms. We further explain in detail how instanton effects in string theory are captured by the Fourier coefficients of Eisenstein series. Our treatise concludes with a detailed list of interesting open questions and pointers to additional topics which go beyond the scope of this book.}

\vspace{5mm}

\hrule

\end{center}

\newpage

\pagestyle{fancy}
\renewcommand{\headrulewidth}{0pt}
\fancyhead{}

\mbox{ }

\vspace{-5mm}

\begin{center}
{\Large \bf Note to the reader}

\end{center}

\noindent This book has grown out of our endeavour to understand the theory of automorphic representations and the structure of Fourier expansions of automorphic forms with a particular emphasis on adelic methods and Eisenstein series. Our intention is also to open a path of communication between mathematicians and physicists, in particular string theorists, interested in these topics. Most of the results presented herein exist already in the literature and we benefitted greatly from~\cite{LanglandsFE,Bump,MR2807433,MR2808915,ShahidiBook,Deitmar,Goldfeld,Terras1,Terras2}; our exposition differs, however, at places from the standard one. A few new results and examples are included as well; in particular, we provide many techniques for working out aspects of the Fourier expansion of Eisenstein series. We would be very grateful to learn of any omissions and mistakes that we have made unintentionally.

Sections that are more advanced or explore topics beyond the main focus of this treatise are marked with an asterisk. 

\vspace{20mm}

\begin{center}
{\Large \bf Acknowledgements}

\end{center}

\noindent Our understanding of the material presented here was greatly facilitated by numerous discussions with colleagues both from the mathematics community and from the string theory community. We are especially indebted to Guillaume Bossard, David Ginzburg, Stephen D. Miller, Hermann Nicolai, Bengt E.W. Nilsson and Boris Pioline for many clarifying and stimulating discussions over the past years. In addition, we gratefully acknowledge the exchanges with Olof Ahl\'en, Sergei Alexandrov, Marcus Berg, Benjamin Brubaker, Lisa Carbone, Martin Cederwall, Brian Conrad, Thibault Damour, Eric D'Hoker, Dennis Eriksson, Alex J. Feingold, Matthias Gaberdiel, Terry Gannon, Ori Ganor, Dmitri Gourevitch, Michael~B.~Green, Murat G\"unaydin, Stefan Hohenegger, Shamit Kachru, Ralf K\"ohl (n\'e Gr\"amlich), Kyu-Hwan Lee, Carlos~R.~Mafra, Gregory~W.~Moore, Jakob Palmkvist, Christoffer Petersson, Siddhartha Sahi, Per Salberger, Gordan Savin, Oliver Schlotterer, Philippe Spindel, Stefan Theisen, Pierre Vanhove, Roberto Volpato, Peter West and Martin Westerholt-Raum. We are also grateful to Olof Ahl\'en, Marcus Berg, Guillaume Bossard and Brian Conrad for valuable comments on the first version of this treatise.

\vspace{10mm}

\renewcommand{\headrulewidth}{.5pt}
\fancyhead[LO]{}
\fancyhead[RO]{}
\fancyhead[CO]{\textit{\titleString}}
\fancyhead[CE]{\textit{\nouppercase\leftmark}}
\fancyhead[LE]{}
\fancyhead[RE]{}

\tableofcontents

\newpage
\setcounter{page}{0}
\thispagestyle{empty}

\mainmatter

\pagestyle{fancy}
\renewcommand{\headrulewidth}{.5pt}
\fancyhead[LO]{}
\fancyhead[RO]{}
\fancyhead[CO]{\textit{\titleString}}
\fancyhead[CE]{\textit{\nouppercase\leftmark}}
\fancyhead[LE]{}
\fancyhead[RE]{}
\fancyfoot{}
\fancyfoot[CO]{\arabic{page}}
\fancyfoot[CE]{\arabic{page}}

\chapter{Motivation and background}
\label{ch:intro}

\vspace{0.5cm}

\epigraph{\emph{An efficient, but abstract, way to approach the subject of automorphic forms is by the introduction of adeles, \\
rather ungainly objects that nevertheless, once familiar, spare much unnecessary thought and many useless calculations.}}{--- Robert P. Langlands${}^*$}

\nonumberfootnote{${}^*$Representation theory - its rise and its role in number theory, Proceedings of the Gibbs symposium (1989)}

\vspace{0.75cm}

\noindent This text grew out of our endeavour to learn the adelic techniques used in the analysis of Eisenstein series in many mathematical works. Part of our motivation came from research problems in string theory were we faced the challenge of calculating certain Fourier coefficients of automorphic forms on exceptional Lie groups. The present text can be viewed as the culmination of the resulting journey through the world of automorphic forms. Even though none of the results that we present here are new, we felt that there might be an interest in such a survey since many of the original sources and textbooks do not make an easy first reading, especially for theoretical physicists like ourselves. Therefore we strove to be as pedagogical \emp{and} precise as possible but will sometimes sacrifice rigour or generality for conveying ideas and explicit examples. The reader is referred to the many sources quoted if he wishes more details on a particular point.

\newpage

\section{Automorphic forms and Eisenstein series}
\label{sec:AutoIntro}

\emp{Automorphic forms} are functions $f(g)$ on a Lie group $G$ that 
\begin{enumerate}
\item[(1)] are invariant under the action of a discrete subgroup $\Gamma\subset G$: $f(\gamma\cdot g) = f(g)$ for all $\gamma\in\Gamma$,
\item[(2)] satisfy eigenvalue differential equations under the action of the ring of $G$-invariant differential operators and
\item[(3)] have well-behaved growth conditions.
\end{enumerate}
A more explicit and refined form of these conditions will be given in chapter~\ref{ch:autforms} when we properly define automorphic forms; here we content ourselves with a qualitative description based on examples. We will mainly be interested in automorphic forms $f(g)$ that are invariant under the action of the maximal compact subgroup $K$ of $G$ when acting from the right: $f(gk)=f(g)$ for all $k\in K$; such forms are called \emp{$K$-spherical}. The automorphic forms are then functions on the coset space $G/K$.

The prime example of an automorphic form is obtained when considering $G=SL(2,\reals)$ and $\Gamma=SL(2,\ints)\subset SL(2,\reals)$. The maximal compact subgroup is $K=SO(2,\reals)$ and the coset space $G/K$ is a constant negative curvature space isomorphic to the Poincar\'e upper half plane $\UHP=\{z=x+iy \st x\in \reals\,\,\mathrm{and}\,\,y>0\}$. A function satisfying the three criteria above is then given by the non-holomorphic function
\begin{align}
\label{Eisenintro}
f_s(z) = \sum_{(c,d)\in\ints^2\atop (c,d)\neq (0,0)} \frac{y^s}{|cz+d|^{2s}}.
\end{align}
The sum converges absolutely for $\Re(s)>1$. The action of an element $\gamma\in SL(2,\ints)$ on $z\in\UHP$ is given by the standard fractional linear form (see section~\ref{sec:SL2})
\begin{align}
\label{SL2ACTintro}
\gamma\cdot z = \frac{az+b}{cz+d} \quad\mathrm{for} \quad 
\gamma=\begin{pmatrix}a&b\\c&d\end{pmatrix}\in SL(2,\ints).
\end{align}
Property (1) is then verified by noting that the integral lattice $(c,d)\in\ints^2$ is preserved by the action of $SL(2,\ints)$ and acting with $\gamma\in SL(2,\ints)$ in (\ref{Eisenintro}) merely reorders the terms in the absolutely convergent sum. Property (2) in this case reduces to a single equation since there is only a single primitive $SL(2, \reals)$-invariant differential operator on $SL(2,\reals)/SO(2, \reals)$. This operator is given by
\begin{align}
\label{DeltaUHP}
\Delta = y^2 \left( \partial_x^2 + \partial_y^2\right)
\end{align}
and corresponds to the Laplace--Beltrami operator on the upper half plane $\UHP$. In group theoretical terms it is the quadratic Casimir operator. Acting with it on the function (\ref{Eisenintro}) one finds 
\begin{align}
    \label{eq:SL2Eisenstein-eigenvalue}
\Delta f_s(z) = s(s-1) f_s(z)
\end{align}
and hence $f_s(z)$ is an eigenfunction of $\Delta$ (and therefore of the full ring of differential operators generated by $\Delta$). Condition (3) relating to the growth of the function here corresponds to the behavior of $f_s(z)$ near the boundary of the upper half plane, more particularly near the so-called \emp{cusp at infinity} when $y \to\infty$.\footnote{For $\Gamma=SL(2,\ints)$ this is the only cusp up to equivalence. With this one means that the fundamental domain of the action of $\Gamma$ on $\UHP$ only touches the boundary of the upper half plane at a single point. See section~\ref{sec:SL2} for pictures and~\cite{LangSL2,KoblitzElliptic,BorelSL2} more details on discrete subgroups of $SL(2,\reals)$.} 
The growth condition requires $f_s(y)$ to grow at most as a power law as $y\to\infty$. To verify this point it is easiest to consider the \emp{Fourier expansion} of $f_s(y)$. This requires a bit more machinery and also paves the way to the general theory. We will introduce it heuristically in section~\ref{intro:Fourier} and in detail in chapter~\ref{ch:fourier}.

The form of the function $f_s(z)$ is very specific to the action of   $SL(2,\ints)$ on the upper half plane $\UHP$. To prepare the ground for the more general theory of automorphic forms for higher rank Lie groups we shall now rewrite (\ref{Eisenintro}) in a more suggestive way. In fact, $f_s(z)$ is (almost) an example of an \emphindex{Eisenstein series} on $G=SL(2,\reals)$. To see this, we first extract the greatest common divisor of the coordinates of the lattice point $(c,d)\in\ints^2$:
\begin{align}
f_s(z) = \left(\sum_{k>0} k^{-2s}\right)\sum_{(c,d)\in\ints^2\atop \gcd(c,d)=1} \frac{y^s}{|cz +d|^{2s}}
= \zeta(2s) \sum_{(c,d)\in\ints^2\atop \gcd(c,d)=1} \frac{y^s}{|c z+d|^{2s}}
\end{align}
where we have evaluated the sum over the common divisor $k$ using the \emp{Riemann zeta function}~\cite{Riemann}
\begin{align}
\label{RZintro}
\zeta(s) = \sum_{n>0} n^{-s}.
\end{align}
Referring back to (\ref{SL2ACTintro}), we can rewrite the summand using an element of the group $SL(2,\ints)$:
\begin{align}
\frac{y^s}{|cz+d|^{2s}} = \lb\Im \left(\gamma\cdot z\right)\rb^s\quad\text{for}\quad 
\gamma=\begin{pmatrix}a&b\\c&d\end{pmatrix}.
\label{cosetSL2}
\end{align}
For this to be possible, two things have to occur: $(i)$ For any co-prime pair $(c,d)$ such a matrix $\gamma\in SL(2,\ints)$ must exist, and $(ii)$ if several matrices exist we must form equivalence classes such that the sum over co-prime pairs $(c,d)$ corresponds exactly to the sum over equivalence classes.
For $(i)$, we note that the condition that $c$ and $d$ be co-prime is necessary since it would otherwise be impossible to satisfy the determinant condition $ad-bc=1$ over $\ints$. At the same time, co-primality is sufficient to guarantee existence of integers $a_0$ and $b_0$ that complete $c$ and $d$ to a matrix $\gamma\in SL(2,\ints)$. In fact, there is a one-parameter family of solutions for $\gamma$ that can be written as
\begin{align}
\label{allSL2}
\begin{pmatrix} a_0+m c & b_0+md\\c &d\end{pmatrix}
=\begin{pmatrix} 1&m \\0 &1\end{pmatrix}
\begin{pmatrix} a_0 & b_0\\c &d\end{pmatrix}
\end{align}
for any integer $m\in\ints$. (That these are all solutions to the determinant condition over $\ints$ is an elementary lemma of number theory, sometimes called \emp{B\'ezout's lemma}~\cite{Bezout}.) The form (\ref{allSL2}) tells us also how to resolve point $(ii)$: We identify matrices that are obtained from each other by left multiplication by a matrix belonging to the \emp{Borel subgroup}
\begin{align}
\label{Borelintro}
B(\ints) = \left\{\begin{pmatrix}\pm 1&m\\0&\pm 1\end{pmatrix}\,\middle|\,m\in\ints\right\}\subset SL(2,\ints).
\end{align}
The interpretation of this group is that it is the stabiliser of the $x$-axis. 

Summarising the steps we have performed, we find that we can write the function (\ref{Eisenintro}) as
\begin{align}
\label{SL2lattice}
f_s(z) = 2\zeta(2s)\sum_{\gamma\in B(\ints)\setminus SL(2,\ints)} \lb \Im\left(\gamma\cdot z\right)\rb^s.
\end{align}
Since we had included the matrix $-\idm$ in the definition of $B(\ints)$, an extra factor of $2$ arises in this formula.

Dropping the multiplicative $\zeta$-factor, we obtain the function
\begin{align}
\label{Eisenintro2}
E(\chi_s,z) = \sum_{\gamma\in B(\ints)\setminus SL(2,\ints)} \chi_s(\gamma\cdot z),
\end{align}
where we have also introduced the notation $\chi_s(z)=\lb\Im(z)\rb^s=y^s$. The reason for this notation is that $\chi_s$ is actually induced from a \emphindex[character!on Borel subgroup]{character} on the real Borel subgroup. We will explain this in more detail below in chapter~\ref{ch:autforms}. Note that this way of writing the automorphic form makes the invariance under $SL(2,\ints)$ completely manifest because it is a sum over images.

The form (\ref{Eisenintro2}) is what we will call an \emphindex{Eisenstein series} on $SL(2,\reals)$ and it is this form that generalises straight-forwardly to Lie groups $G(\reals)$ other than $SL(2,\reals)$ (whereas the form with the sum over a lattice does not, as we discuss in more detail in section~\ref{sec:latticesums}). In complete analogy with (\ref{Eisenintro2}) we define the (minimal parabolic) Eisenstein series on $G(\reals)$ invariant under the discrete group $G(\ints)$ by\footnote{We will always take $G(\ints)$ as the Chevalley group that is defined as the stabiliser (in $G(\reals)$) of a preferred integral basis (Chevalley basis) of the Lie algebra of $G(\reals)$; see section \ref{sec:discgps} below more details.}
\begin{align}
\label{EIntro}
E(\chi,g) = \sum_{\gamma \in B(\ints)\setminus G(\ints)} \chi(\gamma g)
\end{align}
where $\chi$ is (induced from) a character on the Borel subgroup $B(\reals)$ and $g\in G(\reals)$. Eisenstein series are the protagonists of the story we will develop.

\section{Why Eisenstein series and automorphic forms?}

Before delving into the further analysis of Eisenstein series, let us briefly step back and provide some motivation for their study.

\subsection{A mathematician's possible answer} 
Automorphic forms are of great importance in many mathematical fields such as number theory, representation theory and algebraic geometry. The various ways in which automorphic forms enter these seemingly disparate fields are connected by a web of conjectures collectively referred to as the Langlands program \cite{LanglandsWeilLetter, LanglandsProb, GelbartLanglands,KnappLanglandsProgram, KnappFirstSteps,FrenkelLanglands}. 

Much of the arithmetic information is contained in the Fourier coefficients of automorphic forms. The standard examples correspond to modular forms on $G(\reals) = SL(2, \reals)$, where these coefficients  yield eigenvalues of Hecke operators (covered in chapter~\ref{ch:Hecke}) and the counting of points on elliptic curves. 

For arbitrary Lie groups $G(\reals)$ one considers the Hilbert space $L^2(\Gamma\backslash G(\reals))$ of square-integrable functions that are invariant under a left action by a discrete subgroup $\Gamma\subset G(\reals)$. This space carries a natural action of $g\in G(\reals)$, called the right-regular action, through 
\beq
\big[\pi(g)f\big] (x) = f(xg)
\eeq
where $f\in L^2\big(\Gamma\backslash G(\reals)\big)$, $g,x\in G(\reals)$ and $\pi: G(\reals)\to \mathrm{Aut}\big(L^2\big(\Gamma\backslash G(\reals)\big)\big)$ is the right-regular representation map. Since the functions are square-integrable the representation is unitary. This representation-theoretic viewpoint on automorphic forms was first proposed by Gelfand, Graev and Piatetski-Shapiro \cite{GGPS} later developed considerably by Jacquet and Langlands \cite{JL}. This perspective provides the key to generalising the classical theory of modular forms on the complex upper half plane to higher rank Lie groups. 

It is an immediate, important and difficult question as to what the decomposition of the space $L^2(\Gamma\backslash G(\reals))$ into irreducible representations of $G(\reals)$ looks like. The irreducible constituents in this decomposition are called \emph{automorphic representations}. This spectral problem was tackled and solved by Langlands~\cite{LanglandsFE}. The Eisenstein series (and their analytic continuations) form an integral part in the resolution although they themselves are not square-integrable.\footnote{A passing physicist might note that this is very similar to using non-normalisable plane waves as a `basis' for wave functions in quantum mechanics. Indeed the piece $\chi(\gamma g)$ in (\ref{EIntro}) is exactly like a plane wave; the $\gamma$-sum is there to make it invariant under the discrete group by the method of images so that $E(\chi,g)$ are the simplest $\Gamma$-invariant plane waves. The decomposition of an automorphic function (`wave-packet')\index{wave-packet} in this basis (extended by the discrete spectrum of cusp forms and residues) is the content of various trace formulas discussed in section~\ref{classautorep}.} Although we will not describe the full resolution of this problem in these notes, automorphic representations will play a prominent role in our discussion.

\subsection{A physicist's possible answer}
\label{sec_phys}

Many problems in quantum mechanics are characterised by discrete symmetries. At the heart of many of them lies Dirac quantisation where charges (e.g. electric or magnetic) of physical states are restricted to lie in certain lattices rather than in continuous spaces. The discrete symmetries preserving the lattice are often called dualities and can give very interesting different angles on a physical problem. This happens in particular in string theory, where such dualities mix perturbative and non-perturbative effects. 

For the discrete symmetry to be a true symmetry of a physical theory, all observable quantities must be given by functions that are invariant under the discrete symmetry, corresponding to property (1) discussed at the beginning of section~\ref{sec:AutoIntro}. Similarly, the dynamics or other symmetries of the theory impose differential equations on the observables, corresponding to property (2), and the growth condition (3) is typically associated with having well-defined perturbative regimes of the theory. The main example we have in mind here comes from string theory and the construction of scattering amplitudes of type II strings in various maximally supersymmetric backgrounds~\cite{Green:1997tv,Green:1997tn,Green:1998by,Kiritsis:1997em,Obers:1999um,Green:2010wi}. However, the logic is not necessarily restricted to this, see also~\cite{Takhtajan,Petropoulos:2012ne,Bogomolny:1992cj,Aurich:2004ke} for some other uses of automorphic forms in physics.

For these reasons one is naturally led to the study of automorphic forms in physical systems with discrete symmetries. Via this route one is also led to the same spectral problem posed by the mathematician since one needs to determine which automorphic representation a given physical observable belongs to. Again, the Eisenstein series and their properties are the building blocks of such spaces and it is important to understand them well. Furthermore, in a number of examples from string theory it was actually possible to show that the observable is given by an Eisenstein series itself~\cite{Green:1997tv,Green:2010kv}.\footnote{That Eisenstein series are mostly not square-integrable is no problem in these cases since the object computed (part of a scattering amplitude) is not a wavefunction and not required to be normalisable.}

\section{Analysing automorphic forms and adelisation}
\label{intro:Fourier}

We now return to the study of Eisenstein series defined by (\ref{EIntro}) and their properties, starting again with the very explicit example (\ref{Eisenintro}) for $SL(2,\reals)$.

\subsection{Fourier expansion of the \texorpdfstring{$SL(2,\reals)$}{SL(2, R)} series}

The discrete Borel subgroup $B(\ints)$ of (\ref{Borelintro}) acts on the variable $z=x+i y$ as translations by
\begin{align}
\begin{pmatrix}\pm 1&m\\0&\pm 1\end{pmatrix} \cdot z =z\pm m\quad\textrm{for $m\in\ints$}
\end{align}
and therefore any automorphic function (that is by definition invariant under any discrete transformation) is periodic in the $x$ direction with period equal to $1$ corresponding to the smallest non-trivial $m=1$. This means that we can Fourier expand it in modes $e^{2\pi i n x}$. Applying this to (\ref{Eisenintro2}) leads to
\begin{align}
\label{SL2Fourierintro}
E(\chi_s,z)  = \underbrace{\vphantom{\sum_{n\neq 0}} \, C(y) \, }_{\textrm{constant term}\atop\textrm{zero mode}} + \underbrace{\sum_{n\neq 0} a_n (y) e^{2\pi i n x}}_{\textrm{non-zero mode}}.
\end{align}
As we indicated, it is natural to divide the Fourier expansion into two parts depending on whether one deals with the zero Fourier mode (a.k.a. constant term) or with a non-zero mode. Since the Fourier expansion was only in the $x$ direction, the Fourier coefficients still depend on the second variable $y$.\footnote{If one dealt with an automorphic form holomorphic in $z$ (called modular forms in chapter~\ref{ch:autforms} below) this would not be true since the holomorphicity condition links the $x$ and $y$ dependence. The Fourier coefficients in an expansion in $q=e^{2\pi i (x+iy)}=e^{2\pi i z}$ would be pure numbers. This is the origin of the name \emphindex{constant term} for the zero mode in (\ref{SL2Fourierintro}).}

Determining the explicit form of the Fourier coefficients is one of the key problems in the study of Eisenstein series. In the example of $SL(2,\reals)$ this can for instance be done by making recourse to the formulation in terms of a lattice sum that was given in (\ref{Eisenintro}) and using the technique of Poisson resummation. The calculation is reviewed in appendix~\ref{app:SL2Fourier} and leads to the following explicit expression
\begin{align}
\label{SL2FC2}
E(\chi_s,z) = y^s + \frac{\xi(2s-1)}{\xi(2s)} y^{1-s} + \frac{2y^{1/2} }{\xi(2s)} \sum_{m\neq 0} |m|^{s-1/2} \sigma_{1-2s}(m) K_{s-1/2} (2\pi |m| y) e^{2\pi i m x}, 
\end{align}
where 
\begin{align}
\label{CRZ}
\xi(s)=\pi^{-s/2} \Gamma(s/2) \zeta(s)
\end{align}
is the \emphindex[Riemann zeta function!completed]{completion of the Riemann zeta function} \eqref{RZintro}, $K_s(z)$ is the modified Bessel function of the second kind (that decreases exponentially for $z\to\infty$ in accordance with the growth condition) and 
\begin{align}
\sigma_{1-2s}(n) =\sum_{d|n} d^{1-2s}
\label{intro:divisorsum}
\end{align}
is called a divisor sum (or the instanton measure in physics; see chapter~\ref{ch:intro-strings} below) and given by a sum over the positive divisors of $n\neq 0$.

As is evident from (\ref{SL2FC2}), the explicit form of the Fourier expansion can appear quite complicated and involves special functions as well as number theoretic objects. For the case of more general groups $G(\reals)$ the method of Poisson resummation is not necessarily available as there is not always a form of the Eisenstein series as a lattice sum. \emp{It is therefore desirable to develop alternative techniques for obtaining (parts of) the Fourier expansion under more general assumptions.}\footnote{Additional care has to be taken for the Fourier expansion for general $G(\reals)$ also because the translation group $B(\ints)$ is in general not abelian. One can still define (abelian) Fourier coefficients as we will see, however, they fail to capture the full Eisenstein series. There are also non-abelian parts to the Fourier expansion.} This is achieved by lifting the Eisenstein series into an adelic context which we now sketch and explain in more detail in section~\ref{sec_defauto}.

\subsection{Adelisation of Eisenstein series}
\label{sec:adelisation-of-Eisenstein-series}

A standard elementary technique in number theory for analysing equations over $\ints$ is by analysing them instead as congruences for every prime (and its powers) separately (sometimes known as the Hasse principle or the local-global principle based on the Chinese remainder theorem)~\cite{Apostol1,Neukirch} (see also Appendix \ref{sec_Hasse} for some examples). One way of writing all the terms together is by using the \emp{ring of adeles} $\mathbb{A}$. The adeles can formally be thought of as an infinite tuple
\begin{align}
\label{adelesintro}
a = (a_\infty; a_2,a_3,a_5,a_7,\ldots) \in \mathbb{A} = \reals \times \prodp_{p<\infty} \rats_p,
\end{align}
where  $\rats_p$ denotes the \emp{$p$-adic numbers} that are a completion of the rational number $\rats$ that is inequivalent to the standard one (leading to $\reals$) and that is parametrised by a prime number $p$ and defined properly in section~\ref{sec:padic}. The product is over all prime numbers and the prime on the product symbol indicates that the entries $a_p$ in the tuple are restricted in a certain way (see~(\ref{adeleseq}) below for the exact statement). The real numbers $\reals$ can be written as $\rats_\infty$ in this context and interpreted as the completion of $\rats$ at the `prime' $p=\infty$. \emp{Very} crudely, an adele can be thought of as summarising the information of an object modulo all primes.

\emp{Strong approximation} is a similar method that lifts a general automorphic form from being defined on the space $G(\ints)\backslash G(\reals)/K(\reals)$ to the space $G(\rats)\backslash G(\mathbb{A})/K(\mathbb{A})$ so that $G(\rats)$ plays the role of the discrete subgroup that was played by $G(\ints)$ before. However, $G(\rats)$ is a nicer group than $G(\ints)$ since $\rats$ is a field whereas $\ints$ is only a ring. This facilitates the analysis and allows the application of many theorems for algebraic groups. 

A consequence of using strong approximation and adeles is that the result of the calculation factorises according to (\ref{adelesintro}) and one can do the calculation for all primes and $p=\infty$ separately. Indeed, the explicit form (\ref{SL2FC2}) for the Fourier expansion of the $SL(2,\reals)$ Eisenstein series already secretly had this form. This can be seen for example in the constant term since
\begin{align}
\label{Constintro}
\frac{\xi(2s-1)}{\xi(2s)} = \pi^{1/2}\frac{\Gamma(s-1/2)}{\Gamma(s)} \prod_{p<\infty} \frac{1-p^{-2s}}{1-p^{1-2s}}
\end{align}
where we have used the definition of the completed zeta function from~\eqref{CRZ} and the \emp{Euler product formula} for the Riemann zeta function~\cite{Riemann}
\begin{align}
\label{RiemannProd}
\zeta(s) = \sum_{n>0} n^{-s} =\prod_{p<\infty} \frac{1}{1-p^{-s}}.
\end{align}
In (\ref{Constintro}) we clearly recognise a factorised form that is very similar to (\ref{adelesintro}). That this is not an accident will be demonstrated in section~\ref{sec:Sl2const} for $SL(2,\reals)$. For the other Fourier modes in \eqref{SL2FC2} we get a similar factorisation with the modified Bessel function belonging to the $p=\infty$ factor and
\begin{equation}
    \sigma_{1-2s}(m) = \prod_{p < \infty} \gamma_p(m) \frac{1 - p^{-(2s-1)} \abs{m}_p^{2s-1}}{1 - p^{-(2s-1)}}
\end{equation}
where $\abs{m}_p$ is the $p$-adic norm of $m$ defined in section~\ref{sec:padic} and $\gamma_p(m)$ selects all factors with $\abs{m}_p \leq 1$ as shown in section~\ref{sec:p-adic-special-functions}. The complete derivation for the non-constant terms can be found in section~\ref{sec:SL2FC} for the $SL(2,\reals)$ Eisenstein series.

The adelic methods are so powerful that one can obtain a closed, simple and group-theoretic formula for the constant term of Eisenstein series on any (split real) Lie group $G(\reals)$. This formula, known as the \emphindex{Langlands constant term formula}, will be the topic of chapter~\ref{ch:CTF}.

For the (abelian) Fourier coefficients, the adelic methods also help to obtain fairly general results, in particular for the part that involves the finite primes $p<\infty$. For the contribution coming from the $\reals$ in (\ref{adelesintro}) the results are not quite as general; already for $SL(2,\reals)$ this is what gives the modified Bessel function. We discuss the Fourier coefficients in chapter~\ref{ch:Whittaker-Eisenstein}.

\section{Reader's guide and main theorems}

The following is a brief outline of the contents of this treatise. Chapters~\ref{ch:intro-strings},~\ref{ch:p-adic-and-adelic} and~\ref{ch:Lie-groups} are introductory and provide some background material for subsequent chapters. Chapter~\ref{ch:intro-strings} gives a general overview of how automorphic forms enter in computing scattering amplitudes in string theory. This is not intended as a comprehensive introduction to string theory, but its aim is rather to act as a first glimpse, primarily directed towards mathematicians, of a vast and fascinating topic that is closely tied to automorphic forms and representation theory. Throughout the main text we also offer remarks and pointers that indicate physical interpretations of various mathematical notions and results. Chapter~\ref{ch:p-adic-and-adelic}, on the other hand, introduces the basic machinery of $p$-adic and adelic analysis which will be crucial for everything we do later. The main thrust of the chapter is provided by the numerous examples of computing $p$-adic integrals that will be used extensively in proving Langlands constant term formula, and computing Fourier coefficients of Einstein series. In chapter~\ref{ch:Lie-groups} we introduce some basic features of Lie algebras and Lie groups that will be used throughout the book. We first discuss Lie groups and Lie algebras over $\mathbb{R}$ and then move on to algebraic groups over $\mathbb{Q}_p$ as well as adelic groups. 

We shall now discuss the structure of the remainder of the text in a little more detail, with emphasis on the central results in each chapter. 
\begin{itemize}
\item In chapters~\ref{ch:autforms} and~\ref{ch:autoreps} we introduce the general theory of automorphic forms and automorphic representations. We start out gently by discussing how to pass from modular forms on the upper half plane to automorphic forms on the adelic group $SL(2,\mathbb{A})$. We then move on to the general case of arbitrary Lie groups. We define \emp{Eisenstein series} $E(\chi,g)$ for general split real Lie groups $G(\reals)$ that are invariant under the discrete Chevalley subgroup $G(\ints)$. The definition (\ref{EIntro}) requires the choice of a character $\chi$ of a parabolic subgroup $P$ of $G$; alternatively, we can think of $\chi$ as being defined by a choice of weight vector $\lambda$ of the (split real) Lie algebra of $G(\reals)$. We explain how this can be understood from the point of view of the representation theory of $G(\mathbb{R})$, and we show how to lift the function from being defined on $G(\reals)$ to a function defined on $G(\mathbb{A})$ where $\mathbb{A}$ are the \emp{adeles} of the rational number field $\rats$. 
\item A major part of this book is devoted to analysing Fourier expansions of automorphic forms. This is a highly non-trivial subject with many interesting connections to representation theory as well as to physics. In chapter~\ref{ch:fourier} we introduce the general theory of Fourier coefficients and Whittaker functions, with emphasis on the representation-theoretic viewpoint. Toward the end of the chapter we also introduce some more advanced topics, such as nilpotent orbits and wavefront sets, as well as the Piatetski-Shapiro--Shalika formula. 

\item In chapter~\ref{ch:SL2-fourier} we illustrate all the general techniques in the context of Eisenstein series on $SL(2)$. Specifically, using adelic techniques, we provide a detailed proof of the  following classic theorem:
\begin{theorem} The complete Fourier expansion of the Eisenstein series $E(\chi_s,g)$ for $g\in SL(2,\reals)\subset SL(2,\ads)$ is given by:
\begin{align}
E(\chi_s,g)= y^{s}+\frac{\xi(2s-1)}{\xi(2s)}y^{1-s}
+\sum_{m\neq 0}\frac{2 y^{1/2}}{\xi(2s)}  |m|^{s-1/2}\sigma_{1-2s}(m)K_{s-1/2}(2\pi|m|y) e^{2\pi i m x}.
\label{SL2exp_outline}
\end{align}

Furthermore, the Eisenstein series satisfies the functional relation\index{functional relation!for SL(2,R)@for $SL(2,\reals)$}
\begin{align}
E(\chi_s,g) = \frac{\xi(2s-1)}{\xi(2s)} E(\chi_{1-s},g),
\end{align}
\end{theorem}
where $\xi$ is the (completed) Riemann zeta function. The rest of the notation  is explained in chapter~\ref{ch:SL2-fourier}. 

\item The first two terms in the Fourier expansion above correspond to the zeroth Fourier coefficients. These are often collectively referred to as the \emph{constant term} of the Eisenstein series, A very important and general result in this context is provided by the so-called Langlands constant term formula, which yields a remarkably simple expression for the complete constant term of Eisenstein series on arbitrary semi-simple Lie groups. In chapter~\ref{ch:CTF} we give a complete proof of the following theorem of Langlands:
\begin{theorem}[Langlands constant term formula] Let $G(\rats)$ be a split semi-simple algebraic group with $G(\reals)$ a real semi-simple Lie group and $G(\mathbb{A})$ the corresponding adelic group. Let $\lambda$ be a weight of the Lie algebra $\mathfrak{g}$, $\Weyl$ the associated Weyl group, and $N$ a maximal unipotent radical of $G$. We then have
\begin{align}
\lint_{N(\mathbb{Z})\backslash N(\mathbb{R})} E(\lambda,ng) dn = \sum_{w\in \Weyl} a^{w\lambda+\rho} \prod_{\alpha>0\,:\, w\alpha<0} \frac{\xi(\langle \lambda| \alpha\rangle)}{\xi(1+\langle \lambda| \alpha\rangle)},
\label{Langlandsconst_outline}
\end{align}
where $a$ belongs to the Cartan torus $A\subset G$,  and the product runs over positive roots $\alpha$ of $\mathfrak{g}$.
\end{theorem}

\item The infinite sum in the Fourier expansion (\ref{SL2exp_outline}) corresponds to the non-zero coefficients and this is generally referred to as the \emph{non-constant term}. In chapter~\ref{ch:Whittaker-Eisenstein} we discuss the general structure of Fourier coefficients of Eisenstein series on reductive groups $G$. For this part of the expansion much less is known explicitly. However, there exists a beautiful formula due to Kato--Shintani--Casselman--Shalika, commonly known as the \emph{Casselman--Shalika formula}, which gives an explicit expression for the so-called $p$-adic \emphindex[Whittaker function!$p$-adic]{Whittaker function}. This corresponds to a local version of the Fourier coefficient of the Eisenstein series, which can be used to reassemble the full (global) coefficient. A large part of chapter~\ref{ch:Whittaker-Eisenstein} is therefore devoted to proving the following theorem:
\begin{theorem}[Casselman--Shalika formula]
Let $G(\mathbb{Q}_p)$ $(p<\infty)$ be a $p$-adic semi-simple Lie group, $N(\mathbb{Q}_p)$ a maximal unipotent radical of $G(\mathbb{Q}_p)$ and $\psi$ an unramified unitary character on $N(\mathbb{Q}_p)$. The Casselman--Shalika formula is given by:
\begin{align}
\int_{N(\mathbb{Q}_p)} \chi(w_0 na)\overline{\psi(n)}dn=\frac{\eps(\lambda)}{\xi(\lambda)} \sum_{w\in W} (\det(w)) |a^{w\lambda+\rho}|  \prod_{\alpha>0  \atop w\alpha<0} p^{\langle\lambda|\alpha\rangle}
\label{CSformula_outline}
\end{align}
\end{theorem}

\item For certain special types of Fourier coefficients, so-called \emph{degenerate Whittaker functions}, one can take one step further and compute the full global coefficient (and not just the $p$-adic version). In chapter~\ref{ch:Whittaker-Eisenstein} we also prove the following theorem which gives such a formula:
\begin{theorem}
Let $\psi: N(\rats)\backslash N(\ads) \to U(1)$ be a degenerate character with associated subgroup $G'(\mathbb{A})\subset G(\mathbb{A})$ on which it is supported. Then the degenerate Whittaker function on $G(\mathbb{A})$ is given by
\begin{align}
W^\circ_\psi (\chi,a) = \sum_{w_c\wlong'\in \Weyl/\Weyl'} a^{(w_c\wlong')^{-1}\lambda+\rho}  M(w_c^{-1},\lambda) W'^\circ_{\psi^a}(w_c^{-1}\lambda,\id),
\end{align}
where $W'^\circ_{\psi}$ denotes a Whittaker function on the $G'(\mathbb{A})$ subgroup of $G(\mathbb{A})$. The weight $w_c^{-1}\lambda$ is given as a weight of $G'(\mathbb{A})$ by orthogonal projection. 
\end{theorem} 
A more complete formulation is provided in section~\ref{sec:degpsi}. In section~\ref{sec:SL3ex} we also provide an extensive example of how to calculate Whittaker functions for Eisenstein series on $SL(3,\mathbb{A})$.

\item  In chapter~\ref{ch:working} we illustrate how to perform calculations with Eisenstein series in practice. More specifically we explain how to evaluate Langlands constant term formula in concrete  examples, which, in particular, involves a detailed analysis of the functional equation.  We also show how to perform similar evaluations of the Whittaker functions that appear in the non-constant Fourier coefficients. We provide some explicit examples for exceptional Lie groups. 
    
\item It is interesting to note that both the Langlands constant term formula (\ref{Langlandsconst_outline}) and the Casselman--Shalika formula (\ref{CSformula_outline}) have cunning similarities to the Weyl character formula. That this is not an accident is a central insight of Langlands. To understand this requires the additional machinery of \emph{Hecke theory}, which is the topic of chapter~\ref{ch:Hecke}. Here we explain how to pass from the classical notion of Hecke operators acting on modular forms to the general notion of \emph{spherical Hecke algebras} on adelic groups. This analysis leads us to a reformulation of the Casselman--Shalika formula that clearly illustrates the intimate connection with the Weyl character formula. In this context we are naturally lead to the notion of the \emph{Langlands dual group} ${}^L G$ and to the notion of \emph{automorphic $L$-functions}, which form a central ingredient in the Langlands program. The chapter concludes with a discussion of the \emph{Langlands--Shahidi method}, which is a powerful way to construct $L$-functions from automorphic representations.

\item In the concluding chapter~\ref{ch:outlook}, we present various interesting directions which we have not been able to cover in detail. We emphasise open questions and conjectures, many of which have sprung out of problems in string theory, and we have tried to formalise and phrase them in purely mathematical terms. We also briefly discuss some of the key ingredients in the Langlands program and we make various comments and conjectures regarding its extension to Kac--Moody groups. 

\end{itemize}

We end by giving some disclaimers: All groups $G(\reals)$ that will be considered here are associated with split real forms and we also restrict to simple groups. The only base field that we will use for adelisation are the rational numbers $\rats$. Often we will perform formal manipulations of infinite sums and integrals without paying attention to whether the expressions are (absolutely) convergent or not. The expressions typically depend on a set of parameters and for some range of parameters convergence can be established. In many cases, the results can be extended by analytic continuation.

\chapter{String theory scattering and automorphic forms}
\label{ch:intro-strings}

In this chapter we will introduce some basic concepts of string theory with emphasis on scattering amplitudes of closed strings. Our main purpose is to illustrate the deep connection between string theory and automorphic forms, which may come as a surprise to a mathematician. The topic is far too vast for us to do it justice in just a single chapter, but hopefully this will provide a sufficient glimpse to spark the motivation for further studies. The main point we wish to convey  is that the Fourier coefficients of  automorphic forms on higher rank Lie groups capture essential information about string theory scattering amplitudes. We should stress that scattering amplitudes are but one aspect of the relation between string theory and automorphic forms. Other connections are discussed in the outlook chapter, such as to BPS-state counting (section~\ref{sec:BH-counting}), theta correspondences (section~\ref{sec:ThetaCorr}) and moonshine (section~\ref{sec:Moonshine}). For more information about string theory we recommend the books \cite{GSW, Zwiebach:2004tj, Polchinski, BLT} and, for a brief introduction, the lecture notes \cite{DHoker-strings,Tong:2009np}.

\section{String theory concepts}
String theory is a theory of one-dimensional extended objects propagating in a (Lorentzian) \emphindex[target space]{target space-time} $M$.
During their propagation, strings sweep out a two-dimensional \emphindex[world-sheet]{world-sheet} $\Sigma$ and string theory can therefore be thought of as the dynamics of the embedding maps $X:\Sigma\to M$, where both $\Sigma$ and $M$ are endowed with additional structure (like a metric) that enter in the definition of the dynamics.
In \emphindex{superstring theory} ---that we exclusively consider here and refer to generally as string theory--- this additional structure includes \emphindex[supersymmetry!on string world-sheet]{world-sheet supersymmetry} and the space of allowed world-sheets $\Sigma$ is then the space of all closed, orientable super Riemann surfaces. 
We here also specialise to so-called type IIB superstrings; otherwise one might have to include boundaries and non-orientable surfaces as well.
Riemann surfaces are classified by their \emphindex[Riemann surface!genus]{genus} $h\in\ints_{\geq 0}$, a topological invariant.

A fundamental parameter of string theory is the characteristic \emph{length} $\ells$ \emph{of a string}\index{string length}.
More commonly, one uses the parameter $\alpha'=\ells^2$ which can be thought of as the scale of area of the string world-sheet $\Sigma$.
For very small string lengths, $\ells\to 0$, the strings look effectively like point particles.
The (bosonic) spectrum of string excitations in flat ten-dimensional Minkowski space $M=\reals^{1,9}$ is then given by states of mass
\begin{align}
\label{eq:stringspectrum}
m^2 = \frac{2}{\alpha'} N,
\end{align}
where $N\in\ints_{\geq 0}$ is the excitation number of the string states.
The lightest string states are massless and there is an infinite sequence of massive states with quantised masses and the separation between the masses is set by the~\emphindex[string theory!string scale]{string scale $\alpha'$}.
At every mass level one has a finite number of degrees of freedom.

\begin{figure}[t!]
\centering
\includegraphics[scale=.6]{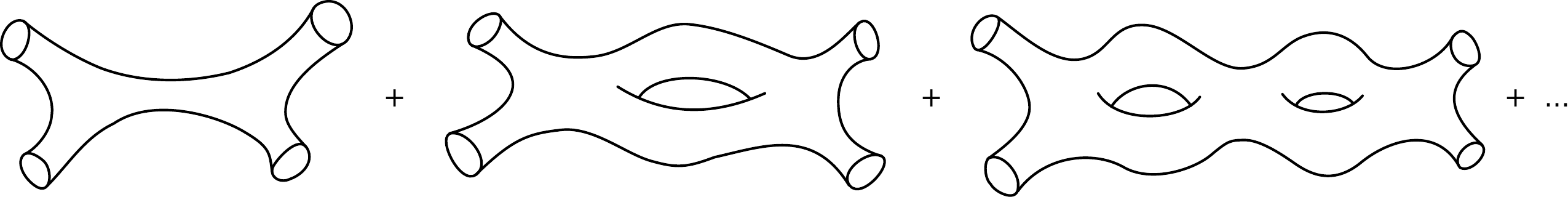}
\caption{\label{fig:LowGenus}\sl String world-sheets as they appear in the scattering of four closed strings. Ignoring the asymptotic boundary states, the diagrams correspond to genera $h=0,1,2$, respectively.}
\end{figure}

Strings can interact by various joining and splitting processes.
Such interactions are studied by string world-sheets with asymptotic states coming together for a scattering process and then separating again.
Examples with few splittings, corresponding to low genus world-sheets, are depicted in figure~\ref{fig:LowGenus}.
In the limit $\alpha'\to 0$ the diagrams lose their waists and reduce to standard~\emphindex[Feynman diagram]{Feynman diagrams} commonly used in the~\emphindex[quantum field theory]{quantum field theory} of~\emphindex[point particles]{point particles}.
At the same time, the massive states ($N>0$) in~\eqref{eq:stringspectrum} become infinitely heavy compared to the massless states ($N=0$) in the limit $\alpha'\to 0$. 

Among the quantities one wants to compute in string theory are~\emphindex[scattering amplitude]{scattering amplitudes}. They provide information about the likelihood of a certain scattering process of strings to take place.
Computing a scattering amplitude in \emphindex[string theory!perturbation theory]{string perturbation theory} requires summing over possible world-sheets of all genera $h$. 
A given topology is then weighted by the so-called \emphindex[string coupling]{string coupling} $\gs$ with a weight $\gs^{2(h-1)}$.
The string coupling is a measure of the strength of string-string interaction, i.e., the joining and splitting of strings. Note that in string theory it is conventional to use the term~\emphindex[string theory!perturbation theory!loop]{loop} when referring to the genus of a world-sheet. The probability of a certain string scattering process is given by the modulus square of the scattering amplitude. The perturbative expansion is a power expansion in $\gs$.

Besides the data of the asymptotic states, the scattering amplitude depends on $\alpha'$, the string coupling and potentially other so-called \emphindex[moduli!in string theory]{moduli fields}. These can be thought of as aspects of the target space $M$ in the form of additional (scalar) fields living on them. Only their vacuum expectation values matter for the present discussion and we will denote these by 
\begin{align}
\label{stringmoduli}
g \in \mathcal{M}\quad\quad\textrm{(moduli expectation values).}
\end{align}
Here, $\mathcal{M}$ is the so-called \emphindex[moduli!space]{moduli space} of string theory. The string coupling constant $\gs$ turns out to be related to one of the moduli fields called the \emphindex[string theory!dilaton]{dilaton}, but in general there are many more moduli fields.

\begin{figure}[t!]
\vspace{2em}
\centering
\input{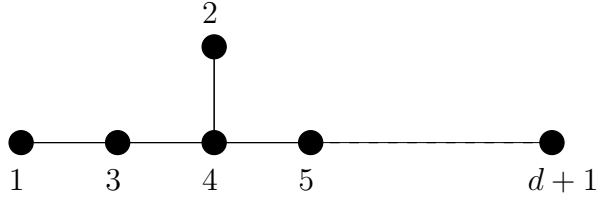}
\caption{\label{fig:CJ}\sl The Dynkin diagram of the Cremmer--Julia symmetry  group $E_{d+1}$ with labelling of nodes in the `Bourbaki convention'.}
\end{figure}

The structure of moduli space is of central importance in understanding the possible forms of string scattering amplitudes. Much is known for flat target spaces of the type $M=\reals^{1,9}$ (\emphindex[Minkowski space]{flat Minkowski space}\index{flat space}) or $M=\reals^{1,9-d}\times T^d$ (\emphindex{toroidal compactification}). In both cases one retains \emphindex[supersymmetry]{maximal supersymmetry} strongly constraining the moduli space. The classical low energy moduli space is a symmetric space of the form
\begin{align}
\label{moduliclass}
\mathcal{M}_{\mathrm{class.}} = G(\reals)/K(\reals) = E_{d+1}(\reals) / K(E_{d+1}(\reals)),
\end{align}
where $E_{d+1}$ is the Cremmer--Julia sequence of duality groups~\cite{Cremmer:1978ds,Julia:1980gr,Cremmer:1997ct} that are listed in table~\ref{tab:CJ} and their Dynkin diagrams shown in figure~\ref{fig:CJ}, and $K(E_{d+1}(\reals))$ are their maximal compact subgroups. Up to $d\leq 7$, these groups are finite-dimensional reductive groups and we restrict to this range first. We will come back to $d\geq 8$ in section~\ref{sec:KM}. For other internal manifolds one can get a large variety of different Lie groups.

\begin{table}[t!]
\centering
\caption{\label{tab:CJ}\sl Table of Cremmer--Julia symmetry groups $G(\reals)$ with compact subgroup $K(\reals)$ and U-duality groups $G(\ints)$ for compactifications of type IIB string theory on a $d$-dimensional torus $T^d$ to $D = 10 - d$ dimensions.}
\begin{tabular}{ccccc}
\toprule
$d$ & $G(\reals)$ & $K(\reals)$ & $G(\ints)$ & $D$\\
\midrule
$0$ & $SL(2,\reals)$ & $SO(2,\reals)$ & $SL(2,\ints)$ & $10$\\
$1$ & $GL(2,\reals)$ & $SO(2,\reals)$ & $SL(2,\ints)$ & $9$\\
$2$ & $SL(2,\reals)\times SL(3,\reals)$ & $SO(2,\reals)\times SO(2,\reals)$ & $SL(3,\ints)\times SL(2,\ints)$ & $8$\\
$3$ & $SL(5,\reals)$ & $SO(5,\reals)$ & $SL(5)$ & $7$\\
$4$ & $SO(5,5,\reals)$ & $(SO(5,\reals)\times SO(5,\reals))/\ints_2$ & $SO(5,5,\ints)$ & $6$\\
$5$ & $E_6(\reals)$ & $USp(8,\reals)/\ints_2$ & $E_6(\ints)$ & $5$\\
$6$ & $E_7(\reals)$ & $SU(8,\reals)/\ints_2$ & $E_7(\ints)$ & $4$\\
$7$ & $E_8(\reals)$ & $SO(16,\reals)/\ints_2$ & $E_8(\ints)$ & $3$\\
\bottomrule
\end{tabular}
\end{table}

The classical low-energy effective theory in $D = 10 - d$ dimensions, which is described by supergravity, has a symmetry given by the non-compact real Lie group ${G}(\mathbb{R})$ \cite{Hull:1994ys,Obers:1998fb}.
However, as mentioned in section \ref{sec_phys}, when passing to the quantum theory, the classical symmetries are generically broken because the (generalised) electro-magnetic charges of physical states become quantised according to the 
\emphindex[Dirac--Schwinger--Zwanziger quantisation]{Dirac--Schwinger--Zwanziger quantisation condition}, and take values in some integral lattice $\Gamma$~\cite{Dirac:1931kp,Teitelboim:1985ya,Teitelboim:1985yc,Nepomechie:1984wu}.
Although the classical symmetry group $G(\reals)$ is broken, there is a discrete subgroup of $G(\reals)$ that survives and remains a symmetry of the full quantum theory.
This quantum symmetry is defined as the subgroup of $G(\reals)$ that preserves the lattice $\Gamma$~\cite{Hull:1994ys}
\begin{equation}
    \{ g\in G(\reals) \st g\Gamma=\Gamma\} .
\end{equation}
This quantum symmetry group is generally referred to as a \emph{U-duality group} which unifies the previously known existing dualities called S- and T-duality and agrees with the Chevalley subgroup $G(\ints)$~\cite{Soule}.
The U-duality groups $G(\ints)$ for toroidal compactifications of type IIB string theory on a torus $T^d$ are also listed in table~\ref{tab:CJ}.

Points in moduli space related by U-duality transformations give rise to equivalent string theories. This implies that the correct moduli space of quantum string theory is not the classical symmetric space~\eqref{moduliclass} but 
\begin{align}
\label{Mexact}
\mathcal{M} \equiv  \mathcal{M}_{\mathrm{quantum}} =  G(\ints) \backslash G(\reals)/ K(\reals) \, ,
\end{align}
and all observables, like string scattering amplitudes, that are functions of the expectation values of the moduli are functions on this space.

Put differently, physical observables are $G(\ints)$-invariant functions on $G(\reals)/K$. In addition, physical constraints, such as supersymmetry, typically force these observables to satisfy differential equations and have a prescribed asymptotic behaviour at infinity, thereby satisfying the conditions (1)--(3) in section \ref{sec:AutoIntro}, characterising automorphic forms. We will explain the origin of these additional constraints in more detail in section~\ref{sec:SUSYcons} but first we study a particularly relevant example.

\section{Four-graviton scattering amplitudes}
\label{sec:intro-four-graviton}

In ten-dimensional type IIB string theory the scattering of four massless string states called \emphindex{gravitons} gives rise to a quantum correction to the standard gravitational interactions in general relativity computed without the inclusion of string theory degrees of freedom. 

The exact form of the full~\emphindex[string theory!perturbation theory!four-graviton scattering]{string theory four-graviton scattering} amplitude is not known but one can attempt to perform series expansions of the amplitude with respect to some of its arguments. There are two common expansions of the amplitude. 

The first expansion is \emphindex[string perturbation theory|see{string theory}]{string perturbation theory}\index{string theory!perturbation theory}, discussed above, in which one treats the string coupling constant $\gs$ as small and computes the contributions to the amplitudes from Riemann surfaces of increasing genus.
This involves an integral over the moduli space of all Riemann surfaces of a given genus and with a number of punctures corresponding to the number of asymptotic scattering states.
These integrals have been studied up to two loops, see for instance~\cite{DHoker:2002gw,DHoker:2005jc,DHoker:2005ht}, and become increasingly hard for increasing genus.
One complication arises from the \emphindex{Schottky problem} of parametrising the moduli space for large genus $h$.
Another serious complication is that one should actually integrate over the moduli space of \emph{super} Riemann surfaces since one is dealing with superstring theory and it is known that this integral cannot be reduced in a simple manner to an integral over ordinary Riemann surfaces for $h\geq 5$~\cite{Donagi:2013dua}.
Finally, the amplitude is not expected to be a convergent series in $\gs$, meaning that there are \emphindex[string theory!non-perturbative effects]{non-perturbative effects} arising from instanton configurations~\cite{Shenker:1990uf,Green:1997tv,Pioline:2009ia}.
These are roughly of the form $e^{-1/\gs}$ and do not admit a Taylor series expansion around weak coupling $\gs=0$ and therefore cannot be captured by string perturbation theory.
The string coupling $\gs$ is one of the coordinates on the moduli space $\mathcal{M}$ and the limit $\gs\to 0$ corresponds to approaching a cusp on $\mathcal{M}$. 

The second expansion of the amplitude is the~\emphindex[string theory!low energy expansion]{low energy expansion} in which one considers the momenta of the scattering particles to be small, leading to an expansion in derivatives of the fields. 
Dimensionless expansion parameters are formed by multiplying momenta squared with $\alpha'$, which is why the expansion is also called the~\emphindex[string theory!$\alpha'$-expansion]{$\alpha'$-expansion}. It is this expansion that makes contact to automorphic forms and we will now study it in detail.

As asymptotic states, gravitons are characterised by their momenta $k_i\in \reals^{1,9}$ ($i=1,2,3,4$) and their polarisations $\epsilon_i\in S^2(\reals^{1,9})$, which are symmetric second rank tensors subject to some constraints whose detail we do not require for the present discussion. Since gravitons are massless, the momenta satisfy $k_i^2=0$, where the norm-squared is computed using the Lorentzian metric on $\reals^{1,9}$. Out of the four momenta $k_i$ one forms the dimensionless Lorentz invariant \emphindex{Mandelstam variables} 
\begin{align}
\label{Mandeldef}
s=-\frac{\alpha'}{4} (k_1+k_2)^2,\quad t=-\frac{\alpha'}{4} (k_1+k_3)^2 \quad\textrm{and}\quad  u=-\frac{\alpha'}{4} (k_1+k_4)^2.
\end{align}
(Identical constructions also apply in the toroidally compactified theory.) Momentum conservation ($k_1+k_2+k_3+k_4=0$) and masslessness imply that $s+t+u=0$. Any symmetric polynomial in $s,t,u$ is then a polynomial in
\begin{align}
\label{symMandel}
\sigma_2 = s^2+t^2+u^2 \quad \textrm{and}\quad \sigma_3=s^3+t^3+u^3.
\end{align}
The string scattering amplitude will therefore be a function of the momenta only through $\sigma_2$ and $\sigma_3$. Similar simplifications arise for the polarisation tensors $\epsilon_i$ that enter the final answer only in a particular combination that we will denote by $\mathcal{R}^4$. It can be expressed as the contraction of two copies of a standard rank 8 tensor $t_8$, whose precise form is for example given in~\cite{Gross:1986iv, GSW}, and four copies of the linearised curvature tensor $\mathcal{R}_{\mu\nu\rho\sigma} \propto k_\mu \epsilon_{\nu\rho} k_\sigma$ (with permutations), see section~\ref{sec:four-graviton-tree-amplitude} for some more details.

Our four graviton amplitude in $D=10-d$ dimensions is therefore of the form
\begin{align}
\mathcal{A}^{(D)} (s,t,u,\epsilon_i; g)\,,
\end{align}
with $g\in\mathcal M$. We recall that the string scale $\alpha'$ was absorbed into the Mandelstam variables $s,t,u$. 

The (analytic part of the) $\alpha'$-expansion of the four-graviton amplitude (in Einstein frame)  takes the form~\cite{Green:2010wi}
\begin{align}
\label{fourgrav}
\mathcal{A}^{(D)}(s,t,u,\epsilon_i;g) = 
\left[\mathcal{E}^{(D)}_{(0,-1)}(g)\frac{1}{\sigma_3}+ 
\sum_{p\geq 0} \sum_{q \ge 0} \mathcal{E}^{\scriptstyle{(D)}}_{\scriptstyle{(p,q)}}(g) \sigma_2^p \sigma_3^q
\right] \mathcal{R}^4 \, .
\end{align}
The interesting objects in this expression are the coefficient functions $\mathcal{E}^{\ordd{D}}_{\gra{p}{q}}(g)$ that are functions on the quantum moduli space $\mathcal{M_\text{quantum}} = G(\ints)\bs G(\reals)/K(\reals)$ in~\eqref{Mexact}.

The first term in~\eqref{fourgrav} plays a special role in that it is the only term that is not polynomial in $\sigma_2$ and $\sigma_3$. It is the lowest order term in the $\alpha'$-expansion and it agrees with what one would calculate in a standard theory of gravity with Lagrangian given by the Ricci scalar only (referred to as the~\emphindex{Einstein--Hilbert term}). The coefficient function $\mathcal{E}^{\ordd{D}}_{\gra{0}{-1}}(g)=3$ is constant. By contrast, the infinite series of terms in $p$ and $q$ come with higher powers of $\alpha'$ and they reflect the contribution of \emph{massive} string states to the graviton scattering process~\cite{Gross:1986iv}.

The \emphindex{low energy effective theory}\index{effective action} is obtained by writing the field theory action whose classical interactions give rise to the same quantum corrected amplitudes obtained from string theory, order by order in $\alpha'$. It gets corrections from the four-graviton amplitudes of the form \cite{Green:2011vz} 
\begin{equation}
    \mathcal{E}^{(D)}_{(p,q)}(g)  D^{2p+3q} R^4 
\end{equation}
where $D$ denotes a covariant derivative and $R^4$ (not to be confused with the Ricci curvature scalar) is a contraction of two $t_8$ tensors and four Riemann curvature tensors like the linearised version for the polarisation term $\mathcal{R}^4$ in the amplitude.

In other words, for the first few orders in $\alpha'$ we get the corrections
\begin{equation}
    \label{eq:four-graviton-effective-action}
    S = S_\text{class.} + \int \! d^Dx \sqrt{-G} \Big( (\alpha')^3 \mathcal{E}^{(D)}_{(0,0)}(g) R^4 + (\alpha')^5 \mathcal{E}^{(D)}_{(1,0)}(g) D^4 R^4 + (\alpha')^6 \mathcal{E}^{(D)}_{(0,1)}(g) D^6 R^4 + \ldots \Big) \, ,
\end{equation}
where $S_\text{class.}$ is the classical, zeroth order low energy effective action described by supergravity. We see why this expansion is also called the \emphindex[string theory!derivative expansion]{derivative expansion}.

As the functions $\mathcal{E}^{\ordd{D}}_{\gra{p}{q}}(g)$ depend on the moduli $g\in \mathcal{M}$, they in particular depend on the string dilaton and thus on the string coupling $\gs$ that controls string perturbation theory in terms of Riemann surfaces as discussed above. However, there is no reason that the dependence on $\gs$ be analytic. Non-analytic terms in $\gs$ are known as \emphindex{non-perturbative effects} and they appear in $\mathcal{E}^{\ordd{D}}_{\gra{p}{q}}(g)$ through so-called \emphindex{instanton contributions}. Their direct determination in terms of a string theory calculation is typically very hard. The action of the U-duality group $G(\ints)$ also includes a transformation that mixes perturbative and non-perturbative effects and therefore using U-duality opens up the opportunity to access non-perturbative effects indirectly.

The fact that the coefficients $\mathcal{E}^{\ordd{D}}_{\gra{p}{q}}$ are functions on the (quantum) moduli space $\mathcal{M}$ can be rephrased as follows: \emph{They are functions on the symmetric space $\mathcal{M}_\textnormal{class.} = G(\reals)/K(\reals)$ invariant under the left action of the discrete (U-duality) group $G(\ints)$}.

This means that $\mathcal{E}^{\ordd{D}}_{\gra{p}{q}}$ satisfy~(1) of the definition of automorphic forms in section~\ref{sec:AutoIntro} and, since they should also have perturbative expansions in the weak string coupling limit $\gs \to 0$ and other similar limits of the moduli space corresponding to cusps in $G/K$; they also satisfy the growth condition~(3). 

The remaining condition (2) from the introduction requires that an automorphic form satisfies appropriate differential equations under the action of $G$-invariant differential operators. It turns out that  supersymmetry in string theory imposes precisely such differential conditions on the coefficient functions. This was analysed first by Green and Sethi in the case of ten-dimensional ($D=10$) type IIB string theory and $p=q=0$~\cite{Green:1998by} and we review aspects of their arguments in section~\ref{sec:SUSYcons}. More recent work on supersymmetry constraints includes~\cite{Green:2010wi,Basu:2011he,Bossard:2014lra,Bossard:2014aea,Bossard:2015uga,Wang:2015aua}. They found that $\mathcal{E}^{\ordd{10}}_{\gra{0}{0}}(g)$ has to satisfy a Laplace equation with an eigenvalue determined by supersymmetry considerations. This, together with the known value for the string tree level ($h=0$ topology) scattering result, uniquely determined the coefficient function to be a non-holomorphic Eisenstein series on $SL(2,\reals)$~\cite{Pioline:1998mn} as discussed by Green and Gutperle in~\cite{Green:1997tv,Green:1997tn}.

In other dimensions $D$ and for small values of $p$ and $q$, there are strong arguments that the coefficient functions $\mathcal{E}^{\ordd{D}}_{\gra{p}{q}}(g)$ satisfy the differential equations~\cite{Green:2010wi,Pioline:2015yea}
\begin{subequations}
    \label{eq:intro-fourgravLap}
    \begin{align}
        \label{eq:intro-LapR4}
        R^4:& &\left(\Delta_{G/K} - \frac{3(11-D)(D-8)}{D-2} \right) \mathcal{E}^{\scriptstyle{(D)}}_{\scriptstyle{(0,0)}}(g) &= 6\pi\delta_{D,8},&\\
        \label{eq:intro-LapD4R4}
        D^4R^4:& &\left(\Delta_{G/K} - \frac{5(12-D)(D-7)}{D-2} \right) \mathcal{E}^{\scriptstyle{(D)}}_{\scriptstyle{(1,0)}}(g) &= 40\zeta(2)\delta_{D,7}+7\mathcal{E}_{\scriptstyle{(0,0)}}^{\scriptstyle{(6)}} \delta_{D,6},&\\
        \label{eq:intro-LapD6R4}
        D^6R^4:&\quad&\left(\Delta_{G/K} - \frac{6(14-D)(D-6)}{D-2} \right) \mathcal{E}^{\scriptstyle{(D)}}_{\scriptstyle{(0,1)}}(g) &= -\left(\mathcal{E}_{\scriptstyle{(0,0)}}^{\scriptstyle{(D)}}\right)^2+ 40\zeta(3) \delta_{D,6}&\\
      &&  &\quad +\frac{55}3\mathcal{E}_{\scriptstyle{(0,0)}}^{\scriptstyle{(5)}} \delta_{D,5} + \frac{85}{2\pi} \mathcal{E}_{\scriptstyle{(1,0)}}^{\scriptstyle{(4)}} \delta_{D,4},&\nn
    \end{align}
\end{subequations}
where $\Delta_{G/K}$ is the \emphindex{Laplace--Beltrami operator} on $G/K$.

We see that the third equation is qualitatively very different from the first two since it has the square of a non-constant function as a source on the right-hand side. This would take us out of the standard domain of automorphic forms and we will discuss this case in more detail in section~\ref{sec:D6R4}.

The Kronecker delta contributions in all three equations in~\eqref{eq:intro-fourgravLap} are related to the existence of ultraviolet and infrared divergences in the underlying supergravity theory and the existence of supersymmetric counterterms.
The ultraviolet divergences arise in those dimensions where also the eigenvalue vanishes and signal logarithmic terms in the coefficient function $\mathcal{E}^{\ordd{D}}_{\gra{p}{q}}$. We refer the reader to~\cite{Green:2010sp} for further discussions of this point. The additional Kronecker delta contributions proportional to coefficient functions associated with fewer derivatives are related to form factor divergences and these are discussed in~\cite{Pioline:2015yea,Bossard:2015oxa}.

Besides these special cases, equations~\eqref{eq:intro-LapR4} and~\eqref{eq:intro-LapD4R4} correspond to eigenfunction conditions from (2) in section~\ref{sec:AutoIntro}. For dimensions lower than ten, there are additional $G$-invariant differential operators other than $\Delta_{G/K}$ but the corresponding conditions are not fully known from string theory. A superspace analysis that generates the other differential equations was pioneered in~\cite{Bossard:2014lra,Bossard:2014aea}. It is expected that the coefficient functions for $R^4$ and $D^4R^4$ satisfy all the required differential equations and hence are standard automorphic forms.

As seen in section~\ref{sec:AutoIntro} in the case of $G(\reals) = SL(2, \reals)$, Eisenstein series are eigenfunctions of the Laplace--Beltrami operator, and comparing with computed scattering amplitudes in string theory one has been able to conjecture the exact forms of the coefficients $\mathcal{E}^{\ordd{D}}_{\gra{0}{0}}$ and $\mathcal{E}^{\ordd{D}}_{\gra{1}{0}}$ in terms of maximal parabolic Eisenstein series which will be defined in chapter~\ref{ch:autforms}. 
Parabolic subgroups, denoted by $P$, are introduced in section~\ref{sec:parsubgp}. 

More precisely, in five, four and three dimensions, with symmetry groups $E_6$, $E_7$ and $E_8$ according to table~\ref{tab:CJ}, if one considers the maximal parabolic subgroups $P$ that have semi-simple Levi parts $SO(5,5)$, $SO(6,6)$ and $SO(7,7)$, respectively, then the solutions
\begin{subequations}
    \label{eq:intro-R4D4R4}
    \begin{align}
        \label{eq:intro-R4fn}
        R^4: & & \mathcal{E}^{(D)}_{(0,0)}(g) &= 2\zeta(3) E(\lambda_{s = 3/2}, P, g),&\\
        \label{eq:intro-D4R4fn}
        D^4 R^4:& & \mathcal{E}^{(D)}_{(1,0)}(g) &= \zeta(5)E(\lambda_{s = 5/2}, P, g).&
    \end{align}
\end{subequations}
to equations~\eqref{eq:intro-LapR4} and~\eqref{eq:intro-LapD4R4} are the conjectured coefficient functions appearing in the four-graviton amplitudes, or equivalently, as corrections to the effective action \eqref{eq:four-graviton-effective-action}.

The weight $\lambda_s$ specifies the character $\chi_s$ on $P$, which defines the Eisenstein series similar to \eqref{EIntro}. It is given by
\begin{align}
    \label{eq:intro-lambda_s}
\lambda_s = 2s \Lambda_P -\rho,
\end{align}
where $\Lambda_P$ denotes the fundamental weight orthogonal to the Levi subgroup $L$ of $P=LU$ and $\rho$ the Weyl vector.

\begin{remark}
    Looking at the coefficients of \eqref{eq:intro-R4D4R4} we recognise the corresponding values from the tree level amplitudes which are computed in section~\ref{sec:four-graviton-tree-amplitude}. This means that the above functions are nothing but the (single) U-duality orbit of the tree level results. This is no longer true for the higher functions $\mathcal{E}_{\gra{p}{q}}^{\ordd{D}}$ as we will discuss in more detail in section~\ref{sec:D6R4}.
\end{remark}

These conjectures have been subjected to numerous consistency checks~\cite{Green:2010kv,Pioline:2010kb,Green:2011vz} and, particularly, capture the known results of scattering amplitudes in the weak coupling limit $\gs \to 0$ which we will discuss in the following section. 

\begin{remark}
We would also like to point out that recent investigations of superstring scattering amplitudes at tree-level and one-loop for more than four particles have revealed very interesting different connections to number theory. Instead of single $\zeta$-values like $\zeta(3)$ one will typically have so-called (elliptic) \emphindex{multiple zeta values} governed by \emphindex[Dringeld associator]{Drinfeld associators}~\cite{Schlotterer:2012ny,Broedel:2013aza,Broedel:2014vla,EnriquezEA}. We note that this structure is at fixed order in string perturbation theory whereas the U-duality invariant functions we are discussing here include all perturbative and non-perturbative effects.
\end{remark}

\section{Physical interpretation of the Fourier expansion}
\label{sec:physical-interp-10D}

We will now study the functions $\mathcal{E}^{\ordd{D}}_{\gra{p}{q}}$ which were found above as the quantum corrections to the low energy effective action in type IIB string theory on tori $T^d$ with $d=10-D$. Since $\mathcal{E}^{\ordd{D}}_{\gra{p}{q}}$ are invariant under the discrete subgroup $G(\ints)$, they are periodic functions and we can extract physical information from their Fourier expansions.

For concreteness, let us consider the $R^4$ and $D^4 R^4$ coefficients $\mathcal{E}^{\ordd{10}}_\gra{0}{0}$ and $\mathcal{E}^{\ordd{10}}_{\gra{1}{0}}$ in ten dimensions, where $G(\reals) = SL(2, \reals)$ --- although the physical interpretations hold for general dimensions and coefficients. 

As stated in section~\ref{sec:AutoIntro}, the classical moduli space $G(\reals)/K(\reals) = SL(2, \reals) / SO(2, \reals)$ is isomorphic to the Poincar\'e upper half plane $\UHP = \{z = x + i y \in \cmplx \st y = \Im z > 0\}$ parametrised by a complex scalar field called the \emphindex{axio-dilaton} here denoted by $z = \chi + i e^{-\phi}$ with $\chi$ being the axion and $\phi$ the dilaton. When $\phi$ is constant we have that $\gs = e^{\phi} = y^{-1}$. More generally, the string coupling constant is the asymptotic value at infinity denoted by $e^{\phi_0}$. The U-duality group $G(\ints) = SL(2, \ints)$ acts on $z$ by \eqref{SL2ACTintro}, including the translation $z \to z + 1$.

To reproduce the right perturbative behaviour, the leading order term in the weak coupling limit $\gs \to 0$, i.e. $y\to \infty$, should be (in Einstein frame): 
\begin{equation}
    \label{eq:SL2-weak-coupling-limit}
    \begin{split}
        \mathcal{E}^{(10)}_{(0,0)} &\sim 2\zeta(3)y^{3/2} \\
        \mathcal{E}^{(10)}_{(1,0)} &\sim \zeta(5)y^{5/2}
    \end{split}
    \qquad \qquad \text{as}\quad y\to \infty.
\end{equation}
These weak coupling limits correspond to the tree-level contribution to the scattering amplitudes and are computed in section~\ref{sec:four-graviton-tree-amplitude}.

The eigenvalue equations \eqref{eq:intro-fourgravLap} for the coefficient functions are
\begin{equation}
    \label{eq:10D-eigeneq}
    \begin{split}
        \Delta \mathcal{E}^{(10)}_{(0,0)}(z) &= \frac{3}{4} \mathcal{E}^{(10)}_{(0,0)}(z) \\
        \Delta \mathcal{E}^{(10)}_{(1,0)}(z) &= \frac{15}{4} \mathcal{E}^{(10)}_{(1,0)}(z)
    \end{split}
\end{equation}
where $\Delta = y^2(\partial_x^2 + \partial_y^2)$ is the Laplace--Beltrami operator on $\UHP$ from \eqref{DeltaUHP}. The coefficient functions are thus automorphic forms as defined in section~\ref{sec:AutoIntro}.

It was first realised by Green et al. in \cite{Green:1997tv} and \cite{Green:1999pu}, that these conditions are solved by
\begin{equation}
\label{eq:sols10}
    \begin{split}
        \mathcal{E}^{(10)}_{(0,0)}(z) &= f_{3/2}(z) = 2 \zeta(3) E(s = 3/2, z) \\
        \mathcal{E}^{(10)}_{(1,0)}(z) &= \frac{1}{2} f_{5/2}(z) = \zeta(5) E(s = 5/2, z)
    \end{split}
\end{equation}
as seen from \eqref{SL2FC2} and \eqref{eq:SL2Eisenstein-eigenvalue} with $f_s(z)$ defined in \eqref{Eisenintro}. This is exactly the $SL(2, \reals)$ variant of \eqref{eq:intro-R4D4R4}. 

\begin{remark}
A priori, the solutions~\eqref{eq:sols10} are only unique up to the addition of cusp forms but they were subsequently ruled out in \cite{Pioline:1998mn} for the $R^4$ coupling.
\end{remark}

The Fourier expansion \eqref{SL2FC2} then has a direct physical interpretation: the first two terms (constant terms) correspond to the perturbative quantum corrections (tree-level and one-loop), while the infinite series of the remaining Fourier modes encode non-perturbative effects. To see this we can expand the Bessel function $K_{s-1/2} (2\pi |n| y)$ in the limit $y \to \infty$ which, for the $R^4$ coefficient, yields 
\begin{equation}
    \label{eq:intro-SL2-fourier-expansion}
    \mathcal{E}^{(D)}_{(0,0)}(z) = 
    \overbrace{ 
        \underbrace{\vphantom{\Bigg(}2\zeta(3)y^{3/2}}_{\substack{\text{tree-level} \\[0.5cm] \includegraphics[height=1.0cm]{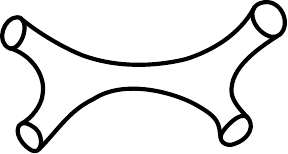}}} + 
        \underbrace{\vphantom{\Bigg(}4\zeta(2)y^{-1/2}}_{\substack{\text{one-loop} \\[0.5cm] \includegraphics[height=1.0cm]{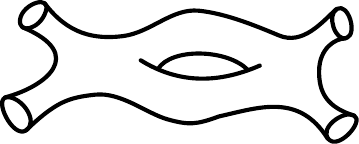}}}
    }^\text{perturbative terms} +
    \overbrace{\vphantom{\Bigg(}2\pi \sum_{m\neq 0} \sqrt{|m|} \sigma_{-2}(m)e^{-S_\text{inst}(z)}\left[1+\mathcal{O}(y^{-1})\right]}^\text{non-perturbative terms}_{\substack{\\[1.3cm]\text{amplitudes in the presence of instantons}}},
\end{equation}
where we have defined the \emphindex[instanton!action]{instanton action}
\begin{equation}
    \label{eq:Instaction}
    S_\text{inst}(z):=2\pi \abs{m} y - 2\pi i mx. 
\end{equation}
It is clear from this expression that the infinite series is exponentially suppressed by $e^{-y}$ in the limit $y = \gs^{-1} \to \infty$ corresponding to the weak coupling limit of the theory and this exponential suppression is characteristic for non-perturbative, instanton effects.

In string theory, these corrections come from so-called \emphindex[D-instantons]{\textnormal{D}-instantons}~\cite{Green:1995my, Green:1997tv} where D stands for Dirichlet as in the Dirichlet boundary conditions that are imposed on strings attached to them. These are special cases of extended objects in string theory called \emphindex[D-branes]{\textnormal{D}-branes}, localised to a single point in space-time and will be discussed in more detail in section \ref{sec:instantons}.

In this setting the divisor sum $\sigma_{-2}(m)=\sum_{d|m} d^{-2}$ is called the \emphindex[instanton!measure]{instanton measure} which weights the contributions of different modes $m$ called \emphindex[instanton!charge]{instanton charges}. As described in \cite{Green:1997tv} and shown in section \ref{sec:instantons}, the summation over divisors counts the degeneracy of D-instanton states of charge $m$. When $m$ is a negative integer it corresponds to an \emph{anti-instanton}.

The higher order corrections in $\gs = y^{-1}$ in the non-perturbative terms in \eqref{eq:intro-SL2-fourier-expansion} are higher genus corrections to the scattering amplitude in the presence of instantons.

In summary, the discrete symmetry leading us to the study of automorphic forms gives a lot of information about string theory: it tells us that there are no perturbative corrections to the $R^4$ term from world-sheets of genus larger than one, and gives an indirect way of computing scattering amplitudes in the presence of instantons in string theory, which is otherwise very difficult to do.

\vspace{5mm}
The above example provided the first hint of an intriguing relation between quantum corrections in string theory and automorphic forms. There is by now a vast literature on this subject; see \cite{Bao:2007er,Bao:2007fx,Bao:2009fg,Bao:2010cc,Basu:2006cs,Basu:2007ru,Green:1997di,Green:2010wi,GreenDualHighDeriv,Gubay:2010nd,Gunaydin:2005mx,Kazhdan:2001nx,Lambert:2006ny,Lambert:2010pj,Obers:1999um,ObersESeries,Pioline:1998mn,Pioline:2001jn,Pioline:2003bk,Pioline:2004xq,Pioline:2005vi,Pioline:2009qt} for a sample. In recent years also the representation theoretic aspects have proven to play an increasingly important role \cite{PerssonAuto,Green:2010kv,Green:2011vz,Pioline:2010kb,FK2012,FK2012Summary,FKP2013}, thus providing ample motivation also from physics for the emphasis on automorphic representations in these notes. 

Section~\ref{sec:outlook-strings}, continues the discussions of this chapter with the topic of automorphic representations for the coefficients $\mathcal{E}^{\ordd{D}}_{\gra{p}{q}}$ and the $D^6R^4$ correction which differs from the lower order terms by also requiring a non-constant source term in the differential equation \eqref{eq:intro-LapD6R4}.

\section{Supersymmetry constraints*}
\label{sec:SUSYcons}

We will now explain the arguments for obtaining the eigenvalue equation \eqref{eq:intro-LapR4} for the coefficient function $\mathcal{E}^{\ordd{10}}_{\gra{0}{0}}$ in ten dimensions using constraints from supersymmetry following \cite{Green:1998by, Green:1997me, Green:1999qt, Sinha:2002zr}. This requires us to include, not only the four-graviton corrections \eqref{eq:four-graviton-effective-action} to the effective action, but all interactions at the same order in $\alpha'$ related by supersymmetry. We denote the different corrections to the action as
\begin{equation}
    \label{eq:SUSY-constraints:full-action}
    S = S^{(0)} + (\alpha')^3 S^{(3)} + (\alpha')^4 S^{(4)} + \ldots
\end{equation}
where the four-graviton term $R^4$ of \eqref{eq:four-graviton-effective-action} is included in $S^{(3)}$ and where we have used that there are no corrections of order $\alpha'$ and $(\alpha')^2$. In fact, the term $S^{(4)}$ also vanishes but this does not affect our discussion.

All the different interaction terms at order $(\alpha')^3$ can be obtained by studying the so-called rigid (or global) limit of the linearised theory \cite{Green:1997me}. This means that we are, at first, only considering \emphindex[supersymmetry!global]{global supersymmetry} transformations for which the transformation parameters are constant spinors $\epsilon$ as opposed to the local case where they are space-time dependent. Moreover, the transformations are of linear order in the field fluctuations around some fixed background. From the linearised theory we will also find relations between the coefficient functions of different interactions. The sought for eigenvalue equation can then be found using the full non-linear symmetry.

The physical fields of the ten-dimensional type IIB supergravity theory are conveniently packaged into a generating function as the coefficients of an expansion in formal, complex Grassmann variables $\theta^\alpha$, $\alpha = 1, \ldots, 16$ together transforming as a Weyl spinor under $Spin(9,1)$. Their complexity means that they also transform with non-trivial weight under the $U(1)=SO(2)$ \emphindex{R-symmetry} of ten-dimensional type IIB symmetry. This $SO(2)$ symmetry is the denominator of the moduli space $SL(2,\reals)/SO(2)$. Together with the space-time coordinates, the variables $\theta^\alpha$ parametrise a \emphindex{superspace} and the generating function is called a \emphindex{superfield}. This construction only makes sense here at the linearised level. Supersymmetry transformations can then be described as coordinate transformations on this space. The corresponding differential operators acting on the superfields are
\begin{equation}
    \label{eq:superspace-operators}
    Q_\alpha = \frac{\partial}{\partial \theta^\alpha}, \qquad
    Q_\alpha^* = -\frac{\partial}{\partial \theta^{*\,\alpha}} + 2 i (\overline \theta \gamma^\mu)_\alpha \frac{\partial}{\partial x^\mu} \\ 
\end{equation}
which (by virtue of acting from the left and the right, respectively) anticommute with the covariant derivatives
\begin{equation}
    D_\alpha = \frac{\partial}{\partial \theta^\alpha} + 2 i (\gamma^\mu \theta^*)_\alpha \frac{\partial}{\partial x^\mu}, \qquad
    D_\alpha^* = -\frac{\partial}{\partial \theta^{*\,\alpha}},
\end{equation}
where $\overline \theta = \theta^\dagger \gamma^0$ is the Dirac conjugate, with $\theta^\dagger$ being the Hermitian conjugate and $\theta^*$ the complex conjugate. We use the Einstein summation convention where repeated indices are summed over. The $\gamma$-matrices are real matrices acting on spinors and satisfy the \emphindex{Clifford algebra} $\{\gamma^\mu, \gamma^\nu\} = 2 \eta^{\mu\nu}$ with the constant Minkowski background metric $\eta^{\mu\nu}=\textrm{diag}(-++\ldots+)$ where $\mu, \nu = 0, \ldots, 9$ are Lorentz indices for the ten directions of space-time. We also use the standard notation that a $\gamma$-matrix with multiple space-time indices $\gamma^{\mu_1 \cdots \mu_n}$ is the anti-symmetrised product of $n$ $\gamma$-matrices with the normalisation $\gamma^{\mu_1 \cdots \mu_n} = \gamma^{\mu_1} \cdots \gamma^{\mu_n}$ for $\mu_1 \neq \ldots \neq \mu_n$. 

Comparing the coefficients in the superfield expansion before and after an infinitesimal coordinate transformation in superspace gives the supersymmetry transformations for the physical fields. We see from \eqref{eq:superspace-operators} (acting on, for example, \eqref{eq:SUSY-constraints:superfield} below) that the supersymmetry transformation relates bosons (with an even number of free spinor indices) and fermions (with an odd number of spinor indices).

The superfield $\Phi$ of ten-dimensional type IIB supergravity satisfies the holomorphic constraint $D^*\Phi = 0$ and the on-shell condition $D^4 \Phi = D^{*4}\Phi^*$ \cite{Green:1997me,Howe:1983sra} which narrows down the expansion of $\Phi$ into its component fields --- the physical fields --- into the following \cite{Green:1997me, Green:1999qt}
\begin{equation}
    \label{eq:SUSY-constraints:superfield}
    \begin{split}
        \Phi &= z_0  + \delta\Phi \\
        &= z_0 + \delta z + \frac{1}{\gs}(i \overline \theta^* \delta\lambda + \delta \hat{G}_{\mu\nu\rho} \overline \theta^* \gamma^{\mu\nu\rho} \theta + \cdots + \mathcal {R}_{\mu\sigma\nu\tau} \overline \theta^* \gamma^{\mu\nu\rho} \theta \overline \theta^* \tensor{\gamma}{^{\sigma\tau}_{\rho}} \theta +
        \cdots + \theta^8 \partial^4 \overline {\delta z}) \, .
    \end{split}
\end{equation}
where $\delta\Phi$ is the linearised fluctuation around a flat background with constant axio-dilaton $z_0 = \chi_0 + i e^{-\phi_0} = \chi_0 + i \gs^{-1}$ and where $\delta \lambda$ is the fluctuation of the dilatino, $\mathcal{R}_{\mu\sigma\nu\tau}$ is the linearised curvature tensor of the metric fluctuation, and $\hat G$ is a supercovariant combination of fermion bilinears and field strengths of Ramond and Neveu--Schwarz two-form potentials. The exact expression for $\hat G$ can be found in Appendix C of \cite{Green:1999qt}. We will continue to use the notation $z = x + iy = z_0 + \delta z$ for the full axio-dilaton.

\begin{remark}
    If we let the metric fluctuation be a plane wave with polarisation $\epsilon_{\mu\nu}$ and momentum $k_\rho$ in the linearised curvature we recover the expression $\mathcal{R}_{\mu\nu\rho\sigma} \propto k_\mu \epsilon_{\nu\rho} k_\sigma$ (with permutations) from section~\ref{sec:intro-four-graviton} together with an exponential that later gives the momentum conservation in the amplitude.
\end{remark}

The linearised action at order $(\alpha')^3$ can be obtained as an integral over half of superspace (that is, no integrations over $\theta^*$) of the form \cite{Green:1997me, Green:1998by}
\begin{equation}
    \label{eq:SUSY-constraints:S3}
    S^{(3)}_\text{linear} = \Re \int d^{10}x \, d^{16}\theta \, F[\Phi]
\end{equation}
where the Grassmann integration amounts to taking the $\theta^{16}$ term of the expansion of $F[\Phi]$. 

From the expansion~\eqref{eq:SUSY-constraints:superfield}, we see that the $\theta^4$ term of $\Phi$ contains the linearised curvature $\mathcal{R}_{\mu\sigma\nu\tau}$ and the $\theta^{16}$ term of $S^{(3)}$ will therefore contain a contraction of four such curvature tensors. In fact, they form the $\mathcal{R}^4$ interaction of section \ref{sec:intro-four-graviton} which can be expressed as the superspace integral \cite{Nilsson:1986, Green:1999qt}
\begin{equation}
    \mathcal{R}^4 = \int d^{16}\theta \, (
    \mathcal R_{\mu\sigma\nu\tau} \overline \theta^* \gamma^{\mu\nu\rho} \theta \overline \theta^* \tensor{\gamma}{^{\sigma\tau}_{\rho}} \theta 
    )^4 \, .
\end{equation}
In the expansion, this interaction comes with a coefficient $\gs^{-4} \partial_{z}^4 F$ ignoring numerical prefactors. Inspecting~\eqref{eq:SUSY-constraints:superfield}, we see that another term in the linearised action $S^{(3)}_\text{linear}$ is the $\delta\lambda^{16}$ interaction with coefficient $\gs^{-16} \partial_{z}^{16} F$ where
\begin{equation}
    (\lambda^n)_{\alpha_{n+1} \cdots \alpha_{16}} = \frac{1}{n!} \epsilon_{\alpha_1 \cdots \alpha_{16}} \lambda^{\alpha_1} \cdots \lambda^{\alpha_n}
\end{equation}
with $\alpha_i$ being spinor indices. Similar terms can be constructed from the other pieces of the expansion~\eqref{eq:SUSY-constraints:superfield}.

The full (no longer linearised) supersymmetric action $S^{(3)}$ can, in this manner, be obtained as \cite{Green:1998by}
\begin{multline}
    \label{eq:SUSY-constraints:full-S3}
    S^{(3)} = \int d^{10}x \sqrt{g} \Big( f^{(12,-12)}(z) \lambda^{16} + f^{(11,-11)}(z) \hat G \lambda^{14} + \ldots + f^{(0,0)}(z) R^4 + \ldots \\ + f^{(-12,12)}(z) \lambda^{*\,16} \Big)
\end{multline}
where the coefficients $f^{(w, -w)}(z)$ are to be determined shortly. Note that we have renamed the $R^4$ coefficient which was labeled as $\mathcal{E}^{\ordd{10}}_{\gra{0}{0}}$ in section \ref{sec:physical-interp-10D}, and that, unlike the indices for $\mathcal{E}^{\ordd{D}}_{\gra{p}{q}}$, the indices $(w, -w)$ denote the coefficient's transformation properties under $SL(2, \ints)$ as described below. We let \eqref{eq:SUSY-constraints:full-S3} define the relative normalisations between the coefficients.

Taking the weak coupling limit, we can relate different terms in the full action using the linearised theory. By reading off the coefficients in the expansion of $F[\Phi]$ in \eqref{eq:SUSY-constraints:S3}, letting $y^{-1} \sim \gs \to 0$ we see that \cite{Green:1997me, Green:1999qt}
\begin{equation}
    \label{eq:SUSY-constraints:linear-relation}
    f^{(12,-12)}(z) \sim y^{12} \partial_z^{12} f^{(0,0)}(z) \, ,
\end{equation}
up to a numerical factor.

However, the coefficient $f^{(0,0)}(z) = \mathcal{E}^{\ordd{10}}_{\gra{0}{0}}(z)$ is invariant under $SL(2, \ints)$ transformations, and the same U-duality arguments give that each of the terms in \eqref{eq:SUSY-constraints:full-S3} is $SL(2, \ints)$ invariant. Due to the transformation properties of $\lambda^{16}$, $\hat G \lambda^{14}$ and the remaining interaction terms, the coefficient functions $f^{(w, \hat w)}(z)$ have to be modular forms with holomorphic weight $w$ and anti-holormorphic weight $\hat w$ transforming as
\begin{equation}
    f^{(w, \hat w)}\!\left(\frac{a z + b}{c z + d}\right) = 
    (c z + d)^{w} (c \overline z + d)^{\hat w} f^{(w, \hat w)}(z) \qquad
    \begin{pmatrix}
        a & b \\
        c & d 
    \end{pmatrix} \in SL(2, \ints).
\end{equation}

This means that \eqref{eq:SUSY-constraints:linear-relation} cannot be the complete relation since $y^{12} \partial_z^{12} f^{(0,0)}(z)$ does not transform as a modular form of holomorphic and anti-holormorphic weights $(12, -12)$. In physics parlance, the derivative is not covariant with respect to R-symmetry. Instead, we need the modular covariant derivative
\begin{equation}
    \mathcal{D}_{(w)} = \frac{\partial}{\partial z} - i \frac{w}{2y}  
\end{equation}
which maps a modular form of weights $(w, \hat w)$ to a modular form of weights $(w+2,\hat w)$. Multiplying with $y = \Im z$ which transforms as
\begin{equation}
    \Im \frac{az+b}{cz+d} = \frac{\Im z}{|cz+d|^2}
\end{equation}
we have that $y \mathcal{D}_{(w)} f^{(w, \hat w)}(z)$ transforms with weights $(w+1, \hat w -1)$.

Then, we may covariantise \eqref{eq:SUSY-constraints:linear-relation} to \cite{Green:1997me} 
\begin{equation}
    \label{eq:SUSY-constraints:old-covariant-relation}
    f^{(12, -12)}(z) \propto y^{12} \mathcal{D}_{(22)} \mathcal{D}_{(20)} \cdots \mathcal{D}_{(0)} f^{(0,0)}(z) \, ,
\end{equation}
where we do not need the numerical proportionality constant.

\begin{remark}
We note that the covariant derivative $\mathcal{D}_{(w)}$ reduces to the ordinary derivative for the terms in $f^{(w,-w)}$ which satisfy
\begin{equation}
    \frac{\partial}{\partial z} f(z) \gg \frac{w}{y} f(z) \qquad \text{as } y \to \infty \, .
\end{equation}
This is not satisfied for the perturbative terms in $f^{(w,-w)}$, which are highlighted in \eqref{eq:intro-SL2-fourier-expansion} for $f^{(0,0)}$. It is, however, satisfied for the non-perturbative terms in \eqref{eq:intro-SL2-fourier-expansion} meaning that the linearised theory 
\eqref{eq:SUSY-constraints:S3} can only capture the leading instanton contributions to the effective action \cite{Green:1999qt}.
\end{remark}

To simplify the notation of \eqref{eq:SUSY-constraints:old-covariant-relation} we also use that $\mathcal{D}_{(w)}(y^n f^{(w+n,\hat w)}) = y^n \mathcal{D}_{(w+n)} f^{(w+n,\hat w)}$ which allows us to write \eqref{eq:SUSY-constraints:linear-relation} as
\begin{equation}
    \label{eq:SUSY-constraints:covariant-relation}
    \begin{split}
        f^{(12, -12)}(z) &\propto y \mathcal{D}_{(11)} \, y \mathcal{D}_{(10)} \cdots y \mathcal{D}_{(0)} f^{(0,0)}(z)  \\
        &= D_{(11)} D_{(10)} \cdots D_{(0)} f^{(0,0)}(z)
    \end{split}
\end{equation}
with
\begin{equation}
    D_{(w)} = i\Big(y \frac{\partial}{\partial z} - i \frac{w}{2} \Big) = i y \mathcal{D}_{(w)}
\end{equation}
where we have introduced an extra factor of $i$ to use the same normalisation as in the existing literature. This covariant derivative takes a modular form of weights $(w, \hat w)$ to a modular form of weights $(w+1, \hat w-1)$.

The coefficient functions with negative holomorphic weights (such as $f^{(-12,12)}$) can be related to $f^{(0,0)}$ using the conjugated covariant derivative
\begin{equation}
    \overline D_{(\hat w)} := \overline{D_{(\hat w)}} = -i \Big( y \frac{\partial}{\partial \overline z} + i \frac{\hat w}{2}\Big)
\end{equation}
which maps a modular form of weights $(w, \hat w)$ to a modular form of weights $(w-1, \hat w+1)$. 
\begin{remark}
    The modular forms above transforming with holomorphic and anti-holomorphic weights $(w, -w)$ are so-called \emphindex[Maass form]{Maass forms} of weight $2w$ which are discussed in sections \ref{sec:Maasswt} and \ref{standardsectionSL2} where the modular covariant derivatives are closely related to the raising and lowering operators of \eqref{eq:Maass-lowering-raising}.
\end{remark}

We will now use the full non-linear, local supersymmetry of the theory to determine the coefficient functions. When adding the $(\alpha')^3$ terms \eqref{eq:SUSY-constraints:full-S3} to our effective action we also need to correct the supersymmetry transformation $\delta_\epsilon$, with local spinor parameter $\epsilon$, acting on the different fields. Schematically, for some generic field $\Psi$ we would then have the expansion
\begin{equation}
    \label{eq:SUSY-constraints:delta}
    \delta_\epsilon \Psi = \Big( \delta^{(0)} + \alpha' \delta^{(1)} + (\alpha')^2 \delta^{(2)} + \ldots \Big) \Psi \, .
\end{equation}

Requiring that the action should be invariant under this supersymmetry transformation (up to total derivatives) we get the following first few conditions
\begin{equation}
    \delta^{(0)} S^{(0)} = 0 \qquad \delta^{(0)} S^{(3)} + \delta^{(3)} S^{(0)} = 0 \qquad \ldots
\end{equation}
where the first equation simply states that the classical supergravity action should be invariant under the ordinary local supersymmetry transformations.

To determine the coefficient functions $f^{(w, -w)}(z)$ it turns out that we do not have to analyse the supersymmetry transformations of every term in the action --- we only need to focus on a few interaction terms with high powers of the dilatino $\lambda$. Expanding $\hat G \lambda^{14}$ in terms of the physical fields, we get, among others, a term $\lambda^{15} \gamma^\mu \psi_\mu^*$ \cite{Green:1998by}
\begin{equation}
    \label{eq:SUSY-constraints:S3-highlighted}
    S^{(3)} = \int d^{10}x \sqrt{g} \Big( f^{(12,-12)}(z) \lambda^{16} - 3 \cdot 144 f^{(11, -11)}(z) \lambda^{15} \gamma^\mu \psi_\mu^* + \ldots \Big) \, .
\end{equation}

When collecting the terms of $\delta_{\epsilon_1}^{(0)} S^{(3)}$ proportional to $\lambda^{16} \overline \epsilon^*_1 \gamma^\mu \psi_\mu^*$ the two highlighted terms in \eqref{eq:SUSY-constraints:S3-highlighted} are the only contributing interactions giving \cite{Green:1998by}
\begin{equation}
    \delta_{\epsilon_1}^{(0)} S^{(3)} = - i \int d^{10}x \sqrt{g} (\lambda^{16} \overline \epsilon^*_1 \gamma^\mu \psi_\mu^*) \Big(8 f^{(12, -12)}(z) + 6 \cdot 144 D_{(11)} f^{(11,-11)}(z) \Big) + \ldots
\end{equation}

There are no possible $\delta^{(3)}$ variations which, when acting on $S^{(0)}$, would give a contribution proportional to $\lambda^{16} \overline \epsilon^*_1 \gamma^\mu \psi_\mu^*$ \cite{Green:1998by}. Hence,
\begin{equation}
    \label{eq:SUSY-constraints:f12f11}
    \begin{gathered}
    \delta_{\epsilon_1}^{(0)} S^{(3)} +  \delta_{\epsilon_1}^{(3)} S^{(0)} = 0 \\
\quad \implies \quad     D_{(11)} f^{(11, -11)}(z) = - \frac{4}{3 \cdot 144} f^{(12, -12)}(z) \, .
    \end{gathered}
\end{equation}

If we instead collect the terms of $\delta_{\epsilon_2^*}^{(0)} S^{(3)}$ proportional to $\lambda^{16} \overline \epsilon_2 \lambda^*$ (again with contributions coming only from the two highlighted terms in \eqref{eq:SUSY-constraints:S3-highlighted}) we get \cite{Green:1998by}
\begin{equation}
    \delta^{(0)}_{\epsilon_2^*} S^{(3)} = - 2i \int d^{10}x \sqrt{g} (\lambda^{16} \overline \epsilon_2 \lambda^*) \Big( \overline D_{(12)} f^{(12, -12)}(z) - 3 \cdot 144 \cdot \frac{15}{2} f^{(11, -11)}(z) \Big) + \ldots
\end{equation}

However, in this case we might also get a contribution from $\delta_{\epsilon_2^*}^{(3)} S^{(0)}$ if we assume that the $(\alpha')^3$ variation of $\lambda^*$ is of the form
\begin{equation}
    \label{eq:SUSY-constraints:delta-3}
    \delta_{\epsilon_2^*}^{(3)} \lambda^*_a = - \frac{i}{6} h(z) (\lambda^{14})_{cd} (\gamma^{\mu\nu\rho}\gamma^0)_{dc} (\gamma_{\mu\nu\rho} \epsilon_2^*)_a
\end{equation}
where $h(z)$ is an unkown function of the axio-dilaton. This gives a contribution \cite{Green:1998by}
\begin{equation}
    \delta_{\epsilon_2^*}^{(3)} S^{(0)} = 180i \int d^{10}x \sqrt{g} (\lambda^{16} \overline \epsilon_2 \lambda^*) h(z) + \ldots
\end{equation}

From the $\epsilon_2^*$ variation at order $(\alpha')^3$ we then get 
\begin{equation}
    \label{eq:SUSY-constraints:f12f11h}
    \begin{gathered}
    \delta_{\epsilon_2^*}^{(0)} S^{(3)} + \delta_{\epsilon_2^*}^{(3)} S^{(0)} = 0 \\ 
\quad \implies \quad     \overline D_{(-12)} f^{(12,-12)}(z) - 3 \cdot 144 \cdot \frac{15}{2} f^{(11,-11)}(z) - 90 h(z) = 0
    \end{gathered}
\end{equation}

The two equations \eqref{eq:SUSY-constraints:f12f11} and \eqref{eq:SUSY-constraints:f12f11h} give us a system of ordinary differential equations for $f^{(12,-12)}$ and $f^{(11,-11)}$ together with the unknown $h(z)$ which we will now determine by requiring that the $(\alpha')$-corrected supersymmetry algebra 
\eqref{eq:SUSY-constraints:delta} closes. That is, we require that the commutator of two local supersymmetry transformations acting on some field $\Psi$ gives a sum of local symmetry transformations that are present in the theory up to equations of motion:
\begin{equation}
    [\delta_{\epsilon_1}, \delta_{\epsilon_2^*}] \Psi =
    \delta_\text{local translation}
    \Psi + \delta_\text{local symmetries} \Psi + (\text{equations of motion})
\end{equation}

Computing the commutator acting on $\lambda^*$ to order $(\alpha')^3$ using our ansatz
\eqref{eq:SUSY-constraints:delta-3} for $\delta^{(3)}$ and comparing with the equations of motion for $\lambda^*$ from the corrected action \eqref{eq:SUSY-constraints:full-S3} one finds that \cite{Green:1998by, Sinha:2002zr}
\begin{equation}
    \label{eq:SUSY-constraints:closure}
    -32 D_{(11)} h(z) = f^{(12,-12)}(z) \, .
\end{equation}
Inserting this in \eqref{eq:SUSY-constraints:f12f11h} multiplied by $4 D_{(11)}$ on the left and using \eqref{eq:SUSY-constraints:f12f11} to eliminate $D_{(11)} f^{(11,-11)}(z)$ we obtain the following eigenvalue equation for $f^{(12,-12)}(z)$
\begin{equation}
    \begin{split}
        \MoveEqLeft
        4 D_{(11)} \overline{D}_{(-12)} f^{(12, -12)}(z) = \\
        &= -4\Big(-3 \cdot 144 \cdot \frac{15}{2} D_{(11)} f^{(11,-11)}(z) - 90 D_{(11)} h(z) \Big)\\
        &= -4\left(-3 \cdot 144 \cdot \frac{15}{2} \Big(-\frac{4}{3 \cdot 144}  f^{(12,-12)}(z) \Big) - 90 \Big( -\frac{1}{32} f^{(12,-12)}(z) \Big) \right) \\
        &= \Big(-132 + \frac{3}{4} \Big) f^{(12,-12)}(z)
    \end{split}
\end{equation}

Then, using the relation \eqref{eq:SUSY-constraints:covariant-relation} between $f^{(12,-12)}(z)$ and $f^{(0,0)}(z)$ we get that
\begin{equation}
    4 y^2 \partial_{z} \partial_{\overline z} f^{(0,0)}(z) = \Delta_{\UHP} f^{(0,0)}(z) = \frac{3}{4} f^{(0,0)}(z)
\end{equation}
which is exactly the eigenvalue equation \eqref{eq:10D-eigeneq} for $\mathcal{E}_{\gra{0}{0}}^{\ordd{10}}$ with $\Delta_{\UHP}$ being the Laplace--Beltrami operator \eqref{DeltaUHP} on the upper half plane. Note that we can now obtain eigenvalue equations for all coefficient functions $f^{(w, -w)}$ at order $(\alpha')^3$ by using the modular covariant derivatives $D_{(w)}$ and $\overline D_{(-w)}$.

With similar arguments requiring invariance of the action and supersymmetry algebra closure at order $(\alpha')^5$ one can find the eigenvalue equation \eqref{eq:10D-eigeneq} for the $D^4 R^4$ term in ten dimensions as was done in \cite{Sinha:2002zr}. The eigenvalue equations \eqref{eq:intro-fourgravLap} for the $R^4$, $D^4 R^4$ and $D^6 R^4$ terms in arbitrary dimensions were given in \cite{Green:2010wi}. 
In~\cite{Bossard:2014lra,Bossard:2014aea,Bossard:2015uga}, an approach using linearised harmonic superspace for maximal supersymmetry was used to investigate possible supersymmetry invariants in arbitrary dimensions. The resulting equations for the coefficient functions are expressed in terms of tensorial differential operators on moduli space and specialise to the scalar Laplace equations presented here in equation~\eqref{eq:intro-fourgravLap}. The tensorial equations are more constraining than the scalar ones and are closely related to the associated varieties of small automorphic representations as discussed in section~\ref{smallreps}.

\section[Computing the four-graviton tree level amplitude*]{Computing the four-graviton tree level \\ amplitude*}
\label{sec:four-graviton-tree-amplitude}
    
A comprehensive discussion of string theory and its scattering amplitudes is beyond the scope of this work. We only give one indicative and, hopefully, illustrative example and refer to the string theory literature~\cite{GSW,Polchinski} for more information. 

\begin{figure}[t]
    \centering

    \includegraphics{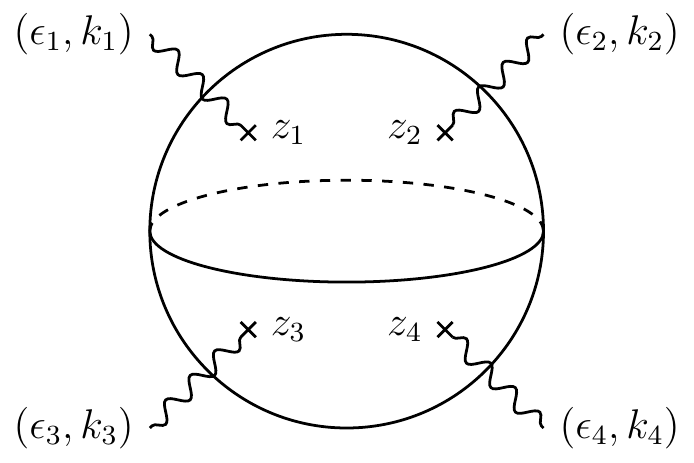}

    \caption{\label{fig:FourGravitons}\sl Riemann sphere with four punctures as it appears in the four graviton scattering amplitude. An external massless graviton with polarisation $\epsilon_i$ and momentum $k_i$ is located at each of the four punctures $z_i\in \cx\mathbb{P}^1$.}
\end{figure}

The example is the four-graviton amplitude at string~\emphindex[string theory!perturbation theory!tree level]{tree level}. The closed string tree level topology is that of a sphere, and the four asymptotic graviton states correspond to four punctures in this sphere as pictured in figure~\ref{fig:FourGravitons}. This configuration can be obtained by a homeomorphism of the left-most diagram in figure~\ref{fig:LowGenus}. By definition, the string scattering amplitude is given by an integral over the moduli space of all Riemann spheres and all possible insertion points for four punctures.  The discussion below uses the \emphindex[Ramond--Neveu--Schwarz formalism]{Ramond--Neveu--Schwarz (RNS) formalism}.

In terms of complex geometry, we can describe the sphere $S^2=\cx \mathbb{P}^1$ by one complex variable $z\in\cx$ everywhere except at the `north-pole'.  An asymptotic graviton state at a puncture $z_i$ corresponds to a \emphindex{vertex operator}. Due to a complication called \emphindex{ghost picture}~\cite{Friedan:1985ge} that arises for the superstring, we actually require it in two related forms, one called ghost picture $0$ and ghost picture $-1$:
\begin{subequations}
\label{VOgrav}
\begin{align}
\label{eq:Pic0}
V_0(z_i; k_i, \epsilon_i) &= \frac{2}{\alpha'} \gs \normord{\epsilon_{i,\mu\nu}\left( \partial X^\mu +  \frac{\alpha'}{2}k_\rho \psi^\rho \psi^\mu\right)\left(\bar{\partial} X^\nu + \frac{\alpha'}{2} k_\sigma \bar{\psi}^\sigma \bar{\psi}^\nu \right) e^{i k_{i,\mu} X^\mu}}\,,\\
\label{eq:Pic-1}
V_{-1}(z_i; k_i, \epsilon_i) &= \gs \, \normord{\epsilon_{i,\mu\nu}\psi^\mu \bar{\psi}^\nu  e^{i k_{i,\mu} X^\mu}}\,.
\end{align}
\end{subequations}
As before, we are using the Einstein summation convention for repeated Lorentz indices $\mu=0,\ldots,9$ that label the ten directions of the Minkoswki target space in which the sphere is embedded via the embedding coordinates $X^\mu\equiv X^\mu(z,\bar{z})$. The polarisation tensor $\epsilon_i\equiv \epsilon_{i,\mu\nu}$ is a second-rank symmetric tensor of the Lorentz group $SO(1,9)$ and the external momentum of the graviton $k_i\equiv k_{i,\mu}$ is light-like in the Minkowski metric: $k_i^2 = k_{i,\mu} k_{i,\nu} \eta^{\mu\nu}=0$ as the scattering gravitons are massless. The colons surrounding this expression indicate a specific normal ordering procedure necessary for the vertex operator to be well-defined on a Fock space vacuum. The field $\psi^\mu$ and its (Dirac) conjugate $\bar{\psi}^\mu$ correspond to the fermionic coordinates that accompany the bosonic $X^\mu$ in superstring theory. The derivatives $\partial$ and $\bar{\partial}$ are with respect to the world-sheet coordinate $z$.

The desired expression for the scattering amplitude is then roughly of the following form for $D=10$ (in string frame)
\begin{align}
\label{treeampl}
\mathcal{A}^{(10)}_{\mathrm{tree}}(s,t,u,\epsilon_i; g) = \lint_{\mathclap{\mathcal{M}_{0,4}}} d\mu \, \langle \prod_{i=1}^4 V(z_i; k_i,\epsilon_i)\rangle_{S^2},
\end{align}
where the angled brackets denote the \emphindex{correlation function} of the vertex operators on the given sphere $S^2$  that we detail below. For reasons of ghost number saturation, two vertex operators have to be taken in ghost picture $0$ (of the form $V_0$ in~\eqref{VOgrav}) and two in ghost picture $-1$ (of the form $V_{-1}$ in~\eqref{VOgrav}). The integral $d\mu$ is over all Riemann spheres (genus $h=0$) with four punctures, so that the positions $z_i$ of the punctures is also integrated over. More concretely, the measure is given by the integral over all metrics that can be put on topological spheres up to diffeomorphisms and Weyl rescalings (local dilatations). To make this well-defined in the path-integral sense one has to divide by the volume of this gauge group so that
\begin{align}
d\mu = \frac{D\gamma_{\cx\mathbb{P}^1}}{\textrm{Vol}(\textrm{diff}\times\textrm{Weyl})}.
\end{align}
Here, $D\gamma_{\cx\mathbb{P}^1}$ indicates all possible metrics on the sphere $\cx\mathbb{P}^1$. According to the \emphindex{Riemann--Roch theorem}, the sphere as a Riemann surface has no metric moduli and any metric $\gamma_{\cx\mathbb{P}^1}$ can be brought into the form of that of the round sphere
\begin{align}
\label{roundmetric}
ds^2 = \frac{ dz d\bar{z}}{(1+|z|^2)^2}
\end{align}
by diffeomorphisms and Weyl rescalings. There is no modulus in this expression and even this form is left invariant by the \emphindex[Killing!group, conformal]{conformal Killing group} of the sphere $PSL(2,\cx)$ that acts by
\begin{align}
\begin{pmatrix}\alpha & \beta \\\gamma &\delta \end{pmatrix}\in PSL(2,\cx) : \quad\quad
z \mapsto \frac{\alpha z+ \beta}{\gamma z + \delta}
\end{align}
on the coordinate $z$. So, even after fixing the $\textrm{diff}\times\textrm{Weyl}$ gauge-freedom to bring the metric to the above form~\eqref{roundmetric}, one has still the freedom to perform transformations from the conformal Killing group $PSL(2,\cx)$ that represent the residual gauge freedom. Without fixing it the integral in~\eqref{treeampl} is of the form
\begin{align}
\lint_{\mathclap{\mathcal{M}_{0,4}}} d\mu = \prod_{j=1}^4 \,\,\lint_{\mathrlap{\cx \mathbb{P}^1}} d^2z_j \, \Delta_{\textrm{FP}},
\end{align}
where $\Delta_{\textrm{FP}}$ is a \emphindex{Faddeev--Popov determinant} arising from the gauge-fixing and that is treated by introducing (super-)ghost systems. We will not be more specific on it here (see for example~\cite{Polchinski}) and only mention its effect on the calculation below. We see that, due to the absence of metric moduli for the sphere, the integral over the moduli space of four-punctured sphere reduces to an integral over the locations $z_i$ of the four punctures. The complex three-dimensional conformal Killing group can be used to fix three of the puncture positions to one's favourite values; a standard choice being $0,1,\infty$ and only a single integral over a single puncture position remains.

Let us return to the correlation function appearing in the expression~\eqref{treeampl} above. It is formally given by a path-integral over all (super-)embeddings $X^\mu$ of the (super-)Riemann sphere into the ten-dimensional target space. Schematically, one has
\begin{align}
 \langle \prod_{i=1}^4 V(z_i; k_i,\epsilon_i)\rangle_{S^2} = \gs^{-\chi(S^2)} \int DX D\psi D\bar{\psi}\, e^{-S[X,\psi,\bar{\psi}] }  \prod_{i=1}^4 V(z_i; k_i,\epsilon_i),
\end{align}
where $S[X,\psi,\bar{\psi}]$ denotes the two-dimensional \emphindex[sigma@$\sigma$-model]{$\sigma$-model} action that is basically the induced volume under the embedding. The \emphindex{Euler number} $\chi(S^2)=2$ of the sphere provides the standard topological weighting of different string diagrams that was mentioned above (cf. figure~\ref{fig:LowGenus}) and here evaluates to $\gs^{-2}$. The string coupling $\gs$ is the only string theory modulus $g$ of equation~\eqref{stringmoduli} of relevance in the present discussion. As our goal here is to give a heuristic derivation of the final formula~\eqref{fourtree} below, we are not displaying or discussing in any more detail aspects related to so-called \emph{pictures} of vertex operators associated with the ghosts arising from gauge-fixing. A proper treatment would modify the above equation~\cite{GSW,Polchinski}; this modification will be taken into account in the final expression below.

The correlation function can be evaluated on the sphere explicitly in terms of the Green's function on the sphere. First, we note that for vertex operators in ghost picture $0$ the fermionic integrals (over $\psi$ and $\bar{\psi}$) pick out the contribution 
\begin{align}
\alpha' \epsilon_{i,\mu\nu} k_\rho k_\sigma e^{i k_{i,\tau} X^\tau}\, \propto  \alpha' \mathcal{R}_{\rho\mu\nu\sigma} e^{ik_{i,\tau} X^{\tau}}
\end{align}
from $V(z_i; k_i, \epsilon_i)$ in~\eqref{eq:Pic0} and we recognise the linearised Riemann curvature tensor mentioned in section~\ref{sec:intro-four-graviton}, and the full integral also provides the necessary contractions of the four Riemann tensors. For the vertex operator in picture $-1$ of~\eqref{eq:Pic-1}, the integral yield contributions that are roughly of the form
\begin{align}
\label{m1contr}
\epsilon_{i,\mu\nu}  e^{i k_{i,\tau} X^\tau}\,\,\propto \frac{\mathcal{R}_{\rho\mu\nu\sigma}}{k_\rho k_\sigma} e^{ik_{i,\tau} X^{\tau}}.
\end{align}
The Green's function on the sphere now evaluates the product of two normal ordered exponentials to 
\begin{align}
\normord{e^{i k_{i,\mu} X^\mu(z_i)}} \; \normord{e^{i k_{j,\nu} X^\nu(z_j)}} \,\,\propto |z_i-z_j|^{\alpha' k_i\cdot k_j} \normord{e^{i (k_1+k_2)_\mu X^{\mu}(z_i)}} + \ldots
\end{align}
which is sometimes called the \emphindex{operator product expansion} (or \emphindex{Wick contraction}) and the dots denote subdominant terms in an expansion around $z_i\to z_j$. In this expression we have used the Lorentz inner product $k_i\cdot k_j= k_{i,\mu} k_j^\mu$ between the two external momenta $k_i$ and $k_j$. Since the external momenta are massless, we can rewrite this product as $k_i\cdot k_j = \frac12 \left(k_i+k_j\right)^2$ and we already begin to recognise a connection to the Mandelstam variables $s$, $t$ and $u$ defined in~\eqref{Mandeldef}. One has to perform the above operator product expansion of $e^{ikX}$ for all four factors appearing in the correlation function and this will lead to a permutation invariant expression of the three Mandelstam variables. 

We can now put all the pieces together: $(i)$ The integral over the moduli space reduces to an integral over the four punctures, $(ii)$ the conformal Killing groups allows to fix three of four punctures to fixed values that we choose to be $0$, $1$ and $\infty$, $(iii)$ the integral over the fermionic variables in the correlation function produces the linearised curvature tensors times exponentials of the form $e^{ikX}$, $(iv)$ these exponentials get converted into factors of the form $|z_i-z_j|^{\alpha' k_i\cdot k_j}$ in all possible ways and $(v)$ include additional contributions from the ghost sector. The ultimate integral is of the form 
\begin{align}
\label{fourtree}
\mathcal{A}^{(10)}_{\mathrm{tree}}(s,t,u,\epsilon_i; g) &=  \frac{4(\alpha')^2\gs^{2}}{\pi} \frac{\mathcal{R}^4}{(k_1\cdot k_3)^2} \delta(k_1+k_2+k_3+k_4) \!\int_{\cx}  \!\! d^2z{|z|^{\alpha' k_1\cdot k_2-2}} |1-z|^{\alpha' k_2\cdot k_3-2}\nn\\
&=  (\alpha')^4\gs^{2}\delta(k_1+k_2+k_3+k_4) \frac{1}{stu} \frac{\Gamma(1-s) \Gamma(1-t)\Gamma(1-u)}{\Gamma(1+s)\Gamma(1+u)\Gamma(1+t)} \mathcal{R}^4 \, .
\end{align}
We have introduced a few factors related to the normalisations of the various measures and vertex operators introduced above.  The prefactor $\gs^2$ is the combination of the tree level value $\gs^{-2}$ and the $\gs^4$ from the four vertex operators of the external graviton states. The shift in the $\alpha'$ power is due to the conformal ghost sector that we have not discussed explicitly and the denominator $(k_1\cdot k_3)^2= \tfrac{4}{(\alpha')^2} t^2$ is due to~\eqref{m1contr}.  We note that the factor $\delta(k_1+k_2+k_3+k_4)$ expresses momentum conservation. As stated in section~\ref{sec:intro-four-graviton}, the linearised curvature tensor $\mathcal{R}_{\mu\nu\rho\sigma} \propto k_{\mu}\epsilon_{\nu\rho} k_\sigma$ is determined by the graviton's momentum and polarisation, and $\mathcal{R}^4$ is a specific contraction of the linearised curvatures of the four gravitons. More precisely, it is given by $t_8 t_8 \mathcal{R}^4$ where the $t_8$ tensor contracts four powers of an antisymmetric matrix $M_{\mu\nu}$ according to $t^{\mu_1\ldots\mu_8} M_{\mu_1\mu_2}\cdots M_{\mu_7\mu_8} = 4\mathrm{Tr} (M^4) - \left( \mathrm{Tr} (M^2)\right)^2$. 

The amplitude~\eqref{fourtree} can be expanded for small values of $s$, $t$ and $u$ which we recall from~\eqref{Mandeldef} contain the string length $\alpha'=\ells^2$. The result is (in string frame)
\begin{align}
\label{fourgravtree}
\mathcal{A}^{(10)}_{\mathrm{tree}}(s,t,u,\epsilon_i; g) 
=\gs^{2}(\alpha')^4\delta(k_1+\ldots+k_4)\left(\frac{3}{\sigma_3}
 +2\zeta(3)  + \zeta(5)\sigma_2 +\frac{2 \zeta(3)^2}{3} \sigma_3 + \ldots \right) \mathcal{R}^4,
\end{align}
where we have used $stu=\frac{1}{3} \sigma_3$ up to momentum conservation implying $s+t+u=0$ for massless states. The above expression provides the tree level contributions to the functions $\mathcal{E}_{\gra{0}{-1}}^{\ordd{10}}$, $\mathcal{E}_{\gra{0}{0}}^{\ordd{10}}$, $\mathcal{E}_{\gra{1}{0}}^{\ordd{10}}$  and $\mathcal{E}_{\gra{0}{1}}^{\ordd{10}}$ appearing in~\eqref{fourgrav} in ten dimensions. We note also that $\gs^2(\alpha')^4= \ell_{10}^8 =\kappa_{10}^2$ in terms of the ten-dimensional Planck scale $\ell_{10}$ and the ten-dimensional gravitational coupling $\kappa_{10}$. The fact that this amplitude is proportional to the square of the coupling is characteristic for a tree level scattering of four particles (in a theory that has cubic vertices).

We note also upon toroidal compactification to  $D<10$ dimensions integrals such as~\eqref{fourtree} receive additional contributions related to the structure of states on the torus $T^d$ ($d=10-D$) that is used in the compactification. In the context of loop amplitudes, this leads naturally to a relation of string amplitudes to \emphindex[theta correspondence]{theta correspondences} which are discussed in section \ref{sec:ThetaCorr}.

\section{Non-perturbative corrections from instantons*}
\label{sec:instantons}

As stated in section \ref{sec:physical-interp-10D}, D-instantons give rise to non-perturbative corrections to the effective action of the form $e^{-1/\gs}$. We will now motivate this by studying string amplitudes in ten dimensions in the presence of D-instantons following \cite{Polchinski:1994fq, Green:1995my}, but first let us introduce the concept of instantons in field theory.

Instantons in a (Wick rotated, Euclidean) field theory are solutions to the equations of motion with finite action. In the semi-classical approximation the \emphindex{partition function} $Z$, which is (loosely speaking) the generating function of scattering amplitudes, is approximated as a sum over local extrema to the action. 
\begin{equation}
\label{eq:Zapprox}
    Z \approx \sum_\text{extrema} e^{-S[\text{extremum}]} \times (\text{quantum fluctuation effects around extremum}) \, .
\end{equation}
The global minimum corresponds to the true vacuum with the usual perturbative corrections while other local minima give instantons corrections.

In ten-dimensional, type IIB supergravity, the following configurations are spherical instanton solutions to the equations of motion \cite{Gibbons:1995vg}
\begin{equation}
    g_{\mu\nu} = \delta_{\mu\nu} \qquad e^{\phi} - e^{\phi_0} = \frac{\abs{q}}{8 \Vol(S^9)} \frac{1}{r^8} \qquad \chi - \chi_0 = \sgn(q) (e^{-\phi} - e^{-\phi_0}) \, ,
\end{equation}
where $\phi_0$ is the value of $\phi$ at $r = \infty$ with $\gs = e^{\phi_0}$ and similarly for $\chi_0$. $\Vol(S^9)$ is the volume of $S^9$ with unit radius, and $q$ is the Ramond--Ramond charge of the solution. For negative $q$ the solution is called an \emphindex[instanton!anti-instanton]{anti-instanton} while for positive $q$ it is called an \emphindex{instanton}. This $q$ is the Noether-charge with respect to translations in $\chi$ and is obtained by a closed hypersurface integral in the Euclidean space-time, enveloping the origin. The integral is invariant under deformations of the surface as long as we do not cross the origin and the instanton is in this way localised to a single point in space-time. Due to \emphindex[Dirac--Schwinger--Zwanziger quantisation]{Dirac--Schwinger--Zwanziger quantisation} this charge becomes quantised as $q = 2 \pi m$ where $m$ is an integer which we colloquially also call the instanton charge.

The value of the action for these solutions are \cite{Green:1997tv} (compare with \eqref{eq:Instaction})
\begin{equation}
    \label{eq:SUGRA-instanton-action}
    S_\text{inst} = 2\pi \abs{m} \gs^{-1} - 2\pi i m \chi_0 \, , 
\end{equation}
which for an instanton and anti-instanton becomes $-2\pi i \abs{m} z_0$ and $2\pi i  \abs{m} \overline{z}_0$ respectively, where $z_0 = \chi_0 + i \gs^{-1}$. This is the value for $S[\text{extremum}]$ that enters in~\eqref{eq:Zapprox}.

In string theory, the objects that correspond to these instanton solutions in supergravity are D-instantons. As discussed in section \ref{sec:physical-interp-10D}, D-instantons are special cases of D-branes on which the endpoints of open strings attach, and which are localised to a single point in space-time.

The contributions to the four-graviton scattering amplitude in string theory shown schematically in figure \ref{fig:LowGenus} only include closed string world-sheets, that is, world-sheets without boundaries (but with punctures for external states). However, when adding D-branes into the picture, the sum over topologies should also include open strings whose world-sheets have boundaries on the D-branes. While the closed string world-sheets give perturbative corrections in $\gs$ to the amplitudes, the open string world-sheets give non-perturbative corrections.

Both the perturbative and non-perturbative contributions are present in the ten-dimensional four-graviton scattering amplitude $\mathcal{A}^{(10)}$ of \eqref{fourgrav} for which we computed the tree level, perturbative term in section \ref{sec:four-graviton-tree-amplitude}. To go beyond the perturbative corrections, we write the amplitude as a sum over the number $n$ of D-instantons at positions $\{y_i\}_{i=1}^{n}$ in space-time
\begin{equation}
    \mathcal{A}^{(10)} = \sum_{n = 0}^{\infty} \mathcal{A}_n \, ,
\end{equation}
where $\mathcal{A}_0$ contains all the perturbative contributions shown in figure \ref{fig:LowGenus} including the tree level term. The positions $y_i$ are later integrated over to restore translational invariance. This integration also imposes momentum conservation of the four graviton states for the non-perturbative corrections \cite{Green:1997tv}. 

Each $\mathcal{A}_n$ is, in turn, a sum over different world-sheets which we now allow to be disconnected and have boundaries which are attached to the points $y_i$ in space-time. We will call the disjoint parts \emphindex[world-sheet!components]{world-sheet components} which are, by themselves, connected but can have holes and handles.

\begin{remark}
We may have a world-sheet with boundaries attached to the same D-instantons that seems disconnected in the view of the world-sheet but that is connected from the space-time point of view --- the different components being connected by the D-instantons. As in ordinary field theory, the scattering amplitudes can be shown to only involve diagrams which are connected when viewed in space-time \cite{Green:1995my}.
\end{remark}

\begin{figure}[t]
\centering
\begin{tikzpicture}[thick]
    \def\diameter{1.5cm}
    \def\sd{0.4cm}
    \node at (0, 0) [draw, circle, minimum size=\diameter, inner sep=0pt, label=left:{$\gs^{-1}$}] {};

    \node at (3, 0) [draw, circle, minimum size=\diameter, inner sep=0pt, label=left:{$\dfrac{1}{2!}$}] {};
    \node at (3, 0) [draw, circle, minimum size=\sd, inner sep=0pt] {};

    \node at (6, 1) [draw, circle, minimum size=\diameter, inner sep=0pt, label=left:{$\dfrac{\gs}{3!}$}] {};
    \node at (5.7, 1) [draw, circle, minimum size=\sd, inner sep=0pt] {};
    \node at (6.3, 1) [draw, circle, minimum size=\sd, inner sep=0pt] {};

    \begin{scope}[even odd rule]
        \clip [domain=-90:90] (3,1) rectangle (7,-2)
        plot ({6.3+1.0*cos(\x)}, {-1+0.4*sin(\x)}) -- plot ({6.3+0.8*cos(-\x)}, {-1+0.2*sin(-\x)});
        \node at (6, -1) [draw, circle, minimum size=\diameter, inner sep=0pt, label=left:{$\gs$}] {}; 
    \end{scope}

    \draw [domain=-90:90] plot ({6.3+1.0*cos(\x)}, {-1+0.4*sin(\x)});
    \draw [domain=-90:90] plot ({6.3+0.8*cos(0.95*\x)}, {-1+0.2*sin(0.95*\x)});

\end{tikzpicture}
\caption{Leading order open world-sheet components with topological weights and symmetry factors for identical boundaries. Each boundary is contracted to the same D-instanton at a point $y_1$ in space-time. The last, lower diagram is a disk with a single boundary and a handle.}
\label{fig:open-worldsheets} 
\end{figure}
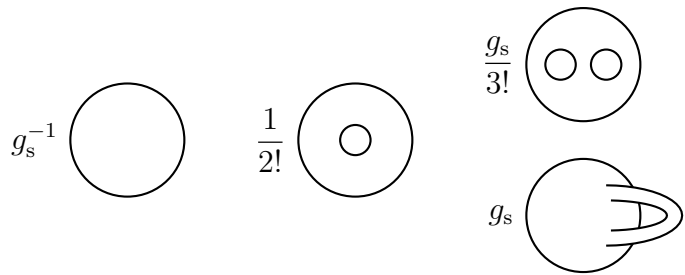

Focussing on the single D-instanton case ($n=1$) at position $y_1$, the contribution $\mathcal{A}_1$ is a sum over disconnected world-sheets with boundaries at the point $y_1$ together with appropriate symmetry factors for exchanging identical world-sheet components and boundaries. Similarly to closed strings, the world-sheet components are weighted by the Euler characteristic as $\gs^{2(h-1)+b}$ where $h$ is the genus of the surface and $b$ is the number of boundaries. The leading order world-sheet components together with their symmetry factors for exchanging identical boundaries are shown in figure \ref{fig:open-worldsheets}. Besides these components we will also have surfaces with punctures corresponding to insertions of external states.

\newcommand{\disk}{\tikz [baseline] \node at (0,0.6ex) [draw, thick, circle, minimum size=1em, inner sep=0pt] {};}
\newcommand{\annulus}{\tikz[baseline]{\node at (0,0.6ex) [draw, thick, circle, minimum size=1em, inner sep=0pt] {}; \node at (0,0.6ex) [draw, thick, circle, minimum size=0.4em, inner sep=0pt] {}; }}

For a given disconnected world-sheet, each component gives a contributing factor to the amplitude $\mathcal{A}_1$ which we will denote by $\langle\disk\rangle$ for the disk and $\langle\annulus\rangle$ for the annulus. Summing over the number of possible disks $d_1$, annuli $d_2$ etc. together with symmetry factors for exchanging these identical components we obtain an exponentiation \cite{Polchinski:1994fq,Green:1995my}
\begin{equation}
    \begin{split}
        \mathcal{A}_1 &= \int d^{10}y_1 \sum_{d_1 = 0}^\infty \frac{1}{d_1!} \Big( \gs^{-1} \langle \disk \rangle \Big)^{d_1} \sum_{d_2 = 0}^\infty \frac{1}{d_2!} \Big( \frac{1}{2!} \langle \annulus \rangle \Big)^{d_2} \cdots \, \mathcal{A}_1^{(\text{external})} \\[0.5em]
        &= \int d^{10}y_1 \, \exp\!\big(\gs^{-1} \langle \disk \rangle + \frac{1}{2!} \langle \annulus \rangle + \ldots\big) \, \mathcal{A}_1^{(\text{external})} \, ,
    \end{split}
\end{equation}
where $\mathcal{A}_1^{(\text{external})}$ contains the contributions from world-sheet components with punctures for external states. With $\langle \disk \rangle$ negative, we see that this gives the expected non-perturbative contribution of the form $e^{-1/\gs}$.

The leading order contribution to $\mathcal{A}_1^{(\text{external})}$ for the four-graviton scattering amplitude is obtained from four disks with one puncture on each corresponding to an external state. This gives exactly the $\mathcal{R}^4$ combination of the momenta and polarisations as in \eqref{fourtree} \cite{Green:1997tv}. Also including a non-zero axion, Green and Gutperle \cite{Green:1997tv} found the single, unit charged instanton leading order correction as
\begin{equation}
    \label{eq:single-instanton-correction}
    \mathcal{A}_1 = C e^{2\pi i z_0} \int d^{10}y_1 e^{i y_1 \cdot(k_1 + k_2 + k_3 + k_4)} \mathcal{R}^4
        + \ldots 
\end{equation}
where $z_0 = \chi_0 + i \gs^{-1}$, $k_j$ are external momenta and $C$ an unspecified constant. This matches our expectations from the field theory arguments above with corrections of the form $e^{-S_\text{inst}}$ using $m=1$ in \eqref{eq:SUGRA-instanton-action}. 
The contribution \eqref{eq:single-instanton-correction} also agrees with the correpsonding term in the Fourier expansion of the Eisenstein series $\mathcal{E}^{\ordd{10}}_{\gra{0}{0}}(z)$ in \eqref{eq:intro-SL2-fourier-expansion} which is the coefficient function to the $\mathcal{R}^4$ interaction. This, together with the $SL(2, \ints)$-invariance and the perturbative corrections, motivated, in fact, the expression for $\mathcal{E}^{\ordd{10}}_{\gra{0}{0}}(z)$ as an Eisenstein series in \cite{Green:1997tv}.  
    (Note that when computing the effective action one also includes the fluctuations in the dilaton and the axion, and the exponent in the instanton contribution should then become $2 \pi i z$ \cite{Green:1997tv} with the full axio-dilaton $z$ in contrast to the constant background $z_0$ used in the amplitude calculation above.)

We will now find the degeneracy for D-instantons of a certain charge $m$ by studying the T-dual picture. This degeneracy affects how each instanton contribution is weighted in the amplitude, more specifically, it explains the summation over divisors in the instanton measure discussed at the end of section \ref{sec:physical-interp-10D}. 

Using T-duality, where we have compactified the Euclidean time direction in space-time on a circle, the worldlines of D-particles become D-instantons in ten dimensions in the dual limit. The T-dual action for a D-particle of integer Ramond--Ramond charge $n$ whose worldline is wrapping the circle $d$ times can be shown to be \cite{Green:1997tv}
\begin{equation}
    S = 2\pi \abs{nd} e^{-\phi} - 2\pi i nd \chi \, ,
\end{equation}
which is exactly $S_\text{inst}$ in \eqref{eq:SUGRA-instanton-action} with $m = nd$ and $z = z_0$. Thus, the number of D-instanton states with charge $m$ amounts to the number of divisors of $m$.

\chapter[Preliminaries on \texorpdfstring{$p$}{p}-adic and adelic technology]{Preliminaries on \texorpdfstring{$p$}{p}-adic and\\ adelic technology}
\label{ch:p-adic-and-adelic}
As seen in section \ref{sec:adelisation-of-Eisenstein-series} the Fourier expansion of the $SL(2, \reals)$ Eisenstein series factorises into an Euler product over all primes $p$. This number theoretic information is best captured by introducing the $p$-adic numbers which, for any prime $p$, are an extension of the rational numbers, and furthermore the ring of adeles, which encapsulates all the different $p$-adic extensions in a single product.

This chapter is intended as an introduction to these objects as well as providing sample calculations that will be used throughout the remaining text. Additional reading can for example be found in~\cite{Neukirch,Deitmar,Apostol1,Gelbart}.
 Readers familiar with the subject are welcome to proceed to the next chapter and come back to the explicit examples when needed later on in the text. Further reading can be found in \cite{Brekke:1993gf} and \cite{Deitmar}.

For the whole of this chapter and most of the remaining text, let $p$ be a prime number. 

\section{\texorpdfstring{$p$}{p}-adic numbers}
\label{sec:padic}
We start by providing the basic definitions and discussing some of the properties of $p$-adic numbers.
\begin{definition}[Integers $\ints_p$]
    The \emphindex[p-adic@$p$-adic!integers]{$p$-adic integers} $\ints_p$ are formal power series in $p$ with coefficients between $0$ and $p-1$
    \begin{align}
        \label{padicint}
        x \in \ints_p \quad \iff \quad x= x_0 p^0 + x_1 p^1 +\ldots\quad \text{ with } x_i\in \ints/p\ints\cong\left\{0,1,\ldots,p-1\right\}\,.
    \end{align}
The $p$-adic integers form a ring. 
\end{definition}
Arithmetic operations on the $p$-adic integers work in the usual manner. However, since all coefficients in the expansion are positive it may not be immediately obvious how the additive inverse (i.e. subtraction) works. As an example, consider the equation $x+1=0$ that should have a solution over $\ints_p$. The inverse is given by the infinite power series in $p$ with all coefficients are equal to $p-1$:
\begin{align}
x=\sum_{i=0}^\infty(p-1)p^i
\end{align}
This is a bit like evaluating the (non-converging) sum $x=\sum_{k=0}^\infty 10^k=1+10+100+\ldots$ in decimal notation to be an infinite string of $1$s. Multiplying by $9$ and then adding $1$ creates a zero for every decimal place. Hence $9x+1=0\Leftrightarrow x=-1/9= 1/(1-10)$ in agreement with a naive application of the geometric series definition.

Next we define the $p$-adic number field.
\begin{definition}[Number field $\rats_p$]
The associated number field is given by the \emphindex[p-adic@$p$-adic!numbers]{$p$-adic numbers} $\rats_p$ that are formal Laurent series in $p$ with a finite number of terms of degree less than zero, i.e. finite polar part
\begin{equation}
    \label{padicrat}
    x = x_k p^k + x_{k+1} p^{k+1} + \ldots \quad \text{with } x_k\neq 0\,,
\end{equation}
where $k$ is some integer not necessarily positive.
\end{definition}
The $p$-adic numbers $\rats_p$ can be thought of as the completion of rational numbers $\rats$ with respect to the following norm.
\begin{definition}[$p$-adic norm $|\cdot|_p$]
    The \emphindex[p-adic@$p$-adic!norm]{$p$-adic norm} on $\rats_p$ is given by
    \begin{align}\label{pnorm}
        |x|_p = p ^{-k} \quad\Leftrightarrow \quad \text{with } x_k\neq 0 \,.
    \end{align}
\end{definition}

The $p$-adic norm is multiplicative
\begin{align}
|x\cdot y|_p = |x|_p |y|_p
\end{align}
and satisfies a stronger triangle inequality than generic norms, namely 
\begin{align}
|x+y|_p \leq \max (|x|_p, |y|_p)\,,
\end{align}
for $x,y\in\rats_p$.
This second property is called \emphindex[p-adic@$p$-adic!ultrametric property]{ultrametric property} and a space with a norm of this type is called \emphindex{non-archimedean} in contrast with archimedean spaces satisfying the usual archimedean triangle inequality. The $p$-adic norm of $0$ is $|0|_p=0$.

The integer $k$ in~\eqref{pnorm} is called the $p$-adic \emphindex[p-adic@$p$-adic!valuation]{valuation} of $\rats$ or $\rats_p$ and is often also denoted by $\nu_p(x)$. Two properties of the $p$-adic valuation, equivalent to the ones above for the $p$-adic norm, are
\begin{align}
\nu_p(x\cdot y)= \nu_p(x) + \nu_p(y)
\end{align}
and
\begin{align}\label{valprop}
\nu_p(x+y)\geq \min( \nu_p(x), \nu_p(y))\,,
\end{align}
where in the last property equality is achieved if $\nu_p(x)\ne \nu_p(y)$.

The integers in the normed space $\rats_p$ can then be expressed as
\begin{align}
\ints_p = \left\{ x\in\rats_p \,\middle|\, |x|_p \leq 1\right\}\,,
\end{align}
i.e. they have an exponent $k\geq0$ of $p$. This shows that the $p$-adic integers are compactly embedded in $\mathbb{Q}_p$.
The complementary set to $\ints_p$ in $\rats_p$ is given by
 \begin{align}\label{QpZp}
\mathbb{Q}_p\setminus \mathbb{Z}_p=\{x\in \mathbb{Q}_p\, |\, |x|_p> 1\}\,.
 \end{align}
Let us provide two simple examples illustrating the $p$-adic expansion of a rational number.
\begin{example}[$p$-adic expansions]
We consider the $p$-adic expansion of the rational number $x=\frac12\in\rats$ for $p=2$ and $p=3$. 

For $p=2$ one has $|x|_2 = 2^1=2$ or $\nu_2(x)=-1$ and hence $\frac12$ is not a $2$-adic integer. As an element of $\rats_2$ one finds $\frac12 = 1\cdot 2^{-1}$ as the expansion of the form (\ref{padicrat}).

For $p=3$ one has $|x|_3 = 3^0=1$ or $\nu_3(x)=0$ and hence $\frac12$ is a $3$-adic integer. Its expansion of the form (\ref{padicrat}) is $\frac12 = 2\cdot 3^0 +\sum_{k>0} 3^k$.
\end{example}

Another useful property for the $p$-adic norm of the greatest common divisor of two integers which will be used in section \ref{sec:SL3ex} is introduced in the following example.
\begin{example}[Norm of a greatest common divisor]
Let $m$ and $n$ be two integers, $d = \gcd(m, n)$, $m' = m/d$ and $n' = n/d$. Then $1 = \gcd(m', n')$ which, for a prime $p$, means that if $\abs{m'}_p < 1$ (that is, $p \mid m'$) then $\abs{n'}_p = 1$ (that is, $p \nmid n'$) and vice versa. Thus, $1 = \abs{\gcd(m',n')}_p = \max(\abs{m'}_p, \abs{n'}_p)$. Hence,
    \begin{equation}
        \label{eq:p-adic-abs-gcd}
        \abs{d}_p = \abs{d}_p \max(\abs{m'}_p, \abs{n'}_p) = \max(\abs{m'd}_p, \abs{n'd}_p) = \max(\abs{m}_p, \abs{n}_p) \, .
    \end{equation}
\end{example}
We also define the multiplicatively invertible $p$-adic numbers of $\ints_p$ and $\rats_p$.
\begin{definition}[Multiplicatively invertible numbers $\ints_p^\times$ and $\rats_p^\times$]
The set of \emphindex[p-adic@$p$-adic!multiplicatively invertible element]{multiplicatively invertible elements} in $\ints_p$ will be denoted by 
\begin{align}
\ints_p^\times = \left\{ x\in \ints_p \,\middle|\, x^{-1} \text{ exists in $\ints_p$}\right\} = \left\{ x\in \ints_p \,\middle|\, |x|_p=1\right\} = \left\{ x\in \rats_p \,\middle|\, |x|_p=1\right\}.
\end{align}
They correspond to those $x$ in (\ref{padicint}) for which $x_0\neq0$. The set of multiplicatively invertible elements $\rats_p^\times$ in $\rats_p$ is defined as
\begin{align}
\rats^\times_p=\left\{x\in\rats_p\,\middle|\, |x|_p\neq0\right\}\,.
\end{align}
\end{definition}
For $p$-adic numbers the case when $p$ is the prime at infinity, i.e. $p=\infty$, is typically associated with standard calculus via
\begin{align}
\rats_\infty = \reals\,.
\end{align}
In accord with the terminology used for more general number fields the case of a finite prime, i.e. $p<\infty$, is sometimes referred to as the~\emphindex[place!non-archimedean]{non-archimedean place}, while $p=\infty$ is called the~\emphindex[place!archimedean]{archimedean place}.

The $p$-adic numbers were introduced in number theory by Hensel with the intention of transferring the powerful tools of complex analysis to power and Laurent series. A theorem by Ostrowski~\cite{Koblitz1984} states that any non-trivial norm on $\rats$ is either the standard Euclidean norm (leading to the real numbers upon completion) or one of the $p$-adic norms.
\begin{remark}[Alternative construction of $\rats_p$]
Another way of defining the $p$-adic numbers is through the following definition of the $p$-adic norm of an ordinary rational number $x\in\rats$:
\begin{align}
\label{eqn:pnorm-on-Q}
|x|_p = p^{-k}\,,
\end{align}
where $k\in\ints$ is the largest integer such that $x=p^k y$ with $y\in\rats$ not containing any powers of $p$ in its numerator or denominator (in cancelled form); this is often stated as $p^k$ divides $x$. It is from this construction that one obtains $\rats_p$ as the completion of $\rats$ and one obtains an embedding of $\rats$ into $\rats_p$. The definition implies that for a prime  $q$ and $k\in\ints$
\begin{align}
\label{qp}
|q^k|_p = \left\{\begin{array}{cl} p^{-k} &\text{if $p=q$}\\
1 &\text{otherwise}\,.
\end{array}\right.
\end{align}
\end{remark}

\section{\texorpdfstring{$p$}{p}-adic integration}

Integration on $\rats_p$ can be defined with respect to the \emphindex[p-adic@$p$-adic!additive measure]{additive measure} $dx$ that is invariant under translation and has a simple scaling transformation 
\begin{align}
    \label{eq:p-adic-measure}
d(x+a)=dx\,,\quad\quad d(ax) = |a|_p dx\,.
\end{align}
The measure is by convention normalised as to give the $p$-adic integers unit volume:
\begin{align}
\label{padicintmesnorm}
\lint_{\ints_p} dx = 1\,.
\end{align}
We will now provide a series of examples of basic $p$-adic integrals. When evaluating such integrals it is often useful to employ different decompositions of $\ints_p$.
One such decomposition is to write $\ints_p$ as a disjoint union 
\begin{align}\label{ZpCk}
\ints_p = \bigsqcup_{k=0}^{p-1} C_k\,,
\end{align}
where $C_k$ denotes the set of those $p$-adic integers with `constant' coefficient (the coefficient of $p^0$ in~\eqref{padicint}) equal to $k$.
Another decomposition of $\ints_p$ employed is to write it as
\begin{align}
\label{intdec}
\ints_p =\bigsqcup_{k=0}^{\infty} p^k \ints_p^\times\,,
\end{align}
\begin{example}[Volume of invertible integers $\ints_p$]
    \label{ex1}
\begin{align}
\int_{\ints_p^\times} dx = \frac{p-1}{p}.
\end{align}
The integral is a simple consequence of the definition~\eqref{padicintmesnorm} and can be understood intuitively by noting that only $p-1$ out of the $p$ choices for the constant coefficient of $x\in\ints_p$ correspond to elements in $\ints_p^\times$.
For a more formal derivation we use decomposition~\eqref{ZpCk} of $\ints_p$ and integrate over each $C_k$ separately. By translation invariance of the measure~\eqref{eq:p-adic-measure} all $C_k$ have the same volume $1/p$. Integrating over all $C_k$ except for the one with $k=0$ one thus obtains the above formula.
\end{example}
The following two examples explore the integration of the $p$-adic norm.
\begin{example}[Integration of the norm over $\ints_p$]
\label{simpleintegral}
\begin{align}
\label{intxs}
\int_{\ints_p} |x|_p^s dx &= \frac{p-1}{p} \frac1{1-p^{-s-1}}\,,
\end{align}
with convergence for $\Re(s)>-1$. This is derived in a few steps:
\begin{align}
\int_{\ints_p} |x|_p^s dx &= \sum_{k=0}^{\infty}\, \int_{p^k\ints_p^\times} |x|_p^s dx 
= \sum_{k=0}^\infty p^{-ks} \int_{p^k\ints_p^\times}  dx
= \sum_{k=0}^\infty p^{-ks} \int_{\ints_p^\times} p^{-k} dy\nn\\
&= \frac{p-1}{p}\sum_{k=0} p^{-k(s+1)} = \frac{p-1}{p} \frac1{1-p^{-s-1}}.
\end{align}
In the first step we have used the decomposition~\eqref{intdec} of the $p$-adic integers.
Then we have used the fact that for $x\in p^k\ints_p^\times$ the norm is $|x|_p=p^{-k}$. After that we have changed variables to $x=p^k y$ with $y\in\ints_p^\times$, used the resulting volume of $\ints_p^\times$ computed in exmaple~\ref{ex1} and carried out the geometric sum. 
\end{example}

Using the identity from the previous example we can also evaluate the following integral which will be used in chapter~\ref{ch:SL2-fourier} and section~\ref{pfin}. 
\begin{example}[Integration of the norm over $\rats_p\setminus\ints_p$]
    \label{ex:p-adic-x-s-integral}
 \begin{align}
 \int_{\mathbb{Q}_p\setminus \mathbb{Z}_p} |x|_p^{s} dx= \frac{p-1}{p} \frac{p^{s+1}}{1-p^{s+1}}\,,
 \label{nonintegerintegral}
 \end{align}
 with $\Re(s)<-1$ and the domain of integration as defined in~\eqref{QpZp}. The integral is then evaluated to be
 \begin{align}
 \int_{\mathbb{Q}_p\setminus\mathbb{Z}_p} |x|_p^{s} dx&= \int_{|x|_p>1}  |x|_p^{s} dx=\sum_{k=1}^{\infty} p^{ks} \int_{p^{-k}\ints_p^\times} dx = \sum_{k=1}^\infty p^{k(s+1)} \int_{\ints_p^\times} dx\nonumber\\
&= \frac{p-1}{p} \sum_{k=1}^{\infty}  p^{k(s+1)} = \frac{p-1}{p} \frac{p^{s+1}}{1-p^{s+1}}.
 \label{usefulintegral}
 \end{align}
The integral converges for $\Re(s)<-1$. Note that the same integral over all of $\rats_p$ does not exist.
\end{example}

\begin{remark}[Multiplicative measure $dx^\times$]
\label{rk-multmeasure}
We denote the \emphindex[p-adic@$p$-adic!multiplicative measure]{multiplicative measure} on $\rats_p^\times$ by $d^\times x$ with its defining relation
\begin{align}
d^\times x = \frac{p}{p-1} \frac{dx}{|x|_p}\,.
\label{multiplicativemeasure}
\end{align}
It satisfies $\int_{\ints_p^\times} d^\times x =1$. It transforms as $d^\times(ax) = d^\times x$. Integrating the function $|x|_p^s$ against the multiplicative measure $d^{\times}x$ the result~\eqref{intxs} simplifies to
\begin{align}
 \int_{\mathbb{Z}_p} |x|_p^{s}\,d^{\times} x= \sum_{k=0}^{\infty}p^{-ks}\int_{\mathbb{Z}^{\times}_p}d^\times x=\sum_{k=0}^{\infty}p^{-ks}=\frac{1}{1-p^{-s}}\, ,
 \end{align}
where in the first step we  used the property~\eqref{intdec}. Note that the same result is obtained if we restrict the integration domain to $\mathbb{Z}_p\backslash \{0\}$, which will be useful in the proof of proposition~\ref{prop-Tateglobal}.  
\end{remark}

\section{Characters and the Fourier transform}
\label{sec:chars}
In this section we introduce the concept of a character which is then used to define the $p$-adic Fourier transform. As before we provide explicit computations of various integrals serving as prototypical examples for later calculations. 
\begin{definition}[Fractional part of a $p$-adic number]
    \label{def:fracpart}
The \emphindex[p-adic@$p$-adic!number, fractional part]{fractional part} $[y]_p$ of a $p$-adic number $y\in\rats_p$ is given by its class in $\rats_p/\ints_p$, or more concretely by the terms in its series expansion with negative powers of $p$:
\begin{align}
\label{fracpart}
\left[ x_{k} p^{k}+\ldots+ x_{-1}p^{-1} +x_0 p^0 +x_1 p^1 +\ldots \right]_p 
= \left\{ \begin{array}{cl} x_{k} p^{k} +\ldots x_{-1}p^{-1} & \text{if $k<0$,}\\ 0 &\text{otherwise}.\end{array}
\right.
\end{align}
Note that we will often suppress the subscript $p$ when there is no risk of confusion.
\end{definition}
We will now show that given a rational number $x$, subtracting all the fractional parts of $x$ with respect to all $\ints_p$ from $x$ leaves a normal integer. This will, for example, be used in sections~\ref{sec_adeles} and \ref{sec:SL3ex}.
\begin{proposition}[Normal integer of a $p$-adic number]
    \label{prop:p-adic-fractional}
    Let $x \in \rats$. Then,
    \begin{equation}
        \label{eq:p-adic-fractional-property}
        x - \sum_{p<\infty} [x]_p \in \ints \, .
    \end{equation}
\end{proposition}
\begin{proof}
    By design, $x - [x]_p \in \ints_p$ and for any prime $q \neq p$
    \begin{equation}
        \abs{[x]_q}_p \leq \max(\abs{x_k q^k}_p, \ldots, \abs{x_{-1}q^{-1}}_p) \leq 1
    \end{equation}
    if $k < 0$ (since $x_i \in \ints$) and otherwise $[x]_q = 0$ which means that $[x]_q \in \ints_p$. Hence, for any prime $p$
    \begin{equation}
        x - \sum_{q<\infty} [x]_q = (x - [x]_p) - \sum_{q\neq p} [x]_q \in \ints_p
    \end{equation}
    which proves the statement. 
\end{proof}
With the definition of the fractional part of a $p$-adic number introduced we can now provide the definition of an additive character.
\begin{definition}[Additive characters]
\emphindex[p-adic@$p$-adic!additive character]{Additive characters}  on $\rats_p$ are defined by
\begin{align}
\label{eq:p-adic-char}
\psi_p\equiv\psi_{p,u} \,:\, \rats_p \to U(1)\,,\quad
\psi_{p,u}(x) = e^{-2\pi i [ux]_p} \quad\quad x,u\in\rats_p\,.
\end{align}
The additive characters of \eqref{eq:p-adic-char} satisfy the relations $\psi_{p,u}(x)\psi_{p,u}(y)=\psi_{p,u}(x+y)$ and $\psi_{p,u}\psi_{p,v}=\psi_{p,u+v}$, as well as $\overline{\psi_{p,u}(x)}=\psi_{p,-u}(x)=\psi_{p,u}(-x)$. The \emphindex[p-adic@$p$-adic!conductor]{conductor} of the character is its kernel $|u|_p\ints_p$, but often we simply call $|u|_p$ the conductor. 
\end{definition}

Note that in the following, we shorten the notation to $\psi_p\equiv \psi_{p,u}$ since it will be more important to keep track of the dependence on the prime $p$; the `mode number' $u$ will be given explicitly where needed. Also in the interest of simplicity of notation we will  often drop the prime $p$ script on the symbol for the fractional part $[\cdot]_p$ when writing out characters explicitly.

\begin{example}[Integration of a character over $p^k\ints_p$]
\label{intchar}
For $k\in\ints$ one has
\begin{align}
\label{charintpk}
\int_{p^k\ints_p} e^{-2\pi i [ux]} dx 
=p^{-k} \gamma_p(up^k)\,,
\end{align}
where the characteristic function $\gamma_p(u)$ of $\ints_p$ in $\rats_p$ is defined as
\begin{align}
\label{charintp}
\gamma_p(u):=\int_{\ints_p} e^{-2\pi i [ux]} dx = \left\{ \begin{array}{cl} 1&\text{if $u\in\ints_p$}\,,\\0&\text{otherwise}\,.
\end{array}\right.
\end{align}
The function $\gamma_p(u)$ is also called the $p$-adic \emphindex[p-adic@$p$-adic!Gaussian]{Gaussian} which will be discussed in more detail in section~\ref{sec:p-adic-special-functions}.

In order to derive this result we start with the case when $k=0$:
\begin{align}
\int_{\ints_p} \psi_u(x) dx =  \int_{\ints_p} e^{-2\pi i [ux]} dx
\end{align}
and the integral only depends on the conductor $|u|_p$. We distinguish two cases: (i) $u\in\ints_p$ and (ii) $u\notin\ints_p$:

(i) If $u\in\ints_p$ then $[ux]_p=0$ for $x\in\ints_p$ and hence the integral equals $\int_{\ints_p}dx=1$.

(ii) If $u\notin\ints_p$ then we are effectively integrating a periodic function over a full period and hence the integral gives zero. More concretely, consider the example when $u=p^{-1}$; then 
\begin{align}
\int_{\ints_p} e^{-2\pi i [p^{-1} x]} dx = \sum_{k=0}^{p-1}\,e^{-2\pi i k/p}\, \int_{C_k} dx = \frac1p \sum_{k=0}^{p-1} e^{-2\pi i k/p} = 0
\end{align}
with $C_k$ defined as in~\eqref{ZpCk} and where we have used the fact that $\int_{C_k}dx=1/p$, c.f. also example~\ref{ex1}.
If $u$ is `more rational' one has to refine the summation region more but will always encounter sums that average to zero.
We have thus derived~\eqref{charintpk} for the case of $k=0$.

The result for the integral in the case when $k\neq0$ then follows by a simple change of variables:
\begin{align}
\label{charintpk2}
\int_{p^k\ints_p} e^{-2\pi i [ux]} dx 
= p^{-k}\int_{\ints_p} e^{-2\pi i [up^kx]} dx 
=p^{-k} \gamma_p(up^k).
\end{align}
\end{example}
\newpage

We will also require the integral over `shells' of $p$-adic numbers.

\begin{example}[Integration of a character over $p^k\ints_p^\times$]
\label{intcharunits}
For $k\in\ints$ we have
\begin{align}
\lint_{p^k\ints_p^\times} e^{-2\pi i [ux]} dx  
& = \left\{\begin{array}{cl}
\frac{p-1}{p}p^{-k} & \text{for $|u|_p\leq p^k$}\\
-p^{-(k+1)} & \text{for $|u|_p=p^{k+1}$}\\
0 & \text{for $|u|_p>p^{k+1}$}
\end{array}\right..
\end{align}
Starting as before with the case $k=0$, this can be related to the preceding example by noting that $\ints_p^\times=\ints_p\setminus(p\ints_p)$:
\begin{align}
\label{charintpstar}
\lint_{\ints_p^\times} e^{-2\pi i [ux]} dx &= \lint_{\ints_p} e^{-2\pi i[ux]} dx -\lint_{p\ints_p} e^{-2\pi i [ux]} dx\nn\\
&= \gamma_p(u) -p^{-1} \lint_{\ints_p} e^{-2\pi i[upx]}dx \nn\\
& = \gamma_p(u) -p^{-1} \gamma_p (pu)\nn\\
& = \left\{\begin{array}{cl}
\frac{p-1}{p} & \text{for $|u|_p\leq 1$, i.e., $u\in\ints_p$}\\
-p^{-1} & \text{for $|u|_p=p$}\\
0 & \text{for $|u|_p>p$}
\end{array}\right..
\end{align}
The result of $p^k\ints_p^\times$ for $k\neq0$ then follows by a change of variables
\begin{align}
\lint_{p^k\ints_p^\times} e^{-2\pi i [ux]} dx  
&=p^{-k}\lint_{\ints_p^\times} e^{-2\pi i [up^kx]} dx = p^{-k}\gamma_p(up^k)-p^{-(k+1)}\gamma_p(up^{k+1})\nn\\
& = \left\{\begin{array}{cl}
\frac{p-1}{p}p^{-k} & \text{for $|u|_p\leq p^k$}\\
-p^{-(k+1)} & \text{for $|u|_p=p^{k+1}$}\\
0 & \text{for $|u|_p>p^{k+1}$}
\end{array}\right..
\end{align}
Note also that this implies that
\begin{align}
\lint_{\rats_p\backslash\ints_p} e^{-2\pi i[ux]}dx = \sum_{k=-1}^{-\infty} \lint_{p^{k}\ints_p^\times} e^{-2\pi i[ux]}dx = -\gamma_p(u)\,.
\end{align}
\end{example}
An important comment here concerns the integral of a character over all of $\rats_p$: Since $\rats_p$ is formally the sum of $p^k\ints_p^\times$ over all $k\in\ints$, we see from the above result that for any $u\in\rats_p$ we obtain formally
\begin{align}
\label{wronglimit}
\lint_{\rats_p} e^{-2\pi i [ux]}d x = 0\quad\quad\text{(not well-defined!)}
\end{align}
which is the analogue of the incorrect equation $\int_\reals e^{2\pi i u x}dx=0$ which could be derived by splitting up $\reals$ into an infinite number of  intervals of length $1/u$ on each of which the integral vanishes. As is well-known, the integral over the whole real line of $e^{2\pi i u x}$ is not well-defined but rather yields a $\delta$-distribution. We will now see that something similar is true for the $p$-adic character $e^{-2\pi i [ux]}$ integrated over $\rats_p$.

Before introducing the concept of a $p$-adic Fourier transform let us make a short comment about function spaces used. The functions which we will be integrating are elements of $\mathcal S(\rats_p)$ which is the \emph{Schwartz-Bruhat} space. These functions generalise the \emph{Schwartz functions} which are infinitely differentiable, with rapidly decreasing derivatives. 

\begin{definition}[Fourier transform]
One defines the \emphindex[p-adic@$p$-adic!Fourier transform]{Fourier transform} over $\rats_p$ by integrating a function $f_p$ on $\rats_p$ against the additive character $\psi_p(x)\equiv\psi_{p,u}(x)$:
\begin{align}
\tilde{f}_p(u) =\lint_{\rats_p} f_p(x) \psi_p(x) dx=\lint_{\rats_p} f_p(x) e^{-2\pi i[ux]}dx\,.
\end{align}
The inverse transform uses the conjugate character
\begin{align}
f_p(x) = \lint_{\rats_p} \tilde{f}_p(u)\overline{\psi_p(x)} du = \lint_{\rats_p} \tilde{f}_p(u) e^{2\pi i[ux]}du\,.
\end{align}
\end{definition}

One can now ask for which functions $f_p$ the transform is well-defined and can actually be inverted. As a first step we calculate the composition of the transforms of the characteristic function of a ball $p^k \ints_p\subset\rats_p$, i.e., $f_p(x)=\gamma_p(p^{-k}x)$
\begin{align}
&\lint_{\rats_p} \overline{\psi_{p}(x)} \lint_{\rats_p}\psi_p(y) \gamma_p (p^{-k}y) dy du
=\lint_{\rats_p} \overline{\psi_{p}(x)} \lint_{p^k\ints_p}\psi_p(y)  dy du\nn\\
&\quad= \lint_{\rats_p} e^{2\pi i[ux]} p^{-k} \gamma_p(up^k)du 
= \lint_{\ints_p} e^{2\pi i [up^{-k}x]}du = \gamma_p(p^{-k}x)\,.
\end{align}
From this calculation we see that restricting to compactly supported (and bounded) functions makes the integrals well-defined. We can relax the assumption of compact support if the function decreases sufficiently fast for larger and larger balls. This is for instance the case when $f_p(x)=|x|_p^s$ with $\Re(s)$ sufficiently negative. However, since in this case $|x|_p^s$ blows up for $|x|_p\to 0$ one has to cut out that region or replace $f_p(x)$ by a different function there. In summary, the $p$-adic Fourier transform is only well-defined on functions that are locally constant (i.e. constant on each $p^k\ints_p^\times$) and have compact support or a sufficiently fast decrease when $|x|_p\to\infty$. 

Let us go through a series of interesting examples.
\begin{example}[Fourier transform of $|x|_p^s$ over $\rats_p\setminus\ints_p$]
\label{BesselFT}
Here, we cut out the compact region of the integers and consider the effect of the damping function $|x|_p^s$ with $\Re(s)<-1$. The result is
\begin{align}
\lint_{\rats_p\backslash\ints_p} |x|_p^s \psi_p(x) d x 
= \gamma_p(u) \left((1-p^{s})\frac{1-p^{s+1} |u|_p^{-s-1}}{1-p^{s+1}}-1\right).
\end{align}
To show this, we denote the integral by $I$ and distinguish two cases: $(i)$ $u$ integral and $(ii)$ $u$ non-integral:

$(i)$: If $u\in \ints_p$ and has conductor $p^k$ with $k\geq 0$, then we evaluate the integral as
\begin{align}
I &=  \int_{\rats_p\backslash\ints_p} |x|_p^s e^{-2\pi i [p^k x]} dx
=  \sum_{\ell=1}^\infty p^{s\ell} \int_{p^{-\ell}\ints_p^\times} e^{-2\pi i [p^k x]} dx\nonumber\\
&= \sum_{\ell=1}^\infty p^{(s+1)\ell} \int_{\ints_p^\times} e^{-2\pi i [p^{k-\ell} x]} d x
\underset{\text{Ex. \ref{intcharunits}}}=  \frac{p-1}{p} \sum_{\ell=1}^{k} p^{(s+1)\ell} -\frac1p\, p^{(k+1)(s+1)}\nonumber\\
&= (1-p^{s}) \frac{1-p^{s+1} |u|_p^{-s-1}}{1-p^{s+1}} -1
\end{align}

$(ii)$: If $u\notin \ints_p$, so that the character has conductor $p^k$ with $k<0$ we find\begin{align}
I=  \int_{\rats_p\backslash\ints_p} |x|_p^s e^{-2\pi i [p^k x]} d x = \sum_{\ell=1}^\infty p^{(s+1)\ell} \int_{\ints_p^\times} e^{-2\pi i [p^{k-\ell}x]} dx=0
\end{align}
by example~\ref{intcharunits} since $k-\ell<-1$ for all $\ell\geq1$ and $k<0$.
\end{example}

The result for the Fourier integral over $\ints_p$ is slightly more complicated.

\begin{example}[Fourier transform of $|x|_p^s$ over $\ints_p$]
Here we cut out the region $\rats_p\backslash\ints_p$. In order to have a bounded function for $|x|_p\to 0$ we now require $\Re(s)>-1$. The Fourier transform now evaluates to
\begin{align}
\int_{\ints_p} |x|^s_p e^{-2\pi i[x u]} d x
&= \gamma_p(u) \frac{p-1}{p} \frac{1}{1-p^{-s-1}} +(1-\gamma_p(u)) |u|_p^{-s-1}\frac{1-p^s}{1-p^{-s-1}}\,.
\end{align}
This can be derived in a few steps
\begin{align}
\int_{\ints_p} |x|^s_p e^{-2\pi i[x u]} d x
&= \sum_{\ell=0}^\infty p^{-s \ell} \int_{p^\ell \ints_p^\times} e^{-2\pi i [ux]} d x
= \sum_{\ell=0}^\infty p^{-(s+1)\ell} \int_{\ints_p^\times} e^{-2\pi i [u p^\ell x]} d x\nn\\
&= \gamma_p(u) \frac{p-1}{p} \frac{1}{1-p^{-s-1}} +(1-\gamma_p(u)) |u|_p^{-s-1}\frac{1-p^s}{1-p^{-s-1}},
\end{align}
where we have treated the cases $u\in\ints_p$ and $u\notin\ints_p$ separately and used equation~(\ref{intxs}) and example~\ref{intcharunits}.
\end{example}

\section{\texorpdfstring{$p$-adic}{p-adic} Gaussian and Bessel function}
\label{sec:p-adic-special-functions}
In this section we will discuss two special functions: the \emphindex[p-adic@$p$-adic!Gaussian]{$p$-adic Gaussian} and the \emphindex[p-adic@$p$-adic!Bessel function]{$p$-adic Bessel function}, which will play a role later on in the text.

The $p$-adic analogue of the Gaussian $e^{-\pi x^2}$ is given by the function
\begin{align}
\label{padicGauss}
    \gamma_p(x) =
    \begin{cases}
        1 & \text{if } x \in \ints_p, \text{ i.e. } |x|_p\leq 1 \\
        0 & \text{if } x \notin \ints_p, \text{ i.e. } |x|_p>1\,,
    \end{cases}
\end{align}
which we have already encountered in example~\ref{intchar} of the previous section. In order to see why it is the generalisation of the real Gaussian $e^{-\pi x^2}$ we recall that the real Gaussian is invariant under Fourier transformation.
Using~\eqref{charintpk} this property is then also easily checked for the $p$-adic version:
\begin{align}\label{GaussFourier}
\tilde{\gamma}_p(u) &= \int_{\rats_p} \psi_u(x) \gamma_p (x) dx = \int_{\ints_p} e^{-2\pi i [ux]} d x = \gamma_p(u)\,.
\end{align}
Let us also note the following useful property of the finite product of the $p$-adic Gaussian which will be used in chapter~\ref{ch:SL2-fourier} in the computation of the Fourier coefficients of Eisenstein series on $SL(2,\mathbb A)$, c.f. equation~\eqref{finitSL2coeff}:
\begin{equation}
    \label{eq:prod-gamma-p}
    \prod_{p < \infty} \gamma_p(m) =
    \begin{cases}
        1 & \text{if } m \in \ints \\
        0 & \text{otherwise}\,.
    \end{cases}
\end{equation}

In order to introduce the $p$-adic version of the Bessel function, recall that the real (modified) Bessel function $K_s$ can be written as the (inverse) Fourier transform of the function 
$||(1,u)||^{-2s}=(1+u^2)^{-s}$ via
\begin{align}
\label{Besselreal}
\int_\reals (1+u^2)^{-s} e^{-2\pi i m u} du = \frac{2\pi^s}{\Gamma(s)} |m|^{s-1/2} K_{s-1/2}\left(2\pi |m|\right),
\end{align}
where $\Gamma(s)$ is the standard Gamma function.
The $p$-adic generalisation of the integrand is through $||(1,u)||_p^{-2s}= (\max(1,|u|_p))^{-2s}$. The normalisation to be chosen is~\cite{Kazhdan:2001nx,KazhdanPolishchuk}
\begin{align}
\label{psphereU}
\tilde{f}(u) = \frac1{1-p^{-2s}} (\max(1,|u|_p))^{-2s}\,.
\end{align}
The Fourier transform of this function is
\begin{align}
f(x) = \gamma_p(x)  \frac{1-p^{-2s+1} |x|_p^{2s-1}}{1-p^{-2s+1}}\,.
\end{align}
which we call the $p$-adic Bessel function and will be used, for example, in section \ref{sec:SL2FC}. 

To demonstrate these properties of the $p$-adic Bessel function we perform the Fourier transform in the following example.
\begin{example}[Fourier transform of $p$-adic Bessel function]
\label{besselfull} Consider the following calculation that converges for $\Re(s)>1/2$
\begin{align}
&\quad\quad\int_{\rats_p} e^{2\pi i [ux]} \frac1{1-p^{-2s}} (\max(1,|u|_p))^{-2s} du\nn\\
&= \frac1{1-p^{-2s}}\int_{\ints_p} e^{2\pi i [ux]} du + \frac{1}{1-p^{-2s}} \int_{\rats_p\setminus\ints_p} |u|_p^{-2s} e^{2\pi i [ux]}d u\,.
\end{align}
We have separated the integral according to the two possible cases of the $\max$ function. The first integral is given in example~\ref{intchar} and the second one in example~\ref{BesselFT}. 
\pagebreak

Combining the results we obtain
\begin{align}
\quad\quad\int_{\rats_p} e^{2\pi i [ux]} \frac1{1-p^{-2s}} (\max(1,|u|_p))^{-2s} du
= \gamma_p(x) \frac{1-p^{-2s+1} |x|_p^{2s-1}}{1-p^{-2s+1}}\,.
\end{align}
\end{example}

\section{Adeles}
\label{sec_adeles}
In the previous sections we have introduced the concept of the $p$-adic completions $\rats_p$ of $\rats$ and we have shown in a number of examples how integration can be carried out locally and also that the real Gaussian and Bessel function have $p$-adic counterparts. The next step will be to organise the completions $\mathbb{Q}_p$ of the rational numbers $\mathbb{Q}$ into a global field called the adeles of $\rats$, denoted by $\mathbb A$, which comprises the $p$-adic completions at all primes, including the prime at infinity, at the same time. The introduction of the adeles as a global number field is in line with the so-called local-to-global principle. For a brief summary highlighting the power of this principle see appendix~\ref{appendixLocGlob}.

\begin{definition}[Adeles $\mathbb A$]
The \emphindex{adeles} $\ads = \ads_\rats$ of $\rats$ are defined as a \emph{restricted direct product}
\beq
\label{adeleseq}
\mathbb{A} = \mathbb{R}\times \prodp_{p<\infty}\mathbb{Q}_p\,,
\eeq
where the restriction on the product (signified by the prime) means that $\mathbb{A}$ consists of those elements 
\beq
a=(a_p)=(a_\infty; a_2, a_3, a_5, a_7, \dots )
\eeq
 such that for almost all finite primes $p$ one has $a_p\in \mathbb{Z}_p$.
\end{definition}

The restriction on the direct product in the definition of the adeles makes them locally compact which is needed for the existence of a Haar measure.
Also as a consequence of the definition, the adeles are endowed with a natural topology, and they are in fact a locally compact ring. We refer the reader to~\cite{Gelbart} for more details on these issues. It will sometimes be useful to talk about the \emphindex[adeles!finite]{finite} adeles $\mathbb{A}_f$ which are defined as the restricted direct product over the finite primes:
\beq
\mathbb{A}_f=\prodp_{p<\infty} \mathbb{Q}_p,.
\eeq
We also define the set of invertible elements of the adeles.
\begin{definition}[Ideles $\mathbb A^\times$]
The~\emphindex{ideles} $\mathbb A^\times$ are the set of invertible elements in $\mathbb{A}$. They are defined as:
\beq
\mathbb{A}^{\times}=\mathbb{R}^{\times} \times \prodp_{p<\infty}\mathbb{Q}^{\times}_p.
\eeq
\end{definition}
The norm for the adeles is induced directly from the local norms.
\begin{definition}[Global norm $|\cdot |_\mathbb{A}$]
The~\emphindex[adeles!norm]{global norm} $|\cdot |_\mathbb{A}$ on $\mathbb{A}$ is induced from the norm $|\cdot|_p$ on the local factors $\mathbb{Q}_p$ according to the formula
\beq
|a|_\mathbb{A} = \prod_{p\leq \infty} |a_p|_p.
\eeq
This is in fact a convergent product since almost all $a_p\in \mathbb{Z}_p$ and hence satisfy $|a_p|_p\leq 1$. Furthermore, the adelic norm  $|a|_\ads$ is only non-zero if $a\in \ads^\times$ and then is a finite product.
\end{definition}  
The \emph{strong approximation principle} states that the set 
\beq
J=\mathbb{R}_+\times \prod_{p<\infty}\mathbb{Z}_p^\times
\label{funddomain}
\eeq
is a fundamental domain for $\mathbb{Q}^{\times}\backslash \mathbb{A}^{\times}$. Hence we can write the ideles as the (disjoint) union
\beq
\mathbb{A}^{\times} = \bigcup_{k\in \mathbb{Q}^\times}k\cdot J.
\label{StrongapproxAbelian}
\eeq 
(A higher rank version of strong approximation is proven in section~\ref{sec_adelisation}.)

The rational numbers \emphindex[adeles!diagonal embedding of $\rats$]{embed diagonally} into the adeles, i.e. $\mathbb{Q}\hookrightarrow \mathbb{A}$, by simply taking  
\beq
\mathbb{Q}\,  \ni\,  x \, \longmapsto \, (x;x, x, x, \dots ) \, \in \, \mathbb{A}\,.
\eeq
One can see that this is indeed and element of the adeles since for $x \in \mathbb{Q}$ the norm $|x|_p$ is non-trivial only for the finite number of $p$'s which divide $x$ in the sense of \eqref{eqn:pnorm-on-Q}. In other words, the prime factorisations of the coprime numerator and denominator of $x$ contain only a finite number of primes. By factorising $x \in \rats$ into its prime factors we see that 
\begin{equation}
    \label{eq:adelic-norm-rationals}
    \abs{x}_\ads = \abs{x}_\infty \prod_{p<\infty} \abs{x}_p = \abs{x}_\infty \abs{x}_\infty^{-1} = 1 \, .
\end{equation}

Following~\cite{Deitmar}, we will now show that with this embedding $\mathbb{Q}$ sits discretely inside $\mathbb{A}$, mimicking the way the integers $\mathbb{Z}$ are embedded as a lattice inside $\mathbb{R}$. As we will see, this fact lies at the heart of the analysis in subsequent sections. 

\begin{proposition}[Discrete embedding of $\rats$ in $\mathbb A$]
\label{prop:discrete-embedding}
$\mathbb{Q}$ sits discretely inside $\mathbb{A}$.
\begin{proof}
Let us first consider $0 \in \mathbb{Q}$ and construct
\begin{equation}
    V = \left(-\frac{1}{2}, \frac{1}{2} \right) \times \prod_{p < \infty}  \mathbb{Z}_p \subset \mathbb{A} \, .
\end{equation}

The subgroup $\mathbb{Z}_p$ is an open ball in $\mathbb{Q}_p$ since $|x|_p$ takes only a discrete set of values, that is, $\mathbb{Z}_p = \{ x \in \mathbb{Q}_p \mid |x|_p \leq 1 \} = \{ x \in \mathbb{Q}_p \mid |x|_p < \alpha \}$ for any $1 < \alpha < p$.

Thus, $V$ is an open neighbourhood of $0$ in $\mathbb{A}$ and for any $x \in V \cap \mathbb{Q}$ we have that $|x|_p \leq 1$ for all $p < \infty$ which means that $x \in \mathbb{Z}$, and $|x|_\infty < \tfrac{1}{2}$ which then gives that $x = 0$. Hence, we have found an open neighbourhood $V$ to $0$ in $\mathbb{A}$ such that $V \cap \mathbb{Q} = \{0\}$. For a general point $r \in \rats$ these arguments generalise by instead considering $r + V$, which makes $\mathbb{Q}$ discrete in $\mathbb{A}$.
\end{proof}
\end{proposition}

With the definition of the adeles as the collection of all local factors at hand, we will now see how to turn a set of local functions into a global one. 

\section{Adelisation}
\label{sec:adelisation}
One can extend a collection of local functions $f_p$ on $\mathbb{Q}_p$ to a global function $f_\mathbb{A}$ on $\mathbb{A}$:
\beq
f_\mathbb{A}(a) = f_\mathbb{A}(a_\infty; a_f)\,, \qquad a_f=(a_2,a_3,a_5, \dots ) \in \mathbb{A}_f\,,
\eeq
via an Euler product
\beq f_\mathbb{A}(a) = \prod_{p\leq \infty} f_p(a_p)\,.
\eeq
Starting from $f_\mathbb{A}$ we can recover a function on $\mathbb{R}$ by setting
\beq 
f_\infty(a_\infty)= f_\mathbb{A}(a_\infty; 1, 1, 1, \dots )\,.
\eeq
One says that $f_\mathbb{A}$ is the \emphindex[adeles!adelisation]{adelisation} of $f_{\mathbb{R}}$. Similarly we can extend to $\mathbb{A}$ the notion of local additive characters $\psi$ on $\mathbb{Q}_p$.

Let $u = (u_\infty, u_2, u_3, \ldots) \in \ads$ and $\psi_p : \mathbb{Q}_p \to U(1)$ be an additive character, such that for finite $p$ this 
coincides with the character
\begin{equation}
    \psi_{p}(x_p) = e^{-2 \pi i [u_p x_p]_p}, \qquad u_p, x_p \in \rats_p \, ,
\end{equation}
defined in section \ref{sec:chars}, while for $p=\infty$ this is the standard character on $\reals$:

\begin{equation}
    \psi_\infty(x_\infty)=e^{2 \pi i u_\infty x_\infty}, \qquad u_\infty, x_\infty \in \mathbb{R}.
\end{equation}
We can then consider a global character 
\begin{equation}
    \psi_\mathbb{A} \, :\, \mathbb{A} \, \to \, U(1)
\end{equation}
as the adelisation of $\psi_\mathbb{R}$, i.e., as the Euler product
\begin{equation}
\label{Ach_fact}
    \psi_\ads(x) = \prod_{p\leq\infty} \psi_p(x_p) =  e^{2\pi i u_\infty x_\infty} \prod_{p < \infty} e^{-2 \pi i [u_p x_p]_p }  \, , 
\end{equation}
which we will denote as $\psi_\ads(x) = e^{2 \pi i u x}$ for short.

The sign difference in the exponentials of the characters at the archimedean and non-archimedean places have been introduced for the following reason. For $u = m \in \rats$ diagonally in $\ads$, the character $\psi_\ads$ is periodic in $\rats$ since for $x \in \ads$ and $r \in \rats$
\begin{equation}
    \psi_\ads(x + r) = \psi_\ads(x) \psi_\ads(r)
\end{equation}
with
\begin{equation}
    \psi_\ads(r) = \prod_{p \leq \infty} \psi_p(r) = \exp\Big( 2\pi i \big( m r  - \sum_{p < \infty} [m r]_p  \big) \Big) = 1
\end{equation}
using proposition~\ref{prop:p-adic-fractional}. Thus, for rational $u$, $\psi_\ads$ is a character on $\rats \bs \ads$. That these are all the characters on $\rats \bs \ads$ is shown in \cite{Deitmar}.

Integration over the adeles is similarly defined using Euler products~\cite{LanglandsEP}. For instance, the integral over an adelic function $f_\mathbb{A}(x)$ 
can be written as
\beq
\int_{\mathbb{A}} f_\mathbb{A}(x)\, dx= \left(\int_\mathbb{R} f_\mathbb{R}(x)\, dx\right)\,\left( \prod_{p< \infty} \int_{\mathbb{Q}_p} f_p(x)\, dx\right).
\eeq

\begin{definition}[Adelic Fourier transform]
The adelic Fourier transform is defined using the global character $\psi_{\mathbb{A}}$ as follows:
\beq
\tilde{f}_{\mathbb{A}}(u)=\int_{\mathbb{A}} f_{\mathbb{A}}(x)\overline{\psi}_{\mathbb{A}}(x) dx.
\eeq
\end{definition}

We will perform several integrals of this type in subsequent sections.

In the following section we will illustrate the usefulness of the adelic framework in the 
context of the Riemann zeta function.

\section{Adelic analysis of the Riemann zeta function}
\label{RiemannZetaEx}

In this section we will illustrate the power of the adelic formalism by analysing the Riemann zeta function from this point of view. This was one of the main points of the celebrated thesis of Tate \cite{Tate} which first introduced the notion of Fourier analysis over the adeles. 

\subsection{The completed Riemann zeta function} 
\label{sec-complRiem}

The first task will be to illustrate how the \emphindex[Riemann zeta function!completed]{completed Riemann zeta function} is a much more  natural object from an adelic perspective, than the ordinary zeta function. Recall first that the completed Riemann zeta function takes the form:
\begin{align}
\xi(s)=\pi^{-s/2} \Gamma(s/2) \zeta(s).
\end{align}
We now have:
\begin{proposition}[Tate's global Riemann integral \cite{Tate}]
\label{prop-Tateglobal}
The completed Riemann zeta function $\xi(s)$ can be written in the following global form: 
\beq
\xi(s)=\int_{\mathbb{A}^{\times}}\gamma_\mathbb{A}(x) |x|_\mathbb{A}^{s} d^{\times}x,
\eeq
where $s\in \mathbb{C}$ and $\gamma_\mathbb{A}=\prod_{p\leq \infty} \gamma_p$ with $\gamma_p$ the $p$-adic Gaussian (\ref{padicGauss}) and $\gamma_\infty=e^{-\pi x^2}$.
\end{proposition}
\begin{proof}
Splitting the integral into an Euler product yields
\beq
\int_{\mathbb{A}^{\times}}\gamma_\mathbb{A}(x) |x|_\mathbb{A}^{s} d^{\times}x=\left(\int_{\mathbb{R}^{\times}} e^{-\pi x ^2} |x|^{s}_\infty d^{\times} x\right)\, \prod_{p< \infty} \int_{\mathbb{Q}^{\times}_p}\gamma_p(x) |x|_p^{s}d^{\times}x.
\eeq
The archimedean integral can be evaluated in terms of a Gamma function:
\begin{align}
\int_{\mathbb{R}^{\times}}e^{-\pi x^2} |x|_\infty^{s} d^{\times}x = \int_{\mathbb{R}} e^{-\pi x^2} |x|_\infty^{s-1} dx = \pi^{-s/2}\Gamma(s/2),
\end{align}
where we made use of \eqref{multiplicativemeasure}. 

Due to the $\gamma_p$-factor, the $p$-adic integrals localise on the $p$-adic integers
\beq
\int_{\mathbb{Q}^{\times}_p}\gamma_p(x) |x|_p^{s}d^{\times}x=\int_{\mathbb{Z}_p\backslash \{0\}}|x|_p^{s}d^{\times}x.
\eeq
By remark~\ref{rk-multmeasure}, this yields
\beq
\int_{\mathbb{Z}_p\backslash \{0\}}|x|_p^{s}d^{\times}x=\int_{\mathbb{Z}_p} |x|_p^{s}d^{\times}x= \frac{1}{1-p^{-s}}.
\eeq
Combining everything and performing the product over primes we obtain
\beq
\int_{\mathbb{A}^{\times}}\gamma_\mathbb{A}(x) |x|_\mathbb{A}^{s} d^{\times}x= \pi^{-s/2}\Gamma(s/2)\prod_{p< \infty} \frac{1}{1-p^{-s}} = \pi^{-s/2}\Gamma(s/2)\zeta(s)= \xi(s).
\eeq
\end{proof} 

\begin{remark} The above result illustrates that the adelic approach gives an elegant integral representation of the 
completed zeta function, where the normalisation factor corresponds to the contribution from the archimedean place $p=\infty$. Such 
integrals were first considered in the thesis of Tate \cite{Tate}, and then developed further by Jacquet and Langlands \cite{JL}. 
\end{remark}

\begin{remark}
It is common to define the archimedean zeta factor by $\zeta_\infty(s)=\pi^{-s/2}\Gamma(s/2)$ and write the global Euler product form of the 
completed Riemann zeta function as
\beq
\xi(s)=\prod_{p\leq \infty} \zeta_p(s).
\eeq
Anticipating later notions, the completed Riemann zeta function can be thought of as an automorphic form on the group $GL(1,\ads)$.
\end{remark}

\subsection{The functional relation}
We shall now take the analysis one step further and prove the following famous theorem using the adelic framework.
\begin{theorem}[Functional relation for the completed Riemann zeta function]
\label{thm:Riemann-func-rel}
The completed Riemann zeta function satisfies the functional relation \begin{align}
\label{CRZfuncrel}
\xi(s) = \xi(1-s).
\end{align}
\end{theorem}
To prove the theorem using the approach of Tate \cite{Tate} we first need the following lemma:
\begin{lemma}[Adelic Poisson resummation] For any (sufficiently nice) function $f_{\mathbb{A}}$ with adelic Fourier transform $\tilde f_\ads$ we have the Poisson summation formula
\beq
\sum_{\gamma\in \mathbb{Q}}f_{\mathbb{A}}(\gamma)=\sum_{\gamma\in \mathbb{Q}} \tilde{f}_\mathbb{A}(\gamma).
\eeq
\end{lemma}
\begin{proof} 
The proof is similar to the proof of the ordinary Poisson summation formula so we will be brief. Define
\beq
F_\mathbb{A}(x)=\sum_{\gamma\in\mathbb{Q}}f_\mathbb{A}(x\gamma).
\label{defF}
\eeq
This function is periodic by construction and so has a Fourier expansion. The Fourier coefficients $F_{\psi_\gamma}$ of $F_\mathbb{A}$ with respect a unitary character $\psi_\gamma$ precisely 
equals the Fourier transform $\tilde{f}_\mathbb{A}(\gamma)$ of the seed function $f_\mathbb{A}(\gamma)$ and so we can write 
\beq
F_\mathbb{A}(x)=\sum_{\gamma\in \mathbb{Q}} \tilde{f}_\mathbb{A}(\gamma) \psi_\gamma(x).
\eeq
Putting $x=1$ in this formula equating it with $F_\mathbb{A}(1)$ from the definition (\ref{defF}) then establishes the result.
\end{proof}
To complete the proof of the theorem we need also the following lemma:
\begin{lemma} 
\label{lemmatheta}
The global theta function 
\beq 
\Theta_\mathbb{A}(x)=\sum_{k\in \mathbb{Q}} \gamma_\mathbb{A}(kx)
\eeq 
satisfies the functional relation
\beq
\Theta_\mathbb{A}(x) = \frac{1}{|x|_\mathbb{A}} \Theta_\mathbb{A}(1/x), \qquad \forall x\in \mathbb{A}^\times.
\eeq
\end{lemma} 
\begin{proof} 
This follows from applying the Poisson summation formula and the fact that the global Gaussian $\gamma_{\mathbb{A}}(x)$ is invariant under Fourier transform.
\end{proof}

\begin{proof}[Proof of theorem \ref{thm:Riemann-func-rel}]
Now let $J$ be the fundamental domain defined in \eqref{funddomain} for $\rats^\times \bs \ads^\times$. By lemma \ref{lemmatheta} we then have
\beq
\int_J \Theta_\mathbb{A}(x) |x|_{\mathbb{A}}^s d^\times x=\int_J \Theta_\mathbb{A}(1/x) |x|_{\mathbb{A}}^{s-1} d^\times x = \int_J \Theta_\mathbb{A}(x) |x|_{\mathbb{A}}^{1-s} d^\times x,
\label{integralJ}
\eeq 
where in the last step we used the fact that the multiplicative measure is invariant under $x\rightarrow x^{-1}$. Finally, using the factorisation 
$\mathbb{A}^{\times} = \bigcup_{k\in \mathbb{Q}^\times} k\cdot J$ (see \eqref{StrongapproxAbelian}), and the fact that $|x|_{\mathbb{A}}=1$ for $x\in \mathbb{Q}$, we can rewrite (\ref{integralJ}) as
\beq
\int_{\mathbb{A}^\times} \gamma_\mathbb{A}(x) |x|_{\mathbb{A}}^s d^\times x = \int_{\mathbb{A}^\times} \gamma_{\mathbb{A}}(x) |x|_{\mathbb{A}}^{1-s} d^\times x,
\eeq
thus establishing the functional relation (\ref{CRZfuncrel}) for the completed Riemann zeta function using proposition \ref{prop-Tateglobal}.
\end{proof}

\chapter{Basic notions from Lie algebras and Lie groups}
\label{ch:Lie-groups}

We will make use of some standard terminology from the theory of Lie groups and Lie algebras that we  summarise for definiteness. \index{Lie group}\index{Lie algebra}We first address real Lie algebras and groups and their complexifications before we turn to the adelic setting with emphasis on the strong approximation theorem.

\section{Real Lie algebras and real Lie groups}

The material reviewed in this section can be found for example in~\cite{Humphreys,FultonHarris,Helgason,Kac}.

\subsection{Split real simple Lie algebras and root systems}
\label{sec:simple-lie-alg}

Let $\mf{g}(\cx)$ be a finite-dimensional and simple complex Lie algebra from the Cartan--Killing classification. We will  consider here only the split real form $\mf{g}\equiv\mf{g}(\reals)$ \index{Lie algebra!split real} of the Lie algebra. We choose a \emphindex{Cartan subalgebra} $\mf{h}\subset\mf{g}$, that is, a maximal abelian subalgebra of semi-simple elements. This means that we can decompose $\mf{g}$ into eigenspaces of $\mf{h}$ in what is called the \emphindex[root!space decomposition]{root space decomposition}:
\begin{align}
\mf{g}=\mf{h}\oplus \bigoplus_{\alpha\in\Delta} \mf{g}_\alpha,
\end{align}
where the root space $\mf{g}_\alpha$ for a generalised eigenvalue $\alpha:\,\mf{h}\to\reals$ is given by
\begin{align}\label{galpha}
\mf{g}_\alpha = \left\{ x\in\mf{g}\,\middle|\, \lb h,x \rb = \alpha(h) x\quad\textrm{for all $h\in\mf{h}$}\right\}.
\end{align}
The set of $\alpha\neq 0$ for which $\mf{g}_\alpha\neq\{0\}$ is called the set of \emphindex[root system!of a Lie algebra]{roots} $\Delta$. By our assumption on the Lie algebra $\mf{g}$ we have that $\dim(\mf{g}_\alpha)=1$ for all $\alpha\in\Delta$. 
Since $\lie g_\alpha$ is one-dimensional there is, for each root $\alpha \in \Delta$, a unique element $H_\alpha \in [\lie g_\alpha, \lie g_{-\alpha}] \subset \lie h$ such that $\alpha(H_\alpha) = 2$.

In the set of roots $\Delta\subset\mf{h}^*$ we choose a system of \emphindex[simple root]{simple roots}
\begin{align}
\Pi=\{\alpha_1,\ldots,\alpha_r\},
\end{align}
where $r=\dim(\mf{h})$ is the \emphindex[Lie algebra!rank]{rank} of the Lie algebra. Then any root $\alpha\in\Delta$ can be written as an integral linear combination of the simple roots
\begin{align}
\alpha =\sum_{i=1}^r m_i\alpha_i,
\end{align}
where either all $m_i\geq0$ (and $\alpha$ is called a \emphindex[root!positive]{positive root}: $\alpha>0$) or all $m_i\leq0$ (and $\alpha$ is called a \emphindex[root!negative]{negative root}: $\alpha<0$). The set of positive/negative roots is denoted by $\Delta_{\pm}$ and they satisfy $\Delta_-=-\Delta_+$. There is a unique \emphindex[root!highest]{highest root} $\theta\in\Delta$ for which the \emphindex[root!height]{height} $\mathrm{ht}(\alpha)=\sum_i m_i$ is maximal. Another important element is the \emphindex{Weyl vector} (which is not necessarily an element of $\Delta$)
\begin{align}
\label{Weylvec}
\rho = \frac12 \sum_{\alpha\in\Delta_+}\alpha.
\end{align}

We define the spaces of \emphindex[step operator]{positive/negative step operators} by
\begin{align}
\mf{n}\equiv\mf{n}_+ = \bigoplus_{\alpha\in\Delta_+} \mf{g}_\alpha\quad\textrm{and}\quad
\mf{n}_- = \bigoplus_{\alpha\in\Delta_-} \mf{g}_\alpha,
\end{align}
as well as the  \emphindex[subalgebra!Borel]{(upper) Borel subalgebra}\index{Borel subalgebra}
\begin{align}
\mf{b} = \mf{h}\oplus\mf{n}.
\end{align}
The spaces $\mf{n}_\pm$ are nilpotent subalgebras of $\mf{g}$; the Borel subalgebra $\mf{b}$ is solvable. One can think of $\mf{n}_\pm$ as strictly upper/lower triangular matrices and $\mf{h}$ as diagonal matrices.

More formally, the notions of nilpotency and solvability for Lie algebras are defined as follows. A \index{Lie algebra!nilpotent} nilpotent Lie algebra is one whose \emphindex{lower central series} $D_{k}(\mf{g}) := \lb \mf{g},D_{k-1}(\mf{g})\rb$ vanishes for some finite $k$. A \index{Lie algebra!solvable} solvable Lie algebra is one whose \emphindex{derived series} $D^{k}(\mf{g}) := \lb D^{k-1}(\mf{g}),D^{k-1}(\mf{g})\rb$ vanishes for some finite $k$. The Borel subalgebra $\mf{b}$ includes semi-simple elements whence the lower central series does not vanish. The semi-simple elements disappear in $D^1(\mf{b})=\lb\mf{b},\mf{b}\rb$ and thus the derived series vanishes, rendering $\mf{b}$ solvable. The derived series will play an role when discussing Fourier expansions of automorphic forms in chapter~\ref{ch:fourier}.

On $\mf{g}$ one can define an \emphindex[invariant bilinear form]{invariant bilinear form}\index{Killing!metric} that we will write as $\langle x|y\rangle$ for $x,y\in\mf{g}$. Invariance means compatibility with the Lie bracket:
\begin{align}
    \label{eq:killing-invariance}
\langle  \lb x,y\rb | z\rangle = \langle x|\lb y,z\rb\rangle.
\end{align}
This form is proportional to the Killing metric. We have that $H_\alpha \in [\lie g_\alpha, \lie g_{-\alpha}]$ implies that $H_\alpha = [X_\alpha, Y_\alpha]$ for some $X_\alpha \in \lie g_\alpha$ and $Y_\alpha \in \lie g_{-\alpha}$.

Then, for any $h \in \lie h$
\begin{equation}
    \label{eq:H-alpha-killing}
    \langle H_\alpha | h \rangle = \langle [X_\alpha, Y_\alpha] | h \rangle = \langle X_\alpha | [Y_\alpha, h] \rangle = \alpha(h) \langle X_\alpha | Y_\alpha \rangle \, .
\end{equation}

With $h = H_\alpha$ this becomes
\begin{equation}
    \langle H_\alpha | H_\alpha \rangle = \alpha(H_\alpha) \langle X_\alpha | Y_\alpha \rangle  = 2 \langle X_\alpha | Y_\alpha \rangle
\end{equation}

Thus, by insertion into \eqref{eq:H-alpha-killing}
\begin{equation}
    \alpha(h) = \frac{2 \langle H_\alpha | h \rangle}{\langle H_\alpha | H_\alpha \rangle}
\end{equation}
Sometimes we will also use the notation $\langle \alpha | h \rangle$ for $\alpha(h)$.

The Cartan element $T_\alpha = 2 H_\alpha / \langle H_\alpha | H_\alpha \rangle$ can then be used to define an inner product on $\lie h^*$ by
\begin{equation}
    \langle \alpha | \beta \rangle = \langle T_\alpha | T_\beta \rangle = \alpha(T_\beta) = \beta(T_\alpha) \, .
\end{equation}

Since $\mf{g}$ is finite-dimensional and simple this bilinear form on $\mf{h}^*$ is positive definite and can be used to define the lengths of root vectors $\alpha$. We  normalise it such that the highest root $\theta$ has length $\theta^2:=\langle\theta|\theta\rangle=2$.

The bilinear form on $\mf{h}^*$ (spanned over $\reals$ by the simple roots) can be used to define a basis of $\mf{h}^*$ dual to the simple roots. The corresponding basis elements are called the \emphindex[fundamental weights]{fundamental weights} $\Lambda_i$ and satisfy 
\begin{align}
\label{fwbasis}
\langle \Lambda_i|\alpha_j\rangle = \frac12\langle\alpha_i|\alpha_i\rangle \delta_{ij}\quad\text{for $i,j=1,\ldots,r$}.
\end{align}
In terms of the fundamental weights one can re-express the Weyl vector of equation (\ref{Weylvec}) as $\rho=\sum_{i=1}^r \Lambda_i$. A general element of $\mf{h}^*$ will be called a \emp{weight} and denoted by $\lambda$.

Associated with the choice of simple roots $\alpha_i$ is also a realisation of the \emphindex{Weyl group} of $\mf{g}$. This is a finite Coxeter group that is generated by the \emphindex[fundamental reflection]{fundamental reflections} $w_i$ ($i=1,\ldots,r$) that are defined through their action on weights $\lambda$ by
\begin{align}
\label{funrefl}
w_i(\lambda) = \lambda - \frac{2\langle \lambda|\alpha_i\rangle}{\langle\alpha_i|\alpha_i\rangle}\alpha_i,
\end{align}
so that in particular $w_i(\rho) = \rho-\alpha_i$. A general word of the Weyl group is given by a succession of fundamental reflections $w=w_{i_1}\cdots w_{i_\ell}$ and we call $\ell=\ell(w)$ the \emphindex[Weyl word!length]{length} of the Weyl word $w$. This assumes that the expression is in reduced form, i.e., that the relations between the generating fundamental $w_i$ have been used to make the word as short as possible. We denote the Weyl group by $\Weyl\equiv \Weyl(\mf{g})$ and its distinguished longest element by $\wlong$. The \emphindex[Weyl word!longest]{longest Weyl word} has the property $\wlong(\Delta_+)=\Delta_-$; all other Weyl words map some positive roots to other positive roots.

Since $\dim(\mf{g}_\alpha)=1$ for all roots $\alpha\in\Delta$ and $\Delta_-$ is opposite to $\Delta_+$ we can define for any $\alpha>0$ a triplet
\begin{align}
\label{SL2triplets}
(E_\alpha,H_\alpha,F_{\alpha})\quad\in \mf{g}_{\alpha}\times \mf{h}\times \mf{g}_{-\alpha}
\end{align}
such that the triplet forms a \index{subalgebra!SL2@$\mf{sl}(2,\reals)$}standard $\mf{sl}(2,\reals)$ subalgebra of $\mf{g}$. The relations of one such $\mf{sl}(2,\reals)$ algebra are 
\begin{align}
\left[ H_\alpha, E_\alpha \right] = 2 E_\alpha,\quad
\left[ H_\alpha, F_\alpha \right] = -2 F_\alpha,\quad
\left[ E_\alpha, F_\alpha \right] = H_\alpha.
\end{align}
We also use the notation $E_{-\alpha}=F_{\alpha}$.

Furthermore, we introduce the following notation for the $\mf{sl}(2,\reals)$ triples associated with the simple roots $\alpha_i$ for $i=1,\ldots, r$:
\begin{align}
\label{simpleChev}
e_i \equiv E_{\alpha_i},\quad
f_i \equiv F_{\alpha_i},\quad
h_i \equiv H_{\alpha_i}.
\end{align}
The $h_i$ form a basis of the Cartan subalgebra $\mf{h}$. The $r$ triples $(e_i,h_i,f_i)$ are sometimes referred to as the \emphindex[Chevalley generators!simple]{simple Chevalley generators}.

The \emphindex{Cartan matrix} $A$ is an $r \times r$ matrix defined by the elements
\begin{equation}
\label{CartanMatrix} 
    A_{ij} = \frac{2\langle \alpha_i | \alpha_j \rangle}{\langle \alpha_i | \alpha_i \rangle}  
    = \frac{2 \alpha_j(h_i)}{\alpha_i(h_i)} = \alpha_j(h_i) .
\end{equation}

We then have that
\begin{equation}
    \begin{split}
        [h_i, e_j] &= \phantom{-} \alpha_j(h_i) e_j = \phantom{-} A_{ij} e_j \\
        [h_i, f_j] &= - \alpha_j(h_i) f_j = - A_{ij} f_j \, .
    \end{split}
\end{equation}
as well as the \emphindex{Serre relations}
\begin{align}
\left( \mathrm{ad}\, e_i\right)^{1-A_{ij}} e_j =0\,,\quad
\left( \mathrm{ad}\, f_i\right)^{1-A_{ij}} f_j =0
\end{align}
for $i\neq j$ and where the adjoint action of $\mf{g}$ on itself is $(\mathrm{ad}\, x) y = [x,y]$.

The Lie algebra $\mf{g}$ has a \emphindex[subalgebra!compact]{compact subalgebra} $\mf{k}$ that is spanned by $E_{\alpha}-E_{-\alpha}$. It is of dimension equal to the number of positive roots. All its elements have negative norm in the invariant bilinear form discussed above.

\subsection{Split real Lie groups and highest weight representations}

Many of the notions just introduced carry over to the group level. Let $G(\reals)$ be a \index{Lie group!split real} real Lie group with Lie algebra $\mf{g}$ of the type just discussed. The link between the Lie algebra and Lie group is given by the standard exponential map (in the identity component of $G(\reals)$).

The Cartan subalgebra $\mf{h}$ of commuting elements is the Lie algebra of an \index{subgroup!abelian} abelian subgroup $A(\reals)\subset G(\reals)$ that we take to be the exponential of $\mf{h}$. Topologically, $A(\reals)\cong (GL(1,\reals)_+)^r$, where the $+$ subscript indicates that we restrict to positive elements. An important remark here is that there is a larger abelian subgroup, sometimes called the (split) \emphindex{Cartan torus} that is of the form $(GL(1,\reals))^r$ and covers $A(\reals)$. We will sometimes abuse notation and refer to $A(\reals)$ as the Cartan torus or even refer to the Cartan torus as $A(\reals)$ as it should always be clear from the context which abelian subgroup is meant. 

The space of nilpotent elements $\mf{n}\equiv\mf{n}_+$ is the Lie algebra of a \emphindex[subgroup!unipotent]{unipotent} subgroup $N(\reals)\subset G(\reals)$. 
The compact subalgebra $\mf{k}\subset\mf{g}$ is the Lie algebra of a (maximal) \emphindex[subgroup!compact]{compact} subgroup $K(\reals)\subset G(\reals)$. 

The \emphindex{Iwasawa decomposition} states that one can write any element $g\in G(\reals)$ uniquely as the product of elements of the three subgroups just introduced, i.e.,
\begin{align}
\label{IwasawaR}
G(\reals) = N(\reals) A(\reals) K(\reals)
\end{align}
with uniqueness of decomposition~\cite{Helgason}. 

The split real Lie algebras $\mf{g}(\reals)$ have (non-unitary) irreducible finite-dimensional representations\index{representation!highest weight} labelled by a \emphindex[dominant weight]{dominant highest weight} $\Lambda$. This is an element of $\mf{h}^*$ that has integral non-negative coefficients when expanded in the basis of fundamental weights $\Lambda_i$ that was introduced in (\ref{fwbasis}). In other words, a dominant highest weight $\Lambda$ satisfies
\begin{align}
\langle \Lambda| \alpha_i\rangle \in \mathbb{N}_0\quad \textrm{for all $i=1,\ldots,r$.}
\end{align}
We denote the \emphindex[highest weight representation]{highest weight representation} of a dominant highest weight $\Lambda$ by $V_\Lambda$. The notion of highest weight implies that there is a vector $v_\Lambda\in V_\Lambda$ that satisfies
\begin{subequations}
\label{hstwt}
\begin{align}
h\cdot v_\Lambda &= \Lambda(h)v_\Lambda\quad&&\textrm{for all $h\in \mf{h}$},\\
E_\alpha \cdot v_\Lambda &=0\quad&&\textrm{for all positive roots $\alpha\in\Delta_+$}.
\end{align}
\end{subequations}
The first condition reflects that the vector $v_\Lambda$ is in the $\Lambda$-eigenspace of the action of $\mf{g}$ (hence it is a weight vector)  and the second condition shows that it is annihilated by all raising operators (hence at highest weight). Here, we have denoted the action of $\mf{g}$ on the representation space $V_\Lambda$ by $\cdot$ for brevity.

The structure of highest weight representations $V_\Lambda$ can be conveniently summarised in terms of its (formal) \emphindex[character!of a highest weight representation]{character} 
\begin{align}
\label{charactersum}
\mathrm{ch}_\Lambda = \sum_{\mu \in \mf{h}^*} \textrm{mult}_{V_\Lambda}(\mu) e^\mu,
\end{align}
where $\textrm{mult}_{V_\Lambda}(\mu)$ denotes the \emphindex[weight multiplicity]{weight multiplicity} of a weight $\mu\in\mf{h}^*$ in the representation $V_\Lambda$, i.e., the dimension of the $\mu$-eigenspace of the action of $\mf{g}$ on $V_\Lambda$. The expression $e^\mu$ denotes an element of the \emphindex{group algebra} of $\mf{h}^*$ and satisfies $e^{\mu_1} e^{\mu_2}=e^{\mu_1+\mu_2}$ for two weights $\mu_1$ and $\mu_2$. Any representation has a character but the advantage of highest weight representations is that there is a nice compact formula that determines the character $\mathrm{ch}_\Lambda$ in terms of $\Lambda$, the root structure of $\mf{g}$ and its Weyl group. This formula is the \emphindex{Weyl character formula}~\cite{Humphreys,FultonHarris}:
\begin{align}
\label{WCF}
\mathrm{ch}_\Lambda = \frac{\sum_{w\in\Weyl}\eps(w) e^{w(\Lambda+\rho)-\rho}}{\prod_{\alpha>0} (1-e^{-\alpha})}.
\end{align}
The product in the denominator is over all positive roots $\alpha \in \Delta_+$ of the algebra $\mf{g}$ and $\rho$ is the Weyl vector defined in (\ref{Weylvec}). The sign $\epsilon(w)=(-1)^{\ell(w)}$ gives the signature of $w$ as an even or odd element in $\Weyl$. As a special case for $\Lambda=0$ one obtains the one-dimensional \emphindex[representation!trivial]{trivial representation} with $\mathrm{ch}_{\Lambda=0} = 1$. This implies the \emphindex{denominator formula}
\begin{align}
\sum_{w\in\Weyl}\eps(w) e^{w(\rho)-\rho}=\prod_{\alpha>0} (1-e^{-\alpha})
\label{denominatorformula}
\end{align}
that ties the structure of the Weyl group to the structure of the root system. There is an alternative form of the character formula that will play a r\^ole later on. This based on observing that
\begin{align}
w \left( e^{\rho} \prod_{\alpha>0} (1-e^{-\alpha})\right) = \epsilon(w) e^\rho \prod_{\alpha>0} (1-e^{-\alpha})
\end{align}
is $\Weyl$ skew-invariant, as follows for example from the denominator identity. This implies that one can write the character $\mathrm{ch}_\Lambda$ alternatively as
\begin{align}
\label{WCF2}
\mathrm{ch}_\Lambda = \sum_{w\in\mathcal{W}} w\left( \frac{e^{\Lambda+\rho}}{e^\rho\prod_{\alpha>0} (1-e^{-\alpha})}\right) = \sum_{w\in\mathcal{W}} w\left( \frac{e^{\Lambda}}{\prod_{\alpha>0} (1-e^{-\alpha})}\right) .
\end{align}

The character $\mathrm{ch}_\Lambda$ is not only a formal object but can actually be interpreted as a function $\mathrm{ch}_\Lambda \,:\, \mf{h}(\cx) \to \cx$ on the Cartan subalgebra $\mf{h}(\cx)$ by replacing $e^{\Lambda}(h) =e^{\Lambda(h)}$ etc. everywhere. Then one has $\mathrm{ch}_\Lambda(h)=\mathrm{Tr}_{V_\Lambda} \exp(h)$ for $h\in \mf{h}$; the trace in the highest weight representation $V_\Lambda$. The resulting expression converges everywhere on the complexified Cartan subalgebra. We can also evaluate the character on elements of the maximal torus by the exponential map. Let $a\in A$ then
\begin{align}
\label{charfn}
\mathrm{ch}_\Lambda(a) 
= \frac{\sum_{w\in\mathcal{W}} \eps(w) a^{w(\Lambda+\rho)}a^{-\rho}}{ \prod_{\alpha>0}(1-a^{-\alpha})}
= \sum_{w\in\mathcal{W}} w\left(\frac{a^{\Lambda}}{ \prod_{\alpha>0} (1-a^{-\alpha})}\right).
\end{align}

\begin{remark} 
The highest weight representations $V_\Lambda$ for split real $G(\reals)$ are finite-dimensional, but not unitary. For complex $G(\cx)$ the representation $V_\Lambda$ is irreducible and \emphindex[representation!unitarizable]{unitarizable} for dominant highest weights. 
\end{remark}

\begin{remark}
For Kac--Moody algebras with symmetrizable Cartan matrix, convergence is restricted to the interior of the complexified \emphindex{Tits cone}~\cite[\S 10.6]{Kac}. Since we will not be dealing with this case, we refer the reader to literature.
\end{remark}

\subsection{Borel and parabolic subgroups}
\label{sec:parsubgp}

An important notion for the development of automorphic representations will be that of Borel and parabolic subgroups. The \emphindex[subgroup!Borel]{(upper) Borel subgroup}\index{Borel subgroup} is given by
\begin{align}
B(\reals) = A(\reals)  N(\reals) =  N(\reals)  A(\reals).
\end{align}
Here, the abelian group $A(\reals)$ denotes the full Cartan torus that covers the exponential of the Cartan subalgebra $\mf{h}$. 

A \emphindex[subgroup!parabolic]{(standard) parabolic subgroup} $P(\reals)$ of $G(\reals)$ is a proper subgroup that contains the standard Borel subgroup $B(\reals)$. If we think of $B(\reals)$ as consisting of upper triangular matrices (in $G(\reals)$) then a parabolic subgroup $P(\reals)$ contains all upper triangular matrices as well as some lower triangular ones. The discussion of this section is valid for both $\reals$ and $\cmplx$ and from here on we will suppress the notation of the underlying field.

Standard parabolic subgroups can be described by choosing a subset $\Sigma$ of the simple roots $\Pi$ of $\lie g$ \cite{CollingwoodMcGovern}. The subset $\Sigma\subset\Pi$ generates a root system $\langle \Sigma \rangle$ which defines a \emphindex[subalgebra!parabolic]{parabolic subalgebra} as follows
\begin{equation}
    \lie p = \lie h \oplus \bigoplus_{\mathclap{\alpha \in \Delta(\lie p)}} \lie g_\alpha \qquad \text{ where } \Delta(\lie p) = \Delta_+ \cup \langle \Sigma \rangle \, .
\end{equation}
For clarity of notation, we suppress typically the dependence on the subset $\Sigma$.

The parabolic subalgebra can be decomposed into semi-simple \emphindex{Levi subalgebra} $\lie l$ and a \emphindex[subalgebra!nilpotent]{nilpotent subalgebra} $\lie u$
\begin{equation}
    \lie p = \lie l \oplus \lie u
\end{equation}
which is called a \emphindex[Levi decomposition!of parabolic subalgebra]{Levi decomposition}.

Explicitly,
\begin{equation}
        \lie l = \lie h \oplus \bigoplus_{\mathclap{\alpha \in \langle \Sigma \rangle}} \lie g_\alpha \qquad
        \lie u = \bigoplus_{\mathclap{\alpha \in \Delta_+ \setminus \langle \Sigma \rangle_+}} \lie g_\alpha \, ,
\end{equation}
where $\langle \Sigma \rangle_+ = \Delta_+ \cap \langle \Sigma \rangle$. Henceforth we will often denote the set difference $\Delta_+ \setminus \langle \Sigma \rangle_+$ as $\Delta(\lie u)$. We also note that $\lie l$ has the same rank as $\lie g$.

The reductive Levi subalgebra is often decomposed further into
\begin{equation}
\label{LanglandsLA}
\lie l = \lie m \oplus \lie a_P
\end{equation}
with $\lie m = [ \lie l, \lie l ]$ being semi-simple and $\mf{a}_P \subset \lie h$ being abelian. The decomposition
\begin{equation}
    \lie p = \lie m \oplus \lie a_P \oplus \lie u 
\end{equation}
of the parabolic subalgebra is referred to as the \emphindex{Langlands decomposition}.

Note that we have decorated $\lie a_P$ with a subscript $P$ to distinguish its corresponding group $A_P$ from the $A$ in the Iwasawa decomposition. Recall that we use $\lie h$ (and not $\lie a$) for the Cartan subalgebra of $\lie g$.

Explicitly we have that 
\begin{equation}
    \begin{split}
        \lie a_P &= \{h \in \lie h \mid \alpha(h) = 0 \; \textrm{ for all } \alpha \in \Sigma \} \\ 
        \lie m &= [\lie l, \lie l] = \lie a_P^\perp \oplus \bigoplus_{\mathclap{\alpha \in \langle \Sigma \rangle}} \lie g_\alpha \, ,
    \end{split}
\end{equation}
where the orthogonal complement $\lie a_P^\perp$ is taken within $\lie h$ with respect to the invariant bilinear form $\langle \cdot | \cdot \rangle$.

\begin{example}[Parabolic subgroups of $\mf{sl}(3,\reals)$]
\label{SL3pars}
As an example we consider the Lie algebra $\mf{g}(\reals)=\mf{sl}(3,\reals)$ of type $A_2$. It has two simple roots $\Pi(\mf{g})= \left\{ \alpha_1, \alpha_2\right\}$ and positive roots given by
\begin{align}
\Delta_+(\mf{g}) =  \left\{ \alpha_1, \alpha_2, \alpha_1+\alpha_2\right\}.
\end{align}
Choosing the subset $\Sigma=\left\{\alpha_1\right\}$ defines a parabolic subalgebra $\mf{p}(\reals)\subset\mf{sl}(3,\reals)$ with root system
\begin{align}
\Delta(\mf{p}) = \underbrace{\left\{ \alpha_1,-\alpha_1\right\}}_{\Delta(\mf{l})}  \cup \underbrace{\left\{ \alpha_2, \alpha_1+\alpha_2\right\}}_{\Delta(\mf{u})}.
\end{align}
The Levi subalgebra $\mf{l}(\reals)$ consists of the embedded $\mf{sl}(2,\reals)$ associated with the simple root $\alpha_1$, together with an additional abelian element:
\begin{align}
\mf{l}(\reals) = \underbrace{\mf{sl}(2,\reals)}_{\mf{m}(\reals)} \oplus \underbrace{\phantom{(}\reals\phantom{)}}_{\mf{a}(\reals)}
\end{align}
The nilpotent part $\mf{u}(\reals)$ is a two-dimensional abelian Lie algebra and transforms in the two-dimensional representation of $\mf{l}(\reals)$.

As (traceless) $(3\times 3)$-matrices the elements of $\mf{p}(\reals)$, $\mf{l}(\reals)$ and $\mf{u}(\reals)$ take the forms
\begin{align}
\mf{p}(\reals):\,
\begin{pmatrix}
* & * &*\\
* & * &*\\
0 &0 &*
\end{pmatrix},\quad
\mf{l}(\reals):\,
\begin{pmatrix}
* & * &0\\
* & * &0\\
0 &0 &*
\end{pmatrix},\quad
\mf{u}(\reals):\,
\begin{pmatrix}
0 & 0 &*\\
0 & 0 &*\\
0 &0 &0
\end{pmatrix}.
\end{align}
\end{example}

At the level of Lie groups there are corresponding notions. Let $P$ be a connected group having $\mf{p}$ as its Lie algebra. Then there are (unique) decompositions
\begin{align}
\label{LanglandsLG}
P = L U = M A_P U \, ,
\end{align}
also called the \emphindex{Levi decomposition} and \emphindex{Langlands decomposition}. The subgroup $L$ is called the \emphindex[subgroup!Levi]{Levi subgroup} and $U$ the \emphindex[subgroup!unipotent|textbf]{unipotent subgroup} or \emphindex[unipotent radical]{unipotent radical} of the parabolic subgroup $P \subset G$. 

A particularly important class of parabolic subgroups is furnished by the so-called \emphindex[subgroup!parabolic!maximal|textbf]{maximal parabolic subgroups}. These are in a sense the largest (proper) parabolic subgroups and are characterised by choosing as a defining set $\Sigma$ all simple roots of $G$ but one: $\Sigma=\Pi\setminus\{\alpha_{i_*}\}$, where we denoted the simple root that is left out by $\alpha_{i_*}$. We will use the notation $P_{i_*}$ to denote the maximal parabolic subgroup associated with such a choice. For maximal parabolic subgroups one has that
\begin{align}
\label{maxparLevi}
L = GL(1)\times M
\end{align}
where $M$ is a semi-simple Lie group. The Dynkin diagram of its Lie algebra is obtained by removing the node $i_*$ from the Dynkin diagram of $\mf{g}$. The parabolic subgroup of example~\ref{SL3pars} is maximal and corresponds to the choice $i_*=2$.

\subsection{Chevalley group notation and discrete subgroups}
\label{sec:discgps}

Using the exponential map, we will often parametrise group elements in terms of some basic elements. Concretely, we define
for roots $\alpha\in\Delta$ and $u\in\reals$ (or another base field of the split Lie algebra)
\begin{align}\label{chevalleydefn}
x_\alpha(u) = \exp( u E_\alpha),
\end{align}
where $E_\alpha$ is the distinguished element of the root space $\mf{g}_\alpha$ that appears in the Chevalley basis constructed in (\ref{SL2triplets}). The one-parameter group generated by $x_\alpha(u)$ for $u\in \reals$ will be denoted by $N_\alpha(\reals)$. 

Furthermore, let
\begin{align}
w_\alpha(u) = x_\alpha(u) x_{-\alpha}(-u^{-1}) x_\alpha(u)\quad\text{and}\quad 
h_\alpha(u) = w_\alpha(u) w_{\alpha}(1)^{-1}.
\end{align}
The notation $w_\alpha(u)$ is connected to the Weyl group defined above by noting that the $w_{\alpha_i}(1)$ (for simple $\alpha_i$) generate a cover of the Weyl group~\cite{Kac}. For $u\approx 1$, the element $h_\alpha(u)$ yields $H_\alpha$.

These elements so defined satisfy
\begin{align}
x_\alpha(u) x_\alpha(v) = x_\alpha(u+v)\quad\text{and}\quad 
h_{\alpha}(u) h_\alpha(v) = h_\alpha(uv)
\end{align}
and for $\alpha\neq-\beta$
\begin{align}
x_\alpha(u)x_\beta(v)x_\alpha(u)^{-1} x_\beta(v)^{-1} 
= \prod_{m,n>0\atop m\alpha+n\beta\in\Delta} x_{m\alpha+n\beta} (c^{\alpha\beta}_{mn} u^m v^n),
\end{align}
which is the exponentiation of the relation $\lb E_\alpha,E_\beta\rb \propto E_{\alpha+\beta}$. The constants $c^{\alpha\beta}_{mn}$ depend on the chosen order in the product and the structure constants of the basis $\{E_{\alpha}\}$. If $\alpha=-\beta$ we obtain instead
\begin{align}
w_\alpha(u)x_\alpha(v) w_\alpha(-u) = x_{-\alpha}(-u^{-2}v).
\end{align}
This is related to the commutator $\lb E_\alpha, E_{-\alpha}\rb=H_{\alpha}$.

We will always take the discrete subgroup $G(\ints)$ that is generated by the $x_\alpha(u)$ and $h_\alpha(u)$ for integer $u$. This group is called the (adjoint) Chevalley group. Another way of obtaining the group is to consider the integer lattice of elements spanned by the generators $E_{\alpha}$, $H_\alpha$ and $F_{\alpha}$ (for $\alpha>0$). Since the structure constants of $\mf{g}$ are integral in this basis (whence Chevalley basis), this lattice actually defines a Lie algebra $\mf{g}(\ints)$ over the ring $\ints$. 
The group $G(\reals)$ acts on $\mf{g}(\reals)$ via the adjoint action. The alternative definition of $G(\ints)$ is as the stabiliser of $\mf{g}(\ints)$ under this action. The group $G(\ints)$ contains representatives of the Weyl group.

\section{\texorpdfstring{$p$}{p}-adic and adelic groups}
\label{sec:adelic-groups}

In this section we introduce some basic properties of linear algebraic groups defined over a number 
field $\field$, which can be either local or global. In our treatment we shall always take $\field$ to be either 
$\mathbb{Q}$, $\mathbb{Q}_p$ or the ring of adeles $\mathbb{A}$. For more details and proofs, see for instance \cite{Borel,GGPS,Gelbart,Bump}.

Recall that a Lie group $G=G(\mathbb{C})$ defined over $\mathbb{C}$ is a differentiable manifold with a compatible group structue. More generally, one can consider \emphindex[algebraic group]{algebraic groups} $G(\field)$ over any number field $\field$. Formally, the group $G(\field)$ is an (affine) algebraic variety equipped with a group structure given by polynomial operations with values in $\field$. We will be interested in linear algebraic groups over $\field$, which are subgroups $G(\field)$ of the group $GL(n, \field)$ of invertible $n\times n$ matrices with entries in $\field$. As we shall see, the notion of algebraic group extends to local fields, like $\field=\mathbb{Q}_p$, or global fields, like $\field=\mathbb{Q}$, or even the adeles $\field=\ads$.

\subsection{\texorpdfstring{$p$}{p}-adic groups}

We shall now take a closer look at algebraic groups $G$ over the local field of $p$-adic numbers $\mathbb{Q}_p$. At the infinite place, $\mathbb{Q}_\infty=\mathbb{R}$, this is just a real Lie group $G(\mathbb{R})$ corresponding to a real form of a complex Lie group $G(\mathbb{C})$. Unless otherwise specified, we will always take this to be the split real form discussed in the preceding section.

Let us focus on the non-archimedean completions $\mathbb{Q}_{p<\infty}$ of $\mathbb{Q}$, comparing with the more familiar situation of real Lie groups $G(\mathbb{R})$ where it is appropriate. If $G(\mathbb{Q})$ is a \emphindex[group!linear algebraic]{linear algebraic subgroup} of $GL(n,\mathbb{C})$ defined by polynomial conditions with coefficients in $\mathbb{Q}$, then we can also speak of the local linear algebraic group $G(\mathbb{Q}_p)$ defined by the same polynomial conditions, but now taken over $\mathbb{Q}_p$. The typical example of a linear algebraic group is $SL(n,\mathbb{C})$ that is defined as the subgroup of $GL(n,\cx)$ such that the polynomial equation $\det(g)=1$ is satisfied. As this equation has rational coefficients, one can define the local linear algebraic groups $GL(n,\mathbb{Q}_p)$ and $SL(n, \mathbb{Q}_p)$, which are simply the corresponding groups of $n\times n$ matrices with entries in $\mathbb{Q}_p$. 

An important fact is that the notion of maximal compact subgroup carries over to the local setting. Recall that for a real Lie group $G(\mathbb{R})$ in its split real form the maximal compact subgroup $K(G)$ is defined as the fixed point set of $G$ under the Chevalley involution. For example, in the case of $G(\mathbb{R})=SL(n,\mathbb{R})$ we have $K(G)=SO(n)$. For real Lie groups, the maximal compact subgroup is unique up to conjugacy. To understand the analogous notion of maximal compact subgroup of $G(\mathbb{Q}_p)$, recall that the $p$-adic integers $\mathbb{Z}_p$ form a compact ring inside $\mathbb{Q}_p$. It follows that the subgroup of integer points 
\beq
G(\mathbb{Z}_p)=G\cap GL(n,\mathbb{Z}_p),
\eeq
sits compactly inside $G(\mathbb{Q}_p)$. For finite primes $p<\infty$ a maximal compact subgroup of $G(\mathbb{Q}_p)$ is $K_p = G(\mathbb{Z}_p)$. 

\begin{remark}
In contrast to the archimedean case, maximal compact subgroups of $G(\rats_p)$ for $p<\infty$ are not all conjugate. The different conjugacy classes of maximal compact subgroups of $G(\rats_p)$ are most easily classified by \emphindex{Bruhat--Tits theory}~\cite{BruhatTits1,BruhatTits2,HumphreysArithmetic}. For example, for $SL(2,\rats_p)$ there are two conjugacy classes of maximal subgroups. One is represented by $SL(2,\ints_p)$ and the other is  represented by $\left(\begin{smallmatrix}p&\\&1\end{smallmatrix}\right) SL(2,\ints_p) \left(\begin{smallmatrix}p^{-1}&\\&1\end{smallmatrix}\right)$~\cite{HumphreysArithmetic}. For $SL(n,\rats_p)$ there are $n$ classes of maximal compact subgroups in general. (However, there is a unique \emphindex{Iwahori subgroup} that generalises the notion of Borel subgroup to local fields.)
\end{remark}

For real Lie groups $G(\mathbb{R})$ we always have a unique \emphindex{Iwasawa decomposition}
\beq
G(\mathbb{R})= N(\mathbb{R}) A(\mathbb{R}) K(\mathbb{R}),
\eeq
where $K(\mathbb{R})$ is the maximal compact subgroup, $A(\mathbb{R})$ is the Cartan torus and $N(\mathbb{R})$ is the nilpotent subgroup generated by the positive Chevalley generators of the Lie algebra of $G$. The notion of Iwasawa decomposition carries over to the local situation, where we have a decomposition of the form
\beq
\label{IwasawaQp}
G(\mathbb{Q}_p)= N(\mathbb{Q}_p) A(\mathbb{Q}_p) G(\mathbb{Z}_p).
\eeq
In contrast to the case of real groups, the local Iwasawa decomposition is \emph{not} unique, however its restriction to the norm on $A$ is, and this fact will play a crucial role later.

\begin{example}[Iwasawa decompositions in $SL(2,\rats_p)$ for $p\leq \infty$]
\label{Iwex}  We now consider in more detail the example of $G(\mathbb{Q}_p)=SL(2,\mathbb{Q}_p)$. The maximal compact subgroup is $K_p=SL(2,\mathbb{Z}_p)$ and the Iwasawa decomposition reads
\beq 
SL(2, \mathbb{Q}_p)=N(\mathbb{Q}_p)A(\mathbb{Q}_p) SL(2,\mathbb{Z}_p),
\eeq
where
\beq
N(\mathbb{Q}_p) =\left\{\left(\begin{array}{cc} 1 & x \\ & 1 \\ \end{array}\right) \, \Big|\, x \in \mathbb{Q}_p\right\}, \qquad A(\mathbb{Q}_p) =\left\{\left(\begin{array}{cc} a &  \\ & a^{-1} \\ \end{array}\right) \, \Big|\, a \in \mathbb{Q}_p^{\times}\right\}.
\eeq
To illustrate this further, let us consider the explicit Iwasawa decomposition of a specific element 
\beq
g=\left(\begin{array}{cc} 1 & \\ u & 1 \\ \end{array} \right) \in SL(2,\mathbb{Q}_p),
\eeq
which will  be of relevance for the analysis in subsequent sections. First notice that if $u\in \mathbb{Z}_p$ then 
$g$ is already in $SL(2,\mathbb{Z}_p)$ and the decomposition is trivial. Consider therefore the case when $u\in \mathbb{Q}_p\backslash \mathbb{Z}_p$, 
for which one could write a $g=nak$ decomposition as follows
\beq
\label{eq:SL2-p-adic-Iwasawa}
\begin{pmatrix}
1&\\u&1 
\end{pmatrix}
=\begin{pmatrix}
1&u^{-1}\\&1 
\end{pmatrix}
\begin{pmatrix}
u^{-1}&\\&u 
\end{pmatrix}
\begin{pmatrix}
&-1\\
1&u^{-1}
\end{pmatrix}.
\eeq
Notice that for $u\in\rats_p\backslash \ints_p$, the element $u^{-1}\in \ints_p$ and therefore the matrix on the right is in $K_p=SL(2,\ints_p)$ such that this represents \emp{a} valid Iwasawa decomposition. 

As emphasised above, the Iwasawa decomposition for groups over non-archimedean fields is not unique. In the present example, all possible Iwasawa decompositions are of the form
\begin{align}
\begin{pmatrix}
1&\\u&1\end{pmatrix}
= \begin{pmatrix}
1&u^{-1}-ke^{-1} u^{-2}\\&1 
\end{pmatrix}
\begin{pmatrix}
(eu)^{-1}&\\&eu 
\end{pmatrix}
\begin{pmatrix}
k&ku^{-1}-e\\
e^{-1}&e^{-1}u^{-1}
\end{pmatrix}
\end{align}
for arbitrary $k\in \ints_p$ and $e\in\ints_p^\times$. We note that since $|e|_p=1$, the norms of the entries of the element $a\in A(\rats_p)$ is unambiguously defined even though the full Iwasawa decomposition is not unique. The relation~\eqref{eq:SL2-p-adic-Iwasawa} corresponds to $k=0$ and $e=1$. (One can render the $p$-adic Iwasawa decomposition unique by imposing further restrictions on the individual elements~\cite{MR2807433}. We will not use this here.)

It is illuminating to compare \eqref{eq:SL2-p-adic-Iwasawa} with the decomposition of the analogous element in $SL(2,\mathbb{R})$. Thus, take
\beq
\left(\begin{array}{cc} 1 & \\ x & 1 \\ \end{array} \right) \in SL(2,\mathbb{R}),
\eeq
so that in this case $x\in \mathbb{R}$. The \emp{unique} Iwasawa decomposition of this element is
\beq
\label{IwaLowerSL2}
\begin{pmatrix}
1&\\x&1 
\end{pmatrix}
=\begin{pmatrix}
1&\frac{x}{1+x^2}\\&1 
\end{pmatrix}
\begin{pmatrix}
1/\sqrt{1+x^2}&\\&\sqrt{1+x^2} 
\end{pmatrix}
k,
\eeq
with
\begin{equation}
    k = \frac{1}{\sqrt{1+x^2}}
    \begin{pmatrix}
        1 & -x \\
        x & 1
    \end{pmatrix}
    \in SO(2, \reals) . 
\end{equation}
Hence, the component along the Cartan torus in~\eqref{eq:SL2-p-adic-Iwasawa} is in fact simpler in the Iwasawa decomposition of $SL(2,\mathbb{Q}_p)$ compared with that  of  $SL(2,\mathbb{R})$. 
\end{example}

\subsection{Adelisation and strong approximation}
\label{sec_adelisation}

We now discuss the central notion of strong approximation that allows the reformulation of many questions concerning $G(\reals)$ and its automorphic forms in terms of questions on the adelic group $G(\ads)$. The description in this section is general; the following section~\ref{sec:SL2SA} gives more details for the case of $G=SL(2)$.

Given an algebraic group $G$ defined over $\mathbb{Q}$ we can consider the adelic group $G(\mathbb{A})$ as the restricted direct product 
\beq
G(\mathbb{A})=G(\mathbb{R})\times G_f ,
\eeq
where
\begin{equation}
    G_f = \prodp_{p<\infty} G(\mathbb{Q}_p) ,
\end{equation}
consisting of elements $g=(g_p)=(g_\infty; g_2, g_3, g_5, \dots)$ such that all but finitely many $g_p\in G(\mathbb{Z}_p)$. 

\begin{remark}
The adelic group $G(\ads)$ is attached to the algebraic group $G$ over $\rats$ (or more precisely over $\mathrm{Spec}(\rats)$). We typically start from a Chevalley group as in the preceding section~\ref{sec:discgps} that provides a clear basis for the construction of the adelic group $G(\ads)$.
\end{remark}

We further set 
\beq
K_f=\prod_{p< \infty} G(\mathbb{Z}_p)
\eeq
and we then have the notion of maximal compact subgroup $K_{\mathbb{A}} $ of $G(\mathbb{A})$ defined as
\beq
K_{\mathbb{A}} = K_\infty \times K_f,
\eeq
where $K_\infty$ is the maximal compact subgroup of $G(\mathbb{R})$. The adelic version of the Iwasawa decomposition thus reads
\beq
G(\mathbb{A})= N(\mathbb{A})A(\mathbb{A})K_{\mathbb{A}} .
\eeq
When $G$ is split of rank $r$, the adelic Cartan torus is given by 
\beq
A(\mathbb{A})= GL(1,\mathbb{A})\times \cdots \times GL(1,\mathbb{A}) \cong (\mathbb{A}^{\times})^{r} .
\eeq 

Since $\mathbb{Q}$ is discrete in $\mathbb{A}$ according to proposition~\ref{prop:discrete-embedding}, it follows that $G(\mathbb{Q})$ is a discrete subgroup of $G(\mathbb{A})$. This implies that the arithmetic coset space $G(\mathbb{Q})\backslash G(\mathbb{A})$ corresponds to the adelisation of $G(\mathbb{Z})\backslash G(\mathbb{R})$. In fact, topologically $G(\mathbb{Q})\backslash G(\mathbb{A})$ is the total space of a fiber bundle over $G(\mathbb{Z})\backslash G(\mathbb{R})$ \cite{GGPS}:
\beq 
\begin{array}{ccc}
  K_f & \hookrightarrow & G(\mathbb{Q})\backslash G(\mathbb{A})  \\
   &  & \downarrow \\
   &  & G(\mathbb{Z})\backslash G(\mathbb{R})
\end{array}
\eeq
One way of stating \emphindex{strong approximation} then asserts that 
\begin{equation}
    \label{eqn:strong-approx}
    G(\ints) \bs G(\reals) \iso G(\rats) \bs G(\adeles) / K_f .
\end{equation}
This has the very useful consequence that any function $\phi_\reals$ on $G(\ints) \bs G(\reals)$ can be lifted to a function $\phi_\adeles$ on the adelisation $G(\rats)\bs G(\adeles)$ where $\phi_\adeles$ is characterized by being right-invariant under $K_f$. The consequences of this for automorphic forms will be discussed in chapter~\ref{ch:autforms}.

The strong approximation theorem \eqref{eqn:strong-approx} can be stated even more generally for open subgroups $K_\Gamma$ of $K_f$ according to \cite{Deitmar} (see also \cite{Platonov:1993, Prasad:1977}):

\begin{theorem}[Strong approximation theorem]\label{thm:SA}
Let $G$ be a topological group with $G(\rats)$ dense in $G_f$, let $K_\Gamma$ be an open subgroup of $K_f$ and $\Gamma = K_\Gamma \cap G(\rats)$. Then 
\begin{equation}
    \begin{split}
        \phi : \Gamma \bs G(\reals) & \to G(\rats) \bs G(\adeles) / K_\Gamma \\
        \Gamma x_\infty & \mapsto G(\rats) (x_\infty; \id) K_\Gamma
    \end{split}
\end{equation}
is a homeomorphism. Here,  $G(\rats)$ is diagonally embedded in $G(\adeles)$; $G(\reals)$ is embedded as $(x_\infty; \id)$ and $K_\Gamma$ as $(\id; k_p)$.
\end{theorem}

\begin{remark}
An assumption in the theorem is that $G(\rats)$ should be dense in $G_f$ and this is equivalent to the statement that for all open subsets $U$ of $G_f$ we have that $U \cap G(\rats) \neq \emptyset$. An example of such a group $G$ that will be useful for us is $SL(n)$ \cite{Bump}. The assumption of denseness can also be reformulated in terms of the semi-simple group $G$ being simply-connected~\cite{PlatonovRapinchuk}.
\end{remark}

\begin{proof} \mbox{}
(of theorem~\ref{thm:SA})
    \begin{itemize}
        \item $\phi$ is well-defined (independent of coset representative)

            Let $x_\infty, y_\infty \in G(\reals)$ such that $\Gamma x_\infty = \Gamma y_\infty$, that is, there exists a $\gamma \in \Gamma$ such that $x_\infty = \gamma y_\infty$.
            
            We have that $\Gamma = K_\Gamma \cap G(\rats)$. 
            Denoting a double coset in $G(\rats) \bs G(\adeles) / K_\Gamma$ by square brackets (with the real and finite places given separately) this leads to
            \begin{equation}
                \phi(\Gamma x_\infty) = [x_\infty; \id] = [\gamma y_\infty; \id] = [ \gamma y_\infty; \gamma \gamma^{-1}] \stackrel{(a)}{=} [y_\infty; \gamma^{-1}] \stackrel{(b)}{=} [y_\infty; \id] = \phi(\Gamma y_\infty)
            \end{equation}
            where we have used (a) that $\gamma \in G(\rats)$ and (b) that $\gamma^{-1} \in K_\Gamma$

        \item $\phi$ is injective
            
            Assume $\phi(\Gamma x_\infty) = \phi(\Gamma y_\infty)$. Then $G(\rats) (x_\infty; \id) K_\Gamma = G(\rats) (y_\infty; \id) K_\Gamma$, that is, there exists a $\gamma \in G(\rats)$ and $k \in K_\Gamma$ such that $(x_\infty; \id) = \gamma (y_\infty; \id) k = (\gamma y_\infty; \gamma k)$.  

          This means that $x_\infty = \gamma y_\infty$ and $\gamma = k^{-1}$. Since $\gamma \in G(\rats)$ and $k \in K_\Gamma$ we then have that $\gamma = k^{-1} \in K_\Gamma \cap G(\rats) = \Gamma$. Thus, $x_\infty = \gamma y_\infty$ implies that $\Gamma x_\infty = \Gamma y_\infty$.

      \item $\phi$ is surjective \big($G(\adeles) = G(\rats) G(\reals) K_\Gamma$\big)

    Let $x = (x_\infty; x_f)$ be an arbitrary element in $G(\adeles) = G(\reals) \times G_f $. 
    We will now show that since $G(\rats)$ is dense in $G_f$, there exists a $\gamma \in G(\rats)$ such that $\gamma^{-1} x_f \in K_\Gamma$. 

    Consider the continuous map $f : G_f \to G_f$ $g\mapsto g x_f$. $K_\Gamma$ is an open set around $\id$ in $G_f$. Since $f$ is continuous $U = f^{-1}(K_\Gamma)$ is an open set in $G_f$ around $x_f^{-1}$. Let $\gamma^{-1} \in U \cap G(\rats)$ which is non-empty as $G(\rats)$ is dense in $G_f$. Then $\gamma^{-1} x_f = f(\gamma^{-1}) \in f(U) = K_\Gamma$.

    Let $k = \gamma^{-1} x_f \in K_\Gamma$. Then,
    \begin{equation}
        x = (x_\infty; x_f) = (x_\infty; \gamma k) = \gamma(\gamma^{-1} x_\infty; k) = \gamma(\gamma^{-1} x_\infty; \id) k \in G(\rats) G(\reals) K_\Gamma
    \end{equation}
    \end{itemize}
    
\end{proof}

\begin{remark}
The generalisation to open subgroups $K_\Gamma$ is important since it allows to treat different discrete subgroups $\Gamma$ in a uniform way. Typically these subgroups are associated with arithmetically defined congruence subgroups.
\end{remark}

\subsection{Strong approximation for \texorpdfstring{$SL(2,\reals)$}{SL(2,R)}}
\label{sec:SL2SA}

In this section, we illustrate the concepts of the preceding section in some examples involving $G=GL(2)$ and $G=SL(2)$. 

\begin{example}[Discreteness of $GL(2,\rats)$ in $GL(2,\ads)$]
    We will now show that $GL(2, \rats)$ is discrete in $GL(2, \adeles)$ by first considering the identity element. The line of reasoning is analogous to the case of $\rats$ being discretely embedded in $\adeles$ that was treated in proposition~\ref{prop:discrete-embedding}.
    
    Let $U \subset GL(2,\reals)$ be an open neighbourhood of $\id$ such that $U \cap GL(2,\ints) = \{\id\}$. Then 
    \begin{equation}
        V = U \times \prod_{p < \infty} GL(2, \ints_p) \subset GL(2, \adeles) 
    \end{equation}
    is an open neighbourhood of $\id$ in $GL(2, \adeles)$.
        Since, with the diagonal embedding, 
    \begin{equation}
        GL(2, \rats) \cap \prod_{p<\infty} GL(2, \ints_p) = GL(2, \ints)
    \end{equation}
    we then have that $GL(2, \rats) \cap V = \{\id\}$.
    For an arbitrary element $g \in GL(2, \rats)$ these arguments generalise directly by instead considering $g V$. We have then that $GL(2, \rats)$ is discrete in $GL(2, \adeles)$. 

    Perhaps, more importantly for our further calculations, it can similarly be shown that $SL(2, \rats)$ is discrete in $SL(2, \adeles)$.
\end{example}

The next example discusses how the strong approximation theorem~\ref{thm:SA} works for $SL(2)$ and different choices of subgroup $\Gamma$.

\begin{example}[Strong approximation for $SL(2)$]
\label{strongapproxcongruence}
Let $G=SL(2)$. Firstly, let $K_\Gamma = K_f = \prod_{p < \infty} G(\ints_p)$. Then $\Gamma = K_\Gamma \cap G(\rats) = G(\ints)=SL(2,\ints)$, the standard modular group. From the above theorem we then get that
    \begin{equation}
     G(\ints) \bs G(\reals) \iso G(\rats) \bs G(\adeles) / K_f.
    \end{equation}

The second example addresses the principal congruence subgroup $\Gamma_0(N)$. Let locally
    \begin{equation}
        \Gamma_0(N)_p = \left\{ 
            \begin{pmatrix}
                a & b \\
                c & d
            \end{pmatrix} \in SL(2, \ints_p) : c \equiv 0 \mod N\ints_p
        \right\}
    \end{equation}
    and $K_\Gamma(N) =K_0(N):= \prod_{p < \infty} K_\Gamma(N)_p$ where
    \begin{equation}
        K_\Gamma(N)_p = 
        \begin{cases}
            SL(2, \ints_p) & p \nmid N \\
            \Gamma_0(N)_p & p \mid N
        \end{cases}
    \end{equation}
    Since $K_\Gamma \subset K_f = \prod_{p < \infty} SL(2,\ints_p)$ we know that $\Gamma = K_\Gamma \cap SL(2, \rats) \subset SL(2, \ints)$.
    
    That $c \equiv 0 \mod N \ints_p$ for all divisors $p$ of $N$ means that (with $c \in \ints$)
    \begin{equation}
        \begin{split}
            c \in N \ints_p \quad \forall p \mid N & \iff \left| \frac{c}{N} \right|_p \leq 1 \quad \forall p \mid N \\
           & \iff \frac{c}{N} \text{ has no $p$ in the denominator } \quad \forall p \mid N \\
           &\iff c \equiv 0 \mod N
        \end{split}
    \end{equation}

    Thus $\Gamma = \Gamma_0(N)$ and from the above theorem
    \begin{equation}
\Gamma_0(N) \bs SL(2,\reals) \iso SL(2, \rats) \bs SL(2, \adeles) / K_\Gamma(N).
    \end{equation}
\end{example}

We finally exhibit an isomorphism of cosets of the discrete subgroups in $G(\ads)$ with cosets of discrete subgroups in $G(\reals)$ that will be central in section \ref{sec:Eisenstein-adelisation}.

\begin{example}[Bijection of Borel cosets]
\label{bijection}
In this example we will (based on the notes of \cite{Garrett:2014}) show that
\begin{equation}
    \begin{aligned}
        \phi : B(\ints)\bs SL(2, \ints) & \to  B(\rats)\bs SL(2, \rats) \\
        B(\ints) \gamma & \mapsto B(\rats) \gamma
    \end{aligned}
\end{equation}
is an isomorphism, where
\begin{equation}
    B(\field) = \left\{
        \begin{pmatrix}
            * & * \\
            0 & *
        \end{pmatrix} \right\} \cap SL(2, \field) \, .
\end{equation}

The mapping is well-defined since if $B(\ints) \gamma' = B(\ints) \gamma$ then $B(\rats) \gamma' = B(\rats) \gamma$ as $B(\ints) \subset B(\rats)$.

It is injective because if $B(\rats) \gamma' = B(\rats) \gamma$ then there exists a $b$ in $B(\rats)$ such that $\gamma' = b \gamma$, but then $b = \gamma' \gamma^{-1} \in SL(2, \ints)$ which means that $b \in B(\rats) \cap SL(2, \ints) = B(\ints)$. Thus, $B(\ints) \gamma' = B(\ints) \gamma$.

For the surjectivity we need to show that every $B(\rats) g$ with $g \in SL(2, \rats)$ has a representative in $SL(2, \ints)$. Let
\begin{equation}
    g = 
    \begin{pmatrix}
        a & b \\
        c & d 
    \end{pmatrix} \in SL(2, \rats)
    \qquad
    b = 
    \begin{pmatrix}
        q & m \\
        0 & q^{-1}
    \end{pmatrix} \in SL(2, \rats)
    \qquad
    bg =
    \begin{pmatrix}
        q a + mc & qb + md \\
        q^{-1} c & q^{-1} d
    \end{pmatrix}
\end{equation}
where $c = c_1 / c_2$ and $d = d_1 / d_2$ with $c_i, d_i \in \ints$ in shortened form with positive denominators. Now set $q = \gcd(c_1d_2, c_2d_1)/(c_2d_2)$ which makes $q^{-1}c$ and $q^{-1}d$ coprime integers, and thus there exist integers $\alpha$ and $\beta$ such that $\alpha q^{-1} d - \beta q^{-1} c = 1$ by B\'ezout's lemma.

If $c = 0$ then $d \neq 0$, $a = 1/d$ and $q = \gcd(0, c_2d_1)/(c_2 d_2) = \abs{c_2 d_1}/(c_2 d_2) = \abs{d}$ meaning that $qa = q^{-1}d = \pm 1$ and we may choose $m$ such that $qb + md$ is integer. On the other hand, if $c \neq 0$ we may choose $m = (\alpha - qa)/c$ giving $qa + mc = \alpha$ and $q b + md = \beta$ which are both integers. This completes the proof.

\end{example}

\chapter{Automorphic forms}

\label{ch:autforms}

In this chapter we explain how to think about automorphic forms as functions on 
$G(\mathbb{Q})\backslash G(\mathbb{A})$, as opposed to the more familiar concept of $G(\mathbb{Z})$-invariant functions on a real Lie group $G(\mathbb{R})$. After introducing some standard terminology for $SL(2,\reals)$ in section~\ref{sec:SL2}, we begin in section~\ref{sec_frommodtoaut} by discussing the passage from classical modular forms on the upper-half plane $\mathbb{H}$ to automorphic forms on the adelic group $SL(2,\mathbb{A})$. This will serve to illustrate some of the main points in a simple and explicit setting. The general adelic picture of automorphic forms on $G(\ads)$ is introduced in section~\ref{sec_defauto} and Eisenstein series on $G(\ads)$ are defined in section~\ref{sec_Eisenstein}. These notions will lead naturally to automorphic representations and the close connection 
with studying the unitary action of $G(\mathbb{A})$ on the Hilbert space $L^2(G(\mathbb{Q})\backslash G(\mathbb{A}))$, which will be discussed in the subsequent chapter~\ref{ch:autoreps}.

\section{Preliminaries on \texorpdfstring{$SL(2,\reals)$}{SL(2,R)}}
\label{sec:SL2}

In this section we introduce our conventions on notation related to $SL(2,\reals)$  that will be used throughout this book.

\subsection{\texorpdfstring{$SL(2,\reals)$}{SL(2,R)} Lie group and \texorpdfstring{$\mf{sl}(2,\reals)$}{sl(2,R)}  Lie algebra}

We take $SL(2,\reals)$ to be the real Lie group defined (in its fundamental representation) by
\begin{align}
SL(2,\reals) = \left\{ g=\begin{pmatrix} a&b\\c&d\end{pmatrix}\,\middle|\, a,b,c,d\in\reals\,\mathrm{and}\,\, \det(g)=ad-bc=1\right\}.
\end{align}
The maximal compact subgroup is $K=SO(2,\reals)$ corresponding to the orthogonal matrices within $SL(2,\reals)$.

The Lie algebra $\mathfrak{sl}(2,\reals)$ has the standard Chevalley basis
\begin{align}
\label{eq:sl2basis}
e = \begin{pmatrix}0&1\\0&0\end{pmatrix}\,,\quad
h = \begin{pmatrix}1&0\\0&-1\end{pmatrix}\,,\quad
f = \begin{pmatrix}0&0\\1&0\end{pmatrix}
\end{align}
with commutation relations
\begin{align}
\label{sl2rels}
\lb h,e\rb =2e\,,\quad
\lb h,f\rb = -2f\,,\quad
\lb e,f\rb =h\,.
\end{align}
The generator $h$ acts diagonally and is called the Cartan generator;  $e$ is a positive step operator and $f$ a negative step operator. The compact subgroup $SO(2,\reals)$ is generated by the combination $e-f$. 

The universal enveloping algebra $\mathcal{U}(\mf{sl}(2,\cx))$ has a distinguished second order element, called the \emphindex[Casimir operator!for $\mf{sl}(2,\reals)$]{Casimir operator} and that we define by
\begin{align}
\label{eq:SL2cas}
\Omega = \frac14 h^2 + \frac12 ef + \frac12 fe = \frac14 h^2-\frac12 h +ef.
\end{align}
This definition is unique up to normalisation. The Casimir operator commutes with all Lie algebra elements.

The Iwasawa decomposition of $SL(2,\reals)$ can be chosen in the form $SL(2,\reals)=NAK$; where $N$ is in the image of the exponential map $\exp$ applied to $e$; the maximal torus is in the image of $\exp$ applied to $h$ and $K$ is the compact subgroup $SO(2,\reals)$ whose identity component is the exponential of $e-f$. Concretely that means that we can write any element $g$ of $SL(2,\reals)$ as
\begin{align}
\label{SL2IWA:app}
    g = nak &= \exp\left(x e\right) \exp\left(\frac{1}{2}\log(y) h\right) \exp\left(\theta (e-f)\right)\nonumber\\ 
    & = 
    \begin{pmatrix}
        1 & x \\
        0 & 1
    \end{pmatrix}
    \begin{pmatrix}
        y^{1/2} & 0 \\
        0 & y^{-1/2}
    \end{pmatrix} 
    \begin{pmatrix}
        \cos\theta & \sin\theta \\
        -\sin\theta & \cos\theta
    \end{pmatrix}
\end{align}
with $k\in K=SO(2,\reals)$ and $y>0$.  

\subsection{The upper half plane \texorpdfstring{$\UHP$}{UHP} and \texorpdfstring{$SL(2,\mathbb{Z})$}{SL(2,Z)}}

A main object of interest to us is the two-dimensional coset space $G/K=SL(2,\reals)/SO(2,\reals)$; a representative for any point of this space is given by the first two factors in (\ref{SL2IWA:app}). The coset space can therefore be parametrised by elements of the \emphindex[upper half plane]{upper half plane}
\begin{align}
\label{eq:UHPn}
\UHP = \left\{ z=x + iy \st x,y \in \reals \text{ and } y > 0 \right\}
\cong G/K.
\end{align}

The coset space $G/K$ (or, equivalently, the upper half plane $\UHP$) carries an action of $SL(2,\reals)$ by left multiplication: An element $\gamma \in G$ transforms a $g$ into $g'=\gamma g$. The action on the explicit parameters $z \in \UHP$ can be read off from writing the new element in Iwasawa form $g'= n'a'k'$. Performing this calculation one finds
\begin{align}
\label{SL2act:app}
z' = \gamma \cdot z = \frac{az + b}{cz + d} \quad \text{for} \quad 
\gamma=\begin{pmatrix}a&b\\c&d\end{pmatrix}\in SL(2,\reals).
\end{align}
Using the Iwasawa decomposition~\eqref{SL2IWA:app}, we see that the point $i$ is left invariant by the maximal compact subgroup $K=SO(2,\reals)$ and that
\begin{align}
g\cdot i = x + iy = z.
\end{align}

For $SL(2,\reals)$, automorphic forms are functions $f(g)$ that are invariant under the action of a discrete subgroup $\Gamma\subset SL(2,\reals)$. Taking $\Gamma=SL(2,\ints)$ to consist of the $SL(2,\reals)$ matrices with integral entries, the invariance $f(\gamma g)=f(g)$ for all $\gamma\in SL(2,\ints)$ means that $f$ is, in fact, a function on the double quotient $SL(2,\ints)\backslash SL(2,\reals)/SO(2,\reals)$ where $SL(2,\ints)$-equivalent points are identified. Using the upper half plane $\UHP$ presentation of $SL(2,\reals)/SO(2,\reals)$ one can give a very explicit description of the double quotient.

\begin{figure}
\centering
{\scalebox{.8}{\input{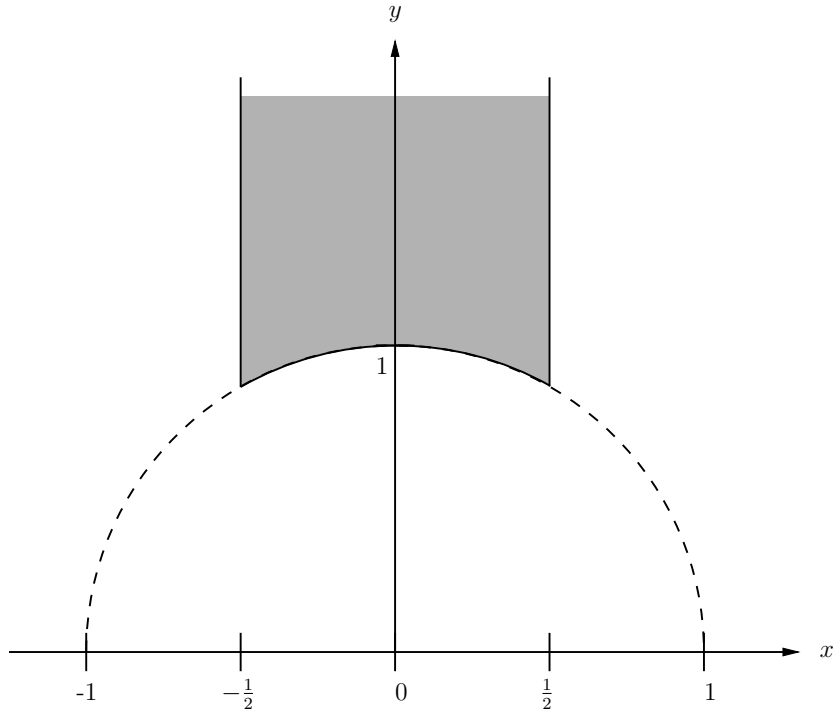}}}
\caption{\label{fig:SL2fund} \textit{A fundamental domain for the action of $SL(2,\ints)$ acting on the upper half plane (grey region). The cusp is at $y \to \infty$.}}
\end{figure}

The group $SL(2,\ints)$ is well-known to be generated by~\cite{MR2112196}:
\begin{align}
T = \begin{pmatrix}1&1\\0&1\end{pmatrix}\quad\mathrm{and}\quad
S = \begin{pmatrix}0&-1\\1&0\end{pmatrix}.
\end{align}
When acting on $z \in \UHP$, they generate
\begin{align}
    T \cdot z = z + 1, \quad S \cdot z = -\frac{1}{z}.
\end{align}
Therefore, $T$ is a translation by one unit and $S$ is inversion in the unit circle combined with a reflection in the $y$-axis. A fundamental domain for the action of $SL(2,\ints)$ on the upper half plane (parametrising the double quotient above) is depicted in figure~\ref{fig:SL2fund}. The fundamental domain clearly displays a single \emphindex{cusp} where it touches the boundary of the space. This cusp corresponds to the limit $y \to \infty$. For a congruence subgroup $\Gamma \subset SL(2,\ints)$ (discussed more in section \ref{LevelNforms}) the cusps are the $\Gamma$-equivalence classes of points in $\rats \cup \{\infty\}$ and they can all be reached by $SL(2,\ints)$-transformations of the $y \to \infty$ cusp~\cite{MR2112196}. The corresponding fundamental domain, similar to figure \ref{fig:SL2fund}, touches the real axis at these points forming `cusps'.

\begin{remark}
What we are dealing with is effectively $PSL(2,\ints)$ rather than $SL(2,\ints)$, where the `$P$' indicates that a matrix has to be identified with minus itself. The reason is that the two matrices have identical action on the upper half plane as easily verified from (\ref{SL2act:app}).
\end{remark}

\section[From classical modular forms to (adelic) automorphic forms]{From classical modular forms to (adelic) \\ automorphic forms}
\label{sec_frommodtoaut}

In this section we show how to pass from the classical notion of a modular form as a function on the complex upper-half plane $\mathbb{H}$ to an automorphic 
form as a function on a Lie group $G$. Here we focus on the example of $SL(2)$, leaving the generalisation to higher rank groups to subsequent sections. We begin by 
going from $\mathbb{H}$ to $SL(2,\mathbb{R})$ and then further to the adelic group $SL(2,\mathbb{A})$. 

\subsection{Holomorphic modular forms}
A \emphindex[modular form!holomorphic]{holomorphic modular form} of weight $w\geq 0$ is a holomorphic function $ f\, :\, \mathbb{H}\rightarrow \mathbb{C}$ which transforms according to 
\beq
f\left(\frac{az + d}{cz+d}\right) = (cz +d)^w f(z),
\label{holmodform}
\eeq
under the discrete action of 
\beq
\begin{pmatrix} a & b \\ c & d \\ \end{pmatrix} \in SL(2,\mathbb{Z}).
\eeq
If $f(z)$ has zero weight, $w=0$, we call it a \emphindex{modular function}. The prefactor $(cz+d)^w$ in~\eqref{holmodform} is often referred to as \emphindex{factor of automorphy}. The defining equation (\ref{holmodform}) implies that $f$ is periodic $f(z+1)=f(z)$ (for any weight $w$) and thus has a Fourier expansion of the form
\beq
f(z)=\sum_{n\in \mathbb{Z}} a(n) q^{n}, \qquad q:=e^{2\pi i z}.
\label{holFourierexp}
\eeq
Decomposing $q=e^{2\pi i z}=e^{2\pi i x}e^{-2\pi y}$ the Fourier coefficients can be computed from  the standard Fourier transform
\beq
a(n) e^{-2\pi ny} = \int_0^1 e^{-2\pi i nx} f(x+iy)dx.
\label{integralformulacoeff}
\eeq
This formula (and its generalisations) will play a key role in subsequent chapters. 

The moderate growth condition mentioned in section \ref{sec:AutoIntro} can be formulated as the statement that 
\beq
|f(x+iy)|\leq C \cdot y^N
\label{moderategrowthSL2}
\eeq
for some constants $C, N$ as $y\to \infty$ for any $x\in \mathbb{R}$. For holomorphic modular forms this is in fact equivalent to the statement that all negative Fourier coefficients $a(n), n<0,$ in (\ref{holFourierexp}) vanish. To see this we simply use the integral representation (\ref{integralformulacoeff}) for the Fourier coefficient and calculate its norm:
\beq
|a(n) e^{-2\pi ny}|=\left|\int_0^1 e^{-2\pi i nx} f(x+iy)dx\right|
\eeq
Removing the oscillating exponential we obtain a sequence of inequalities
\beq
\left|\int_0^1 e^{-2\pi i nx} f(x+iy)dx\right|\leq \int_0^1 \left|f(x+iy)\right| dx\leq \int_0^1 C\cdot y^N dx=C\cdot y^N. 
\eeq
Thus we arrive at the  inequality
\beq
|a(n) e^{-2\pi ny}|\leq C\cdot y^N.
\eeq
and when $n<0$ the exponential $e^{-2\pi ny}$ blows up as $y\to\infty$ so therefore we must have $a(n)=0$ for $n<0$ as claimed.

\begin{example}[Classical holomorphic Eisenstein series]
Classic examples of holomorphic modular forms on $\mathbb{H}$ are provided by the \emphindex[Eisenstein series!holomorphic for $SL(2,\ints)$]{holomorphic Eisenstein series} defined by
\begin{align}
\label{holEisen}
E_{2w}(z) = \frac{1}{2}\sum_{(c,d)\in\ints^2 \atop (c,d)=1} \frac{1}{(cz+d)^{2w}}.
\end{align}
One can check that this satisfies all the criteria stated above for integral $w\geq 2$. The (finite-dimensional) space $\mathcal{M}_{2w}(SL(2,\mathbb{Z}))$ of weight $2w$ holomorphic modular forms is a ring, famously generated by the Eisenstein series $E_4(z)$ and $E_6(z)$~(see, e.g., \cite{MR2409678} for a proof). The Fourier expansions of $E_4(z)$ and $E_6(z)$ are given by
\beqa
E_4(z)&=&1+240\sum_{n=1}^{\infty} \sigma_3(n)q^n=1+240q+2160q^2+\cdots 
\nn \\
E_6(z)&=& 1-504\sum_{n=1}^{\infty} \sigma_5(n)q^n=1-504q-16632q^2+\cdots
\eqa
where 
\beq
\sigma_s(n)=\sum_{d|n} d^s
\eeq
is the divisor function. For proofs see for example the classic book by Serre \cite{MR0344216}.
\end{example}

\subsection{Modular forms for congruence subgroups*}
\label{LevelNforms}
It is often of interest in number theory to consider holomorphic modular forms for congruence subgroups $\Gamma\subset SL(2,\mathbb{Z})$ (a good reference is the book by Diamond and Shurman \cite{MR2112196}). These satisfy an analogous relation to (\ref{holmodform}) but with extra restrictions on the transformation matrix and possibly with a character appearing on the right hand side. Consider for example the congruence subgroups 
\beqa
\Gamma_1(N)&=&\left\{\begin{pmatrix} a & b \\ c & d \\ \end{pmatrix} \in SL(2,\mathbb{Z})\, \bigg| \, N|c, a-1, d-1\right\}
\nn \\ 
\Gamma_0(N)&=&\left\{\begin{pmatrix} a & b \\ c & d \\ \end{pmatrix} \in SL(2,\mathbb{Z})\, \bigg| \, N|c \right\},
\eqa
where $\Gamma_0(N)$ contains $\Gamma_1(N)$ as a normal subgroup of finite index $\phi(N)$ (the Euler totient function). The space of weight $w$ modular forms for $\Gamma_1(N)$ (resp. $\Gamma_0(N)$) is then denoted by $\mathcal{M}_w(\Gamma_1(N))$ (resp. $\mathcal{M}_w(\Gamma_0(N))$). Since the  quotient $\Gamma_0(N)/\Gamma_1(N)$ is isomorphic to the multiplicative group $(\mathbb{Z}/N\mathbb{Z})^{\times}$ of order $\phi(N)$, one can relate modular forms on these different congruence subgroups through the introduction of Dirichlet characters. A \emph{Dirichlet character} $\chi$ is a group homomorphism 
\beq
\chi : (\mathbb{Z}/N\mathbb{Z})^{\times} \to \mathbb{C}^{\times}
\eeq
where the product between two Dirichlet characters $\chi_1$ and $\chi_2$ is defined by $(\chi_1\chi_2)(g)=\chi_1(g)\chi_2(g)$ for $g\in (\mathbb{Z}/N\mathbb{Z})^{\times}$. One can then decompose the space $\mathcal{M}_w(\Gamma_1(N))$ in terms of modular forms for the larger group $\Gamma_0(N)$, at the expense of introducing Dirichlet characters:
\beq
\mathcal{M}_w(\Gamma_1(N))=\bigoplus_\chi \mathcal{M}_w(\Gamma_0(N), \chi),
\eeq
where functions in the $\chi$-eigenspace $\mathcal{M}_w(\Gamma_0(N), \chi)$ obey a generalisation of (\ref{holmodform}):
\beq
f\left(\frac{az + d}{cz+d}\right) = \chi(d)(cz +d)^w f(z), \qquad f\in \mathcal{M}_w(\Gamma_0(N), \chi).
\eeq
Functions in $\mathcal{M}_w(\Gamma_0(N), \chi)$ are said to be of \emph{level N}. 
\begin{example}[Classical level $N$ Eisenstein series]
An example of a modular form for $\Gamma_0(N)$ is provided by the level $N$, weight $2w$ Eisenstein series 
\begin{align}
\label{holEisenLevelN}
E_{2w}(z; \chi) = \sum_{(m,n)\in\ints^2 \atop (m,n)=1} \frac{\chi(n)}{(mz+n)^{2w}},
\end{align}
generalising the classical series (\ref{holEisen}).
\end{example}

\subsection{From holomorphic modular forms to automorphic forms on \texorpdfstring{$SL(2,\mathbb{R})$}{SL(2, R)}}
\label{sec:fromholtoaut}

We shall now see how to adapt the theory of holomorphic modular forms on $\mathbb{H}=SL(2,\reals)/SO(2,\reals)$ to the more general framework of automorphic forms defined on $G=SL(2,\reals)$, invariant under the left action of $SL(2,\ints)$.

Given a weight $w$ holomorphic modular form $f : \mathbb{H}\to \mathbb{C}$ we define a new (complex)  function $\varphi_f$ on $ SL(2,\mathbb{R})$ through the assignment
\beq
f\longmapsto \varphi_f(g) = (c i + d)^{-w} f(g \cdot i),
\label{phifromf}
\eeq
where $g=\left(\begin{smallmatrix} a & b \\ c & d \\ \end{smallmatrix}\right) \in SL(2,\mathbb{R})$. The prefactor here is chosen in such away as to cancel the factor of automorphy in~\eqref{holmodform} in order for the function $\varphi_f$ to be \emph{invariant} under $SL(2,\ints)$:
\beq
\varphi_f(\gamma g)=\varphi_f(g), \qquad \gamma \in SL(2,\mathbb{Z}).
\eeq
According to our definition in section \ref{sec:AutoIntro}, $\varphi_f$ is thus an automorphic function on $SL(2,\reals)$. Note that the condition of moderate growth is satisfied automatically since the seed function $f$ is holomorphic. 

We can make the \emphindex[lift!from $\mathbb{H}$ to $SL(2,\reals)$]{lift~\eqref{phifromf} from $\mathbb{H}$ to $SL(2,\reals)$} more explicit by making use of the Iwasawa decomposition of an element $g\in SL(2,\mathbb{R})$:
\beq
g=nak = \begin{pmatrix} 1 & x \\& 1 \end{pmatrix} \begin{pmatrix} y^{1/2} & \\ & y^{-1/2} \end{pmatrix} \begin{pmatrix} \cos \theta & \sin \theta \\ -\sin\theta & \cos \theta  \end{pmatrix}, 
\label{IwasawaSL2R1}
\eeq
with $n\in N(\mathbb{R}), a\in A(\mathbb{R}), k\in SO(2,\mathbb{R})$. Acting with $g$ on the point $i$ one finds then
\begin{align}
g \cdot i = x+iy \equiv z.
\end{align}
We recall that $K=SO(2,\mathbb{R})$ leaves the point $i$ invariant. Plugging the Iwasawa decomposition of $g$ into the right-hand side of~\eqref{phifromf} we can write $\varphi_f$ as a function of the three variables $(x,y,\theta)$:
\beq
\label{varphi}
\varphi_f(g)=\varphi_f(x, y, \theta)=e^{iw\theta} y^{w/2} f(x+iy).
\eeq

Moreover, under the right-action of 
\beq
k=\begin{pmatrix} \cos \vartheta & \sin \vartheta \\ -\sin\vartheta & \cos \vartheta  \end{pmatrix} \in SO(2,\mathbb{R})
\eeq
 it transforms by a phase: 
 \beq
 \varphi_f(gk)= e^{iw\vartheta}\varphi_f(g).
 \label{varphiktransf}
 \eeq
This implies that the original transformation property (\ref{holmodform}) of $f$ under $SL(2,\mathbb{Z})$ has been traded for the above phase transformation of $\varphi_f(g)$ under $K=SO(2,\mathbb{R})$. While $f$ itself was invariant under $SO(2,\mathbb{R})$ one instead says that $\varphi_f$ is \emphindex[K-finite@$K$-finite]{$K$-finite}, implying that the action of $K$ on $f$ generates a finite-dimensional vector space; in the present example this is represented by the one-dimensional space of characters $\sigma : k\mapsto  e^{iw\vartheta}$ through 
\begin{align}
\label{modphik}
\varphi_f(gk)=\sigma(k) \varphi_f(g).
\end{align}

Next, we address the question how the automorphic form $\varphi_f$ incorporates the holomorphy of $f$ on $\mathbb{H}$:
\beq
\frac{\partial}{\partial \bar{z}} f = \frac{1}{2}\left(\frac{\partial}{\partial x} +i \frac{\partial}{\partial y}\right)f=0,
\label{holomorphy}
\eeq
where $z=x+iy$. The corresponding statement for $\varphi_f$ is that it satisfies
\beq
F\varphi_f = -2i  e^{-2i\theta} \left(y\frac{\partial}{\partial \bar{z}} -\frac{1}{4} \frac{\partial}{\partial \theta}\right) \varphi_f=0.
\label{Fdiffop}
\eeq
We will now give a group-theoretic interpretation to this differential condition. 

The group $SL(2,\reals)$ acts on smooth functions on $SL(2,\reals)$ via the \emphindex[action!right regular on functions]{right-regular action}, where we recall that `regular' here refers to the fact that the \emphindex{Radon--Nikodym derivative} is trivial. Let $g'$ be an element of $SL(2,\reals)$ and $\varphi(g)$ a function on $SL(2,\reals)$. The right-regular action is defined by:
\begin{align}
\label{eq:rra}
\Big(\pi(g') \varphi\Big) (g) = \varphi(gg').
\end{align}
The action of the Lie algebra $\mf{sl}(2,\reals)$ is then given by differential operators acting on \emp{smooth} functions. Using~\eqref{eq:rra} one finds the following differential operators corresponding to the Chevalley basis generators:
\begin{subequations}
\begin{align}
h &= -2\sin(2\theta) y \partial_x+ 2\cos(2\theta) y \partial_y + \sin(2\theta) \partial_\theta,\nn\\
e &= \cos(2\theta) y \partial_x +  \sin(2\theta) y \partial_y + \sin^2\theta \partial_\theta,\nn\\
f&= \cos(2\theta) y \partial_x +  \sin(2\theta)y\partial_y -\cos^2\theta \partial_\theta. 
\end{align}
\end{subequations}
The compact generator $e-f$ of $SO(2,\reals)$ acts by $\partial_\theta$.  We record also the inverse relations
\begin{subequations}
\begin{align}
y \partial_x &= \frac{1}{2} ((e+f) \cos (2 \theta )+e-f-h \sin (2 \theta )),\\
y \partial_y &=\frac{1}{2} ((e+f)  \sin (2 \theta )+h \cos (2 \theta )),\\
\partial_\theta &=e-f.
\end{align}
\end{subequations}

The Casimir operator~\eqref{eq:SL2cas} then becomes a second order differential operator, namely the Laplacian
\begin{align}
\label{eq:SL2Lapapp}
\Delta_{SL(2,\mathbb{R})} = y^2 \left( \partial_x^2 + \partial_y^2 \right) - y \partial_x \partial_\theta.
\end{align}

We shall also make use of the so-called \emphindex[compact basis!for $\mf{sl}(2,\reals)$]{compact basis}. This is a representation of $\mf{sl}(2,\reals)$ in terms of $(2\times 2)$-matrices different from~\eqref{eq:sl2basis} and given explicitly by
\begin{align}
H=-i(e-f),\quad
E= \frac12\left( h +i (e+f) \right),\quad
F= \frac12\left( h -i (e+f) \right),
\end{align}
that is,
\begin{align}
\label{eq:sl2cpt}
H = i\begin{pmatrix}0&-1\\1&0\end{pmatrix},\quad
E = \frac12\begin{pmatrix}1&i\\i&-1\end{pmatrix},\quad
F = \frac12\begin{pmatrix}1&-i\\-i&-1\end{pmatrix}.
\end{align}
The generators satisfy the standard $\mf{sl}(2,\reals)$ commutation relations
\begin{align}
\left[ H, E \right] =2E,\quad
\left[ H, F \right] =-2F,\quad
\left[ E,F \right] =H.
\end{align}
The Cartan generator $H$ is Hermitian in this basis and this is the reason for the name compact basis.

The representation~\eqref{eq:sl2cpt} of $\mf{sl}(2,\reals)$ is unitarily equivalent to the standard representation~\eqref{eq:sl2basis} through the transformation
\begin{align}
U H U^\dagger = h, \textrm{etc.} \quad\quad
\textrm{for}\quad
U=\frac12
\begin{pmatrix} -1+i &1+i \\
-1+i & -1-i
\end{pmatrix}.
\end{align}
The differential operators associated with this basis are then given by
\begin{subequations}
\begin{align}
H &= -i\partial_\theta,\\
E &= 2i e^{2i \theta} \left(y \partial_z -\frac{1}{4} \partial_\theta\right),\\
F &= -2i e^{-2i \theta} \left(y \partial_{\bar{z}} -\frac{1}{4} \partial_\theta\right),
\end{align}
\end{subequations}
where we have used standard holomorphic and anti-holormorphic derivatives:
\begin{align}
\partial_z = \frac{1}{2} \left(\partial_{x} - i \partial_{y}\right),\quad
\partial_{\bar{z}} = \frac{1}{2} \left(\partial_{x} + i \partial_{y}\right).
\end{align}
Because the compact basis is unitarily equivalent, the Casimir operator does not change.

\begin{remark}
The change of basis is basically that induced by the \emphindex{Sekiguchi isomorphism}~\cite{Sekiguchi,KostantRallis} that enters in the description of real nilpotent orbits.
\end{remark}

The above discussion implies that  the differential operator $F$ in~\eqref{Fdiffop} may in fact be identified with the lowering operator in the basis $(E, F, H)$ of the Lie algebra  $\mathfrak{sl}(2,\mathbb{R})$. This implies that $\varphi_f$ may be viewed as a lowest weight state of a representation of $\mathfrak{sl}(2,\mathbb{R})$.  Furthermore, we note that the (Hermitian) generator $H$ in this basis corresponds to
\beq
e^{i \theta H}=\begin{pmatrix} \cos \theta & \sin \theta \\ -\sin\theta & \cos \theta  \end{pmatrix}  \in SO(2,\mathbb{R}),
\label{Hkrelation}
\eeq
and hence corresponds to the differential operator $H=-i\partial_\theta$. Therefore, $H$ is diagonal on $\varphi_f$ with eigenvalue $w$:
\beq 
H\varphi_f = w \varphi_f.
\eeq
{}From the commutation relations we further deduce that $E$ raises the $H$-eigenvalue $w$ by $+2$,  while  $F$ lowers it by the same amount. This implies that the holomorphic Eisenstein series $E_{2w}$ can be viewed as lowest weight vectors in the holomorphic discrete series of $SL(2,\mathbb{R})$, providing our first glimpse of the general connection between automorphic forms and representation theory, a topic that will be discussed in more generality in chapter \ref{ch:autoreps} and onwards. See also section~\ref{standardsectionSL2} for some more details on the specific case of $SL(2)$ treated above.

Before we proceed we shall mention one final important property of $\varphi_f$, namely that it is an eigenfunction of the Laplacian \eqref{eq:SL2Lapapp} on $SL(2,\mathbb{R})$. Acting on $\varphi_f$ one obtains that
\beq
\Delta_{SL(2,\mathbb{R})} \varphi_f = \frac{w}{2}\left(\frac{w}{2}-1\right) \varphi_f.
\label{Laplacevarphif}
\eeq
As we will see, all the properties of $\varphi_f$ discussed above will  have counterparts in the general theory of automorphic forms.  

The automorphic lift of weight $w$, level $N$ modular forms for $\Gamma_0(N)$ will be treated in section~\ref{ex_adelisationholf}. 

\begin{example}[Lift of a holomorphic Eisenstein series]
\label{LiftholEis}
The lift of the holomorphic Eisenstein series $f(z) = E_{2w}(z)$ to an automorphic form on $SL(2,\mathbb{R})$ can be written using~\eqref{varphi} as
\begin{align}
\varphi_f(g) = e^{2iw\theta} y^w E_{2w}(x+iy).
\end{align}
Another common way of writing the lift is by using the co-called \emphindex{slash operator} that is defined as follows. For $f$ a weight $w$ holomorphic modular form, $g\in SL(2,\mathbb{R})$, define the \emph{slash operator} $f|_wg  : \mathbb{H} \to \mathbb{C}$ by
\beq
\left(f|_w g\right)(z):=(cz+d)^{-w}f\left(\frac{az+b}{cz+d}\right), \qquad g=\begin{pmatrix} a & b \\ c & d \\ \end{pmatrix} \in SL(2,\mathbb{R}).
\eeq
Using the slash operator the defining relation (\ref{holmodform}) can be written simply as 
\beq
\left(f|_w \gamma\right)(z)=f(z), \qquad \gamma=\begin{pmatrix} a & b \\ c & d \\ \end{pmatrix} \in SL(2,\mathbb{Z}).
\eeq
By a calculation similar to the one used in proving equation (\ref{cosetSL2}) we can now rewrite the Eisenstein series $E_{2w}(z)$ in (\ref{holEisen}) directly as a function on $SL(2,\mathbb{R})$. Parametrising an element $g\in SL(2,\mathbb{R})$ in Iwasawa form as $g=nak$ (see~\ref{IwasawaSL2R1}) we obtain
\beq
\varphi_f(nak)=\left(f|_w nak\right)(i)=e^{iw\theta}y^{w/2} \, \sum_{\gamma\in N(\mathbb{Z})\backslash SL(2,\mathbb{Z})} \left(1|_w \gamma\right)(nak \cdot i).
\eeq
\end{example}
 
\begin{remark}
    We would like to make a cautionary remark regarding the generalisation of the above discussion to arbitrary groups $G(\reals)$. Modular forms are holomorphic functions on  $\UHP\cong SL(2,\mathbb{R})/SO(2,\mathbb{R})$ with simple transformation properties under $SL(2,\mathbb{Z})$, and it seems natural to try and generalise this construction to higher rank real Lie groups $G(\mathbb{R})$. One might suspect a generalisation to holomorphic functions  $f\, :\, G(\mathbb{R})/K \, \to \mathbb{C}$, where $K$ is the maximal compact subgroup of $G(\mathbb{R}$), transforming with some weight under the action of a discrete subgroup $G(\mathbb{Z})\subset G(\mathbb{R})$. However, this only works whenever the coset $G(\mathbb{R})/K$ carries a complex structure. In the case above this complex structure is provided by the fact that the maximal compact subgroup $K=SO(2,\mathbb{R})\cong U(1)$. In general, the maximal subgroup $K$ of some $G(\reals)$ does not have an isolated $U(1)$ factor that can provide a complex structure on $G(\reals)/K$ and therefore we could not expect to have a general theory of holomorphic modular forms on $G/K$.  A standard example with a complex structure is provided by $G=Sp(2n, \mathbb{R}), \, K=U(n),$ in which case $Sp(2n;\mathbb{R})/U(n)$ is a hermitian symmetric domain known as the Siegel upper half space. This leads to the theory of holomorphic Siegel modular forms (see, e.g., \cite{deGeer} for a review). 
\end{remark}

 \subsection{Maass forms and non-holomorphic Eisenstein series}
 \label{Maass}

As just discussed, it is in general too restrictive (and often impossible) to consider holomorphic modular forms. It is therefore called for to look for a theory of arbitrary (\emp{non-holomorphic})  functions $f \, :\, G(\mathbb{R})/K \to \mathbb{C}$ which  transform nicely under the action of some discrete subgroup $G(\mathbb{Z})\subset G(\mathbb{R})$. This leads to the  notion of an \emphindex{automorphic form} that we will now discuss for $SL(2,\reals)$. 
 
In addition to the holomorphic modular forms, the classical theory also contains an interesting class of \emp{non-holomorphic functions} $f : SL(2,\mathbb{Z})\backslash \mathbb{H} \to \mathbb{R}$ . These non-holomorphic functions are eigenfunctions of the Laplacian $\Delta_{\mathbb{H}}$ on $\mathbb{H}=SL(2,\mathbb{R})/SO(2,\mathbb{R})$ (that is simply obtained from~\eqref{eq:SL2Lapapp} since $\partial_\theta=0$ on $\UHP$):
\beq
\label{LapSL2}
\Delta_{\mathbb{H}} f = y^2\left(\frac{\partial^2}{\partial x^2}+\frac{\partial^2}{\partial y^2}\right) f=\lambda f
\eeq
and by definition are \emp{invariant} under $SL(2,\ints)$:
\begin{align}
f(\gamma\cdot z) = f(z);
\end{align}
there is no non-trivial weight $w$ compared to~\eqref{holmodform}. Similarly, to the holomorphic case we require that $f(z)$ is of moderate growth, i.e. that it grows at most polynomially for $y\to \infty$ (see (\ref{moderategrowthSL2})).

Functions on $SL(2,\mathbb{R})$ satisfying these conditions are called \emphindex[Maass form]{Maass (wave) forms}, and they can also be fit into the general framework of automorphic forms, with even less effort than for the holomorphic modular forms. Given a Maass form $f$ on $\mathbb{H}$ we lift this to a function $\varphi_f$ on $SL(2,\mathbb{R})$ according to~\eqref{varphi}
\beq
f \longmapsto \varphi_f(g)=\varphi_f \left( \begin{pmatrix} 1 & x \\& 1 \end{pmatrix} \begin{pmatrix} y^{1/2} & \\ & y^{-1/2} \end{pmatrix} k\right) =f(x+iy), 
\label{liftMaass}
\eeq
where we used the Iwasawa decomposition $g=nak\in SL(2,\mathbb{R})$  given in equation (\ref{IwasawaSL2R}). The lift in this case is trivial since $w=0$.

The associated function $\varphi_f(g)$ then satisfies 
\beq
\varphi_f(\gamma g k) = \varphi_f(g), \qquad \qquad \gamma \in SL(2,\mathbb{Z}), \quad k\in SO(2,\mathbb{R}),
\eeq
and so is indeed an automorphic form on $SL(2,\mathbb{R})$.

Important examples of Maass forms are 
provided by the \emphindex[Eisenstein series!non-holomorphic]{non-holomorphic Eisenstein series} with parameter $s\in \cx$
\beq
E(s, z) = \frac{1}{2} \sum_{(c,d)\in \mathbb{Z}^2\atop (c,d)=1} \frac{y^s}{|cz+d|^{2s}}.
\label{nonholeis}
\eeq
This converges absolutely for $\Re(s) > 1$, but according to Langlands general theory \cite{LanglandsFE} it can  be analytically continued to a meromorphic function of $s\in \mathbb{C}\backslash \{0,1\}$. This crucial fact relies on the \emph{functional relation}\index{functional relation!for SL(2,R)@for $SL(2,\reals)$}
\begin{align}
\label{SL2func}
\xi(s) E(s,z) = \xi(1-s) E(1-s,z),
\end{align}
where $\xi(s)$ is the completed Riemann zeta function (\ref{CRZ}). 

One can verify that the Eisenstein series $E(s,z)$ indeed defines an $SL(2,\mathbb{Z})$-invariant eigenfunction of the Laplacian $\Delta_{\mathbb{H}}$ with eigenvalue $\lambda=s(s-1)$. The non-holomorphic Eisenstein series $E(s, z)$ provides the simplest example of the class of Eisenstein series on a group $G(\reals)$ that will be our main concern in the following. 
 
It is instructive to rewrite $E(s,z)$ as defined in~\eqref{nonholeis}. We parametrise an arbitrary group element $g\in SL(2,\mathbb{R})$ according to the same Iwasawa decomposition $g=nak$ as in (\ref{IwasawaSL2R1}). Then introduce a character $\chi_s : B=NA \to \mathbb{C}^{\times}$ defined by
\beq
\chi_s(na)=y^s, \qquad \qquad n\in N, \, a \in A,
\eeq
and extend it to all of $SL(2,\mathbb{R})$ by requiring it to be trivial on $SO(2,\mathbb{R})$: $\chi_s(nak)=\chi_s(na)$. The Eisenstein series $E(s,z)$ can now be equivalently 
written as
\beq
E(s,g) = \sum_{\gamma \in B(\mathbb{Z})\backslash SL(2,\mathbb{Z})} \chi_s(\gamma  g),
\label{PoincareSL2}
\eeq
where the quotient by the discrete Borel subgroup $B(\mathbb{Z})=\{ \left(\begin{smallmatrix} \pm1 & m \\& \pm1 \end{smallmatrix}\right) \, |\, m\in \mathbb{Z}\}$ is needed since it leaves $\chi_s$ invariant. (This was also explained in the introduction.) It should be apparent that this reformulation of the Eisenstein series is well suited for generalisations to higher rank Lie groups $G(\mathbb{R})$. This will be discussed in section \ref{sec_Eisenstein} below.

\subsection{Maass forms of non-zero weight*}
\label{sec:Maasswt}

One can generalise the definition of Maass form given above to include non-holomorphic functions which transform with a weight. We define a \emph{weight $w$ Maass form} to be a non-holomorphic function $f:\mathbb{H}\to \mathbb{C}$ satisfying 
\beq
f\left(\frac{az + d}{cz+d}\right) =\left(\frac{cz+d}{|cz+d|}\right)^w f(z), \qquad w\in \mathbb{Z}.
\eeq
A weight $w$ Maass form is furthermore an eigenfunction of the weight $w$ Laplacian $\Delta_w$ which is a modification of (\ref{LapSL2}): 
\beq 
\Delta_{\mathbb{H}}^{w}:= y^2\left(\frac{\partial^2}{\partial x^2}+\frac{\partial^2}{\partial y^2}\right) -i w y \frac{\partial}{\partial x}.
\label{weightLaplacian}
\eeq
We can elucidate the meaning of this differential operator by lifting the weight $w$ Maass form $f$ to an automorphic form $\varphi_f$ on $SL(2,\mathbb{R})$ through a straightforward generalisation of (\ref{liftMaass}):
\beq
f\longmapsto \varphi_f(g):= \left(\frac{ci+d}{|ci+d|}\right)^{-w} f(g\cdot i), \qquad g=\begin{pmatrix} a & b \\ c & d \\ \end{pmatrix} \in SL(2,\mathbb{R}).
\eeq
Rewriting this in Iwasawa form (\ref{IwasawaSL2R1}) yields 
\beq
\varphi_f(g)=\varphi_f \left( \begin{pmatrix} 1 & x \\& 1 \end{pmatrix} \begin{pmatrix} y^{1/2} & \\ & y^{-1/2} \end{pmatrix} k\right)=e^{iw\theta} f(x+iy). 
\eeq
We then recognise the weight $w$ Laplacian $\Delta_{\mathbb{H}}^{w}$ in (\ref{weightLaplacian}) as nothing but the full Laplacian on $SL(2, \mathbb{R})$ (\ref{eq:SL2Lapapp}) after evaluating the 
derivative on $\theta$:
\beq
\Delta_{SL(2,\mathbb{R})}\varphi_f(g)=e^{iw\theta}\left[y^2\left(\frac{\partial^2}{\partial x^2}+\frac{\partial^2}{\partial y^2}\right) -i w y \frac{\partial}{\partial x}\right]f(x+iy)=e^{iw\theta} \Delta_{\mathbb{H}}^{w} f(z).
\eeq

\begin{example}[Non-holomorphic Eisenstein series of weight $w$]
A classic example of a weight $w$ Maass form is the following generalisation of the non-holomorphic Eisenstein series (\ref{nonholeis}):
\beq
E_w(s, z)=\sum_{(c,d)\in \mathbb{Z}^2\atop (c,d)=1} \frac{y^s}{|cz+d|^{2s}}\left(\frac{cz+d}{c\bar{z}+d}\right)^{w},
\eeq
which transforms as 
\beq
E_w\left(s, \frac{az+b}{cz+d}\right)= \left(\frac{cz+d}{c\bar{z}+d}\right)^{w/2} E_w(s, z)=\left(\frac{cz+d}{|cz+d|}\right)^{w} E_w(s, z).
\eeq
We will come back to this Eisenstein series in section \ref{standardsectionSL2}.
\end{example}

\subsection{Adelisation of non-holomorphic Eisenstein series} 
\label{sec:Eisenstein-adelisation}

As we discussed in section \ref{sec_adelisation}, strong approximation ensures that we can always lift a function on $SL(2, \mathbb{Z})\backslash SL(2,\mathbb{R})$ to an adelic function on $SL(2,\mathbb{Q})\backslash SL(2,\mathbb{A})$, where the role of the discrete subgroup is now played by $SL(2,\mathbb{Q})$. Recall that this lift also requires that the resulting function is right-invariant under $K_f=\prod_{p< \infty} SL(2,\mathbb{Z}_p)$. It is now a simple matter to generalise the Eisenstein series $E(s,g)$ to such an adelic function. First extend the definition of $\chi_s$ to a function $\chi_s : B(\mathbb{A}) \to \mathbb{C} $, which is invariant under the left action of $B(\mathbb{Q})$. We extend it to all of $SL(2,\mathbb{A})$ using the global Iwasawa decomposition $SL(2,\mathbb{A}) = B(\mathbb{A})K_{\mathbb{A}} $ and demanding it to be trivial on $K_{\mathbb{A}} =SO(2,\mathbb{R}) \times K_f$. Note that this automatically takes care of the required condition of $K_f$-invariance on the right. The adelic Eisenstein series then takes the form 
\beq
E(s, g_\mathbb{A})=\sum_{\gamma \in B(\mathbb{Q})\backslash SL(2,\mathbb{Q})} \chi_s(\gamma g_\mathbb{A}),
\label{adelic_EisensteinSL2}
\eeq
which is a function on $SL(2,\mathbb{A})$ satisfying
\beq
E(s, \gamma g_\mathbb{A} k_\mathbb{A})=E(s, g_\mathbb{A}), \qquad \qquad \gamma\in SL(2,\mathbb{Q}), \quad k_\mathbb{A}\in K_{\mathbb{A}}.
\eeq
As shown in example~\ref{bijection}, the range of the sum in (\ref{adelic_EisensteinSL2}) is in fact in bijection with the range of summation in (\ref{PoincareSL2}):
\beq
B(\mathbb{Q})\backslash SL(2,\mathbb{Q})\cong B(\mathbb{Z})\backslash SL(2,\mathbb{Z}).
\label{Borelbijection}
\eeq
Therefore, if we restrict to elements $g_\mathbb{A} = (g_\infty;  1 , 1 ,\dots) \in SL(2,\mathbb{A})$, with $g_\infty \in SL(2,\mathbb{R})$, then the adelic Eisenstein series reduces  to the real Eisenstein series (\ref{PoincareSL2}). More details of this procedure can be found in section \ref{standardsectionSL2}.

With a little more effort one can also obtain the adelisation of the function $\varphi_f(g)$ in (\ref{phifromf}), for $f$ a weight $w$ holomorphic modular form on $\mathbb{H}$. This analysis is done in example~\ref{ex_adelisationholf}. 

Even though our main interest often lies with automorphic forms on real Lie groups $G(\mathbb{R})$, the adelic reformulation turns out to be extremely convenient for 
many purposes, not the least of which being the calculational advantages that it brings when computing Fourier expansions of automorphic forms, a topic which we will be concerned with in the second half of this treatise.

\section{Adelic automorphic forms}
\label{sec_defauto}

We shall now treat automorphic forms in the adelic framework. There are various degrees of generality here; for instance, one can define the theory of automorphic forms over the adeles $\mathbb{A}_\field$ of an arbitrary number field $\field$. However, we shall continue to assume that $\field=\mathbb{Q}$ in what follows. The framework of adelic automorphic forms was originally developed in the books by Gelfand--Graev--Piatetski-Shapiro \cite{GGPS}, and Jacquet--Langlands \cite{JL}. Good introductions can be found in the books by Gelbart \cite{Gelbart}, Bump \cite{Bump} and Goldfeld-Hundley \cite{MR2807433,MR2808915}.

\subsection{Main definition}

In what follows we let $G$ be a split algebraic group defined over $\mathbb{Q}$ and $G(\mathbb{A})$ its adelisation as in section~\ref{sec_adelisation}.  
The typical example we have in mind is $G(\mathbb{A})=SL(n, \mathbb{A})$. Let us now state our definition of an automorphic form:

\begin{definition}[Automorphic form]
\label{defauto} An \emphindex[automorphic form|textbf]{automorphic form} is a smooth function $\varphi \, :\, G(\mathbb{Q})\backslash G(\mathbb{A})\to \mathbb{C}$ 
satisfying the following conditions:
\begin{enumerate}
\item {\it left $G(\mathbb{Q})$-invariance:} $\varphi(\gamma g)=\varphi(g), \qquad \gamma \in G(\mathbb{Q})$,
\item {\it right $K$-finiteness:} $\text{dim}_{\mathbb{C}}\langle \varphi(gk)|k\in K_\mathbb{A}\rangle<\infty$,
\item {\it $\mathcal{Z}(\mathfrak{g_{\mathbb{C}}})$-finiteness:} $\text{dim}\langle X\varphi(g)|X\in\mathcal{Z}(\mathfrak{g}_{\mathbb{C}})\rangle<\infty$,
\item {\it $\varphi$ is of moderate growth:} for any norm $||\cdot ||$ on $G(\mathbb{A})$ there exists a positive integer $n$ and a constant $C$ such that 
$|\varphi(g)| \leq C ||g||^n$.
\end{enumerate}
\end{definition}

\begin{remark}
We denote by $\mathcal{A}(G(\mathbb{Q})\backslash G(\mathbb{A}))$ the space of automorphic forms satisfying definition \ref{defauto}. This  is a subspace of the space $C^{\infty}( G(\mathbb{Q})\backslash G(\mathbb{A}))$ of~\emph{smooth} functions on $G(\mathbb{Q})\backslash G(\mathbb{A})$. An adelic function $\varphi(g)$, with $g=(g_\infty; g_f)\in G(\mathbb{A})=G(\mathbb{R}) \times G(\mathbb{A}_f)$, is said to be smooth if it is $C ^{\infty}$ with respect to the archimedean variables $g_\infty \in G(\mathbb{R})$ and locally constant with respect to the non-archimedean variables $g_f\in G(\mathbb{A}_f)$.
\end{remark}

Let us now elaborate a little on the definition \ref{defauto} of an adelic automorphic form:
\begin{itemize}
\item Condition $(1)$ ensues as a straightforward generalisation of invariance of the function under a discrete subgroup of $G(\mathbb{A})$.
\item The condition of right $K$-finiteness means that the vector space $V$ spanned by the functions $k\mapsto \varphi(gk)$, $k\in K_{\mathbb{A}}$, is finite-dimensional. We have already seen an example of a non-trivial $K$ representation in~\eqref{modphik}. When $\varphi$ is $K$-invariant on the right, it would be more appropriate to refer to $\varphi$ as an automorphic function rather than form but we will use the more general term. 
\item In condition $(3)$, $\mathfrak{g}$ is the Lie algebra associated with the group $G$ and $\mathcal{Z}(\mathfrak{g})$ is the center of its universal enveloping algebra $\mathcal{U}(\mathfrak{g}_\mathbb{C})$. The center $\mathcal{Z}(\mathfrak{g})$ is the space of bi-invariant differential operators on $G$, i.e. the quadratic Casimir and  higher-order operators.
The condition of $\mathcal{Z}(\mathfrak{g})$-finiteness then implies that $\varphi$ is contained in a $\mathcal{Z}(\mathfrak{g})$-invariant finite-dimensional subspace of $C^{\infty}( G(\mathbb{Q})\backslash G(\mathbb{A}))$. Equivalently, if $X\in \mathcal{Z}(\mathfrak{g})$ then $\mathcal{Z}(\mathfrak{g})$-finiteness implies that there exists a polynomial $R(X)$ such that $R(X)\varphi=0$. 
\end{itemize}

\begin{remark}\label{Condition5} It is sometimes useful to specify the transformation properties of an automorphic form with respect to the 
center $Z(\mathbb{A})$ of $G(\mathbb{A})$. To this end, let $\omega$ be a \emphindex{central character}, i.e. a  homomorphism $\omega  : Z(\mathbb{A}) \to \mathbb{C}^{\times}$, 
which is trivial on $Z(\mathbb{Q})$. An automorphic form $f$ is then said to have central character $\omega$ if it satisfies conditions (1)--(4) along with the additional condition
\begin{enumerate}
\item[5.] $f(zg) =\omega(z)f(g)$.
\end{enumerate}
\end{remark}

We shall now illustrate these defining properties of an automorphic form, by giving two examples. First we look at the non-holomorphic Eisenstein series~\eqref{adelic_EisensteinSL2} and we will verify its properties according to the above definition.

\begin{example}[Verification of automorphic properties of an Eisenstein series on $SL(2,\mathbb{A})$]
 \label{ex_definingrelations} Consider the case when $G(\mathbb{A})=SL(2,\mathbb{A})$. We now verify the conditions (1)--(4) in definition \ref{defauto} for  the non-holomorphic Eisenstein series defined in (\ref{adelic_EisensteinSL2}):
 \beq
 E(s, g)=\sum_{\gamma \in B(\mathbb{Q})\backslash SL(2,\mathbb{Q})} \chi_s(\gamma g), \qquad g\in SL(2,\mathbb{A}).
 \eeq
 \begin{itemize} 
\item By construction, $E(s, g)$ is left $SL(2,\mathbb{Q})$-invariant and so satisfies condition $(1)$.
\item Moreover, by definition the function $\chi_s$ is invariant under any $k_p\in SL(2,\mathbb{Z}_p)$, $\chi_s(gk_p)=\chi_s(g)$, and hence condition $(2)$  is also satisfied. 
\item To understand condition $(3)$ we recall that $E(s, g)$ is an eigenfunction of the $\mathfrak{g}$-invariant Laplacian $\Delta_{\mathbb{H}}$ on $SL(2,\mathbb{R})/SO(2,\mathbb{R})$ with eigenvalue $\lambda = s(s-1)$. Hence, $E(s, g)$ is in the kernel of the operator $(\Delta_{\mathbb{H}} -\lambda)\in \mathcal{Z}(\mathfrak{g})$, and since for $\mathfrak{g}=\mathfrak{sl}(2,\mathbb{A})$, $\mathcal{Z}(\mathfrak{g})=\mathbb{C}[\Delta_{\mathbb{H}}]$, we have that condition $(3)$ is satisfied. 
\item The final part consists in verifying the moderate growth condition $(4)$. To this end one must translate the classical moderate growth condition (\ref{moderategrowthSL2}) to the adelic picture. A norm $||\, ||$ on $SL(2,\mathbb{A})$ can be defined as follows (see, e.g., \cite{MR2807433}): 
\beq
||g||:= \prod_{p\leq \infty} \text{max}\left\{ |a|_p, |b|_p, |c|_p, |d|_p, |ad-bc|_p^{-1}\right\}, \qquad g=\begin{pmatrix} a & b \\ c & d \\ \end{pmatrix},
\eeq
where it is understood that $|a|_p=|a_p|_p$ etc., with $a_p$ the $p$:th component of the adele $a=(a_\infty, a_2, a_3, \dots)\in \mathbb{A}$. For a proof that the moderate growth condition of $E(s,g)$ with respect to this norm follows from the classical moderate growth on $SL(2,\mathbb{R})$, see p. 122-123 of \cite{MR2807433}.
\end{itemize} 
\end{example}
Before we move on to analysing Eisenstein series on arbitrary reductive groups we shall give an additional important definition: 

\begin{definition}[Cusp form]
\label{autodef}An automorphic form $f\in \mathcal{A}(G(\mathbb{Q})\backslash G(\mathbb{A}))$ is a \emphindex{cusp form} if for all parabolic subgroups $P(\mathbb{A})\subset G(\mathbb{A})$ it satisfies
\beq
\int_{U(\mathbb{Q})\backslash U(\mathbb{A})} f(ug) du=0,
\label{cuspcond}
\eeq
where $U$ is the unipotent radical in the Levi decomposition $P(\mathbb{A})=L(\mathbb{A})U(\mathbb{A})$, and 
$du$ is the left-invariant Haar measure on $U$. The subspace of cusp forms will be denoted by $\mathcal{A}_0(G(\mathbb{Q})\backslash G(\mathbb{A}))\subset \mathcal{A}(G(\mathbb{Q})\backslash G(\mathbb{A}))$.
\end{definition}

\begin{remark}
This definition generalises the notion of cusp form found in the classical theory, namely holomorphic modular forms $f(z)$ whose Fourier expansion in $q=e^{2\pi iz}$ contains no term of order $q^0$. An example is provided by Ramanujan's discriminant $\Delta(z)$ of weight $w=12$.
\end{remark}

\begin{remark}
The integral in (\ref{cuspcond}) can be thought of as the zeroth Fourier coefficient of $f(g)$ with respect to $U$; by analogy with the classical theory it is  called the ``constant term'' of $f(g)$, although in general it is by no means constant. From this perspective a cusp form is simply an automorphic form with vanishing constant term. Constant terms are analysed in detail for $SL(2,\mathbb{A})$ in chapter~\ref{ch:SL2-fourier} and in full generality in chapter~\ref{ch:CTF}. 
\end{remark}

\subsection{Adelic lift with Hecke character*}
 \label{ex_adelisationholf} We shall now construct the adelic lift of a holomorphic modular form $f$. To illustrate the power of the adelic formalism we will consider the general case addressed in section~\ref{LevelNforms}, namely let $f\in \mathcal{M}_w(\Gamma_0(N), \chi)$, i.e a level $N$ holomorphic modular form for $\Gamma_0(N)$ with Dirichlet character $\chi$. We can now use strong approximation (see sec.~\ref{sec_adelisation}) to lift $f$ to a function on $SL(2,\mathbb{A})$. Recall from section~\ref{sec_adelisation} (see in particular example~\ref{strongapproxcongruence}) that strong approximation implies that any $g\in SL(2,\mathbb{A})$ can be (non-uniquely) written as 
\beq
g=\gamma g_\infty k_f, \qquad \gamma\in SL(2,\mathbb{Q}), \, g_\infty \in SL(2,\mathbb{R}), \, k_f \in K_0(N),
\eeq
where $K_0(N)\subset K_f=\prod_{p<\infty} SL(2,\mathbb{Z}_p)$ was defined in example~\ref{strongapproxcongruence}. In order to define a lift to $SL(2,\mathbb{A})$ we must first lift the Dirichlet character $\chi$ to the adelic setting. This can be done using the $GL(1, \mathbb{A})=\mathbb{A}^\times$-version of strong approximation:
\beq
\mathbb{A}^{\times}=\mathbb{Q}^{\times}\, \mathbb{R}_+ \prod_{p<\infty}\mathbb{Z}_p^\times.
\eeq 
This implies that any Dirichlet character $\chi : (\mathbb{Z}/N\mathbb{Z})^{*} \to \mathbb{C}^{\times}$ has a canonical lift to an adelic (Hecke) character 
\beq 
\omega_\chi : \mathbb{Q}^{\times}\backslash \mathbb{A}^{\times} \to \mathbb{C}^{\times}.
\eeq
 Indeed, such a character can be decomposed as 
 \beq
 \omega_\chi=\omega_{\chi, \infty} \, \prod_{p<\infty} \omega_{\chi, p},
 \label{factorisationHecke}
 \eeq
 where the archimedean factor $\omega_{\chi, \infty}$ is taken to be trivial, and each local factor $\omega_{\chi, p}$ equals the Hecke character $\chi : (\mathbb{Z}/N\mathbb{Z})^{\times}\to \mathbb{C}^{\times}$ for $N$ a power of the prime $p$.
 
Next, we lift the local character $\omega_{\chi, p}$ to a character on $SL(2,\mathbb{Z}_p)$ via the map $\left(\begin{smallmatrix} a & b \\ c & d \end{smallmatrix} \right) \mapsto \omega_{\chi, p}(d)$. The adelic lift of the holomorphic modular form $f$ is then defined by:
\beq
\varphi_f(g):=(ci+d)^{-w} f(g_\infty \cdot i)\omega_{\chi}(k_f), 
\label{adelicliftholform}
\eeq
where $g_\infty=\left(\begin{smallmatrix} a & b \\ c & d \\ \end{smallmatrix}\right) \in SL(2,\mathbb{R})$. We can also write this in terms of the slash operator used in example~\ref{LiftholEis}:
\beq
\varphi_f(g)=\left(f|_w g_\infty\right)(i)\, \omega_\chi(k_f).
\eeq

Having extended the definition of $\varphi_f$ to an adelic automorphic form we wish to verify the conditions (1)--(4) of the definition:
\begin{itemize}
\item  Condition $(1)$ is  satisfied by construction: $\varphi_f(\gamma g)=\varphi_f(g)$, for any $g\in SL(2,\mathbb{A})$ and $\gamma\in SL(2,\mathbb{Q})$.
\item Condition $(2)$, concerning right $K$-finiteness can be seen as follows. Finiteness under the non-archimedean $K_f$ is a consequence of the relation 
\beq 
\varphi_f(g k_f)=\varphi_f(g) \omega_\chi(k_f)
\eeq
while at the archimedean place we have  
\beq
\varphi_f(g k_\infty)=\varphi_f(g) e^{iw\theta},
\eeq
 where $k_\infty=k_\infty(\theta) \in SO(2,\mathbb{R})$ as in (\ref{Hkrelation}). 
\item $\mathcal{Z}(\mathfrak{g})$-finiteness (Condition $3$) again follows from the fact that $\varphi_f$ is an eigenfunction of the Laplacian:
\beq
\Delta \varphi_f = \frac{w}{2}\left(\frac{w}{2}-1\right) \varphi_f.
\eeq
\item Finally the condition of moderate growth (Condition $4$) is satisfied if the coefficients $a(n)$ in the $q$-expansion of $f(z)$ satisfy $a(n)=0$ whenever $n<0$ which holds since $f$ is holomorphic.
\end{itemize}

Finally, we shall see that $\varphi_f$ is in fact an example of an automorphic form \emph{with central character}, as in the supplementary Condition $(5)$ mentioned in remark~\ref{Condition5}.  To this end we must first view $\varphi_f$ as a function on $GL(2,\mathbb{A})$ as opposed to $SL(2,\mathbb{A})$. The defining relation (\ref{adelicliftholform}) is still valid for $g\in GL(2,\mathbb{A})$ and conditions $(1)$--$(4)$ go through without change. Our aim is now to check how $\varphi_f(g)$ transforms under the non-trivial centre $Z(GL(2,\mathbb{A}))=\mathbb{A}^{\times}$. An element $z\in Z(GL(2,\mathbb{A}))$  can be represented by the diagonal matrix
\beq
z=\begin{pmatrix} r & \\ & r\\ \end{pmatrix}, \qquad r\in \mathbb{A}.
\eeq
Strong approximation then yields the decompositions
\beqa
g&=&\gamma g_\infty k_f, \qquad \gamma\in GL(2,\mathbb{Q}), \, g_\infty \in GL(2,\mathbb{R})^{+}, \, k_f\in K_0(N), 
\nn \\
r&=& \alpha\, r_\infty \, r_f, \qquad \alpha\in \mathbb{Q}^\times, \, r_\infty\in \mathbb{R}_+, \, r_f\in \prod_{p<\infty} \mathbb{Z}_p^{\times},
\eqa
and consequently
\beq
zg=\begin{pmatrix} \alpha & \\ & \alpha \\ \end{pmatrix} \gamma \, \begin{pmatrix} r_\infty& \\ & r_\infty \\ \end{pmatrix} g_\infty \, \begin{pmatrix} r_f & \\ & r_f \\ \end{pmatrix} k_f \in GL(2,\mathbb{A}).
\eeq
We can now proceed to calculate the action of $Z$ on the automorphic form $\varphi_f$:
\beq
\varphi_f(zg)=\left(f|_w \begin{pmatrix} r_\infty& \\ & r_\infty \\ \end{pmatrix} g_\infty\right)(i) \, \omega_\chi \left(\begin{pmatrix} r_f & \\ & r_f \\ \end{pmatrix} k_f\right).
\eeq
To evaluate this we first notice that 
\beq
\left(f|_w \begin{pmatrix} r_\infty& \\ & r_\infty \\ \end{pmatrix}\right)=f, 
\eeq
and hence 
\beq
\left(f|_w \begin{pmatrix} r_\infty& \\ & r_\infty \\ \end{pmatrix} g_\infty\right)(i)=\left(f|_w  g_\infty\right)(i).
\eeq
Using the multiplicative property of the Hecke character we further have
\beq
\omega_\chi \left(\begin{pmatrix} r_f & \\ & r_f \\ \end{pmatrix} k_f\right)=\omega_\chi \left(\begin{pmatrix} r_f & \\ & r_f \\ \end{pmatrix} \right)\, \omega_\chi(k_f).
\eeq
By definition the Hecke character $\omega_\chi$ is trivial on $\mathbb{Q}^{\times}$ and at the archimedean place. Thus, using strong approximation we can write 
\beq 
\omega_\chi(r_f)=\omega_\chi(\alpha r_\infty r_f)=\omega_\chi(z), \qquad z\in Z(GL(2,\mathbb{A}).
\eeq
Combining everything we then find 
\beq
\varphi_f(zg)=\omega_\chi(z) \varphi_f(g),
\eeq
verifying that $\varphi_f$ is an automorphic form with central character $\omega=\omega_\chi$ as in remark~\ref{Condition5}.

\section{Eisenstein series}
\label{sec_Eisenstein}

We now want to generalise the construction of adelic Eisenstein series given in section \ref{sec_frommodtoaut} to arbitrary reductive  groups $G(\mathbb{A})$. To this end we must first recall the process of constructing representations of $G$ via induction from a standard parabolic subgroup $P \supset B$. In this section we shall take $P=B$, the Borel subgroup which is the minimal parabolic subgroup. The case of arbitrary (standard) parabolic subgroups will be treated in the subsection~\ref{nonminEis}.

\subsection{Adelic multiplicative characters}
\label{sec_multcharborel}

Fix a Borel subgroup $B(\mathbb{A})\subset G(\mathbb{A})$ with Levi decomposition $B(\mathbb{A})=N(\mathbb{A})A(\mathbb{A})$. Recall that since $G(\mathbb{A})$ is split, $A(\mathbb{A})\cong (\mathbb{A}^{\times})^{\text{rank}\, \mathfrak{g}}$. Introduce a multiplicative character 
\beq
\chi \, :\, B(\mathbb{Q})\backslash B(\mathbb{A}) \to \mathbb{C}^{\times}, 
\eeq
defined by 
\beq
\chi(na)=\chi(a), \qquad \qquad n\in N(\mathbb{A}), \, a\in A(\mathbb{A}). 
\eeq
Using the Iwasawa decomposition we can extend $\chi$ to all of $G(\mathbb{A})$ by demanding it to be trivial on $K_\mathbb{A}$:
\beq
\chi(g)=\chi(nak)=\chi(na)=\chi(an)=\chi(a), \qquad k\in K_\mathbb{A}.
\eeq
Although we extend the character to all of $G(\mathbb{A})$ it is only multiplicative on  $B(\mathbb{A})$:
\beq\label{charmult}
\chi(b b')=\chi(b)\chi(b') = \chi(a)\chi(a'), \qquad b, b'\in B(\mathbb{A}).
\eeq
On the other hand, to evaluate it on a product of two elements $g, g'\in G(\mathbb{A})$ we  have
\beq 
\chi(g g') = \chi(bkb' k')=\chi(bkb')=\chi (b \tilde{b}\tilde{k})=\chi(b \tilde{b})=\chi(b)\chi(\tilde{b}),
\eeq
where $\tilde{b}\tilde{k}$ is the Iwasawa decomposition of $kb'$. From this we see also
\begin{align}
\chi(bg)=\chi(b)\chi(g),\quad\quad b\in B(\mathbb{A}), g\in G(\mathbb{A}).
\end{align}

The global character splits  into an Euler product over local factors:
\beq
\chi(g)=\prod_{p\leq \infty} \chi_p(g_p), \qquad \qquad g_p\in G(\mathbb{Q}_p).
\eeq

There is a one-to-one correspondence between such characters and weights of the Lie algebra $\mathfrak{g}(\mathbb{R})$, or, more precisely, complex linear functionals $\lambda\in \mathfrak{h}_\mathbb{C}^{\star}=\mathfrak{h}(\mathbb{R})^{\star}\otimes_\mathbb{R} \mathbb{C}$, where $\mathfrak{h}(\mathbb{R})$ is the Cartan subalgebra of $\mathfrak{g}(\mathbb{R})$. We  define a  logarithm map $H$  as follows:
\beq
\label{eq:Hlog}
H : G(\mathbb{A}) \to \mathfrak{h}(\mathbb{R}), 
\eeq
defined by 
\beq
\label{eq:Borel-logarithm-map}
H(g)=H(nak)=H(a)=\log |a|.
\eeq
The absolute value is defined as follows. Parametrise the group element $a\in A(\mathbb{A})$ by  
\beq
a=\exp\left(\sum_{\alpha\in \Pi} u_\alpha  H_\alpha \right), \qquad H_\alpha \in \mathfrak{h}(\mathbb{R}), \quad u_\alpha \in \mathbb{A},
\eeq
where $\Pi$ denotes the set of simple roots of $\mathfrak{g}(\mathbb{R})$.

 Then we define
\beq
\label{eq:A-absolute-value}
\log |a| := \log \exp\left( \sum_{\alpha\in \Pi} |u_\alpha |  H_\alpha \right) = \sum_{\alpha\in \Pi} |u_\alpha|  H_\alpha = \sum_{\alpha\in \Pi}\left(\prod_{p\leq \infty} |u_{\alpha, p}|_p\right) H_\alpha ,
\eeq
where each $u_{\alpha, p}\in \mathbb{Q}_p$.

The choice of character $\chi$ can now be parametrised by the choice of linear functional $\lambda$ according to the formula:
\beq
\chi(g)=e^{\left<\lambda + \rho|H(g)\right>}=|a^{\lambda+\rho}|.
\label{characterBorel}
\eeq
Here, we introduced a convenient short-hand notation. The translation by the Weyl vector $\rho$ constitutes a convenient choice of normalisation. 

\begin{remark}[Modulus character]
The map 
\beq
b\mapsto e^{\left< 2 \rho |H(b)\right>}\equiv \delta_B(b), \qquad b\in B(\mathbb{A}),
\eeq
is often called the \emphindex[modulus character]{modular function} (or `\emphindex{modulus character}') of $B$. It is defined by 
\beq
\delta_B(b)= \Big|\det \text{ad}(b) \big|_\mathfrak{n} \Big|.
\eeq
In words, it is the modulus of the determinant of the adjoint representation of $b\in B(\mathbb{A})$, restricted to the Lie algebra $\mathfrak{n}$ of the unipotent radical $N$. By virtue of the properties (\ref{eq:Borel-logarithm-map}) of the logarithm map we have 
\beq
\delta_B(b)=\delta_B(na)=\delta_B(a).
\eeq
The modulus character corresponds to the Jacobian that relates the left- and right-invariant Haar measures on $B$. This implies in particular that under conjugation by $a\in A(\mathbb{A})$, i.e. 
\beq
n\mapsto ana^{-1}, 
\eeq
the Haar measure $dn$ on $N(\mathbb{A})$ transforms by 
\beq
dn \mapsto \delta_B(a) dn.
\eeq
This fact will play a crucial role in our calculations in chapter~\ref{ch:SL2-fourier} and onwards. 

See example \ref{ex_HaarSL2} for an explicit description for $SL(2,\mathbb{A})$. Using the modulus character we can write $\chi$ in the alternative form
\beq
\label{eq:sepmod}
\chi(g)=e^{\left<\lambda |H(g)\right>} \delta^{1/2}_B(g).
\eeq
This form of the character is common in the mathematical literature. 
\end{remark}

\begin{example}[Haar measure and modulus character for the Borel subgroup of $SL(2,\mathbb{A})$]
\label{ex_HaarSL2} For $SL(2,\mathbb{A})$ we can take 
\beq
b=na= \left(\begin{array}{cc}
1 & u \\
& 1\\
\end{array} \right) 
\left(\begin{array}{cc}
v &  \\
& v^{-1}\\
\end{array} \right),
\eeq
in which case the right-invariant Haar measure is 
\beq
dn da = \frac{du dv}{v},
\eeq
and the modulus character is  given by 
\beq
\delta_B(na) = \delta_B\left(\left(\begin{array}{cc}
1 & u \\
& 1\\
\end{array} \right) 
\left(\begin{array}{cc}
v &  \\
& v^{-1}\\
\end{array} \right)\right)= |v|^{2}.
\eeq

\end{example}

\subsection{Eisenstein series}
\index{Eisenstein series|textbf}

With the definition of the character $\chi$ on the Borel subgroup, we are now in a position to state Langlands' definition of an Eisenstein series for an arbitrary reductive group $G(\mathbb{A})$. 

\begin{definition}[Eisenstein series] Let $G(\mathbb{A})$ be an adelic group with discrete subgroup $G(\rats)$ and Borel subgroup $B(\rats)$. Let $\chi$ be a character on $B(\ads)$ that is trivial on $B(\rats)$. The (minimal) Eisenstein series is defined as the sum over images of the coset $B(\mathbb{Q})\backslash G(\mathbb{Q})$ by
\beq
E(\chi, g)=\sum_{\gamma\in B(\mathbb{Q})\backslash G(\mathbb{Q})} \chi(\gamma g)\,,
\eeq
and using the explicit parametrisation~\eqref{characterBorel}, the definition reads
\beq
E(\lambda, g)=\sum_{\gamma\in B(\mathbb{Q})\backslash G(\mathbb{Q})} e^{\left<\lambda + \rho|H(\gamma g) \right>}.
\label{generalEisenstein}
\eeq
\end{definition}
The series defined here is not the only possible type of Eisenstein series that one can define, although it is the one that we will be most interested in. However, in section~\ref{ESInducedRep} we will treat Eisenstein series in the context of automorphic representations. This will then provide us with a way of deriving different types of Eisenstein series, including the one just defined.\\
 
For the above series Godement~\cite{Godement,BorelBoulder} proved that the sum converges absolutely whenever $\lambda$ lies in the open subset
\beq
\{ \lambda \in \mathfrak{h}_\mathbb{C}^{\star} \, |\, \Re(\lambda)\in \rho + (\mathfrak{h}^{\star})^{+}\},
\label{absoluteconvergence}
\eeq
where the positive chamber $(\mathfrak{h}^{\star})^{+}$ is defined by 
\beq
(\mathfrak{h}^{\star})^{+}=\{ \Lambda \in \mathfrak{h}^{\star}\, |\, \left< \Lambda, H_\alpha\right> > 0, \, \forall \, \alpha\in \Pi\},
\eeq
so that we require $\langle\lambda, H_\alpha\rangle >1$ for all simple roots $\alpha$.
For discussing the spectral decomposition of $\mathcal{A}(G(\mathbb{Q})\backslash G(\mathbb{A}))$, see also section~\ref{classautorep}, one is mainly interested in the case when $\lambda \in i\mathfrak{h}^{\star}$. This choice for $\lambda$ is motivated in section~\ref{indautorep} by the fact that the inducing representation is unitary in this case. However, such values lie outside the domain of absolute convergence (\ref{absoluteconvergence}) of $E(\lambda, g)$, which seems worrisome for the spectral decomposition. This puzzle is resolved by the remarkable result of Langlands that the Eisenstein series $E(\lambda, g)$ can in fact be analytically continued outside of the domain (\ref{absoluteconvergence}) to a meromorphic function on all of $\mathfrak{h}_\mathbb{C}^{\star}$. To establish the analytic continuation a crucial property of $E(\lambda, g)$ is its \emph{functional relation} which relates its value at $\lambda$ to its value at the Weyl-transform of $\lambda$:
\begin{align}
\label{funrel5}
E(\lambda,g) = M(w,\lambda) E(w\lambda,g)\,, \qquad w\in \mathcal{W}(\mathfrak{g}),
\end{align}
where $M(w, \lambda)$ is a known function. The functional relation will be discussed in more detail in chapter~\ref{ch:CTF}. Another important property is that $E(\lambda,g)$ is an eigenfunction of the Laplace operator $\Delta_{G/K}$ \index{Eisenstein series!Laplace eigenfunction}:
\begin{align}
\label{eq:EisenLapl}
\Delta_{G/K} E(\lambda,g) = \frac12\left(\langle\lambda|\lambda\rangle - \langle \rho|\rho\rangle\right) E(\lambda,g).
\end{align}
This formula is derived in appendix~\ref{App:Laplace}. In fact, the Eisenstein series $E(\lambda,g)$ is a common eigenfunction of all $G(\mathbb{A})$-invariant differential operators which is a reflection of its $\mathcal{Z}(\mf{g})$-finiteness. The following is a useful property of Eisenstein series:
\begin{proposition}
In the special case when $\lambda = -\rho$ we have 
\beq
E(-\rho, g)=1.
\eeq
\end{proposition}
\begin{proof}
We note first that the value $\lambda=-\rho$ is outside the region of absolute convergence of the Eisenstein and the function must therefore be defined by analytic continuation using the functional equation~\eqref{funrel5}. One also observes that by~\eqref{eq:EisenLapl}, $E(-\rho, g)$ is an eigenfunction of the Laplacian $\Delta_{G/K}$  with eigenvalue zero; hence it must be a constant function. To fix the constant to unity, we note that by Langlands constant term formula (see chapter~\ref{ch:CTF}), the constant term of $E(-\rho, g)$ with respect to the maximal unipotent radical $N$ is 
\beq
\int_{N(\mathbb{Q})\backslash N(\mathbb{A})}E(-\rho, ng)dn=1,
\eeq
from which the claim follows.
\end{proof}

\chapter{Automorphic representations and Eisenstein series}
\label{ch:autoreps}
In this chapter we introduce the concept of an automorphic representation associated with an adelic group $G(\mathbb{A})$. Since this is a rather difficult concept to grasp at first sight, we shall begin with a heuristic discussion before we delve into the technical definition. Our main focus will then lie with the so-called principal series which is the representation relevant for  general Eisenstein series. We  provide some remarks on the problem of classifying all automorphic representations of $\mathcal{A}(G(\mathbb{Q})\backslash G(\mathbb{A}))$. The theory of automorphic representations is illustrated through a detailed discussion of Eisenstein series from the point of view of representations induced from parabolic subgroups $P(\mathbb{A})\subset G(\mathbb{A})$. We  conclude the chapter with an extensive example that explains the representation-theoretic relation between holomorphic and non-holomorphic Eisenstein series via the embedding of the holomorphic discrete series inside the principal series of the non-compact group $SL(2,\mathbb{R})$. 					

\section[A first glimpse]{A first glimpse}
We have already seen hints in section~\ref{sec_frommodtoaut} that automorphic forms on $SL(2,\mathbb{Z})\backslash SL(2,\mathbb{R})$ are intimately related to the representation theory of $SL(2,\mathbb{R})$. Here we will further develop this point of view and also generalise it to the adelic framework. 

The main idea is that the space of smooth functions $\varphi: SL(2,\mathbb{Z})\backslash SL(2,\mathbb{R})\to \mathbb{C}$ carries several actions:
\begin{itemize}
\item First, we have the action $\pi$ of $SL(2,\mathbb{R})$ by \emph{right-translation}:
\beq
\label{eq:rra5}
[\pi(h)\varphi](g):=\varphi(gh), \qquad g,h\in SL(2,\mathbb{R}).
\eeq
\item Second, we have the action of the universal enveloping algebra $\mathcal{U}(\mathfrak{sl}(2,\mathbb{C}))$ by differential operators:
\beq
(D_X \cdot \varphi)(g):=\frac{d}{dt}\varphi(g\cdot e^{tX})\Big|_{t=0}, \qquad X\in \mathcal{U}(\mathfrak{sl}(2,\mathbb{C})).
\label{diffop}
\eeq
\end{itemize}
Whenever one has a group action on a space it is natural to look for a decomposition into irreducible representations of the group. Moreover, since the centre $\mathcal{Z}$ of the universal enveloping algebra $\mathcal{U}(\mathfrak{g}_{\mathbb{C}})$ commutes with $SL(2,\mathbb{R})$ it is also natural to distinguish the irreducible components in terms of their eigenvalues with respect to differential operators in $\mathcal{Z}$. For these reasons the theory of automorphic forms on $SL(2,\mathbb{R})$ is closely related to the decomposition of the space $\mathcal{A}(SL(2,\mathbb{Z})\backslash SL(2, \mathbb{R}))$ into irreducible representations with respect to the right-regular action of $SL(2,\mathbb{R})$, compatible with the action by $\mathcal{U}$. To get an idea of what this entails, let us now look at an extremely simplified, though still enlightening, example. 

\begin{example}[Fourier analysis on $\mathbb{Z}\backslash \mathbb{R}$]
\label{example:FourierZR}
In this example we will look at the abelian situation where the space $SL(2,\mathbb{Z})\backslash SL(2,\mathbb{R})$ is replaced by $\mathbb{Z}\backslash \mathbb{R}$, the circle group. This is the setting of classical Fourier analysis. The space of smooth functions $C^{\infty}(\mathbb{Z}\backslash \mathbb{R})$ is then just  Fourier series where the coefficients are constrained to decay rapidly with increasing Fourier number. Let us now formalise this a little and try to analyse it in the spirit of automorphic forms. Consider the unitary character $\chi_k : \mathbb{Z}\backslash \mathbb{R}\to U(1)$, defined by $\chi_k(x)=e^{2\pi i kx}$ for $x\in \mathbb{R}, k\in \mathbb{Z}$. Any function $f\in C^{\infty}(\mathbb{Z}\backslash \mathbb{R})$ can then be expanded in a Fourier series in terms of these characters
\beq
f(x)=\sum_{k\in \mathbb{Z}}c_k\chi_k(x), 
\eeq
Recall that an automorphic form is also required to satisfy a moderate growth condition. In the present setting we can choose square-integrability of $f(x)$ as a suitable condition for moderate growth. Thus the space of automorphic forms on $\mathbb{Z}\backslash \mathbb{R}$ can be taken to be the Hilbert space $L^2(\mathbb{Z}\backslash \mathbb{R})\subset C^{\infty}(\mathbb{Z}\backslash \mathbb{R})$ where the Fourier coefficients satisfy the square-integrability condition
\beq
\sum_{k\in \mathbb{Z}} |c_k|^2<\infty. 
\label{squareintegrabillity}
\eeq
Each character $\chi_k$ generates a one-dimensional irreducible subspace $V_k=\mathbb{C}\chi_k\subset L^2(\mathbb{Z}\backslash \mathbb{R})$ and the \emph{regular representation} $\pi$ of $\mathbb{R}$ defined by 
\beq
(\pi(y)f)(x):=f(x+y), \qquad x, y\in \mathbb{R}, 
\label{regularaction}
\eeq
is diagonalized by the subspaces $V_k$:
\beq
\pi(y)\cdot v=\chi_k(y) v, \qquad v\in V_k, \, y\in \mathbb{R}.
\label{diagonalactionVk}
\eeq
The set of equivalence classes of unitary representations of a group $G$ is called the \emph{unitary dual}, usually denoted $\widehat{G}$. In our example, the unitary dual $\widehat{\mathbb{R}}$ is simply the space of Fourier coefficients subject to the condition (\ref{squareintegrabillity}):
\beq
\widehat{\mathbb{R}}=L^2(\mathbb{Z})=\{(c_k)\, | \, \sum_k |c_k|^2<\infty\}.
\eeq
This gives the spectral decomposition  of the Hilbert space $L^2(\mathbb{Z}\backslash \mathbb{R})$. The fact that the spectrum is discrete is a general feature of spectral theory on compact spaces like $S^1\cong \mathbb{Z}\backslash \mathbb{R}$. 
 
\end{example}

Before we proceed with the adelic perspective we shall consider one more simple example that illustrates another feature that has a counterpart on the general theory of automorphic forms. 

\begin{example}[Fourier analysis on $\mathbb{R}$]
\label{example:FourierR}
Consider the same setting as in the previous example, namely $G=\mathbb{R}$, but we now take the discrete subgroup $G(\mathbb{Z})$ to be trivial. In other words, we are interested in the 
regular representation of $\mathbb{R}$ on the Hilbert space $L^2(\mathbb{R})$. The regular action $\pi$ is defined in the same way as in (\ref{regularaction}), but now this action is diagonalised on a continuous family of characters $\chi_\zeta : \mathbb{R}\to U(1)$ defined by $\chi_\zeta(x)=e^{2\pi i \zeta x}$ with $\zeta, x\in \mathbb{R}$. On the irreducible subspaces $V_\zeta = \mathbb{C}\chi_\zeta$ we then have 
\beq 
\pi(x)\cdot v = \chi_\zeta(x) v, \qquad \zeta, x\in \mathbb{R}, \, v\in V_\zeta.
\eeq
The unitary dual in this case is a ``continuous direct sum'', or \emph{direct integral} of irreducible representations, meaning that any function $f\in L^2(\mathbb{R})$ can be written as a continuous version of a Fourier series:
\beq
f(x) =\int_{\mathbb{R}}\hat{f}(\zeta)\chi_\zeta(x) d\zeta, 
\eeq
where $\hat{f}(\zeta)$ is the standard Fourier transform of $f$ with respect to the character $\chi_\zeta$. 

From the above analysis we conclude that the spectral decomposition of $L^2(\mathbb{R})$ with respect to the regular action of $\mathbb{R}$ has only a continuous part, in stark contrast with the situation in example~\ref{example:FourierZR} above. The appearance of a continuous spectrum is a general feature of spectral analysis on non-compact spaces, just like the discrete spectrum always appears for compact spaces. We also note a curious feature, namely that although the characters $\chi_k$ are used to decompose the spectrum $L^2(\mathbb{R})$ they are in fact \emph{not} square-integrable. This is not a problem since the Fourier transform always preserves square integrability. An elaborate version of this phenomenon will reappear later in this chapter. 
\end{example} 

In the previous examples, we have illustrated how the spectral analysis on compact and non-compact spaces reveals very different properties. When we generalise this to the non-abelian setting of $SL(2,\mathbb{Z})\backslash SL(2,\mathbb{R})$ we actually combine these properties in the following sense. Consider for a moment the space of square-integrable automorphic forms $L^2(SL(2,\mathbb{Z})\backslash SL(2,\mathbb{R}))$ which is a subspace of all automorphic forms where the moderate growth condition is replaced by the square-integrability condition. In contrast to the abelian case $\mathbb{Z}\backslash \mathbb{R}$, the space $SL(2,\mathbb{Z})\backslash SL(2,\mathbb{R})$ is certainly non-compact and we therefore expect that a spectral analysis would give rise to a continuous spectrum. In addition, the quotient $SL(2,\mathbb{Z})\backslash SL(2,\mathbb{R})$ has \emph{finite volume} and therefore also gives rise to a discrete spectrum; $SL(2,\ints)$ is what is known as a \emphindex[subgroup!co-compact]{co-compact subgroup}. Indeed, it was proven by Selberg that 
\beq
L^2(SL(2,\mathbb{Z})\backslash SL(2,\mathbb{R}))=\mathbb{C}\oplus L^2_{\mathrm{cusp}}(SL(2,\mathbb{Z})\backslash SL(2,\mathbb{R}))\oplus L^2_{\mathrm{cont}}(SL(2,\mathbb{Z})\backslash SL(2,\mathbb{R})),
\label{spectralSL2}
\eeq
where: 
\begin{itemize}
\item the first factor $\mathbb{C}$ represents the constant functions (these are considered to be part of the discrete spectrum and arise also as the residue of the non-holomorphic Eisenstein series $E(s,z)$ at $s=1$, see also section~\ref{SL2limits}), 
\item the remainder of the discrete spectrum is $L^2_{\mathrm{cusp}}(SL(2,\mathbb{Z})\backslash SL(2,\mathbb{R}))$ which is spanned by Maass cusp forms $\varphi$ with a discrete set of eigenvalues $\lambda_n , n=1, 2, 3, \dots, $ with respect to the Laplacian $\Delta_\mathbb{H}$, 
\item the continuous spectrum $L^2_{\mathrm{cont}}(SL(2,\mathbb{Z})\backslash SL(2,\mathbb{R}))$ is a direct integral of non-holomorphic Eisenstein series $E(s, z)$ (\ref{nonholeis}). 
\end{itemize}

\begin{remark} Note that the non-holomorphic Eisenstein series $E(s,z)$ are not-square integrable and so are not themselves part of $L^2(SL(2,\mathbb{Z})\backslash SL(2,\mathbb{R}))$ they nevertheless play a key role in parametrising the unitary dual, in a very similar vein as the non-square integrable, continuous characters $\chi_\zeta$ occurred in the spectral decomposition of $L^2(\mathbb{R})$ in example~\ref{example:FourierR}. The spectral decomposition \eqref{spectralSL2} forms a crucial ingredient in the \emphindex{Selberg trace formula} (see~\cite{Arthur:2005} for a nice introduction). 
\end{remark}

Although Eisenstein series are not square integrable they are still important for the representation theoretic aspects of automorphic forms and therefore it is natural to enhance the 
space of automorphic forms to include non-square integrable objects. One then replaces the square-integrable condition with a more general moderate growth condition, leading to the full space of automorphic forms $\mathcal{A}(SL(2,\mathbb{Z})\backslash SL(2,\mathbb{R}))$. We thus have the inclusions of function spaces
\beq
L^2_{\mathrm{cusp}}(SL(2,\mathbb{Z})\backslash SL(2,\mathbb{R}))\subset L^2(SL(2,\mathbb{Z})\backslash SL(2,\mathbb{R}))\subset \mathcal{A}(SL(2,\mathbb{Z})\backslash SL(2,\mathbb{R})).
\eeq

The representation-theoretic aspects of automorphic forms is however not yet complete, as the above treatment is missing an important ingredient. The space $\mathcal{A}(SL(2,\mathbb{Z})\backslash SL(2,\mathbb{R}))$ carries an additional action by \emph{Hecke operators}, which has not yet been taken into account. We will treat Hecke operators in detail in chapter~\ref{ch:Hecke} so here we shall merely offer some qualitative remarks. A Hecke operator is an operation $T_p : \mathcal{A}\to \mathcal{A}$, parametrised by a prime number $p<\infty$. The set of all Hecke operators $\{T_p\}_{p<\infty}$ forms a commutative ring, called the \emph{Hecke algebra}. An element $\varphi \in \mathcal{A}(SL(2,\mathbb{Z})\backslash SL(2,\mathbb{R}))$ which is an eigenvector for all Hecke operators
\beq
T_p \varphi = \lambda_p\varphi,
\eeq
is called a \emph{Hecke eigenform}. Here the eigenvalues $\lambda_p$ carry the arithmetic information contained in $\varphi$. However, the right-regular action of $SL(2,\mathbb{R})$ on $\mathcal{A}$ cannot be used to accommodate the action of the Hecke algebra and so the analysis of automorphic forms in terms of the representation theory of $SL(2,\mathbb{R})$ is incomplete. One of the reasons for passing to the adelic picture is precisely to remedy this problem. The basic idea is this: if we consider the right-regular action of $SL(2,\mathbb{A})$ on the space $\mathcal{A}(SL(2,\mathbb{Q})\backslash SL(2,\mathbb{A}))$, then the Hecke eigenvalues $\{\lambda_p\}$ parametrise the irreducible representations of the right regular action of the local subgroups $SL(2,\mathbb{Q}_p)$ on $\mathcal{A}(SL(2,\mathbb{Q})\backslash SL(2,\mathbb{A}))$. Thus, from the adelic perspective the Hecke algebra plays the same role at the non-archimedean places of $SL(2,\mathbb{A})$, as the universal enveloping algebra $\mathcal{U}(\mathfrak{sl}(2,\mathbb{C}))$ does at the archimedean place. The Hecke action is thus implicitly already taken into account in the general definition~\ref{defauto} of an adelic automorphic form. 

\newpage

\begin{example}[The action of Hecke operators on non-holomorphic Eisenstein series]
For illustration we consider here a simple example of how the Hecke operators act on the non-holomorphic Eisenstein series $E(s,z)$. For any integer $n>0$ we define the operator $T_n$ as follows:
\beq 
(T_n E)(s,z):=\frac{1}{n}\sum_{d|n}\sum_{b=0}^{d-1}E\left(s, \frac{nz+bd}{d^2}\right).
\eeq
In chapter~\ref{ch:Hecke} we will show that 
\beq
(T_n E)(s,z)=\lambda_n E(s, z), 
\eeq
with 
\beq
\lambda_n =n^{s-1/2}\sigma_{1-2s}(n).
\eeq
This is precisely the numerical Fourier coefficient in the Fourier expansion of $E(s,z)$; see equation~(\ref{SL2FC2}). The operators $T_n$ further satisfy the following basic relation
\beq
T_m T_n= \sum_{d|(m,n)}\frac{1}{d}T_{mn/d^2},
\eeq
characterising the Hecke algebra. This algebra is generated by the subset of Hecke operators $T_p$ for $p$ a (finite) prime, hence the fundamental information is contained in the prime eigenvalues $\lambda_p$, as claimed in the main text. In chapter~\ref{ch:Hecke} we provide much more details on Hecke operators and explain their link with the representation theory of $SL(2, \mathbb{Q}_p)$.
\end{example}

We now want to make sense of combined action on $\mathcal{A}(SL(2,\mathbb{Q})\backslash SL(2,\mathbb{A}))$ of $SL(2,\mathbb{A})$ by right-translation as well as the action of the universal enveloping algebra $\mathcal{U}(\mathfrak{sl}(2,\mathbb{C}))$ by differential operators. To this end it is useful to distinguish between the action of the finite part $SL(2,\mathbb{A}_f)=\prod_{p<\infty}SL(2,\mathbb{Q}_p)$ and the archimedean part $SL(2,\mathbb{R})$. For any $\varphi\in \mathcal{A}(SL(2,\mathbb{Q})\backslash SL(2,\mathbb{A}))$ we then have 
\beqa
(\pi(h_f)\varphi)(g)&=&\varphi(g h_f), \qquad g\in SL(2,\mathbb{A}), \, h_f\in SL(2,\mathbb{A}_f).
\nn \\
(\pi(h_\infty)\varphi)(g)&=&\varphi(g h_\infty), \qquad g\in SL(2,\mathbb{A}), \, h_\infty\in SL(2,\mathbb{R}).
\eqa
Here it is understood that the elements $h_f$ and $h_\infty$ are embedded in the canonical way into the adelic group. For instance, when we write $gh_\infty$ we really mean 
\beq
g\cdot \left(h_\infty, \begin{pmatrix} 1 & \\ & 1  \end{pmatrix}, \begin{pmatrix} 1 & \\ & 1  \end{pmatrix}, \dots, \begin{pmatrix} 1 & \\ & 1  \end{pmatrix}\right), \qquad g\in SL(2,\mathbb{A}), \,  h_\infty \in SL(2,\mathbb{R}).
\eeq
These two actions of course commute with each other. The action by $SL(2,\mathbb{A}_f)$ at the non-archimedean places also commutes with the $\mathcal{U}(\mathfrak{sl}(2,\mathbb{C}))$-action at the archimedean place, and so this gives a well-defined representation. On the other hand, the right-regular action of $K_\infty=SO(2, \mathbb{R})\subset SL(2,\mathbb{R})$ does \emph{not} commute with $\mathcal{U}(\mathfrak{sl}(2,\mathbb{C}))$. Rather, for $X\in \mathcal{U}(\mathfrak{sl}(2,\mathbb{C}))$ and $k_\infty\in SO(2,\mathbb{R})$ one has 
\beq
D_X \cdot \pi(k_\infty)=\pi(k_\infty) \cdot D_{k_\infty^{-1}Xk_\infty},
\label{gKmoduleproperty}
\eeq
where $D_X$ is the differential operator~(\ref{diffop}). One can check this by a  direct calculation:
\beqa 
D_X\cdot \left(\pi(k_\infty)\varphi\right)(g)&=& D_X\cdot \varphi(gk_\infty)
\nn \\ 
&=&  D_X \cdot \varphi\left(g\cdot \left(k_\infty, \begin{pmatrix} 1 & \\ & 1  \end{pmatrix}, \begin{pmatrix} 1 & \\ & 1  \end{pmatrix}, \dots, \begin{pmatrix} 1 & \\ & 1  \end{pmatrix}\right)\right) 
\eqa
where $k_\infty \in SO(2,\mathbb{R})$. Now using the definition of $D_X$ we find that the right hand side can be written as
\beq
\frac{d}{dt}\varphi \left(g\cdot \left( e^{Xt} k_\infty, \begin{pmatrix} 1 & \\ & 1  \end{pmatrix}, \begin{pmatrix} 1 & \\ & 1  \end{pmatrix}, \dots, \begin{pmatrix} 1 & \\ & 1  \end{pmatrix}\right)\right)\bigg|_{t=0}.
\label{gKverification}
\eeq
Inserting the identity $k_\infty k_\infty^{-1}$ and using the following property of the matrix exponential
\beq
k_\infty^{-1} e^{Xt} k_\infty=e^{k_\infty^{-1}X k_\infty t},
\eeq
we can rewrite equation~(\ref{gKverification}) as 
\beq
\frac{d}{dt}\varphi \left(gk_\infty \cdot \left( k_\infty^{-1} e^{Xt} k_\infty, \begin{pmatrix} 1 & \\ & 1  \end{pmatrix}, \dots, \begin{pmatrix} 1 & \\ & 1  \end{pmatrix}\right)\right)\bigg|_{t=0}=\pi(k_\infty) \cdot \left( D_{k_\infty^{-1}Xk_\infty}\cdot \varphi(g)\right),
\eeq
which is the right hand side of~(\ref{gKmoduleproperty}). This turns out to be the characteristic property of a so-called \emph{$(\mathfrak{g}, K_\infty)$-module}, a notion which will be properly defined in the next section.

\begin{remark} To ensure that the space $\mathcal{A}(SL(2,\mathbb{Q})\backslash SL(2,\mathbb{A}))$ is preserved under all three actions defined above, one must of course verify that they are compatible with definition~\ref{defauto}. In other words one should check that the three functions 
\beq
D_X\cdot \varphi(g), \qquad (\pi(h_\infty)\varphi)(g), \qquad  (\pi(h_f)\varphi)(g),
\eeq
all satisfy conditions $(1)$--$(4)$ in definition~\ref{defauto}. See, e.g., section 5.1 of \cite{MR2807433} for a detailed check of this.
\end{remark}

\begin{remark}
When speaking about an \emph{automorphic representation} of $SL(2,\mathbb{A})$ one really refers to a structure that carries a standard group representation with respect to the finite part  $SL(2,\mathbb{A}_f)$, and a $(\mathfrak{g}, K_\infty)$-module structure at the archimedean place. In the following section we will give the precise definition for an arbitrary reductive group $G(\mathbb{A})$ and discuss some central features of automorphic representations. 
\end{remark}

\section{Automorphic representations}
\label{sec:AR}

In this section we shall give the precise definition of an automorphic representation of an adelic group $G(\mathbb{A})$ and present some of the key features that will be important in subsequent chapters. Just as in the $SL(2,\mathbb{A})$-discussion of the previous section, we are interested in the combined actions of $G(\mathbb{A}_f)$ and $K_\infty$ by right-translation and the action of $\mathcal{U}(\mathfrak{g}_\mathbb{C})$ by differential operators. The general analysis goes through in a similar vein as above and the conclusion is that the space $\mathcal{A}(G(\mathbb{Q})\backslash G(\mathbb{A}))$ does not carry a group representation with respect the whole group $G(\mathbb{A})$, but only with respect to the finite part $G(\mathbb{A}_f)$. At the real place one has instead that $\mathcal{A}(G(\mathbb{Q})\backslash G(\mathbb{A}))$ carries the structure of a $(\mathfrak{g}, K_\infty)$-module, whose definition we will now recall. 

\begin{definition}[$(\mathfrak{g}, K)$-module] 
A $(\mathfrak{g}, K_\infty)${\it -module} is a complex vector space $V$ which carries an action of both the Lie algebra $\mathfrak{g}$ and $K_\infty $, such that all vectors $v\in V$ are $K$-finite, i.e. $\text{dim}\left< k_\infty\cdot v \, |\, k_\infty \in K_\infty\right><\infty$. The actions of $\mathfrak{g}$ and $K_\infty$ are furthermore required to be compatible in the following sense 
\beq
\label{compgK}
X\cdot k_\infty \cdot v = k_\infty\cdot \text{Ad}_{k_\infty^{-1}}(X)  \cdot v, \qquad k_\infty\in K_\infty, \, X\in \mathfrak{g}.
\eeq
\end{definition}
 
\begin{remark} In our context the complex vector space $V$ is $\mathcal{A}(G(\mathbb{Q})\backslash G(\mathbb{A}))$, the action by $X\in \mathfrak{g}$ is by the differential operator $D_X$ and the action by $k_\infty\in K_\infty$ is by right-translation. In this setting $k_\infty \cdot \text{Ad}_{k_\infty^{-1}}(X)$ means $\pi(k_\infty)\cdot D_{k_\infty^{-1} X k_\infty}$ and hence equation (\ref{gKmoduleproperty}) is precisely the compatibility condition (\ref{compgK}) for a $(\mathfrak{g}, K_\infty)$-module.
\end{remark}
 
\begin{remark} Let us  offer some remarks on the usefulness of $(\mf{g}, K_\infty)$-modules. The notion of $(\mf{g}, K_\infty)$-module was introduced by Harish-Chandra in his efforts on ``algebraisation'' of representations. Function spaces on groups are themselves typically not specific enough and there can be many function spaces that share the same algebraic features. For example, one can consider continuous functions on a group manifold and they are a perfectly nice representation of $G$. However, unless the functions are differentiable, this representations does not give rise to a representation of the Lie algebra $\mathfrak{g}$ that would be represented by an algebra of differential operators. There are many different types of differentiable functions on $G$ and the notion of $(\mf{g}, K_\infty)$-module mainly serves to eliminate the ambiguities related to choosing a type of differentiable function. 
\end{remark}

With the above concepts introduced, we are now ready to state the definition of an automorphic representation.

\begin{definition}[Automorphic representation]
\label{autorepdef}
A representation $\pi$ of $G(\mathbb{A})$ is called an \emphindex{automorphic representation} if it occurs as an irreducible constituent  in the 
decomposition of $\mathcal{A}(G(\mathbb{Q})\backslash G(\mathbb{A}))$ with respect to the  simultaneous action by 
\beq
(\mathfrak{g}_\infty, K_\infty)\times G(\mathbb{A}_f),
\eeq
where $K_\infty$ and $G(\mathbb{A}_f)$ act by right-translation and $\mathfrak{g}_\infty$ by differential operators at the archimedean place.
\end{definition}

We shall for short denote by $V$ the complex vector space on which an automorphic representation $\pi$ acts. Then $V$ is simultaneously a $(\mathfrak{g}_\infty, K_\infty)$-module and a $G(\mathbb{A}_f)$-module and the automorphic representation is also often denoted by the pair $(\pi, V)$.\\

Let $K_f\subset G(\mathbb{A}_f)$ be a compact open subgroup (not necessarily maximal) and $\sigma $  an irreducible representation of $K_\infty \times K_f$. Denote by $V[\sigma]$ the space of vectors in $V$ that transform according to $\sigma$ under the action of $K_\infty \times K_f$. We then have the following important definition:

\begin{definition}[Admissible representation]
A $(\mathfrak{g}_\infty, K_\infty)\times G(\mathbb{A}_f)$-module $V$ is called \emphindex{admissible} if the subspace $V[\sigma] \subset V$ is finite-dimensional for all $\sigma$.
\end{definition}

One  then has the following central result due to Flath.

\begin{theorem}[Tensor decomposition theorem]
\label{thm:Flath}\index{tensor decomposition theorem}
Let $(\pi, V)$ be an admissible automorphic representation. Then there exists an Euler product decomposition into local factors~
\beq
(\pi, V)=\bigotimes_{p\leq \infty} (\pi_p, V_p),
\label{autorepseuler}
\eeq
where the archimedean component $(\pi_\infty, V_\infty)$ is a $(\mathfrak{g}_\infty, K_\infty)$-module according to the discussion above, while the non-archimedean components $(\pi_p, V_p)$ furnish representations of $G(\mathbb{Q}_p)$.
\end{theorem}

\begin{proof} See \cite{Flath}.
\end{proof}

Let us finally also introduce the notion of an unramified (or spherical) representation and associated spherical vector.

\begin{definition}[Unramified representation] 
{\label{def_unramifiedreps}}An admissible automorphic representation $\pi_p$ is called \emphindex{unramified} (or \emph{spherical}) if $V_p$ contains a non-zero vector ${\sf f}_p$ which is invariant under $K_p=G(\mathbb{Z}_p)$. We then call such an ${\sf f}_p$ a \emphindex{spherical vector}. Globally, one has the important notion that $(\pi, V)$ is a spherical automorphic representation if $\pi_p$ is spherical for almost all $p$. 
\end{definition}

\section{Principal series representation}
\label{indautorep}

Fix a Borel subgroup $B$ and a quasi-character $\chi: B\to \mathbb{C}^{\times}$ defined as in \eqref{characterBorel}:
\beq
\chi=e^{\left<\lambda+\rho|H\right>}.
\eeq
Consider now the following space of smooth functions on $G(\mathbb{A})$:
\beq
I(\chi)=\big\{ f : G(\mathbb{A})\to \mathbb{C} \, \big|\, f(bg)=\chi(b) f(g), \, b\in B(\mathbb{A})\big\}\,.
\label{principal}
\eeq
This is the function space of an induced  representation of $G(\mathbb{A})$ called the \emphindex[principal series!representation]{principal series representation}; it is also often denoted by $\text{Ind}_{B(\mathbb{A})}^{G(\mathbb{A})}\chi$. The principal series $I(\chi)$ provides an important example of an automorphic representation thanks to the theory of Eisenstein series which will be discussed in section~\ref{ESInducedRep} below. In general $I(\chi)$ is a reducible representation. However, one can show that when $\chi=\otimes_p\chi_p $ is an unramified character, $I(\chi)$ is an irreducible, admissible representation and so affords a decomposition into local factors:
\beq
I(\lambda)=\bigotimes_{p\leq \infty} I_p(\lambda)=\bigotimes_{p\leq \infty} \text{Ind}_{B(\mathbb{Q}_p)}^{G(\mathbb{Q}_p)}\chi_p.
\label{decprincipal}
\eeq

\begin{remark} The space $I(\lambda)$ can be viewed as the total space of a fiber bundle $I(\lambda)\to \mathfrak{h}_\mathbb{C}^{\star}$, with the fiber over each point $\lambda\in \mathfrak{h}_\mathbb{C}^{\star}$ consisting of the space of functions on $G(\mathbb{A})$ which transform by the character $e^{\left<\lambda+\rho|H\right>}$ under the left action of $B(\mathbb{A})$. 
\end{remark}

\begin{definition}[Standard section]
An element ${\sf f}_\lambda \in I(\lambda)$ is called a \emphindex{standard section} if it is $K_\mathbb{A}$-finite and its restriction to $K_\mathbb{A}$ is independent of $\lambda$.
\end{definition}

By virtue of (\ref{decprincipal}), any standard section ${\sf f}_\lambda \in I(\lambda)$  splits into a product of local 
factors 
\beq
{\sf f}_\lambda =\bigotimes_{p\leq \infty} {\sf f}_{\lambda, p}.
\eeq

\begin{definition}[Gelfand--Kirillov dimension]
Although the principal series representations $I_p(\lambda)$ are infinite-dimensional, one 
can still attach to them a notion of ``size'', which is called the \emphindex[dimension!Gelfand--Kirillov]{functional, or Gelfand--Kirillov, dimension} \index{dimension!functional} and denoted by $\text{GKdim}$. This is defined as the smallest number of variables on which we can realize the 
functions in $I_p(\lambda)$. 
\end{definition}

\begin{example}[Gelfand--Kirillov dimension]
The Gelfand-Kirillov dimension of the Hilbert space $L^2(\mathbb{R}^n)$ is $n$. Similarly, by the Iwasawa decomposition $G(\mathbb{Q}_p)=B(\mathbb{Q}_p) G(\mathbb{Z}_p)$ the sections ${\sf f}_{\lambda, p}\in I_p(\lambda)$ are determined by their restriction to $B(\mathbb{Q}_p)\backslash G(\mathbb{Q}_p)=G(\mathbb{Z}_p)$ and hence the functional dimension is (for generic $\lambda$)
\begin{align}
\label{GK-dim}
\text{GKdim}(I_p(\lambda))=\text{dim}_{\mathbb{Q}_p}\, B(\mathbb{Q}_p)\backslash G(\mathbb{Q}_p).
\end{align}
\end{example}

\section{Eisenstein series and induced representations}\label{ESInducedRep}
Let us now discuss the definition of Eisenstein series from the perspective of induced representations. 
One can think of an Eisenstein series as providing a $G(\mathbb{A})$-equivariant map from the induced representation $I(\lambda)$ of~\eqref{principal} into the space of automorphic forms:
\beq
E \, :\, I(\lambda) \,  \to\,  \mathcal{A}(G(\mathbb{Q})\backslash G(\mathbb{A})).
\eeq
For any standard section ${\sf f}_\lambda \in I(\lambda)$ the construction of the corresponding Eisenstein series is given by
\beq
E({\sf f}_\lambda, g)=\sum_{\gamma\in B(\mathbb{Q})\backslash G(\mathbb{Q})} {\sf f}_\lambda(\gamma g).
\label{sectionEisenstein}
\eeq
As ${\sf f}_\lambda$ varies in the fiber of $I(\lambda) \to \mathfrak{h}_\mathbb{C}^{\star}$ we thus obtain a family of Eisenstein series $E({\sf f}_\lambda, g)$ that satisfy all the conditions of definition~\ref{defauto} for an automorphic form; in particular, $K_\mathbb{A}$-finiteness follows from the fact that ${\sf f}_\lambda$ is a standard section. By virtue of the decomposition ${\sf f}_\lambda=\otimes_{p\leq \infty} {\sf f}_{\lambda, p}$, the Eisenstein series $E({\sf f}_\lambda, g)$ can be defined by choosing all the local factors ${\sf f}_{\lambda, p}$ separately. This gives a lot of freedom in defining the  Eisenstein series and is one of the main reasons why the adelic formalism is so powerful (see section~\ref{standardsectionSL2} for a detailed demonstration in the case of $SL(2,\mathbb{A})$). 

\begin{remark}
In order to recover the particular Eisenstein series of definition~\eqref{generalEisenstein}, one chooses the standard section ${\sf f}_\lambda$ to be equal to the inducing character, ${\sf f}_\lambda=e^{\langle\lambda + \rho|H\rangle}=\chi$.
\end{remark}

\begin{remark}
In section~\ref{standardsectionSL2}  we will illustrate for $G(\mathbb{A})=SL(2,\mathbb{A})$ how the more general construction in~\eqref{sectionEisenstein} interpolates between holomorphic and non-holomorphic Eisenstein series on the upper-half plane $\mathbb{H}$.
\end{remark}

\section{Classifying automorphic representations}
\label{classautorep}

It is one of the central unsolved problems in the theory of automorphic forms to classify all automorphic representations. In fact, all admissible automorphic representations have been classified (see, e.g., \cite{Bump}). This includes in particular the spherical, or unramified, representations.

The task of decomposing $\mathcal{A}(G(\mathbb{Q})\backslash G(\mathbb{A}))$ into irreducible representations is closely connected to the problem of decomposing the Hilbert space $L^2(G(\mathbb{Q})\backslash G(\mathbb{A}))$ under the unitary action of $G$. A priori this might seem a little surprising since an automorphic form need not be square-integrable; indeed the Eisenstein series $E(s, g)$ considered in section \ref{sec_frommodtoaut} provide an example of a non-square integrable automorphic form. The decomposition of 
$L^2(G(\mathbb{Q})\backslash G(\mathbb{A}))$ splits into two orthogonal spaces 
\beq
L^2(G(\mathbb{Q})\backslash G(\mathbb{A})) = L_{\text{discrete}}^2(G(\mathbb{Q})\backslash G(\mathbb{A}))\oplus L_{\text{continuous}}^2(G(\mathbb{Q})\backslash G(\mathbb{A})),
\label{spectraldecomposition}
\eeq
corresponding respectively to the discrete and continuous parts of the spectrum. It turns out that the discrete spectrum is spanned by cusp forms and residues of Eisenstein series \cite{LanglandsFE,Moeglin}.

\vspace{.5cm}

 It is a fundamental result in the spectral theory of automorphic forms that the space $\mathcal{A}_0(G(\mathbb{Q})\backslash G(\mathbb{A}))$ is  the subspace of $L_{\text{discrete}}^2(G(\mathbb{Q})\backslash G(\mathbb{A}))$ corresponding to smooth, cuspidal, $K$-finite and $\mathcal{Z}(\mathfrak{g})$-finite vectors occurring in the decomposition of the unitary representation $L^2(G(\mathbb{Q})\backslash G(\mathbb{A}))$. This is the reason that cusp forms constitute an essential part in the theory of automorphic forms.

While the discrete spectrum can be understood in this way as a direct sum of invariant subspaces spanned by cusp forms, the space $L_{\text{continuous}}^2(G(\mathbb{Q})\backslash G(\mathbb{A}))$ rather decomposes into a \emph{direct integral} over principal series representations of $G(\mathbb{R})$. Such integrals turn out to be parametrised by Eisenstein series, even though these by themselves are not square integrable (see the following section and also~\cite{Gelbart} for more on the continuous spectrum and the relation with Eisenstein series). This situation is a generalisation of the problem of decomposing the Hilbert space $L^2( \mathbb{R})$ via Fourier analysis in terms of the non-square integrable characters (Fourier modes) $e^{2\pi i xy}$, $x,y\in \mathbb{R}$, as discussed in example~\ref{example:FourierR}. The construction of Eisenstein series on $G(\mathbb{Q})\backslash G(\mathbb{A})$, generalising the function $E(g, s)$ of section \ref{sec_frommodtoaut}, therefore constitutes an equally important part of the theory of automorphic forms as that of analyzing the space of cusp forms. Moreover, the complement of $\mathcal{A}_0(G(\mathbb{Q})\backslash G(\mathbb{A}))$ inside the discrete spectrum $L_{\text{discrete}}^2(G(\mathbb{Q})\backslash G(\mathbb{A}))$ is conjecturally spanned by \emph{residues} of Eisenstein series $E(\lambda, g)$ for special values of the weight $\lambda$. Thus, one expects that the discrete spectrum decomposes according to:
\beq
L_{\text{discrete}}^2(G(\mathbb{Q})\backslash G(\mathbb{A}))=L_{\text{cusp}}^2(G(\mathbb{Q})\backslash G(\mathbb{A}))\oplus L_{\text{res}}^2(G(\mathbb{Q})\backslash G(\mathbb{A})).
\eeq

Arthur has outlined a set of conjectures that characterise precisely the weights $\lambda$ for which the representation becomes square-integrable \cite{Arthur89unipotentautomorphic} (for proofs of Arthur's conjectures in some cases, see \cite{MR748508,MR1426903,MR1300285,MR1026752,MR3034297}). See also section~\ref{sec-Square} for further discussions of square-integrability of Eisenstein series.

\begin{example}[A representation-theoretic viewpoint on Eisenstein series on $SL(2,\mathbb{A})$]
 \label{EisensteinSL2example}We now analyse the general Eisenstein series $E(\lambda,g)$ more explicitly for $G(\mathbb{A})=SL(2,\mathbb{A})$. In this case the space of (complex) weights $\mathfrak{h}^{\star}_\mathbb{C}$ is one-dimensional and spanned by the fundamental weight $\Lambda$ dual to the unique simple root $\alpha$ of the Lie algebra $\mathfrak{sl}(2,\mathbb{A})$. The Weyl vector $\rho$ is also identical to $\Lambda=\alpha/2$. Therefore, we can parametrise the weight appearing in (\ref{generalEisenstein}) with a single parameter $s\in\cx$ according to 
\begin{align}
\lambda= 2s \Lambda -\rho=(2s-1)\Lambda\quad\Longrightarrow \quad \lambda+\rho = 2s\Lambda\,.
\end{align}
The character $\chi : B(\mathbb{Q})\backslash B(\mathbb{A})\to \mathbb{C}^{\times}$ in (\ref{characterBorel}) can now be written  as
\begin{align}
\chi_s(g)\equiv e^{\left<\lambda+\rho|H(g)\right>}=e^{\left<2s\Lambda|H(a)\right>}\, .
\end{align}
We can write all these objects explicitly in the fundamental representation of $SL(2,\mathbb{A})$:
\begin{align}
\label{Iwasawa3}
g = nak = \left(\begin{array}{cc}
1 & u \\
& 1\\
\end{array} \right) 
\left(\begin{array}{cc}
v &  \\
& v^{-1}\\
\end{array} \right)
\, k
\end{align}
with $k\in K_{\mathbb{A}} $. Evaluated on the group element (\ref{Iwasawa3}) we then find
\begin{align}
\chi_s(g) = e^{2s\left<\Lambda|H(a)\right>} = |v|^{2s}
\end{align}
since $H(a) = \log |v|\cdot H_\alpha$ where $H_\alpha$ is the Cartan generator of $\mathfrak{sl}(2,\mathbb{A})$. The general Eisenstein series $E(\lambda,g)$ in (\ref{generalEisenstein}) now becomes
\begin{align}
\label{generalSL2AEisenstein}
E(s,g)=\sum_{\gamma\in B(\mathbb{Q})\backslash SL(2,\mathbb{Q})} \chi_s(\gamma g),
\end{align}
which is indeed equivalent to (\ref{adelic_EisensteinSL2}). This Eisenstein series is  attached to the induced representation 
\beq
I(s)=\text{Ind}_{B(\mathbb{A})}^{SL(2,\mathbb{A})} \chi_s.
\eeq
This representation is unitary when $s=\tfrac{1}{2}+it\in \tfrac{1}{2} + i\mathbb{R}_+$. In this simple example one can also give a more explicit description of the spectral problem of decomposing the space $L^2(SL(2,\mathbb{Q})\backslash SL(2,\mathbb{A}))$ with respect to the right-regular action of $SL(2,\mathbb{A})$. The decomposition (\ref{spectraldecomposition}) becomes in this case
\beq
L^2(SL(2,\mathbb{Q})\backslash SL(2,\mathbb{A}))=L_{\text{cusp}}^2(SL(2,\mathbb{Q})\backslash SL(2,\mathbb{A}))\oplus \mathbb{C}\oplus \int_0^{\infty} I\left(\tfrac{1}{2}+it\right) dt,
\eeq
where the discrete spectrum $L_{\text{discrete}}^2$ is represented by the space of cusp forms $L_{\text{cusp}}^2$ together with the space $\mathbb{C}$ of constant functions that correspond to the residual spectrum from the pole of $E(s,g)$ at $s=1$.  The continuous spectrum $L_{\text{continuous}}^2$ corresponds to the integral over the principal series $I(s)$, restricted to the unitary domain $s\in \tfrac{1}{2} + i\mathbb{R}_+$. For a proof of this statement, see the book by Gelbart \cite{Gelbart}. 
\end{example}

\section[Embedding of the discrete series in the principal series]{Embedding of the discrete series\\ in the principal series}
\label{standardsectionSL2}

Our aim in this section is to illustrate in great detail the construction of the general  Eisenstein series $E({\sf f}_\lambda, g)$ in (\ref{sectionEisenstein}) for the special case of $SL(2,\mathbb{A})$. We  will in particular demonstrate that when restricted to a function on $SL(2, \mathbb{Z})\backslash SL(2,\mathbb{R})$ this yields a generalisation of the classical non-holomorphic Eisenstein series, which in fact interpolates between the non-holomorphic function $E(s, z)$, $z\in \mathbb{H}$,  and the weight $w$ holomorphic Eisenstein series $E_w(z)$. We explain how to understand this representation-theoretically in terms of the embedding of the holomorphic discrete series of $SL(2,\mathbb{R})$ into the principal series.

\subsection{Eisenstein series for arbitrary standard sections}
\label{sec:nonspherical}

Let $I(\lambda)=\mathrm{Ind}_{B(\mathbb{A})}^{SL(2,\mathbb{A})}\chi_s$ be the induced representation \eqref{principal} for $SL(2,\mathbb{A})$. As in example \ref{EisensteinSL2example} we take the inducing character $\chi_s : B(\mathbb{Q})\backslash B(\mathbb{A}) \to \mathbb{C}^{\times}$ (extended to all of $SL(2,\mathbb{A})$) to be defined by 
\beq
\chi_s(b k)=\chi_s(b)=\chi_s\begin{pmatrix} v & \star \\  & v^{-1} \\ \end{pmatrix} = \left|v \right|^{2s}, \qquad s\in \mathbb{C},
\label{ex_defSL2char}
\eeq
where $b\in B(\mathbb{A})$ and $k\in K_{\mathbb{A}}=SO(2, \mathbb{R})\times \prod_{p<\infty} SL(2,\mathbb{Z}_p)$. 

Let ${\sf f}_\lambda = \otimes_p {\sf f}_{\lambda, p} \in {I}(\lambda)$ with each local factor
\beq
{\sf f}_{\lambda, p}\in {I}_p(\lambda)=\mathrm{Ind}_{B(\mathbb{Q}_p)}^{SL(2,\mathbb{Q}_p)}\chi_p, 
\eeq
determined by its restriction to $SL(2,\mathbb{Z}_p)=B(\mathbb{Q}_p)\backslash SL(2,\mathbb{Q}_p)$.  For the purposes of this example we shall now fix these local sections as follows.
\begin{itemize}
\item For the non-archimedean  places $p<\infty$ we choose the section ${\sf f}_{\lambda, p}$ to be the unique (normalised) spherical vector ${\sf f}^{\circ}_{\lambda, p}$ in ${I}_p(\lambda)$ defined by (see also section \ref{sec_sphericalvector})
\beq {\sf f}^{\circ}_{\lambda, p}(g_p)={\sf f}^{\circ}_{\lambda, p}(b_pk_p)=\chi_{s,p}(b_p), \qquad \qquad {\sf f}^{\circ}_{\lambda, p}(k_p)={\sf f}^{\circ}_{\lambda, p}(1)=1,
\label{ex_finite_spherical}
\eeq
where $b_p\in B(\mathbb{Q}_p)$ and $k_p\in SL(2,\mathbb{Z}_p)$. 
\item For the archimedean place $p=\infty $ we define ${\sf f}_{\lambda, \infty}\in I_\infty(\lambda)$ according to 
\beq
{\sf f}_{\lambda, \infty}(b_\infty)=\chi_{s,\infty}(b_\infty), \qquad {\sf f}_{\lambda, \infty}(g_\infty k_\infty)=e^{iw\theta} {\sf f}_{\lambda, \infty}(g_\infty), \qquad {\sf f}_{\lambda, \infty}(1)=1,
\label{eq_archsec}
\eeq
where $w\in \mathbb{Z}, \, g_\infty\in SL(2,\mathbb{R}),\,  b_\infty\in B(\mathbb{R})$ and 
\beq
k_\infty=\left(\begin{array}{cc}
\phantom{-}\cos\theta & \phantom{-}\sin \theta \\
-\sin \theta & \phantom{-}\cos \theta \\
\end{array} \right)\in SO(2, \mathbb{R}). 
\eeq
Notice that with this definition ${\sf f}_{\lambda, \infty}$ is a $K_\infty=SO(2,\mathbb{R})$-finite, but \emph{non-spherical} section of $I_\infty(\lambda)$. 
\end{itemize}
With these definitions of the local factors, the product 
\beq
{\sf f}_\lambda = {\sf f}_{\lambda, \infty} \otimes \bigotimes_{p<\infty} {\sf f}_{\lambda, p}^{\circ}
\eeq 
becomes a standard section (because $\lambda$ and $w$ are independent parameters) of the global representation $I(\lambda)$. 

With this choice of section we now construct the Eisenstein series 
\beq
E({\sf f}_\lambda, g)=\sum_{\gamma\in B(\mathbb{Q})\backslash G(\mathbb{Q})} \left({\sf f}_{\lambda, \infty}(\gamma g_\infty) \times \prod_{p< \infty} {\sf f}^{\circ}_{\lambda, p}(\gamma g_p)\right),
\label{nonsphericaleisensteinSL2}
\eeq
with $g=(g_\infty; g_2, g_3, \dots)\in SL(2,\mathbb{A})=SL(2,\mathbb{R})\times \prodp_{p<\infty} SL(2,\mathbb{Q}_p)$.
We now want to analyse the restriction of this adelic Eisenstein series to a function on $SL(2,\mathbb{R})$. To this end we fix the adelic group element to be the identity at 
all finite places:
\beq 
g=(g_\infty; 1, 1, \cdots )\in SL(2,\mathbb{A}). 
\eeq
In example~\ref{bijection} we showed the bijection of cosets $B(\mathbb{Q})\backslash SL(2,\mathbb{Q})\cong B(\mathbb{Z})\backslash SL(2,\mathbb{Z})$ with each $B(\rats)g$ coset having a representative in $SL(2, \ints)$ for $g \in SL(2, \rats)$. We will now use this write the Eisenstein series as a sum over $\gamma\in SL(2,\mathbb{Z})$. At the finite places $SL(2,\mathbb{Z})$ embeds into $SL(2,\mathbb{Z}_p)$, and hence, by (\ref{ex_finite_spherical}), we have 
\beq
{\sf f}^{\circ}_{\lambda, p}(\gamma)=1 \qquad (\text{for $\gamma\in SL(2,\mathbb{Z})$ and $p<\infty$}).
\eeq
The Eisenstein series $E({\sf f}_\lambda, g)$  therefore 
restricts to 
\beqa
E({\sf f}_{\lambda, \infty}, g_\infty)&=& \sum_{B(\mathbb{Q})\backslash SL(2,\mathbb{Q})}{\sf f}_{\lambda, \infty}(\gamma g_\infty)\cdot \prod_{p<\infty} {\sf f}^{\circ}_{\lambda, p}(\gamma)
\nn \\
&=&\sum_{B(\mathbb{Z})\backslash SL(2,\mathbb{Z})}{\sf f}_{\lambda, \infty}(\gamma g_\infty)\cdot \prod_{p<\infty} {\sf f}^{\circ}_{\lambda, p}(\gamma)
\nn \\
&=&\sum_{B(\mathbb{Z})\backslash SL(2,\mathbb{Z})}{\sf f}_{\lambda, \infty}(\gamma g_\infty).
\eqa
To relate this   to a function on the upper-half plane $\mathbb{H}$, we use the Iwasawa decomposition 
\beq
g_\infty = b_\infty k_\infty = n_\infty a_\infty k_\infty = \begin{pmatrix} 1 & x \\ & 1  \end{pmatrix} \begin{pmatrix} y^{1/2} & \\ & y^{-1/2} \end{pmatrix} \left(\begin{array}{cc}
\phantom{-}\cos\theta & \phantom{-}\sin \theta \\
-\sin \theta & \phantom{-}\cos \theta \\
\end{array} \right),
\label{IwasawaSL2R}
\eeq
which yields
\beq
E({\sf f}_{\lambda, \infty}, g_\infty)=e^{iw\theta} \sum_{B(\mathbb{Z})\backslash SL(2,\mathbb{Z})}{\sf f}_{\lambda, \infty}(\gamma b_\infty).
\label{ex_eisres}
\eeq
It remains to evaluate ${\sf f}_{\lambda, \infty}(\gamma b_\infty)$. To this end we perform an additional Iwasawa decomposition of $\gamma b_\infty$ with the result:
\beq
\gamma b_\infty = b'_\infty k'_\infty, \qquad \qquad \gamma= \begin{pmatrix} a & b \\ c & d  \end{pmatrix}\in SL(2,\mathbb{Z}), 
\eeq
with 
\beq
b'_\infty = \begin{pmatrix} \frac{y^{1/2}}{|cz+d|} & \star \\ 0 & \frac{|cz+d|}{y^{1/2}}  \end{pmatrix}, \qquad k'_\infty =\begin{pmatrix} \cos\theta' & \sin\theta' \\ -\sin\theta' & \cos\theta'  \end{pmatrix}=\frac{1}{|cz+d|} \begin{pmatrix} cx+d & -cy \\ cy & cx+d  \end{pmatrix} ,
\eeq
and $z=x+iy = b_\infty \cdot i$. Further using that $e^{i\theta}=\cos \theta + i \sin \theta$ we find 
\beq
e^{i\theta'}=\frac{|cz+d|}{cz+d},
\eeq
and hence, by (\ref{ex_defSL2char}) and (\ref{eq_archsec}), the section in (\ref{ex_eisres}) evaluates to
\beq
{\sf f}_{\lambda, \infty}(\gamma b_\infty)={\sf f}_{\lambda, \infty}(b'_\infty k'_\infty)= \left(\frac{|cz+d|}{cz+d} \right)^{w}\chi_{s,\infty}(b'_\infty)=\left( \frac{|cz+d|}{cz+d}\right)^w \frac{y^s}{|cz+d|^{2s}}.
\eeq
We thereby arrive at the following explicit expression for the Eisenstein series 
\beq
E({\sf f}_{(s,w), \infty}, g_\infty) = e^{iw\theta} \sum_{(c, d)=1} \frac{y^{s}}{(cz+d)^w |cz+d|^{2s-w}}.
\label{ex_interpolating}
\eeq
This is a non-holomorphic function on $SL(2,\mathbb{Z})\backslash SL(2,\mathbb{R})$ with weight $e^{iw\theta}$ under the right action of $k_\infty \in SO(2,\mathbb{R})$. 

\subsection{Representation theoretic interpretation}

Let us now analyse the Eisenstein series (\ref{ex_interpolating}) a little closer. First observe that $E({\sf f}_{(s,w), \infty}, g_\infty) $ interpolates between the classical holomorphic and non-holomorphic Eisenstein series on the upper-half plane. Indeed, restricting the value of $s$ to 
$s=w/2$ we obtain
\beq
E({\sf f}_{(w/2,w), \infty}, g_\infty)= e^{iw\theta} y^{w/2} \sum_{(c,d)=1} \frac{1}{(cz+d)^w}=e^{iw\theta} y^{w/2} E_w(z), 
\label{restriction}
\eeq
which we recognize as $\varphi_{f}(g_\infty)$ in the terminology of section \ref{sec_frommodtoaut} (see equation (\ref{varphi})) with $f=E_w(z)$ being the classical 
weight $w$ holomorphic Eisenstein series on $\mathbb{H}$. Similarly, restricting to $w=0$ in (\ref{ex_interpolating}) we obtain the classical non-holomorphic 
Eisenstein series
\beq
E({\sf f}_{(s,0), \infty}, g_\infty)=\sum_{(c,d)=1} \frac{y^{s}}{|cz+d|^{2s}}=E(s,z).
\eeq
Note that this is  compatible with the fact that fixing $w=0$ is equivalent to choosing the local section ${\sf f}_{\lambda, p}$ to be spherical also at the archimedean place, ${\sf f}_{\lambda, \infty}={\sf f}^{\circ}_{\lambda, \infty}$. The Eisenstein series $E({\sf f}_\lambda, g)$ in (\ref{sectionEisenstein}) then reduces to (\ref{adelic_EisensteinSL2}) which is indeed  the adelisation of $E(s, z)$.

While the non-holomorphic Eisenstein series $E(s,z)$  is naturally associated with the principal series $\text{Ind}_{B(\mathbb{R})}^{SL(2,\mathbb{R})}\chi_s$, the holomorphic Eisenstein series $E_w(z)$  is  rather associated with the so-called \emphindex{holomorphic discrete series} $\mathcal{D}(w)$ of $SL(2,\mathbb{R})$ for $w\in 2\ints_+$. 

In order to exhibit this, we make use of some of the notions of section~\ref{sec:Maasswt} on Maass forms of non-zero weight. To this end, let $\mathcal{H}(w)$ be the Hilbert space of smooth, complex square-integrable functions $g(z)$ on the upper-half plane $\UHP\cong SL(2, \mathbb{R})/SO(2,\mathbb{R})$ that are Maass forms of weight $w$. \index{Maass form!arbitrary weight}Specifically, they transform under $SL(2,\ints)$ as
\beq
g(z)\longmapsto \left(\frac{cz+d}{|cz+d|}\right)^{w}g\left(\frac{az+b}{cz+d}\right).
\label{defMaassw}
\eeq
If $f(z)$ is a \emph{holomorphic} modular form of weight $w$, then $g(z) = y^{w/2} f(z)\in \mathcal{H}(w)$. This is in particular true for $g(z)=y^{w/2} E_w(z)$ with the holomorphic Eisenstein series appearing in~\eqref{restriction}. This transition from modular holomorphic forms to Maass forms is closely related to the lift~\eqref{varphi} to $SL(2,\reals)$.

For Maass forms of weight $w$ one can define the following differential operators, called \emphindex[Maass operators]{Maass lowering and raising operators}:
\begin{subequations}
    \label{eq:Maass-lowering-raising}
\begin{align}
L_w &:=-(z-\bar{z})\frac{\partial}{\partial \bar{z}}-\frac{w}{2},\\
R_w &:=(z-\bar{z})\frac{\partial}{\partial {z}}+\frac{w}{2}.
\end{align}
\end{subequations}
These have the property that they lower and raise the weights of a Maass form by two units:
\begin{align}
L_w :\mathcal{H}_w \to \mathcal{H}_{w-2},\quad \quad 
R_w :\mathcal{H}_w \to \mathcal{H}_{w+2}.
\end{align}
If $f(z)$ is a holomorphic modular form of weight $w$, then one finds that the associated Maass form $g(z) = y^{w/2} f(z)$ is annihilated by the lowering operator:
\begin{align}
L_w g(z) = L_w \left( y^{w/2} f(z) \right) =0.
\end{align}

The Maass operators can be understood representation-theoretically as follows. Let 
\begin{align}
\varphi_f(x,y,\theta) = e^{iw\theta} y^{w/2} f(x+iy) = e^{iw\theta} y^{w/2} f(z)
\end{align}
be the  lift of $f\in\mathcal{M}_w\equiv \mathcal{M}_w(SL(2,\ints))$ to $SL(2,\mathbb{R})$ as in section \ref{sec:fromholtoaut}. (We see that the Maass form $g(z) = y^{w/2}f(z)$ on the upper half-plane is in some sense half of the lift.) Recall from~\eqref{Fdiffop} that $\varphi_f$ satisfies the differential equation
\beq
F \varphi_f = -2i  e^{-2i\theta} \left(y\frac{\partial}{\partial \bar{z}} -\frac{1}{4} \frac{\partial}{\partial \theta}\right) \varphi_f=0,
\label{Fdiffop2}
\eeq
where $F$ is the differential operator-realisation of the negative Chevalley generator $F$ of $SL(2,\mathbb{R})$ in the compact basis \eqref{eq:sl2cpt} (see also section~\ref{sec:SL2} for more details). Using \eqref{restriction} we can 
rewrite this as
\beq
F \varphi_f = e^{iw\theta}  L_w\left( y^{w/2} f(z)\right)=0,
\eeq
revealing that the Maass operator $L_w$ is nothing but the Chevalley generator $F$ after evaluating the derivative on $\theta$. Similarly, the Maass operator $R_w$ is the positive Chevalley operator $E$. We conclude from this that the holomorphic Eisenstein series $f(z)=E_w(z)$ gives rise to the lowest weight vector in the holomorphic discrete series of $SL(2,\mathbb{R})$ via~\eqref{restriction}. The weight $w$ is measured by the Cartan generator $H$ (c.f. \eqref{eq:sl2cpt}):
\beq
H \varphi_f=w\varphi_f.
\eeq
One can also  check that 
\beq
H(E  \varphi_f)= [H, E]   \varphi_f +E   (H  \varphi_f) =2 E  \varphi_f+ w E  \varphi_f= (w+2) E \varphi_f,
\eeq
This means that $\varphi_f$ for $f(z)=E_w(z)$ is the lowest weight state in a representation of $\mathfrak{sl}(2,\mathbb{R})$, whose states are obtained by acting successively with the raising operator 
$E$:
\beq
\{\varphi_f, E \varphi_f, E^2 \varphi_f, \dots \}.
\label{states}
\eeq
Here each vector $E^{n}  \varphi_f$, $n\geq 0$, is an  automorphic form on $SL(2,\mathbb{A})$, and hence belongs to the space $\mathcal{A}(SL(2,\mathbb{Q})\backslash SL(2,\mathbb{A}))$. The span of the states (\ref{states}) is a subspace $V$ of $\mathcal{A}(SL(2,\mathbb{Q})\backslash SL(2,\mathbb{A}))$ that is clearly invariant under the $\mathfrak{sl}(2,\mathbb{R})$-action. It is furthermore preserved by $K=SO(2,\mathbb{R})$, since each vector $E^{n} \varphi_f\in \mathcal{A}(SL(2,\mathbb{Q})\backslash SL(2,\mathbb{A}))$ is $K$-finite by definition \ref{defauto}. Thus, the vector space $V$ spanned by (\ref{states}) is a $(\mathfrak{g}, K)$-module. This is the $(\mathfrak{g}, K)$-module underlying the holomorphic discrete series $\mathcal{D}(w)$. It is conventional to think of the holomorphic Eisenstein series $f=E_{w}$ as the explicit lowest weight state.

In general  the principal series $\text{Ind}_{B(\mathbb{R})}^{SL(2,\mathbb{R})}\chi_s$ is not a lowest (or highest) weight representation; indeed the general Eisenstein series \eqref{ex_interpolating} is not annihilated by either of $E$ or $F$. However, as we restrict to the integer points $s=w/2$ of the complex weight space where $\chi_s$ lives, we  land on an irreducible submodule of $\text{Ind}_{B(\mathbb{R})}^{SL(2,\mathbb{R})}\chi_s$ which can be identified with the holomorphic discrete series $\mathcal{D}(w)$. In other words, we have discovered the well-known fact that the holomorphic discrete series can be embedded into the principal series for special values of the inducing character: 
\beq
\mathcal{D}(w) \subset \text{Ind}_{B(\mathbb{R})}^{SL(2,\mathbb{R})}\chi_s\big|_{s=w/2}.
\eeq

We should in fact be a little more careful. In \eqref{restriction} it is understood that the weight is restricted to be a positive  integer $w>0$. We should therefore distinguish between positive and negative weights in the spherical vector \eqref{eq_archsec}. The case $w>0$ leads to \eqref{restriction} as we just discussed. The negative weight case $w<0$ leads to the same conclusion, except that the restriction \eqref{restriction} now corresponds to the \emph{anti-holomorphic} Eisenstein series $\overline{E}_w(\bar{z})$. This Eisenstein series is then naturally associated with the anti-holomorphic discrete series $\overline{\mathcal{D}}(w)$ of $SL(2,\mathbb{R})$ which is defined analogously to \eqref{defMaassw} for anti-holormorphic functions $\overline{f}(\bar{z})$. The anti-holomorphic Eisenstein series $\overline{E}_w(\bar{z})$ lifts to a function $\varphi_{\overline{f}}$ which is annihilated by $E$, rather than $F$ and can therefore be interpreted as a \emph{highest weight vector} of $\overline{\mathcal{D}}(w)$ with weight $-w$. The negative Chevalley generator $F$ then lowers the weight by 2. 

\begin{figure}[t]
\begin{center}
\includegraphics[width=5cm]{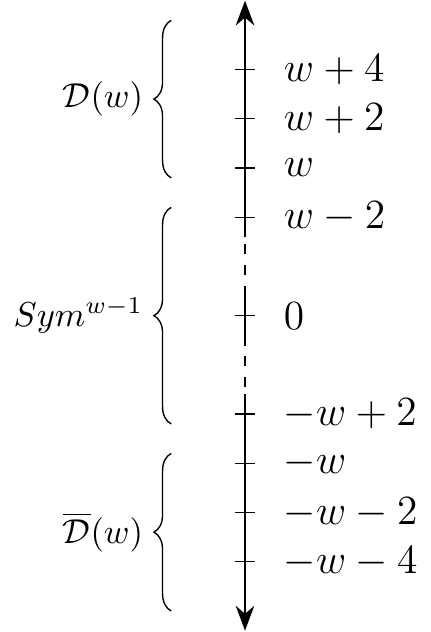}
\caption{Weight diagram for $SL(2,\mathbb{R})$.}
\label{fig:weightSL2}
\end{center}
\end{figure}

The above discussion shows that both the holomorphic and anti-holomorphic discrete series can be embedded into the principal series. The complement is a finite-dimensional representation of $SL(2,\mathbb{R})$, known as $Sym^{w-1}$; this is the $w$-dimensional symmetric power representation of $SL(2,\mathbb{R})$ acting on homogeneous degree $w-1$ polynomials in two real variables. Indeed, it is easy to see from the weight diagram in figure~\ref{fig:weightSL2} that the number of weights that are excluded from the holomorphic and anti-holomorphic discrete series are precisely equal to $w$. This implies that $Sym^{w-1}$ is the following quotient of the principal series by the discrete series for $w\in\ints_+$
\beq
\text{Ind}_{B(\mathbb{R})}^{SL(2,\mathbb{R})} \chi_{w/2} \Big/ \Big(\mathcal{D}(w)\oplus \overline{\mathcal{D}}(w)\Big)=Sym^{w-1}.
\eeq
The symmetric power representation is non-unitary.

\section{Eisenstein series for non-minimal parabolics*}
\label{nonminEis}
With a little effort one can generalise the construction of Eisenstein series to any parabolic subgroup. In what follows we restrict to standard parabolics, that is, those that contain the Borel subgroup $B(\adeles)=A(\adeles)N(\adeles)$ as discussed in section \ref{sec:parsubgp}. Fix such a standard parabolic subgroup $P(\adeles) \subset G(\adeles)$ with Langlands decomposition as in \eqref{LanglandsLG}
\begin{equation}
    P(\adeles) = L(\adeles) U(\adeles) = M(\adeles) A_P(\adeles) U(\adeles) 
\end{equation}

The full group $G(\adeles)$  factorises (non-uniquely) as 
\begin{equation}
    \label{decompGP2}
    G(\adeles)=M(\adeles)A_P(\adeles)U(\adeles)K_\adeles.
\end{equation}

For an arbitrary element of $G(\adeles)$ we write
\begin{equation}
    g = luk = ma_Puk, \qquad l \in L(\adeles), \, m\in M(\adeles), \, a_P\in A_P(\adeles), \, u\in U(\adeles), \, k\in K_\adeles.
\end{equation}
See example \ref{ex_parabolicsSLn} for some details of these decompositions in the cases $G(\adeles)=GL(n, \adeles)$ and $G(\ads) = SL(n, \ads)$.

\subsection{Multiplicative characters}

We now want to define multiplicative characters on $P(\mathbb{A})$ analogously to what we did for the Borel subgroup in section \ref{sec_multcharborel}. These will be homomorphisms 
\begin{equation}
    \chi_P\, :\, P(\rats)\backslash P(\adeles) \, \rightarrow \mathbb{C}^{\times}.
\end{equation}
determined by their restriction to the Levi subgroup
\begin{equation}
    \chi_P(lu)=\chi_P(l) \qquad l\in L(\adeles), \, u\in U(\adeles) \, .
\end{equation}

As in the minimal parabolic case the characters can be described by roots, but now the image in root space will be $\lie a_P^\ast$, the Lie algebra of $A_P$, instead of $\lie h^\ast$ and we will use a generalisation of the logarithm map $H$ from \eqref{eq:Borel-logarithm-map}.

Let $H_P : P(\adeles) \to \lie a_P(\reals)$ be defined by
\begin{equation}
    \label{eq:HP-def}
    H_P(p) = H_P(lu) = H_P(ma_Pu) = H_P(a_P) = \log |a_P| \qquad a_P \in A_P \subseteq A
\end{equation}
where the absolute value is defined as in \eqref{eq:A-absolute-value}.

A character can then be defined using a weight $\lambda_P \in \lie a_P^\ast$ as
\begin{equation}
    \chi_P(l) = e^{\langle \lambda_P + \rho_P | H_P(l) \rangle}
\end{equation}
where $\rho_P$ is now the restriction of the full Weyl vector to the positive roots $\Delta(\lie u) = \Delta_+ \setminus \langle \Sigma \rangle_+$ of $\lie g$ from section \ref{sec:parsubgp}
\begin{equation}
    \label{eq:rho_P}
    \rho_P=\frac{1}{2}\sum_{\alpha \in \Delta(\lie u)}\alpha.
\end{equation}

In example \ref{ex_parabolicsSLn} we give some details of the above construction for the case of $GL(n, \adeles)$.

\begin{example}[Parabolic subgroups and characters for $GL(n, \mathbb{A})$]
    \label{ex_parabolicsSLn} For $G(\adeles)=GL(n, \adeles)$ there is a bijection between standard parabolic subgroups $P(\adeles)$ and ordered partitions $(n_1, \dots, n_q)$ of $n$. It is then sometimes useful to start from this point of view when parametrising the subgroup $P$ instead of the one based on subsets $\Sigma \subset \Pi$ from section \ref{sec:parsubgp}.

For a given such partition we then have that $P(\adeles)=L(\adeles) U(\adeles) = M(\ads) A_P(\ads) U(\ads)$ can be expressed explicitly as
\begin{equation}
    L(\adeles)=\left\{ 
    \begin{pmatrix}
        l_1 &  & 0 \\ & \ddots & \\ 0 & & l_q 
    \end{pmatrix} \middle|\mathrel{} l_i\in GL(n_i, \adeles)\right\} ,
    \quad
    U(\adeles)=\left\{ 
    \begin{pmatrix}
        \id_{n_1} & \star & \star \\ & \ddots & \star\\ 0 & & \id_{n_q} 
    \end{pmatrix}
    \right\}
\end{equation}

\begin{equation}
    M(\adeles)=\left\{ 
    \begin{pmatrix}
        m_1 &  & 0 \\ & \ddots & \\ 0 & & m_q 
    \end{pmatrix} \middle|\mathrel{} m_i\in SL(n_i, \adeles)\right\} ,
    \quad
    A_P(\adeles)=\left\{ 
    \begin{pmatrix}
        a_1 \id_{n_1} &  & 0 \\ & \ddots & \\ 0 & & a_q \id_{n_q} 
    \end{pmatrix}
    \middle|\mathrel{} a_i\in GL(1, \adeles)\right\},
\end{equation}
where $\id_{n}$ denotes the $n\times n$ identity matrix. 

Similarly, instead of working with the Chevalley basis for $\lie a_P$ it is useful to choose a basis that reflects the block form in the parametrisation above.

We choose a basis $\tilde H_i$ for $\lie a_P$ such that $H_P$ from \eqref{eq:HP-def} becomes
\begin{equation}
    \label{eq:HP-GLn}
    H_P\Big(
    \begin{psmallmatrix}
        a_1 \id_{n_1} &  & 0 \\ & \ddots & \\ 0 & & a_q \id_{n_q} 
    \end{psmallmatrix} \Big) = 
    \sum_{i = 1}^q n_i \log\abs{a_i} \tilde H_i \, .
\end{equation}

For
\begin{equation}
    l = 
    \begin{psmallmatrix}
        l_1 &  & 0 \\ & \ddots & \\ 0 & & l_q 
    \end{psmallmatrix}
\end{equation}
we then obtain that
\begin{equation}
    H_P(l) = \sum_{i = 1}^q \log\abs{\det l_i} \tilde H_i \, .
\end{equation}

Now introduce a basis $\tilde \Lambda_i$ for $\lie a_P^*$ dual to $\tilde H_i$, that is, $\langle \tilde \Lambda_i | \tilde H_j \rangle = \delta_{ij}$, and let $\lambda_P$ and $\rho_P$ in $\lie a_P^*$ be parametrised by $\lambda_P = \sum_{i = 1}^q s_i \tilde \Lambda_i$ and $\rho_P = \sum_{i = 1}^q \rho_i \tilde \Lambda_i$ with $s_i, \rho_i \in \cmplx$. Note that since $\tilde H_i$ is not the Chevalley basis, the $\tilde \Lambda_i$ are not the standard fundamental weights.

Any character on $P(\rats) \bs P(\ads)$ can then be constructed by
\begin{equation}
    \label{eq:GLn-character}
    \chi_P(lu) = \chi_P(l) = e^{\langle \lambda_P + \rho_P | H_P(l) \rangle} = \prod_{i=1}^q \abs{\det l_i}^{s_i + \rho_i} \, .
\end{equation}

For $G = SL(n, \ads)$ we have the restriction $\prod_{i=1}^q a_i = 1$ which reduces the number of independent elements in the sum \eqref{eq:HP-GLn} spanning $\lie a_P$. In the same way, the parameters in $\lambda_P = \sum_{i=1}^q s_i \tilde \Lambda_i$ are also restricted. A general character on $P(\rats) \bs P(\ads)$ for $SL(n, \ads)$ can thus be seen as special cases of \eqref{eq:GLn-character}. Explicitly, we require that $\langle \lambda_P | \id_n \rangle = 0$ and, since $\id_n = \sum_i n_i \tilde H_i$, we get the restriction that $\sum_i n_i s_i = 0$.

\end{example}

\subsection{Parabolically induced representations}\label{parindrep}
Associated with the parabolic subgroup $P(\mathbb{A})$ we now consider the following space of functions
\beq
I_P(\lambda)=\{ f\, :\, G(\mathbb{A})\, \rightarrow \mathbb{C} \, |\, f(gp)=e^{\left<\lambda+\rho_P| H_P(p)\right>} f(g), \, g\in G(\mathbb{A}), \, p\in P(\mathbb{A})\}
\eeq
where $\lambda \in \lie a_P^*$. Note that we have suppressed the subscript $P$ on $\lambda$ for brevity. 

This is the function space of the induced representation $\text{Ind}_{P(\mathbb{A})}^{G(\mathbb{A})} \chi_P = \text{Ind}_{P(\mathbb{A})}^{G(\mathbb{A})} e^{\left< \lambda+\rho_P | H_P\right>}$. We will also refer to this as the \emphindex{principal series}, although strictly speaking that name should be reserved for the case $P(\mathbb{A})=B(\mathbb{A})$, the Borel subgroup. In that case we have $I_B(\lambda)=I(\lambda)$ from (\ref{principal}).  The generic functional dimension of the representation $\mathrm{Ind}_{P(\mathbb{A})}^{G(\mathbb{A})} \chi_P$ is, similarly to~\eqref{GK-dim}, given by
\begin{align}
\label{GK-dim2}
\mathrm{GKdim} (I_P(\lambda)) = \dim(G)-\dim(P).
\end{align}

It is now straightforward to construct an Eisenstein series associated with the induced representation $I_P(\lambda)$. It takes the form
\beq\label{ESeriesMaxParabDefn}
E(\lambda, P, g)=\sum_{\gamma\in P(\mathbb{Q})\backslash G(\mathbb{Q})} e^{\left<\lambda +\rho_P|H_P(\gamma g)\right>}.
\eeq
In case $P$ is the Borel subgroup we write $E(\lambda, B, g)=E(\lambda, g)$ and we recover the Eisenstein series in (\ref{sectionEisenstein}). 

Just as for the Borel subgroup, we can of course also start from any standard section ${\sf f}_\lambda\in I_P(\lambda)$ and obtain another Eisenstein series
\beq
E({\sf f}_\lambda, P, g)=\sum_{\gamma\in P(\mathbb{Q})\backslash G(\mathbb{Q})} {\sf f}_\lambda(\gamma g).
\eeq

    One can generalise the Eisenstein series $E(\lambda, P, g)$ even further by modifying the induced representation $I_P(\lambda)$ as follows. For $P(\ads) = L(\ads) U(\ads)$ let $\sigma$ be a representation of $L(\mathbb{A})$ and define
\beq
\sigma_\lambda(l)=\sigma(l) e^{\left<\lambda +\rho_P| H_P(l)\right>}, \qquad l\in L(\mathbb{A}), \, \, \lambda \in \mathfrak{a}_P^{\star}(\mathbb{C}).
\eeq
We then consider the associated induced representation $\text{Ind}_{P(\mathbb{A})}^{G(\mathbb{A})}\sigma_\lambda$. This corresponds to automorphic forms on $L(\mathbb{Q})\backslash L(\mathbb{A})$, extended to $P(\mathbb{A})$ by triviality on $U(\mathbb{A})$. More specifically, it is the space of functions
\beq
\phi \, :\,  \big(L(\mathbb{Q})U(\mathbb{A})  \big) \backslash G(\mathbb{A}) \rightarrow \mathbb{C}, 
\eeq
such that for each fixed $g\in G(\mathbb{A})$ the function 
\beq
\phi_g \, :\, l \rightarrow \phi(lg)
\eeq
is a vector in the finite-dimensional  space $V$ of automorphic forms on $L(\mathbb{Q})\backslash L(\mathbb{A})$ transforming according to the representation $\sigma$.

For each $ \phi\in \text{Ind}_{P(\mathbb{A})}^{G(\mathbb{A})}\sigma_\lambda, \, \, \lambda \in \mathfrak{a}_P^{\star}(\mathbb{C})$ we now have the Eisenstein series 
\beq
E(\lambda, P, \phi, g)=\sum_{\gamma\in P(\mathbb{Q})\backslash G(\mathbb{Q})} \phi(\gamma g) e^{\left< \lambda+\rho_P |H_P(\gamma g)\right>}.
\eeq 
In the mathematical literature, one often takes $\phi\in \mathcal{A}_0(L(\mathbb{Q})\backslash L(\mathbb{A}))$, the space of cusp forms on $L(\mathbb{A})$. In this case, Langlands has proven the analytic continuation and functional relation for $E(\lambda, P, \phi, g)$ \cite{LanglandsFE}.

\begin{proposition}
\label{prop-minimaltomaximal}
The Eisenstein series $E(\lambda, g)$, induced from the Borel subgroup $B$, is a special case of the Eisenstein series $E(\lambda, P_{i_*}, \phi, g)$, where $P_{i_*}$  is a maximal parabolic subgroup associated with the simple root $\alpha_{i_*}$ (see section~\ref{sec:parsubgp}). 
\end{proposition}
\begin{proof}
To see this, we follow the argument in \cite{Green:2010kv}. First note that $B_{L_{i_*}}=L_{i_*}\cap B$ is a Borel subgroup of the Levi $L_{i_*}\subset P_{i_*}$. This implies that any $\gamma\in B(\mathbb{Q})\backslash G(\mathbb{Q})$ can be uniquely decomposed as $\gamma=\gamma_1\gamma_2$, with 
\beq
\gamma_1\in B_{L_{i_*}}(\mathbb{Q})\backslash L_{i_*}(\mathbb{Q}), \qquad \gamma_2\in P_{i_*}(\mathbb{Q})\backslash G(\mathbb{Q}).
\eeq
We can thus rewrite the Eisenstein series $E(\lambda, g)$ as follows:
\beq
E(\lambda, g)=\sum_{\gamma\in B(\mathbb{Q})\backslash G(\mathbb{Q})}e^{\left<\lambda+\rho | H(\gamma g)\right>}=\sum_{\gamma_1\in B_{L_{i_*}}(\mathbb{Q})\backslash L_{i_*}(\mathbb{Q})}\sum_{\gamma_2\in P_{i_*}(\mathbb{Q})\backslash G(\mathbb{Q})}e^{\left<\lambda+\rho| H(\gamma_1\gamma_2g)\right>}. 
\eeq
Any complex weight $\lambda\in \mathfrak{a}^{\star}(\mathbb{C})$ can be decomposed by orthogonal projection on the plane with normal vector $\Lambda_{i_*}$, the fundamental weight associated with node $i_*$. Explicitly, 
\begin{align}
\lambda& =\lambda_{\parallel i_*}+\lambda_{\perp i_*},\quad\textrm{with}\quad
\lambda_{\parallel i_*} = \frac{\langle \lambda|\Lambda_{i_*}\rangle}{\langle \Lambda_{i_*}|\Lambda_{i_*}\rangle} \Lambda_{i_*},\quad
\lambda_{\perp i_*} = \lambda - \lambda_{\parallel i_*}.
\end{align}
By the defining properties of $\Lambda_{i_*}$, we find that $\lambda_{\perp i_*}$ is a complex linear combination of all simple roots different from $\alpha_{i_*}$ and $\lambda_{\parallel i_*}$ is proportional to $\Lambda_{i_*}$. Furthermore, for any $l\in L_{i_*}$ generated by only \emph{positive} roots  in the Lie algebra $\mathfrak{l}_{i_*}$ we have for any $g\in G$
\beq
\left<\lambda_{\parallel i_*}+\rho_{\parallel i_*}| H(lg)\right>=\left<\lambda_{\parallel i_*}+\rho_{\parallel i_*}| H(g)\right>.
\eeq
Hence, in the domain of absolute convergence, we may decompose the summation as
\beqa
E(\lambda, g)&= & \sum_{\gamma_1\in B_{L_{i_*}}(\mathbb{Q})\backslash L_{i_*}(\mathbb{Q})}\sum_{\gamma_2\in P_{i_*}(\mathbb{Q})\backslash G(\mathbb{Q})}e^{\left<\lambda+\rho| H(\gamma_1\gamma_2g)\right>}
\nn \\
&=&
\sum_{\gamma_2\in P_{i_*}(\mathbb{Q})\backslash G(\mathbb{Q})}\bigg[\sum_{\gamma_1\in B_{L_{i_*}}(\mathbb{Q})\backslash L_{i_*}(\mathbb{Q})}e^{\left<\lambda_{\perp i_*}+\rho_{\perp i_*}| H(\gamma_1\gamma_2g)\right>}\bigg]e^{\left<\lambda_{\parallel i_*}+\rho_{\parallel i_*}| H(\gamma_2g)\right>}
\nn \\ 
&=& \sum_{\gamma_2\in P_{i_*}(\mathbb{Q})\backslash G(\mathbb{Q})}e^{\left<\lambda_{\parallel i_*}+\rho_{\parallel i_*}| H(\gamma_2g)\right>}
\phi(\gamma_2g),
\eqa
where the  function 
\beq
\phi(g)=\sum_{\gamma_1\in B_{L_{i_*}}(\mathbb{Q})\backslash L_{i_*}(\mathbb{Q})}e^{\left<\lambda_{\perp i_*}+\rho_{\perp i_*}| H(\gamma_1 g)\right>}
\eeq
 is an Eisenstein series on the Levi $L_{i_*}$, induced from the Borel subgroup $B_{L_{i_*}}$.
\end{proof}
\begin{remark}
Proposition~\ref{prop-minimaltomaximal} can be straightforwardly generalised to give a relation between Eisenstein series $E(\lambda, g)$, induced from the Borel subgroup $B$, and Eisenstein series $E(\lambda, P, \phi, g)$ induced from an arbitrary parabolic subgroup $P$ (not necessarily maximal). 
\end{remark}

To illustrate the general analysis of this section, we shall conclude with two explicit examples dealing with the case of maximal parabolic subgroups $P(\mathbb{A})$. This is the opposite extreme compared to the Borel subgroup, which we recall is a minimal parabolic. 

\begin{example}[Eisenstein series on $SL(n, \mathbb{A})$ induced from a maximal parabolic $P$]
\label{ex_maxparabolicSLn} Again we consider $SL(n, \mathbb{A})$ and in this example we take $P(\mathbb{A})$ to be a \emphindex[subgroup!parabolic!maximal]{maximal parabolic}. Maximal parabolic subgroups are simply classified by partitions $n\mapsto (n_1, n_2)$. The Levi decomposition is therefore $P(\mathbb{A})=L(\mathbb{A})U(\mathbb{A})$ with Levi subgroup given by 
\beq
L(\mathbb{A})=\left\{ \left(\begin{array}{cc} l_1 & \\ & l_2  \end{array} \right) \, |\, l_i\in GL(n_i, \mathbb{A})\right\}.
\eeq
The character $\chi_P$  evaluates to 
\beq
\chi_P(luk)=\chi_P(l)=|\det l_1|^{s_1+\rho_1}|\det l_2|^{s_2+\rho_2},
\eeq
for $l \in L(\mathbb{A}), u \in U(\mathbb{A}), k \in K_\mathbb{A}$. Since the restriction to $SL(n, \ads)$ requires $n_1s_1+n_2s_2=0$ we effectively only have one independent parameter $s\in \mathbb{C}$.

Let now $\phi \in \text{Ind}_{P(\mathbb{A})}^{SL(n, \mathbb{A})} \sigma_s$ such that 
\beq
\phi(luk)=\phi(l). 
\eeq
The associated Eisenstein series is 
\beq
E(\lambda, P, s, g)=\sum_{\gamma\in  P(\mathbb{Q})\backslash G(\mathbb{Q})} \phi(\gamma g) \chi_P(\gamma g).
\eeq 
\end{example}

The final example discusses a construction using a 5-grading of $\mathfrak{g}$ that is possible for all simply-laced ADE groups but $SL(2,\ads)$. 

\begin{example}[Eisenstein series on $E_6$, $E_7$, $E_8$ induced from Heisenberg parabolic subgroups]
\label{ex_maximalparabolic}
Let now $G(\mathbb{A})$ be the adelisation of either $E_6$, $E_7$ or $E_8$ with Lie algebra $\mathfrak{g}$. We shall analyse the above construction for a very special type of maximal parabolic subgroup of $G$, known as the \emphindex{Heisenberg parabolic}, henceforth denoted by $P_{\text{Heis}}$. This parabolic subgroup is associated with the highest root $\theta$ of $\mathfrak{g}$. Similar arguments can be made for the ADE-series of simple Lie algebras \cite{kazhdan1990smallest,Kazhdan:2001nx} but then not necessarily resulting in $P_\text{Heis}$ being maximal. 

Associated with $\theta$  the Lie algebra exhibits a canonical 5-grading 
\beq
\mathfrak{g}=\mathfrak{g}_{-2} \oplus \mathfrak{g}_{-1}\oplus \mathfrak{g}_0 \oplus \mathfrak{g}_1 \oplus \mathfrak{g}_2,
\eeq
where the subscript indicates the eigenvalue under the Cartan generator $H_\theta$ associated with $\theta$, and $\mathfrak{g}_{\pm 2}$ are one-dimensional subspaces spanned by the corresponding Chevalley generators $E_{\pm \theta}$. The triple $(H_\theta, E_\theta, E_{-\theta})$ generates an $\mathfrak{sl}(2,\mathbb{R})$ subalgebra:
\beq
[H_\theta, E_\theta]=2E_\theta, \qquad [H_\theta, F_\theta]=-2F_\theta, \qquad [E_\theta, F_\theta]=H_\theta.
\label{ex_sl2theta}
\eeq
 The zeroth subspace $\mathfrak{g}_0$ is of the form $\mathfrak{m}_{\text{Heis}}\oplus \mathbb{C} H_\theta$, where $\mathfrak{m}_{\text{Heis}}\subset \mathfrak{g}$ is a reductive Lie algebra corresponding to  the commutant of the $\mathfrak{sl}(2,\mathbb{R})$-algebra (\ref{ex_sl2theta}) inside $\mathfrak{g}$. The nilpotent subspace $\mathfrak{g}_1\oplus \mathfrak{g}_2$ is a Heisenberg  algebra of dimension $\text{dim}_{\mathbb{R}}\, \mathfrak{g}_1+1\equiv 2d+1$, with commutator
\beq
[\mathfrak{g}_1, \mathfrak{g}_1]\subseteq \mathfrak{g}_2.
\label{HeisLiealgebra}
\eeq
We set 
\beq 
\mathfrak{p}=\mathfrak{g}_0 \oplus \mathfrak{g}_1 \oplus \mathfrak{g}_2,
\eeq
which is the Lie algebra of a maximal parabolic subgroup $P_{\text{Heis}}\subset G$,  the \emphindex{Heisenberg parabolic}. Its Levi and Langlands decompositions are 
\beq
P_{\text{Heis}}=L_{\text{Heis}}U_{\text{Heis}}=M_{\text{Heis}}A_{\text{Heis}} U_{\text{Heis}},
\eeq
where the Levi subgroup $L_{\text{Heis}}=M_{\text{Heis}} A_{\text{Heis}}$ is the exponentiation of $\mathfrak{g}_0$ further decomposing into $\lie a_\text{Heis} = \cmplx H_\theta$ and $\lie m_\text{Heis}$ above, and the unipotent radical $U_{\text{Heis}}$ is the Heisenberg group whose Lie algebra is $\mathfrak{g}_1\oplus \mathfrak{g}_2$.

For $P_\text{Heis}$ we define a logarithm map $H_P : P_\text{Heis}(\ads) \to \lie a_\text{Heis}(\reals) = \reals H_\theta$ according to \eqref{eq:HP-def}. Where for $a_P = \exp(v H_\theta)$ with $v \in \ads$
\begin{equation}
    H_P(p) = H_P(m a_P u) = H_P(a_P) = |v| H_\theta \, .
\end{equation}

Let $\Lambda_\theta$ be the weight dual to $H_\theta$, i.e. defined by 
\begin{equation}
\langle \Lambda_\theta \mid H_\theta \rangle = 1, 
\end{equation}
and parametrise an arbitrary linear functional $\lambda \in \lie a^*
_\text{Heis}(\cmplx)$ by $\lambda = 2s \Lambda_\theta - \rho_{P_\text{Heis}}$ with $s \in \cmplx$.
The Weyl vector is  $\rho_{P_{\text{Heis}}}=\Lambda_\theta$ so we have
\begin{equation}
    \lambda = (2s-1)\Lambda_\theta, \qquad \lambda+\rho_{P_\text{Heis}}=2s\Lambda_\theta. 
\end{equation}
 
Putting this together we obtain a character 
 \beq 
 \chi_{P_{\text{Heis}}}\equiv \chi_s : P_{\text{Heis}}(\mathbb{Q})\backslash P_{\text{Heis}}(\mathbb{A}) \to \mathbb{C}^{\times}
 \eeq
defined by 
\beq
\chi_s=e^{\left< 2s\Lambda_\theta | H_P\right>} .
\eeq

We extend it to all of $G(\mathbb{A})$ by demanding that it is trivial on $K_\mathbb{A}$ by virtue of the decomposition \eqref{decompGP2}. Explictly, we have 
\beq
\chi_s(g)=\chi_s(ma_Puk)=\chi_s(a_P) = e^{\left<2s\Lambda_\theta |H_P(a_P)\right>}=|v|^{2s}.
\eeq

The associated induced representation $\text{Ind}_{P_{\text{Heis}}(\mathbb{A})}^{G(\mathbb{A})} \chi_s$ is called the \emphindex{degenerate principal series}. At the infinite place it has functional dimension 
\beq
\GKdim \text{Ind}_{P_{\text{Heis}}(\mathbb{R})}^{G(\mathbb{R})} \chi_s=\text{dim}\, P_{\text{Heis}}(\mathbb{R})\backslash G(\mathbb{R})=\text{dim}\, \mathfrak{g}_1\oplus \mathfrak{g}_2= 2d+1, 
\eeq 
and depends on a single complex parameter $s$. In contrast, the generic principal series induced from the Borel subgroup $B$ depends on $r=\text{rank}\, \mathfrak{g}$ parameters 
$(s_1, \dots, s_r)\in \mathbb{C}^{r}$. Formally one can view $\text{Ind}_{P_{\text{Heis}}(\mathbb{A})}^{G(\mathbb{A})} \chi_s$ as the limit of $\text{Ind}_{B(\mathbb{A})}^{G(\mathbb{A})} e^{\left<\sum_{i=1}^{r} s_i \Lambda_i | H\right>}$ when projecting onto the complement of a complex co-dimension one locus in $\mathbb{C}^r$. 

For any standard section ${\sf f}_s\in \text{Ind}_{P_{\text{Heis}}(\mathbb{A})}^{G(\mathbb{A})} \chi_s$ the Eisenstein series attached to the degenerate principal series is 
\beq
E({\sf f}_s, P_{\text{Heis}}, g)= \sum_{\gamma \in P_{\text{Heis}}(\mathbb{Q})\backslash G(\mathbb{Q})} {\sf f}_s(\gamma g).
\eeq
 
 This Eisenstein series has interesting properties because its residues at the poles in the complex $s$-plane
give rise to automorphic forms attached to special types of (unipotent) representations of $G$ which have very small functional dimensions (typically of dimension less than $2d +1)$. 
The smallest such representation is known as the \emphindex{minimal representation} of $G$ and it has functional dimension $d+1$. Automorphic forms attached to minimal 
representations were analysed from this point of view in \cite{GRS}, and has also played an important role in physical applications \cite{Gunaydin:2001bt,Kazhdan:2001nx,Pioline:2005vi,Gunaydin:2005mx,Pioline:2009qt,Pioline:2010kb,Green:2011vz}. See also sections \ref{sec:outlook-strings} and \ref{sec:BH-counting}.
\end{example}

\section{Induced representations and spherical vectors*}
\label{sec:IndRep}

In this section, we elaborate on some aspects of induced representations and present a different construction of Eisenstein series based on global spherical vectors.

\subsection{Analytic induction}

We first recall briefly the notion of analytic induction of a representation from a parabolic subgroup $P$ of a group $G$ that was already discussed in section~\ref{indautorep}. Generalising the set-up there, let
\begin{align}
\sigma \,:\, P \to  GL(V)
\end{align}
be a (smooth) representation of the parabolic subgroup $P$ on a (finite-dimensional) vector space $V$. In the case of the principal series, the vector space $V$ is one-dimensional, $P$ is a Borel subgroup and the representation $\sigma$ is by a character $\chi$. 
\begin{definition}[Induced representation]
The \emphindex{induced representation} of $\sigma$ from $P$ to $G$ is defined to be the space of smooth functions
\begin{align}
\Ind_P^G\sigma = \left\{ f \,:\, G\to V \st f(pg) = \sigma(p) f(g)\right\} .
\end{align}
The group $G$ acts on this function space by right translation.
\end{definition}

Functions in the induced representation space are by construction fully determined by their values on the coset $P\backslash G$ (and the representation $\sigma$). The coset typically has a complicated cell decomposition but often one can learn many important aspects by focussing on the biggest cell. As was already mentioned in~\eqref{GK-dim2}, the functional dimension of the this induced representation is given by the dimension of the coset $P\backslash G$, see for instance~\cite{Borho,MR0404366}.

\begin{example}[Principal series of $SL(2,\reals)$]
Consider the case $G=SL(2,\reals)$ and $P=\bar{B}=\left\{\begin{pmatrix}A&0\\Z&A^{-1}\end{pmatrix}\st A\in \reals^\times, Z\in \reals\right\}$, the standard lower parabolic subgroup. In the \emphindex{line model} of $\bar{B}\backslash G$ one decomposes a matrix $g\in SL(2,\reals)$ as 
\begin{align}
g=\bar{b}n:\quad
\begin{pmatrix}
a& b\\c&d
\end{pmatrix}
= \begin{pmatrix}
A & \\
Z & A^{-1}
\end{pmatrix}
\begin{pmatrix}
1 & x\\
&1
\end{pmatrix}
=\begin{pmatrix}
A & Ax\\
Z & Zx+A^{-1}
\end{pmatrix}
\end{align}
which is always possible as long as $a\neq 0$. The representative of $g$ is then
\begin{align}
\label{eq:line}
n = \begin{pmatrix}
1 & x\\
&1
\end{pmatrix}
\quad\textrm{with}
\quad
x= \frac{b}{a} \in \reals.
\end{align}
If $a=0$, this \emphindex[LU-decomposition]{LU-type decomposition} breaks down. One can check that there is only a single class $\bar{B}\backslash G$ with $a=0$ and that this can be represented by
\begin{align}
\begin{pmatrix}
0&-1\\1&0
\end{pmatrix}.
\end{align}
This point should be thought of as the `point at infinity' that compactifies the real line $x\in\reals$ to a compact space isomorphic to a circle. (The Iwasawa decomposition implies that $\bar{B}\backslash G$ is a quotient of the compact space $K$ and this way of looking at the coset is called the \emphindex{circle model} or \emphindex{ball model}.) The set of matrices with $a\neq 0$ corresponds to the big cell whereas the `point at infinity' is the single smaller cell in this example.

Suppose that $\bar{B}$ is represented on the one-dimensional vector space $V=\cx$ by 
\begin{align}
\sigma\left( \begin{pmatrix} A&\\Z& A^{-1} \end{pmatrix}\right) = |A|^{-2s}.
\end{align}
This is a character on $\bar{B}$ and so we are dealing with a principal series representation as discussed in section~\ref{indautorep}. A given representative~\eqref{eq:line} in the big cell transforms as follows under right multiplication by $SL(2,\reals)$:
\begin{align}
n g = \bar{b}' n'
\end{align} 
with 
\begin{align}
n' &= \begin{pmatrix}
1& x'\\&1
\end{pmatrix}
\quad\textrm{with}\quad
x' = \frac{dx+b}{cx+a},\nn\\
\bar{b}' &=  \begin{pmatrix}1&\\ \star&1\end{pmatrix} 
\begin{pmatrix}
A&\\
&A
\end{pmatrix}
\end{align}
with $A=c x+a$. The $\star$ entry is easy to compute but of no further importance to us. If we think of $\Ind_{\bar{B}}^G \sigma$ as the space of functions $f(x)$ on the real line since $x$ is a good coordinate on $\bar{B}\backslash G$ almost everywhere, the functions in $\Ind_{\bar{B}}^G \sigma$ satisfy
\begin{align}
f( x g) = |cx+a|^{-2s} f\left(\frac{dx+b}{cx+a}\right)
\end{align}
under $g=\begin{pmatrix}a&b\\c&d\end{pmatrix}\in SL(2,\reals)$. 

The action above can be linearised to obtain a realisation of the Lie algebra $\mf{sl}(2,\reals)$ on functions on the real line corresponding to the principal series. The result is
\begin{align}
\label{eq:sl2diffs}
h = -2x\partial_x -2s,\quad
e = \partial_x,\quad
f = -x^2 \partial_x -2sx,
\end{align}
where the standard generators~\eqref{eq:sl2basis} have been used. The above action can be formally extended by including the point at infinity corresponding to the small cell.
\end{example}

One can also investigate induced representations for $p$-adic groups.
\begin{example}[Principal series of $SL(2,\rats_p)$]
We use the same type of decomposition of parabolic subgroup and line model as in the example above, so that the induced function space becomes
\begin{align}
\label{eq:SL2pInd}
\Ind_{\bar{B}(\rats_p)}^{SL(2,\rats_p)} \sigma \cong \left\{ f \,:\, \rats_p \to \cx \st f(xg) = |cx+a|_p^{-2s} f\left(\frac{x+b}{cx+a}\right)\right\}
\end{align}
almost everywhere.
\end{example}

The next example deals with the \emphindex{conformal realisation} of the orthogonal group. 

\begin{example}[Conformal realisation of $SO(d,d)$]
\label{ex:CFT}
Our discussion will be based on the Lie algebra $\mf{so}(d,d)$ that can be presented by generators $M^{IJ}=-M^{JI}$ for $I,J=1,\ldots, 2d$ that satisfy
\begin{align}
\lb M^{IJ}, M^{KL} \rb = 2\eta^{K[J} M^{I]L}-2\eta^{L[J} M^{I]K} = \eta^{JK} M^{IL} -\eta^{IK} M^{JL}-\eta^{JL} M^{IK} +\eta^{IL} M^{JK}.
\end{align}
The $SO(d,d)$ metric $\eta=\textrm{diag}(+-+-\ldots+-)$ is fully split.

We consider the Lie algebra $\mf{so}(d,d)$ (for $d\geq 3$) in the graded decomposition with respect to $\mf{so}(d-1,d-1)$:
\begin{align}
\label{eq:so3grad}
\mf{so}(d,d) = {\underbrace{(2(d-1))}_{K^i}}^{(-1)} \oplus \big( \underbrace{\mf{so}(d-1,d-1)}_{M^{ij}} \oplus \underbrace{\mf{gl}(1)}_{D}\big)^{(0)} \oplus {\underbrace{(2(d-1))}_{P^i}}^{(+1)}.
\end{align}
Indices can be raised and lowered freely with the $\mf{so}(d-1,d-1)$ metric $\eta$ that is obtained by truncating the $\mf{so}(d,d)$ metric. The superscripts denote the eigenvalues under the $\mf{gl}(1)$ generator $D$. The commutation relations in this decomposition are
\begin{align}
\label{eq:conf}
\lb M^{ij} , M^{kl} \rb & = \eta^{jk} M^{il}+\eta^{ik} M^{jl}-\eta^{jl} M^{ik}+\eta^{il} M^{jk},\nn\\
\lb M^{ij}, P_k \rb &= 2\delta^{[j}_{k} P^{i]},\nn\\
\lb M^{ij}, K_k \rb &= 2\delta^{[j}_{k} K^{i]},\nn\\
\lb M^{ij}, D \rb &=0,\nn\\
\lb D, P^i \rb &= P^i,\nn\\
\lb D, K^i \rb &= -K^i,\nn\\
\lb P^i, K^j\rb &= -2M^{ij} -2\eta^{ij} D.
\end{align}
This is the standard conformal algebra. The index range is $i=1,\ldots, 2d-2$.

We choose the parabolic subgroup $P$ to be generated by all generators with $D$-eigenvalue $\leq 0$, i.e., $P$ has generators $K^i$, $M^{ij}$ and $D$. For this example we take the representation $\sigma$ of $P$ to be again a character, namely such that $\sigma(e^{\log(a) D}) = |a|^{4s}$. The big cell in the coset space $\bar{P}\backslash G$ is then represented by the unipotent elements 
\begin{align}
n= \exp(x^i P_i) 
\end{align}
and has dimension $2d-2$, corresponding to the Gelfand--Kirillov dimension of the induced representation. The Lie algebra $\mf{so}(d,d)$ is then represented on functions of the $2d-2$ variables $x^i$ by the differential operators
\begin{align}
\label{eq:soddDS}
M_{ij} &= x_i \partial_j - x_j \partial_i,\nn\\
P_i &= \partial_i,\nn\\
D &= -x^i\partial_i +4s,\nn\\
K_i &=  2x_i x^k\partial_k - x^2 \partial_i -8sx_i.
\end{align}
We recognise the familiar form of generators of the conformal group on scalar fields defined in $(2d-2)$-dimensional space-time~\cite{DiFrancesco:1997nk}, albeit here in a Wick-rotated form. The generators $P_i$ act as standard translation operators, the generators $M^{ij}$ as geometrical rotation generators, the generator $D$ acts as a dilatation operator and $K_i$ as a special conformal transformation. The value $s$ that appears in the inducing character is related to the conformal weight of the field. It is also possible to extend this representation to non-scalar fields by including a non-trivial representation of the $M^{ij}$ generator in the induction process; this leads to bosonic and fermionic higher spin analogues of the scalar field representation, see for example~\cite{Fernando:2015tiu} for a discussion.

In the representation~\eqref{eq:conf}, the quadratic Casimir operator evaluates to
\begin{align}
-M^{IJ} M_{IJ} &= - M^{ij} M_{ij} - P^i K_i - K^i P_i +2D^2 \nn\\
&=8s(2s-2d+2).
\end{align}
In language that will be introduced in section~\ref{sec:orbits}, the representation~\eqref{eq:conf} is associated with the next-to-minimal nilpotent orbit of Bala--Carter type $2A_1$.

The degenerate series representation~\eqref{eq:conf} admits a subquotient for the value $s=\frac{d-2}2$ that corresponds to the minimal unitary representation of $SO(d,d)$ with Gelfand--Kirillov dimension $2d-3$. In conformal field theory, this representation is known as the \emphindex{singleton} and plays a distinguished role in the context of the \emphindex{AdS/CFT correspondence}, see for example~\cite{Bekaert:2011js} for a discussion.
\end{example}

\subsection{Spherical vector constructions}

The importance of spherical elements in a representation has already been emphasised repeatedly in this chapter. Using explicit realisations of induced representations of the type just discussed, one can try to find spherical elements very concretely. This knowledge can also be used in the construction of Eisenstein series and Fourier expansions as will become more clear in the subsequent chapters. We begin with a few examples.

\begin{example}[Spherical vector in $SL(2,\reals)$ and $SL(2,\rats_p)$ principal series]
\label{ex:SL2sph}
The compact subgroup $SO(2)\subset SL(2,\reals)$ is generated by the element
\begin{align}
e-f = (1+x^2) \partial_x + 2sx,
\end{align}
where we have used the explicit differential operators~\eqref{eq:sl2diffs} in the line model. We are therefore asked to look for a function $f^\circ_\infty(x)$ that is annihilated by this differential operator:
\begin{align}
(1+x^2) \partial_x f^\circ_\infty(x) +2sx f_\infty^\circ(x) =0.
\end{align}
This ordinary first order differential equation is easy to solve with the result
\begin{align}
f^\circ_\infty(x) = c (1+x^2)^{-s},
\end{align}
where $c\in \reals$ is an integration constant. It is useful to normalise by putting $c=1$. The definition can also be extended to include the point at infinity $x=\infty$ by letting $f^\circ_\infty(\infty)=1$.

For $SL(2,\rats_p)$ one has to use the action~\eqref{eq:SL2pInd} and study sphericality at the group level. Invariance under $SL(2,\ints_p)$ is given by the $p$-adic spherical vector
\begin{align}
f_p^\circ(x) = \left\{\begin{array}{cl}
1&\textrm{if $x\in\ints_p$}\\
|x|_p^{-2s} &\textrm{if $x\notin \ints_p$}
\end{array}\right.
\end{align}
This corresponds to the inducing character.

A global spherical vector can then be defined by
\begin{align}
\label{eq:sphA}
f^\circ(x) = \left[\prod_{p<\infty} f_p^\circ(x)\right] f^\circ_\infty(x).
\end{align}

It is also interesting to consider the spherical vector in the Fourier transformed version in the line model. According to section~\ref{sec:adelisation}, the Fourier transformation is defined generally by
\begin{align}
\tilde{f}(u) = \int_{\bar{\mathbb{A}}} e^{-2\pi i u x} f(x) dx.
\end{align}
The Fourier transform of the individual pieces occurring in the global spherical vector~\eqref{eq:sphA} are
\begin{subequations}
\label{FourierSpherical}
\begin{align}
\tilde{f}_\infty^\circ(u) &= \frac{2\pi^s}{\Gamma(s)} |u|^{s-1/2} K_{s-1/2}(2\pi|u|),\\
\tilde{f}_p^\circ(u) &= (1-p^{-2s}) \gamma_p(u) \frac{1-p^{-2s+1} |u|^{2s-1}}{1-p^{-2s+1}},
\end{align}
\end{subequations}
where we have used the results from example~\ref{besselfull} for $p<\infty$ and the standard real Fourier transform of the Bessel function. 

\end{example}

We also determine the spherical vector in the conformal realisation of the orthogonal group that was introduced in example~\ref{ex:CFT}.

\begin{example}[Spherical vector for $SO(d,d)$ in conformal realisation]

In the big cell, functions in the induced representation~\eqref{eq:conf} depend on the $2d-2$ variables $x^i$ that are in the vector representation of $SO(d-1,d-1)$. It is useful to split the $2d-2$ variables as $x^a,x^\ba$ where $a=1,\ldots,d-1$ and $\ba=1,\ldots,d-1$ and both span a Euclidean space. Sphericality then implies that the spherical vector $f^\circ_\infty$ only depends on the Euclidean norms $r=||x^a||$ and $\br=||x^\ba||$.  This makes it sufficient to consider $f^\circ_\infty(r,\br)$ that satisfies
\begin{align}
x_a (P^a - K^a) f^\circ_\infty(r,\br) =0,
\quad  x_\ba (P^\ba + K^\ba)  f^\circ_\infty(r,\br) = 0.
\end{align}
The spherical vector is
\begin{align}
f^\circ_\infty(r,\br) = \left(1 + r^4 + \br^4 - 2 r^2 \br^2 + 2 r^2 + 2 \br^2\right)^s = \left((1+\br^2-r^2)^2+4.r^2\right)^s
\end{align}
It does not depend on $d$ and the argument of the power is manifestly positive. We have again chosen a convenient normalisation. This spherical vector corresponds to the inducing character. 
\end{example}

We will now outline an approach to Eisenstein series (or automorphic forms more generally) that was advocated in~\cite{Kazhdan:2001nx,KazhdanPolishchuk,Pioline:2003bk}. The setting is global and requires a group $G(\ads)$ with a discrete subgroup $G(\rats)$ as well as a induced representation $\Ind_P^G\sigma$. The approach for constructing an automorphic form is to form the following pairing
\begin{align}
\varphi(g)  = \langle f_\rats , \rho(g) f^\circ \rangle,
\end{align}
where the objects involved in this expression are the following
\begin{itemize}
\item $\langle\cdot,\cdot\rangle$ is an inner product on the function space $\Ind_P^G\sigma$ of an induced representation. In a given model of the coset $P\backslash G$, this is an integral over the product of two functions.
\item $f_\rats$ represents a distribution on the function space $\Ind_P^G\sigma$ that is invariant under the discrete subgroup $G(\rats)$: $f_\rats(\gamma x) = f_\rats(x)$ for all $\gamma\in G(\rats)$.
\item $f^\circ$ is a spherical vector in the representation space. It is acted upon by the group $G(\ads)$ by right translation $\rho(g)$, cf.~\eqref{eq:rra5}. 
\end{itemize}
The function $\varphi(g)$ thus defined is manifestly right $K$-invariant and left $G(\rats)$-invariant.

In more mundane terms, if $x$ are the coordinates on $P\backslash G$ and $d\mu(x)$ an appropriate invariant measure on the induced function space, then
\begin{align}
\label{eq:sphc}
\varphi(g) = \lint_{P\bs G} f_\rats(x) f^\circ(xg) d\mu(x).
\end{align}
It is important here that one uses the full coset space and not only the big cell.

\begin{example}[Reconstructing the $SL(2,\ads)$ Eisenstein series from spherical vectors]
Let us carry this out for $SL(2,\ads)$ using the ingredients already developed in example~\ref{ex:SL2sph} above. Additionally, we require the $SL(2,\rats)$ invariant distribution. Formally, this is given by~\cite{Kazhdan:2001nx}
\begin{align}
f_\rats(x) = \sum_{q\in\bar{\rats}} \delta(x-q)
\end{align}
in terms of the standard $\delta$-distribution. The summation here is over the compactified $\bar{\rats}$ that also includes the point at infinity corresponding to the small cell. The measure in the line model is simply $d\mu(x) = dx$. Plugging this into the definition~\eqref{eq:sphc} leads to
\begin{align}
\label{eq:phiSV}
\varphi(g) = \sum_{q\in \bar{\rats}} f^\circ(qg).
\end{align}
Restricting to $g\in SL(2,\reals)$ and using~\eqref{eq:sphA} leads to
\begin{align}
\varphi(g) = \sum_{q\in \bar{\mathbb{Q}}} \left[\prod_{p<\infty} f_p^\circ (q)\right] \frac{y^s}{ ((x+q)^2 +y^2)^{s}}
\end{align}
Here, we have used the coordinate $z=x+iy$ from section~\ref{sec_frommodtoaut} on $SL(2,\reals)/SO(2)$.
We can next split the sum into $q=\frac{m}{n}\in \rats$ (with $m$ and $n$ coprime and $n>0$) and $q=\infty$ related to the small cell. Using that
\begin{align}
\prod_{p<\infty} f_p^\circ\left(\frac{m}{n}\right) = n^{-2s}
\end{align}
we then obtain
\begin{align}
\label{ToEisen}
\varphi(g) &= \sum_{(m,n)=1\atop n>0} n^{-2s} \frac{y^s}{((x+m/n)^2+y^2)^{s}} + y^{s} 
= \sum_{(m,n)=1\atop n>0} \frac{y^s}{( (n x+m)^2+(ny)^2)^{s}} + y^{s} \nn\\
&= \frac12 \sum_{(m,n)=1} \frac{y^s}{|nz+m|^{2s}} = E(\chi_s,g).
\end{align}
This agrees precisely with the standard definition~\eqref{Eisenintro2} of the $SL(2,\reals)$ Eisenstein series, including all constant terms. The extra term $y^s$ is exactly the contribution from the point at infinity in the line model (the small Bruhat cell) since $f_\infty^\circ(\infty g) = y^{s}f_\infty^\circ(\infty) = y^s$ by using the transformation of the point at infinity and the value of the spherical vector there.
\end{example}

\chapter[Whittaker functions and Fourier coefficients]{Whittaker functions and\\ Fourier coefficients}
\label{ch:fourier}

In this chapter, we analyse the general structure of the Fourier expansions of automorphic forms, with particular emphasis on Eisenstein series and the associated theory of Whittaker functions. We will discuss both local and global aspects. As advanced topics we introduce the useful notion of wave-front set~\cite{Matumoto,MoeglinWaldspurger,MoeglinWaveFront,MoeglinHoweWaveFront} and discuss the method of Piatetski-Shapiro and Shalika~\cite{PiatetskiShapiro,Shalika}.
General references are~\cite{HJacquet,Bump,Goldfeld} and we also found the discussions in~\cite{Green:2011vz,MillerSahi,GinzburgHundley,KobayashiSavin,JiangLiuSavin} very useful.

\section{Preliminary example: \texorpdfstring{$SL(2, \mathbb{R})$}{SL(2, R)} Whittaker functions}
\label{sec_FCsl2ex}

In section \ref{intro:Fourier}, we discussed the Fourier expansion of the non-holomorphic Eisenstein series $E(s, z)$ where $z=x+iy$ is on the upper half plane $\UHP$. Invariance under 
$SL(2,\mathbb{Z})$ implies the periodicity of the series in the real $x$-direction: 
\begin{align}
\label{SL2period}
E(s, x+1+iy)=E(s, x+iy),
\end{align}
and hence we have a Fourier expansion of the form 
\beq
E(s, x+iy)=\sum_{m\in \mathbb{Z}} a_m(y) e^{2\pi i mx},
\label{EisExpSec4}
\eeq
where the $y$-dependent Fourier coefficients $a_m(y)$ can be  extracted from the explicit expansion stated in~\eqref{SL2FC2} and will be derived in detail in chapter~\ref{ch:SL2-fourier} using adelic methods. Let us now reinterpret $E(s,z)$ as a function on $SL(2,\mathbb{R})=N(\mathbb{R})A(\mathbb{R})K(\mathbb{R})$ according to the prescription in section \ref{Maass}. To this end we define
\beq
\varphi_E(g) = \varphi_E(nak)=\varphi_E \left( \begin{pmatrix} 1 & x \\ & 1  \end{pmatrix} \begin{pmatrix} y^{1/2} & \\ & y^{-1/2} \end{pmatrix} \left(\begin{array}{cc}
\cos\theta & \sin \theta \\
-\sin \theta & \cos \theta \\
\end{array} \right)\right) =E(s,x+iy),
\eeq
where $n\in N(\mathbb{R})$, $a\in A(\mathbb{R})$, $k\in K(\mathbb{R})=SO(2,\mathbb{R})$. From this point of view, the periodicity~\eqref{SL2period} of $E(s,z)$ in the variable $x$ is equivalent to the invariance of $\varphi_E(g)$ under discrete left-translations: $\varphi_E(ng)=\varphi_E(g)$, $n\in N(\mathbb{Z})$. This follows from the simple calculation for $n=\begin{psmallmatrix}1&1\\&1\end{psmallmatrix}$:
\begin{align}
\varphi_E\left(\begin{pmatrix}1&1\\&1\end{pmatrix}g\right) &= \varphi_E \left( \begin{pmatrix} 1 & x+1 \\ & 1  \end{pmatrix} \begin{pmatrix} y^{1/2} & \\ & y^{-1/2} \end{pmatrix} \left(\begin{array}{cc}
\cos\theta & \sin \theta \\
-\sin \theta & \cos \theta \\
\end{array} \right)\right) \nn\\
& =E(s,x+1+iy),
\end{align}
which equals $\varphi_E(g)= E(s,x+i y)$ by left $N(\ints)$-invariance.

More generally, we can consider an automorphic form $\varphi$ on $SL(2,\mathbb{Z})\backslash SL(2,\mathbb{R})$, satisfying 
\beq
\varphi(\gamma g k)=\sigma(k) \varphi(g), \qquad \gamma\in SL(2,\mathbb{Z}), \, k\in K(\reals)=SO(2,\reals),
\eeq
where $\sigma$ can be  a non-trivial finite-dimensional representation of $K(\reals)$. When $\sigma$ is non-trivial, the function $\varphi $ depends on all three coordinates $(x,y, \theta)$. When $\sigma$ is trivial and hence $\varphi$ independent of $k$, the function is spherical.

The automorphy of $\varphi$ includes invariance under $N(\ints)$ and therefore $\varphi(g)=\varphi(x, y, \theta) $ will have a Fourier expansion of the same form as the one for $E(s,z)$, although the precise coefficients will of course be different depending on the choice of $\varphi$. To pave the way for higher rank groups, we now wish to recast this expansion in a form that can be easily generalised.

To this end, let $\psi : N(\mathbb{Z})\backslash N(\mathbb{R}) \to U(1)$ be a \emphindex[character!unitary]{unitary multiplicative character} on $N(\reals)$ which is trivial on $N(\mathbb{Z})$. The space of such characters is $\text{Hom}\left(N(\mathbb{Z}) \backslash N(\mathbb{R}), U(1)\right)\cong \mathbb{Z}$ and we can  parametrise the choice of character  by a single integer $m$ via 
\beq
\psi\left(\begin{pmatrix} 1 & x \\  & 1 \end{pmatrix} \right)= e^{2\pi i m x}, \qquad m\in \mathbb{Z}, \, x\in \mathbb{R}.
\label{sl2ch}
\eeq 
This is therefore nothing but a set of Fourier modes. If $\psi$ is \emphindex[character!non-trivial]{non-trivial}, i.e. $m\neq 0$, we say that $\psi$ is \emphindex[character!generic]{generic}. For higher rank groups, if a character is non-trivial, it does not necessarily mean that it is also generic. In definition~\ref{def:GC} we will extend our concept of this notion to the case of higher rank groups by introducing a more refined notion of generic vs. non-generic (or degenerate) characters.

Then, due to the periodicity of the automorphic form, $\varphi(ng)=\varphi(g)$, $n\in N(\mathbb{Z})$, we can write $\varphi(g)$ as a Fourier expansion along $N(\mathbb{R})$:
\beq
\varphi(g)=\sum_{\psi\in \textrm{Hom}(N(\mathbb{Z})\backslash N(\reals),U(1))} W_\psi(g)\,,
\label{SL2WhittakerExpansion}
\eeq
where the sum runs over all possible characters $\psi$ and hence over $m\in\ints$. We have also defined the \emphindex{Whittaker--Fourier coefficient} 
\beq
W_\psi(g)=\int_{N(\mathbb{Z})\backslash N(\mathbb{R})} \varphi(ng)\overline{\psi(n)}dn,
\label{SL(2)Whittaker}
\eeq
with $dn$ the \emphindex{Haar measure} on $N(\mathbb{Z})\backslash N(\mathbb{R})$. The Haar measure is normalised such that $\int_{N(\mathbb{Z})\bs N(\reals)}dn=1$. 

\begin{remark}
For brevity we shall often refer to $W_\psi(g)$ simply as a \emph{Whittaker coefficient}, in place of the more awkward \emph{Whittaker--Fourier coefficient}.
\end{remark}

The expansion~\eqref{SL2WhittakerExpansion} is a reformulation of~\eqref{EisExpSec4} as we will now illustrate.
By the Iwasawa decomposition $g=nak$ it follows that $W_\psi(g)$ is determined by its restriction to $A(\mathbb{R})$: 
\begin{align}
W_\psi (nak)&=\int_{N(\mathbb{Z})\backslash N(\mathbb{R})} \varphi(n'nak)\overline{\psi(n')}dn'\nn \\
&= \sigma(k)\int_{N(\mathbb{Z})\backslash N(\mathbb{R})} \varphi(\tilde{n}a)\overline{\psi(\tilde{n} n^{-1})}d\tilde{n}\nn \\
&= \psi(n) \sigma(k) W_\psi(a), 
\label{WhittakerRestSL2}
\end{align}
where we used the multiplicativity of $\psi$  as well as the invariance of the Haar measure under translations by $N(\mathbb{R})$.  In particular, this allows us to rewrite the expansion in a way that is more akin to the classical form~\eqref{EisExpSec4}: 
\beq
\varphi(g)=\sum_{\psi} W_\psi(ak) \psi(n).
\label{eq:SL2expansion}
\eeq
Note that contrary to standard \emphindex[harmonic analysis]{harmonic analysis} the function $W_\psi(g)$ is not a numerical coefficient, but also contains explicitly the Fourier variable(s) that one is expanding in. This is made explicit by the factor $\psi(n)$ appearing in~\eqref{WhittakerRestSL2} and~\eqref{eq:SL2expansion}.

By an explicit Iwasawa parametrisation of $g$ in terms of the variables $(x, y, \theta)$ as in (\ref{IwasawaSL2R}) and the character $\psi$ in terms of an integer~\eqref{sl2ch}, the integral~\eqref{SL(2)Whittaker} takes the more familiar form
\beq
W_\psi(ak)=W_m(y, \theta)=\int_{0}^{1} \varphi(x, y, \theta) e^{-2\pi i mx} dx.
\eeq 
(For trivial $\sigma(k)$ the integral is independent of $\theta$ and equal to $a_m(y)$ of (\ref{EisExpSec4}).)

The general $SL(2,\reals)$ expansion~\eqref{eq:SL2expansion} contains two types of terms, corresponding to $m=0$ ($\psi=1$) and $m\neq0$ ($\psi\neq1$) that are useful to distinguish:

\begin{definition}[Constant terms and Whittaker coefficients for $SL(2,\reals)$]
The sum~\eqref{eq:SL2expansion} can be split into
\beq
\varphi(g)=W_1(ak)+\sum_{\psi\neq1} W_\psi(g),
\eeq
where the first term is independent of $n$ and called the \emphindex{constant term}. It is defined by
\beq
W_1(ak)=\int_{N(\mathbb{Z})\backslash N(\mathbb{R})} \varphi(nak) dn.
\eeq
and we will sometimes also denote it by $C(ak)\equiv W_1(ak)$. The functions $W_\psi(g)$ for non-trivial characters $\psi$ ($m\neq 0$) are the proper Whittaker coefficients.
\end{definition}
 
\begin{remark}
The functions $W_\psi(g)$ were termed Whittaker functions by Jacquet~\cite{HJacquet} because they reduce to the classical Whittaker function $W_{k,m}(y)$ for the group $GL(2,\reals)$. For $SL(2,\reals)$ they are given basically by modified Bessel functions that arise for the special case $k=0$; see also appendix~\ref{app:SL2Laps}. For higher rank groups $G(\reals)$, Whittaker functions define more complicated special functions. For example, for $GL(n,\reals)$ it was shown in~\cite{Stade1,Stade2} how these generalised Whittaker functions can be obtained as nested integrals over standard Whittaker functions.  We will encounter an instance of this for $SL(3,\reals)$ in sections~\ref{sec:FCsec} and~\ref{sec:SL3ex} (see, e.g.,~\eqref{eq:SL3-generic-archimedean}). The general asymptotics of these functions near a cusp will be discussed in section~\ref{sec:asympt} where also their relation to string theory effects is mentioned.
\end{remark}

\begin{example}[Fourier and $q$ expansion of holomorphic Eisenstein series]
Consider now the example when $\varphi=\varphi_f$ with $f(z)=E_{2w}(z)$ being a weight $2w$ holomorphic Eisenstein series
\beq
E_{2w}(z)=\frac12\sum_{(c,d)=1} \frac{1}{(cz+d)^{2w}}, \qquad \quad z\in \mathbb{H}, \, \, w>1, \,\,w\in\mathbb{Z}.
\eeq
This function is spherical ($\theta$-independent) and has a well-known Fourier expansion (sometimes called \emphindex[q-expansion@$q$-expansion]{$q$-expansion})
\beq
E_{2w}(z)=1+\sum_{m=1}^{\infty} a_m q^{m}, \qquad q=e^{2\pi i z},
\eeq
where the coefficients are given by
\beq
a_m=\frac{2}{\zeta(1-2w)} \sigma_{2w-1}(m), 
\eeq
with $\sigma_{2w-1}(m)$ the \index{divisor sum} sum over positive divisors as in~\eqref{intro:divisorsum} 
\beq
\sigma_s(m)=\sum_{d|m} d^s.
\eeq
In this case the constant term and Whittaker coefficients are given by 
\begin{align}
C(a) &\equiv W_1(a)= 1,
\nn \\
W_\psi(z) &\equiv W_\psi(na)= W_m(z)=\frac{2}{\zeta(1-2w)} \sigma_{2w-1}(m)q^m, \qquad m>0.
\end{align}
The coefficients can be alternatively be expressed in terms of \emphindex{Bernoulli numbers}, see for example~\cite{Apostol2}.
Notice that the holomorphicity of $E_{2w}(z)$ requires that $W_m(z)$ vanishes unless $m>0$. As mentioned in section~\ref{standardsectionSL2}, this is due to the holomorphic Eisenstein series' $E_{2w}$ being associated with the discrete series representation of $SL(2,\reals)$. (There are other common normalisations of holomorphic Eisenstein series where the constant term is not given by $1$ but by $2\zeta(2w)$.)
\end{example}

For completeness, we also recall the constant terms and Whittaker coefficients for the non-holomorphic Eisenstein series $E(s,z)$ on $SL(2,\reals)$ from the introduction.

\begin{example}[Fourier expansion of non-holomorphic Eisenstein series]
In the case when $\varphi=\varphi_E$, with  $E(s, z)$ the non-holomorphic Eisenstein series on $\mathbb{H}$, the constant term $W_1(a)$ and Whittaker coefficient $W_\psi(na) $ will be derived in chapter~\ref{ch:SL2-fourier} with the result
\begin{equation} 
    \begin{split}
        W_1(a) &= W_1(y)=y^s + \frac{\xi(2s-1)}{\xi(2s)} y^{1-s} \\
        W_\psi(na) &= W_m(x, y)= \frac{2y^{1/2} }{\xi(2s)} |m|^{s-1/2} \sigma_{1-2s}(m) K_{s-1/2} (2\pi |m| y) e^{2\pi i m x}, \quad \text{ with } m\neq 0. \\
    \end{split}
\end{equation}
In contrast to the holomorphic case,  the `constant term' here is not really constant; it is a function on the Cartan torus that is parametrised by the imaginary part $y$ of $z=x+iy$. As we will see below, this is in fact a general feature, namely the constant term of a spherical automorphic function is only constant with respect to the coordinates along the unipotent radical $N(\mathbb{R})$ of the Borel subgroup $B(\mathbb{R})\subset G(\mathbb{R})$.
\end{example}

\section{Fourier expansions and unitary characters}
\label{sec:FCsec}

We now turn to the general analysis of Fourier coefficients of automorphic forms on semi-simple Lie groups $G$, and we also switch to the adelic framework. For this we first require the notion of a unitary character $\psi$ on a unipotent subgroup $U\subset G$ that generalises the Fourier mode $e^{2\pi i m x}$ in~\eqref{sl2ch}. This is discussed in detail in section~\ref{sec:UC}. We will then discuss the notion of Fourier expansion for different types of unipotent groups $U$ in the sequel. 

\subsection{Unitary characters}
\label{sec:UC}

\begin{definition}
\label{def:UC}
Let $U(\ads)$ be a unipotent subgroup of the adelic group $G(\ads)$. A \emphindex[character!unitary|textbf]{unitary character} on $U(\ads)$ is a group homomorphism
\begin{align}
\psi : U(\rats)\bs U(\ads) \to U(1)
\end{align}
and we also require it to be trivial on the discrete subgroup $U(\rats)=U(\ads)\cap G(\rats)$ since we will study in the context of automorphic forms on $G(\ads)$ that are invariant under the discrete subgroup $G(\rats)$. The space of all unitary characters on $U(\ads)$ that are trivial on $U(\rats)$ are called the integral points of the \emphindex[character!variety]{character variety}.
\end{definition}

\begin{remark}
Unipotent groups are required if one wants to have non-trivial unitary characters. On the simple group $G(\ads)$ there are no non-trivial unitary characters.
\end{remark}

Definition~\ref{def:UC} generalises~\eqref{sl2ch}. As $\psi$ is a group homomorphism to the abelian group $U(1)$, it is trivial on the \emphindex{commutator subgroup}
\begin{align}
[U,U] = \left\{ u_1u_2u_1^{-1} u_2^{-1} \st u_1, u_2 \in U\right\}.
\end{align}
In other words, 
\begin{align}
\label{psitriv}
\psi([U,U]) =1,
\end{align}
such that $\psi$ is sensitive only to the \emphindex[abelianisation of a unipotent group]{abelianisation} $[U,U]\bs U$. We note that $[U,U]$ equals the second member of the derived series of $U$ defined in section~\ref{sec:simple-lie-alg}. We will discuss the relevance of the derived series for Fourier expansions in more detail below in section~\ref{sec:ANA}.

It is convenient to have a more explicit parametrisation of possible unitary characters $\psi$. To this end we restrict to the case where $U$ is the unipotent of a standard parabolic subgroup $P=LU$ as defined in section~\ref{sec:parsubgp}. As always we are working with a fixed choice of split Cartan torus $A\subset G$. Such unipotent groups $U$  can be generated from the product of \emphindex{one-parameter subgroups}
\begin{align}
U_\alpha = \big\{ x_\alpha(u_\alpha) = \exp(u_\alpha E_\alpha) \st u_\alpha \in \ads\big\},
\end{align}
with $\alpha$ ranging over the subset $\Delta(\mf{u})$ of positive roots of $\mf{g}$ corresponding to the Lie algebra $\mf{u}$ of $U$:
\begin{align}
U = \prod_{\alpha\in \Delta(\mf{u})} U_\alpha.
\end{align}
The restriction of $\psi$ to any of the one-parameter subgroups $U_\alpha$ then yields a unitary character
\begin{align}
\psi_\alpha : U_\alpha(\rats) \bs U_\alpha(\ads) \to U(1) .
\end{align}
As any one-parameter subgroup $U_\alpha$ is abelian and satisfies the isomorphism
\begin{align}
U_\alpha(\rats) \bs U_\alpha(\ads) \cong \rats\bs \ads,
\end{align}
the unitary character $\psi_\alpha$ can therefore be parametrised by a rational number $m_\alpha\in\rats$ as discussed in section~\ref{sec_adeles}, see also~\cite[Thm 5.4.3]{Deitmar}, and can be thought of as the global function
\begin{align}
\label{adsUC}
\psi_\alpha\left( x_\alpha(u_\alpha) \right) = e^{2\pi i m_\alpha u_\alpha}
\end{align}
and we will sometimes refer to the $m_\alpha$ as \emphindex{mode numbers} or \emphindex[instanton!charge|textbf]{instanton charges} as this is their interpretation in a string theory context, see chapter~\ref{ch:intro-strings}. 

The triviality~\eqref{psitriv} of $\psi$ can then be restated as
\begin{align}
\psi \left( \prod_{\alpha\in \Delta([\mf{u},\mf{u}])} U_\alpha \right) =1
\end{align}
and the non-trivial unitary characters are therefore sensitive only to the one-parameters subgroups $U_\alpha$ such that $\alpha$ is a `root' of $\mf{u}$ but not of $[\mf{u},\mf{u}]$. This means that the parametrisation of different unitary characters $\psi$ on $U$ only requires the knowledge of the mode numbers $m_\alpha$ for the positive roots $\alpha$ that belong to $\Delta(\mf{u})$ but not to $\Delta([\mf{u},\mf{u}])$. We define
\begin{align}
\Delta^{(1)}(\mf{u}) := \Delta(\mf{u}) \setminus \Delta([\mf{u},\mf{u}])
\end{align}
to be these roots. 
\begin{remark}
The notation $\Delta^{(1)}(\mf{u})$ indicates that these are the `roots' of the abelianisation $[U,U]\bs U$ of the degree one piece $U=U^{(1)}$ in the \emphindex{derived series} of $U$. See section~\ref{sec:ANA} for a more detailed discussion of the relevance of the derived series of $U$ for Fourier expansions and section~\ref{sec:simple-lie-alg} for the notion of derived series.
\end{remark}

The above considerations lead to 
\begin{proposition}[Parametrisation of unitary characters]
\label{prop:UC}
Let $U(\ads)$ be a unipotent subgroup of $G(\ads)$. Unitary characters $\psi: U(\rats)\bs U(\ads) \to U(1)$ can be parametrised uniquely by a set of mode numbers $\left\{ m_\alpha \in \rats \st \alpha \in \Delta^{(1)}(\mf{u})\right\}$. The unitary character is then given by
\begin{align}
\label{eq:UCIC}
\psi\left( \prod_{\alpha\in\Delta^{(1)}(\mf{u})} x_\alpha(u_\alpha)\right) = \exp\left( 2\pi i \sum_{\alpha \in\Delta^{(1)}(\mf{u})} m_\alpha u_\alpha\right).
\end{align}
It factorises into local places as in~\eqref{Ach_fact}.
\end{proposition}

\begin{proof}
The triviality of $\psi$ on the commutator subgroup $[U,U]$ shows that it suffices to define $\psi$ on the abelianisation that is constructed from the one-parameter subgroups $U_\alpha$ with $\alpha\in\Delta^{(1)}(\mf{u})$ for which the characters were determined in~\eqref{adsUC} above. The group homomorphism property of $\psi$ then yields the proposition.
\end{proof}

The following notions will be important in the sequel.
\begin{definition}[Generic and degenerate characters]
\label{def:GC}
Let $\psi : U(\mathbb{Q})\backslash U(\mathbb{A}) \to U(1)$ be a  \emph{global} character as in~\eqref{eq:UCIC}. 
\begin{enumerate}
\item[$(i)$] $\psi$ is called \emphindex[character!generic|textbf]{generic} if $m_\alpha\neq0$ for all $\alpha\in\Delta^{(1)}(\mf{u})$, i.e. if the character is non-trivial on each one-parameter subgroups $U_\alpha(\mathbb{A})$ for $\alpha\in\Delta^{(1)}(\mf{u})$. 
\item[$(ii)$] If $m_\alpha=0$ for all $\alpha\in\Delta^{(1)}(\mf{u})$, the character $\psi$ is called \emphindex[character!trivial|textbf]{trivial}.
\item[$(iii)$] Furthermore, if $m_\alpha\neq0$ for at least one, but not all, $\alpha\in\Delta^{(1)}(\mf{u})$, the character $\psi$ is called \emphindex[character!non-generic|textbf]{non-generic} or \emphindex[character!degenerate|textbf]{degenerate}.
\end{enumerate}
\end{definition} 

We illustrate these notions by the following example.
\begin{example}[Unitary characters on the maximal unipotent of $SL(n,\ads)$]
\label{SLnUC}
Consider the case $G=SL(n,\ads)$ and $U(\ads)=N(\ads)$ to be the (maximal) unipotent subgroup of the Borel subgroup $B(\ads)$, implying $\mf{n}=\mf{u}$. The set $\Delta(\mf{n})$ is given by all positive roots $\Delta_+$ of $\mf{sl}(n)$ and the set $\Delta^{(1)}(\mf{n})$ equals the ($n-1$) simple roots $\Pi\subset \Delta_+$. In the fundamental representation we can write elements of $n\in N$ as $(n\times n)$-matrices of the form
\begin{align}
n = \begin{pmatrix}
1 & u_1 & * & *& \cdots\\
& 1 & u_2 &* &\cdots\\
&& \dots &&\\
&&&1&u_{n-1}\\
&&&&1
\end{pmatrix}.
\end{align}
The starred entries are of no relevance for the discussion of unitary characters as they are associated with the commutator subgroup $[N,N]$. A character $\psi$ on $N$ is determined by $n-1$ rational numbers $m_i$ ($i=1,\ldots,n-1$) such that
\begin{align}
\psi(n) = \exp(2\pi i \sum_{i=1}^{n-1} m_i u_i).
\end{align}
The character $\psi$ is generic when all $m_i\neq 0$. It is degenerate when some $m_i$ vanish and then it does not depend on the corresponding one-parameter subgroups.
\end{example}

 We recall from section~\ref{sec_adeles} that a global unitary character $\psi_\alpha$ on $\rats\bs\ads$ as in ~\eqref{adsUC} factorises as 
\beq
\psi_\alpha=\prod_{p\leq \infty} \psi_{\alpha,p}\,,
\eeq
where for $p<\infty$
\begin{align} 
\psi_{\alpha,p} &: U(\mathbb{Z}_p)\backslash U(\mathbb{Q}_p)  \rightarrow U(1),& \psi_{\alpha,p}(x_\alpha(u)) = e^{-2\pi i [m_\alpha u]}
\end{align}
in terms of the fractional part~\eqref{fracpart} of a $p$-adic number, and for $p=\infty$
\begin{align}
\psi_\infty &: U(\mathbb{Z})\backslash U(\mathbb{R})\rightarrow U(1),& \psi_{\alpha,p}(x_\alpha(u)) = e^{2\pi i m_\alpha u}.
\end{align}
This factorisation extends to characters $\psi$ on unipotent groups $U$:
\begin{align}
\psi = \prod_{p\leq \infty} \psi_p.
\end{align}
Definition~\ref{def:GC} extends to all local characters $\psi_p$. Moreover, we have the following notion:

\begin{definition}[Unramified unitary character]
\label{def_unram}A generic \emph{local} character $\psi_{p}$ for $p<\infty$ is called \emphindex[character!unramified]{unramified} if for all  $\alpha\in\Delta^{(0)}(\mf{u})$ one has 
\beq
\psi_{\alpha,p}\left( e^{u E_\alpha}\right)=e^{-2\pi i [ u]}, \qquad   u\in \mathbb{Q}_p.
\eeq 
Equivalently, this means that all instanton charges $|m_\alpha|_p=1$ in (\ref{eq:UCIC}). We call a global character unramified if $m_\alpha=1$ for all $\alpha$.
\end{definition}

\subsection{General Fourier coefficients vs. Whittaker coefficients}
\label{sec:WhF}

Now that we have the Fourier modes in terms of characters $\psi$ on unipotent subgroups $U$, it is possible to define Fourier coefficients of automorphic forms. 

\begin{definition}[Fourier coefficient]
\label{def:F}
Let $\varphi$ be an automorphic form on $G(\ads)$, i.e., an element of the space $\mathcal{A}(G(\rats)\bs G(\ads))$, and $U(\ads)$ a unipotent subgroup of $G(\ads)$. The \emphindex[Fourier coefficient!of an automorphic form]{Fourier coefficient} of $\varphi$ with respect to the unitary character $\psi$ on $U$ is given by:
\begin{align}
\label{eq:FCU}
F_\psi(\varphi,g) = \lint_{U(\rats)\bs U(\ads)} \varphi(ug) \overline{\psi(u)} du,
\end{align}
 where $du$ denotes the invariant Haar measure on $U$.
\end{definition}
\begin{remark}
The Fourier coefficient $F_\psi$ can be viewed either as a function $F_\psi(g)$ on $G(\ads)$ for fixed $\varphi$, or as a functional $F_\psi(\varphi)$ on $\mathcal{A}(G(\rats)\bs G(\ads))$. When it is clear from the context which fixed $\varphi$ is meant, we may write simply $F_\psi(g)$ for conciseness.
\label{FourierCoefficientRemark}
\end{remark}

A short calculation similar to~\eqref{WhittakerRestSL2} shows that Fourier coefficients satisfy
\begin{align}
F_\psi (\varphi,ug) = \psi(u) F_\psi(\varphi,g) \quad\quad \textrm{for all $u\in U$.}
\end{align}

We make the additional definitions for the case $U(\ads)=N(\ads)$.
\begin{definition}[Whittaker coefficient]
\label{def:Wh}
Let $\varphi$ be an automorphic form on $G(\ads)$, $N(\ads)$ be the maximal unipotent subgroup of a fixed Borel $B(\ads)$ and $\psi$ be a unitary character on $N(\ads)$. 
\begin{itemize}
\item[$(i)$]The integral
\begin{align}
\label{abelianFC}
W_\psi(\varphi,g) = \lint_{N(\rats)\bs N(\ads)} \varphi(ug) \overline{\psi(u)} du.
\end{align}
is called the \emphindex{Whittaker coefficient} of $\varphi$ with respect to $\psi$. 
\item[$(ii)$]
If $\varphi$ is $K_\ads$ invariant, the Whittaker coefficient $W_\psi(\varphi,g)$ is right-invariant under $K_\ads$ and the Whittaker coefficient is then called \emphindex[Whittaker coefficient!spherical]{spherical}. We denote it by $W^\circ_\psi(\varphi,g)$. The spherical Whittaker coefficient is completely determined by its values on the Cartan torus $A(\ads)$: Writing $g=nak$ in Iwasawa decomposed form one has
\begin{align}
W^\circ_\psi(\varphi,nak) = \psi(n) W^\circ_\psi(\varphi,a).
\end{align}
This is the case for Eisenstein series.
\end{itemize}
\end{definition}

\begin{remark}
Even though definition~\ref{def:Wh} is a special case of~\ref{def:F}, it is useful to distinguish this case notationally. Throughout this work, we will denote Whittaker coefficients (i.e., Fourier coefficients along the maximal unipotent $N$)  by $W_\psi$ and reserve the notation $F_\psi$ for the case when the unipotent $U$ is different from $N$. Whittaker coefficients will be the main focus of this work and studied in detail in chapter~\ref{ch:Whittaker-Eisenstein} for Eisenstein series.
\end{remark}

We note that if $U$ is the unipotent of some standard parabolic subgroup $P=LU$ and $\varphi$ $K_\ads$-invariant, the general Fourier coefficient $F_\psi(\varphi,g)$ is determined by its values on the Levi subgroup $L$ and one could define a spherical Fourier coefficient $F^\circ_\psi$ but we will not make use of this notion. 

\begin{definition}[Constant term]
\label{CTdef}
\mbox{ }

\begin{enumerate}
\item[$(i)$]
The Fourier coefficient of an automorphic form $\varphi$ with respect to the trivial character $\psi=\id$ on $U$ is called the \emphindex[constant term!along $U$]{constant term along $U$}:
\begin{align}
C_U(\varphi,g) = \lint_{U(\rats)\bs U(\ads)} \varphi(ug) du.
\end{align}
It is independent of $u\in U$: $C_U(\varphi,ug)=C_U(\varphi,g)$. 
\item[$(ii)$]
For the case $U=N$, we will call it simply the \emphindex{constant term} and denote it by
\begin{align}
C(\varphi,g) \equiv C_N(\varphi,g) = \lint_{N(\rats)\bs N(\ads)} \varphi(ng) dn.
\end{align}
If $\varphi$ is spherical, the constant term is a function only of the Cartan torus $A(\ads)$: Using Iwasawa decomposition $C(\varphi,nak) = C(\varphi,a)$.
\end{enumerate}
\end{definition}

\subsection{Abelian vs. non-abelian Fourier expansions}
\label{sec:ANA}
\index{Fourier expansion!abelian}\index{Fourier expansion!non-abelian}

In the $SL(2)$ example of section~\ref{sec_FCsl2ex}, the Whittaker coefficient $W_\psi$ were used in~\eqref{SL2WhittakerExpansion} to give a \emph{complete} Fourier expansion of an automorphic form $\varphi$ by summing over all possible unitary characters $\psi$.  
\label{eq:non-abelian-Whittaker}

It is a natural question how this carries over to higher rank groups $G(\ads)$. In view of proposition~\ref{prop:UC}, we can already anticipate that the Fourier expansion with unitary characters $\psi$ on a unipotent group $U$ will in general be incomplete since the characters $\psi$ only depend on the abelianisation $[U,U]\bs U$, see~\eqref{psitriv}. For $SL(2)$ the (maximal) unipotent group $N$ is abelian and we did not have to consider this subtlety. The general statement is
\begin{proposition}[Partial Fourier sum]
\label{prop:NAFC}
Let $U(\ads)$ be a unipotent subgroup of $G(\ads)$ and $\varphi$ be an automorphic form on $G(\ads)$. Then the sum of Fourier coefficients over all unitary characters $\psi$ on $U$ yields
\begin{align}
\sum_\psi F_\psi(\varphi,g) = \lint_{[U,U](\rats)\bs [U,U](\ads)} \varphi(ug) du.
\end{align}
In other words, the sum of the Fourier coefficients reconstitutes only the average of the automorphic form over the commutator subgroup $[U,U]$. If $U$ is abelian, the Fourier expansion is complete.
\end{proposition}

\begin{proof}
See~\cite{MillerSahi}.
\end{proof}

In order to obtain a complete Fourier expansion when the unipotent $U$ is non-abelian, one has to consider the \emphindex{derived series} of $U$ (cf. also section~\ref{sec:simple-lie-alg}):
\begin{align}
\label{eq:DSeries}
U^{(i+1)} = [U^{(i)} , U^{(i)}] \,,\quad \quad U^{(1)} = U.
\end{align}
Since $U$ is unipotent, the derived series trivializes after finitely many steps: $U^{(i_0)}=\{\id\}$ for some $i_0 \geq 1$ and we assume $i_0$ to be the smallest integer for which $U^{(i_0)}=\{\id\}$. If $U$ is abelian, one has $i_0=2$. The successive quotients
$U^{(i+1)}\bs U^{(i)}$
are the abelianisations of the unipotent groups $U^{(i)}$ for any integer $i\geq 1$. A unitary character $\psi^{(i)}$ on $U^{(i)}$ is trivial on $U^{(i+1)}$. One can define Fourier coefficients for any of the $U^{(i)}$ by the same formula as in definition~\ref{def:F}:
\begin{align}
F_{\psi^{(i)}}(\varphi,g) = \lint_{U^{(i)}(\rats) \bs U^{(i)}(\ads)} \varphi(ug) \overline{\psi^{(i)}(u)} du.
\end{align}
As an immediate analogue of proposition~\ref{prop:NAFC} one has that
\begin{align}
\label{eq:CTUder}
\sum_{\psi^{(i)}} F_{\psi^{(i)}}(\varphi,g) = \lint_{U^{(i+1)}(\rats) \bs U^{(i+1)}(\ads)} \varphi(ug) du.
\end{align}
We observe that the right-hand side is nothing but the constant term of $\varphi$ along $U^{(i+1)}$, corresponding to $\psi^{(i+1)}=1$. It is therefore natural that the complete \emphindex[Fourier expansion!non-abelian]{non-abelian Fourier expansion} of $\varphi$ along $U$ is given by
\begin{align}
\label{eq:NAFCU}
\varphi(g) = C_U(\varphi,g) +\sum_{\psi^{(1)}\neq 1} F_{\psi^{(1)}}(\varphi,g) +\sum_{\psi^{(2)}\neq 1} F_{\psi^{(2)}}(\varphi,g) + \ldots +\sum_{\psi^{(i_0)}\neq 1} F_{\psi^{(i_0)}}(\varphi,g).
\end{align}
The trivial character $\psi^{(i)}=1$ is always excluded because the sum of the preceding terms reconstitutes the constant term along $U^{(i)}$ by~\eqref{eq:CTUder}. Note that unitary characters $\psi^{(1)}$ are characters on $U^{(1)}=U$ and therefore equal the unitary characters we have been discussing in definition~\ref{def:F}. We will sometimes refer to the Fourier coefficients in~\eqref{eq:NAFCU} associated with $U^{(i)}$ and $i\geq 2$ as \emphindex[Fourier coefficient!non-abelian]{non-abelian Fourier coefficients} and the ones associated with $U^{(1)}=U$  as the \emphindex[Fourier coefficient!abelian]{abelian Fourier coefficient}.

The same structure of the expansion and terminology arises for the case when the unipotent $U$ is given by the maximal unipotent $N$. Then we have
\begin{align}
\varphi(g)=\underbrace{\phantom{\sum_{\psi^{(1)}}}\!\!\!\!\!\!\!\!C(g)}_{\text{constant term}}+\underbrace{\sum_{\psi^{(1)} \neq1}W_{\psi^{(1)}}(g)}_{\text{abelian term}}+\underbrace{\sum_{\psi^{(2)}\neq 1} W_{\psi^{(2)}}(g)}_{\text{non-abelian term}}+\cdots\,,
\label{generalFourierExpansion}
\end{align}
where we have suppressed the fixed automorphic function $\varphi$ on the right-hand side.

\begin{remark}
Our main interest in this work lies with the abelian Whittaker coefficients $W_{\psi^{(1)}}$ and we will discuss them in more detail in the following sections and in particular in chapter~\ref{ch:Whittaker-Eisenstein}.
Non-abelian Fourier expansions have been carried out in detail for $SL(3,\reals)$ in~\cite{Bump,VT,Narita,Pioline:2009qt} and this will be reviewed in section~\ref{sec:SL3ex}. Non-abelian Fourier expansions for the non-split real group $SU(2,1)$ can be found in~\cite{Ishikawa,Bao:2009fg} and  some further comments on the non-abelian coefficients will be offered in chapter~\ref{ch:outlook}. 

\end{remark}

\section{Induced representations and Whittaker models}

We now specialise to the case when the automorphic form  $\varphi\in \mathcal{A}(G(\mathbb{Q})\backslash G(\mathbb{A}))$ is an Eisenstein series
\beq
\label{eq:ESWh}
E({\sf f}_\lambda, g)=\sum_{\gamma\in B(\mathbb{Q})\backslash G(\mathbb{Q})} {\sf f}_\lambda (\gamma g), \qquad g\in G(\mathbb{A}),
\eeq
constructed from a standard section ${\sf f}_\lambda$ of the (in general, non-unitary) principal series $\text{Ind}_{B(\mathbb{A})}^{G(\mathbb{A})} \chi$, cf. section~\ref{indautorep}. Here $\chi=e^{\langle \lambda +\rho|H\rangle}$ is the inducing character on the Borel subgroup $B(\mathbb{A})=N(\mathbb{A})A(\mathbb{A})$, as defined in section \ref{sec_multcharborel}. For the constant term of $E({\sf f}_\lambda, g)$ one can derive an explicit formula; this is done in great detail for $SL(2, \mathbb{A})$ in chapter~\ref{ch:SL2-fourier}. The formula for arbitrary split groups $G(\mathbb{A})$, due to Langlands, will be derived in chapter~\ref{ch:CTF}. Here we are   interested in the representation theoretic properties of the non-constant (abelian) Whittaker coefficients of $E({\sf f}_\lambda, g)$.

\subsection{Global considerations}

For a character $\psi$ on $N(\mathbb{Q})\backslash N(\mathbb{A})$ the abelian coefficients of $E({\sf f}_\lambda, g)$ are given by the  Whittaker coefficient $W_\psi$ of the type (\ref{abelianFC}). Plugging $E({\sf f}_\lambda, g)$ from~\eqref{eq:ESWh} into  (\ref{abelianFC}) and exchanging the order of summation and integration we obtain the formula
\beq
\label{WhES}
W_\psi ({\sf f}_\lambda, g)=\sum_{\gamma \in B(\mathbb{Q})\backslash G(\mathbb{Q})} \int_{N(\mathbb{Q})\backslash N(\mathbb{A})} {\sf f}_\lambda(\gamma n g ) \overline{\psi(n)} dn.
\eeq

Representation theoretically, $W_\psi({\sf f}_\lambda, g)$ belongs to the \emphindex[representation!induced]{induced representation}
\beq
\text{Ind}_{N(\mathbb{A})}^{G(\mathbb{A})}\psi = \left\{W_\psi:G(\mathbb{A})\to \mathbb{C} \, \Big|\, W_\psi(ng)=\psi(n)W_\psi(g), \, \, n \in N(\mathbb{A})\right\}.
\eeq
Equation (\ref{WhES}) thus gives an embedding
\beq
I(\lambda)=\text{Ind}_{B(\mathbb{A})}^{G(\mathbb{A})}\chi \hookrightarrow \text{Ind}_{N(\mathbb{A})}^{G(\mathbb{A})}\psi .
\label{embedding}
\eeq

\begin{definition}[Whittaker model]
\label{def-Whittakermodel}
The space 
\beq
Wh_\psi(\lambda) =\{ W_\psi({\sf f}_\lambda) | {\sf f}_\lambda \in I(\lambda)\}\subset \text{Ind}_{N(\mathbb{A})}^{G(\mathbb{A})}\psi
\eeq
 is called a  \emphindex{Whittaker model} of $I(\lambda)$, and its elements  \emphindex[Whittaker functionals]{Whittaker functionals} or \emphindex{Whittaker vector} (c.f.~Remark~\ref{FourierCoefficientRemark}). The associated map 
\beq
{\sf f}_\lambda \mapsto W_\psi({\sf f}_\lambda),
\eeq
is an \emphindex[intertwiner]{intertwiner} between the principal series $I(\lambda)$ and its Whittaker model $Wh_\psi(\lambda)$.
\end{definition}

\begin{remark}
An important result about Whittaker models is their uniqueness: for each fixed section ${\sf f}_\lambda\in I(\lambda)$ and fixed generic character $\psi$ there exists a \emph{unique} Whittaker function $W_\psi({\sf f}_\lambda)$ (see, e.g., \cite{Bump,Cogdell}). This property is known \emphindex{multiplicity one} and was shown originally for $GL(n)$ locally for archimedean and non-archimedean fields in~\cite{JL,Shalika}. We note that it does not hold for $SL(n)$ if $n>2$~\cite{Blasius}.
\end{remark}

In chapter~\ref{ch:Whittaker-Eisenstein} we will show that, \emp{for generic $\psi$}, the Whittaker coefficient can be written as a single integral rather than a sum. The argument relies on the Bruhat decomposition of $G(\mathbb{Q})$, which allows one to trade the sum over $\gamma \in B(\mathbb{Q})\backslash G(\mathbb{Q})$ for a sum over the Weyl group $\Weyl(\mathfrak{g})$. The end result is that the Whittaker coefficient may be written as
\beq
W_\psi({\sf f}_\lambda, g)=\lint_{N(\mathbb{A})} {\sf f}_\lambda(\wlong n g) \overline{\psi(n)} dn.
\label{WhittakerLongestWeyl}
\eeq
The sum over $\gamma$ has reduced to a single contribution represented by $\wlong$, the longest element in the Weyl group $\Weyl(\mathfrak{g})$ (for the details see chapter~\ref{ch:Whittaker-Eisenstein}).
\begin{remark}
The expression \eqref{WhittakerLongestWeyl} is sometimes known as a \emphindex[Jacquet--Whittaker integral]{Jacquet--Whittaker integral}~\cite{HJacquet}. We will often refer to it simply as the (global) Whittaker function. 
\end{remark}

\subsection{Local considerations}

Recall from section~\ref{sec:AR} that by Flath's tensor product theorem \index{Flath's tensor decomposition theorem}the principal series decomposes into a product over all places~\cite{Flath}
\beq
\text{Ind}_{B(\mathbb{A})}^{G(\mathbb{A})} \chi = \bigotimes_{p\leq \infty} \text{Ind}_{B(\mathbb{Q}_p)}^{G(\mathbb{Q}_p)} \chi_p,
\eeq
and we have a similar decomposition for $\text{Ind}_{N(\mathbb{A})}^{G(\mathbb{A})}\psi$:
\beq
\text{Ind}_{N(\mathbb{A})}^{G(\mathbb{A})}\psi=\bigotimes_{p\leq \infty} \text{Ind}_{N(\mathbb{Q}_p)}^{G(\mathbb{Q}_p)}\psi_p.
\eeq
To each standard section ${\sf f}_{\lambda, p}\in \text{Ind}_{B(\mathbb{Q}_p)}^{G(\mathbb{Q}_p)} \chi_p$ we then have a local ($p$-adic) Whittaker function
\beq
W_{\psi_p}({\sf f}_{\lambda, p}, g)=\int_{N(\mathbb{Z}_p)\backslash N(\mathbb{Q}_p)} {\sf f}_{\lambda, p}(\wlong n g) \overline{\psi_p(n)}dn, \quad g\in G(\mathbb{Q}_p)
\eeq
and the global Whittaker model $Wh_\psi (\chi)$ splits accordingly
\beq
Wh_{\psi}(\chi)=\bigotimes_{p\leq \infty} Wh_{\psi_p}(\chi_p).
\eeq
In chapter~\ref{ch:Whittaker-Eisenstein} we will derive an explicit formula (\emphindex{Casselman--Shalika formula}) for the $p$-adic Whittaker function $W_{\psi_p}({\sf f}_{\lambda, p}, g)$, $p<\infty$, in the special case when  $W_{\psi_p}({\sf f}_{\lambda, p}, g)$ is \emp{spherical}
and $\psi$ unramified, notions that were defined in definitions~\ref{def:Wh} and~\ref{def_unram}, respectively.

For generic characters $\psi$, the global Whittaker function can then be recovered as an Euler product over all places
\beq
W_\psi({\sf f}_{\lambda}, g)=\prod_{p\leq \infty} W_{\psi_p}({\sf f}_{\lambda, p},g_p), \qquad g\in G(\mathbb{A}), \, \, g_p\in G(\mathbb{Q}_p).
\eeq
It is sometimes useful to separate the finite places $p<\infty$ from the infinite place $p=\infty$ and make the following definition:

\begin{definition}[finite Whittaker function]
Consider the Whittaker function $W_\psi^{\text{fin}}$ obtained by taking the product over all the \emph{finite places}:
\beq
W_{\psi}^{\text{fin}}({\sf f}_{\lambda}^{\text{fin}}, g_f)=\prod_{p<\infty} W_{\psi_p}({\sf f}_{\lambda, p},g_p), \qquad g_f=(1; g_2, g_3, \dots )\in G(\mathbb{A}_f).
\label{def_finiteWhittaker}
\eeq 
We call this the \emphindex[Whittaker function!finite]{finite Whittaker function}. 
\end{definition}

\begin{remark}
    The finite Whittaker function plays an important role in string theory where it contributes to the \emphindex[instanton!measure]{instanton measure}, as we illustrate in example \ref{ex_localWhittakerSL2} below and as was discussed in chapter~\ref{ch:intro-strings}.
\end{remark}

\subsection{Spherical  Whittaker functions}
\label{sec_sphericalvector}

Here we will introduce a special class of Whittaker functions which are spherical in an appropriate sense. Assume that $\textrm{Ind}_{B(\mathbb{A})}^{G(\mathbb{A})} \chi$ is \emphindex[character!on Borel subgroup!unramified]{unramified}, i.e.  for almost all places $p$ the local component $\text{Ind}_{B(\mathbb{Q}_p)}^{G(\mathbb{Q}_p)} \chi_p$ is \emphindex[representation!spherical]{spherical}. This implies that there exists a \emph{unique} (up to normalisation) section ${\sf f}^{\circ}_{\lambda, p}\in \text{Ind}_{B(\mathbb{Q}_p)}^{G(\mathbb{Q}_p)} \chi_p$ that satisfies
\beq
{\sf f}^{\circ}_{\lambda, p}(bk) = \chi_p(b), \qquad {\sf f}^{\circ}_{\lambda, p}(k)={\sf f}^{\circ}_{\lambda, p}(1)=1,
\label{defsphericalvector}
\eeq
where $b\in B(\mathbb{Q}_p)$ and $k\in G(\mathbb{Z}_p)$.

\begin{definition}[spherical vector]
We call ${\sf f}^{\circ}_{\lambda, p}\in \text{Ind}_{B(\mathbb{Q}_p)}^{G(\mathbb{Q}_p)} \chi_p$, defined by (\ref{defsphericalvector}), a \emphindex[spherical vector]{spherical vector}.
\end{definition}

\begin{definition}[spherical Whittaker function] 
To each spherical vector ${\sf f}^{\circ}_{\lambda, p}$ and generic character $\psi_p$ we can associate a \emphindex[Whittaker function!spherical]{spherical Whittaker function} $W_{\psi_p}^{\circ}\in Wh_{\psi_p}(\chi_p)$, defined by 
\beq
W_{\psi_p}^{\circ}(\lambda,g)=\int_{N(\mathbb{Z}_p)\backslash N(\mathbb{Q}_p)} {\sf f}^{\circ}_{\lambda, p}(\wlong ng) \overline{\psi_p(n)}dn, \qquad g\in G(\mathbb{Q}_p).
\eeq
\end{definition}

As before, the spherical Whittaker function satisfies the relation
\beq
W^{\circ}_{\psi_p}(\lambda,nak)=\psi_p(n)W^{\circ}_{\psi_p}(\lambda,a),
\label{sphericalWhittakerRelation}
\eeq
where $n\in N(\mathbb{Q}_p), \, a\in A(\mathbb{Q}_p), \, k\in G(\mathbb{Z}_p)$. This again implies that $W^{\circ}_{\psi_p}(\lambda,g)$  is completely determined by its restriction to the Cartan torus $A(\mathbb{Q}_p)$, where it equals
\beq
\label{Wspherical}
W_{\psi_p}^{\circ}(\lambda,a)=\int_{N(\mathbb{Q}_p)}{\sf f}_{\lambda,p}^\circ(\wlong na) \overline{\psi_p(n)}dn.
\eeq

\begin{example}[Spherical Whittaker function for $SL(2,\ads)$]
\label{ex_localWhittakerSL2} We now illustrate the discussion for the Eisenstein series $E(s, g)$ on $SL(2,\mathbb{A})$. The results below are all derived in section \ref{sec:SL2FC}. Recall from example \ref{EisensteinSL2example} that the Eisenstein series  is obtained by choosing the standard section ${\sf f}_\lambda$ to be the spherical vector ${\sf f}^{\circ}_\lambda={\sf f}^{\circ}_s$, such that 
\beq
E({\sf f}^{\circ}_s,g)=\sum_{\gamma \in B(\mathbb{Q})\backslash SL(2,\mathbb{Q})} {\sf f}^{\circ}_s(\gamma g)=\sum_{\gamma \in B(\mathbb{Q})\backslash SL(2,\mathbb{Q})} \chi_s(\gamma na),
\eeq
where $\chi_s = e^{\left< 2s\Lambda |H\right>}, \, \Lambda=\alpha/2$ with $\alpha$ the simple root of  $\mathfrak{sl}(2, \mathbb{R})$. The local spherical Whittaker function (\ref{Wspherical}) is 
\beq
W^{\circ}_{\psi_p}(s, a)=\int_{N(\mathbb{A})} \chi_s(\wlong n a)\overline{\psi(n)} dn, \qquad a\in A(\mathbb{Q}_p).
\eeq

As will be shown in detail in section~\ref{sec:SL2FC}, the integral equals 
\beq
W^{\circ}_{\psi_\infty}(s, y)= \frac{2\pi^s}{\Gamma(s)}  y^{1/2}  |m|^{s-1/2} K_{s-1/2}(2\pi |m| y), 
\label{realWhittakerSL2}
\eeq
at the archimedean place $p=\infty$ (see~\eqref{SL2Whittinf}).
Here, $m\in \mathbb{Z}^\times$, and the Cartan torus $A(\reals)$ is given by
\beq
\begin{pmatrix} y^{1/2} & \\ & y^{-1/2} \end{pmatrix}, \qquad y\in \mathbb{R}_{>0}.
\eeq

At the non-archimedean places, the integral becomes (cf.~\eqref{SL2Wloc})
\beq
W^{\circ}_{\psi_p}(s, v)=|v|_p^{-2s+2}\gamma_p(mv^2)(1-p^{-2s})\frac{1-p^{-2s+1}|mv^2|_p^{2s-1}}{1-p^{-2s+1}},
\label{ex_localWhittakersl2}
\eeq
with $p<\infty$ and $m \in \mathbb{Q}^\times$, and we parametrised the torus $A(\mathbb{Q}_p)$ by 
\beq
\begin{pmatrix} v & \\ & v^{-1} \end{pmatrix}, \qquad v\in \mathbb{Q}_p^\times
\eeq
for all $p<\infty$.

The associated finite Whittaker function (\ref{def_finiteWhittaker}), evaluated at the identity $v=1$, is only non-vanishing for $m\in \mathbb{Z}^\times$ because of the $\gamma_p$ factors as is seen in section \ref{sec:p-adic-special-functions}.  For $m \in \mathbb{Z}^\times$
\begin{equation}
    \label{ex_finiteWhittakersl2}
    W^{\circ, \text{fin}}_{\psi}(s, 1) = \prod_{p<\infty} W^{\circ}_{\psi_p}(s, 1) = 
    \left( \prod_{p<\infty} (1 - p^{-2s}) \right) \left( \prod_{p<\infty} \frac{1 - p^{-(2s-1)} \abs{m}_p^{2s-1}}{1 - p^{-(2s-1)}} \right) \, .
\end{equation}

The first factor is simply the Euler product \eqref{RiemannProd} of the (inverse of the) Riemann zeta function $\zeta(2s)^{-1}$. We will now show that the second factor is actually the divisor sum $\sigma_{t}(m)$ defined in \eqref{intro:divisorsum} denoting $t = 1-2s$ for brevity.

Assume first that $m = p^{a}$ for some prime $p$ and positive integer $a$. Then
\begin{equation}
    \sigma_{t}(m) = \sum_{d | m} d^{t} = 1 + p^{t} + p^{2t} + \ldots + p^{at} = \frac{1 - p^{(a+1)t}}{1 - p^t} \, .
\end{equation}

For $m = p^a q^b$ we get
\begin{equation}
    \begin{split}
        \sigma_{t}(m) &= 1 + p^{t} + q^{t} + p^{2t} + q^{2t} + p^{t} q^{t} + \ldots + p^{at} q^{bt}\\
        &= (1 + p^{t} + \ldots + p^{t})(1 + q^{t} + \ldots + q^{bt}) = \sigma_{t}(p^a) \sigma_{t}(q^b) \, .
    \end{split}
\end{equation}

Similarly, for the general case with $m$ having the prime factorisation $m = p_1^{a_1} \cdots p_r^{a_r}$,
\begin{equation}
    \sigma_t(m) = \sigma_t(p_1^{a_1}) \cdots \sigma_t(p_r^{a_r}) = \prod_{i=1}^r \frac{1 - p_i^t p_i^{a_i t} }{1 - p_i^t} = \prod_{p<\infty} \frac{1 - p^t \abs{m}_p^{-t} }{1 - p^t} 
\end{equation}
since $\abs{m}_p = p_j^{-a_j}$ for $p = p_j$ (some $j$) and otherwise $\abs{m}_p = 1$. In other words, the finite spherical Whittaker function for $SL(2,\ads$) (and the divisor sum $\sigma_t$) are \emphindex[series!multiplicative]{multiplicative}.
Thus, for non-zero integer $m$
\begin{equation}
    W^{\circ, \text{fin}}_{\psi}(s, 1) = \frac{1}{\zeta(2s)} \sigma_{1-2s}(m) \, .
\end{equation}

Comparing with the discussion in section \ref{intro:Fourier} we conclude that the finite Whittaker function $W^{\circ, \text{fin}}_\psi$, defined in (\ref{def_finiteWhittaker}), is closely related to the \emphindex[instanton!measure]{instanton measure} in string theory. More precisely, when evaluating the finite Whittaker function at the identify in $SL(2,\mathbb{A}_f)$ we  obtain the divisor sum which is characteristic for so-called D$(-1)$-instanton effects in string theory (see \cite{Green:1997tv}). This in fact also holds for more general groups $G(\mathbb{A})$ and gives a strong physics motivation for the detailed analysis of the Casselman--Shalika formula presented in section \ref{sec_CSformula}.
\end{example}

\section{Fourier coefficients and nilpotent orbits*}
\label{sec:orbits}

When considering the Fourier expansion along a unipotent radical $U(\mathbb{A})$ that is part of a standard parabolic subgroup $P(\mathbb{A})=L(\mathbb{A}) U(\mathbb{A})$, one can group the \index{Fourier coefficient!orbit under Levi}Fourier integrals (\ref{eq:FCU}) into orbits of the Levi factor $L(\mathbb{Q})$, see for example~\cite{MillerSahi,Green:2011vz}. There is a close connection to the theory of nilpotent orbits of the adjoint action of $G(\cx)$ on its Lie algebra $\mf{g}(\cx)$ and the notion of wavefront sets through the work of M\oe{}glin--Waldspurger~\cite{MoeglinWaldspurger,MoeglinWaveFront}, Matumoto~\cite{Matumoto}, Ginzburg--Rallis--Soudry~\cite{GinzburgRallisSoudryFourier,GinzburgConjectures,GinzburgTowards}, Gomez--Gourevitch--Sahi~\cite{GomezGourevitchSahi,GourevitchSahi1,GourevitchSahi2}, Jiang--Liu--Savin~\cite{JiangLiuSavin} and many others.

\begin{remark}
The discussion of the present section only applies to Fourier expansions along unipotent radicals $U$ of non-minimal parabolic subgroups; for expansions along $N(\mathbb{A})$ contained in the (minimal parabolic) Borel subgroup $B(\mathbb{A})$ the orbits under the abelian Levi factor become single points.
\end{remark}

\subsection{Character variety orbits}
\label{sec:character-variety-orbits}

Let $\psi$ denote a unitary character on $U(\mathbb{A})$ that is trivial on $U(\mathbb{Q})$ and consider the Fourier integral $F_{\psi}(g)\equiv F_\psi(\varphi,g)$ of an automorphic form $\varphi$ as defined in definition~\ref{def:F}. We consider $\varphi$ fixed for the following discussion and will suppress it in the notation $F_\psi(g)$. Under the action of an element $\gamma\in L(\mathbb{Q})$, that is an element $\gamma$ of the intersection of the discrete subgroup with the Levi factor, the Fourier coefficient changes as follows
\begin{align}
\label{Whittakerorbit}
F_{\psi}(\gamma g) &= \lint_{U(\mathbb{Q})\backslash U(\mathbb{A})} \varphi(u \gamma g) \overline{\psi(u)} du
=  \lint_{U(\mathbb{Q})\backslash U(\mathbb{A})} \varphi(\gamma^{-1} u \gamma g) \overline{\psi(u)} du\nn\\
&=  \lint_{U(\mathbb{Q})\backslash U(\mathbb{A})} \varphi(u g) \overline{\psi(\gamma u\gamma^{-1})} du
= F_{ \psi^\gamma}( g) 
\end{align}
where we have used the fact that $\varphi$ is invariant under discrete transformations as well as the fact that the change of coordinates $u\to \gamma^{-1}u \gamma$ is uni-modular since $\gamma$ is in the discrete subgroup. In the last step, we have defined the transformed character 
\begin{align}
\psi^{\gamma}(u) := \psi(\gamma u \gamma^{-1})
\end{align}
and identified its Fourier coefficient. The transformed character $\psi^\gamma$ is well-defined since the Levi component $L(\mathbb{Q})$ acts on $U(\mathbb{Q})$ by conjugation. In view of the terminology introduced in definition~\ref{def:UC}, the orbits thus produced are called \emphindex[character!variety orbit]{character variety orbits}. 
We also introduce the following notion:
\begin{definition}
Let $\psi$ be a unitary character on the unipotent subgroup $U(\ads)$ of a standard parabolic subgroups $P(\ads)= L(\ads) U(\ads)$. The set
\begin{align}
C_\psi = \left\{ \gamma \in L(\rats) \st \psi^\gamma=\psi \right\}
\end{align}
is called the \emphindex[character!stabiliser]{stabiliser} of the character $\psi$. We will sometimes use the same terminology when referring to the action of $L(\reals)$ or $L(\cx)$ on the corresponding character variety.
\end{definition}

The calculation~\eqref{Whittakerorbit} shows that the Fourier coefficient $F_\psi$ is invariant (automorphic) under the stabiliser subgroup $C_\psi$. 

The adjoint action of $L(\mathbb{Q})$ on $U(\mathbb{Q})$ can be described more explicitly by realizing the original character $\psi$ in terms of a weight vector similar to proposition~\ref{prop:UC}. The Lie algebra $\mf{u}$ consists of nilpotent elements $X\in \mf{u}$ and we can write an element $u\in U$ as $u= e^X$. A unitary character $\psi$ on $U$ is then given by an element $\omega$ of the dual space $\mf{u}^*$ via
\begin{align}
\label{Uchar}
\psi(e^X) = \exp\left(2\pi i \omega(X) \right)
\end{align}
and the triviality~\eqref{psitriv} of $\psi$ on the commutator subgroup $[U,U]$ enforces that 
\begin{align}
\omega\left([\mf{u},\mf{u}]\right) = 0 ,
\end{align}
so that $\omega$ is not an arbitrary element of $\mf{u}^*$ but one associated with the Lie algebra of the abelianisation $[U,U]\bs U$.  Clearly, the abelianisation $[U,U]\bs U$ is preserved by the adjoint action of $L(\rats)$ on $U(\ads)$ and $L(\rats)$ therefore acts dually on the space of allowed $\omega$. 
By virtue of (\ref{Whittakerorbit}), the Fourier coefficients for all characters in one orbit are related and it suffices to calculate the Fourier coefficient of one representative of an orbit. In practice, it is more convenient to take the dual of $\omega$ and study the adjoint nilpotent orbits of the action of $L(\mathbb{Q})$ on $\mathfrak{u}(\mathbb{Q})$, where one can also restrict to the abelian quotient $[\mf{u},\mf{u}]\bs \mf{u}$.

\begin{remark}
Let $\Sigma$ be the subset of the simple roots $\Pi$ that defines a standard parabolic subgroup $P=LU$, cf. section~\ref{sec:parsubgp}. The nilpotent Lie algebra $\mf{u}$ of $U$ has a (finite) graded decomposition 
\begin{align}
\mf{u} = \bigoplus_{j\in\ints} \mf{u}_j,\quad
\textrm{with}\quad
\mf{u}_j = \left\langle
E_\alpha \st \alpha = \sum_{\beta \in \Pi} n_\beta \beta\in\Delta_+ \quad \textrm{and} \quad
\sum_{\alpha\in \Pi\setminus \Sigma} n_\alpha =j
\right\rangle.
\end{align}
Each space $\mf{u}_j$ is preserved by the adjoint action of $L$ and the space of characters $\psi$ on $U$ is dual to $\mf{u}_1$. The character variety orbits can therefore be viewed dually in $\mf{u}_1$.  
The space $\mf{u}_1$ is isomorphic (as a vector space) to $[\mf{u},\mf{u}]\bs \mf{u}$.
\end{remark}

\begin{example}[Mirabolic subgroups of $GL(n,\reals)$]
\label{ex:mirabolic}
Consider $G(\reals)=SL(n,\reals)$ that can be represented by $(n\times n)$-matrices. A maximal parabolic subgroup can be chosen with Levi factor $L(\reals)=GL(n-1,\reals)$ through the following matrices
\begin{align}
\label{miraL}
L(\reals)=\left\{\begin{pmatrix}
* &0&0&\cdots &0&0\\
0 &* &*&\cdots &*&*\\
0 &* &*&\cdots &*&*\\
\vdots & \vdots &&&&\vdots\\
0 &* &*&\cdots &*&*
\end{pmatrix}\right\}
=\left\{\begin{pmatrix}
r &0\\
0 &  m 
\end{pmatrix}\st r\in \reals ,\, m\in GL(n-1,\reals) \textrm{ such that } \det(m) = r^{-1}\right\}
\end{align}
and associated $(n-1)$-dimensional unipotent radical
\begin{align}
U(\reals)=\left\{\begin{pmatrix}
1 &*&*&\cdots &*&*\\
0 &1 &0&\cdots &0&0\\
0 &0 &1&\cdots &0&0\\
\vdots & \vdots &&&&\vdots\\
0 &0 &0&\cdots &0&1
\end{pmatrix}\right\}
=\left\{\begin{pmatrix}
1 & u^T\\
0 & \idm_{n-1} 
\end{pmatrix}\st u \in \reals^{n-1} \right\}.
\end{align}
The unipotent radical is abelian in the present case and acted upon by $L=GL(n-1,\reals)$. Characters $\psi$ can be thought of as being given by $(n-1)$-column vectors $\omega$ that contract into $X\in \mathrm{Lie}(U)=\mf{u}$ and define the character via (\ref{Uchar}). These parabolic subgroups are sometimes referred to as \emphindex[subgroup!mirabolic]{mirabolic subgroups}.

For the local transformation of the Fourier coefficients~\eqref{Whittakerorbit} at the archimedean place one needs to restrict to orbits under $L(\ints)$ that force $r=1$ and $m\in SL(n-1,\ints)$ in~\eqref{miraL} (or $r=-1$ and $\det(m)=-1$ but this does not influence the discussion below). The action of the Levi subgroup $L(\ints)$ on $U(\reals)$ is by
\begin{align}
\begin{pmatrix}
1 & \\
& m
\end{pmatrix}
\begin{pmatrix}
1& u^T \\
&1\end{pmatrix}
\begin{pmatrix}
1& \\&m^{-1}
\end{pmatrix}
=\begin{pmatrix}
1 & u^T m^{-1}\\
&1 
\end{pmatrix}.
\end{align}
The group $L(\ints)$ then acts on the character variety $\mf{u}^*$ modelled by a vector $\omega\in\reals^n$ by
\begin{align}
\omega \mapsto m^{-1} \omega,
\end{align}
that is, simply by left multiplication of the column vector. The character variety $\{ \omega\in \reals^{n-1}\}$ decomposes into infinitely many orbits with representatives
\begin{align}
\begin{pmatrix} \sigma \\0\\\vdots\\0\end{pmatrix} \quad \textrm{for} \quad \sigma\in \reals_{\geq 0}
\end{align}
under this action. If the character is trivial on integral points (as is the case for unitary characters trivial on $U(\ints)$) the representatives are labelled by $\sigma\in \ints_{\geq 0}$:
\begin{align}
\ints^{n-1} = \bigcup_{\sigma\in \ints_{\geq 0}} \left\{SL(n-1,\ints)\cdot \begin{pmatrix} \sigma \\0\\\vdots\\0\end{pmatrix}\right\}.
\end{align}
An arbitrary vector $\omega\in \ints^{n-1}$ belongs to the orbit with $\sigma=\gcd(\omega)$.
\end{example}

\begin{remark}
Classifying the orbits of the action of $L$ on $U$ over $\mathbb{Z}$ or $\mathbb{Q}$ is in general a difficult task, see~\cite{JiangRallis,BhargavaI,SavinWoodbury,Krutelevich} for some results.
A slightly coarser description can be obtained by complexifying the Levi subgroup to $L(\cx)$ and studying the complex orbits. All such complex orbits have been determined in the literature~\cite{Littelmann,deGraaf,MillerSahi}, using the methods of Dynkin~\cite{Dynkin}, Kostant~\cite{Kostant,KostantRallis}, Bala--Carter~\cite{BalaCarterI,BalaCarterII}, Vinberg~\cite{Vinberg1975} and Kac~\cite{KacVinberg}. 
\end{remark}

\subsection[Wavefront sets and vanishing theorems for Fourier coefficients]{Wavefront sets and vanishing theorems for\\ Fourier coefficients}
\label{sec:wavefront-sets}

There are many different choices of parabolic subgroup $P(\ads)=L(\ads)U(\ads)$ and associated Fourier expansions along their unipotents $U$. All the different character variety orbits of the action of $L(\rats)$ on unitary characters on $U(\ads)$ are associated with nilpotent elements $\omega\in \mf{u}^*\subset \mf{g}^*$. A given character variety orbit therefore lies in a \emphindex{coadjoint nilpotent orbit} of the action of $G(\reals)$ on elements of $\mf{g}^*$ that are dual to nilpotent elements. Properties of automorphic representations of $G(\ads)$ are only associated with structures arising from $G$, implying that the character variety orbits (under the action of $L$) are less fundamental than the \emphindex[nilpotent orbit]{nilpotent orbit} they embed into. We will not fully develop the theory of nilpotent orbits here but refer the reader to the books~\cite{CollingwoodMcGovern,CarterBook,Spaltenstein} for a detailed exposition. Below we will mention only some aspects that are of relevance to our discussion.

The approach using nilpotent orbits is useful because it sometimes allows to determine that certain Fourier coefficients must vanish identically without actually calculating them. At the heart of this is the notion of the (complexified) \emph{wavefront set} of an irreducible automorphic representation $\pi$ (cf.~definition~\ref{autorepdef}).

\begin{definition}[Wavefront set]
Let $\pi_p$ be an automorphic representation of $G(\rats_p)$ at a local place $p$. The \emphindex[wavefront set|textbf]{wavefront set} of $\pi_p$ is given by
\begin{align}
\WF(\pi_p) = \bigcup_{i\in I} \overline{\cO}_i,
\end{align}
where the $\cO_i$ are a finite collection of complex nilpotent orbits and the closure is with respect to the Zariski topology naturally defined on the set of nilpotent elements of $\mf{g}(\cx)$. The $\cO_i$ appearing in the sum are characterized by admitting a non-trivial Fourier coefficient~\cite{MoeglinWaveFront,GinzburgTowards}.
\end{definition}

The wavefront set therefore is the closure of a (set of) nilpotent orbits~\cite{Joseph,BorhoBrylinski}. Originally, it is defined as the \emphindex[annihilator ideal]{annihilator ideal} associated with the representation $\pi_p$; in the case of the so-called minimal representation it is also referred to as the \emphindex{Joseph ideal}.

\begin{remark}
We will also use the notion of a \emph{global wavefront set} of an adelic representation $\pi=\otimes_{p\leq \infty} \pi_p$ of $G(\ads)$. It is a priori not guaranteed that the local wavefront set $\WF(\pi_p)$ does not vary as $p$ varies and therefore one has to treat this notion with care. For Eisenstein series induced by characters of the form~\eqref{characterBorel} this does not happen. Global wavefront sets have been discussed for example in~\cite{JiangLiuSavin} where it was also shown that the maximal orbits in wavefront sets have to be so-called \emph{special} orbits.\index{nilpotent orbit!special} This property was known for local wavefront sets due to~\cite{MoeglinWaveFront}.
\end{remark}

A  nilpotent orbit for a Lie algebra $\mf{g}$ is the orbit of a nilpotent element $X\in\mf{g}$ under the action of the adjoint group $G$ with Lie algebra $\mf{g}$, see for example~\cite{CollingwoodMcGovern,Spaltenstein} for an introduction. The theorems of~\cite{Joseph,BorhoBrylinski} show that one can associate (the closure of) a unique nilpotent orbit in $\mf{g}$ to any irreducible automorphic representation $\pi$ of $G$, meaning that the wavefront set of irreducible automorphic representation is given by the closure of unique maximal orbit (w.r.t. the partial closure ordering).
One can also consider the action of the adjoint group $G$ on the dual Lie algebra $\mf{g}^*$ and study coadjoint nilpotent orbits. Using the non-degenerate Killing form, we can identify adjoint and coadjoint nilpotent orbits. By the correspondence~\eqref{Uchar} one can view characters $\psi$ on some unipotent $U$ as elements of $\mf{g}^*$ and the character variety orbits lie therefore in coadjoint nilpotent orbits.

The link to the $L(\mathbb{C})$-orbits of Fourier coefficients $F_{\psi_U}$ of an automorphic function $\varphi$ is provided\index{wavefront set!Moeglin--Waldspurger theorem}\index{wavefront set!Matumoto theorem} by the theorems of M\oe glin--Waldspurger and Matumoto~\cite{MoeglinWaldspurger,Matumoto,MillerSahi} that assert that a Fourier coefficient can only be non-zero if its associated character variety orbit in $\mf{u}^*$ (under the action of $L(\mathbb{C})$) intersects a coadjoint nilpotent orbit in $\mf{g}^*\supset\mf{u}^*$ (under the action of $G(\mathbb{C})$) that belongs to the wavefront set associated with the automorphic representation to which $\varphi$ belongs.

\begin{example} [Minimal representation of $E_6$]
Suppose $\varphi$ belongs to the minimal representation of the exceptional Lie group $E_6(\reals)$ of dimension $78$. Then its associated wavefront set is the closure of the minimal orbit (of dimension $22$).  The minimal representation of $E_6$ can be realized as a special point in the degenerate principal series representation that is associated with a maximal parabolic subgroup $P=LU$ with Levi factor $L=SO(5,5)\times GL(1)$. The unipotent $U$ in this case is a Heisenberg group of dimensions $21$. The (dualized) character variety $\mf{u}_1=[\mf{u},\mf{u}]\bs \mf{u}$ has dimension $20$ and is acted upon by $L(\reals)$. After complexification one finds that $\mf{u}_1$ breaks up into five different character variety orbits under $GL(6,\mathbb{C})$~\cite{MillerSahi}. Of these only the trivial and the smallest non-trivial one intersect the closure of the minimal coadjoint nilpotent orbit. One concludes that the Fourier coefficients in the remaining three character variety orbits must vanish in the minimal representation. 

We also note that the Gelfand--Kirillov dimension of the degenerate principal series in this case is $21=20+1$, corresponding to the dimension of the Heisenberg group. At the special point where the minimal representation is realized the $20$-dimensional space can be polarized into $10$ `coordinates' and $10$ `momenta' and the Heisenberg algebra is realized on functions of the $10$ coordinate variables on which the momenta act as derivative operators. This action of the Heisenberg group extends to all of $E_6(\reals)$ and can also be given an oscillator realization~\cite{Gunaydin:2001bt}. 

This example is based on \cite{GRS,Pioline:2010kb,MillerSahi,Green:2011vz} where more information can be found. The minimal representation discussed here is an example of a small representation that we will discuss in more detail in sections~\ref{sec:Slambda} and~\ref{smallreps}.
\end{example}

The connection between nilpotent orbits and Fourier coefficients is made more concrete in the work of Ginzburg~\cite{GinzburgConjectures}. We follow~\cite{Gustafsson:2014iva} in the following discussion. To the nilpotent orbit $\cO\equiv \cO_X$ of a nilpotent element $X\in\mf{g}$ one can associate a \emphindex{Jacobson--Morozov triple} $H,X,Y\in \mf{g}$ that satisfies the standard $\mf{sl}(2)$ Lie algebra relations. The orbit is uniquely characterised by the (unique) Weyl chamber image of $H$ under the action of the Weyl group. This leads to a labelling of nilpotent orbits by weighted Dynkin diagrams, where the weights are non-negative integers. (These integers are less than or equal to two but this does not matter for our discussion.) Any integrally weighted Dynkin diagram gives rise to a graded decomposition 
\begin{align}
\mf{g} = \bigoplus_{i\in \ints} \mf{g}_i,
\end{align}
where $\mf{g}_i$ is the space of elements in $\mf{g}$ with eigenvalue $i$ under the adjoint action of the $H$ that lies in the Weyl chamber. All $\mf{g}_i$ are of finite dimension and there are only finitely many non-trivial $\mf{g}_i$ since $\mf{g}$ is finite-dimensional. We define
\begin{align}
\mf{l}_\cO = \mf{g}_0,\quad \mf{u}_\cO = \bigoplus_{i\geq 1} \mf{g}_i \quad\textrm{and}\quad 
\mf{v}_\cO = \bigoplus_{i\geq 2} \mf{g}_i.
\end{align}
Let $L_\cO$, $U_\cO$ and $V_\cO$ be corresponding subgroups of $G$. A nilpotent orbit $\cO$ has a unique stabiliser $C_\cO\subset L_\cO$ that is a reductive group. 

\begin{definition}[orbit Fourier coefficient]
\label{def:orbit-coefficient}
Let $\cO$ be a non-trivial nilpotent orbit and let $\psi_V:  V_\cO(\rats) \bs V_\cO(\ads)\to U(1)$ be a unitary character on $V_\cO$. We require $\psi_V$ to have the same stabiliser type under the action of $L_\cO$ as the stabiliser $C_\cO$ of the orbit $\cO$. Then the \emphindex[Fourier coefficient!orbit]{orbit Fourier coefficient} of an automorphic form $\varphi$ belonging to some automorphic representation $\pi$ is defined as
\begin{align}
F_\cO(\varphi, \psi_V,g) = \lint_{V_\cO(\rats)\bs V_\cO(\ads)} \varphi(vg) \overline{\psi_V(v)}dv.
\end{align}
For the trivial orbit $\cO=\{0\}$ we define the orbit Fourier coefficient to be the constant term along the maximal unipotent $N(\ads)$ as in definition~\ref{CTdef}.
\end{definition}

The orbit Fourier coefficients vanish when the orbit does not belong to the wavefront set and allow a rewriting of the Fourier expansion of an automorphic function in terms of a sum over nilpotent orbits. This is similar to the expansion of the \emphindex{Howe--Harish-Chandra expansion} of the \emphindex[character!distribution]{character distribution} of an automorphic representation~\cite{Howe,HC}:
\begin{align}
\mu(\pi) = \sum_{\cO \in \WF(\pi)} c_\cO \mu_\cO
\end{align}
For local automorphic representations $\pi_p$, the numbers $c_\cO$ are computed by the M\oe{}glin--Waldspurger theorem~\cite{Rodier1,MoeglinWaldspurger}.

\begin{remark} 
It is often possible to relate the orbit Fourier coefficients to (degenerate) Whittaker coefficients and this was done for example in~\cite{GinzburgConjectures,MR2123125,Gustafsson:2014iva}. Turning the argument around, one might suspect that the wavefront set can be effectively computed by studying the degenerate Whittaker coefficients with charges defining the parabolic subgroups defining a nilpotent orbit in the Bala--Carter classification. This is borne out for minimal representations~\cite{MR2123125,MillerSahi} and also well supported for some other small representations relevant for string theory~\cite{Bossard:2014lra,Gustafsson:2014iva,Bossard:2015oxa}, see also the discussion in sections~\ref{smallreps} and~\ref{sec:WFred}.
\end{remark}

\begin{remark}
The connection between the \emphindex{associated variety} of automorphic representations and degenerate Whittaker coefficients was also studied in detail recently in work of Gourevitch and Sahi together with Gomez~\cite{GourevitchSahi1,GourevitchSahi2,GomezGourevitchSahi}. Their main emphasis was on local results for the real and $p$-adic place but some results in the global case can be found in~\cite{GomezGourevitchSahi}. An important refinement contained in their work is the introduction of generalised Whittaker pairs that go beyond Jacobson--Morozov triples for nilpotent orbits. In this more general language it is often possible to more easily determine the vanishing (or not) of certain generalised degenerate Whittaker coefficients and from this obtain information about the wavefront set of a representation.

\end{remark}

\section{Method of Piatetski-Shapiro and Shalika*}
\index{Piatetski-Shapiro!method of}
\index{Shalika!method of}
\label{sec:Piatetski-Shapiro-Shalika}

The grouping of Fourier coefficients into orbits under a Levi subgroup $L$ discussed in the previous section is a powerful tool for analyzing automorphic forms. This is at the heart of the method of Piatetski-Shapiro and Shalika that expresses an automorphic form on $GL(n,\reals)$ completely in terms of its Whittaker coefficients (with respect to $N$)~\cite{PiatetskiShapiro,Shalika}. We briefly explain how this connection between Fourier coefficients along $U$ and Whittaker coefficients along $N$ comes about in the case of $GL(n,\reals)$ following~\cite{GanLectures,MillerSahi}. Generalisations to some other groups have been discussed by Miller and Sahi~\cite{MillerSahi}.

Let $P=LU$ be a parabolic subgroup of $G$. According to proposition~\ref{prop:NAFC} we have for a spherical automorphic form $\varphi(g)=\varphi(gk)$ in an automorphic representation $\pi$ that
\begin{align}
\label{PSFourier}
\sum_{\psi} F_{\psi}^\circ(g) = \lint_{U^{(2)}(\mathbb{Q})\bs U^{(2)}(\mathbb{A})} \varphi(ug) du,
\end{align}
where $U^{(2)} = [U,U]$ is the derived group of $U$ and 
\begin{align}
F_{\psi} (g)= \lint_{U(\mathbb{Q})\bs U(\mathbb{A})} \varphi(ug) \overline{\psi(u)} du
\end{align}
is the Fourier coefficient of $\varphi$ along $U$ for the character $\psi$.

Now we want to group the sum over the characters $\psi$ into complex orbits thanks to (\ref{Whittakerorbit}). A given character $\psi$ can have a stabiliser $C_{\psi}(\mathbb{A})\subset L(\mathbb{A})$ under the action of $L(\mathbb{A})$ and $ C_{\psi}(\mathbb{Q})=C(\mathbb{Q})\cap L_{\psi_U}(\mathbb{A})$ is a discrete subgroup of it.  Writing the set of complex character orbits as $\WF(\pi)$ one can write (\ref{PSFourier}) as~\cite{MillerSahi} 
\begin{align}
\label{PSformula}
\sum_{\psi} F_{\psi}^\circ(g) = \sum_{\cO \in\WF(\pi)} \sum_{\psi\in \cO}\sum_{\gamma\in  C_{\psi}(\mathbb{Q})\backslash L(\mathbb{Q})} F_{\psi} (\gamma g),
\end{align}
where the sum over $\psi\in \cO$ denotes \emph{single representatives} of the different integral orbits contained in the complex orbit $\cO\in \WF(\pi)$. The method of Piatetski-Shapiro and Shalika then uses the fact that $F_{\psi}(\gamma g)$ is a function on the reductive stabiliser $C_{\psi}(\mathbb{A})$ and automorphic under $C_{\psi}(\mathbb{Q})$ and so can be expanded in the same manner, yielding an iterative procedure for determining the Fourier expansion.

In the case of $GL(n,\reals)$ this can be done very successfully in terms of iterations of parabolic subgroups of the type discussed in example~\ref{ex:mirabolic}. As was explained there, the unipotent subgroup is abelian and therefore the Fourier expansion (\ref{PSFourier}) recovers the whole automorphic function $\varphi$. Moreover, there is a unique non-trivial complex orbit of $GL(n-1,\cx)$ acting on the $(n-1)$-dimensional $U(\cx)$. The trivial orbit corresponds to trivial $\psi=1$ and corresponds to the constant term in the expansion along $U$, cf. definition~\ref{CTdef}. In order to avoid having to include this term at every iteration step, we now assume until the end of this section that $\varphi$ is a cusp form. Then the sum over complex orbits in $\WF(\pi)$ has only a single element.

The integral orbits contained in the single complex orbit can also be identified easily in this case. They are represented by (non-zero) integers $m_{\alpha_1}$ and the representatives can be chosen to be such that
\begin{align}
\psi(u) = e^{2\pi i m_{\alpha_1} u_{\alpha_1}}
\end{align}
where $\alpha_1$ is the first simple root. (Allowing $m_{\alpha_1}$ to be integral instead of integral and positive actually overcounts the integral orbits by a factor of two but this has no impact on the final result.) In terms of matrices as in example~\ref{ex:mirabolic} this can be written as 
\begin{align}
u = \exp\begin{pmatrix}
0 &u_{\alpha_1}& u_{\alpha_1+\alpha_2}&\cdots &u_{\alpha_1+\ldots+ \alpha_{n-1}}\\
0 &0 &0&\cdots &0\\
0 &0 &0&\cdots &0\\
\vdots & \vdots &&&\vdots\\
0 &0 &0&\cdots &0
\end{pmatrix}\,,\quad
\omega = (m_{\alpha_1},0,0,\ldots,0)^T \in \mf{u}^*
\end{align}
The stabiliser of such a character in $L(\mathbb{Q})=GL(n-1,\mathbb{Q})$ is given by $GL(n-2,\mathbb{Q})\subset GL(n-1,\mathbb{Q})$.  At the present stage we have therefore from (\ref{PSformula})
\begin{align}
\label{PSstep1}
\varphi (g) = \sum_{m_{\alpha_1}\in \ints} \sum_{\gamma\in GL(n-2,\mathbb{Q})\bs GL(n-1,\mathbb{Q})} F_{\psi} (\gamma g).
\end{align}

The Fourier coefficient $F_{\psi}(\gamma g)$ appearing in (\ref{PSstep1}) is therefore an automorphic form on $GL(n-1,\mathbb{Q})$ automorphic under $GL(n-2,\mathbb{Q})$. The iteration now consists in repeating the same process for this smaller subgroup. What this will produce is a sum over $m_{\alpha_2}\in\ints$ and an automorphic form on $GL(n-2,\mathbb{A})$ and so on. At the end of the iteration we obtain
\begin{align}
\label{PSfull}
\varphi (g) = \sum_{m_{\alpha_1}, \ldots m_{\alpha_{n-1}} \in \mathbb{Z}} \sum_{\gamma\in N(n-1,\mathbb{Q})\backslash GL(n-1,\mathbb{Q})} W_{\psi}(\gamma g)
\end{align}
where 
\begin{align}
W_{\psi}(g) = \lint_{N(\mathbb{Q})\bs N(\mathbb{A})} \varphi(ng) \overline{\psi(n)}dn
\end{align}
is a standard Whittaker coefficient on $N$ for the character with instanton charges $m_\alpha$ for the simple roots $\alpha\in\Pi$ and as defined in~\eqref{abelianFC}. (The integration domain is enlarged from $U$ to $N$ by combining some of the intermediate sums over cosets.) Reassembling the sum over all these characters we therefore can also also write
\begin{align}
\varphi(g) = \sum_{\gamma\in N(n-1,\mathbb{Z})\bs GL(n-1,\mathbb{Z})} \lint_{N^{(2)}(\mathbb{Z})\bs N^{(2)}(\mathbb{R})} \varphi(n\gamma g) dn
\end{align}
where $N^{(2)}$ is the derived group of $N$, cf.~(\ref{eq:DSeries}), and we have projected back down to $\reals$. The power of the formula (\ref{PSfull}) is that it allows us to reconstruct the whole (cuspidal) automorphic form from its standard Whittaker coefficients by taking suitable translates of them. This is important since, according to (\ref{generalFourierExpansion}), an automorphic function also contains terms beyond the standard Whittaker coefficients in its expansion. The result of Piatetski-Shapiro and Shalika tells us how to compute these non-abelian terms as translates of abelian terms. We will see a similar structure later when we study the case of $SL(3)$ in detail in section~\ref{sec:SL3ex}.

\chapter[Fourier coefficients of Eisenstein series on \texorpdfstring{$SL(2,\mathbb{A})$}{SL(2, A)}]{Fourier coefficients of\\ Eisenstein series on \texorpdfstring{$SL(2,\mathbb{A})$}{SL(2, A)}}
\label{ch:SL2-fourier}

In this chapter we apply the formalism developed in chapter~\ref{ch:fourier} to the classical theory of non-holomorphic Eisenstein series $E(s,z)$ on the double coset $SL(2,\mathbb{Z})\backslash \mathbb{H}$, with $z\in \mathbb{H}=SL(2,\mathbb{R})/SO(2,\mathbb{R})$ which was already presented as a canonical example in the introduction. Following the analysis in the previous section, we will consider the adelic treatment of this Eisenstein series. The purpose of this chapter is to give an example for the Fourier expansion of an Eisenstein series, where the method can be made explicit and is still fully tractable. We try to carefully introduce every step in the calculation, however the explanation of some of the underlying theory is postponed to the next chapter, where Langlands constant term formula is derived in full detail. Where appropriate in this chapter, we refer the reader to the next chapter for more detailed explanations. In appendix~\ref{app:PSKS}, we present another derivation of the Fourier expansion based on Kloosterman sums.

\section{Statement of theorem}
Before we state the theorem, let us introduce some of the necessary terminology. Recall from chapter~\ref{ch:Lie-groups} that the adelic group $SL(2,\mathbb{A})$ has the maximal compact subgroup $K_{\mathbb{A}}$, defined by
\begin{align}
\label{KA}
K_{\mathbb{A}}  = K(\mathbb{R})\times \prod_{p<\infty} K(\mathbb{Q}_p)=SO(2,\mathbb{R})\times \prod_{p<\infty} SL(2,\mathbb{Z}_p)\,,
\end{align}
see for example~\cite{LanglandsEP}. We then have, by strong approximation~\eqref{eqn:strong-approx}, that\index{strong approximation!for $SL(2,\reals)$}
\begin{align}
SL(2, \mathbb{A})=SL(2,\mathbb{Q})SL(2,\mathbb{R})K_{\mathbb{A}} .
\end{align}
This decomposition ensures that any automorphic form $\varphi$ on $SL(2,\mathbb{Z})\backslash SL(2,\mathbb{R})/SO(2,\mathbb{R})$ corresponds to an automorphic form on $SL(2, \mathbb{Q})\backslash SL(2,\mathbb{A})/K_{\mathbb{A}} $. 

For the adelic group $SL(2,\mathbb{A})$ we have the Iwasawa decomposition
\begin{align}
SL(2,\mathbb{A})=N(\mathbb{A})A(\mathbb{A})K_{\mathbb{A}} \,.
\end{align}
Given a generic group element $g=nak$ of $SL(2,\mathbb{A})$, we define the character $\chi$ on $SL(2,\mathbb{A})$ in analogy with the general definition~\eqref{characterBorel} such that
\begin{align}
\label{sl2char}
\chi(nak)=|a^{\lambda+\rho}|\,,
\end{align}
where $\lambda$ is some weight vector of $\mathfrak{sl}(2,\mathbb{A})$ and $\rho$ is the Weyl vector. In the case of $SL(2,\mathbb{A})$ the space of (complex) weights is one-dimensional and spanned by the fundamental weight $\Lambda_1$ dual to the unique simple root $\alpha_1$ of $\mathfrak{sl}(2,\mathbb{A})$. The Weyl vector $\rho$ is also identical to $\Lambda_1$. Therefore, we can parametrise the weight appearing in (\ref{sl2char}) with a single parameter $s\in\cx$ as
\begin{align}
\label{charparam}
\lambda= 2s \Lambda_1 -\rho=(2s-1)\Lambda_1\quad\Rightarrow \quad \lambda+\rho = 2s\Lambda_1\,.
\end{align}

Furthermore, we make use of  the function $H(g)$ of~\eqref{eq:Hlog}, which denotes the Lie algebra element associated with the abelian part $a$ in the Iwasawa decomposition of $g=nak$, such that $a=\exp(H(g))$. With this function the character can now be written as
\begin{align}
\label{chis}
\chi_s(g)\equiv |a^{\lambda+\rho}|=e^{(\lambda+\rho)(H(g))}=e^{(\lambda+\rho)(H(a))}= e^{2s\Lambda_1(H(a))}\,,
\end{align}
where we have introduced the notation $\chi_s$ for the character parametrised by $s\in\cx$. In the case of $SL(2,\mathbb{A})$ we can write all these objects explicitly (in the fundamental representation) as $(2\times 2)$-matrices as follows:
\begin{align}
\label{Iwasawa}
g = nak = \left(\begin{array}{cc}
1 & u \\
& 1\\
\end{array} \right) 
\left(\begin{array}{cc}
v &  \\
& v^{-1}\\
\end{array} \right)
\, k
\end{align}
with $k\in K_{\mathbb{A}} $ of (\ref{KA}). Here, $u$ and $v$ are adelic numbers. In keeping with the identification of the upper half plane (cf. also section~\ref{sec:SL2})
\begin{align}
\UHP=\left\{ z=x+iy \in \cx \st \Im(z)>0\right\}=SL(2,\reals)/SO(2,\reals),
\end{align}
this requires that at the archimedean place $p=\infty$ we have to use the following parametrisation
\begin{align}
\label{IwasawaArch}
g_\infty = n_\infty a_\infty k_\infty  = \left(\begin{array}{cc}
1 & x \\
& 1\\
\end{array} \right) 
\left(\begin{array}{cc}
y^{1/2} &  \\
& y^{-1/2}\\
\end{array} \right) k_\infty
\,,
\end{align}
with $y>0$ and $k_\infty\in SO(2)$.
Evaluated on the group element (\ref{Iwasawa}), the character (\ref{chis}) yields
\begin{align}
\label{SL2chis}
\chi_s(g) = e^{2s\Lambda_1(H(a))} = |v|^{2s}
\end{align}
since $H(a) = \log |v|\cdot H_1$ where $H_1$ is the Cartan generator of $SL(2,\mathbb{A})$ and the norm is the adelic one. For the archimedean place this implies with~\eqref{IwasawaArch} that $\chi_s(g_\infty)= y^s$, where we have embedded $g_\infty$ into $G(\ads)$ as $g=(g_\infty,1,1,\ldots)$.

The adelic Eisenstein series $E(\chi,g)$ is then defined by summing the character over a coset according to ($g\in SL(2,\mathbb{A})$)
\begin{align}
E(\chi_s,g)=\sum_{\gamma\in B(\mathbb{Q})\backslash SL(2,\mathbb{Q})} e^{(  \lambda +\rho)(H(\gamma g))}\,,
\end{align}
where the Borel subgroup $B(\mathbb{Q})=N(\mathbb{Q})A(\mathbb{Q})$. Recall that this is the definition of the $SL(2,\mathbb{A})$ Eisenstein series attached to the induced representation which was given in equation~\eqref{generalSL2AEisenstein} of example~\ref{EisensteinSL2example}. The sum converges absolutely for $\Re(s)>1$.

We are now ready to state the theorem:
\begin{theorem}[Fourier expansion $SL(2,\mathbb{A})$ Eisenstein series]\label{SL2AESexp} The expansion of the $SL(2,\ads)$ Eisenstein series $E(\chi_s,g)$ with respect to the unipotent radical $N$ of $SL(2,\mathbb{A})$ is given for $g=g_\infty\in SL(2,\reals)\subset SL(2,\ads)$ by:
\begin{align}
E(\chi_s,g)= \sum_{\psi} W_\psi(s,g)&= y^{s}+\frac{\xi(2s-1)}{\xi(2s)}y^{1-s}\nn\\
&+\sum_{m\neq 0}\frac{2}{\xi(2s)} y^{1/2} |m|^{s-1/2}\sigma_{1-2s}(m)K_{s-1/2}(2\pi|m|y) e^{2\pi i m x},
\end{align}
where there terms on the right-hand side of the first line constitute the constant term and the second line provides the non-constant terms. Here, we have used the parametrisation~\eqref{IwasawaArch} for $g_\infty\in SL(2,\reals)$. 

Furthermore, the Eisenstein series satisfies the functional relation\index{functional relation!for SL(2,R)@for $SL(2,\reals)$}
\begin{align}
\label{SL2funcrel}
E(\chi_s,g) = \frac{\xi(2s-1)}{\xi(2s)} E(\chi_{1-s},g).
\end{align}
\end{theorem}

\begin{proof} 
The proof of this theorem constitutes the rest of the present chapter.
\end{proof}

To prove the theorem we now wish to analyse the Fourier expansion of $E(\chi_s,g)$ along the unipotent radical $N$. This was already outlined in section~\ref{sec_FCsl2ex}. According to the general discussion of the previous chapter, we have the following expansion
\begin{align}
E(\chi_s,g)=\sum_{\psi \in \textrm{Hom}(N(\mathbb{Q})\backslash N(\mathbb{A}),U(1))} W^\circ_\psi(s,g).
\end{align}
We recall from section~\ref{sec:FCsec} that the superscript indicates that the Whittaker coefficient  $W^\circ_\psi$ is spherical, i.e., $K$-independent: $W^\circ_\psi(nak) = W^\circ_\psi(na)$. \index{Whittaker coefficient!for $SL(2,\reals)$}
We shall distinguish the `constant' Whittaker coefficient $W^\circ_1(s,g)$ corresponding to the special case of a trivial character $\psi\equiv 1$
\begin{align}
W^\circ_1(s,g)=\int\limits_{N(\mathbb{Q})\backslash {N}(\mathbb{A})}E(\chi_s,ng)dn\,,
\end{align}
and the remaining `non-constant' coefficients given by 
\begin{align}
W^\circ_\psi(s,g)=\int\limits_{N(\mathbb{Q})\backslash {N}(\mathbb{A})}E(\chi_s,ng)\overline{\psi(n)}dn.
\end{align}
The expressions $W_1(s,g)$ and $W_\psi(s,g)$ are sometimes simply referred to as the `constant term' and the `Fourier coefficients', respectively.

Plugging-in the definition of the Eisenstein series, and interchanging the sum and integration, we can rewrite the coefficients 
in the following form
\begin{subequations}
\begin{align}
\label{SL2cst}
W_1^\circ(s,g)&=\sum_{\gamma\in B(\mathbb{Q})\backslash SL(2,\mathbb{Q})}\; \int\limits_{N(\mathbb{Q})\backslash {N}(\mathbb{A})} \chi_s(\gamma ng) dn,\\
\label{SL2noncst}
W_\psi^\circ(s,g)&= \sum_{\gamma\in B(\mathbb{Q})\backslash SL(2,\mathbb{Q})}\; \int\limits_{N(\mathbb{Q})\backslash {N}(\mathbb{A})} \chi_s(\gamma ng) \overline{\psi(n)} dn,
\end{align}
\end{subequations}
where we recall that the sums converge absolutely for $\Re(s)>1$.

We now proceed with the analysis of the constant and non-constant terms, starting with the constant term.

\section{Constant term}
\label{sec:Sl2const}
As seen in more detail in section~\ref{sec:Bruhat}, the constant term~\eqref{SL2cst} can be re-written as
\begin{align}
W_1^\circ(s,g)&=\sum_{\gamma\in B(\mathbb{Q})\backslash SL(2,\mathbb{Q})}\; \int\limits_{N(\mathbb{Q})\backslash {N}(\mathbb{A})} \chi_s(\gamma ng) dn
\nn \\
&=\sum_{\gamma \in B(\mathbb{Q})\backslash G(\mathbb{Q})/B(\mathbb{Q})}\;\sum_{\delta\in \gamma^{-1}B(\mathbb{Q})\gamma\cap B(\mathbb{Q})\backslash B(\mathbb{Q})}\; \int\limits_{N(\mathbb{Q})\backslash N(\mathbb{A})} \chi_s(\gamma \delta na) dn.
\label{trick}
\end{align}
Because of the quotient by $B(\mathbb{Q})$ on the left in the original $\gamma$ sum, we must make sure to not overcount the coset representatives $\delta$ and this is achieved by the restriction on the $\delta$ sum. To simplify the integral further we shall need  the following result:
\begin{proposition}[Bruhat decomposition]
\begin{align}
\label{BruhatSL2}
SL(2,\mathbb{Q}) = \bigcup_{w\in \Weyl}B(\mathbb{Q})w B(\mathbb{Q}).
\end{align}
\end{proposition}
\begin{proof} To establish (\ref{BruhatSL2}) we begin by noting that for the first coset representative $w$, we have the double coset $Bw B=B$ and for the second coset representative\index{Bruhat decomposition!for $SL(2,\reals)$} we get
\begin{align}
B\begin{pmatrix}&1\\-1&\end{pmatrix} B
&=\left\{ \begin{pmatrix}a&b\\&a^{-1}\end{pmatrix}\begin{pmatrix}&1\\-1&\end{pmatrix} 
\begin{pmatrix}\tilde{a}&\tilde{b}\\&\tilde{a}^{-1}\end{pmatrix}\,:\,b,\tilde{b}\in\mathbb{Q},\,a,\tilde{a}\in\mathbb{Q}^\times\right\}\nonumber\\
&=\left\{ \begin{pmatrix}-b\tilde{a}&a\tilde{a}^{-1}-b\tilde{b}\\-a^{-1}\tilde{a}&-a^{-1}\tilde{b}\end{pmatrix}\,:\,b,\tilde{b}\in\mathbb{Q},\,a,\tilde{a}\in\mathbb{Q}^\times\right\}\nonumber\\
&=\left\{ \begin{pmatrix}a&\frac{ad-1}{c}\\c&d\end{pmatrix}\,:\,a,d\in\mathbb{Q},\,c\in\mathbb{Q}^\times\right\}
\end{align}
and hence the Bruhat decomposition (\ref{BruhatSL2}) corresponds to the division of $SL(2,\mathbb{Q})$ into those matrices with lower left entry equal to zero and those where it is non-zero.
\end{proof}

Using this we can now unfold  the $\delta$ sum in \eqref{trick} to the integration domain by enlarging it, which yields
\begin{align}
W_1^\circ(s,g)= \sum_{\gamma\in B(\mathbb{Q})\backslash SL(2,\mathbb{A})/B(\mathbb{Q})}\; \int\limits_{\gamma^{-1} B(\mathbb{Q}) \gamma\cap N(\mathbb{Q})\backslash N(\mathbb{A})} \chi_s(\gamma ng )dn.
\end{align}
 The measure on this larger space is induced from the embedding $N(\mathbb{Q})\to N(\mathbb{A})$.

We can simplify the summation even further by using the embedding of the Weyl group $\mathcal{W}$ into $K(\mathbb{Q})$; the sum over cosets has only two contributions arising from the trivial and non-trivial coset representatives that can be chosen as
\begin{align}
w=\left(\begin{array}{cc}
1 & \\
& 1\\
\end{array} \right)
\quad \text{or}\quad  
\left(\begin{array}{cc}
 & 1\\
-1 & \\
\end{array} \right).
\label{Weyl}
\end{align}
These correspond precisely to the fundamental Weyl reflections of the Weyl group $\mathcal{W}$ of the Lie algebra $\mathfrak{sl}(2,\mathbb{Q})$. 
Denoting the coset representatives by Weyl words $w$,  we can therefore write the constant term as 
\begin{align}
W^\circ_1(s,g)=\sum_{w\in \Weyl}C_w=\sum_{w\in \Weyl}\; \int\limits_{w^{-1} B(\mathbb{Q}) w\cap N(\mathbb{Q})\backslash N(\mathbb{A})} \chi_s(w n g )dn\,,
\end{align}
where we have defined individual contributions $C_w$ to the constant term that are labelled by elements of the Weyl group.

\subsection{Trivial Weyl word}

In the case when the Weyl word is the trivial Weyl reflection, i.e. $w=\id$, the integral reduces to 
\begin{align}
C_\id= \int\limits_{B(\mathbb{Q}) \cap N(\mathbb{Q})\backslash N(\mathbb{A})} \chi_s( n g )dn&=\int\limits_{N(\mathbb{Q}) \backslash N(\mathbb{A})} \chi_s( ng )dn\nn\\&=|v|^{2s} \int\limits_{N(\mathbb{Q}) \backslash N(\mathbb{A})} dn=|v|^{2s},
 \label{trivialorbit}
\end{align}
where we have used the fact that the Haar measure on $N(\mathbb{Q}) \backslash N(\mathbb{A})$ is normalised to 1,
and applied the definition~\eqref{SL2chis} of the character $\chi_s$ for the Iwasawa decomposed group element $g$.

\subsection{Non-trivial Weyl word}

When $w$ is the non-trivial Weyl reflection in (\ref{Weyl}), it is clear that we have a trivial intersection
\begin{align}
w^{-1} B(\mathbb{Q}) w\cap N(\mathbb{Q}) =\left\{ \id \right\}
\end{align}
and hence the integral for the non-trivial Weyl word simplifies to
\begin{align}
C_w=\int\limits_{w^{-1} B(\mathbb{Q}) w\cap N(\mathbb{Q})\backslash N(\mathbb{A})} \chi_s(w n g )dn=\int\limits_{{N}(\mathbb{A})} \chi_s(wng)dn.
\end{align}
To evaluate the integral we first note that we can restrict the argument to $\chi_s(wng)=\chi_s(wna)$ since we integrate over $N(\mathbb{A})$ and $\chi_s$ is trivial on $K_{\mathbb{A}} $.  Therefore, we have to evaluate
\begin{align}
\int\limits_{{N}(\mathbb{A})} \chi_s(wna)dn.
\end{align}
Now we choose a parametrisation of $N(\mathbb{A})$ by
\begin{align}
\label{npar}
N(\mathbb{A})=\left\{ \left(\begin{array}{cc}
1 & u\\
 & 1\\
\end{array} \right)\, | \, u\in \mathbb{A}\right\},
\end{align}
and of $a$ as in (\ref{Iwasawa}) to write the integral explicitly as
\begin{align}
\int\limits_{{N}(\mathbb{A})} \chi_s(wna)dn 
= \int\limits_{A} \chi_s \Bigg(
\underbrace{\begin{pmatrix}&1\\-1\end{pmatrix}}_{w}
\underbrace{\begin{pmatrix}1&u\\&1\end{pmatrix}}_{n}
\underbrace{\begin{pmatrix}v&\\&v^{-1}\end{pmatrix}}_{a}
\Bigg)d u
\end{align}
We now want to separate out how the integral depends on $a$. This is done by writing 
\begin{align}
\label{amove}
wna = wa a^{-1}na =(waw^{-1}) w(a^{-1}na).
\end{align}
The $a$-dependence comes from both parentheses in this relation. The factor in the first parenthesis is 
\begin{align}
waw^{-1} = \begin{pmatrix}&1\\-1\end{pmatrix}
\begin{pmatrix}v&\\&v^{-1}\end{pmatrix}
\begin{pmatrix}&1\\-1\end{pmatrix}=
\begin{pmatrix}-v^{-1}\\&-v\end{pmatrix}
\end{align}
and lies in $A(\mathbb{A})$. It can therefore be extracted from the character $\chi_s$ using $\chi_s(waw^{-1})=|v|^{-2s}$ by the definition (\ref{SL2chis}) of $\chi_s$ and using the multiplicative properties of $\chi_s$. 

The factor in the second parenthesis in (\ref{amove}) is a conjugation of $N(\mathbb{A})$ by a diagonal element $a$ and can be undone by a change of integration variable. Explicitly, we have
\begin{align}
a^{-1} n a = \begin{pmatrix}v^{-1} &\\&v\end{pmatrix}
\begin{pmatrix}1&u\\&1\end{pmatrix}
\begin{pmatrix}v &\\&v^{-1}\end{pmatrix}
= \begin{pmatrix}1 &v^{-2}u &\\&1\end{pmatrix}
\end{align}
Making the change of variables $u\rightarrow v^{-2}u$ that maps the (Haar) measure $du\rightarrow |v|^2 du$, we can combine the contributions from the two parentheses in (\ref{amove}) to obtain the $a$-dependence
\begin{align}
\label{nontrivWadep}
\int\limits_{\mathbb{A}} \chi_s(wna)dn = \underbrace{|v|^{-2s}}_{\chi_{s}(waw^{-1})} \underbrace{|v|^{2}}_{\text{change of $du$}} \int\limits_{\mathbb{A}} 
\chi_s\Bigg(
\underbrace{\begin{pmatrix}&1\\-1\end{pmatrix}}_{w}
\underbrace{\begin{pmatrix}1&u\\&1\end{pmatrix}}_{n}
\Bigg)du.
\end{align}
In order to evaluate the remaining integral, we rewrite the character according to
\begin{align}
\chi_s(wn)=\chi_s(wnw^{-1}w)=\chi(wnw^{-1}),
\end{align}
where we have used that we have embedded the Weyl group in $K_{\mathbb{A}} $ and the fact that $\chi_s$ is right invariant under $K_{\mathbb{A}} $. Inserting the explicit parametrisations for $w$ and $n$ we find
\begin{align}
wnw^{-1} = \begin{pmatrix}1&\\-u&1\end{pmatrix}.
\end{align}
We see that the Weyl transformation $w$ maps the upper triangular element into a lower triangular element as expected since the non-trivial $w$ maps the (unique) positive root of $SL(2,\mathbb{A})$ to the unique negative root. To evaluate the character $\chi_s$ we will need to perform an Iwasawa decomposition of its argument. By Langlands' theory, see~\cite{LanglandsEP}, the  integral (\ref{nontrivWadep} enjoys complete factorisation into a product
\begin{align}
\int\limits_{\mathbb{A}} 
\chi_s\Bigg(
\begin{pmatrix}&1\\-1\end{pmatrix}
\begin{pmatrix}1&u\\&1\end{pmatrix}
\Bigg)du=\prod_{p\leq \infty} \int_{\mathbb{Q}_p} \chi_{s,p}\Bigg(
 \left(\begin{array}{cc}
1 & 0\\
-u & 1\\
\end{array} \right) \Bigg)du
\label{factorisation}
\end{align}
such that one can analyse the integrals for each prime $p$ separately. 

\vspace{5mm}

\noindent{\bf Archimedean place $p=\infty$}.
We first prove the following result for the archimedean integral corresponding to the real prime at infinity $\mathbb{Q}_\infty=\mathbb{R}$. 
\begin{lemma}[Gindikin--Karpelevich formula for $SL(2,\mathbb{R})$]
At the archimedean place the integral (\ref{factorisation}) evaluates to:
\beq
\int_{\mathbb{R}} \chi_{s,p}\Bigg(
 \left(\begin{array}{cc}
1 & 0\\
-u & 1\\
\end{array} \right) \Bigg)du=\sqrt{\pi}\frac{\Gamma(s-1/2)}{\Gamma(s)}\,.
\eeq
\end{lemma}
\begin{proof}

At the archimedean place, we denote the parameters of $SL(2,\reals)/SO(2,\reals)$ by $x$ and $y^{1/2}$ rather than $u$ and $v$ as shown in~\eqref{IwasawaArch}. The integral then becomes
\begin{align}
\int_{-\infty}^{\infty} \chi_{s,\infty}\Bigg( \left(\begin{array}{cc}
1 & 0\\
-x & 1\\
\end{array} \right) \Bigg)dx.
\end{align}
In order to evaluate this we must bring the argument of the character into Iwasawa form, i.e. we must find $n\in N$ and $a\in A$ such that 
\begin{align}
 \left(\begin{array}{cc}
1 & 0\\
-x & 1\\
\end{array} \right)  =na k, 
\label{archargument}
\end{align}
for some $k\in K_\infty=SO(2,\mathbb{R})$. This was done in example~\ref{Iwex} with the result (cf.~\eqref{IwaLowerSL2}) that
\begin{align}
\left(\begin{array}{cc}
1 & 0\\
-x & 1\\
\end{array} \right)
= \left(\begin{array}{cc}
1 & \frac{-x}{1+x^2}\\
 & 1\\
\end{array} \right)  \left(\begin{array}{cc}
\sqrt{1+x^2}^{-1}  & \\
 & \sqrt{1+x^2}\\
\end{array} \right)k
\end{align}
and the character therefore evaluates to
\begin{align}
\label{Iwreal}
\chi_{s,\infty}\Bigg( \left(\begin{array}{cc}
1 & 0\\
-x & 1\\
\end{array} \right) \Bigg)
=\chi_{s,\infty}\Bigg(
 \left(\begin{array}{cc}
\sqrt{1+x^2}^{-1}  & \\
 & \sqrt{1+x^2}\\
\end{array} \right) \Bigg)
=(\sqrt{1+x^2})^{-2s}.
\end{align}
We then find for the integral 
\begin{align}
\int_{-\infty}^{\infty} \chi_{s,\infty}\Bigg( \left(\begin{array}{cc}
1 & 0\\
-x & 1\\
\end{array} \right)\Bigg)dx=\int_{-\infty}^{\infty} (\sqrt{1+x^2})^{-2s} dx=\sqrt{\pi}\frac{\Gamma(s-1/2)}{\Gamma(s)}\,.
\label{archimedean}
\end{align}
\end{proof}
\vspace{5mm}

\noindent{\bf Non-archimedean places $p<\infty$}.  We shall now prove the corresponding result for finite primes: 
\begin{lemma}[Gindikin--Karpelevich formula for $SL(2,\mathbb{Q}_p)$ \cite{LanglandsEP}]
 \begin{align}
 \int_{\mathbb{Q}_p} \chi_{s,p}\Bigg( \left(\begin{array}{cc}
1 & 0\\
-u & 1\\
\end{array} \right) \Bigg)du&=1+\frac{p-1}{p}\,\frac{p^{-2s+1}}{1-p^{-2s+1}}
=\frac{1-p^{-2s}}{1-p^{-2s+1}}.
\end{align}
\end{lemma}

\begin{proof}

We now consider the terms in (\ref{factorisation}) for which $p$ is a finite prime:
\begin{align}
\int_{\mathbb{Q}_p} \chi_{s,p}\Bigg( \left(\begin{array}{cc}
1 & 0\\
-u & 1\\
\end{array} \right) \Bigg)du.
\label{nonarchimedeanintegral1}
\end{align}
The Iwasawa decomposition of this element was also discussed in example~\ref{Iwex}. When $u\in \mathbb{Z}_p$ we have
\begin{align}
 \left(\begin{array}{cc}
1 & 0\\
-u & 1\\
\end{array} \right)\in SL(2,\mathbb{Z}_p)\,.
\end{align}
The compact part of $SL(2,\mathbb{Q}_p)$ is $K_{p}=SL(2,\mathbb{Z}_p)$ and hence, since $\chi_s$ is trivial on $K$, we have
\begin{align}
\chi_{s,p}\left( \Bigg(\begin{array}{cc}
1 & 0\\
-u & 1\\
\end{array} \right)  \Bigg)=1\,.
\end{align}
We may thus split the integral into
\begin{align}
\int_{\mathbb{Z}_p}du+\int_{\mathbb{Q}_p\backslash\mathbb{Z}_p} \chi_{s,p}\Bigg( \left(\begin{array}{cc}
1 & 0\\
-u & 1\\
\end{array} \right)  \Bigg)du\,,
\label{nonarchimedeanintegral2}
\end{align}
where the first term is  unity by choice of normalisation~\eqref{padicintmesnorm} for the measure $du$. When $u\in \mathbb{Q}_p$ but not in $\mathbb{Z}_p$, i.e. $|u|_p>1$, we write the matrix in an Iwasawa decomposition (which is not unique, cf. example~\ref{Iwex})
\begin{align}
 \left(\begin{array}{cc}
1 & 0\\
-u & 1\\
\end{array} \right) = \left(\begin{array}{cc}
 1 & *\\
 & 1\\
\end{array} \right)\left(\begin{array}{cc}
u^{-1} & \\
 & u\\
\end{array} \right)k.
\end{align}
The remaining integral becomes
\begin{align}
\int_{\mathbb{Q}_p\backslash\mathbb{Z}_p} \chi_{s,p}\Bigg(
 \left(\begin{array}{cc}
 1 & *\\
 & 1\\
\end{array} \right)\left(\begin{array}{cc}
u^{-1} & \\
 & u\\
\end{array} \right)k  \Bigg)du=\int_{\mathbb{Q}_p\backslash\mathbb{Z}_p} |u|_p^{-2s} du.
\end{align}
We recognize this as an integral of the type (\ref{usefulintegral}) that we already evaluated, so the result is
\begin{align}
 \int_{\mathbb{Q}_p\backslash\mathbb{Z}_p} |u|_p^{-2s} du=\frac{p-1}{p}\,\frac{p^{-2s+1}}{1-p^{-2s+1}}.
 \end{align}
 Combining this with the first term in (\ref{nonarchimedeanintegral2}) proves the claim.
 \end{proof}
 
 \begin{remark}
Integrals of this type will be evaluated more generally inequation~\eqref{GindiKarp}, leading to a more general \emphindex[Gindikin--Karpelevich formula!for $SL(2,\rats_p)$]{Gindikin--Karpelevich formula}.
\end{remark}

\subsection{The global form of the full constant term}

We are now ready to assemble all the pieces and write down the complete result for the constant term $W_1(s,g)$. The only remaining step is to compute the product over all finite primes in (\ref{factorisation}). Recalling the Euler product representation of the Riemann zeta function (\ref{RiemannProd}), we find
\begin{align}
\prod_{p< \infty} \int_{\mathbb{Q}_p} \chi_{s,p}\left( \left(\begin{array}{cc}
1 & 0\\
-u & 1\\
\end{array} \right) \right)du=\prod_{p<\infty} \frac{1-p^{-2s}}{1-p^{-2s+1}}=\frac{\zeta(2s-1)}{\zeta(2s)}.
\end{align}
Combining this with the result from the archimedean integral (\ref{archimedean}), including the overall pre-factor from~\eqref{nontrivWadep}, as well as the contribution from the trivial Weyl word in (\ref{trivialorbit}) we finally find
\begin{align}
W_1^\circ(s,g)= |v|^{2s}+\sqrt{\pi}\frac{\Gamma(s-1/2)}{\Gamma(s)}\frac{\zeta(2s-1)}{\zeta(2s)}|v|^{-2s+2}.
\end{align}
Here, we have left $v\in \ads$. Restricting to $g\in SL(2,\reals)$ we can write this in terms of $v=y^{1/2}$ instead, 
with the result that the first term 
scales like $y^{s}$ while the second term scales as $y^{1-s}$. The relation between the exponents is that induced from the non-trivial Weyl reflection:
\begin{align}
 s\to 1-s\,.
 \end{align}
Referring back to our particular parametrisation~\eqref{sl2char} and~\eqref{charparam} of the character $\chi_s$, we recall that for $w$ being trivial, we have $\lambda+\rho=2s\Lambda_1$, while for $w$ being the non-trivial Weyl word we obtain $w\lambda+\rho=2(1-s)\Lambda_1$. Hence the parameter $s$ is seen to be related by the above transformation. Recall from section~\ref{sec-complRiem} that the completed zeta function is given by
\begin{align}
 \xi(s)=\pi^{-s/2}\Gamma(s/2)\zeta(s),
 \label{completedZeta}
 \end{align}
and satisfies the functional relation
 \begin{align}
 \xi(s)= \xi(1-s).
 \end{align}
 Using this we can  write the constant term in the following compact way
 \begin{align}
 W_1^\circ(s,g)= |v|^{2s}+\frac{\xi(2s-1)}{\xi(2s)}|v|^{-2s+2}.
 \end{align}
For $v=y^{1/2}$ this agrees with the constant term in~\eqref{SL2FC2} and the statement of theorem~\ref{SL2AESexp}. We note that constant terms therefore satisfy the functional relation
\begin{align}
\label{funcSL2cst}
W_1^\circ(s,g) = \frac{\xi(2-2s)}{\xi(1-2s)} C(1-s,g)=\frac{\xi(2s-1)}{\xi(2s)} C(1-s,g),
\end{align}
where the functional relation for the completed Riemann zeta function~\eqref{CRZfuncrel} has been used.

\section{The non-constant Fourier coefficients}
\label{sec:SL2FC}

The Whittaker coefficient $W^\circ_\psi$ of~\eqref{SL2noncst} we want to compute is given by
\begin{align}
W^\circ_\psi(s,g)=\sum_{w\in\Weyl}F_{w,\psi}=\sum_{w\in \Weyl}\; \int\limits_{w^{-1} B(\mathbb{Q}) w\cap N(\mathbb{Q})\backslash N(\mathbb{A})} \chi_s(w n g )\overline{\psi(n)}dn\,.
\end{align}
Using the Iwasawa decomposition of $g=n'ak$ and performing a change of variables this expression can be re-written as
\begin{align}\label{WPsi}
W^\circ_\psi(s,g)=\overline{\psi(n'^{-1})}\sum_{w\in \Weyl}\; \int\limits_{w^{-1} B(\mathbb{Q}) w\cap N(\mathbb{Q})\backslash N(\mathbb{A})} \chi_s(w n a )\overline{\psi(n)}dn\,.
\end{align}

As for the constant term case, the sum over the Weyl group has two contributions, one each, for when the Weyl word is trivial, $\id$, and non-trivial, w,
\begin{align}
W^\circ_\psi=F_{\id,\psi}+F_{w,\psi}\,.
\end{align}
The two contributions will be treated separately below.

Given our standard parametrisation of $n$ by the variable $u$, we define the character $\psi$ (against which we integrate) as a direct product
\begin{align}
\psi = \prod_{p\leq \infty} \psi_{p}
\end{align}
with
\begin{align}
\label{charprod}
\psi_{p}(u)=\begin{cases}
    e^{2\pi imu}& \text{  for $p=\infty$}\,,\\
    e^{-2\pi i[mu]}& \text{  for $p<\infty$}\,.
  \end{cases}
\end{align}
The function $[\,\cdot\,]$ returns the fractional part of a $p$-adic number as defined in (\ref{fracpart}). An important point here is that we are interested in characters $\psi$ of the continuous group $N(\rats)$ embedded diagonally in $N(\mathbb{A})$. Therefore the coefficient $m$ is a rational number and identical in all $\psi_{p}$. 

In equation~(\ref{charprod}), we have not indicated the conductor $m$ as a subscript on the character $\psi_{p}$ in contrast to section~\ref{sec:chars} in order to keep the notation light. It is always understood that the conductor is $m$. Note that for the pre-factor in~\eqref{WPsi} we have $\overline{\psi(n'^{-1})}=\psi(n')$. 

\subsection{Trivial Weyl word}
 
In the trivial case, i.e. when $w$ is equal to the identity matrix, the integral in~\eqref{WPsi} takes the form
\begin{align}
\int\limits_{N(\mathbb{Q})\backslash N(\mathbb{A})}\chi_s(na)\overline{\psi(n)}dn\,.
\end{align}
As before, we use the definition $\chi_s(na)=|v|^{2s}$. 

The complete expression for the `trivial' term of the Fourier coefficient is then given by
\begin{align}
F_{\id,\psi}(s,g) = \psi(n')|v|^{2s} \int\limits_{N(\mathbb{Q})\backslash N(\mathbb{A})} \overline{\psi(n)}dn\,.
\end{align}
We now proceed to write this expression as a product over all primes, including the place at infinity, as
\begin{align}
F_{\id,\psi}(s,g)=\prod\limits_{p\leq\infty}\psi_{p}(n')|v|_p^{2s} \int\limits_{N(\mathbb{Z}_p)\backslash N(\mathbb{Q}_p)}\overline{\psi_p(u)}du\,.
\end{align}
This has to be evaluated separately for each prime $p\leq\infty$. Starting with the archimedean $p=\infty$,
the domain of integration is $\ints\backslash\reals\cong [0,1]$. This leads to the integral
\begin{align}
F_{\id,\psi,\infty}(s,g)=\psi_\infty(n')|v|_\infty^{2s} \int\limits_{0}^{1} e^{-2\pi imu}du =0
\end{align}
since this is the integral of a periodic function (with mean value zero) over a full period. Therefore, the full Fourier coefficient vanishes for the trivial Weyl word:
\begin{align}
F_{\id,\psi}(s,g)= 0.
\end{align}
This is the reflection of a general phenomenon that will be discussed in section~\ref{sec:redlong} below.

\subsection{Non-trivial Weyl word}

In the case when the Weyl word is non-trivial, with representation
\begin{align}
w=\left(
\begin{array}{cc}
0&1\\
-1&0
\end{array}
\right)\,,
\end{align}
the corresponding term in the Fourier coefficient reads
\begin{align}
F_{w,\psi}(s,g)=\psi(n')\int\limits_{N(\mathbb{A})}\chi_s(wna)\overline{\psi(n)}dn\,.
\end{align}
We now perform the same transformation (\ref{amove}) to remove the $a$-dependence from $\chi_s$.
Under a change of variables $n\rightarrow ana^{-1}$, the integration measure transforms as $dn\rightarrow |v|^2dn$, and we obtain
\begin{align}
\psi(n')|v|^2\int\limits_{N(\mathbb{A})}\chi_s(wan)\overline{\psi(ana^{-1})}dn\,.
\end{align}
Inserting $w^{-1}w$ in the argument before and after $n$, we find for the character $\chi_s(wan)=\chi_s(waw^{-1})\chi_s(wnw^{-1}w)=|v|^{-2s}\chi_s(wnw^{-1})$, where we have again used the fact that $\chi$ is right invariant under a Weyl group transformation. The full expression then takes the form
\begin{align}
\psi(n')|v|^{-2s+2}\int\limits_{N(\mathbb{A})}\chi_s(wnw^{-1})\overline{\psi(ana^{-1})}dn\,.
\end{align}
Now we write the expression in the standard way as a product over all places
\begin{align}
F_{w,\psi}(s,g)=\prod\limits_{p\leq\infty}\psi_{p}(n')|v|_p^{-2s+2}\int\limits_{N(\mathbb{A})}\chi_{s,p}(wnw^{-1})\overline{\psi_{p}(ana^{-1})}dn
\end{align}
and evaluate the archimedean and non-archimedean places separately.

\vspace{5mm}

\noindent{\bf The archimedean place $p=\infty$}:
The Iwasawa decomposition of $wnw^{-1}$ is as in (\ref{Iwreal}) and again leads to
$\chi_{s,\infty}(wnw^{-1})=\sqrt{1+x^2}^{-2s}$. Furthermore, the character evaluates to $\psi_{\infty }(ana^{-1})=\exp(2\pi i m y x)$, such that overall we obtain
\begin{align}
\label{SL2Whittinf}
F_{w,\psi,\infty}&=\psi_{\infty }(n')|y|_\infty^{-s+1}\int\limits_{-\infty}^{\infty}(1+x^2)^{-s}e^{-2\pi i m y x}dx\nonumber\\
&=\frac{2\pi^s}{\Gamma(s)}  y^{1/2}  |m|^{s-1/2} K_{s-1/2}(2\pi |m| y) \psi_{\infty }(n'),
\end{align}
where we have used the integral representation of the modified Bessel function given in equation~(\ref{Besselreal}) and $y>0$. 

\vspace{5mm}

\noindent{\bf The non-archimedean places $p<\infty$}:
We have to analyse the integral
\begin{align}
\label{SL2Floc}
F_{w,\psi,p}=\psi_{p}(n')|v|_p^{-2s+2}\int\limits_{\mathbb{Q}_p}\chi_{s,p}(wnw^{-1})\overline{\psi_{p}(ana^{-1})}dn\,.
\end{align}
We will set $a=n'=1$ along the finite primes since we are interested in the Eisenstein series as a function on $SL(2,\reals)$ only.
From the Iwasawa decomposition of $wnw^{-1}$ following~\eqref{nonarchimedeanintegral1} we know that $\chi_{s,p}(wnw^{-1})$ is given by 
\begin{align}
\chi_{s,p}(wnw^{-1}) = \max(1,|u|_p)^{-2s},
\end{align}
where $u$ parametrises $N(\rats_p)$ as in (\ref{npar}) and we have to integrate this against the appropriate character
\begin{align}
\int\limits_{\mathbb{Q}_p}\text{max}(1,|u|_p)^{-2s}e^{2\pi i [m u]}du\,.
\end{align}
Using example~\ref{besselfull}, this integral evaluates to
\begin{align}
\label{SL2Wloc}
\int\limits_{\mathbb{Q}_p}\text{max}(1,|u|_p)^{-2s}e^{2\pi i [m u]}du=\gamma_p(m)(1-p^{-2s})\frac{1-p^{-2s+1}|m|_p^{2s-1}}{1-p^{-2s+1}}\,.
\end{align}
Taking the product over all finite places yields
\begin{align}\label{finitSL2coeff}
\prod\limits_{p<\infty}F_{w,\psi,p}=\left(\prod_{p<\infty} (1-p^{-2s})\right)\left(\prod_{p<\infty} \gamma_p(m) \frac{1-p^{-(2s-1)}|m|_p^{2s-1}}{1-p^{-(2s-1)}}\right).
\end{align}
The first factor is equal to $\zeta(2s)^{-1}$ by virtue of~\eqref{RiemannProd}. We can restrict to $m\in\ints$ due to the occurrence of the $p$-adic Gaussian $\gamma_p(m)$ for all $p<\infty$ as seen in section~\ref{sec:p-adic-special-functions}. Writing $m$ then in terms of its unique prime factorisation $m=\prod q_i^{k_i}$ with $q_i$ primes, we can rewrite the second factor (cf. example \ref{ex_localWhittakerSL2}). Consider first the case when $m=q^k$ for a single prime $q$. Then (\ref{qp}) implies that the second factor can be written as a sum over (positive) divisors of $q^k$
\begin{align}
&\left(\prod_{p<\infty} \frac{1-p^{-(2s-1)}|q^k|_p^{2s-1}}{1-p^{-(2s-1)}}\right)
= \frac{1-q^{-(2s-1)} q^{-k(2s-1)}}{1-q^{-(2s-1)}}= \frac{1-q^{-(k+1)(2s-1)}}{1-q^{-(2s-1)}} \nonumber\\
&\quad = \sum_{d|q^k} d^{-(2s-1)} .
\end{align}
By multiplicativity of the expressions, we therefore obtain for a general integral $m$
\begin{align}
\left(\prod_{p<\infty} \gamma_p(m) \frac{1-p^{-(2s-1)}|m|_p^{2s-1}}{1-p^{-(2s-1)}}\right)=\sum_{d|m} d^{-2s+1}=: \sigma_{1-2s}(m),
\end{align}
where we have used the general divisor sum $\sigma_{1-2s}(m)$ over positive divisor of an integer (cf.~(\ref{intro:divisorsum})).

Putting everything together we therefore obtain a non-vanishing coefficient only for integral $m$ with value
\begin{align}
\label{SL2FourierFull}
W^\circ_\psi (s,g) = F_{w,\psi}(s,g)=\frac{2}{\xi(2s)} y^{1/2} |m|^{s-1/2}\sigma_{1-2s}(m)K_{s-1/2}(2\pi|m|y) \psi_{\infty }(n'),
\end{align}
where we have used the definition $\xi(s)=\pi^{-s/2}\Gamma(s/2)\zeta(s)$ for the completed Riemann $\zeta$-function~\eqref{CRZ}.

Finally, we address the functional relation~\eqref{SL2funcrel} for the non-zero Fourier coefficients.
The modified Bessel function has the property $K_{s-1/2}(w) = K_{1/2 -s} (w)$ for all $w>0$. For the divisor sum one finds similarly 
\begin{align}
\sigma_{1-2s}(m) = \sum_{d|m} d^{1-2s} = |m|^{1-2s} \sum_{d|m} \left(\frac{|m|}{d}\right)^{2s-1} = |m|^{1-2s} \sigma_{2s-1}(m).
\end{align}
Putting this together, we obtain
\begin{align}
W^\circ_\psi(s,g) = \frac{\xi(2s-1)}{\xi(2s)} W^\circ_\psi (1-s,g),
\end{align}
where again the functional relation~\eqref{CRZfuncrel} of the Riemann zeta function was used.
This concludes the proof of theorem~\ref{SL2AESexp}.

\chapter{Langlands constant term formula}
\label{ch:CTF}
In this chapter we shall provide a detailed proof of the Langlands constant term formula for Eisenstein series on an arbitrary reductive group $G$. This generalises the results of 
the previous chapter for $G=SL(2)$. We will also discuss the general functional relation satisfied by Eisenstein series and we explain how to define and 
evaluate constant terms with respect to non-maximal unipotent radicals.
\section{Statement of theorem}

We start from the following definition of the minimal parabolic Eisenstein series
\begin{align}
\label{eisdefA}
E(\chi,g)=\sum_{\gamma\in B(\mathbb{Q})\backslash G(\mathbb{Q})}\chi(\gamma g)\,.
\end{align}
This is a valid rewriting since the cosets $B(\mathbb{Q})\backslash G(\mathbb{Q})$ are in bijection with those of $B(\mathbb{Z})\backslash G(\mathbb{Z})$ (see example~\ref{bijection} for a proof of this for $SL(2)$). When writing (\ref{eisdefA}), we can allow $g\in G(\mathbb{A})$. The real function~\eqref{EIntro} is re-obtained by setting $g=(g_\infty,1,1,\ldots)$, i.e., setting the components along $G(\mathbb{Q}_p)$ equal to the identity for $p\neq \infty$.

As in section~\ref{sec_Eisenstein} we parametrise the character $\chi$ by
\begin{align}
\label{charweight}
\chi(nak) = a^{\lambda+\rho}
\end{align}
in terms of a weight $\lambda$ of the Lie algebra and $\rho$ is the Weyl vector defined in~\eqref{Weylvec}.

Our interest in the present section is to evaluating the so-called constant terms in the minimal parabolic subalgebra $B$ (standard Borel); that is we shall prove:
\begin{theorem}[Langlands constant term formula]\label{LCFthm}
    \index{constant term!formula|see{Langlands constant term formula}}
    \index{Langlands constant term formula}
The constant term of $E(\lambda,g)$ with respect to the unipotent radical $N\subset B$ is given by:
\begin{align}
\lint_{N(\mathbb{Q})\backslash N(\mathbb{\mathbb{A}})} E(\lambda,ng) dn = \sum_{w\in \mathcal W} a^{w\lambda+\rho} \prod_{\alpha>0\,|\, w\alpha<0} \frac{\xi(\langle \lambda| \alpha\rangle)}{\xi(1+\langle \lambda| \alpha\rangle)},
\end{align}
where $dn$ is the Haar measure that is normalised such that $N(\mathbb{Q})\backslash N(\mathbb{A})$ has unit volume. 
\end{theorem}

\begin{proof} \mbox{} The proof of this theorem constitutes the greater part of the present chapter and is contained in sections~\ref{sec:Bruhat} to~\ref{consttermass}.
\end{proof} 

Clearly, the constant term (\ref{constterm}) depends only on $a$: For $g=n'ak$ in Iwasawa form, right $K$-invariance and a change of integration variation reduce the dependence to $a$. In what follows we shall therefore define
\begin{align}
\label{constterm}
C(\chi,a)=\int\limits_{N(\mathbb{Q})\backslash N(\mathbb{A})} E(\chi,ng) dn
= \int\limits_{N(\mathbb{Q})\backslash N(\mathbb{A})} E(\chi,na) dn\,,
\end{align}
and we view the integral as a function on the Cartan torus. 

\section{Bruhat decomposition}
\label{sec:Bruhat}

The first step in evaluating (\ref{constterm}) is to use the \emphindex{Bruhat decomposition}~\cite{Borel}:
\begin{align}
G(\mathbb{Q}) = \bigcup_{w\in \mathcal W} B(\mathbb{Q})w B(\mathbb{Q})
\end{align}
that describes the group $G(\mathbb{Q})$ as a disjoint union of double cosets by the Borel subgroup $B(\mathbb{Q})\subset G(\mathbb{Q})$. The group $\mathcal W$ is the Weyl group of $G(\mathbb{R})$. Clearly, we could restrict the group $B(\mathbb{Q})$ on the right to the subgroup generated by those positive step operators that are mapped to negative step operators by the action of $w$. The ones that stay positive are already contained in the Borel subgroup on the left. One can think of the Bruhat decomposition as the extension of the tessellation of the Cartan subalgebra into Weyl chambers to the full group.

Using the same trick as in section~\ref{sec:Sl2const} we can rewrite the constant term as
\begin{align}
\label{Bruhatapplied}
C(\chi,a) &= \sum_{\gamma \in B(\mathbb{Q})\backslash G(\mathbb{Q})} \; \int\limits_{N(\mathbb{Q})\backslash N(\mathbb{A})} \chi(\gamma na) dn&\nonumber\\
&=\sum_{w\in \mathcal W} \;\int\limits_{w^{-1}B(\mathbb{Q})w\cap N(\mathbb{Q})\backslash N(\mathbb{A})} \chi(w na) dn\,.&
\end{align}
Continuing from (\ref{Bruhatapplied}) we look at the individual terms
\begin{align}
\label{constind}
C_w = \int\limits_{w^{-1}B(\mathbb{Q})w\cap N(\mathbb{Q})\backslash N(\mathbb{A})} \chi(w na) dn
\end{align}
and note that the integration domain can be simplified to 
\begin{align}
\label{constsimpl1}
C_w = \int\limits_{N_w(\mathbb{A})} \chi(w na) dn\,,
\end{align}
where $N_w(\mathbb{A})$ is generated from a product over the positive roots that are mapped to negative roots by the given Weyl word $w$
\begin{align}
\label{Uwdef}
N_w(\mathbb{A}) = \prod_{\alpha>0\,|\,w\alpha<0} N_{\alpha}(\mathbb{A})
\end{align}
and $N_{\alpha}(\mathbb{A})$ is the subgroup of $G(\mathbb{A})$ generated by the step operator $E_{\alpha}$ and its dimension is given by the length $\ell$ of the reduced Weyl word $w$.
This simplification of (\ref{constind}) uses two facts:
\begin{enumerate}
\item[$(i)$]  (Upper) Borel elements that get mapped to lower Borel elements by $w$ have trivial intersection with $N(\mathbb{Q})$ and therefore the quotient becomes trivial and leaves an integral over all of $\mathbb{A}$ in that direction. 
\item[$(ii)$] If an (upper) Borel element is mapped to another upper Borel element by the action of $w$, one is left with the integral over the corresponding quotient. However, since the part $wn$ of the argument is then still a Borel element, the character is insensitive to it and one is left with the volume of corresponding Borel directions which is normalised to unity.
\end{enumerate}
Therefore, in (\ref{constsimpl1}) we have carried out many trivial integrals and are only left with the non-trivial integrals where $wn$ is really a lower Borel element.

\section{Parametrising the integral}

We will eventually evaluate integral (\ref{constsimpl1}) using a recursive method and we start by parametrising it conveniently. First, we need to know something about $N_w(\mathbb{A})$ defined in (\ref{Uwdef}). We fix a reduced expression $w=w_{i_1}w_{i_2}\cdots w_{i_\ell}$ for the Weyl word $w$ of length $\ell$. The subscripts refer to the nodes of the Dynkin diagram of $G(\mathbb{Q})$ and $w_i$ are the fundamental reflections that generate the Weyl group. Then one can explicitly enumerate all positive roots that are mapped to negative roots by the action of $w$ as follows. Define
\begin{align}
\gamma_{k} = w_{i_\ell} w_{i_{\ell-1}}\cdots w_{i_{k+1}}\alpha_{i_k}
\end{align}
where $\alpha_{i_k}$ is the $i_k$th simple root. That this gives a valid description of the positive roots generating $U_w$ can be checked easily by induction. Therefore we have
\begin{align}
\{\alpha>0\,|\, w\alpha<0\} = \{\gamma_i \,|\, i=1,\ldots,\ell\}\,.
\end{align}
We also note that there is a simple expression for the sum of all these roots:
\begin{align}
\label{sumgamma}
\gamma_1+\ldots+\gamma_\ell = \rho - w^{-1} \rho
\end{align}
which can again be checked by induction. We note in particular $\gamma_\ell=\alpha_{i_\ell}$. 

In the next step we use the Chevalley basis notation to write elements $u\in N_w(\mathbb{A})$ as
\begin{align}
u = x_{\gamma_1}(u_1)\cdots x_{\gamma_\ell}(u_\ell)
\end{align}
with the Chevalley generator $x_{\alpha}(v)$ being defined by
\begin{align}
x_\alpha(v) = e^{v E_\alpha}\,,
\end{align}
where $E_\alpha$ is the generator of the $\alpha$ root space normalised to unity (for both short and long roots) and $v\in\mathbb{A}$ is the parameter of the group element. With this parametrisation, we can rewrite our individual term $C_w$ of~\eqref{constsimpl1} as
\begin{align}
\label{constsimpl2}
C_w = \int_{\mathbb{A}^\ell} \chi\left(w x_{\gamma_1} (u_1) \cdots x_{\gamma_\ell}(u_\ell)a\right) du_1 \cdots du_\ell\,.
\end{align}

\section{Obtaining the \texorpdfstring{$a$}{a} dependence of the integral}

It is now possible to extract the dependence on $a$ from the integral (\ref{constsimpl2}).
This is done by conjugating the abelian element $a$ through to the left in the argument of $\chi$. The result is 
\begin{align}
C_w =  \int_{\mathbb{A}^\ell} \chi\left(wa x_{\gamma_1} (a^{-\gamma_1}u_1) \cdots x_{\gamma_\ell}(a^{-\gamma_\ell}u_\ell)\right) du_1 \cdots du_\ell\,,
\end{align}
where we have used the fact that Cartan elements act diagonally on the $E_{\gamma_i}$ root spaces. In the next step we can move $w$ past $a$ in the argument of $\chi$ and employ the multiplicative property~\eqref{charmult} of the character to obtain
\begin{align}
\label{constaexplicit}
C_w &= \chi(waw^{-1})  \int_{\mathbb{A}^\ell} \chi\left(w x_{\gamma_1} (a^{-\gamma_1}u_1) \cdots x_{\gamma_\ell}(a^{-\gamma_\ell}u_\ell)\right) du_1 \cdots du_\ell&\nonumber\\
&= \chi(waw^{-1}) a^{\gamma_1+\ldots \gamma_\ell} \int_{\mathbb{A}^\ell} \chi\left(w x_{\gamma_1} (u_1) \cdots x_{\gamma_\ell}(u_\ell)\right) du_1 \cdots du_\ell\,,&
\end{align}
where we have also rescaled the $u$-variables and moved the Jacobi factor outside of the integral. In the form (\ref{constaexplicit}), one can read off the full $a$-dependence of the constant term. Using (\ref{charweight}) and (\ref{sumgamma}) we can rewrite the $a$-dependence of the constant term, remembering that $\chi$ is $\mathcal W$-invariant from the right, as
\begin{align}
\label{Cw}
C_w &= a^{w^{-1}(\lambda+\rho)} a^{\rho-w^{-1}\rho} I_w&\nonumber\\
&= a^{w^{-1}\lambda+\rho} I_w\,,&
\end{align}
where the remaining, $a$-independent integral is
\begin{align}
\label{Iw}
I_w = \int_{\mathbb{A}^\ell} \chi\left(x_{w\gamma_1} (u_1) \cdots x_{w\gamma_\ell}(u_\ell)\right) du_1 \cdots du_\ell\,,
\end{align}
and we have applied $w$ to all the Chevalley elements and used again $K$-invariance of $\chi$ on the right.

\section{Solving the remaining integral by induction}
\label{sec:CSTInd}

We will now solve (\ref{Iw}) by induction. First, we note that $w\gamma_i$ is a negative root for all $i$ by virtue of the definition of $\gamma_i$. Therefore, the factors $x_{w\gamma_i}(u_i)$ appearing in (\ref{Iw}) are elements of the lower triangular Borel subgroup of $G(\mathbb{A})$. To evaluate the character $\chi$ in (\ref{Iw}) according to (\ref{charweight}), we need to perform an Iwasawa decomposition and isolate the $A(\mathbb{A})$ part of the argument of the character $\chi$. We start by Iwasawa decomposing the last Chevalley factor in the argument of the character according to
\begin{align}
\label{Iwaell}
x_{w\gamma_\ell}(u_\ell) = n(u_\ell) a(u_\ell) k(u_\ell)\,.
\end{align}
The (negative) step operator $E_{w\gamma_\ell}$ that enters in $x_{w\gamma_\ell}(u_\ell)$ is part of an $SL(2,\mathbb{A})$ subgroup of $G(\mathbb{A})$ and the Iwasawa decomposition (\ref{Iwaell}) takes place in that subgroup. We choose to label the $SL(2,\mathbb{A})$ subgroup by its {\em positive} root $-w\gamma_\ell$, so that the corresponding Cartan generator is proportional to $H_{-w\gamma_\ell}$. The problem of Iwasawa decomposing the $SL(2,\mathbb{A})$ associated with $-w\gamma_\ell$ is different for $p=\infty$ and $p<\infty$ and will be treated separately in~\ref{pinf} and~\ref{pfin} below. 

Inserting the Iwasawa decomposed (\ref{Iwaell}) into (\ref{Iw}), we can again drop the compact element on the right. Then we can conjugate $n(u_\ell)$ through to the left. For this we have to pass through the negative step operators $x_{w\gamma_i}(u_i)$ for $i=1,\ldots,\ell-1$. This produces an element that can be arranged as a product of nilpotent elements on the left times negative step operators $x_{w\gamma_i}(u'_i)$ for $i=1,\ldots,\ell-1$ with different $u'_i$. The nilpotent elements disappear in the character and the transformation of the space of parameters $u_1$ to $u_{\ell-1}$ is uni-modular.We tacitly perform the corresponding change of variables $u_i\to u_i'$. We therefore obtain that $n(u_\ell)$ can be completely absorbed and we are left with:
\begin{align}
I_w = \int_{\mathbb{A}^\ell} \chi \left( x_{w\gamma_1}(u_1)\cdots x_{w\gamma_{\ell-1}}(u_{\ell-1}) a(u_\ell)\right) du_1\cdots du_\ell.
\end{align}
In the next step, we conjugate $a(u_\ell)$ to the left. This rescales again the $u$ variables with the result
\begin{align}
\label{Iwrec}
I_w &= \int_{\mathbb{A}^\ell} \chi\left(a(u_\ell) x_{w\gamma_1}(a(u_\ell)^{-w\gamma_1}u_1)\cdots x_{w\gamma_{\ell-1}}(a(u_\ell)^{-w\gamma_{\ell-1}}u_{\ell-1})\right)du_1\cdots du_\ell&\nonumber\\
&= \int_{\mathbb{A}} \chi(a(u_\ell)) a(u_\ell)^{w(\gamma_1+\ldots \gamma_{\ell-1})} du_\ell \cdot I_{w'}\,,&
\end{align}
where we have undone the scaling of the variables at the cost of introducing a Jacobi factor and introduced $w'$ through 
\begin{align}
w=w' w_{i_\ell},
\end{align}
i.e., it is obtained from the Weyl word $w$ by removing the right-most fundamental reflection. This is the recursion formula we are after. All that remains now is to evaluate one integral over $\mathbb{A}$.

\section{The Gindikin--Karpelevich formula}

Using the expression (\ref{charweight}) for the character $\chi$, the desired integral is
\begin{align}
\label{Iell}
I_\ell=
\int_{\mathbb{A}} a(u_\ell)^{\lambda+\rho+w(\gamma_1+\ldots \gamma_{\ell-1})} du_\ell =\prod_{p\leq \infty} \int_{\mathbb{Q}_p} a(u_\ell)^{\lambda+\rho+w(\gamma_1+\ldots +\gamma_{\ell-1})} du_\ell,
\end{align}
and can be evaluated for each finite and infinite prime $\leq\infty$ as follows. For this one needs the explicit Iwasawa decomposition expressions for $a(u_\ell)$ that will be derived below for $p=\infty$ and $p<\infty$. We also introduce the notation
\begin{align}
\label{auell}
a(u_\ell)^{\lambda+\rho+w(\gamma_1+\ldots \gamma_{\ell-1})} = |\phi_\ell|^{z_\ell+1}
\end{align}
with (\ref{charweight}) and
\begin{align}
\label{zell}
z_\ell= (\lambda+w'\rho)(H_{-w\gamma_\ell}) -1 = -\langle\lambda|w\gamma_\ell\rangle,
\end{align}
where we have used $a(u_\ell) = e^{\log(\phi_\ell) H_{-w\gamma_\ell}}$ to introduce the Cartan generator $H_{-w\gamma_\ell}$ of the $SL(2,\mathbb{A})$ associated with the $-w\gamma_\ell$ positive root space. We have also used
\begin{align}
w(\gamma_1+\ldots +\gamma_{\ell-1}) = w'\rho-\rho
\end{align}
and 
\begin{align}
(w'\rho)(H_{-w\gamma_\ell})= -\langle w'\rho | w\gamma_\ell\rangle = -\langle \rho | w_\ell\gamma_\ell\rangle = \langle\rho|\alpha_{i_\ell}\rangle =1\,,
\end{align}
since $\gamma_\ell=\alpha_{i_\ell}$ and we have normalised the symmetric bilinear form such that $\rho$ has unit inner product with all simple roots.  
The precise value of $\phi_\ell$ depends on whether one is at $Q_\infty=\reals\subset\mathbb{A}$ or at $Q_p\subset\mathbb{A}$ for $p<\infty$.

\subsection{Integral over \texorpdfstring{$\mathbb{R}$: $p=\infty$}{R: p=inf}}\label{pinf}
At the archimedean place we have result: 
\begin{lemma}[archimedean Gindikin--Karpelevich formula]
\beq
\int_{\mathbb{Q}_p} a(u_\ell)^{\lambda+\rho+w(\gamma_1+\ldots +\gamma_{\ell-1})} du_\ell=\sqrt{\pi}\frac{\Gamma(z_\ell/2)}{\Gamma\left((z_\ell+1)/2\right)}.
\eeq
\end{lemma}
\begin{proof}

If $u=u_\ell\in \mathbb{R}$, the Iwasawa decomposition (\ref{Iwaell}) is (see example~\ref{Iwex})
\begin{align}
\begin{pmatrix}
1&0\\u&1 
\end{pmatrix}
=\begin{pmatrix}
1&\frac{u}{1+u^2}\\0&1 
\end{pmatrix}
\begin{pmatrix}
1/\sqrt{1+u^2}&0\\0&\sqrt{1+u^2}
\end{pmatrix}
k
\end{align}
with
\begin{equation}
    k = \frac{1}{\sqrt{1+u^2}}
    \begin{pmatrix}
        1 & -u \\
        u & 1
    \end{pmatrix}
    \in SO(2, \reals) . 
\end{equation}

The diagonal matrix is $a(u_\ell)$. Substituting it into the integral (\ref{Iell}) for $u=u_\ell\in\mathbb{R}$, one obtains with (\ref{auell})
\begin{align}
\int_{\mathbb{R}} (1+u^2)^{-(z_\ell+1)/2} du = \sqrt{\pi}\frac{\Gamma(z_\ell/2)}{\Gamma\left((z_\ell+1)/2\right)}.
\end{align}
\end{proof}

\subsection{Integral over \texorpdfstring{$\mathbb{Q}_p$}{Qp} for finite \texorpdfstring{$p$}{p}}\label{pfin}
We now prove the corresponding result for finite primes:
\begin{lemma}[non-archimedean Gindikin--Karpelevich formula]
In~\eqref{Iell}, the non-archimedean integral evaluates to
\beq 
 \int_{\mathbb{Q}_p} a(u_\ell)^{\lambda+\rho+w(\gamma_1+\ldots \gamma_{\ell-1})} du_\ell=\frac{1-p^{-z_\ell-1}}{1-p^{-z_\ell}},
 \eeq
where $z_\ell>0$.
 \end{lemma}
 \begin{proof}
If $u=u_\ell\in \mathbb{Z}_p$, the matrix 
\begin{align}
\begin{pmatrix}
1&0\\u&1 
\end{pmatrix}
\end{align}
is in $SL(2,\mathbb{Z}_p)$ which is the compact part of $SL(2,\mathbb{Q}_p)$. Therefore $a(u_\ell)=1$ in this case and the integral is trivial. 

If $u\in \mathbb{Q}_p\backslash \mathbb{Z}_p$, the (non-unique) Iwasawa decomposition yields (see example~\ref{Iwex})
\begin{align}
\begin{pmatrix}
1&0\\u&1 
\end{pmatrix}
=\begin{pmatrix}
1&*\\0&1 
\end{pmatrix}
\begin{pmatrix}
u^{-1}&0\\0&u 
\end{pmatrix}
k\,.
\end{align}
Even though the Iwasawa decomposition is not unique, the norm of $a(u_\ell)$ is defined uniquely. We therefore obtain 
\begin{align}\label{GindiKarp}
\int_{\mathbb{Q}_p} |\phi_\ell|_p^{z_\ell+1} du &= \int_{\mathbb{Z}_p} dx + \int_{\mathbb{Q}_p\backslash \mathbb{Z}_p} |u|_p^{-z_\ell-1} du&\nonumber\\
&= 1+ \frac{p-1}{p} \frac{p^{-z_\ell}}{1-p^{-z_\ell}} = \frac{1-p^{-z_\ell-1}}{1-p^{-z_\ell}}\,,&
\end{align}
where we have used the integral~\eqref{usefulintegral} and the normalisation of the measure~\eqref{padicintmesnorm}. 
\end{proof}

\subsection{The global formula}

Putting finite and infinite contributions together, we therefore obtain (with $z_\ell$ as in (\ref{zell}))
\begin{align}
\label{Iellres}
I_\ell = \sqrt{\pi}\frac{\Gamma(z_\ell/2)}{\Gamma\left((z_\ell+1)/2\right)}\prod_{p<\infty} \frac{1-p^{-z_\ell-1}}{1-p^{-z_\ell}}
=\frac{\xi(z_\ell)}{\xi(z_\ell+1)}\,,
\end{align}
where we have used the Euler product formula~\eqref{RiemannProd} for the Riemann zeta function and its completion
\begin{align}
\xi(s) = \pi^{-s/2} \Gamma(s/2) \zeta(s)
\end{align}
that satisfies the functional relation (cf.~\eqref{CRZfuncrel})
\begin{align}
\xi(s)=\xi(1-s)\,.
\end{align}

\section{Assembling the constant term}\label{consttermass}

We can now write the final formula for the constant term (\ref{constterm}) by assembling (\ref{Cw}) and the result (\ref{Iellres}) inserted into the recursion relation (\ref{Iwrec}). The answer is 
\begin{align}
\int\limits_{N(\mathbb{Q})\backslash N(\mathbb{A})} E(\chi,ng) dn = \sum_{w\in \mathcal W} a^{w^{-1}\lambda+\rho} \prod_{\alpha>0\,|\, w\alpha<0} \frac{\xi(-\langle \lambda| w\alpha\rangle)}{\xi(1-\langle \lambda| w\alpha\rangle)}\,.
\end{align}
By relabelling the sum by $w\to w^{-1}$ we obtain the standard Langlands formula for the constant term in the minimal parabolic
\begin{align}
\label{LCF}
\int\limits_{N(\mathbb{Q})\backslash N(\mathbb{A})} E(\chi,ng) dn = \sum_{w\in \mathcal W} a^{w\lambda+\rho} \prod_{\alpha>0\,|\, w\alpha<0} \frac{\xi(\langle \lambda| \alpha\rangle)}{\xi(1+\langle \lambda| \alpha\rangle)}\,.
\end{align}
We note that the inner product here is normalised such that $\langle \rho|\alpha_i\rangle=1$ for all simple roots $\alpha_i$. Often one denotes the intertwining coefficient by
\begin{align}
\label{eq:intertwiner}
M(w,\lambda) =\prod_{\alpha>0\,|\, w\alpha<0} \frac{\xi(\langle \lambda| \alpha\rangle)}{\xi(1+\langle \lambda| \alpha\rangle)}\,.
\end{align}
This concludes the proof of theorem~\ref{LCFthm}.

\begin{remark}
The derivation above of the constant term formula made heavy use of the adeles $\mathbb{A}$. This was most noticeable when evaluating the integral (\ref{Iell}) that yielded the completed Riemann $\zeta$-functions in their Euler product form. Still, one may wonder whether this level of abstraction was really necessary. For $SL(2,\reals)$ one can obtain the constant terms alternatively by Poisson resummation techniques, see appendix~\ref{app:SL2Fourier} for a summary, without ever making reference to $p$-adic numbers. What this requires, however, is an explicit understanding of the sum over the cosets $B(\ints)\backslash G(\ints)$ in the definition of the Eisenstein series and their relation to sums over integer lattices. In the general case, the description of these cosets is not easy to characterise and the lattice sum descriptions typically involve representation theoretic constraints on the sum and this is discussed in more detail in section~\ref{sec:latticesums}. None of these details are required for obtaining the constant term formula when using the $p$-adic description and this is where the power of the method lies.
 \end{remark}

\section{Functional relations for Eisenstein series}
The definition of the Eisenstein series (\ref{eisdefA}) is initially restricted to the domain of (absolute) convergence of the defining sum. As we mentioned before this requires that the weight $\lambda$ entering in the definition of the character $\chi$ through the relation (\ref{charweight}) has sufficiently large real parts. More precisely, it is required to lie in the so-called Godemont range 
\begin{align}
\Re \langle\lambda|\alpha_i\rangle  >1\quad\quad\textrm{for all simple roots $\alpha_i$.} 
\end{align}
The Eisenstein series $E(\chi,g)$ can then be defined by analytic continuation in the complexified weight $\lambda$ to almost all values of $\lambda$. Note that a complex weight $\lambda$ corresponds to a character $\chi$ taking values not in $U(1)$ but in $\cx^\times$. We are mainly interested in real weights. As shown in~\cite{LanglandsEisenstein,LanglandsFE}, this continuation is possible everywhere except for certain hyperplanes that are related to the integral weight lattice.

An important property of the Eisenstein series is that they obey functional relations. More precisely one has:
\begin{theorem}[Functional relation for Eisenstein series \cite{LanglandsFE}]
 For each $w\in \mathcal{W}$ the Eisenstein series $E(\lambda, g)$ satisfies the functional relation:
  \begin{align}
\label{funrel}
E(\lambda,g) = M(w,\lambda) E(w\lambda,g)\,.
\end{align}
In other words, the Eisenstein series along the Weyl orbit of a character are all proportional to each other. 
\end{theorem}
\begin{proof}
To prove this, note first that the coefficients $M(w,\lambda)$ given by~\eqref{eq:intertwiner} satisfy the following property:
\begin{lemma}
\label{lem-functintertwine}
\begin{align}
\label{Mproduct}
M(w_1w_2,\lambda) = M(w_1,w_2\lambda) M(w_2,\lambda)\,, \qquad, \forall \, w_1, w_2\in \mathcal{W}.
\end{align} 
\end{lemma}

\begin{proof}
 Assume first that $w_1$ and $w_2$ are fundamental Weyl reflections, say $w_i$ and $w_j$. This yields 
\beqa
M(w_i, w_j\lambda)M(w_j, \lambda)&=&\prod_{\alpha>0\, \st\,  w_i\alpha<0} \frac{\xi(\left<w_j\lambda|\alpha\right>)}{\xi(1+\left<w_j\lambda|\alpha\right>)} \prod_{\alpha>0\, \st \,  w_j\alpha<0}\frac{\xi(\left<\lambda|\alpha\right>)}{\xi(1+\left<\lambda|\alpha\right>)}
\nn \\ 
&=& \prod_{w_j\alpha'>0\, \st\,  w_iw_j\alpha'<0} \frac{\xi(\left<\lambda|\alpha'\right>)}{\xi(1+\left<\lambda|\alpha'\right>)}\prod_{\alpha>0\, \st\,  w_j\alpha<0}\frac{\xi(\left<\lambda|\alpha\right>)}{\xi(1+\left<\lambda|\alpha\right>)},
\eqa
where in the second step we made the substitution $\alpha'=w_j\alpha$ and used the invariance of the bilinear form: $\left<w_j\lambda|\alpha\right>=\left<\lambda|w_j\alpha\right>$. 

We want to show that the two factors combine into the left hand side of~\eqref{Mproduct}.  To this end we rewrite the two disjoint sets of roots
    \begin{equation}
        \begin{split}
            A &= \{\alpha > 0 \land w_j \alpha < 0\} \\
            B &= \{w_j \alpha > 0 \land w_i w_j \alpha < 0\}
        \end{split}
    \end{equation}
    using Lemma 3.7 from \cite{Kac} which states that if any $\alpha \in \Delta_+$ satisfies $w_i \alpha < 0$ for some fundamental reflection $w_i$ then $\alpha = \alpha_i$; the root corresponding to $w_i$.

    Thus, $\alpha > 0 \land w_j \alpha < 0 \implies \alpha = \alpha_j \implies w_i w_j \alpha = - w_i \alpha_j < 0$
    and
    \begin{equation}
        A = \{\alpha > 0 \land w_i w_j \alpha < 0 \land w_j \alpha < 0\} \, .
    \end{equation}
    
    Similarly, $\alpha' = w_j \alpha > 0 \land w_i \alpha' = w_i w_j \alpha < 0 \implies \alpha' = \alpha_i \implies \alpha = w_j \alpha' = w_j \alpha_i > 0$ and
    \begin{equation}
        B = \{\alpha > 0 \land w_i w_j \alpha < 0 \land w_j \alpha > 0 \} \, .
    \end{equation}

    This gives the disjoint union
    \begin{equation}
        A \cup B = \{\alpha > 0 \land w_i w_j \alpha < 0 \} = \{\alpha > 0 \land w \alpha < 0 \}
    \end{equation}
    and it follows that 
\beqa
& &  \prod_{w_j\alpha'>0\, \st\,  w_iw_j\alpha'<0} \frac{\xi(\left<\lambda|\alpha'\right>)}{\xi(1+\left<\lambda|\alpha'\right>)}\prod_{\alpha>0\, \st\,  w_j\alpha<0}\frac{\xi(\left<\lambda|\alpha\right>)}{\xi(1+\left<\lambda|\alpha\right>)}
 \nn \\
 &=&\prod_{\alpha >0\, \st \, w_i w_j\alpha <0} \frac{\xi(\left<\lambda|\alpha\right>)}{\xi(1+\left<\lambda|\alpha\right>)}=M(w_i, w_j\lambda).
 \eqa
 The general formula~\eqref{Mproduct} for arbitrary products of fundamental reflections follows by iterating this procedure.
\end{proof}
The functional relation (\ref{funrel}) for the constant term of the Eisenstein series, now follows  from this result applied to the constant term formula (\ref{LCF}). The fact that this also extends to the non-constant terms was shown by Langlands~\cite{LanglandsFE}.
\end{proof}
\begin{remark}
The functional relation (\ref{funrel})  shows the limitations of the analytic continuation: For weights $\lambda$ and Weyl words $w\in \mathcal W$ for which $M(w,\lambda)$ is not a non-zero finite number, the relation appears ill-defined. This can only happen for $\lambda$ on certain hyperplanes as indicated above. Another apparent limitation of the functional relation is that for choosing a Weyl word $w$ that stabilises the weight $\lambda$ one would require $M(w,\lambda)$. This is not guaranteed to be true. For generic $\lambda$, (\ref{funrel}) provides a valid relation. The remaining cases are those when $E(\lambda,g)$ actually develops poles (as a function of $\lambda$) and one has to consider appropriate normalising factors in the functional relation to make it well-defined.
\end{remark}

\section{Expansion in maximal parabolics*}
\label{sec:CTFmax}

In the previous sections we have explained how to  expand Eisenstein series  along a minimal parabolic subgroup, with unipotent radical $N$. It is also possible to make an expansion along different parabolic subgroups with smaller unipotent radical. The analogue of the constant term (\ref{constterm}) then retains a dependence on the coordinates on $N(\mathbb{A})$ since only a subset of them is integrated out. In this section we state  and prove a theorem giving the formula for the constant terms of $E(\chi,g)$ in an expansion along a maximal parabolic subgroup.

The maximal parabolic subgroup, which we denote by $P_{j_\circ}$, is defined with respect to a particular choice of simple root $\alpha_{j_\circ}$ (i.e node $j_\circ$ in the Dynkin diagram) of $G$. This choice defines a subset of the simple roots, $\Pi_{j_\circ}=\Pi\backslash \{\alpha_{j_\circ}\}$, where $\Pi$ is the set of all simple roots. From this we furthermore define $\Gamma_{j_\circ}$ to be the set of all positive roots of $G$, which are given by linear combinations of the simple roots contained in $\Pi_{j_\circ}$, i.e., those roots of $G$ that do not contain $\alpha_{j_\circ}$. Using the standard Langlands decomposition, we can write a parabolic subgroup as the product of the Levi subgroup $L$ and the unipotent radical $U$. In the present case this reads $P_{j_\circ}=L_{j_\circ} U_{j_\circ}$. At the level of the corresponding Lie algebras one obtains $\mathfrak{p}_{j_\circ}=\mathfrak{l}_{j_\circ}\oplus \mathfrak{u}_{j_\circ}$ (not a direct sum of Lie algebras) with
\begin{align}
\mathfrak{p}_{j_\circ}=\mathfrak{a}\oplus\bigoplus_{\alpha\in\Delta_+\cup(-\Gamma_{j_\circ})}\mathfrak{g}_\alpha 
\end{align}
and
\begin{align}
\mathfrak{l}_{j_\circ}=\mathfrak{a}\oplus\bigoplus_{\alpha\in\Gamma_{j_\circ}\cup-\Gamma_{j_\circ}}\mathfrak{g}_\alpha\quad\text{and}\quad\mathfrak{u}_{j_\circ}=\bigoplus_{\substack{\alpha\in\Delta_+\backslash\Gamma_{j_\circ}}}\mathfrak{g}_\alpha\,,
\end{align}
respectively. Note that the case of the minimal parabolic above can be recovered by setting $\Gamma_{j_\circ}=\emptyset$.

Proceeding in analogy with the definitions of the minimal parabolic expansion we define the constant term part of the maximal parabolic expression as
\begin{align}
C_{j_\circ}=\int\limits_{U_{j_\circ}(\mathbb{Q})\backslash U_{j_\circ}(\mathbb{A})} E(\lambda,ug) du\,.
\end{align}
The subscript $j_\circ$ indicates the restriction to the maximal parabolic subgroup. A similar constant term formula as in the case of the minimal parabolic expansion can also be derived for this case. 
Upon deleting the $j_\circ$th node from the Dynkin diagram of $G$, we will denote the group associated with the Dynkin diagram which is left, by $G'$. We also note that the Levi factor $L_{j_\circ}$ of $P_{j_\circ}$ can then be written as $L_{j_\circ}=GL(1)\times G'$, where the one-parameter group $GL(1)$ is parametrised by a single variable in $\mathbb{R}^\times$.

Let us then state the formula for the constant term in this maximal parabolic:
\begin{theorem}[Constant terms for maximal parabolics \cite{Moeglin}]\label{consttermmaxparab} 
The constant term  of $E(\lambda,g)$ with respect to the unipotent radical $U_{j_\circ}\subset P_{j_\circ}$ is given by:
\begin{align}\label{parabexp}
\int\limits_{U_{j_\circ}(\mathbb{Q})\backslash U_{j_\circ}(\mathbb{A})} E(\lambda,ug) du = \sum_{w\in \Weyl_{j_\circ}\backslash\Weyl} e^{\langle\left(w\lambda+\rho\right)_{\parallel j_\circ}|H(g)\rangle}M(w,\lambda)E^{G'}\left(\chi_w,g\right)\,.
\end{align}
Below we provide a concise explanation of the notation used and give a proof of the formula.
\end{theorem}

In equation~\eqref{parabexp} the Weyl group of $P_{j_\circ}$ is denoted by $\Weyl_{j_\circ}$ and the sum on the right-hand-side is then restricted to a sum over a coset of the Weyl group, in contrast to the minimal parabolic case.
Furthermore the projection operators $(\cdot)_{\parallel j_\circ}$ and $(\cdot)_{\perp j_\circ}$ are defined as follows when acting on a weight $\lambda\in \mathfrak{a}^{*}$:
\begin{subequations}
\begin{align}
\lambda_{\parallel j_\circ}&:=\frac{\left<\Lambda_{j_\circ}|\lambda\right>}{\left<\Lambda_{j_\circ}|\Lambda_{j_\circ}\right>} \Lambda_{j_\circ}, &\\
\lambda_{\perp j_\circ}&:=\lambda - (\lambda)_{\parallel j_\circ}.&
\end{align}
\end{subequations}
These correspond, respectively, to the component of $\lambda$ parallel and orthogonal to the fundamental weight $\Lambda_{j_\circ}$. The orthogonal component is given by a linear combination of simple roots of $G'$. The character $\chi_w$ follows the definition~\eqref{charweight}, however with $\lambda$ now replaced by the weight $\left(w\lambda\right)_{\perp j_\circ}$. The $G'$ invariant Eisenstein Series on the right-hand-side of the equation is independent of the $GL(1)$ factor in the decomposition of the Levi subgroup $L_{j_*}$, as this dependence is projected out using the $(\cdot)_{\perp j_\circ}$ operator and appears solely through the exponential prefactor. Let us also note that for simplicity of notation we have put $g$ in the argument of the Eisenstein series on the right, even though $g$ lies effectively in $G'$.

We also note that the formula (\ref{parabexp}) is well-defined and independent of the choice of coset representative due to the functional relation (\ref{funrel}). 

\begin{proof} \mbox{}
Sources for the analysis presented here are~\cite{Moeglin,GRS}.\index{constant term!in maximal parabolic}
Consider two parabolic subgroups $P_1(\mathbb{A})$ and $P_2(\mathbb{A})$  of $G(\mathbb{A})$. The first one we take to be the one defining an Eisenstein series through
\begin{align}
E(\chi,g) = \sum_{\gamma\in P_1(\mathbb{Q})\backslash G(\mathbb{Q})} \chi(\gamma g),
\end{align}
where $\chi: P_1(\rats)\bs P_1(\mathbb{A}) \to \cx^\times$ is a character on the parabolic $P_1(\mathbb{A})$. The parabolic subgroup $P_2(\mathbb{A})=L_2(\mathbb{A}) U_2(\mathbb{A})$ is used to define the constant term along $U_2(\mathbb{A})$ via
\begin{align}
\label{CTpar}
C_{U_2}(\chi,g) = \lint_{U_2(\mathbb{Q})\backslash U_2(\mathbb{A})} E(\chi,ug)du.
\end{align}
By the definition of the integral, the result is determined by its dependence on $g=l\in L_2(\mathbb{A})$ and we will restrict to the Levi factor now. Most of the steps in the evaluation of (\ref{CTpar}) will be very similar to those in section~\ref{sec:Bruhat}.

Using the Bruhat decomposition, we can rewrite the integral (\ref{CTpar}) as
\begin{align}
C_{U_2}(\chi,l) = \sum_{w\in \Weyl_1 \backslash \Weyl / \Weyl_2} C_{w,U_2}(\chi,l)\,,
\end{align}
where the individual of the double coset of the Weyl group $\Weyl$ is
\begin{align}
C_{w,U_2}(\chi,l)= \sum_{\delta \in w^{-1} P_1(\rats) w \cap P_2(\rats)\backslash P_2(\rats)}\;
\lint_{U_2(\mathbb{Q})\backslash U_2(\mathbb{A})} \chi(w\delta u l)  du.
\end{align}
Now the sum over $\delta\in P_2(\rats)$ can be split into the Levi and unipotent part according to $\delta = \gamma_l \gamma_u$ and then one can unfold the sum over $\gamma_u$ onto the integration domain as in (\ref{Bruhatapplied}). The result is
\begin{align}
C_{w,U_2}(\chi,l) = \sum_{\gamma_l \in w^{-1} P_1(\rats) w \cap L_2(\rats)\backslash L_2(\rats)}\; 
\lint_{w^{-1} P_1(\rats) w\cap U_2(\mathbb{Q})\backslash U_2(\mathbb{A})} \chi(w\gamma_l u l)du.
\end{align}

We now specialize to $P_2(\mathbb{A})$ being maximal parabolic. Then the Levi factor takes the form (cf.~\ref{maxparLevi})
\begin{align}
L_2(\mathbb{A}) = GL(1,\mathbb{A}) \times M_2(\mathbb{A})
\end{align}
with $M_2$ reductive and we parametrise the $L_2(\mathbb{A})$ element as $l=rm$. We next separate out the dependence on the $GL(1,\mathbb{A})$ element $r$ by moving it to the left within $\chi$. This leads to
\begin{align}
\label{CTindpar}
C_{w,U_2}(\chi,rm) = r^{w^{-1}\lambda+\rho}\sum_{\gamma_m \in w^{-1} P_1(\rats) w \cap M_2(\rats)\backslash M_2(\rats)}\; 
\lint_{w^{-1} P_1(\rats) w\cap U_2(\mathbb{Q})\backslash U_2(\mathbb{A})} \chi(w u \gamma_m  m)du,
\end{align}
by combining the contribution from $\chi(wrw^{-1})$ and the change of the measure $du$. Note also that we changed the summation over the $L_2(\rats)$ cosets to one over $M_2(\rats)$ cosets since the two agree. We have also interchanged $\gamma_m$ and $u$ as the corresponding change of variables is uni-modular ($\gamma_m\in M_2(\rats)$ is discrete).

Let us analyse the properties of the integral
\begin{align}
\label{indchar}
\mathcal{I}=\lint_{U_2^w(\mathbb{Q})\backslash U_2(\mathbb{A})} \chi(w u \gamma_m  m)du
\end{align}
that is a function from $M_2(\mathbb{A})\to \cx$. We also defined for simplicity
\begin{align}
M_2^w(\rats) := w^{-1}P_1(\rats)w\cap M_2(\rats)
\quad\textrm{and} \quad U_2^w(\rats) := w^{-1}P_1(\rats)w\cap U_2(\rats)
\end{align}
and we note that $M_2^w(\rats)$ is a parabolic subgroup of $M_2(\rats)$.
The integral (\ref{indchar}) is invariant under $\eps\in M_2^w(\rats)$ according to
\begin{align}
\lint_{U_2^w(\mathbb{Q})\backslash U_2(\mathbb{A})} \chi(w u\eps  \gamma_m  m)du
= \lint_{U_2^w(\mathbb{Q})\backslash U_2(\mathbb{A})} \chi(w \eps u \gamma_m  m)du =\mathcal{I}
\end{align}
since $\eps=w^{-1} p_1 w$ for some $p_1\in P_1(\rats)$ and the character $\chi: P_1(\rats)\bs P_1(\mathbb{A})\to \cx^\times$ is invariant under $P_1(\rats)$.  But this means that
\begin{align}
\label{inducedcharacter}
\chi^w_2( \gamma_m m) = \lint_{U_2^w(\mathbb{Q})\backslash U_2(\mathbb{A})}  \chi(w u \gamma_m  m)du
\end{align}
is a character $\chi^w_2 : M_2^w(\rats)\bs M_2(\mathbb{A}) \to \cx^\times$ if it is non-zero.
(The integral $\mathcal{I}$ serves as an intertwiner \index{intertwiner} from characters on $P_1(\mathbb{A})$ to characters on $M_2^w(\mathbb{A})$.)

Inserting this back into the individual constant term (\ref{CTindpar}) we obtain
\begin{align}
C_{w,U_2}(\chi,rm) = r^{w^{-1}\lambda+\rho}\sum_{\gamma_m \in M_w^2(\rats)\backslash M_2(\rats)}  \chi_2^w(\gamma_m m).
\end{align}
The sum over $\gamma_m$ now produces an Eisenstein series on $M_2(\mathbb{A})$ so that in all
\begin{align}
C_{U_2}(\chi,rm) = \sum_{w\in \Weyl_1\bs \Weyl/\Weyl_2} r^{w^{-1}\lambda+\rho} E^{M_2}(\chi^w_2,m),
\end{align}
where we indicated that the Eisenstein series is on $M_2(\mathbb{A})$. This expression can be simplified a bit more by identifying the character $\chi_2^w$ in terms of a weight of $\mf{m}_2=\mathrm{Lie}(M_2)$. To this end we evaluate (\ref{inducedcharacter}) at a semi-simple element of $M_2(\mathbb{A})$, i.e., $m=a$. This leads to
\begin{align}
\chi_2^w (a) &= \lint_{U_2^w(\rats)\bs U_2(\mathbb{A})} \chi(wua) du 
= \chi(waw^{-1}) \delta_{\bar{U}_2^w}(a)  \lint_{U_2^w(\rats)\bs U_2(\mathbb{A})} \chi(wu) du\nn\\
&=M(w^{-1},\lambda) a^{(w^{-1}\lambda+\rho)_{M_2}},
\end{align}
where the last symbol denotes the orthogonal projection onto the space of $M_2$ weights. The exponent comes about as follows. The character $\chi(waw^{-1})$ evaluates to $a^{w^{-1}(\lambda+\rho)}$ and the modulus character on $\bar{U}_2^w$ is determined by the sum over all roots of $U_2$ that are not mapped to roots of $P_1$, i.e., the total exponent of $a$ 
\begin{align}
w^{-1}(\lambda+\rho) + \sum\limits_{\alpha\in\Delta(\mf{u_2})\st  w\alpha\notin \Delta(\mf{p}_1)}\alpha = (w^{-1}\lambda + \rho)_{M_2}.
\end{align}
We have made furthermore made use of our knowledge of the Gindikin--Karpelevich type integral, cf. the evaluation of~(\ref{Iw}). In summary we arrive at
\begin{align}
C_{U_2}(\chi,rm) = \sum_{w\in \Weyl_1\bs \Weyl/\Weyl_2} r^{w^{-1}\lambda+\rho} M(w^{-1},\lambda) E^{M_2}((w^{-1}\lambda)_{M_2},m),
\end{align}
in agreement with (\ref{parabexp}) if one replaces $w$ by $w^{-1}$ which maps the double coset to $ \Weyl_2\bs \Weyl/\Weyl_1$.
\end{proof}

\chapter[Whittaker coefficients of Eisenstein series]{Whittaker coefficients of\\ Eisenstein series}
\label{ch:Whittaker-Eisenstein}

In this chapter, we derive theorem~\ref{CasselmShalikaform} that states the formula of \index{Casselman--Shalika formula} Casselman--Shalika~\cite{CasselmanShalika} (see also~\cite{ShahidiBook}) for the \emp{local} abelian Fourier coefficients in the minimal parabolic (Borel) subgroup $B(\mathbb{A})\subset G(\mathbb{A})$ for the Eisenstein series $E(\chi, g)$. This formula is used to evaluate Fourier integrals with a generic character $\psi$. By the discussion in chapter~\ref{ch:fourier} the global form of these Fourier coefficients are captured by the \index{Whittaker coefficient!spherical} spherical Whittaker coefficient
\begin{align}
\label{Whittakershort}
W^{\circ}_{\psi}(\chi,g)&=\int\limits_{N(\mathbb{Q})\backslash N(\mathbb{A})}E(\chi,ng)\overline{\psi(n)}dn
\end{align}
for a (quasi-)character $\chi: \,B(\mathbb{A})\to \cx^\times$ and a general unitary character $\psi: N(\mathbb{Q})\backslash N(\mathbb{A})\to U(1)$. For $SL(2,\mathbb{A})$ we have already evaluated this integral in section~\ref{sec:SL2FC}. As there, a useful strategy is to factorise the integral and perform it at all places separately. It will turn out that only for the finite primes $p<\infty$ and so-called generic and \emphindex[character!unramified]{unramified characters} $\psi$ (to be defined below) a nice and compact formula exists. In sections~\ref{sec:genpsi} and~\ref{sec:degpsi}, we will explain how to also evaluate (\ref{Whittakershort}) for arbitrary generic or even degenerate characters $\psi$.

\section[Reduction of the integral and the longest Weyl word]{Reduction of the integral and the longest\\ Weyl word}
\label{sec:redlong}

To begin with, we bring the integral (\ref{Whittakershort}) into a form that is more amenable to evaluation.
As discussed in section~\ref{sec_sphericalvector}, the spherical Whittaker coefficient satisfies
\begin{align}
\label{Wtransform}
W^{\circ}_\psi(\chi,ngk) = \psi(n) W^{\circ}_\psi(\chi,g)
\end{align}
and is therefore determined by its values on $A(\mathbb{A})$ due to the Iwasawa decomposition (\ref{IwasawaR}) and (\ref{IwasawaQp}). Hence, we will only consider it for $g=a\in A(\mathbb{A})$ in the sequel. For the discussion in this subsection we assume $\psi$ to be generic (see definition~\ref{def:GC}), i.e., it does not vanish on any simple root generator.

We start evaluating (\ref{Whittakershort}) by applying the Bruhat decomposition as for the constant term to obtain
\begin{align}
\label{BruhatW}
W^{\circ}_\psi(\chi,a)&=\sum\limits_{\gamma\in B(\mathbb{Q})\backslash G(\mathbb{Q})}\;\int\limits_{N(\mathbb{Q})\backslash N(\mathbb{A})}\chi(\gamma na)\overline{\psi(n)}dn \nn\\
&=\sum_{w\in\Weyl}\int\limits_{w^{-1}B(\mathbb{Q})w\cap N(\mathbb{Q})\backslash N(\mathbb{A})}\chi(wna)\overline{\psi(n)}dn\,.
\end{align}
{}From the last line let us define $W^{\circ}_\psi(\chi,a)=\sum_{w\in\Weyl}F_{w,\psi}$, where
\begin{align}\label{Fw}
F_{w,\psi}=\int\limits_{w^{-1}B(\mathbb{Q})w\cap N(\mathbb{Q})\backslash N(\mathbb{A})}\chi(wna)\overline{\psi(n)}dn
\end{align}
is the contribution from the Weyl word $w$.

Let us start with an analysis of the integration range of the Fourier integral~\eqref{Fw}, given by the coset $w^{-1}B(\mathbb{Q})w\cap N(\mathbb{Q})\backslash N(\mathbb{A})$, and the corresponding contribution to $W_\psi^\circ(\chi,a)$.
It is clear that the intersection in the denominator of this coset consists of those elements of the (upper) Borel subgroup that are mapped to upper Borel elements under the action of the Weyl element $w$.  For the whole denominator we can therefore write
\begin{align}
w^{-1}B(\mathbb{Q})w\cap N(\mathbb{Q}) =  \prod\limits_{\alpha>0\atop w\alpha>0}N_\alpha(\mathbb{Q})\,.
\end{align}
With this, the integration range conveniently splits up in the following way
\begin{align}
w^{-1}B(\mathbb{Q})w\cap N(\mathbb{Q})\backslash N(\mathbb{A})\simeq\left(\prod\limits_{\beta>0\atop w\beta>0}N_\beta(\mathbb{Q})\backslash N_\beta(\mathbb{A})\right)\cdot\left(\prod\limits_{\gamma>0\atop w\gamma<0} N_\gamma(\mathbb{A})\right)\,.
\end{align}
Let us introduce some notation. We denote the union in the first parenthesis as
\begin{align}\label{NBeta}
N^w_{\{\beta\}}:=\left(\prod\limits_{\beta>0\atop w\beta>0}N_\beta(\mathbb{Q})\backslash N_\beta(\mathbb{A})\right)
\end{align}
and the union in the second parenthesis as
\begin{align}\label{NGamma}
N^w_{\{\gamma\}}:=\left(\prod\limits_{\gamma>0\atop w\gamma<0} N_\gamma(\mathbb{A})\right)\,.
\end{align}
Here, the sets of roots $\{\beta\}$ and $\{\gamma\}$ contain precisely those roots, which satisfy the conditions imposed on the products in~\eqref{NBeta} and~\eqref{NGamma}, respectively. It is important to note that there is a qualitative difference in the two sets: $N^w_{\{\gamma\}}$ is non-compact while $N^w_{\{\beta\}}$ is \emp{compact}. 
With this splitting of the integration range the contribution $F_w$ then takes the following form
\begin{align}
F_{w,\psi}&=\int\limits_{N^w_{\{\beta\}} N^w_{\{\gamma\}}}\chi(wna)\overline{\psi(n)}dn\nn\\
&=\int\limits_{N^w_{\{\beta\}}}\int\limits_{N^w_{\{\gamma\}}}\chi(wn_\beta n_\gamma a)\overline{\psi(n_\beta n_\gamma)}dn_\beta dn_\gamma\,.
\end{align}
Inserting $w^{-1}w$ between $n_\beta$ and $n_\gamma$ and splitting the character up into two factors, we obtain
\begin{align}
\int\limits_{N^w_{\{\beta\}}}\int\limits_{N^w_{\{\gamma\}}}\chi(wn_\beta w^{-1}w n_\gamma a)\overline{\psi(n_\beta)}\,\overline{\psi( n_\gamma)}dn_\beta dn_\gamma\,.
\end{align}
Let us note that the character $\chi$ is left invariant under any subgroup that is given by the exponential of positive root generators. In particular this applies to elements $n$ of the (upper) Borel subgroup. Hence by definition of $N^w_{\{\beta\}}$, the character $\chi$ is insensitive to the factor $wn_\beta w^{-1}$ in the argument and we can split-off the integral over $n_\beta$, leaving us with
\begin{align}\label{FwSplit}
F_{w,\psi}=\int\limits_{N^w_{\{\beta\}}}\overline{\psi(n_\beta)}dn_\beta\int\limits_{N^w_{\{\gamma\}}}\chi(w n_\gamma a)\overline{\psi(n_\gamma)}dn_\gamma\,.
\end{align}

Given the form of~\eqref{FwSplit}, we see that in the integral over $n_\beta$, effectively a periodic function is integrated over a full period in the compact space $N^w_{\{\beta\}}$. Provided that the character $\psi$ is non-trivial along at least one simple root contained in $\{\beta\}$, this means that the whole integral will vanish. Since this is always true for generic $\psi$ and we arrive at the conclusion
\begin{align}
F_{w,\psi} = 0 \quad \textrm{unless $w= \wlong$.}
\end{align}
This is true since all Weyl transformations except for the longest Weyl word $\wlong$ leave at least one simple root positive. Therefore we arrive at the following expression for the Whittaker function for a generic $\psi$:
\begin{align}
\label{Wwlong}
W^{\circ}_\psi(\chi,a) = \lint_{N(\mathbb{A})} \chi(\wlong n a) \overline{\psi(n)} dn.
\end{align}
We will also refer to this expression as the~\emp{Jacquet integral}, see~\cite{HJacquet}.

\section{Unramified local Whittaker functions}
\label{sec:unramified-local-Whittaker-coefficients}

The integral (\ref{Wwlong}) should now be evaluated for all places separately, that is
\beqa
 \lint_{N(\mathbb{A})} \chi(\wlong n a) \overline{\psi(n)} dn& =&\int_{N(\mathbb{R})} \chi_\infty(\wlong na)\overline{\psi_\infty(n)}dn \times \prod_{p<\infty} \int_{N(\mathbb{Q}_p)} \chi_p(\wlong na)\overline{\psi_p(n)}dn
 \nn \\
 & =& W^{\circ}_{\psi_\infty} \times \prod_{p<\infty} W^{\circ}_{\psi_p}.
 \label{eulerWhittaker}
 \eqa
However, while we will derive a nice closed formula for the local places $p<\infty$, a general expression for the real place is not known to the best of our knowledge. In the case $SL(2,\mathbb{A})$, the resulting expression was given by the Bessel function~\eqref{SL2Whittinf} and for $SL(3,\mathbb{A})$ it is known that the triple integral over the three unipotent generators gives convoluted integrals of Bessel functions~\cite{VT,BumpSL3}, see also section~\ref{sec:SL3ex}. For general groups with `more non-abelian' unipotent subgroups a proliferation of this nested structure of special functions is to be expected. On the other hand, the Whittaker functions for finite $p<\infty$ contain the essential number theoretic information that is reflected as instanton measures in string theory applications~\cite{Green:1997tv,Moore:1998et}. Therefore we will from now on consider only the group $G(\rats_p)$ for $p<\infty$ and calculate the local Whittaker functions. To ease notation we shall suppress all subscripts involving primes, hence in the remainder of this section we write $\psi$ for $\psi_p$ and $\chi$ for  $\chi_p$.

\subsection{Unramified characters \texorpdfstring{$\psi$}{psi}} 

The formula of Casselman and Shalika for the local Whittaker functions is most conveniently stated when one restricts the character $\psi$ to be \emp{unramified} (see definition \ref{def_unram} and~\cite[p.~219]{CasselmanShalika}). Recall from definition \ref{def_unram} that this means that the character $\psi:N(\rats_p)\to U(1)$ has what a physicist might call `unit instanton charges', i.e., when the element $n$ is expanded in terms of positive step operators in a Chevalley basis as
\begin{align}
    \label{eq:n-param}
    n = \bigg(\prod_{\alpha\in \Delta_+\setminus\Pi} x_{\alpha}(u_\alpha)\bigg)\bigg(\prod_{\alpha\in\Pi}x_\alpha(u_\alpha)\bigg)\in N(\rats_p),
\end{align}
where we have ordered the individual factors in a convenient way. Note that for evaluating $\psi(n)$ the order does not matter since $\psi$ is a homomorphism between abelian groups.
An unramified character $\psi$ is then one that satisfies
\begin{align}
    \label{eq:psi-unramified}
\psi(n) = \exp\left( 2\pi i \lb \sum_{\alpha\in\Pi} m_\alpha u_\alpha\rb\right)
\end{align}
with $m_\alpha =1$ for all simple roots $\alpha\in\Pi$. An unramified character is automatically generic.

The local Whittaker function for an unramified vector will be denoted simply by 
\begin{align}
\label{Wunramified}
W^{\circ}(\chi,a) = \lint_{N(\rats_p)} \chi(\wlong n a ) \overline{\psi(n)} dn,
\end{align}
where the reference to $\psi$ has been suppressed for notational convenience and we do not display the fact that we are using a fixed $p<\infty$. Standard manipulations similar to~\eqref{constaexplicit} on (\ref{Wunramified}) lead to 
\begin{align}
\label{CStoUS}
W^{\circ}(\chi,a) &=  \chi(\wlong a \wlong^{-1} ) \lint_{N(\rats_p)}\chi(\wlong a^{-1} n a) \overline{\psi(n)} dn\nn\\
&= \chi(\wlong a \wlong^{-1})\delta(a) \lint_{N(\rats_p)}\chi(\wlong n) \overline{\psi(ana^{-1})} dn\nn\\
&= |a^{\wlong\lambda+\rho}|\lint_{N(\rats_p)}\chi(\wlong n) \overline{\psi^a(n)} dn,
\end{align}
where we defined $\psi^a(n) := \psi(ana^{-1})$.

\subsection{Vanishing properties}
\label{VP}

One advantage of restricting to unramified characters is that it is very simple to determine the support of $W^\circ(\chi,a)$ (cf. also~\cite[Lemma~5.1]{CasselmanShalika}).  Consider an element $n\in N(\ints_p)\subset K_{p}$; then, by right $K_{p}$-invariance and the transformation properties (\ref{Wtransform}),
\begin{align}
    \label{eq:urnamified-vanishing-condition}
W^{\circ}(\chi,a) = W^{\circ}(\chi, an) = W^{\circ}(\chi,ana^{-1}a) = \psi(ana^{-1}) W^{\circ}(\chi,a).
\end{align}
Therefore, $W^{\circ}(\chi,a)$ can only be non-vanishing if $\psi(ana^{-1})=1$ which requires that $a^\alpha\in \ints_p$ for all positive roots $\alpha$ since $\psi$ is unramified. 

\section{The Casselman--Shalika formula}
\label{sec_CSformula}

We are now ready to state the main theorem of this section: 

\begin{theorem}[Casselman--Shalika formula]\label{CasselmShalikaform}
    \index{Casselman--Shalika formula|textbf}
The local unramified Whittaker function $W^{\circ}(\chi, a)$, defined by the integral (\ref{eulerWhittaker}) for each $p<\infty$, is given by 
\begin{align}
\int_{N(\mathbb{Q}_p)} \chi(\wlong na)\overline{\psi(n)}dn=\frac{\eps(\lambda)}{\zeta(\lambda)} \sum_{w\in \Weyl} (\det(w)) |a^{w\lambda+\rho}|  \prod_{\alpha>0  \atop w\alpha<0} p^{\langle\lambda|\alpha\rangle}
\end{align}
\end{theorem}

\begin{proof} The proof of this theorem will constitute the remainder of this section \ref{sec_CSformula}. \end{proof}

When translated to our notation, the Casselman--Shalika formula found for $W^{\circ}(\chi,a)$ in~\cite[Thm.~5.4]{CasselmanShalika} takes the form 
\begin{align}
    \label{CSus}
    W^{\circ}(\chi,a) = 
        \frac{\eps(\lambda)}{\zeta(\lambda)} \sum_{w\in \Weyl} (\det(w)) \bigg(\prod_{\alpha>0  \atop w\alpha<0} p^{\langle\lambda|\alpha\rangle}\bigg) |a^{w\lambda+\rho}| =
        \frac1{\zeta(\lambda)} \sum_{w\in \Weyl} \eps(w\lambda) |a^{w\lambda+\rho}|
\end{align}
with 
\begin{align}
    \label{zetaeps}
    \zeta(\lambda) = \prod_{\alpha>0} \frac{1}{1-p^{-(\langle\lambda|\alpha\rangle+1)}},\quad
    \eps(\lambda) = \prod_{\alpha>0} \frac{1}{1-p^{\langle\lambda|\alpha\rangle}}.
\end{align}

The latter identity of \eqref{CSus} follows from
\begin{equation}
    \epsilon(w_i \lambda) = \prod_{\alpha > 0} \frac{1}{1 - p^{\langle \lambda | w_i \alpha \rangle}} = \frac{1}{1 - p^{-\langle \lambda | \alpha_i \rangle}} \prod_{\substack{ \alpha > 0 \\ \alpha \neq \alpha_i }} \frac{1}{1 - p^{\langle \lambda | \alpha \rangle}} = -p^{\langle \lambda | \alpha_i\rangle} \epsilon(\lambda)
\end{equation}
where $w_i$ is a fundamental reflection switching the sign of $\alpha_i$ and permuting the remaining positive roots. Recall that $\det w = (-1)^{\ell(w)}$ where $\ell(w)$ is the length of $w$ as introduced in section~\ref{sec:simple-lie-alg}.
Formula \eqref{CSus} is valid only for unramified $\psi$ and we have used $\chi$ and $\lambda$ interchangeably.

Our strategy for proving theorem~\ref{CasselmShalikaform} will be a mixture of the works of Jacquet~\cite{HJacquet} and Casselman--Shalika~\cite{CasselmanShalika}. The argument consists of the following steps:
\begin{enumerate}
\item Derivation of a functional equation for the Whittaker function under Weyl transformations on $\chi$
\item Use this to show that a suitable multiple of the Whittaker function is Weyl invariant and write it as a sum over Weyl images
\item Determine one term in this sum and derive all other terms from it. This will yield formula (\ref{CSus})
\end{enumerate}
Finally we will also show in section~\ref{sec:genpsi} how the formula (\ref{CSus}) can be used to derive the Whittaker functions for all generic characters $\psi$. 

However, as a preparatory `step 0', we first recall and slightly extend some results from chapter~\ref{ch:SL2-fourier} where the Fourier coefficients for Eisenstein series on $SL(2,\mathbb{A})$ were discussed. Namely, after equation~(\ref{SL2Floc}) we derived the Whittaker function at finite places, evaluated at the identity $a=1\in A(\rats_p)$,
for general $\psi$. Repeating the same steps but $(i)$ keeping $a$ arbitrary and $(ii)$ choosing an unramified character ($m=1$) leads to 
\begin{align}
\label{SL2Whitt}
G(\rats_p)=SL(2,\rats_p):\quad
F_{\wlong,\psi,p}&= \lint_{N(\rats_p)} \chi(\wlong n a) \overline{\psi(n)}d n &\nn\\
&= \chi(\wlong a\wlong^{-1})\delta(a)\lint_{N(\rats_p)} \chi(\wlong n) \overline{\psi(ana^{-1})}dn&\nn\\
&= \gamma_p(v^2) (1-p^{-2s}) \frac{|v|^{-2s+2} -p^{-2s+1}|v|^{2s}}{1-p^{-2s+1}}
\end{align}
with $\chi(a) =|a|^{2s}$ and $a=\mathrm{diag} (v,v^{-1})=v^{H_\alpha}$ in terms of the unique positive root $\alpha$ of $\mf{sl}(2,\reals)$. This formula, after dividing by $(1-p^{-2s})$ exhibits invariance under the Weyl reflection $s\leftrightarrow 1-s$. We will see how this feature generalises to arbitrary $G$ and why it is basically a consequence of this $SL(2,\rats_p)$ calculation. Equation (\ref{SL2Whitt}) also manifestly exhibits the vanishing property of section~\ref{VP} since the $p$-adic Gaussian vanishes unless $|v^2| = |a^\alpha| \leq 1$.

Before embarking on the proof proper, we also record the following
\begin{proposition}[Holomorphy of local Whittaker functions~\cite{CasselmanShalika}]
\index{Whittaker function!holomorphy}
\label{prop:hol}
The local Whittaker function $W^\circ(\chi,a)$ depends holomorphically on the quasi-character $\chi$. 
\end{proposition}
\begin{proof} Inspection of formula~\eqref{CSus} immediately reveals holomorphy when $\chi$ is in the Godement domain~\eqref{absoluteconvergence} of absolute convergence. This extends to all $\chi$ by virtue of the functional relation derived below.
\end{proof} 
The holomorphy of the Whittaker function in the case of $SL(2,\rats_p)$ (as a function of $s$) can also be seen from the explicit expression~\eqref{SL2Whitt} above. For $s\to \frac12$, the expression stays finite. We will comment in much more detail on the behaviour of Eisenstein series in chapter~\ref{ch:working}.

\subsection{Functional relation for the local Whittaker function}

We follow Jacquet's thesis~\cite{HJacquet}. First one defines a function associated to the Whittaker function by
\begin{align}
\label{FWrel}
F(\lambda,g) = W^{\circ}_{\psi}(\lambda,\wlong^{-1} g)
\end{align}
for $g\in G(\rats_p)$. This leads to the integral expression
\begin{align}
\label{auxfn1}
F (\lambda,g) = \lint_{N_-(\rats_p)} \chi(n_- g) \overline{\psi_-(n_-)} d n_-
\end{align}
for the associated function. Here, objects with a minus subscript refer to the unipotent opposite to the standard unipotent $N(\rats_p)$. In other words, $N_-(\rats_p)$ designates the subgroup of $G(\rats_p)$ generated by the exponentials of the \emp{negative} roots, whereas the usual $N(\rats_p)$ is associated with the positive roots. The reason that the opposite group arises here is because $\wlong$ maps all positive roots to negative ones (possibly combined with an outer automorphism). We will derive a functional relation for $F$ under Weyl transformations which by (\ref{FWrel}) will imply one for the Whittaker function.

The method for deriving the functional relation will be by reducing to the functional relation for $SL(2,\rats_p)$ that is manifest in (\ref{SL2Whitt}) and then using the fact that $G(\rats_p)$ is made up of $SL(2,\rats_p)$ subgroups.

Let $\alpha_i$ be a simple positive root of $G(\rats_p)$. Then define for $g\in\rats_p$
\begin{align}
\label{auxfn2}
F_{i}(\lambda,g) = \lint_{N_{i,-}(\rats_p)} \chi(n_{i,-} g) \overline{\psi_{i,-}(n_{i,-})} d n_{i,-},
\end{align}
where the integral is now only over the one-dimensional subgroup generated by $x_{-\alpha_i}(u)$ and similarly the character $\psi_{i,-}$ is one of the (lower) unipotent of the $SL(2,\rats_p)$ associated with $\alpha_i$ and can be obtained from $\psi_-$ by restriction to the subgroup $N_{i,-}$. The function $F_{i}$ is useful because for any $\alpha_i$ we can write
\begin{align}
N_-(\rats_p) = N_{i,-}(\rats_p) \hat{N}_-(\rats_p),
\end{align}
where $\hat{N}_-(\rats_p)$ are the lower unipotent elements that are not of the form $x_{-\alpha_i}(u)$ for some $u\in \rats_p$. Associated with the factorisation above is a unique decomposition $n_- = n_{i,-} \hat{n}_-$ and then the integral (\ref{auxfn1}) leads to
\begin{align}
\label{aux1aux2}
F (\lambda,g) = \lint_{\hat{N}_-(\rats_p)} F_{i}(\hat{n}_-g) \overline{\psi_-(\hat{n}_-)} d\hat{n}_-
\end{align}
by carrying out the integral over $dn_{i,-}$.

The $SL(2,\rats_p)$ projected function (\ref{auxfn2}) has the following invariances
\begin{align}
F_{i}(\lambda, \hat{n}gk) = F_{i}(\lambda,g) \quad\textrm{for $\hat{n}\in\hat{N}(\rats_p)$ and $k\in K_{p}=G(\ints_p)$,}
\end{align}
where $\hat{N}(\rats_p)$ is the unipotent subgroup opposite to $\hat{N}_-(\rats_p)$. It is generated by all positive roots but $\alpha_i$. The set of these roots is invariant under the Weyl reflection $w_i$. Let $P_i$ be the next-to-minimal parabolic subgroup defined by the (non-unique) decomposition
\begin{align}
G(\rats_p) = P_i(\rats_p) K_{p} = \hat{N}(\rats_p) L_i(\rats_p) K_{p} = \hat{N}(\rats_p) \hat{A}(\rats_p) SL(2,\rats_p)_{\alpha_i} K_{p}
\end{align}
with $\hat{A}(\rats_p)$ the part of the split torus $A(\rats_p)$ that is not contained in the torus of the embedded $SL(2,\rats_p)_{\alpha_i}$. Using this decomposition and the invariances of $F_i$, one finds that the function $F_{i}(g)$ is determined by its values on elements of the form $g=\hat{a}g_i$ with $\hat{a}\in\hat{A}(\rats_p)$ and $g_i\in SL(2,\rats_p)_{\alpha_i}$. On such values one has that
\begin{align}
\label{Firest}
F_{i} (\lambda,\hat{a}g_i) =| \hat{a}^{\lambda+\rho-\alpha_i}| \lint_{N_{i,-}(\rats_p)} \chi(n_{i,-} g_i) \overline{\psi_{i,-}^{\hat{a}}(n_{i,-})} d n_{i,-}
\end{align}
with $\psi_-^{\hat{a}}( n_{i,-}) = \psi_-(\hat{a} n_{i,-} \hat{a}^{-1})$. The integral is basically the integral we have done in (\ref{SL2Whitt}) with the only change that $\chi$ is now defined on all of $G(\rats_p)$ in which $SL(2,\rats_p)_{\alpha_i}$ is embedded. The result is determined by diagonal $a_i$ and reads for $\hat{a}=1$ 
\begin{align}
 \lint_{N_{i,-}(\rats_p)} \chi(n_{i,-} a_i) \overline{\psi_{i,-}(n_{i,-})} d n_{i,-}
= \gamma_p(a_i^{-\alpha_i} 
)  (1-p^{-(\langle\lambda|\alpha_i\rangle+1)}) 
\frac{1-p^{-\langle\lambda|\alpha_i\rangle}|a_i^{-\alpha_i}
|^{\langle\lambda|\alpha_i\rangle}}{1-p^{-\langle\lambda|\alpha_i\rangle}}|a_i^{\lambda+\rho-\alpha_i}|.
\end{align}
Under the Weyl reflection $w_i$  one has
$w_i\lambda = \lambda-\langle\lambda|\alpha_i\rangle \alpha_i$
and the function $F_{i}$ therefore satisfies
\begin{align}
F_{i} (w_i\lambda,g) = F_{i} (\lambda,g)  
\frac{1-p^{-(1+\langle w_i\lambda|\alpha_i\rangle)}}{1-p^{-(1+\langle\lambda|\alpha_i\rangle)}},
\end{align}
where one also must keep track of the non-trivial $\hat{a}$ given by the prefactor in (\ref{Firest}).
The relation (\ref{aux1aux2}) then gives immediately the same transformation under $w_i$ for $F(\lambda,g)$ and therefore for the unramified Whittaker function:
\begin{align}
    \label{eq:whittaker-special-relation}
W^{\circ} (w_i\lambda,a) =  \frac{\zeta_p(w_i,\lambda)}{\zeta_p(w_i,-\lambda)} W^{\circ}(\lambda,a)
\end{align}
where we defined the local $\zeta$ factor
\begin{align}
\zeta_p(w,\lambda)=\prod_{\alpha>0 \atop w\alpha<0}\frac{1}{1-p^{-(1+\langle\lambda|\alpha\rangle)}}.
\end{align}
For a general Weyl transformation $w\in \Weyl$ we find therefore 
\begin{align}
\label{Whittfuncrel}
W^{\circ} (w\lambda,a) = \frac{\zeta_p(w,\lambda)}{\zeta_p(w,-\lambda)}W^{\circ}(\lambda,a) 
\end{align}
which we check in example \ref{ex:whittaker-weyl-relation}.

This is not surprisingly the same factor that appeared in (the functional relation for) the constant term, see~\eqref{GindiKarp}.

\begin{example}
    \label{ex:whittaker-weyl-relation}
    Let us check \eqref{Whittfuncrel} with $w = w_i w_j$ starting from \eqref{eq:whittaker-special-relation} where $w_i$ and $w_j$ are two (different) fundamental reflections. Using \eqref{eq:whittaker-special-relation} twice we have that
    \begin{equation}
        W^\circ (w \lambda, a) = W^\circ (w_i w_j \lambda, a) = \frac{\zeta_p(w_i, w_j \lambda)}{\zeta_p(w_i, - w_j \lambda)} W^\circ(w_j \lambda, a) = \frac{\zeta_p(w_i, w_j \lambda)}{\zeta_p(w_i, - w_j \lambda)} \frac{\zeta_p(w_j, \lambda)}{\zeta_p(w_j, -\lambda)} W^\circ(\lambda, a) \, .
    \end{equation}

    Consider now the factor
    \begin{equation}
        \begin{split}
            \frac{\zeta_p(w_i, w_j \lambda)}{\zeta_p(w_i, - w_j \lambda)} \frac{\zeta_p(w_j, \lambda)}{\zeta_p(w_j, -\lambda)}
            &= \prod_{\substack{\alpha' > 0 \\ w_i \alpha' < 0}} \frac{1 - p^{-(1 - \langle w_j \lambda | \alpha' \rangle)}}{1 - p^{-(1 + \langle w_j \lambda | \alpha' \rangle)}}
            \prod_{\substack{\alpha > 0 \\ w_j \alpha < 0}}\frac{1 - p^{-(1 - \langle \lambda | \alpha \rangle)}}{1 - p^{-(1 + \langle \lambda | \alpha \rangle)}} \\
            &= \prod_{\substack{w_j \alpha > 0 \\ w_i w_j \alpha < 0}} \frac{1 - p^{-(1 - \langle \lambda | \alpha \rangle)}}{1 - p^{-(1 + \langle \lambda | \alpha \rangle)}}
            \prod_{\substack{\alpha > 0 \\ w_j \alpha < 0}}\frac{1 - p^{-(1 - \langle \lambda | \alpha \rangle)}}{1 - p^{-(1 + \langle \lambda | \alpha \rangle)}}
        \end{split} 
    \end{equation}
    where we have made the substitution $\alpha' = w_j \alpha$ and used the fact that $\langle w_j \lambda | w_j \alpha \rangle = \langle \lambda | \alpha \rangle$. Applying the same argument as in the proof of Lemma~\ref{lem-functintertwine} we can then combine the products into
    \begin{equation}
        W^\circ (w \lambda, a) = \prod_{\substack{\alpha > 0 \\ w \alpha < 0}}\frac{1 - p^{-(1 - \langle \lambda | \alpha \rangle)}}{1 - p^{-(1 + \langle \lambda | \alpha \rangle)}}  W^\circ(\lambda, a) = \frac{\zeta_p(w,\lambda)}{\zeta_p(w,-\lambda)}W^{\circ}(\lambda,a) \, ,
    \end{equation} 
     as claimed.
\end{example}

\subsection{Weyl invariant combination}

As for the constant term (and the full Eisenstein series),  one can obtain a Weyl invariant form by multiplying through by the  denominator of $\zeta$ factors associated with the longest Weyl word. Denoting
\begin{align}
\label{zetalambda}
\zeta (\lambda) \equiv \zeta_p(\wlong,\lambda) = \prod_{\alpha>0} \frac{1}{1-p^{-(1+\langle\lambda|\alpha\rangle)}}
\end{align}
one has that the function
\begin{align}
\label{Whittinv}
\zeta(\lambda) W^{\circ}(\lambda,a)
\end{align}
is Weyl invariant. This is checked simply by combining (\ref{Whittfuncrel}) with the transformation of $\zeta_p(\lambda)$ that can be derived straightforwardly.

Because of the Weyl invariance of (\ref{Whittinv}), we write it as a sum over Weyl images as
\begin{align}
\label{sumimages}
\zeta(\lambda) W^{\circ}(\lambda,a) = \sum_{w\in \Weyl} c(w\lambda) |a^{w\lambda +\rho}|
\end{align}
since the invariant function has to be a polynomial in $a^{\lambda+\rho}$ (and its images). The fact that the local Whittaker function is a single Weyl orbit follows from the considerations in~\cite{CasselmanShalika}.

\subsection{Determining a special coefficient}

Next we determine $c(w\lambda)$ for $w=\wlong$ which is the coefficient of $|a^{\wlong\lambda+\rho}|$ in (\ref{sumimages}). Referring back to (\ref{CStoUS}) we see that the coefficient of $|a^{\wlong\lambda+\rho}|$ in $W^{\circ}(\lambda,a)$ is obtained as the $a$-independent part of the integral
\begin{align}
\lint_{N(\rats_p)}\chi(\wlong n) \overline{\psi^a(n)} dn.
\end{align}
The integral is a polyonmial in $a$ and we can obtain its $a$-independent part formally by sending $a$ to zero. (This is only formal because, of course, $0\notin A(\rats_p)$.) Therefore, the $a$-independent part of this integral can be obtained by removing the character $\psi^a$ from the integral and then one is left with the same integral as in the constant term (\ref{Iw}) for $w=\wlong$. The result then is the same as the local factor for $\rats_p$ in  the constant term formula (\ref{LCF}), viz.
\begin{align}
c(\wlong\lambda)& = \zeta(\lambda)\prod_{\alpha>0}\frac{1-p^{-(\langle\lambda|\alpha\rangle+1)}}{1-p^{-\langle\lambda|\alpha\rangle}}
=\prod_{\alpha>0} \frac1{1-p^{-\langle\lambda|\alpha\rangle}}
= \prod_{\alpha>0} \frac1{1-p^{\langle\wlong\lambda|\alpha\rangle}}
= \eps(\wlong\lambda).
\end{align}
This means that the general coefficient is given by
\begin{align}
c(\lambda)=\eps(\lambda)=\prod_{\alpha>0} \frac1{1-p^{\langle\lambda|\alpha\rangle}}
\end{align}
and the general formula for the Whittaker function for an unramified character is
\begin{align}
\label{CS2}
W^{\circ}(\lambda,a) =  \frac1{\zeta(\lambda)} \sum_{w\in \Weyl} \eps(w\lambda) |a^{w\lambda+\rho}|,
\end{align}
thus demonstrating (\ref{CSus}). This concludes the proof of theorem~\ref{CasselmShalikaform}.

\section{Whittaker functions for generic \texorpdfstring{$\psi$}{psi}}
\label{sec:genpsi}

Theorem~\ref{CasselmShalikaform} is only valid for unramified character $\psi$, but we will now show that it can also be used for generic characters indirectly. Recall from definition~\ref{def_unram} that for an unramified character $m_\alpha = 1$ for all $\alpha \in \Pi$ and that a generic character has $m_\alpha \neq 0$ for all $\alpha$. 

Let us take a closer look at the so-called `twisted' character $\psi^a(n) = \psi(a n a^{-1})$ introduced above where $\psi$ without superscript $a$ is the unramified character. We note that periodicity of $\psi^a$ is of course different from the one of $\psi$, but this will not influence our reasoning.

From \eqref{eq:n-param} and \eqref{eq:psi-unramified} we have that
\begin{equation}
    \psi(n) = \exp \left( -2\pi i \left[ \sum_{\alpha \in \Pi} u_\alpha \right] \right), \qquad n = \bigg(\prod_{\alpha\in \Delta_+\setminus\Pi} x_{\alpha}(u_\alpha)\bigg) \bigg(\prod_{\alpha\in\Pi}x_\alpha(u_\alpha)\bigg)
\end{equation}
where $x_\alpha(u_\alpha) = \exp(u_\alpha E_\alpha)$.

With insertions of $a a^{-1}$, the expression for $a n a^{-1}$ splits into factors of $a x_\alpha(u_\alpha) a^{-1}$. 
Using the Baker-Campbell-Hausdorff formula, these factors can be found as
\begin{equation}
    a x_\alpha(u_\alpha) a^{-1} = \exp(e^t \, u_\alpha E_\alpha), \qquad 
\end{equation}
where $t$ is defined by
\begin{equation}
    [\log a, \log x_\alpha(u_\alpha)] = t \log x_\alpha(u_\alpha) \, .
\end{equation}
Let $a$ be parametrised as
\begin{equation}
    a = \exp\left( \sum_{\beta \in \Pi} \log (v_\beta) H_\beta \right) \, ,
\end{equation}
which gives
\begin{equation}
    [\log a, \log x_\alpha(u_\alpha)] = \sum_{\beta \in \Pi} \log(v_\beta) u_\alpha [H_\beta, E_\alpha] = \underbrace{\sum_{\beta \in \Pi} \alpha(H_\beta) \log(v_\beta)}_{=t} u_\alpha E_\alpha \, .
\end{equation}

Thus,
\begin{equation}
    a n a^{-1} = \bigg(\prod_{\alpha\in \Delta_+\setminus\Pi} x_{\alpha}(u'_\alpha)\bigg) \bigg(\prod_{\alpha\in\Pi}x_\alpha(u'_\alpha)\bigg), \qquad u'_\alpha = e^t u_\alpha = \bigg( \prod_{\beta\in\Pi} (v_\beta)^{\alpha(H_\beta)} \bigg) u_\alpha
\end{equation}
and, finally, by listing the simple roots as $\alpha_i \in \Pi$ for $i = 1, \ldots, r$ and denoting the associated elements $u_{\alpha_i}$ and $v_{\alpha_j}$ as $u_i$ and $v_j$ respectively
\begin{equation}
    \psi^a(n) = \exp\left(-2\pi i \left[ \sum_{i=1}^r u'_{i}  \right] \right) 
    = \exp\left( -2 \pi i \left[ \sum_{i=1}^r \bigg( \prod_{j=1}^r (v_j)^{A_{ji}} \bigg) u_{i}
 \right] \right)
\end{equation}
where we have introduced the Cartan matrix $A_{ij}$ defined in \eqref{CartanMatrix}.

We now note that this is really a generic character with
\begin{equation}
    m_i = m_{\alpha_i} = \prod_{j=1}^r (v_j)^{A_{ji}}
\end{equation}
and that any generic character can be expressed in this way with the inverse relation
\begin{equation}
\label{mtov}
    v_j = \prod_{i=1}^r (m_i)^{A^{-1}_{ij}}
\end{equation}
where $A^{-1}_{ij}$ is the inverse Cartan matrix.

Now that we can express a generic character in terms of the unramified character, we would like to find the Whittaker function for $\psi^a$ using \eqref{CS2} indirectly. More specifically, we ultimatelly want to find $W^{\circ}_{\psi^a}(\chi,a')$ with $a' = 1$ along the finite primes where $g = n a' k$ and $a' \in A(\rats_p)\subset G(\rats_p)$. For each $p$ this gives a contributing factor to the instanton measure as discussed in example \ref{ex_localWhittakerSL2}. 

This will bring us one step closer to finding the Fourier coefficients of the Eisenstein series with general instanton charges $m_\alpha$ in \eqref{Wwlong} and not only the restricted case of an unramified character. 

Using similar steps as taken in \eqref{CStoUS}, but in reverse order, we obtain
\begin{align}
\label{unramified2generic}
W^{\circ}_{\psi^a}(\chi,\id) = \left(\chi(\wlong a \wlong^{-1})\delta(a)\right)^{-1} W^{\circ}(\chi,a)=|a^{-(\wlong\lambda+\rho)}| W^\circ(\chi,a) \, .
\end{align}
Therefore, the local instanton measure for a generic character $\psi^a$ with instanton charges $m_\alpha$ can be expressed through the local instanton measure evaluated for an unramified character at non-trivial $a=\prod_{\alpha\in\Pi} v_\alpha^{H_\alpha}\in A(\rats_p)$.

\begin{example}
We illustrate formula (\ref{unramified2generic}) by recovering the result (\ref{SL2Wloc}) for $SL(2,\mathbb{A})$. In this case, there is only one simple root $\alpha$ and $A_{\alpha\alpha}=2$. The unramified Whittaker function is as given in (\ref{SL2Whitt}). If we want to get the Whittaker function for a character $\psi^a$ with instanton charge $m$, then (\ref{mtov}) tells us that we have $v=m^{1/2}$ and from (\ref{unramified2generic}) we find that 
\begin{align}
W^{\circ}_{\psi^a}(\chi,\id) &= |v|^{2s-2} \gamma_p(v^2)(1-p^{-2s})\frac{|v|^{-2s+2}-p^{-2s+1}|v|^{2s}}{1-p^{-2s+1}}&\nn\\
&=\gamma_p(m) (1-p^{-2s})\frac{1-p^{-2s+1}|m|^{2s-1}}{1-p^{-2s+1}}
\end{align}
in agreement with (\ref{SL2Wloc}).
\end{example}

\section{Degenerate Whittaker coefficients}
\label{sec:degpsi}

While the Casselman--Shalika formula (\ref{CS2}) provides an elegant expression for unramified \emp{local} characters  and, via (\ref{unramified2generic}), also for Fourier coefficients of generic characters, it is desirable to understand also Fourier coefficients for \emph{non}-generic characters $\psi$. These are also sometimes referred to as \emphindex[Whittaker function!degenerate]{degenerate Whittaker functions} in the literature~\cite{Zelevinsky80,MoeglinWaldspurger,GourevitchSahi1} and have the property that they only depend on a subset of the simple roots of $G(\mathbb{R})$ rather than all simple roots.

In this section, we will prove the following theorem that holds for \emp{global} characters $\psi$~\cite{FKP2013,Hashizume}:
\begin{theorem}
\label{DegWhittThm}
Let $\psi: N(\rats)\backslash N(\ads) \to U(1)$ be a degenerate character with $\supp(\psi)=\Pi'\neq \Pi$ with associated subgroup $G'(\mathbb{A})\subset G(\mathbb{A})$. Let $w_c\wlong'$ be the representatives of the coset $\Weyl/\Weyl'$ defined below in~\eqref{Cpsidec}. Then the degenerate Whittaker coefficient on $G(\mathbb{A})$ is given by
\begin{align}
\label{degW1}
W^\circ_\psi (\chi,a) = \sum_{w_c\wlong'\in \Weyl/\Weyl'} a^{(w_c\wlong')^{-1}\lambda+\rho}  M(w_c^{-1},\lambda) W'^\circ_{\psi^a}(w_c^{-1}\lambda,\id),
\end{align}
where $W'^\circ_{\psi}$ denotes a Whittaker function on the $G'(\mathbb{A})$ subgroup of $G(\mathbb{A})$. The weight $w_c^{-1}\lambda$ is given as a weight of $G'(\mathbb{A})$ by orthogonal projection. 
\end{theorem} 

\begin{remark}
In this theorem and in the remainder of the chapter we suppress the adelic absolute value on $|a^\mu|$ in order to ease the notation. 
\end{remark}

Before embarking on the proof, we explain the notation used here. For a global character 
\begin{align}
\psi\left(\prod_{\alpha\in\Pi} x_\alpha(u_\alpha)\right) = \exp\left( 2\pi i \sum_{\alpha\in\Pi} m_\alpha u_\alpha\right),
\end{align}
we call
\begin{align}
\label{supppsi}
\supp(\psi)  = \left\{ \alpha\in\Pi \,|\, m_\alpha\neq 0\right\}\subset \Pi
\end{align}
determined by the non-vanishing $m_\alpha$ the \emphindex[character!support of]{support of the character} $\psi$. With this notion, the definition~\ref{def:GC} becomes
\begin{align}
\supp(\psi) =\Pi \quad&\Longleftrightarrow\quad \textrm{$\psi$ generic},\nn\\
\supp(\psi) \neq\Pi \quad&\Longleftrightarrow \quad\textrm{$\psi$ non-generic or degenerate.}\nn
\end{align}
We note that a degenerate character $\psi:N(\rats)\backslash N(\ads) \to U(1)$ canonically defines a simple proper subgroup $G'\subset G$. This subgroup $G'$ is the one with simple root system $\Pi'=\supp(\psi)$; its Dynkin diagram is the subdiagram of the Dynkin diagram of $G$ obtained by restricting to the nodes corresponding to $\supp(\psi)$. The subgroup $G'$ has a Weyl group $\Weyl'$ with longest Weyl word $\wlong'$.

\begin{proof}
Using the Bruhat decomposition, the spherical Whittaker coefficient $W_\psi^\circ(\chi,a)$ can be written as a sum over the Weyl group $\Weyl$ of $G$ as in~\eqref{BruhatW}
\begin{align}
\label{Wintdeg}
W_\psi^\circ(\chi,a) = \lint_{N(\rats)\backslash N(\mathbb{A})} E(\chi,na) \overline{\psi(n)}dn = \sum_{w\in \Weyl} F_{w,\psi} (\chi,a)
\end{align}
with
\begin{align}
\label{FwSplit2}
F_{w,\psi} (\chi,a) = \lint_{w^{-1}B(\rats)w\cap N(\rats)\backslash N(\mathbb{A})} \chi(wna) \overline{\psi(n)}dn
=\int\limits_{N^w_{\{\beta\}}}\overline{\psi(n_\beta)}dn_\beta\int\limits_{N^w_{\{\gamma\}}}\chi(w n_\gamma a)\overline{\psi(n_\gamma)}dn_\gamma\,.
\end{align}
The various $F_{w,\psi}$ can be analysed as in section~\ref{sec:redlong} and we have used the $w$-dependent split of positive roots of $N$ into two sets of $\{\beta\}$ and $\{\gamma\}$ as in~\eqref{NBeta} and~\eqref{NGamma}. Importantly, for degenerate $\psi$ the integral over the compact domain $N_{\{\beta\}}^w$ can be non-vanishing for Weyl words $w$ different from $\wlong$: if $\psi$ is trivial on all the $n_\beta$ in (\ref{FwSplit2}), the corresponding integral yields unity rather than zero. This means that the Weyl word $w$ must map all elements in $\supp(\psi)$ to negative roots in order to avoid the vanishing of $F_{w,\psi}$ and the sum over $\Weyl$ in (\ref{Wintdeg}) can be restricted to the subset
\begin{align}
\mathcal{C}_\psi = \left\{ w\in\Weyl\,|\, w\alpha<0\quad\textrm{for all $\alpha\in\supp(\psi)$}\right\}.
\end{align}
(If $\psi$ is generic, one recovers $\mathcal{C}_\psi=\left\{\wlong\right\}$ in agreement with the discussion of section~\ref{sec:redlong}.) 

We will now parametrise the set $\mathcal{C}_\psi$ explicitly. Denote by $\Weyl'$ the Weyl subgroup generated by the fundamental reflections associated with $\Pi'=\supp(\psi)$ only. It is the Weyl group of  $G'$ and has its own longest Weyl word that we denote by $\wlong'$. The longest Weyl word $\wlong'$ has the desired property that it maps all elements in $\supp(\psi)$ to negative roots and it is the only Weyl word in $\Weyl'$ with this property. In fact, any element in $\mathcal{C}_\psi$ can be represented in a form that involves the longest word $\wlong'$ of $\Weyl'$:
\begin{align}
\label{Cpsidec}
w\in\mathcal{C}_\psi \quad\Longleftrightarrow \quad
w= w_c\wlong'.
\end{align}
Here, $w_c\in \Weyl$ must satisfy
\begin{align}
\label{wccond}
w_c \alpha>0 \quad \textrm{for all $\alpha\in\supp(\psi)$}
\end{align}
in order for $w=w_c\wlong'$ to belong to $\mathcal{C}_\psi$. 

The words $w_c\in\Weyl$ can be constructed as carefully chosen representatives of the coset $\Weyl/\Weyl'$. Consider the weight
\begin{align}
\Lambda_\psi = \sum_{i\,:\,\alpha_i \notin \supp(\psi)} \Lambda_i,
\end{align}
 i.e., the sum of fundamental weights of $G$ that are not associated with the support of the degenerate character $\psi$. The weight $\Lambda_\psi$ is stabilised by $\Weyl'$ and its $\Weyl$-orbit is in bijection with the coset $\Weyl/\Weyl'$. A standard result for Weyl groups is that if $w(\alpha_i)<0$ for some simple root, then $\ell(ww_i)<\ell(w)$~\cite[Lemma~3.11]{Kac}. Therefore, if $w(\alpha_i)<0$ and $\alpha_i\in\supp(\psi)$ we have
 \begin{align}
w(\Lambda_\psi) = (ww_i)(\Lambda_\psi)
 \end{align}
since $\Lambda_\psi$ is stabilised by the fundamental reflections from $\supp(\psi)$ (these generate $\Weyl'$). This means that if $w(\alpha_i)<0$ there is a shorter Weyl word $ww_i$ leading to the same point in the $\Weyl$-orbit of $\Lambda_\psi$ as the word $w$ does. By induction, the shortest word leading to a given point in the Weyl orbit of $\Lambda_\psi$ must be those $w_c\in\Weyl$ that satisfy $w_c\alpha>0$ for all $\alpha\in\supp(\psi)$. Hence, the words $w_c$ appearing in~\eqref{wccond} are the shortest words leading to the points of the $\Weyl$-orbit of $\Lambda_\psi$. Such shortest words are not necessarily unique; for a given $\Weyl$-orbit point any shortest word $w_c$ will do. An explicit construction of the $w_c$ can be achieved by the same orbit method as in section~\ref{EvalLCF}, see also~\cite{FK2012,FKP2013}.

With the parametrisation $w=w_c\wlong'$ of the elements of $\mathcal{C}_\psi$ we thus arrive at the following expression for the degenerate Whittaker integral (\ref{Wintdeg}):
\begin{align}
W^\circ_\psi (\chi,a) = \sum_{w_c\wlong'\in \Weyl/\Weyl'} F_{w_c\wlong',\psi}(\chi,a),
\end{align}
where it is understood that $w_c\wlong'$ is the specific coset representative described above.

The quantities $F_{w_c\wlong',\psi}(\chi,a)$ can be evaluated by reducing them to Whittaker functions of the subgroup $G'(\mathbb{A})\subset G(\mathbb{A})$ associated with $\supp(\psi)$ as follows. First, we separate out the $a$-dependence as usual by conjugating it to the left and using the multiplicativity of $\chi$
\begin{align}
 F_{w_c\wlong',\psi}(\chi,a) &= \lint_{(w_c\wlong')^{-1} B(\rats)w_c\wlong'\cap N(\rats)\backslash N(\mathbb{A})} \chi(w_c\wlong' na) \overline{\psi(n)}dn&\nonumber\\
 & = a^{(w_c\wlong')^{-1}\lambda+\rho} \lint_{(w_c\wlong')^{-1} B(\rats)w_c\wlong'\cap N(\rats)\backslash N(\mathbb{A})} \chi(w_c\wlong' n) \overline{\psi^a(n)}dn&
\end{align}
with $\psi^a(n) = \psi(ana^{-1})$ as before. We can also rewrite the integration into the two sets $N_{\{\beta\}}^w$ and $N_{\{\gamma\}}^w$ (for $w=w_c\wlong'$) as in~\eqref{FwSplit2} and we know that by construction the integral over $N_{\{\beta\}}^w$ gives unity.

The remaining integral over $N_{\{\gamma\}}^w$ is then over all positive roots $\gamma$ that are mapped to negative roots by the action of $w=w_c\wlong'$ and we drop the $\gamma$ subscript for ease of notation. The particular form of $w$ implies that we can parametrise the unipotent element as $n=n_c n'$ where $n'\in N'(\mathbb{A})$ is the (full) unipotent radical of the standard minimal Borel subgroup $B'(\ads)$ of $G'(\mathbb{A})$ that is determined by $\psi$; and $n_c$ are the remaining elements whose total space we call $N_c(\mathbb{A})$. We note also that $\wlong' n_c (\wlong')^{-1}$ is generated exactly by the positive roots that are mapped to negative roots by $w_c$ alone. The degenerate character $\psi$ only depends on $n'$, i.e., $\psi^a(n_cn') = \psi^a(n')$.

Putting these observations together one obtains
\begin{align}
\label{FintDW}
 F_{w_c\wlong',\psi}(\chi,a) = a^{(w_c\wlong')^{-1}\lambda+\rho} \lint_{N_c(\mathbb{A})} \lint_{N'(\mathbb{A})}  
  \chi(w_c\wlong' n_c n') \overline{\psi^a(n')}dn_c dn'.
\end{align}
As the next step one can rewrite the argument of the character $\chi$ as
\begin{align}
\chi(w_c\wlong' n_c n') = \chi\left(w_c\wlong' n_c (w_c\wlong')^{-1}w_c\wlong' n'\right)= \chi\left(w_c\wlong' n_c (w_c\wlong')^{-1}w_c\tilde{n}\tilde{a}\right),
\end{align}
where we have performed an Iwasawa decomposition (in $G'(\mathbb{A})$) of $\wlong'n' = \tilde{n}\tilde{a}\tilde{k}$ and used left-invariance of $\chi$ under $K'(\mathbb{A})\subset K(\mathbb{A})$ in the last step. In the next step we want to perform another Iwasawa decomposition (now in $G(\mathbb{A})$) of 
\begin{align}
\label{IwaDW}
w_c\tilde{n}\tilde{a} = \hat{n}\hat{a}\hat{k}.
\end{align} 
The important observation now is that $\tilde{n}\in N'(\mathbb{A})$ and $w_c$ satisfies~\eqref{wccond} which implies that $w_c N'(\mathbb{A})w_c^{-1} \subset N(\mathbb{A})$. Therefore, the Iwasawa decomposition~\eqref{IwaDW} has
\begin{align}
\label{IwaDW2}
\hat{n} = w_c\tilde{n} w_c^{-1},\quad
\hat{a} = w_c\tilde{a} w_c^{-1},\quad
\hat{k} = w_c.
\end{align} 
Inserting this back into the integral~\eqref{FintDW} one can bring the element $\hat{n}\in N(\mathbb{A})$ to the left. This will induce a uni-modular change of the integration variables $dn_c$ as in section~\ref{sec:CSTInd}. Conjugating the element $\hat{a}\in A(\mathbb{A})$ to the left will induce a non-trivial change of measure ($w=w_c\wlong'$):
\begin{align}
\lint_{N_c(\mathbb{A})} \chi(w n_c w^{-1} \hat{n} \hat{a}) dn_c
&= \lint_{N_c(\mathbb{A})} \chi(\hat{n} \hat{a} w n_c w^{-1}) \hat{a}^{w_c\rho-\rho} dn_c 
 = \lint_{N_c(\mathbb{A})} \chi(w n_c w^{-1}) \tilde{a}^{w_c^{-1}\lambda-\rho} dn_c \nn\\
&=\tilde{a}^{w_c^{-1}\lambda-\rho} \lint_{N_c(\mathbb{A})} \chi(w n_c w^{-1})  dn_c
= \chi'(\tilde{a}) \lint_{N_c(\mathbb{A})} \chi(w n_c w^{-1})  dn_c.
\end{align}
We have evaluated the character $\chi$ on $\hat{n}\hat{a}$ in the second step according to $\chi(\hat{n}\hat{a}) = \hat{a}^{\lambda+\rho} = \tilde{a}^{w_c^{-1}\lambda +w_c^{-1}\rho}$ due to~\eqref{IwaDW2}. In the last step, we have used that $\tilde{a}$ does not depend on $n_c$ and can therefore be taken out of the integral and defined the character 
\begin{align}
\chi'(\tilde{a}) = \tilde{a}^{w_c^{-1} \lambda+\rho} = \chi'(\tilde{a})  = \chi'(\wlong' n'),
\end{align}
on the group $G'(\mathbb{A})$. In the last step we have used the definition of $\tilde{a}$.

Putting everything together in~\eqref{FintDW} one obtains the factorised expression
\begin{align}
 F_{w_c\wlong',\psi}(\chi,a) = a^{(w_c\wlong')^{-1}\lambda+\rho} \lint_{N_c(\mathbb{A})} \chi(w_c\wlong'n_c) dn_c \cdot \lint_{N'(\mathbb{A})} \chi'(\wlong' n') \overline{\psi^a(n')} dn'.
\end{align}
The two separate integrals are both of types we have encountered before: The $N_c(\mathbb{A})$ integral is precisely the Gindikin--Karpelevich expression (\ref{Iw}) for the Weyl word $w_c\in\Weyl$ and so gives a factor $M(w_c^{-1},\lambda)$ defined in \eqref{eq:intertwiner}, and the second integral is the generic Whittaker function~\eqref{Wwlong} for the subgroup $G'(\mathbb{A})\subset G(\mathbb{A})$ with \emp{generic} Fourier character $\psi^a$, in the representation given by the weight $w_c^{-1}\lambda$, projected orthogonally to $G'(\mathbb{A})$ and evaluated at the identity $\id\in A'(\mathbb{A})$. This completes the proof of theorem~\ref{DegWhittThm}.
\end{proof}

As a consequence of the theorem, Whittaker functions of non-generic characters $\psi$ can be evaluated as sums over Whittaker functions of subgroups on which the character is generic. We stress again that the choice of coset representative of $\Weyl/\Weyl'$ is important here. If the full Whittaker function on the subgroup is known, the above formula provides the explicit expression for any character $\psi$. Thanks to the Casselman--Shalika formula, this means that the \emph{local} Whittaker function ($p<\infty$) can be calculated for \emph{any} character, generic or not. The archimedean part is typically more intricate.

\begin{remark}\label{rmk:degW}
Theorem~\ref{DegWhittThm} of course also remains true in the case of generic $\psi$ since then $\Weyl'=\{\id\}$ is trivial and the sum on the right-hand side is just the decomposition of the generic Whittaker function into Bruhat cells and nothing is gained. The power of the theorem arises in cases where one deals with an Eisenstein series that does not have any generic Whittaker coefficients and one can then use~\eqref{degW1} to determine the degenerate ones. This will be explored in more detail in section~\ref{sec:CSeval} below.
\end{remark}

\section{Whittaker coefficients on~\texorpdfstring{$SL(3,\mathbb{A})$}{SL(3, A)}}
\label{sec:SL3ex}

We illustrate the general considerations above through the explicit example of $SL(3,\mathbb{A})$. The Eisenstein series on $SL(3,\reals)$, $GL(3,\reals)$ and this group have been studied in great detail in the literature~\cite{VT,BumpSL3} by various techniques. 

The split real group $SL(3,\reals)$ has rank two and we denote the two simple roots by $\alpha_1$ and $\alpha_2$. The corresponding Cartan generators will be called $H_1\equiv H_{\alpha_1}$ and $H_2\equiv H_{\alpha_2}$. A general element $a\in A(\mathbb{A})$ will be written as 
\begin{align}
a= v_1^{H_1} v_2^{H_2}.
\end{align}
The Eisenstein series is determined by the weight 
\begin{align}
\lambda = (2s_1-1) \Lambda_1 + (2s_2-1) \Lambda_2
\end{align}
in terms of the fundamental weights dual to the simple roots.

The Weyl group consists of six elements:
\begin{align}
\Weyl = \left\{ \id, w_1, w_2, w_1w_2, w_2w_1, w_1w_2w_1\right\}.
\end{align}

We will first compute the constant term using Langlands constant term formula of chapter~\ref{ch:CTF}. Then, using the results of sections \ref{sec:redlong}--\ref{sec_CSformula}, we find the local part of a Whittaker function with an unramified character, which, with the help of section \ref{sec:genpsi}, can then be used to compute the local part of any Whittaker function with a generic character. The remaining, degenerate, Whittaker functions are then found following the arguments of section \ref{sec:degpsi}. Lastly, we will comment on the non-abelian Whittaker coefficients.

\subsection{Constant terms}

We first evaluate the Langlands constant term formula (\ref{LCF}). This yields a sum of six terms:
\begin{align}
\label{SL3consts}
&\lint_{N(\rats)\backslash N(\mathbb{A})} E(\chi,ng) dn = v_1^{2s_1} v_2^{2s_2}
+\frac{\xi(2s_1-1)}{\xi(2s_1)} v_1^{2-2s_1} v_2^{2s_1+2s_2-1} 
+\frac{\xi(2s_2-1)}{\xi(2s_2)} v_1^{2s_1+2s_2-1} v_2^{2-2s_2} \nonumber\\
&\quad+\frac{\xi(2s_1-1)\xi(2s_1+2s_2-2)}{\xi(2s_1)\xi(2s_1+2s_2-1)} v_1^{2s_2} v_2^{3-2s_1-2s_2}
+\frac{\xi(2s_2-1)\xi(2s_1+2s_2-2)}{\xi(2s_2)\xi(2s_1+2s_2-1)} v_1^{3-2s_1-2s_2} v_2^{2s_1}\nonumber\\
&\quad+\frac{\xi(2s_1-1)\xi(2s_2-1)\xi(2s_1+2s_2-2)}{\xi(2s_1)\xi(2s_2)\xi(2s_1+2s_2-1)} v_1^{2-2s_2}v_2^{2-2s_1}.
\end{align}
Here, $v_1$ and $v_2$ are real positive parameters.

\subsection{Generic Whittaker coefficients}
\label{sec:SL3-generic-Whittaker-coefficients}
We first determine the \emp{local} Whittaker function for an unramified character $\psi$ by using the Casselman--Shalika formula in the form (\ref{CSus}). The quantities $1/\zeta(\lambda)$ and $\epsilon(\lambda)$ of (\ref{zetaeps}) evaluate to 
\begin{subequations}
\begin{align}
\label{zetSL3}
\frac1{\zeta(\lambda)} &= (1- p^{-2s_1})(1-p^{-2s_2})(1-p^{1-2s_1-2s_2}),&\\
\eps(\lambda) &= \frac{1}{(1-p^{2s_1-1})(1-p^{2s_2-1})(1-p^{2s_1+2s_2-2})}
\end{align}
\end{subequations}
and the full unramified local coefficient is then 
\begin{align}
\label{SL3unr}
W^\circ(\lambda,a) &= \frac{\eps(\lambda)}{\zeta(\lambda)} \Big(
|v_1|^{2s_1} |v_2|^{2s_2} 
- p^{2s_1-1} |v_1|^{2-2s_1} |v_2|^{2s_1+2s_2-1} 
- p^{2s_2-1} |v_1|^{2s_1+2s_2-1} |v_2|^{2-2s_2} &\nonumber\\
&\quad+ p^{4s_1+2s_2-3} |v_1|^{2s_2} |v_2|^{3-2s_1-2s_2} 
+ p^{2s_1+4s_2-3} |v_1|^{3-2s_1-2s_2} |v_2|^{2s_1} &\nonumber\\
&\quad-p^{4s_1+4s_2-4} |v_1|^{2-2s_2} |v_2|^{2-2s_1}
\Big).
\end{align}
Here, $v_1$ and $v_2$ are in $\rats_p$. 

{}From (\ref{SL3unr}) we can deduce the Whittaker coefficient for a generic character with non-zero instanton charges $m_1$ and $m_2$, i.e., one that satisfies
\begin{align}
\psi^a\left( x_{\alpha_1} (u_1) x_{\alpha_2} (u_2)\right) = \exp\left(2\pi i [m_1 u_1 + m_2 u_2]\right)
\end{align}
by exploiting (\ref{unramified2generic}). For this we require $v_1 = m_1^{2/3} m_2^{1/3}$ and $v_2=m_1^{1/3} m_2^{2/3}$ in the expression above as well as the prefactor $|a^{-(\wlong\lambda+\rho)}|= |v_1|^{2s_2-2}|v_2|^{2s_1-2}$.
The result is
\begin{align}
\label{SL3gen}
W^\circ_{\psi^a}(\chi,\id) &= \frac{\eps(\lambda)}{\zeta(\lambda)} \Big(
|m_1|^{2s_1+2s_2-2} |m_2|^{2s_1+2s_2-2}
-p^{2s_1-1} |m_1|^{2{s_2}-1}|m_2|^{2s_1+2s_2-2} &\nonumber\\
&\quad-p^{2s_2-1} |m_1|^{2s_1+2s_2-2}|m_2|^{2s_1-1} 
 + p^{4s_1+2s_2-3} |m_1|^{2s_2-1}
  + p^{2s_1+4s_2-3} |m_2|^{2s_1-1}&\nonumber\\
  &\quad - p^{4s_1+4s_2-4}\Big).&
\end{align}
As is well-known~\cite{BumpSL3}, this can also be expressed in terms of a Schur polynomial in $(m_1,m_2)$ which here encodes the character of a highest weight representation of $\mathfrak{sl}(3,\cx)$. Taking the product over all $p<\infty$ produces double divisor sums. 

In this case, we can also work out the archimedean Whittaker function. The Whittaker function at $p=\infty$ can be explicitly written as a convoluted integral of two modified Bessel functions as we will now show.

Starting from \eqref{eulerWhittaker} and using the same standard manipulations as in \eqref{CStoUS} we have that 
\begin{equation}
    W_{\psi_\infty}^\circ(\chi_\infty, a) = \int_{N(\reals)} \chi_\infty(\wlong n a) \overline{\psi_\infty(n)} dn = |a^{\wlong \lambda + \rho}| \int_{N(\reals)} \chi_\infty(\wlong n) \overline{\psi_\infty(ana^{-1})} dn 
\end{equation}
with the generic character $\psi_\infty$ given by two integers $m_1$ and $m_2$ through
\begin{equation}
    \begin{gathered}
        \psi_\infty\big(x_{\alpha_1}(u_1) x_{\alpha_2}(u_2)\big) = \exp\big(2\pi i(m_1 u_1 + m_2 u_2)\big) \\
        \chi_\infty(v_1^{H_1}v_2^{H_2}) = \abs{v_1}^{2 s_1} \abs{v_2}^{2 s_2} \\
        \abs{a^{\wlong \lambda + \rho}} = \abs{v_1}^{2 - 2 s_2} \abs{v_2}^{2 - 2 s_1} \\
        n = 
        \begin{psmallmatrix}
            1 & u_1 & z \\
            0 & 1 & u_2 \\
            0 & 0 & 1
        \end{psmallmatrix} \qquad
        \wlong = 
        \begin{psmallmatrix}
            0 & 0 & 1 \\
            0 & 1 & 0 \\
            1 & 0 & 0
        \end{psmallmatrix} 
        \, .
    \end{gathered}
\end{equation}
Evaluating the integrand we obtain
\begin{multline}
    \int_{N(\reals)} \chi_\infty(\wlong n) \overline{\psi(ana^{-1})} dn = \\
     \int_{\reals^3} \big(1 + (1 + u_1^2)u_2^2 - 2 u_1 u_2 z + z^2\big)^{-s_1} \big(1 + u_1^2 + z^2\big)^{-s_2} \times  \\ \times \exp\left(-2\pi i \frac{v_1^3 m_1 u_1 + v_2^3 m_2 u_2}{v_1 v_2} \right) du_1 du_2 dz \, .
\end{multline}
Using the variable substitution $u_2 \to (u_2 + u_1 z)/(1 + u_1^2)$ and integrating over $u_2$ we get
\begin{multline}
    \frac{2 \pi^{s_1}}{\Gamma(s_1)} \abs{\frac{m_2 v_2^2}{v_1}}^{s_1 - \frac{1}{2}} \int_{\reals^2} \frac{(1 + u_1^2 + z^2)^{\frac{1}{4} - \frac{1}{2} s_1 - s_2}}{\sqrt{1 + u_1^2}} K_{s_1 - \frac{1}{2}}\left(2 \pi \abs{\frac{m_2 v_2^2}{v_1}} \frac{\sqrt{1 + u_1^2 + z^2}}{1 + u_1^2}\right) \times \\ \times \exp\left(-2\pi i \left( \frac{m_1 v_1^2}{v_2} u_1 + \frac{m_2 v_2^2}{v_1} \frac{u_1 z}{1 + u_1^2} \right) \right) du_1 dz \, .
\end{multline}
With standard manipulations (see for example \cite[Lemma~7]{VT}), this integral can be expressed as a convoluted integral of two Bessel functions giving $W_{\psi_\infty}^\circ$ as
\begin{multline}
    \label{eq:SL3-generic-archimedean}
    W_{\psi_\infty}^\circ(\chi_\infty, a) = \frac{4 \pi^{2s_3 + \frac{1}{2}}\abs{v_1 v_2}}{\Gamma(s_1) \Gamma(s_2) \Gamma(s_3)} \abs{m_1 m_2}^{s_3 - \frac{1}{2}}  \abs{\frac{v_1}{v_2}}^{s_1 - s_2} \times \\ \times \int_0^\infty K_{s_3 - \frac{1}{2}}\!\left( 2\pi \abs{\frac{m_1 v_1^2}{v_2}} \sqrt{1 + 1/x} \right) K_{s_3 - \frac{1}{2}}\!\left( 2\pi \abs{\frac{m_2 v_2^2}{v_1}} \sqrt{1 + x} \right) x^{\frac{s_2 - s_1}{2}} \frac{dx}{x}
\end{multline}
where we have introduced $s_3 = s_1 + s_2 - \frac{1}{2}$ for compactness.

\subsection{Degenerate Whittaker coefficients}
\label{sec:SL3-degenerate-Whittaker-coefficients}

We now evaluate the Whittaker coefficients for non-generic characters, i.e., those where either $m_1$ or $m_2$ vanishes. Note that it is not trivially possible to obtain this result from the generic one above by setting some parameters to zero. We will employ theorem~\ref{DegWhittThm} and perform this in the example $m_2=0$. Then the support of the character is only on the first simple root, so that $\wlong'=w_1$ and the subgroup $G'(\mathbb{A})$ is the one associated with the first simple root only.
The possible Weyl words that contribute to~\eqref{degW1} are
\begin{align}
\label{eq:SL3cosW}
w=w_c\wlong' \in \left\{ \id \wlong', w_2\wlong',w_1w_2\wlong'\right\}
= \left\{ w_1, w_2w_1, w_1w_2w_1\right\}.
\end{align}
As a first step, we calulcate the projected weights $w_c^{-1}\lambda$ and $M(w_c^{-1},\lambda)$ factors that appear in (\ref{degW1}) for the three choices:
\begin{subequations}
\begin{align}
w_c&=\id :& \lambda' &= (w_c^{-1}\lambda)_{G'} = (2s_1-1)\Lambda'_1,& M(w_c^{-1},\lambda)&=1,\\
w_c&=w_2 :& \lambda' &=(2s_1+2s_2-2)\Lambda'_1,& M(w_c^{-1},\lambda)&=\frac{\xi(2s_2-1)}{\xi(2s_2)},\\
w_c&=w_1w_2 :& \lambda' &= (2s_2-1)\Lambda'_1,& M(w_c^{-1},\lambda)&=\frac{\xi(2s_1-1)\xi(2s_1+2s_2-2)}{\xi(2s_1)\xi(2s_1+2s_2-1)}.
\end{align}
\end{subequations}
where $\Lambda'_1 = \alpha_1/2$ is the fundamental weight for $G'(\adeles)$.

This will need to be combined with 
\begin{align}
\psi^a(x_{\alpha_1}(u_1)) = \psi(a x_{\alpha_1} (u_1)a^{-1}) = \exp\left(2\pi i a^{\alpha_1} u_1 m_1\right)
= \exp\left(2\pi i v_1^2v_2^{-1} m_1 u_1\right)
\end{align}
and the $SL(2,\mathbb{A})$ Whittaker function for $\lambda'=(2s'-1)\Lambda'_1$ given by (cf.~(\ref{SL2FourierFull})) 
\begin{align}
W'^\circ_{\psi^a}(\lambda',\id) 
= \frac{2 (2\pi)^{1/2-s'}}{\xi(2s')} \sigma_{2s'-1}(m_1) \mathcal{K}_{1/2-s'} (2\pi |m_1| v_1^2v_2^{-1}),
\end{align}
where we introduced the short-hand $\mathcal{K}_t(x) = x^{-t}K_{-t}(x)$ in order to facilitate comparison with~\cite{Pioline:2009qt}. Recall also  the compact notation $s_3=s_1+s_2-\frac12$.
The resulting expression for the $(m_1,0)$ degenerate Whittaker function is then
\begin{align}
    \label{eq:SL3-degenerate}
W^\circ_\psi(\chi,a) &= \frac{2 (2\pi)^{1/2-s_1}}{\xi(2s_1)}  v_1^{2-2s_1} v_2^{2s_1+2s_2-1} \sigma_{2s_1-1}(m_1) \mathcal{K}_{1/2-s_1} (2\pi |m_1| v_1^2v_2^{-1})\nonumber\\
&\quad\!\!\!\!+\frac{2 (2\pi)^{1/2-s_3}}{\xi(2s_3)}\frac{\xi(2s_2-1)}{\xi(2s_2)}  v_1^{3-2s_1-2s_2} v_2^{2s_1} \sigma_{2s_3-1}(m_1) \mathcal{K}_{1/2-s_3} (2\pi |m_1| v_1^2v_2^{-1})\nonumber\\
&\quad\!\!\!\!+\frac{2 (2\pi)^{1/2-s_2}}{\xi(2s_2)}  \frac{\xi(2s_1-1)\xi(2s_3-1)}{\xi(2s_1)\xi(2s_3)}
v_1^{2-2s_2} v_2^{2-2s_1} \sigma_{2s_2-1}(m_1) \mathcal{K}_{1/2-s_2} (2\pi |m_1| v_1^2v_2^{-1}).
\end{align}
This matches also the expressions in~\cite{Pioline:2009qt} if one adapts the conventions and corrects a typo there. More precisely, one uses $v_1=\nu^{-1/6}\tau_2^{1/2}$, $v_2=\nu^{-1/3}$ and exchanges $s_1$ and $s_2$ to find the $\Psi_{0,q}$ coefficient in~\cite[Eq.~(3.45)]{Pioline:2009qt}, if one fixes the third summand there. A similar calculation can be carried out for the degenerate Whittaker coefficient associated with instanton charges $(0,m_2)$; it simply amounts to interchanging the subscripts $1$ and $2$ everywhere thanks to the Dynkin diagram automorphism of $\mathfrak{sl}(3,\reals)$.

\subsection{Non-abelian Fourier coefficients}
\label{sec:SL3-non-abelian}

So far in this section we have only studied Fourier coefficients $N$, but since the characters on $N$ are trivial on the centre $Z = N^{(2)} = [N, N]$ they do not capture the complete Fourier expansion of $E(\chi, g)$ as discussed in section \ref{sec:ANA}. To have a complete expansion we also need Fourier coefficients on $Z$ with non-trivial characters $\psi_Z : Z(\rats) \bs Z(\ads) \to U(1)$ parametrised by $k \in \rats^\times$
\begin{equation}
    \psi_Z(n_{(2)}) = e^{2\pi i k z} \qquad n_{(2)} = 
    \begin{psmallmatrix} 
        1 & 0 & z \\
        0 & 1 & 0 \\
        0 & 0 & 1
    \end{psmallmatrix}
    \in Z(\ads) \, .
\end{equation}
To avoid ambiguities, we will denote the character on $Z$ by $\psi_Z$ and the characters on $N$ by $\psi_N$.

Recalling \eqref{eq:non-abelian-Whittaker}, the Fourier coefficients on $Z$ are defined by
\begin{equation}
    W^\circ_{\psi_Z}(\chi, g) = \intl{Z(\rats)\bs Z(\ads)} E(\chi, n_{(2)} g) \overline{\psi_Z(n_{(2)})} \, dn_{(2)} \, .
\end{equation}
For the remaining parts of this section, we will drop the superscript for the spherical property and write the charge explicitly as $W_{\psi_Z}^{(k)}$ for clarity.

We will now show that these Fourier coefficients on $Z$ are determined by the Whittaker coefficients (i.e. Fourier coefficients on $N$)  but before we can make an exact statement we need to make a few definitions.

Let $k, m_2 \in \rats$ with $k = a_1/b_1$ and $m_2 = a_2/b_2$ in shortened form where $a_i \in \ints$ and $b_i \in \nats$. Define
\begin{equation}
    \label{eq:gcd-rational}
    d = d(k, m_2) \coloneqq \frac{\gcd(a_1 b_2, a_2 b_1)}{b_1 b_2}
\end{equation}
which is then strictly positive since $k \neq 0$ and let $k' \coloneqq k/d = a_1 b_2 / \gcd(a_1 b_2, a_2 b_1) \in \ints$ and $m'_2 \coloneqq m_2 / d = a_2 b_1 / \gcd(a_1 b_2, a_2 b_1) \in \ints$. Then, there exists integers $\alpha$ and $\beta$ such that
\begin{equation}
    \alpha m'_2 - \beta k' = \gcd(k', m'_2) = 1 \, .
\end{equation}

The ambiguity in the definition of $\alpha$ and $\beta$ is discussed in the proof of the following proposition.

\begin{proposition}
    \label{prop:SL3-non-abelian-Whittaker-from-abelian}
    Let $k \in \rats^\times$ with $\alpha, \beta, k'$ and $m'_2$ defined as above. Then
    \begin{equation}
        \label{eq:SL3-non-abelian-Whittaker-from-abelian}
        W_{\psi_Z}^{(k)}(\chi, g) = \sum_{m_1, m_2 \in \rats} W_{\psi_N}^{(m_1, d)}\big(\chi, l g\big) \qquad l = 
        \begin{psmallmatrix}
            \alpha & \beta & 0 \\
            k' & m'_2 & 0 \\
            0 & 0 & 1 
        \end{psmallmatrix} \in SL(3, \ints)
    \end{equation}
    where $g = (g_\infty, g_2, g_3, \ldots)$ is an arbitrary element of $G(\ads)$.
\end{proposition}

We will consider the restriction $g = (g_\infty, \id, \id, \ldots)$ giving integer charges in proposition~\ref{prop:SL3-non-abelian-expansion}. By $W_{\psi_N}^{(m_1, d)}$ we mean the Whittaker coefficients on $N$ given by $\psi_N$ with instanton charges $m_1$ and $d$ for the simple roots, which were calculated in \eqref{SL3gen}, \eqref{eq:SL3-generic-archimedean} and \eqref{eq:SL3-degenerate}. 

\begin{proof}
    
    To show \eqref{eq:SL3-non-abelian-Whittaker-from-abelian}, let first $l$ be defined as in that equation. We can expand $W_{\psi_Z}^{(k)}(\chi, g)$ further as    
        \begin{align}
            W_{\psi_Z}^{(k)}(\chi, g) 
            = \intl_{\rats \bs \ads}\!\! E(\chi,
            \begin{psmallmatrix}
                1 & 0 & z \\
                0 & 1 & 0 \\
                0 & 0 & 1
            \end{psmallmatrix}
            g) e^{-2\pi i k z} \, dz
            = \sum_{m_2 \in \rats} \, \intl_{(\rats \bs \ads)^2} \!\!\!\! E(\chi,
            \begin{psmallmatrix}
                1 & 0 & z \\
                0 & 1 & x_2 \\
                0 & 0 & 1
            \end{psmallmatrix}
            g) e^{-2\pi i (k z + m_2 x_2)} \, dz dx_2 . 
        \end{align}
    Using the automorphic invariance of $E(\chi, g)$ we can make the following conjugation with $l$
    \begin{equation}
        \begin{split}
            W_{\psi_Z}^{(k)}(\chi, g)
            &= \sum_{m_2} \intl_{(\rats \bs \ads)^2} E(\chi, l
            \begin{psmallmatrix}
                1 & 0 & z \\
                0 & 1 & x_2 \\
                0 & 0 & 1
            \end{psmallmatrix}
            l^{-1} lg) e^{-2\pi i (k z + m_2 x_2)} \, dz dx_2 \\
            &= \sum_{m_2} \intl_{(\rats \bs \ads)^2} E(\chi,
            \begin{psmallmatrix}
                1 & 0 & -d(x_2 - \alpha (k z + m_2 x_2)/d)/k \\
                0 & 1 & (k z + m_2 x_2)/d \\
                0 & 0 & 1
            \end{psmallmatrix}
            l g) e^{-2\pi i (k z + m_2 x_2)} \, dz dx_2 \\
            &= \sum_{m_2} \intl_{(\rats \bs \ads)^2} E(\chi,
            \begin{psmallmatrix}
                1 & 0 & x_2 \\
                0 & 1 & x_3 \\
                0 & 0 & 1
            \end{psmallmatrix}
            l g) e^{-2\pi i d x_3} \, dx_2 dx_3 \, ,  
        \end{split}
    \end{equation}
    where we have made the substitution $(k z + m_2 x_2)/d \to x_3$ and then $-d(x_2 - \alpha x_3)/k \to x_2$ leaving the integration domain the same. According to \eqref{eq:p-adic-measure} and \eqref{eq:adelic-norm-rationals} the measure is also unchanged. We note that the ambiguity in $\alpha$ simply results in an extra shift in the periodic variable $x_2$. 
    
    We expand one more time
    \begin{equation}
        \begin{split}
            W_{\psi_Z}^{(k)}(\chi, g) 
            &= \sum_{m_1, m_2 \in \rats} \, \intl_{(\rats \bs \ads)^3} E(\chi,
            \begin{psmallmatrix}
                1 & x_1 & 0 \\
                0 & 1 & 0 \\
                0 & 0 & 1
            \end{psmallmatrix}
            \begin{psmallmatrix}
                1 & 0 & x_2 \\
                0 & 1 & x_3 \\
                0 & 0 & 1
            \end{psmallmatrix}
            l g) e^{-2\pi i (m_1 x_1 + d x_3)} \, d^3x \\
            &= \sum_{m_1, m_2} \intl_{(\rats \bs \ads)^3} E(\chi,
            \begin{psmallmatrix}
                1 & x_1 & x_2 \\
                0 & 1 & x_3 \\
                0 & 0 & 1
            \end{psmallmatrix}
            l g) e^{-2\pi i (m_1 x_1 + d x_3)} \, d^3x \\  
            &=
            \sum_{m_1, m_2} W_{\psi_N}^{(m_1, d)}(\chi, lg)      
        \end{split}
    \end{equation}
    where we, in the second step, have made the substitution $x_2 + x_1 x_3 \to x_2$. 
\end{proof}

Note that when inserting $g = (g_\infty, \id, \id, \ldots) \in G(\ads)$ into \eqref{eq:SL3-non-abelian-Whittaker-from-abelian}, reducing the adelic Eisenstein series on the left hand side to the real Eisenstein series, the arguments on the right hand side become non-trivial at the finite places. To be able to use the expressions for $W_{\psi_N}$ above which require trivial arguments at the finite places, we need to factor out these effects.

\begin{proposition}
    \label{prop:SL3-non-abelian-expansion}
    Let $\tau = u_1 + i v_1^2 / v_2 \in \UHP$ and 
    \begin{equation}
        \label{eq:SL3-non-abelian-corollary-defs}
        \gamma =
        \begin{psmallmatrix}
            \alpha & \beta \\
            k' & m'_2
        \end{psmallmatrix} \in SL(2, \ints) \qquad
        \gamma(\tau) = \frac{\alpha \tau + \beta}{k' \tau + m'_2} \qquad
        a'_{\Im \gamma(\tau)} =
        \begin{psmallmatrix}
            v'_1 & 0 & 0 \\
            0 & v'_2 / v'_1 & 0 \\
            0 & 0 & 1/v'_2 
        \end{psmallmatrix}
    \end{equation}
    with $v'_1 = \sqrt{ v'_2 \Im \gamma(\tau)}$ and $v'_2 = v_2$, and $g = (g_\infty, \id, \ldots) \in G(\ads)$. Then $W_{\psi_Z}^{(k)}$ is non-vanishing only for $k \in \ints$ for which
    \begin{equation}
        \label{eq:SL3-non-abelian-expansion}
        W_{\psi_Z}^{(k)}\big(\chi, (g_\infty, \id, \ldots)\big) = \sum_{m_1, m_2 \in \ints} W_{\psi_N}^{(m_1, d)}\big(\chi, (a'_{\Im \gamma(\tau)}, \id, \ldots)\big) e^{-2\pi i(m_1 \Re \gamma(\tau) + m_2 u_2 + k z)}
    \end{equation}
\end{proposition}

Note that the sums over rationals have collapsed to sums over integers and that
\begin{equation}
    l =
    \begin{psmallmatrix}
        \gamma & 0 \\
        0 & 1
    \end{psmallmatrix} \, .
\end{equation}

This proves the results of \cite{VT} and \cite{Proskurin_1984} reviewed in \cite{Pioline:2009qt, Persson:2010ms} with only a few manipulations using the compact framework of adelic automorphic forms.

\begin{proof}
In~\eqref{eq:SL3-non-abelian-Whittaker-from-abelian} the argument for $W_{\psi_N}^{(m_1, d)}$ is $lg = (lg_\infty; l, l, \ldots)$ and since $W_{\psi_N}(\chi, n' a' k') = \psi_N(n') W_{\psi_N}(\chi, a')$, we factorise $lg$ at the archimedean and non-archimedean places into their respective Iwasawa decompositions. We have that $l \in SL(3, \ints)$ which makes the $p$-adic Iwasawa decomposition trivial with $l \in K_p$. This was the reason for choosing $l$ on this particular form.

We then use the following relation, similar to \eqref{eq:urnamified-vanishing-condition}, to obtain conditions for $m_1$ and $m_2$. For $\hat n = (\id; \hat n_2, \hat n_3, \ldots)$ with $\hat n_p \in N(\ints_p) \subset K_p$ we have that $\hat n \in K_\ads$ and
\begin{equation}
    \label{eq:p-adic-Whittaker-vanishing}
    W_{\psi_N}(\chi, a) = W_{\psi_N}(\chi, a \hat n) = W_{\psi_N}(\chi, a \hat n a^{-1} a) = \psi_{N}(a \hat n a^{-1}) W_{\psi_N}(\chi, a)
\end{equation}
which requires that $\psi_{N}(a \hat n a^{-1}) = 1$ for $W_{\psi_{N}}(\chi, a)$ to be non-vanishing. 

Specifically, for $W_{\psi_N}^{(m_1, d)}(\chi, a)$ with $a = (a_\infty; \id, \id, \ldots)$ we require that
\begin{equation}
    \label{eq:SL3-Whittaker-vanishing}
    1 = \psi_N(a \hat n a^{-1}) = \psi_{N,\infty}(\id) \prod_{p<\infty} \psi_{N,p}(\hat n_p) = \exp\left(-2\pi i \sum_{p<\infty} [m_1 u_1 + d u_2]_p\right) 
\end{equation}
for all $u_1, u_2 \in \ints_p$ where
\begin{equation}
    \hat n_p =
    \begin{psmallmatrix}
        1 & u_1 & z \\
        0 & 1 & u_2 \\
        0 & 0 & 1
    \end{psmallmatrix} \, .
\end{equation}

This implies that $\sum_{p<\infty}[m_1]_p \in \ints$ and $\sum_{p<\infty}[d]_p \in \ints$, which, according to proposition \ref{prop:p-adic-fractional}, gives that $m_1, d \in \ints$. That $d$ is integer means that, for all primes $p$
\begin{equation}
    1 \geq \abs{d}_p = \frac{\max(\abs{a_1b_2}_p, \abs{a_2b_1}_p)}{\abs{b_1b_2}_p} = \max(\abs{k}_p, \abs{m_2}_p)
\end{equation}
according to \eqref{eq:p-adic-abs-gcd}, and hence, that $k$ and $m_2$ are also integers.

For the archimedean place we have the Iwasawa decomposition
\begin{equation}
    \begin{gathered}
        l g_\infty = l n_\infty a_\infty k_\infty = l
        \begin{psmallmatrix}
            1\vphantom{a_1/b_1} & u_1 & z \\
            0 & 1\vphantom{a_2/b_2} & u_2 \\
            0 & 0 & 1\vphantom{a_3/b_3}
        \end{psmallmatrix}
        \begin{psmallmatrix}
            v_1\vphantom{a_1/b_1} & 0 & 0 \\
            0 & v_2 / v_1 & 0 \\
            0 & 0 & 1 / v_2
        \end{psmallmatrix}
        k_\infty
        = n'_\infty a'_\infty k'_\infty \qquad \text{with} \\[1em]
        \begin{aligned}
            u'_1 &= - \frac{d^2 (m_2 + k u_1) v_2^2}{k^3 v_1^4 + k(m_2 + k u_1)^2 + v_2^2}  + \frac{d \alpha}{k} & \qquad u'_2 &= \frac{m_2 u_2 + k z}{d} \\[1em]
            v'_1 &= \frac{v_1 v_2 d}{\sqrt{k^2 v_1^4 + (m_2 + k u_1)^2 v_2^2}} & \qquad v'_2 &= v_2
        \end{aligned}
    \end{gathered}
\end{equation}

We define $\tau = u_1 + i v_1^2/v_2 \in \UHP$, which, under the $l$-translation on $g_\infty$ above, transforms as $\tau \to \tau'$ with
\begin{equation}
    \tau' = u_1' + i \frac{(v'_1)^2}{v'_2} = \gamma(\tau) \qquad \gamma =
    \begin{psmallmatrix}
        \alpha & \beta \\
        k' & m'_2
    \end{psmallmatrix} \in SL(2, \ints) \, .
\end{equation}
Putting it all together we obtain for $k \in \ints^\times$ and $g = (g_\infty; \id, \id, \ldots)$
\begin{equation}
    \begin{split}
        \MoveEqLeft
        W_{\psi_Z}^{(k)}(\chi, (g_\infty; \id, \ldots)) =\\
        &=\sum_{m_1, m_2 \in \ints} W_{\psi_N}^{(m_1, d)}(\chi, a') \psi_N^{(m_1, d)}(n') \\
       &= \sum_{m_1, m_2 \in \ints} W_{\psi_N}^{(m_1, d)}(\chi, (a'_\infty; \id, \ldots)) \Big(\prod_{p<\infty} \psi_{N,p}^{(m_1, d)}(n'_p) \Big) \psi_{N,\infty}^{(m_1,d)}(n'_\infty) \\
       &= \sum_{m_1, m_2 \in \ints} W_{\psi_N}^{(m_1, d)}(\chi, (a'_{\Im \gamma(\tau)}; \id, \ldots)) e^{-2\pi i(m_1 \Re \gamma(\tau) + m_2 u_2 + k z)}
    \end{split}
\end{equation}
where $a'_{\Im \gamma(\tau)}$ is defined in \eqref{eq:SL3-non-abelian-corollary-defs}.
\end{proof}

The remaining Whittaker coefficients on $N$ with trivial arguments at the non-archimedean places were computed in sections \ref{sec:SL3-generic-Whittaker-coefficients} and \ref{sec:SL3-degenerate-Whittaker-coefficients}.

\begin{remark}
    The physical intepretation of the $SL(2, \ints)$ action described by $\gamma$ is described by S-duality of type IIB string theory compactified on a Calabi-Yau threefold \cite{Pioline:2009qt, RoblesLlana:2006is, PerssonAuto}.
    In this setting the parameters $z$ and $u_2$ are scalar fields sourced by D5- and NS5-branes with charges $m_2$ and $k$ (more generally denoted by $p$ and $q$). The branes form bound states that are often reffered to as $(p, q)$ 5-branes. 
    The two scalar fields transform as an $SL(2)$-doublet under S-duality mirrored by their transformation under $g \to lg$ and the charges $p$ and $q$, which appear as $m_2$ and $k$ in \eqref{eq:SL3-non-abelian-expansion}, break the classical $SL(2, \reals)$ symmetry of the supergravity theory to the discrete $SL(2,\ints)$ symmetry of the quantum corrected effective action described by Eisenstein series.
    In short, this tells us that if we can compute the effects from a $(p, 0)$ 5-brane, the results for any $(p, q)$ 5-brane follow from S-duality which is mirrored in the sum over matrices $\gamma$ in \eqref{eq:SL3-non-abelian-expansion}.
\end{remark}

\section[The Casselman--Shalika formula and Langlands duality*]{The Casselman--Shalika formula and\\ Langlands duality*}
\label{sec:CSLD}

In this section we provide an alternative view on the Casselman--Shalika formula~\eqref{CSus}. For the present analysis it is useful to separate out the modulus character contribution $a^\rho=\delta^{1/2}(a)$ in  formula~\eqref{CSus} and write
\begin{align}
a^{w\lambda+\rho} = a^\rho a^{w\lambda}=\delta^{1/2}(a) a^{w\lambda}.
\end{align}
Let $\psi$ be an unramified character on $N$. The Casselman--Shalika formula (\ref{CSus})  for the $p$-adic spherical Whittaker function on $\rats_p$ evaluated at $a\in A(\rats_p)$ is
\begin{align}
\label{CSusR}
W^{\circ}(\lambda,a) = \frac{1}{\zeta(\lambda)} \delta^{1/2}(a) \sum_{w\in \Weyl}  w\left(\frac{|a^{\lambda}|}{\prod_{\alpha>0} (1-p^{\langle\lambda|\alpha\rangle})}\right)
\end{align}
with 
\begin{align}
\label{zetacoeff2}
\frac{1}{\zeta(\lambda)} = \prod_{\alpha>0} \left(1-p^{-1} p^{-\langle\lambda|\alpha\rangle}\right).
\end{align}

The sum over the Weyl group in (\ref{CSusR}) resembles closely the Weyl character formula (\ref{WCF2}) for highest weight modules. In order to make this resemblance exact, we compare with the rewritten character formula in (\ref{charfn}) that we reproduce here for convenience:
\begin{align}
\mathrm{ch}_\Lambda(b) 
= \sum_{w\in\mathcal{W}} w\left(\frac{b^{\Lambda}}{ \prod_{\beta>0} (1-b^{-\beta})}\right),
\end{align}
where $\beta$ runs over the roots of the group whose representation is being constructed and $b$ is an element of its Cartan torus.

 An important first observation now is that because of the way $\lambda$ appears in the numerator and in the denominator of~\eqref{CSusR} the comparison can only work if the character we are trying to match onto is one of the \emphindex[Langlands dual group|textbf]{Langlands dual group} ${}^LG$, or $L$-group for short, which is a complex algebraic group canonically associated to $G$~\cite{LanglandsWeilLetter}. See also sections~\ref{sec_Lgroup} and~\ref{sec:Langlands} for more details. The $L$-group is obtained by interchanging roots and co-roots~\cite{KnappLanglandsProgram}, see also~\cite{Goddard:1976qe} for a realisation in physics. The root systems of $G$ and ${}^LG$ are in bijection and the two groups have isomorphic Weyl groups. 

Denoting the roots of the Langlands dual group by $\alpha^\vee$ instead of $\beta$, we are therefore looking for an element $b$ of the dual torus ${}^LA$ such that $|b^{-\alpha^\vee}| = p^{\langle\alpha|\lambda\rangle}$. This condition fixes uniquely an element $b=a_\lambda\in {}^LA$, where we emphasise that the particular element depends on $\lambda$. To ensure that the numerator matches the character of an irreducible highest weight module $V_\Lambda$ of ${}^LG$ we also need to evaluate~\eqref{CSusR} at a very specific point $a\equiv a_\Lambda$ of $A(\rats_p)$. This element $a_\Lambda$ is fixed by the requirement that the following (duality) relation hold
\begin{align}
\label{LDual}
a^\Lambda_\lambda = a^\lambda_\Lambda,
\end{align}
where the left-hand side derives from evaluating the character $\mathrm{ch}_\Lambda$ at the place $a_\lambda$ and the right-hand side is what one obtains by evaluating the spherical Whittaker function at the special point $a_\Lambda\in A$. 

We observe that $\Lambda$ parametrises points in the space of co-roots $\mf{h}$ of $G$.  By contrast, $\lambda$ is an element of the space of roots (or weights) $\mf{h}^*$ of $G$ from the start, so that one has to consider $a_\lambda$ as an element of the \emphindex[torus!dual]{dual torus} ${}^LA$ of the \emphindex[Langlands dual group|textbf]{Langlands dual group} ${}^LG$. Putting everything together we can write the spherical Whittaker function evaluated at $a_\Lambda$ in terms of the character of the highest weight representation $V_\Lambda$ of ${}^LG$ as
\begin{align}
\label{CSWCF}
W^\circ (\lambda, a_\Lambda)= \left\{\begin{array}{cl} \frac{1}{\zeta(\lambda)}
 \delta^{1/2}(a_\Lambda) \mathrm{ch}_\Lambda(a_\lambda) & \textrm{if $\Lambda$ a dominant integral weight of ${}^LG$},\\
0 & \textrm{otherwise}.
\end{array}\right.
\end{align}
The vanishing for non-dominant weights $\Lambda$ of ${}^LG$ is a consequence of the vanishing properties of Whittaker functions discussed in section~\ref{VP}. 

To summarise the main result of this section: Local spherical Whittaker functions for a principal series representations parametrised by a weight $\lambda$ of $G$ and evaluated at special points $a_\Lambda$ associated with dominant weights $\Lambda$ of the Langlands dual group ${}^LG$ can be evaluated in terms of the character $\mathrm{ch}_\Lambda$ of the irreducible highest weight $V_\Lambda$ of ${}^LG$ evaluated at a point $a_\lambda$ determined by the parameter of the principal series.

The parameter $a_\lambda\in {}^LA$ is called the \emphindex[Satake parameter]{Satake--Langlands parameter} of the principal series representation of $G(\rats_p)$ determined by the weight $\lambda$ and we will come back to it in a slightly different guise in section~\ref{sec_Lgroup}. We also note that the element $a_\Lambda\in A(\rats_p)$  actually corresponds to an equivalence class $A(\rats_p)/A(\ints_p)$ due to sphericality (right $K(\rats_p)$ invariance) of the Whittaker function. 

For the case of $GL(n,\rats_p)$ one has ${}^LG=GL(n,\cx)$. If one considers the case when $\Lambda$ is the highest weight of the fundamental $n$-dimensional representation $\cx^n$, then the character $\mathrm{ch}_\Lambda$ is given by the degree $n$ Schur polynomial $S_n$\index{Schur polynomial}
\beq
W^{\circ}(\lambda, a_\Lambda)=\frac{1}{\zeta(\lambda)} \delta^{1/2}(a_\Lambda) S_n(\alpha_1, \dots, \alpha_n), \qquad G=GL(n, \mathbb{Q}_p).
\eeq
Here, $\lambda$ is thought of as the diagonal matrix $\lambda=\mathrm{diag}(\alpha_1,\ldots,\alpha_n)$.
This formula for $GL(n)$ was first proven by Shintani in 1976 \cite{Shintani}, and it was subsequently generalised by Casselman--Shalika in 1980 to (\ref{CScharacter}) which holds for any $G$. Remarkably the general formula was in fact conjectured by Langlands already in 1967 in a letter to Godement \cite{LanglandsGodement}, a fact that was apparently unknown to Casselman and Shalika at the time of their proof \cite{CasselmanLGroup}.

\begin{example}[$SL(2,\rats_p)$ spherical Whittaker function and $SL(2,\cx)$ characters]
For the case $SL(2,\rats_p)$ the spherical Whittaker function for unramified $\psi$ was given explicitly in (\ref{SL2Whitt}) for $\lambda=(2s-1)\rho$ and general $a=v^{H_\alpha}$ as
\begin{align}
\label{SL2W2}
W^\circ(\lambda,a) = \gamma_p(v^2) (1-p^{-2s}) \frac{|v|^{-2s+2}-p^{-2s+1}|v|^{2s}}{1-p^{-2s+1}}.
\end{align}
In order to verify the expression (\ref{CSWCF}) we need to evaluate them at the special values $a_\Lambda$ where $\Lambda=NH_\alpha/2$ is a dominant integral weight of ${}^LSL(2,\rats_p)=SL(2,\cx)$ for $N\in \ints_{\geq 0}$. This means $v^2=p^N$ and the Whittaker function evaluates to
\begin{align}
\label{SL2WaLambda}
W^\circ(\lambda,a_\Lambda) =(1-p^{-2s}) \frac{p^{\frac{N}{2}(2s-2)}-p^{-2s+1-Ns}}{1-p^{-2s+1}}.
\end{align}
The Whittaker function vanishes if $N$ is not in $\ints_{\geq 0}$ because of the factor $\gamma_p(v^2)$.

Let us now determine the right-hand side of~\eqref{CSWCF}. For $\Lambda=\frac{N}{2} H_\alpha$,
the character of the $(N+1)$-dimensional highest weight representation of ${}^LSL(2,\rats_p)\cong PSL(2,\cx)$ is
\begin{align}
\mathrm{ch}_\Lambda = e^{NH_\alpha/2}+ e^{(N-2)H_\alpha/2} +\ldots + e^{-N H_\alpha/2} = \frac{e^{-NH_\alpha/2}-e^{(N+2)H_\alpha/2}}{1-e^{H_\alpha}}.
\end{align}
This has to be evaluated at $a_\lambda= p^\lambda=p^{(2s-1)\Lambda_\alpha}$ which leads to
\begin{align}
\mathrm{ch}_\Lambda(a_\lambda) = \frac{p^{N(2s-1)/2} - p^{-(N+2)(2s-1)/2}}{1-p^{-2s+1}},
\end{align}
where we recall that the $p$-adic characters are evaluated with the $p$-adic norm such that for instance $e^{H_\alpha}(a_\lambda) = |p^{2s-1}| = p^{-2s+1}$.
For $v^2=p^N$, the modulus character evaluates to $\delta^{1/2}(a_\Lambda) = |p^{N/2}|=p^{-N/2}$ and one also has $\frac{1}{\zeta(\lambda)}=1-p^{-2s}$ from (\ref{zetacoeff2}). Putting everything together in (\ref{CSWCF}) leads to
\begin{align}
\label{CS2Lang}
W^\circ(\lambda,a_\Lambda)=
(1-p^{2s}) p^{-N/2} \frac{p^{N(2s-1)/2} - p^{-(N+2)(2s-1)/2}}{1-p^{-2s+1}} = (1-p^{-2s}) \frac{p^{N(2s-2)/2} - p^{-2s+1-sN}}{1-p^{-2s+1}}
\end{align}
which equals (\ref{SL2WaLambda}).

\end{example}

\begin{remark}
Using formula~\eqref{unramified2generic} we can also reinterpret~\eqref{CSWCF} in terms of a Whittaker function for the twisted character 
\begin{align}
\psi_\Lambda(n):= \psi(a_\Lambda n a_\Lambda^{-1})
\end{align}
as
\begin{align}
W^\circ_{\psi_\Lambda} (\lambda,\id) = a_\Lambda^{-\wlong \lambda -\rho} W^\circ(\lambda, a_\Lambda) 
= \frac{1}{\zeta(\lambda)}a_\lambda^{-\wlong\Lambda}  \mathrm{ch}_\Lambda(a_\lambda),
\end{align}
where we used~\eqref{LDual}.
\end{remark}

\chapter{Working with Eisenstein series}
\label{ch:working}

After having developed the formal theory of Eisenstein series and their Fourier expansion in the previous chapters we would like to discuss Eisenstein series from a more practical point of view in this chapter. In concrete examples this typically means obtaining as much information as possible for a particular Eisenstein series $E(\chi,g)$, that is a particular given $\chi$. Many of the general theorems either simplify for such a $\chi$ or have to be evaluated with much care as $E(\chi,g)$ might be divergent for the chosen $\chi$. This chapter deals with developing methods for addressing these issues. In particular, we exhibit methods for efficiently evaluating the constant term formula~\eqref{LCF} and  formula~\eqref{degW1} for the Whittaker coefficients of a given Eisenstein series. We will also discuss the pole structure of Eisenstein series (as a function of $\chi$) in examples, their residues as well as different normalisations. In this chapter, the emphasis is on illustrating different methods through many examples; for proofs of general statements we will typically refer to the appropriate literature.

Many of the properties of Eisenstein series are controlled by the completed Riemann zeta function whose properties we briefly recall.
\begin{proposition}[Properties of completed Riemann zeta function]
\label{RZprop}
\index{Riemann zeta function!properties}
As a function of $s\in \cx$, the completed Riemann zeta function $\xi(s)=\pi^{-s/2} \Gamma(s/2) \zeta(s)$ has simple poles at $s=0$ and $s=1$ with residues $-1$ and $+1$, respectively. It is non-zero everywhere else. It satisfies the functional relation $\xi(s)=\xi(1-s)$.
\end{proposition}
\begin{proof} 
The first statements follow directly from the definition and the properties of gamma and zeta functions. The functional relation was shown originally by Riemann using analytic continuation~\cite{RiemannZeta}.
\end{proof}

\section{The \texorpdfstring{$SL(2,\reals)$}{SL(2,R)} Eisenstein series as a function of \texorpdfstring{$s$}{s}}

We begin with the $SL(2,\reals)$ Eisenstein series $E(s,z)$ that was analysed in great detail in chapter~\ref{ch:SL2-fourier} with its complete Fourier expansion given in theorem~\ref{SL2AESexp}. We repeat the result here for convenience:
\begin{align}
\label{SL2FE9}
E(s,z) = y^s + \frac{\xi(2s-1)}{\xi(2s)} y^{1-s} + \frac{2}{\xi(2s)} y^{1/2} \sum_{m\neq 0 }  |m|^{s-1/2} \sigma_{1-2s}(m) K_{s-1/2}(2\pi|m| y) e^{2\pi i m x},
\end{align}
where $z=x+i y$ is an element of the upper half plane $\UHP=SL(2,\reals)/SO(2)$. The original definition of $E(s,z)$ only converged for $\Re(s)>1$ but by virtue of the functional relation (cf. theorem~\ref{SL2AESexp})
\begin{align}
\label{funcSL29}
E(s,z) = \frac{\xi(2s-1)}{\xi(2s)} E(1-s,z)
\end{align}
or through analytic continuation of the Fourier expansion~\eqref{SL2FE9} one can define $E(s,z)$ for almost all complex $s$. We restrict our discussion to real $s$ for simplicity. 

\subsection{Limiting values in original normalisation}
\label{SL2limits}
From the explicit form~\eqref{SL2FE9} one sees that special things might happen for the values $s=0$, $s=\tfrac12$ and $s=1$. All of them are outside the original domain of convergence. Let us note that the region $0 \leq \Re(s) \leq 1$ is often called the critical strip.
\begin{enumerate}[leftmargin=5em, style=nextline]
\item[\textbullet\ $s=0$: ] This is the limit where the inducing character $\chi_s(z) = y^s$ becomes trivial. Taking the limit in the expression~\eqref{SL2FE9} for the Fourier expansion one also sees that all terms go to zero except for the first. This is due to the factors $\frac1{\xi(2s)}$ that vanish linearly for $s\to0$ while everything else stays bounded. The proper limiting behavoiur is therefore
\begin{align}
E(s,z) = 1+ O(s).
\end{align}
The constant value $1$ could have been expected from the triviality of the inducing character but the definition in terms of a~\emph{Poincar\'e sum} is ill-defined. Only after analytic continuation of the sum one obtains the constant $E(0,z)=1$.

Representation theoretically, the function $E(s,z)$ in the limit $s\to 0$ belongs to the trivial representation of $SL(2,\reals)$.  

\item[\textbullet\ $s=\tfrac12$: ] Inspection of the Fourier expansion~\eqref{SL2FE9} shows that the non-zero Fourier modes disappear in this limit due to the $\frac{1}{\xi(2s)}$ prefactor. For the constant terms one has to take the limit of the quotient of $\xi$-functions which is found to be $-1$ and the two contributions to the constant term cancel, leading to
\begin{align}
E(s,z) = 0 + O\left(s-\tfrac12\right).
\end{align}
The first order term is a member of the principle series. It is on the critical line and is therefore almost unitary.

\item[\textbullet\ $s=1$: ] This is the most interesting case. The Fourier expansion~\eqref{SL2FE9} shows that the second constant term diverges in the limit $s\to1$ while all other terms remain finite. The residue at the simple pole can be calculated easily from the completed Riemann zeta functions:
\begin{align}
E(s,z) = \frac{3}{\pi(s-1)} + O\left((s-1)^0\right).
\end{align}
The residue is therefore a constant function and is therefore also of the same type as the limit $s\to 0$ discussed above. This is not surprising since the functional relation~\eqref{funcSL29} relates the values $s=0$ and $s=1$ and one sees that the prefactor introduces the additional pole. Representation theoretically, the residue of the series $E(s,z)$ at the simple pole $s=1$ belongs to the trivial representation of $SL(2,\reals)$.

The term at order $(s-1)^0$ can also be evaluated from the Fourier expansion using the fact that the modified Bessel function $K_{1/2}$ has an exact asymptotic expansion in terms of a simple exponential. One finds
\begin{align}
\label{Kron1}
E(s,z) &= \frac{3}{\pi(s-1)} - \frac{6}{\pi} \bigg( -\frac{\pi}{6} y + \log(4\pi\sqrt{y}) -12 \log A \\
&\quad\quad - \sum_{m>0} \sigma_{-1}(m) e^{2\pi i m (x+i y)} - \sum_{m>0} \sigma_{-1}(m) e^{2\pi i m (x-i y)}
\bigg)+ O(s-1).\nn
\end{align}
Here, $A$ is the \emphindex{Glaisher--Kinkelin constant} that satisfies $\log A = \frac1{12}-\zeta'(1)$. The expression can be rewritten by using the \emphindex{Dedekind $\eta$ function}
\begin{align}
\eta(z) = q^{1/24} \prod_{n=1}^\infty (1-q^n),
\end{align}
where $q=e^{2\pi i z} =e^{2\pi i (x+iy)}$ on the right-hand side. From the product formula for $\eta(z)$ one concludes
\begin{align}
\log \eta(z) &= \frac1{24} \log q + \sum_{n>0} \log (1-q^n) = \frac{\pi i}{12} (x+iy) - \sum_{n>0} \sum_{k>0} k^{-1} q^{k n} \nn\\
&= \frac{\pi i}{12} (x+iy) - \sum_{m>0} \sum_{d|m} d^{-1} q^{m} = \frac{\pi i}{12} (x+iy) - \sum_{m>0} \sigma_{-1}(m)  q^{m}.
\end{align}
The $s$-independent term in~\eqref{Kron1} can therefore be written as
\begin{align}
&\quad- \frac{6}{\pi} \left( -\frac{\pi}{6} y + \log(4\pi\sqrt{y}) -12 \log A 
 - \sum_{m>0} \sigma_{-1}(m) q^m - \sum_{m>0} \sigma_{-1}(m) \bar{q}^m \right)\nn\\
&=-\frac{6}{\pi} \Big(-12\log A + \log(4\pi) +\log \left(\sqrt{y} |\eta(z)|^2 \right)\Big),
\end{align}
leading to
\begin{align}
\label{Kron2}
E(s,z) &= \frac{3}{\pi(s-1)} + \frac{6}{\pi} \Big(12\log A - \log(4\pi) -\log \left(\sqrt{y} |\eta(z)|^2 \right)\Big)+ O(s-1).
\end{align}
This formula is known as the (first) \emphindex{Kronecker limit formula}. Even though neither $\eta(z)$ nor $|\eta(z)|^2$ are $SL(2,\ints)$ invariant, the particular combination appearing in this expression is invariant. 
\end{enumerate}

\begin{remark}
The particular combination of constants in the Kronecker limit formula~\eqref{Kron2} depends on the way the Eisenstein series is normalised. The formula is more commonly stated for the $SL(2,\ints)$ invariant lattice sum (cf.~\eqref{Eisenintro}) for which one finds
\begin{align}
\sum_{(c,d)\in \ints^2\atop (c,d)\neq (0,0)} \frac{y^s}{|cz+d|^{2s}} &= 2\zeta(2s) E(s,z)\\
& = \frac{\pi}{s-1} + 2\pi \Big(\gammaE - \log(2) -\log \left(\sqrt{y} |\eta(z)|^2 \right)\Big)+ O(s-1)
\end{align}
if one uses the following relation between the Glaisher--Kinkelin constant $A$ and the \emphindex{Euler--Mascheroni constant} $\gammaE$: $12\log A - \log (4\pi) = \gammaE -\log 2 - \frac{\zeta'(2)}{\zeta(2)}$.
\end{remark}

\subsection{Weyl symmetric normalisation}

The functional relation~\eqref{funcSL29} suggests to define a \emphindex[Eisenstein series!completed]{completed Eisenstein series} in analogy with the completed Riemann zeta function by the definition
\begin{align}
E^\compl(s,z) = \xi(2s) E(s,z).
\end{align}
This then has the simple property that
\begin{align}
E^\compl(s,z) = E^\compl(1-s,z)
\end{align}
and we call this the \emphindex{Weyl symmetric normalisation} as it yields a function invariant under Weyl transformations acting on the character. Indeed, the non-trivial Weyl reflection $w$ of $SL(2,\reals)$ acts on the weight $\lambda_s=(2s-1)\rho$ by
\begin{align}
w\lambda_s = -(2s-1) \rho = (2(1-s)-1)\rho = \lambda_{1-s}
\end{align}
and so exchanges $s$ and $1-s$. This was of course already used and apparent in the constant terms in~\eqref{SL2FE9}.

Since the normalising factor has poles and zeroes of its own, the discussion of the behaviour of $E^\compl(s,z)$ as a function of $s$ is slightly changed from the one above. More precisely, the completed function $E^\compl(s,z)$ has a simple poles at $s=0$ and $s=1$, whereas it has a non-trivial limit for $s=\tfrac12$. Representation theoretically, $E^\compl(\tfrac12, z)$ belongs to the principal series.

\section{Properties of Eisenstein series}

The behaviour of the $SL(2,\reals)$ Eisenstein series at the special values of $s$ above was completely controlled by the constant terms. This is a general feature due to the holomorphy of the Fourier coefficients, see proposition~\ref{prop:hol}. As we have full control of the constant terms thanks to the Langlands constant term formula (theorem~\ref{LCFthm}), we can in principle completely determine the behaviour of an Eisenstein series $E(\lambda,g)$ on a group $G(\reals)$ as a function of $\lambda$. As the number of constant terms is generically equal to the order of the Weyl group $\Weyl$ of $G$ this can be quite tedious due to the large number of terms that have to be considered. In section~\ref{EvalLCF}, we will present a method that makes the problem more tractable for the case of non-generic $\lambda$ when the Eisenstein series $E(\lambda,g)$ is not attached to the full principal series but to a degenerate principal series. The prime example of this is when it becomes a maximal parabolic Eisenstein series as defined in section~\ref{nonminEis}. Before focussing on these cases in section~\ref{EvalLCF}, we offer a few general and cautionary remarks.

\subsection{Validity of functional relation}

As Langlands showed in his seminal work~\cite{LanglandsFE}, the functional equation \eqref{funrel} repeated here for convenience
\begin{align}
\label{eq:FE9}
E(\lambda,g) = M(w,\lambda) E(w\lambda, g)
\end{align}
is valid for almost all $\lambda \in \mathfrak{h}^*(\cx)$. The exceptions are affine hyperplanes in the complex vector space $\mathfrak{h}^*(\cx)$. These affine hyperplanes are associated with poles and zeroes of the intertwining factor $M(w,\lambda)$. Since all poles and zeroes are of finite order, one can make sense of the functional relation even on these planes by treating also the Eisenstein series as meromorphic functions with finite order poles and singularities.

\begin{example}[Functional relation for $SL(2,\reals)$ Eisenstein series at a simple zero and pole]
For $SL(2,\reals)$ and general $\lambda=(2s-1)\rho$, the intertwining factor for the non-trivial Weyl element $w=\wlong$ is
\begin{align}
M(\wlong,\lambda) = \frac{\xi(2s-1)}{\xi(2s)}
\end{align}
and has a simple zero at $s=0$ and a simple pole at $s=1$. The functional relation~\eqref{eq:FE9} remains valid even at these places if one considers 
\begin{align}
E(s,z) &= 1 + s \hat{E}_0(z) + O\left(s^2\right),\nn\\
E(s,z) &= \frac{3}{\pi(s-1)} + \hat{E}_1(z) + O\left( s-1\right)
\end{align}
around $s=0$ and $s=1$, respectively. The expansion of the intertwining factor around these values is
\begin{align}
\frac{\xi(2s-1)}{\xi(2s)} = -\frac{\pi s}{3 } + O(s^2) = \frac{3}{\pi(s-1)} +  O\left((s-1)^0\right),
\end{align}
such that 
\begin{align}
E(s,z) = 1+ s \hat{E}_0(z) +O(s^2) &= \left( -\frac{\pi s}{3} + O(s^2)\right)\left( -\frac{3}{\pi s} + \hat{E}_1 + O(s)\right)\nn\\
&= 1-\frac{\pi s}{3} \hat{E}_1 + O(s^2)
\end{align}
and so the functional relation relates $\hat{E}_0(z)$ and $\hat{E}_1(z)$ (as well as all higher order terms).
\end{example}

Of interest are also fixed planes of the action of the Weyl group action. In these cases, the functional relation~\eqref{eq:FE9} constrains the Eisenstein series on the fixed plane. 

\begin{example}[Functional relation for $SL(2,\reals)$ Eisenstein series with $\lambda=0$]
\label{ex:SL2fixedplane}
For $SL(2,\reals)$ and $E(s,z)$ the fixed plane is $s=\tfrac12$, corresponding to $\lambda=\lambda_{1/2} = 0$. The intertwining factor at this place takes the value $M(w,0) = -1$ such that the functional relation implies
\begin{align}
E(\tfrac12,z) = - E(\tfrac12,z) \quad\Longrightarrow \quad E(\tfrac12,z) = 0,
\end{align}
consistent with the analysis in section~\ref{SL2limits}.
\end{example}

We consider also a few examples of functional relations for higher rank groups.
\begin{example}[Functional relation for $SL(3,\reals)$ Eisenstein series]
The most general Eisenstein series on $G=SL(3,\reals)$ is given by a weight
\begin{align}
\lambda_{s_1,s_2} = 2s_1\Lambda_1 + 2s_2 \Lambda_2 -\rho
\end{align}
that is parametrised by two complex parameters $s_1$ and $s_2$. The $\Lambda_i$ are as always the fundamental weights. We denote the corresponding character by $\chi_{s_1,s_2}(a) = a^{\lambda_{s_1, s_2}+\rho}$ and the Eisenstein series by
\begin{align}
E(s_1,s_2,g) = \sum_{B(\ints)\bs SL(3,\ints)} \chi_{s_1,s_2} (\gamma g).
\end{align}
The sum is absolutely convergent for $\Re(s_1)>1$ and $\Re(s_2)>1$~\cite{BumpSL3}. The Weyl group of $SL(3,\reals)$ is isomorphic to the symmetric group on three letters and hence consists of six elements and the constant terms were already given in~\eqref{SL3consts}. Denoting the fundamental reflections by $w_1$ and $w_2$ one finds that
\begin{subequations}
\begin{align}
w_1\lambda_{s_1,s_2} &= (1-2s_1)\Lambda_1 + 2(s_1+s_2-1)\Lambda_2 = \lambda_{1-s_1,s_1+s_2-\tfrac12},\\
w_2\lambda_{s_1,s_2} &= 
 \lambda_{s_1+s_2-\tfrac12,1-s_2}.
\end{align}
\end{subequations}
The other Weyl images can be obtained similarly. One functional relation is therefore
\begin{align}
E(s_1,s_2,g) &= M(w_1,\lambda_{s_1,s_2}) E(1-s_1,s_1+s_2-\tfrac12,g) 
=\frac{\xi(\langle\alpha_1 | \lambda_{s_1,s_2}\rangle)} {\xi(\langle\alpha_1 | \lambda_{s_1,s_2}\rangle+1)}  E(1-s_1,s_1+s_2-\tfrac12,g) \nn\\
&= \frac{\xi(2s_1-1)}{\xi(2s_1)} E(1-s_1,s_1+s_2-\tfrac12,g).
\end{align}
That this is a valid relation can be checked on the constant terms from~\eqref{SL3consts}. We can consider the limit $s_1\to \tfrac12$ to conclude
\begin{align}
E(\tfrac12,s_2,g) = - E(\tfrac12,s_2,g) \quad\Longrightarrow \quad
E(\tfrac12,s_2,g) = 0.
\end{align}
This is exactly as in the $SL(2,\reals)$ case in example~\ref{ex:SL2fixedplane} above. Again $s_1=\tfrac12$ corresponds to a fixed plane of a fundamental reflection and in this case one always obtains a vanishing Eisenstein series.
\end{example}

More involved examples are obtained for exceptional groups. These will play an important role in section~\ref{sec:outlook-strings} in the context of string theory.

\begin{figure}[t!]
\centering
\input{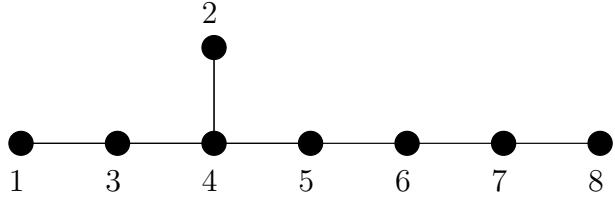}
\caption{\label{fig:E8dynk}\sl The Dynkin diagram of $E_8$ with labelling of nodes in the `Bourbaki convention'.}
\end{figure}

\begin{example}[Functional relation for $E_8(\reals)$ maximal parabolic Eisenstein series]
\label{ex:E8adj}

Consider a (maximal parabolic) Eisenstein series on $E_8$ with Dynkin diagram given in figure~\ref{fig:E8dynk}. For the weight
\begin{align}
\label{E8adwt}
\lambda_s= 2s \Lambda_8 -\rho
\end{align}
the associated character $\chi_s(a)= a^{\lambda_s+\rho} = a^{2s\Lambda_8}$ is invariant under the maximal parabolic subgroup with semi-simple part $E_7$. (Here, $\Lambda_8$ denotes as always the fundamental weight associated with node $8$.) We therefore have a family of maximal parabolic Eisenstein series
\begin{align}
E(s,g) \equiv E(\lambda_s, P, g) = \sum_{\gamma\in P(\ints)\backslash G(\ints)} \chi_s(\gamma g).
\end{align}
There are many functionally related Eisenstein series. The Weyl group of $E_8$ has order $|\Weyl(E_8)|=696\,729\,600$, but for the particular choice of parameter $\lambda_s$ in~\eqref{E8adwt} not all Weyl images of $\lambda_s$ give different Eisenstein series. Instead it suffices to consider elements of the coset $\Weyl(E_8)/\Weyl(E_7)$ and representatives that do not end on an element of $\Weyl(E_7)$. (This will be the main theme of section~\ref{EvalLCF}.) We consider as an example the Weyl word
\begin{align}
w=w_1w_3w_4w_5w_6w_7w_8.
\end{align}
Then
\begin{align}
w(2s\Lambda_8-\rho) = (8-2s)\Lambda_1+(2s-5)\Lambda_2 -\rho.
\end{align}
The intertwining factor is 
\begin{align}
M(w,2s\Lambda_8-\rho) = \frac{\xi(2s-7)}{\xi(2s)}
\end{align}
and therefore
\begin{align}
E(2s\Lambda_8-\rho,g ) = \frac{\xi(2s-7)}{\xi(2s)} E( (8-2s)\Lambda_1+(2s-5)\Lambda_2 -\rho,g).
\end{align}
At $s=\frac52$ this specialises to
\begin{align}
E(5\Lambda_8-\rho,g) = \frac{\xi(3)}{\xi(5)} E(3\Lambda_1-\rho,g)
\end{align}
and therefore relates a specific maximal parabolic Eisenstein series `on node $8$' to another specific maximal parabolic Eisenstein series, this time `on node $1$'. The former appears in the discussion of minimal theta series for $E_8$~\cite{GRS} while the latter version appears commonly in string theory, see chapter~\ref{ch:intro-strings} and section~\ref{sec:outlook-strings}, and both are related to the minimal unitary representation as we will discuss more in sections~\ref{sec:Slambda} and~\ref{smallreps} below.
\end{example}

\subsection{Weyl symmetric normalisation}

For a general (minimal) Eisenstein series $E(\lambda,g)$ on a split real simple group $G$ one can define a completed version according to
\begin{align}
E^\compl(\lambda,g) = \underbrace{\left[\prod_{\alpha>0} \xi(\langle\lambda|\alpha\rangle+1)\right]}_{N_\lambda} E(\lambda,g).
\end{align}
The~\emphindex{normalising factor} $N_\lambda$ is the denominator in $M(\wlong,\lambda)$. This function is completely invariant under the action of the Weyl group $\Weyl$: For any $w\in\Weyl$ one has
\begin{align}
E^\compl(w\lambda,g) = E^\compl(\lambda,g).
\end{align}
To see this it is sufficient to consider the action of a fundamental reflection $w_i\in \Weyl$:
\begin{align}
E^\compl(w_i \lambda, g) = N_{w_i\lambda} E(w_i\lambda,g) = N_{w_i\lambda} M(w_i,\lambda)^{-1} E(\lambda,g) =N_\lambda E(\lambda,g) = E^\compl(\lambda,g).
\end{align}
The prefactor works as follows
\begin{align}
N_{w_i\lambda} M(w_i,\lambda)^{-1} &= \frac{\xi(\langle \lambda|\alpha_i \rangle+1)}{\xi(\langle \lambda|\alpha_i\rangle)}\prod_{\alpha>0} \xi(\langle \lambda|w_i\alpha\rangle+1)\nn\\
&= \frac{\xi(\langle \lambda|\alpha_i \rangle+1)}{\xi(\langle \lambda|\alpha_i\rangle)}\xi(-\langle\lambda|\alpha_i\rangle+1)\prod_{\alpha>0\atop \alpha\neq \alpha_i} \xi(\langle \lambda|\alpha\rangle+1)\nn\\
&= \prod_{\alpha>0} \xi(\langle\lambda| \alpha\rangle +1) = N_\lambda,
\end{align}
where we have used the fact that $w_i$ permutes the set of positive roots $\Delta_+\setminus\{\alpha_i\}$ as well as the functional relation of the completed Riemann zeta function.

The normalising factor for minimal Eisenstein series has as many factors as the order of $\Weyl$. When the character defined by a weight $\lambda$ has a stabiliser that is larger than the Borel subgroup $B$, it is sufficient to use a  normalising factor with fewer factors to obtain a suitably symmetric combination. Consider the case of a non-minimal parabolic Eisenstein series $E(\lambda,g)$ where the stabiliser is given by a parabolic subgroup $P(\ints)$, see section~\ref{nonminEis}. Then the semi-simple part $M$ of the Levi subgroup $L$ of $P=LU$ has a Weyl group $\Weyl(M)$. The normalising factor $N_{P,\lambda}$ in this case can be chosen to be the denominator (after cancelling all factors) of $M(w,\lambda)$ where $w$ is defined by $\wlong(G) = w \wlong(M)$ through the longest words in $\Weyl(G)$ and $\Weyl(M)$. An alternative definition of $w$ is as the longest Weyl word in the $\Weyl(G)$ orbit of $\lambda+\rho$. The~\emphindex[Eisenstein series!normalised]{normalised Eisenstein series}
\begin{align}
E^\compl(\lambda,P,g) = N_{P,\lambda} E(\lambda,P,g).
\end{align}
This completed function then is either invariant under a reflection group isomorphic to the Weyl group generated by the simple reflections in $\Weyl(G)$ that do not belong to $\Weyl(M)$, or it maps to a similar one that is obtained by intertwining additionally by an (outer) Dynkin diagram automorphism. Furthermore, the Weyl normalised series has a different pole structure compared to the Eisenstein series $E(\lambda,P,g)$ with standard normalisation as the normalising factor has zeroes and poles.

In the case of maximal parabolic Eisenstein series with weight
\begin{align}
\lambda=2s\Lambda_{i_*} -\rho
\end{align}
this means a reflection symmetry $s\leftrightarrow \frac{\langle \rho|\Lambda_{i_*}\rangle}{\langle\Lambda_{i_*}|\Lambda_{i_*}\rangle}-s$. We illustrate this by two examples.

\begin{example}[Weyl normalisation of $SL(3,\reals)$ maximal parabolic Eisenstein series]
The first example contains a non-trivial diagram automorphism. Consider the group $G=SL(3,\reals)$ and the weight
\begin{align}
\lambda_s = 2s \Lambda_1 -\rho.
\end{align}
The associated character $\chi_s(a)=a^{\lambda_s+\rho} = a^{2s\Lambda_1}$ is invariant under a maximal parabolic subgroup $P$ with Levi factor $L=GL(1,\reals)\times M$ with $M=SL(2,\reals)$. Denote the associated maximal parabolic Eisenstein series by
\begin{align}
E(s,g) \equiv E(\lambda_s, P, g) = \sum_{\gamma\in P(\ints)\backslash G(\ints)} \chi_s(\gamma g).
\end{align}
The Weyl word $w$ that enters in the definition of the normalising factor $N_{P,\lambda_s}$ is given by the relation
\begin{align}
\underbrace{w_2w_1w_2}_{\wlong(G)} = \underbrace{w_2w_1}_{w} \underbrace{w_2}_{\wlong(M)}
\end{align}
such that
\begin{align}
M(w,\lambda_s) = \frac{\xi(2s-2)\xi(2s-1)}{\xi(2s-1)\xi(2s)}=\frac{\xi(2s-2)}{\xi(2s)}  \quad\Longrightarrow N_{P,\lambda_s} = \xi(2s).
\end{align}
The normalised series then has a reflection symmetry $s\leftrightarrow \frac32-s$ but it maps to the maximal parabolic Eisenstein series associated with the parabolic subgroup $P'$ obtained by the diagram automorphism. In other words, the characters that are being related are
\begin{align}
2s\Lambda_1-1 \quad \leftrightarrow 2\left(\tfrac32-s\right) \Lambda_2 -\rho.
\end{align}

\end{example}

The second example does not have any non-trivial automorphisms.

\begin{example}[Weyl normalisation of $E_8(\reals)$ maximal parabolic Eisenstein series]
We consider the $E_8$ Eisenstein series from example~\ref{ex:E8adj} with weight
\begin{align}
\lambda_s = 2s \Lambda_8-\rho.
\end{align}
The parabolic subgroup leaving the associated character invariant has semi-simple part $E_7$. 
The normalising factor in this case is associated with a Weyl word $w$ of length $\ell(w)=57$ that we do not spell out. The normalising factor turns out to be
\begin{align}
\label{E8adjnorm}
N_{P,\lambda_s} = \xi(2s) \xi(2s-5) \xi(2s-9) \xi(4s-28)
\end{align}
and the thus normalised Eisenstein series $E^\compl(s,g)=N_{P,\lambda} E(s,g)$ is invariant under the reflection $w_8$ that leads to the reflection law
\begin{align}
E^\compl(s,g) = E^\compl\left(\tfrac{29}2-s,g\right).
\end{align}
This example is also discussed in~\cite{GRS} and we will say more about below in example~\ref{ex:E8min}. Here, we note that the normalising factor~\eqref{E8adjnorm} has introduced a pole at $s=\frac52$ (and also for other values). This means that the special Eisenstein series from example~\ref{ex:E8adj} now appears as a~\emphindex[Eisenstein series!residue]{residue of an Eisenstein series}.
\end{example}

\subsection{Square-integrability of Eisenstein series}
\label{sec-Square}

Langlands provided a criterion for Eisenstein series to be~\emphindex[Eisenstein series!square integrable]{square integrable}~\cite[\S{5}]{LanglandsFE}. We state the criterion for Eisenstein series $E(\lambda,g)$ such that they are finite at $\lambda$, meaning that they do not have a zero or pole as a meromorphic function of $\lambda$ for the $\lambda$ chosen. 

\begin{remark}
If $E(\lambda,g)$ has a pole or zero at a given $\lambda$ one has to consider a one-parameter family in the neighbourhood of $\lambda$ and study a suitably normalised version such that the zeroth order term becomes finite~\cite{MR3034297}. In the case of $SL(2,\reals)$ and $E(s,z)$ this means multiplying by $(s-1)$ if one wants to study the square integrability of $E(s,z)$ at $s=1$. Similarly, one would have to multiply by $s^{-1}$ for the $s=0$ case.
\end{remark}

Under the assumption of a finite $E(\lambda,g)$, the constant term formula of theorem~\ref{LCFthm} implies that
\begin{align}
\lint_{N(\rats)\bs N(\ads)} E(\lambda,ng) dn = \sum_{w\in\Weyl} M(w,\lambda) a^{w\lambda+\rho}
\end{align}
is well-defined and non-vanishing function of $a$. 
\begin{proposition}[Square integrability of Eisenstein series~\cite{LanglandsFE}]
A finite Eisenstein series $E(\lambda,g)$ in the sense just described is square integrable if and only if
\begin{align}
\label{L2cond}
\Re \langle w\lambda | \Lambda_i \rangle<0 \quad\quad \textrm{for all $i=1,\ldots,\mathrm{rank}(G)$}
\end{align}
for all $w\in\Weyl$ such that $M(w,\lambda)\neq 0$.
\end{proposition}
The intuition behind this proposition is that the condition ensures that all terms fall off fast enough as one approaches any cusp of $G(\ints)\backslash G(\reals)$. A proof can be found in~\cite{LanglandsFE} and we content ourselves here with some examples.

\begin{example}[Non square integrability of $SL(2,\reals)$ Eisenstein series]
Consider square integrability of Eisenstein series on $SL(2,\reals)$. For $\lambda_s =(2s-1)\rho$ one has to check the condition~\eqref{L2cond} for the Weyl words $w=\id$ and $w=\wlong$. Plugging in the explicit expressions leads to
\begin{align}
\Re \langle\lambda_s | \Lambda_1 \rangle= \Re\frac12(2s-1) <0 
\quad\mathrm{and}\quad
\Re\langle \wlong \lambda_s|\Lambda_1 \rangle =- \Re  \frac12(2s-1)<0.
\end{align}
Clearly, these two conditions cannot be satisfied simultaneously and therefore we recover the well-known result that non-holomorphic Eisenstein series on $SL(2,\reals)$ are never square integrable. There is, however, a limiting case $\Re s=\tfrac12$, where the conditions are almost satisfied. This corresponds to the Eisenstein series on the critical line $s=\tfrac12 +it$ (for $t\in \reals$) that are $\delta$-function normalisable. See for instance~\cite{Gelbart,Terras1} for a discussion of these properties of Eisenstein series on $SL(2,\reals)$.
\end{example}

More interesting is the case when there are non-trivial square integrable functions within a degenerate principal series.
\begin{example}[Square integrability of $E_8$ maximal parabolic Eisenstein series for special $\lambda$]
Consider the maximal parabolic $E_8$ Eisenstein series $E(s,g)$ with weight
\begin{align}
\lambda=2s \Lambda_8 -\rho
\end{align}
that was introduced in example~\ref{ex:E8adj}. Computing the constant term one finds $240$ non-vanishing $M(w,\lambda)$. Checking the criterion~\eqref{L2cond} one finds that it is satisfied for the values
\begin{align}
s=\frac52,\quad s=\frac92,\quad s=7
\end{align}
and hence these are normalisable Eisenstein series for $E_8$ that belong to the discrete spectrum of the Laplacian on $G(\reals)/K(\reals)$ for $G=E_8$. 

The value $s=\tfrac52$ was also discussed in example~\ref{ex:E8adj}  and it was shown there that $E(s,g)$ for this value is functionally related to another known normalisable maximal parabolic Eisenstein series~\cite{Green:2011vz} that appears in string theory for the $R^4$ correction, see also section~\ref{sec:outlook-strings}.

The value $s=\tfrac92$ can be analysed using the functional relation
\begin{align}
E(2s\Lambda_8-\rho) = \frac{\xi(2s-10)\xi(2s-13)}{\xi(2s)\xi(2s-5)} E(2(7-s)\Lambda_1 + 2(s-9/2)\Lambda_2-\rho)
\end{align}
that shows that for $s=\tfrac92$, the adjoint $E_8$ series is connected to the maximal parabolic series on node $1$ with $s=\tfrac52$. This is the case that appears in string theory for the $D^4R^4$ correction and it is known that the function is associated with the next-to-minimal series~\cite{Green:2011vz}.

The value $s=7$ is interesting because it does not represent any simplification in the wavefront set (compared to generic $s$) and so is just part of the residual discrete spectrum with orbit type $A_2$. These cases were also analysed in~\cite{MR3034297}.
\end{example}

\section{Evaluating constant term formulas}
\label{EvalLCF}
\index{Langlands constant term formula!methods for evaluating}
Langlands constant term formula (cf.~theorem~\ref{LCFthm}) 
\begin{align}
\label{CTF9}
\lint_{N(\ints)\bs N(\reals)} E(\lambda,ng) dn = \sum_{w\in\Weyl} M(w,\lambda) a^{w\lambda+\rho}
\end{align}
is nice and compact but evaluating it will a priori produce as many terms as there are different elements in the Weyl group. Since the order of the Weyl group becomes large very quickly as the rank of $G(\reals)$ grows this can render the resulting expressions rather unwieldy. However, by dint of choice of the parameter $\lambda$ of the (degenerate) principal series the sum over Weyl elements may simplify as then some of the coefficients $M(w,\lambda)$ appearing in~\eqref{CTF9} vanish. For convenience we also recall that the definition of the intertwiner
\begin{align}
\label{Mw9}
M(w,\lambda) = \prod_{\alpha>0\atop w\alpha<0} \frac{\xi(\langle\lambda|\alpha\rangle)}{\xi(\langle\lambda|\alpha\rangle+1)}
\end{align}
and its multiplicative property
\begin{align}
\label{Mwprod9}
M(w_1 w_2,\lambda) = M(w_1,w_2\lambda) M(w_2,\lambda)
\quad \textrm{for any $w_1,w_2\in\Weyl$.}
\end{align}
A convenient method for evaluating the Langlands constant term formula can then be developed by exploiting the multiplicative relation~\eqref{Mwprod9}. We first adumbrate this method that we will refer to as the \emphindex[orbit method]{orbit method}. Then we discuss a number of examples and finally mention further simplifications that arise for constant terms in non-maximal unipotent subgroups $U\subset N$. The corresponding constant term formula was given in section~\ref{sec:CTFmax}.

\subsection{The orbit method}
\label{sec:OM}

The factor $M(w,\lambda)$ is, by its definition in~\eqref{Mw9}, given by the product of factors of the form 
\begin{align}
c(k)=\frac{\xi(k)}{\xi(k+1)},
\end{align} 
where $k=\langle \lambda|\alpha\rangle$ and $\alpha$ runs over all positive roots that satisfy $w\alpha<0$. The function $c(k)$ is sometimes referred to as the \emphindex{Harish-Chandra $c$-function}. It has a simple pole at $k=1$ and a simple zero at $k=-1$; otherwise it takes finite non-zero values for real $k$ and satisfies $c(k)c(-k)=1$ as well as $c(0)=-1$. For vanishing $M(w,\lambda)$, we are therefore particularly interested in roots $\alpha$ which satisfy $w\alpha<0$ and $\langle \lambda|\alpha\rangle=-1$.

To characterize these $\alpha$ further, let us define the \emphindex[stabiliser of a weight]{stabiliser of the weight} $\lambda$
\begin{align}
\label{stabdef}
\stab(\lambda) = \left\{\alpha\in\Pi \,|\, \langle\lambda+\rho|\alpha\rangle =  0\right\},
\end{align}
so that it is the subset of the simple roots $\Pi$ for which $\lambda+\rho$ has vanishing Dynkin labels. 

\begin{remark}
If $\stab(\lambda)\neq\{\}$, the corresponding Eisenstein series $E(\lambda,g)$ belongs to a \emph{degenerate} principal series. Let $P\subset G$ be the parabolic subgroup corresponding to $\stab(\lambda)\subset \Pi$ as defined in section~\ref{sec:parsubgp}. Then 
\begin{align}
E(\lambda,g) = E(\lambda,P,g) = \sum_{\gamma\in P(\ints) \bs G(\ints)} e^{\langle \lambda+\rho_P| H_P(\gamma g)\rangle}
\end{align}
as explained in section~\ref{parindrep}.
\end{remark}

\begin{example}[Stabiliser of a maximal parabolic $\lambda$]
As an example and referring back to section~\ref{nonminEis} we note that maximal parabolic Eisenstein series have very large stabilisers, corresponding to $\stab(\lambda)=\Pi\setminus\{\alpha_{i_*}\}$ for the value $i_*$ that determines the maximal parabolic subgroup under which $\chi(a)=a^{\lambda+\rho}$ is left-invariant.
\end{example}

If $w_i$ is the fundamental Weyl reflection in the simple root $\alpha_i$ defined in (\ref{funrefl}), it clearly maps $w_i(\alpha_i)=-\alpha_i$ and this is the only positive root that is mapped to a negative root by the fundamental reflection $w_i$~\cite{Kac}. If  furthermore $\alpha_i\in \stab(\lambda)$, then 
\begin{align}
\alpha_i\in \stab(\lambda)\,\Leftrightarrow \, \langle \lambda|\alpha_i \rangle =-1
\quad\Rightarrow\quad
M(w_i,\lambda)=c(-1)=0.
\end{align}
By the multiplicative property~\eqref{Mwprod9} one can then deduce that all Weyl words $w$ that end (on the right) on a fundamental reflection $w_i$ with $\alpha_i$ in $\stab(\lambda)$ obey $M(w,\lambda)=0$, see also~\cite{Green:2011vz}. Another way of putting this is that only those $w$ can have non-vanishing $M(w,\lambda)$ that lie in
\begin{align}
\mathcal{C}(\lambda) = \left\{ w\in \Weyl\,|\, w\alpha >0 \quad\textrm{for all $\alpha\in\stab(\lambda)$}\right\}.
\end{align}
Depending on $\stab(\lambda)$ this set can be much smaller than $\Weyl$. In fact, its order is given by $|\mathcal{C}(\lambda)|=|\Weyl|/|\Weyl(\stab(\lambda))|$, where $\Weyl(\stab(\lambda))$ is the subgroup of $\Weyl$ that is generated by taking only words in the fundamental reflections associated with $\stab(\lambda)\subset\Pi$. Moreover, one then has the simplified constant term formula
\begin{align}\label{LCF2}
\lint_{N(\rats)\bs N(\mathbb{A})} E(\chi,ng)dn = \sum_{w\in\mathcal{C}(\lambda)} a^{w\lambda+\rho} M(w,\lambda).
\end{align}

The elements in $\mathcal{C}(\lambda)$ can be constructed using the Weyl orbit of a dominant weight $\Lambda$ that is defined as follows:
\begin{definition}
Let $\lambda\in\mathfrak{h}^*$ be a weight with stabiliser $\stab(\lambda)$ as defined in~\eqref{stabdef} and $r=\dim\,\mathfrak{h}^*$ denote the rank of the underlying group. Let $I\subset\{1,\ldots,r\}$ be such that a simple root $\alpha_i$ of $\mathfrak{g}$ belongs to $\stab(\lambda)$ if and only if $i\in I$. Let $\bar{I}$ be the complement of $I$ in $\{1,\ldots,r\}$. Then the dominant weight $\Lambda$ associated to $\lambda$ is defined as a sum over fundamental weights as
\begin{align}
\label{Lambdadef}
\Lambda = \sum_{i\in \bar{I}} \Lambda_i.
\end{align}
In other words, one considers the $\lambda+\rho$ and replaces all non-zero Dynkin labels by $1$ to obtain $\Lambda$.
\end{definition}

Clearly, $\Weyl(\stab(\lambda))$ stabilises $\Lambda$ thus defined and the number of distinct points in the orbit $\Weyl\cdot\Lambda$ equals $|\Weyl|/|\Weyl(\stab(\lambda))|$. Therefore the points in the Weyl orbit are in bijection with the set $\mathcal{C}(\lambda)$. 

In order to establish the bijection, we use the fact that for each element $\mu$ in the Weyl orbit of $\Lambda$ there is a shortest element $w\in\Weyl$ that satisfies $w\Lambda=\mu$. These elements $w$ are exactly the Weyl words that make up the set $\mathcal{C}(\lambda)$. They can also be seen as specific representatives of the coset $\Weyl/\Weyl(\stab(\lambda))$ whose size was already argued above to determine the number of summands in (\ref{LCF2}). The shortest element leading to an element $\mu$ is not necessarily unique but all choices of the same shortest length yield the same factor $M(w,\lambda)$.

\index{Weyl orbit!algorithm}
A standard algorithm for constructing the Weyl orbit $\Weyl\cdot\Lambda$ of a dominant weight $\Lambda$ is as follows:
\begin{enumerate}
\item Define the initial set of orbit points as $\mathcal{O}= \{\Lambda\}$. This is the `highest' element (with respect to the \emphindex{height function} $\height(\mu) = \langle\rho|\mu\rangle$ on $\mathfrak{h}^*$) in the orbit and others will be constructed by using lowering Weyl reflections.
\item For a given $\mu\in\mathcal{O}$ compute the Dynkin labels $p_i = \langle \mu | \alpha_i\rangle$ with respect to all simple roots $\alpha_i$. 
\item If $p_i>0$ for some $i=1,\ldots,\mathrm{rank}(G)$, then construct $\mu'=w_i \mu$ where $w_i$ is the fundamental Weyl reflection in the simple root $\alpha_i$. If $\mu'$ is not already in the orbit $\mathcal{O}$, add it. 
\item For any weight $\mu$ in $\mathcal{O}$ for which steps 2 and 3 have not been carried out go to step 2.
\end{enumerate}

\begin{remark}
As for the initial dominant weight $\Lambda$ the Dynkin labels $p_i$ are zero for all simple roots $\alpha_i\in\stab(\lambda)$, any Weyl word thus constructed will end on a letter (fundamental Weyl reflection) $w_i$ that does not belong to $\Weyl(\stab(\lambda))$.
\end{remark}

\begin{remark}
In practice, it is very advisable to think of the Weyl orbit of $\Lambda$ in terms of a graph where nodes correspond to weights $\mu$ that lie in the orbit $\mathcal{O}$ and links are labelled by the fundamental reflections that relate two such weights. This graph can be constructed algorithmically starting from $\Lambda$ which corresponds to the identity element $\id\in\Weyl$ by the above algorithm and one also keeps track of the corresponding Weyl words in this way. Elements $\mu$ that are farther from the dominant weight in this graph correspond to longer Weyl words.
\end{remark}

With the Weyl orbit $\Weyl\cdot\Lambda$ one has constructed all Weyl words that belong to $\mathcal{C}(\lambda)$ and can therefore evaluate the constant term formula~\eqref{LCF2}. 
We consider an example to illustrate the method.

\newpage

\begin{figure}[t!]
\centering
\begin{picture}(100,100)
\thicklines
\put(50,50){\circle*{10}}
\put(50,50){\line(0,1){30}}
\put(50,50){\line(2,-1){30}}
\put(50,50){\line(-2,-1){30}}
\put(50,80){\circle*{10}}
\put(80,35){\circle*{10}}
\put(20,35){\circle*{10}}
\put(17,19){$1$}
\put(47,35){$2$}
\put(77,19){$3$}
\put(35,76){$4$}
\end{picture}
\caption{\label{fig:D4dynk}\textit{The Dynkin diagram of $SO(4,4)$ with labelling of simple roots.}}
\end{figure}
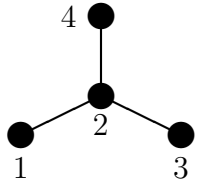

\begin{figure}
\centering
\includegraphics{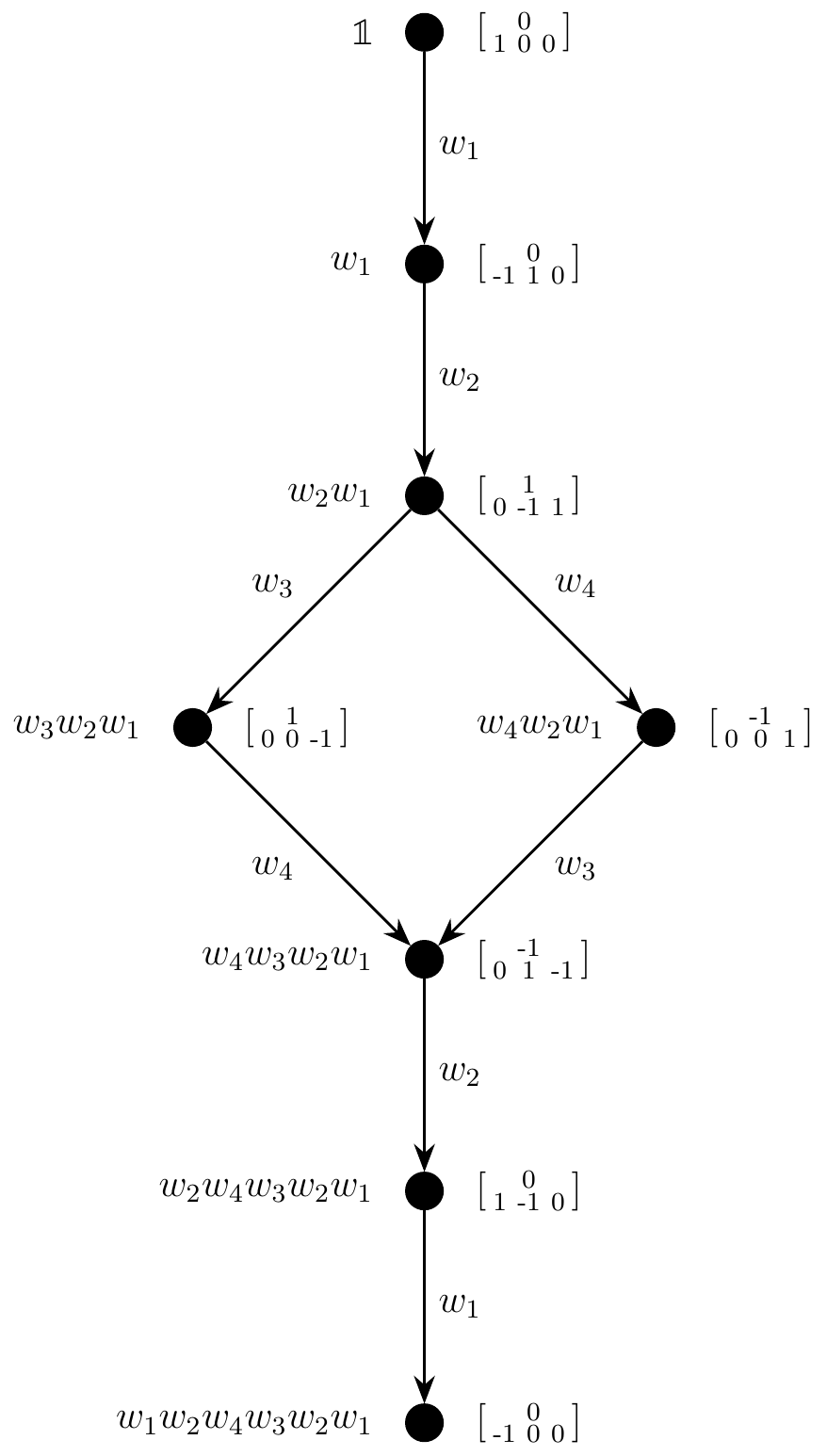}

\caption[]{\label{fig:D4Weyl} \textit{The Weyl orbit of the fundamental weight $\Lambda_1=\DDIV{1}{0}{0}{0}$ under the $D_4$ Weyl group. For each image point in the orbit, we have listed the Dynkin labels and a choice of shortest Weyl word that leads to the given point. The shortest Weyl words are those that make up the set $\mathcal{C}(\lambda)$ for $\lambda=2s\Lambda_1-\rho$ that contribute to the constant terms in (\ref{LCF2}).}}
\end{figure}

\begin{example}[Orbit method for $SO(4,4;\reals)$ maximal parabolic Eisenstein series]
\label{ex:D4series}
We consider the group $SO(4,4;\reals)$ with Dynkin diagram of type $D_4$ shown in figure~\ref{fig:D4dynk}.  The Weyl group $\Weyl$ is of order $192$ in this case. Taking 
\begin{align}
\lambda=2s\Lambda_1-\rho
\end{align}
yields a maximal parabolic Eisenstein series and 
\begin{align}
\stab(\lambda) = \left\{\alpha_2,\alpha_2,\alpha_3\right\},
\end{align}
such that $\Weyl(\stab(\lambda))$ is of type $\Weyl(A_3)$ and order $24$. The dominant weight $\Lambda$ of~\eqref{Lambdadef} equals the sum of all fundamental weights $\Lambda_i$ such that $\alpha_i\notin \stab(\lambda)$ and thus $\Lambda=\Lambda_1$.

The Weyl orbit of $\Lambda=\Lambda_1$ consists of only eight points. Figure~\ref{fig:D4Weyl} shows the graph of this Weyl orbit. 

The Weyl orbit can be calculated by starting from the highest weight $\Lambda_1$ and applying fundamental Weyl reflections in those simple roots whose \index{Dynkin labels} Dynkin labels are positive. This is the implementation of the algorithm above. For the example shown in figure~\ref{fig:D4Weyl} this allows only $w_1$ acting on $\Lambda_1$. 

Considering now the weight $\lambda=2s\Lambda_1-\rho$ that defines the $SO(4,4)$ Eisenstein series $E(\lambda,g)$, we see that the eight Weyl elements potentially contributing to the constant term formula (\ref{LCF2}) are
\begin{align}
\label{eq:D4words}
\mathcal{C}(\lambda) = \left\{ \id, w_1, w_2w_1, w_3w_2w_1,w_4w_2w_1, w_4w_3w_2w_1,w_2w_4w_3w_2w_1,w_1w_2w_4w_3w_2w_1\right\}.
\end{align}
The corresponding factors $M(w,\lambda)$ are:
\begin{align}
\label{D4consts}
\begin{array}{r|l}
w & M(w,\lambda)\\\hline
\id &1\\
w_1 &c(2 s-1)\\
w_2w_1 &c(2 s-1) c(2 s-2)\\
w_3w_2w_1 & c(2 s-1) c(2 s-2) c(2 s-3)\\
w_4w_2w_1 & c(2 s-1) c(2 s-2) c(2 s-3)\\
w_4w_3w_2w_1 & c(2 s-1) c(2 s-2)^2 c(2 s-3)^2\\
w_2w_4w_3w_2w_1 & c(2 s-1) c(2 s-2)^2 c(2 s-3)^2c(2s-4)\\
w_1w_2w_4w_3w_2w_1 & c(2 s-1) c(2 s-2)^2 c(2 s-3)^2c(2s-4)c(2s-5)\\
\end{array}
\end{align}
The table clearly reflects the multiplicative property (\ref{Mwprod9}) of the factors $M(w,\lambda)$: Moving one step down the Weyl orbit adds a single factor $c(k)$ to $M(w,\lambda)$. 

Depending on the value of $s$ some of the factors $M(w,\lambda)$ can vanish leading to a further reduction in the number of constant terms in (\ref{LCF2}).  As already argued based on the multiplicative property (\ref{Mwprod9}) we should start at the top of the Weyl orbit. Let us look at a few examples, keeping in mind that we are looking for values of $s$ where there are more factors $c(-1)$ than $c(+1)$ in the product.

The simplest case is of course $s=0$. Then only $w=\id$ has a non-vanishing $M(w,\lambda)=1$ and this is the whole constant term. This is not surprising since $s=0$ corresponds to $\lambda=-\rho$, yielding the trivial constant automorphic function $E(-\rho,g)\equiv 1$.

The next simplest case is $s=\frac12$. For this choice one the two Weyl words $w=\id$ and $w=w_1$ contribute to the constant term (\ref{LCF2}). Working out their contributions one finds that they cancel (using $c(0)=-1$) and the constant term vanishes. (The same things happens for the $SL(2,\reals)$ series; see section~\ref{SL2limits}.)

For the value $s=1$ the factor $c(2s-3)$ leads to a vanishing contribution but the factor $c(2s-1)$ has a pole so one needs to take the limit carefully. In all there are five non-vanishing contributions to the constant term since the last three orbit points (out of the total eight) contain the factor $c(2s-3)^2$. Summing up the non-vanishing contributions leads to
\begin{align}
\label{D4log}
v_1^2+\frac{6v_2}{\pi}\left( \gamma_E-\log(4\pi)-\log\left( v_1v_2^{-2}v_3v_4\right)\right)+v_3^2+v_4^2
\end{align}
where we parametrised $a= v_1^{h_1} v_2^{h_2} v_3^{h_3} v_4^{h_4}$. The logarithms arise when taking the limit $s\to 1$ and reflect the confluence of the eigenvalues of two polynomial eigenfunctions of the Laplace operator.

Further simplifications occur for $s=\frac32$ and $s=2$ that we leave to the reader to evaluate.

If one had started with a fixed value of $s$ for which simplifications occur, it would have been sufficient to construct the Weyl orbit up to the points where the $M(w,\lambda)=0$.

\end{example}

\begin{figure}
\centering
\includegraphics{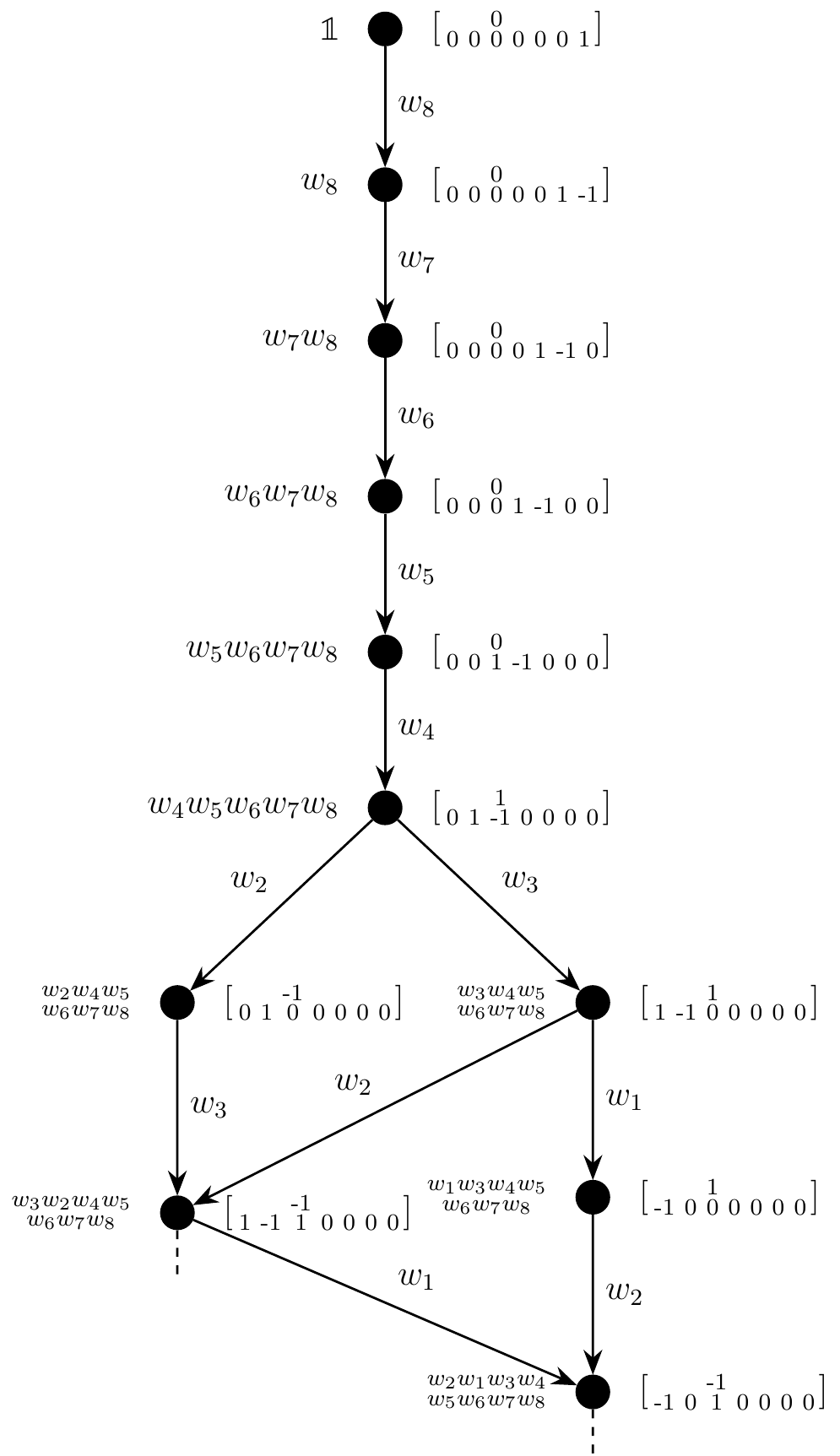}
\caption[]{\label{fig:E8Weyl} \textit{The Weyl orbit of the fundamental weight $\Lambda_1=\DEVIII{0}{0}{0}{0}{0}{0}{0}{1}$ under the $E_8$ Weyl group. For each image point in the orbit, we have listed the Dynkin labels and a choice of shortest Weyl word that leads to the given point. The shortest Weyl words are those that make up the set $\mathcal{C}(\lambda)$ that contribute to the constant terms in example~\ref{ex:E8min}.}}
\end{figure}

\subsection{Special \texorpdfstring{$\lambda$}{lambda}-values and \texorpdfstring{$E(\lambda,g)$}{E(lambda,g)}}
\label{sec:Slambda}

As we have seen in the $SO(4,4)$ example~\ref{ex:D4series} just now and in section~\ref{SL2limits}, there can be special points $\lambda\in\mathfrak{h}^*$ where the constant terms (and the whole Eisenstein series) simplify. Parametrising the weight $\lambda$ in terms of (complex) parameters $s_i$, these special points correspond to specific values for the $s_i$. These simplifications were already observed for the exceptional group $G_2$ by Langlands in his original work~\cite{LanglandsFE}, see also~\cite{MR1426903}.

In order to detect such simplifications, it is not efficient to calculate the whole set $\mathcal{C}(\lambda)$ and then the coefficients $M(w,\lambda)$ as in~\eqref{D4consts}. Due to the partially ordered structure of the Weyl orbit and the multiplicative property~\eqref{Mwprod9}  it suffices to also calculate the factor $M(w,\lambda)$ at the same as one constructs $w$ using the Weyl orbit method of section~\ref{sec:OM}.  One need not construct further any path of the graph (of increasing word length) where one of the intermediate words satisfies $M(w,\lambda)=0$. This simplifies the calculation of the constant term formula considerably~\cite{FK2012}. 

\begin{example}[Orbit method for $E_8$ maximal parabolic Eisenstein series in the minimal representation]
\label{ex:E8min}
Consider again the maximal parabolic $E_8$ Eisenstein series of example~\ref{ex:E8adj} with weight $\lambda_s=2s\Lambda_8-\rho$. For calculating the constant term, we require the Weyl orbit of the dominant weight $\Lambda_8=[0,0,0,0,0,0,0,1]$ in Dynkin label notation. In total, the Weyl orbit of $\Lambda_8$ has $240$ elements. The beginning of the Weyl orbit, computed with the orbit method is depicted in figure~\ref{fig:E8Weyl}. 

The corresponding factors $M(w,\lambda)$ for $\lambda=2s\Lambda-\rho$ are given by the following table:
\begin{align}
\label{E8consts}
\begin{array}{r|l}
w & M(w,\lambda)\\\hline
\id &1\\
w_8 &c(2 s-1)\\
w_7w_8 &c(2 s-1) c(2 s-2)\\
w_6w_7w_8 & c(2 s-1) c(2 s-2) c(2 s-3)\\
w_5w_6w_7w_8 & c(2 s-1) c(2 s-2) c(2 s-3)c(2s-4)\\
w_4w_5w_6w_7w_8 & c(2 s-1) c(2 s-2)^2 c(2 s-3)c(2s-4)c(2s-5)\\
w_2w_4w_5w_6w_7w_8 & c(2 s-1) c(2 s-2)^2 c(2 s-3)c(2s-4)c(2s-5)c(2s-6)\\
w_3w_4w_5w_6w_7w_8 & c(2 s-1) c(2 s-2)^2 c(2 s-3)c(2s-4)c(2s-5)c(2s-6)\\
w_3w_2w_4w_5w_6w_7w_8 & c(2 s-1) c(2 s-2)^2 c(2 s-3)c(2s-4)c(2s-5)c(2s-6)^2\\
w_1w_3w_4w_5w_6w_7w_8 & c(2 s-1) c(2 s-2)^2 c(2 s-3)c(2s-4)c(2s-5)c(2s-6)c(2s-7)\\
w_2w_1w_3w_4w_5w_6w_7w_8 & c(2 s-1) c(2 s-2)^2 c(2 s-3)c(2s-4)c(2s-5)c(2s-6)^2c(2s-7)\\
\vdots & \vdots
\end{array}
\end{align}
Simplifications arise as always for $s=0$ and $s=\tfrac12$. Another interesting case is $s=\tfrac52$. One sees that the nineth and eleventh entry in the table contain a factor $c(2s-6)^2$ that makes the corresponding intertwiner vanish. Since the two Weyl words are the two bottom words in the orbit constructed thus far in figure~\ref{fig:E8Weyl} one knows that all remaining Weyl words coming from the orbit method applied to $\Lambda_8$ will have vanishing $M(w,\lambda)$ and therefore the constant term consists of nine terms. Taking the limit $s\to\tfrac52$ carefully again gives logarithmic terms as in~\eqref{D4log} that we do not display here. The value $s=\tfrac52$ gives the simplest possible constant and the Eisenstein series is attached to the minimal representation as was already mentioned before.
\end{example}

\begin{remark}
Simplifications in the constant term have corresponding simplifications in the Whittaker coefficients as we will see in section~\ref{sec:CSeval}. They are typically associated with subrepresentations in the (degenerate) principal series called \emphindex[small representation|see{representation}]{small representations}\index{representation!small}. At these places the functional dimension of the automorphic representation reduces. This is discussed in more detail in section~\ref{smallreps}.
\end{remark}

\begin{example}[Minimal representation of $SL(3,\reals)$]
The Eisenstein series on $SL(3,\mathbb{A})$ introduced in section~\ref{sec:SL3ex} with weight
\begin{align}
\lambda = 2s_1\Lambda_1 + 2s_2 \Lambda_2 -\rho
\end{align}
simplifies for special values of the parameters $s_i$. Putting 
\begin{align}
s_1=0 \quad\textrm{or}\quad s_2=0 \quad\textrm{or}\quad s_3=s_1+s_2-\frac12=0
\end{align}
or 
\begin{align}
s_1=1 \quad\textrm{or}\quad s_2=1 \quad\textrm{or}\quad s_3=s_1+s_2-\frac12=1
\end{align}
makes the generic Eisenstein series into one on a maximal parabolic subgroup as defined in section~\ref{sec:parsubgp}; the case $s_1=0$ corresponds to inducing from the maximal parabolic subgroup $P_1(\mathbb{A})$. In this case the Fourier expansion simplifies considerably as can be seen be inspecting the expression of section~\ref{sec:SL3ex}. This is already manifest from (\ref{zetSL3}) that appears in the expression of any generic Whittaker coefficient. Precisely for the choices above $1/\zeta(\lambda)$ vanishes identically, implying that all generic Whittaker coefficients vanish.

For the degenerate Whittaker coefficients one also obtains shorter expressions: Out of the three Weyl elements displayed in~\eqref{eq:SL3cosW}  two have a vanishing and one is left with a single modified Bessel function with associated divisor sum. The non-abelian Fourier coefficient also simplifies as is shown in~\cite{Pioline:2009qt}.
\end{example}

\subsection{Constant terms in maximal parabolic subgroups}

Given a unipotent $U\subset N$ and an Eisenstein series $E(\lambda,g)$ on a group $G$ one can define the constant term along $U$, c.f. equation~\eqref{CTdef} and section~\ref{sec:CTFmax} by
\begin{align}
C_U= \lint_{U(\ints)\backslash U(\reals)} E(\lambda,ug) du.
\end{align}
When $U=U_{j_\circ}$ is the unipotent of a maximal parabolic subgroup $P_{j_\circ}=L_{j_\circ}U_{j_\circ}$ associated with node $j_\circ$-th simple root, a general formula for this constant term was given in theorem~\ref{consttermmaxparab}:
\begin{align}
\label{CTMax9}
\lint_{U_{j_\circ}(\ints)\bs U_{j_\circ}(\reals)} E(\lambda,ug) du = \sum_{w\in \Weyl_{j_\circ}\bs \Weyl} e^{\langle w\lambda+\rho)_{\parallel j_\circ} | H(g)\rangle} M(w,\lambda) E^{M_{j_{\circ}}}((w\lambda)_{\perp j_\circ},m).
\end{align}
The Eisenstein series on the right-hand side is one on the semi-simple part $M_{j_\circ}$ of the Levi subgroup $L_{j_\circ}= GL(1,\reals)\times M_{j_\circ}$ and the exponential prefactor is a function only on the $GL(1,\reals)$ factor. We note that it is the same numerical coefficient $M(w,\lambda)$ as in~\eqref{CTF9} that controls this constant term. As we have explained in section~\ref{sec:OM}, one can restrict the Weyl words to the set $\mathcal{C}(\lambda)$ that is in bijection with the Weyl orbit of $\lambda+\rho$ (or an equivalent dominant weight $\Lambda$ defined in~\eqref{Lambdadef}). This bijection implies that it suffices to consider Weyl words $w$ in the left coset $\Weyl/\Weyl(\stab(\lambda))\cong \mathcal{C}(\lambda)$. The addition quotient in formula~\eqref{CTMax9} then allows the restriction to the double coset~\cite{GRS,FK2012}
\begin{align}
\label{DCos}
w\in \Weyl_{j_\circ}\bs \Weyl / \Weyl(\stab(\lambda)).
\end{align}
This double coset typically has very few representatives that allow for a swift evaluation of formula~\eqref{CTMax9}. 

\begin{remark}
The double coset~\eqref{DCos} also depends on $\lambda$ and there can therefore be similar simplifications as those discussed in section~\ref{sec:Slambda}.
\end{remark}

\begin{example}[Constant term of $SO(4,4;\reals)$ Eisenstein series with respect to a maximal parabolic subgroup]
For the $SO(4,4;\reals)$ Eisenstein series considered in example~\ref{ex:D4series} with weight $\lambda=2s\Lambda_1-\rho$ we compute the constant term along the unipotent $U_3$ of the maximal parabolic subgroup $P_3$. The Weyl group of the semi-simple Levi part is again of type $\Weyl(A_3)$ and generated by the fundamental reflections $w_1$, $w_2$ and $w_4$. Inspecting the list~\eqref{eq:D4words} of elements of $\mathcal{C}(\lambda)$ shows that the double coset~\eqref{DCos} in this case has only two representatives, namely
\begin{align}
\id\quad\textrm{and}\quad w_3w_2w_1 .
\end{align}
If we denote the coordinate on $GL(1,\reals)$ by
\begin{align}
r = e^{\langle\Lambda_3|H(g)\rangle},
\end{align}
we obtain for trivial representative the decomposition
\begin{align}
(2s\Lambda_1)_{\parallel 3} =\frac{\langle 2s\Lambda_1|\Lambda_3\rangle}{\langle\Lambda_3|\Lambda_3\rangle}\Lambda_3 =  s \Lambda_3,\quad
(2s\Lambda_1-\rho)_{\perp 3} = (2s-1)\Lambda_1^{M_3} -\Lambda_2^{M_3}-\Lambda_3^{M_3} 
\end{align}
(where $\Lambda_i^{M_3}$ are the three fundamental weights of $M_3(\reals)=SL(4,\reals)$) and a similar decomposition for the other representative the following constant term
\begin{align}
\lint_{U_3(\ints)\bs U_3(\reals)} E(\lambda,ug) du = r^s E([2s-1,-1,-1],m) + r^{3-s} \frac{\xi(2s-3)}{\xi(2s)} E([-1,-1,2(s-1)-1],m),
\end{align}
where we have indicated the weight on the semi-simple subgroup $M_3$ by its Dynkin labels and have evaluated the intertwiner $M(w_3w_2w_1,\lambda)$ using~\eqref{D4consts}.
\end{example}

\section{Evaluating spherical Whittaker coefficients}
\label{sec:CSeval}
\index{Whittaker coefficient!method for evaluating}

We now turn the question of efficiently evaluating degenerate Whittaker coefficients that are given by theorem~\ref{DegWhittThm} whose result we briefly recall. The final formula there was
\begin{align}
\label{degW9}
W^\circ_\psi (\lambda,a) = \sum_{w_c\wlong'\in \Weyl/\Weyl'} a^{(w_c\wlong')^{-1}\lambda+\rho}  M(w_c^{-1},\lambda) W'^\circ_{\psi^a}(w_c^{-1}\lambda,\id).
\end{align}
We briefly recall the notation used in this formula. $\psi$ denotes a degenerate character on the maximal unipotent $N\subset B\subset G$. It has support $\supp(\psi)\subset \Pi$ given by~\eqref{supppsi} and this subset of simple roots defines a semi-simple subgroup $G'\subset G$ with Weyl group $\Weyl'=\Weyl(\supp(\psi))$. The longest element $\Weyl'$ is called $\wlong'$ and $w_c\in\Weyl$ satisfies $w_c\alpha>0$ for all $\alpha\in\supp(\psi)$. The representative $w_c$ of the coset $\Weyl/\Weyl'$ can be constructed using the orbit method of section~\ref{sec:OM}. For any $a\in A(\ads)$, the twisted character $\psi^a(n)=\psi(ana^{-1})$ restricted to the unipotent $N'\subset G'$ is generic and $W'^\circ_{\psi^a}(w_c^{-1}\lambda,\id)$ is the (generic) spherical Whittaker function on the $G'$ of the prinicipal series representation given by the restriction of the weight $w_c^{-1}\lambda$ to $G'$ and evaluated at the identity. An example of this formula was worked out for $SL(3,\ads)$ in section~\ref{sec:SL3ex}.

One sees from formula~\eqref{degW9} that it is again an intertwining coefficient $M(w_c^{-1},\lambda)$ that controls possible simplifications in the degenerate Whittaker coefficients. We know from the discussion of the orbit method in section~\ref{sec:OM} that only those $w_c^{-1}$ give a non-trivial $M(w_c^{-1},\lambda)$ that are the minimal representatives of the coset $\Weyl/\Weyl(\stab(\lambda))$ (constructed by the orbit method). Due to the inverse $w_c^{-1}$ we are therefore again faced with a double coset
\begin{align}
\Weyl(\supp(\psi))\bs \Weyl / \Weyl(\stab(\lambda)).
\end{align}
This has very few representatives one has to consider when evaluating~\eqref{degW9} for a non-generic $\lambda$ and $\psi$. 

For which Eisenstein series $E(\lambda,g)$ does formula~\eqref{degW9} actually offer the prospect of helping find complete information about the Fourier expansion? As mentioned already in remark~\ref{rmk:degW} this will happen when $\lambda$ is such that $E(\lambda,g)$ is not in the generic principal series but in a degenerate one or even at one of the special $\lambda$ values discussed in section~\ref{sec:Slambda}.

\subsection[Degenerate principal series and degenerate Whittaker coefficients]{Degenerate principal series and degenerate\\ Whittaker coefficients}
\index{degenerate principal series!and Whittaker coefficients}

If $\lambda$ is such that $\stab(\lambda)\neq \{\}$, the Eisenstein series $E(\lambda,g)$ is associated to the degenerate principal series, see section~\ref{nonminEis}. Prime example are maximal parabolic Eisenstein series when $\stab(\lambda)=\Pi\setminus\{\alpha_{j_*}\}$ where $\alpha_{j_*}$ is a single simple root that defines a maximal parabolic subgroup $P_{j_*}\subset G$. More generally, we define a parabolic subgroup
\begin{align}
P_\lambda = L_\lambda U_\lambda, 
\end{align}
such that the semi-simple part $M_\lambda\subset L_\lambda$ has the simple root system given by $\stab(\lambda)$. Then the Eisenstein series $E(\lambda,g)$ belongs to the degenerate principal series with Gelfand--Kirillov dimension\index{Gelfand--Kirillov dimension}
\begin{align}
\GKdim( I(\lambda) ) = \dim \left(P_\lambda\bs G\right)
\end{align}
or even a subrepresentation of this in case $\lambda$ sits at a special value. 

For automorphic functions in a degenerate principal series one has that typically not all Whittaker coefficients are non-zero and often the generic ones are absent. In order to determine which Whittaker coefficients are non-zero we recall the notion of a wavefront set\index{wavefront set} introduced in section~\ref{sec:orbits}. The wavefront set is the set of nilpotent orbits of $G(\reals)$ such that there are non-trivial Fourier coefficients (or Whittaker coefficients) associated with it. Characters $\psi$ are associated with nilpotent elements of $\mathfrak{g}=\mathfrak{g}(\reals)$ and one has to consider their $G(\reals)$ orbits in an automorphic representation. 

Nilpotent orbits of $\mathfrak{g}$ under $G(\reals)$ come with a certain (even) dimension and they must be able to `fit into' the automorphic representation of $E(\lambda,g)$ for a non-trivial (non-vanishing) Fourier coefficient to exist. There is a symplectic structure on a nilpotent orbit~\cite{CollingwoodMcGovern} and only a Lagrangian subspace corresponds to the character $\psi$ of a Whittaker coefficient. 
Let $X\in\mathfrak{g}$ be a nilpotent element and $\mathcal{O}_X$ its corresponding nilpotent orbit under $G(\reals)$. The constraint just explained means that there can be non-trivial Fourier coefficients only if 
\begin{align}
\label{eq:orbconstr}
\frac12 \dim_\reals \mathcal{O}_X \leq \dim\left(P_\lambda\bs G\right).
\end{align}

\begin{remark}
If $\lambda$ is generic such that $E(\lambda,g)$ is in the full principal series, then $P_\lambda$ equals the standard Borel subgroup $B\subset G$ and the Gelfand--Kirillov dimension equals $\frac12\left(\dim\mathfrak{g} - \dim\mathfrak{h}\right)$. At the same time, the largest nilpotent orbit (called the principal orbit) has dimension $\dim\mathfrak{g} - \dim\mathfrak{h}$, confirming the fact that such a generic Eisenstein series will have generic Fourier coefficients in general.
\end{remark}

The condition~\eqref{eq:orbconstr} puts strong constraints on the orbits one has to consider for a degenerate principal series representation. Moreover, if a nilpotent orbit $\mathcal{O}_X$ has a representative $X$ that lies completely in $[N,N]\bs N$ such that there is a character $\psi_X: N\to U(1)$ associated with it, one can test whether $\mathcal{O}_X$ belongs to the wavefront set by computing the (degenerate) Whittaker coefficient for $\psi_X$ with formula~\eqref{degW9}. 

\begin{example}[Minimal orbit and $A_1$-type Whittaker coefficients]
Any simple group $G(\reals)$ has a unique minimal non-trivial nilpotent orbit $\mathcal{O}_{\mathrm{min}}$ that is given by the orbit of a generator $E_\theta$ from the root space of the highest root $\theta$. If $G(\reals)$ is simply-laced, one can alternatively choose as a nilpotent representative any simple step operator $X=E_{\alpha_i}$ where $\alpha_i$ is a simple root. The corresponding character $\psi_X$ is maximally degenerate and has $\supp(\psi_X)=\{\alpha_i\}$ and one can compute the associated degenerate Whittaker coefficients using~\eqref{degW9} in terms of Whittaker coefficients on $SL(2,\reals)$ (which are completely known). The minimal orbit is called type $A_1$ in Bala--Carter terminology~\cite{BalaCarterI} and this relates to the subgroup $G'$ that appears in the formula for degenerate Whittaker coefficients.
\end{example}

As another example we consider the consequences of the condition~\eqref{eq:orbconstr} on a degenerate principal series of the group $E_8$.

\begin{example}[Wavefront sets of the adjoint $E_8$ series]
\label{ex:E8adj2}
Consider the maximal parabolic Eisenstein series $E(\lambda,g)$ of $E_8(\reals)$ given by the weight
\begin{align}
\lambda = 2s \Lambda_8 -\rho
\end{align}
as in example~\ref{ex:E8adj}. The degenerate principal series that $E(\lambda,g)$ belongs to is of Gelfand--Kirillov dimension (for generic $s$)
\begin{align}
\GKdim I(\lambda) = \dim(E_8) - \dim(P_\lambda) = 248- (133+1+56+1) = 57.
\end{align}
According to~\eqref{eq:orbconstr}, the largest nilpotent orbit that can contribute to the wavefront set of $E(\lambda,g)$ is therefore of dimension $114$. Here is a list of nilpotent orbits of (split) $E_8(\reals)$ of small dimension~\cite{DjokovicE8,CollingwoodMcGovern}
\begin{center}
\begin{tabular}{ccc}
\toprule
$\dim\mathcal{O}$ & Bala--Carter label & weighted diagram over $\cx$\\
\midrule
$0$ & $0$ & $\DEVIII{0}{0}{0}{0}{0}{0}{0}{0}$\\
$58$ & $A_1$ & $\DEVIII{0}{0}{0}{0}{0}{0}{0}{1}$\\
$92$ & $2A_1$ & $\DEVIII{1}{0}{0}{0}{0}{0}{0}{0}$\\
$112$ & $3A_1$ & $\DEVIII{0}{0}{0}{0}{0}{0}{1}{0}$\\
$114$ & $A_2$ and $(4A_1)''$ & $\DEVIII{0}{0}{0}{0}{0}{0}{0}{2}$\\
$\vdots$&$\vdots$&$\vdots$\\
\bottomrule
\end{tabular}
\end{center}
The last entry corresponds to a single complex orbit of type $A_2$ that splits into two real orbits~\cite{DjokovicE8}. All these orbits have representatives $X$ in $[N,N]\bs N$ and the associated Whittaker coefficients can be calculated using~\eqref{degW9}. More complicated Whittaker coefficients are absent in this degenerate principal series. For special $s$ values not all the orbits in the above table appear in the wavefront set.
\end{example}

\subsection{Whittaker coefficients of maximal parabolic Eisenstein series}

For maximal parabolic Eisenstein series one can also make statements about the vanishing of some Whittaker coefficient. Consider an Eisenstein series on $G(\mathbb{A})$ induced from a character $\chi: P_{i_*}(\mathbb{A})\to \cx^\times$, i.e., one that is in the degenerate principal series and is parametrised by a single complex parameter $s\in \cx$ through the weight $\lambda=2s\Lambda_{i_*}-\rho$, see proposition~\ref{prop-minimaltomaximal}. The Whittaker integral for an arbitrary character $\psi: N(\mathbb{Q})\backslash N(\mathbb{A})\to U(1)$ leads to (cf. (\ref{FwSplit}))
\begin{align}
W_{\psi} (\chi,a) = \sum_{w\in \Weyl(P_{i_*}) \backslash \Weyl} \lint_{N^w_{\{\beta\}}} \overline{\psi(n_\beta)} dn_\beta \cdot \lint_{N_{\{\gamma\}}^w} \chi(w n_\gamma a ) \overline{\psi(n_\gamma)}dn_\gamma,
\end{align}
where the important point is now that the set of contributing Weyl words is restricted to the quotient $\Weyl(P_{i_*})\backslash\Weyl$ from the outset. Again, the integral over $N^w_{\{\beta\}}$ can make the whole expression vanish and imposes constraints on $w$ and $\psi$. Now the set of positive roots $\beta$ that appear in that integral is 
\begin{align}
\left\{ \beta>0 \,|\, w\beta \in \Delta_{P_{i_*}}\right\},
\end{align}
where $\Delta_{P_{i_*}}$ denotes the subset of all roots in $\Delta$ that belong to the maximal parabolic $P_{i_*}(\mathbb{A})$; it involves all positive roots and some negative roots. The set of $\beta$ now always involves a simple root for any $w$, therefore for a generic $\psi$ the integral over $n_\beta$ will vanish and we conclude that the generic Whittaker coefficient for a maximal parabolic Eisenstein series vanishes. 

Another way to see this is by noting that the factor $1/\zeta(\lambda)$ that appears in (\ref{zetaeps}) contains generally factors $(1-p^{-(\langle\lambda|\alpha_i\rangle+1)})=(1-p^{-2s_i})$ for all simple roots $\alpha_i$. For any degenerate principal series Eisenstein series one of these factors vanishes identically, and there is no pole a the same $s_i$ values in the factor $\eps(\lambda)$ (see~(\ref{zetaeps})). This is guaranteed by the holomorphy of the local Whittaker function~\cite{CasselmanShalika}.

We will come back to this in the discussion in section~\ref{sec-FourWhitt}.

\subsection{Examples of degenerate Whittaker coefficients}

We now present some explicit expressions for degenerate Whittaker coefficients calculated with the help of~\eqref{degW9}. The examples are taken mainly from~\cite{FKP2013}. The following notation will be used
\begin{align}
B_{s,m} (v) = \frac{1}{\xi(2s)} \tilde{B}_{s,m}(v) = \frac{2}{\xi(2s)} |v|^{s-1/2}|m|^{s-1/2} \sigma_{1-2s}(m) K_{s-1/2}(2\pi |m| v) 
\end{align}
for a Whittaker coefficients on an $SL(2,\reals)$ subgroup. For an $SL(3,\reals)$ subgroup we write similarly
\begin{align}
B_{s_1,s_2,m_1,m_2}(v_1,v_2) = \frac{1}{\xi(2s_1)\xi(2s_2)\xi(2s_1+2s_2-1)} \tilde{B}_{s_1,s_2,m_1,m_2}(v_1,v_2) .
\end{align}
The explicit expression for $\tilde{B}_{s_1,s_2,m_1,m_2}(v_1,v_2)$ in terms of an integral over two Bessel functions can be found in~\cite{BumpSL3,Pioline:2009qt}, see also the $SL(3,\ads)$ example in section~\ref{sec:SL3ex}. For our purposes, we only need to know that it is finite and non-zero for all values of $s_1$ and $s_2$. The same is true for $\tilde{B}_{s,m}(v)$.

\begin{example}[Degenerate Whittaker coefficients of type $A_1$ for Eisenstein series of $E_{n\geq6}(\reals)$]
\label{ex:degen-whitt-Ed}
We consider maximal parabolic Eisenstein series of the finite-dimensional exceptional groups $E_n(\mathbb R)$ with weight vector $\lambda=2\cdot\frac32\Lambda_1-\rho$ for $n=6,7,8$. As before we use the standard Bourbaki labelling of the nodes of the Dynkin diagram. In these examples all non-vanishing (abelian) Whittaker coefficients turn out to be given by a finite sum of $n$ Whittaker coefficients on the $SL(2,\mathbb R)$ subgroup associated with each node of the Dynkin diagram. The full expression for the Whittaker coefficients will be given by
\begin{align}
\sum_{\psi\neq 0} W^\circ_\psi(\lambda, na) =\sum_{\alpha\in \Pi} \sum_{\psi_\alpha} c_\alpha(a)W'^\circ_{\psi^{a}_\alpha}(\lambda'_\alpha, \id)\psi_\alpha(n)\,,
\end{align} 
where $W^\circ_\psi$ on the left-hand side is given by~\eqref{degW9} and a detailed derivation of the right-hand side can be found in~\cite{FKP2013}. In the language of section~\ref{sec:ANA} this is only the abelian part of the Fourier expansion along the maximal unipotent $N$. For the maximally degenerate character $\psi_\alpha$ associated with the simple roots $\alpha$, $m_\alpha$ is the only non-zero charge and $c_\alpha(a)$ is a function of the variables parametrising the Cartan torus. Furthermore, $W'^\circ_{\psi^{a}_\alpha}$ is a generic Whittaker coefficient on $SL(2,\mathbb R)$ subgroup associated with the simple roots $\alpha$ and $\lambda'_\alpha$ is the projection of the weight $\lambda$ onto this subgroup. 

We can then provide lists of the degenerate Whittaker coefficients for each case where the results are taken from~\cite{FKP2013}:
\begin{itemize}
\item $E_8:$
\begin{center}
\renewcommand{\arraystretch}{1.9}
\begin{tabular}{ | c || c | }
  \hline                       
  $\psi_\alpha$ & $c_\alpha(a)W'^\circ_{\psi^{a}_\alpha}(\chi'_\alpha, \id)$ \\ \hline \hline
 $(m,0,0,0,0,0,0,0)$ & $v_3^2v_1^{-1} B_{3/2,m}\left(v_1^2v_3^{-1}\right)$ \\ \hline
 $(0,m,0,0,0,0,0,0)$ & $\frac{v_2^2\tilde{B}_{0,m}(v_2^2v_4^{-1})}{\xi(3)}$ \\ \hline
 $(0,0,m,0,0,0,0,0)$ & $\frac{\xi (2) v_4 B_{1,m}\left(v_3^2 v_1^{-1}v_4^{-1}\right)}{\xi (3)}$ \\ \hline
 $(0,0,0,m,0,0,0,0)$ & $\frac{v_4\tilde{B}_{1/2,m}(v_4^2v_2^{-1}v_3^{-1}v_5^{-1})}{\xi(3)}$  \\ \hline
 $(0,0,0,0,m,0,0,0)$ & $\frac{v_5^2\tilde{B}_{0,m}(v_5^2v_4^{-1}v_6^{-1})}{\xi(3)v_6}$ \\ \hline
 $(0,0,0,0,0,m,0,0)$ & $\frac{\xi (2) v_6^3v_7^{-2} B_{-1/2,m}\left(v_6^2 v_5^{-1}v_7^{-1}\right)}{\xi (3)}$  \\ \hline
 $(0,0,0,0,0,0,m,0)$ & $\frac{v_7^4 B_{-1,m}\left(v_7^2 v_6^{-1}v_8^{-1}\right)}{r_8^3}$ \\ \hline
  $(0,0,0,0,0,0,0,m)$ & $\frac{\xi (4) v_8^5 B_{-3/2,m}\left(v_8^2 v_7^{-1}\right)}{\xi (3)}$ \\ \hline
 \end{tabular}
 \end{center}
Clearly, there is a simple pattern involving the weight of the Bessel function and the number of Weyl reflections needed to relate a given degenerate character $\psi$ to the one in the first line where the unipotent of the Whittaker expansion is aligned with the unipotent of the maximal parabolic defining the Eisenstein series. In cases where the function $B_{s,m}$ would be ill-defined due to a pole in the $\xi$-prefactor one has that consistently only the differently normalised and finite function $\tilde{B}_{s,m}$ appears. This pattern continues to hold for the other rank cases.

\item $E_7:$
\begin{center}
\renewcommand{\arraystretch}{1.9}
\begin{tabular}{ | c || c | }
  \hline                       
  $\psi_\alpha$ & $c_\alpha(a)W'^\circ_{\psi^{a}_\alpha}(\chi'_\alpha, \id)$ \\ \hline \hline
 $(m,0,0,0,0,0,0)$ & $v_3^2 v_1^{-1}B_{\frac{3}{2},m}\left(v_1^2 v_3^{-1}\right)$ \\ \hline
 $(0,m,0,0,0,0,0)$ & $\frac{v_2^2\tilde{B}_{0,m}(v_2^2v_4^{-1})}{\xi(3)}$ \\ \hline
 $(0,0,m,0,0,0,0)$ & $\frac{\xi (2) v_4 B_{1,m}\left(v_3^2 v_1^{-1}v_4^{-1}\right)}{\xi (3)}$ \\ \hline
 $(0,0,0,m,0,0,0)$ & $\frac{v_4\tilde{B}_{1/2,m}(v_4^2v_2^{-1}v_3^{-1}v_5^{-1})}{\xi(3)}$  \\ \hline
 $(0,0,0,0,m,0,0)$ & $\frac{v_5^2\tilde{B}_{0,m}(v_5^2v_4^{-1}v_6^{-1})}{\xi(3)v_6}$ \\ \hline
 $(0,0,0,0,0,m,0)$ & $\frac{\xi (2) v_6^3 v_7^{-2} B_{-1/2,m}\left(v_6^2 v_5^{-1}v_7^{-1}\right)}{\xi (3)}$  \\ \hline
 $(0,0,0,0,0,0,m)$ & $v_7^4 B_{-1,m}\left(v_7^2 v_6^{-1}\right)$ \\ \hline
 \end{tabular}
 \end{center}

\item $E_6$:
\begin{center}
\renewcommand{\arraystretch}{1.9}
\begin{tabular}{ | c || c | }
  \hline                       
  $\psi_\alpha$ & $c_\alpha(a)W'^\circ_{\psi^{a}_\alpha}(\chi'_\alpha, \id)$ \\ \hline \hline
 $(m,0,0,0,0,0)$ & $v_3^2 v_1^{-1}B_{3/2,m}\left(v_1^2v_3^{-1}\right)$ \\ \hline
 $(0,m,0,0,0,0)$ & $\frac{v_2^2\tilde{B}_{0,m}(v_2^2v_4^{-1})}{\xi(3)}$ \\ \hline
 $(0,0,m,0,0,0)$ & $\frac{\xi (2) v_4 B_{1,m}\left( v_3^2 v_1^{-1}v_4^{-1}\right)}{\xi (3)}$ \\ \hline
 $(0,0,0,m,0,0)$ & $\frac{v_4\tilde{B}_{1/2,m}(v_4^2v_2^{-1}v_3^{-1}v_5^{-1})}{\xi(3)}$  \\ \hline
 $(0,0,0,0,m,0)$ & $\frac{v_5^2\tilde{B}_{0,m}(v_5^2v_4^{-1}v_6^{-1})}{\xi(3)v_6}$ \\ \hline
 $(0,0,0,0,0,m)$ & $\frac{\xi (2) v_6^3 B_{-1/2,m}\left(v_6^2 v_5^{-1}\right)}{\xi (3)}$  \\ \hline
 \end{tabular}
 \end{center}

 \end{itemize}
\end{example}
The following provides an example of a degenerate Whittaker coefficient of type $A_2$.
\begin{example}[Degenerate Whittaker coefficients of type $A_2$ for Eisenstein series of $E_8(\reals)$]
\label{ex:E8-A2-Whittaker}
For the $E_8$ series of example~\ref{ex:E8adj2} we compute the Whittaker coefficient associated with the degenerate character $\psi$ on $N$ with `charges'
\begin{align}
\psi \quad \leftrightarrow \quad \DEVIII{0}{0}{0}{0}{0}{0}{m}{n}.
\end{align}
This choice of character $\psi$ is associated with the $114$-dimensional nilpotent orbit of type $A_2$ from the table in example~\ref{ex:E8adj2}. For simplicity we put the torus element $a=\id$. Then formula~\eqref{degW9} gives
\begin{align}
\label{E8WA2}
W_\psi(\lambda,\id) = \frac{\xi(2s-11)\xi(2s-14)\xi(2s-18)\xi(4s-29)}{\xi(2s)\xi(2s-5)\xi(2s-9)\xi(4s-28)} B_{6-s,\tfrac{19}2-s,m,n}\Big(1,1\Big).
\end{align}
(The contributing Weyl word has length $30$ and we do not spell it out here.) As was argued in example~\ref{ex:E8min}, the value $s=\tfrac52$ corresponds to a simpler Eisenstein series where the constant term simplifies. From the above formula we can see this also in the Whittaker coefficient of type $A_2$. The prefactor tends to zero for $s\to\tfrac52$ while the $SL(3,\reals)$ Whittaker coefficient stays finite and hence the degenerate  coefficient~\eqref{E8WA2} disappears. 

In the case $s=\tfrac52$ one check similarly that all Whittaker coefficients but the ones of type $A_1$ vanish, consistent with the fact that the corresponding Eisenstein series belongs to an automorphic realisation of the minimal representation.
\end{example}

\subsection[Relation between general Fourier coefficients and Whittaker coefficients]{Relation between general Fourier coefficients and\\ Whittaker coefficients}
\label{sec-FourWhitt}

In this section we investigate, based on the methods of Miller--Sahi 
\cite{MillerSahi} and Ginzburg \cite{GinzburgConjectures}, how to compute Fourier coefficients $F_{\psi_U}$ on unipotent subgroups $U$ (from definition~\ref{def:F}) in terms of Whittaker coefficients $W_{\psi_N}$, the latter of which are known using the methods of chapter~\ref{ch:Whittaker-Eisenstein}.  Details can be found in \cite{Gustafsson:2014iva} and similar ideas were discussed in section~\ref{sec:Piatetski-Shapiro-Shalika}. 
As we are discussing characters and Fourier coefficients on different subgroups we will adopt the subscript notation used in section~\ref{sec:SL3-non-abelian} to avoid ambiguities.

Since both Whittaker coefficients on $N$ and the constant term simplifies for autormorphic forms in small representations as seen in the examples above, it is natural to also expect simplifications for more general Fourier coefficients. Using the wavefront set and the arguments in the previous sections one can tell for which representations a Fourier coefficient is non-vanishing, but by rewriting $F_{\psi_U}$ in terms of $W_{\psi_N}$ it is also possible to see how a non-vanishing $F_{\psi_U}$ simplifies for smaller representations.

We have already seen an example of such a computation in proposition \ref{prop:SL3-non-abelian-Whittaker-from-abelian} where we, for $G = SL(3)$, showed how the Fourier coefficients on $Z = [N,N]$ can be expressed as a sum of $G$-translated Whittaker coefficients on $N$. When restricting to Eisenstein series in the minimal representation this expression simplifies as follows.

\begin{example}[Non-abelian $SL(3)$ Fourier coefficient in minimal representation]
    \label{ex:SL3-non-abelian-min-rep}
    For the example $G = SL(3)$ with $\lambda = (2s_1-1) \Lambda_1 + (2s_2-1) \Lambda_2$ the generic principal series has Gelfand--Kirillov dimension $\GKdim(I(\lambda)) = \dim G - \dim B = 3$, but for $(s_1, s_2) = (s, 0)$ or $(s_1, s_2) = (0, s)$ it reduces to $\GKdim(I(\lambda_\text{min})) = 2$ for which the Eisenstein series belongs to the minimal automorphic representation.
    
    The orbits of $SL(3, \cmplx)$ are \cite{CollingwoodMcGovern} 
    \newcommand{\SLIIWeightedDiagram}[2]{%
        \begin{tikzpicture}[]%
            \tikzstyle{dot}=[circle, draw, thick, fill, scale=0.5]%
            \node (n1) at (0, 0) [dot,label=above:$\scriptstyle{#1}$] {};%
            \node (n2) at (0.5, 0) [dot,label=above:$\scriptstyle{#1}$] {};%
            \draw [thick] (n1) -- (n2);%
        \end{tikzpicture}}

    \begin{center}
        \begin{tabular}{ccc}
            \toprule
            $\dim \mathcal{O}$ & Bala-Carter label & weighted diagram \\
            \midrule
            0 & $0$ & \SLIIWeightedDiagram{0}{0} \\
            4 & $A_1$ & \SLIIWeightedDiagram{2}{2} \\
            6 & $A_2$ & \SLIIWeightedDiagram{1}{1} \\
            \bottomrule
        \end{tabular}
    \end{center}

    This means that for $\pi_\text{min}$ we only have Whittaker coefficients of type $A_1$, that is, maximally degenerate Whittaker coefficients charged under a single simple root.

    In proposition \ref{prop:SL3-non-abelian-Whittaker-from-abelian} it was shown that
    \begin{equation}
        \label{eq:SL3-non-abelian-Whittaker-from-abelian-pre-min-rep}
        W_{\psi_Z}^{(k)}(\chi, g) = \sum_{m_1, m_2 \in \rats} W_{\psi_N}^{(m_1, d)} (\chi, l g)
    \end{equation}
    where $d = d(k, m_2)$ as defined in \eqref{eq:gcd-rational} and $l$ depends on $k$ and $m_2$ as described in the proposition. Recall from section \ref{sec:SL3ex} that $W_{\psi_N}^{(m_1, d)}$ for $m_1 \neq 0$ contains a convoluted integral of two Bessel functions, while $W_{\psi_N}^{(0, d)}$ is simpler being proportional to a single Bessel function.
    
    When restricting $\chi$ to $\chi_\text{min}$ (parametrised by $s$), we get that the sum over charges in \eqref{eq:SL3-non-abelian-Whittaker-from-abelian-pre-min-rep} collapses to $m_1 = 0$ since $d \neq 0$, simplifying the expression for $W_{\psi_Z}^{(k)}$ which then only contains single Bessel functions. Note though that the sum over $m_2$ (that is, over $l$-translates) still remains, that is,
\begin{equation}
    W_{\psi_Z}^{(k)}(\chi_\text{min}, g) = \sum_{m_2 \in \rats} W_{\psi_N}^{(0, d)}(\chi_\text{min}, l g) \, .
\end{equation}
When inserting the argument $g = (g_\infty, \id, \id, \ldots)$ the sum over rational charges becomes a sum over integers, similar to what happens in proposition~\ref{prop:SL3-non-abelian-expansion}, since $l \in SL(3, \ints)$ by design.
\end{example}

Let us now consider another example for $G = SL(3)$, but with Fourier coefficients on a maximal parabolic subgroup. Here, the expression simplifies even further in the minimal representation resulting in a single translated maximally degenerate Whittaker coefficient on $N$.

\begin{example}[Maximal parabolic Fourier coefficient in minimal representation]    
    \label{ex:SL3-maximal-parabolic-Fourier}
    Continuing with $G = SL(3)$, we will now see that a Fourier coefficient $F_{\psi_U}$ on the maximal parabolic subgroup corresponding to the first simple root
    \begin{equation}
        P = P_{1} = L U = 
        \begin{psmallmatrix}
            * & * & * \\
            0 & * & * \\
            0 & * & *    
        \end{psmallmatrix} \qquad L = 
        \begin{psmallmatrix}
            * & 0 & 0 \\
            0 & * & * \\
            0 & * & *
        \end{psmallmatrix} \qquad U =
        \begin{psmallmatrix}
            1 & * & * \\
            0 & 1 & 0 \\
            0 & 0 & 1
        \end{psmallmatrix}
    \end{equation}
    can be expressed as a single $L$-translated, maximally degenerate Whittaker coefficient in the minimal representation. Let
    \begin{equation}
        u = 
        \begin{psmallmatrix}
            1 & u_1 & u_2 \\
            0 & 1 & 0 \\
            0 & 0 & 1
        \end{psmallmatrix} \in U
        \qquad
        \psi_U(u) = e^{2\pi i (m_1 u_1 + m_2 u_2)} \qquad m_1, m_2 \in \rats \qquad m_1 m_2 \neq 0 \, .
    \end{equation}
    Then, $d = d(m_1, m_2)$ as defined in \eqref{eq:gcd-rational} is strictly positive and $m'_i \coloneqq m_i / d \in \ints$ with $\gcd(m'_1, m'_2) = 1$ which tells us that there exist integers $\alpha$ and $\beta$ such that
    \begin{equation}
        l =
        \begin{psmallmatrix}
            1 & 0 & 0 \\
            0 & \alpha & \beta \\
            0 & -m'_2 & m'_1
        \end{psmallmatrix} \in L(\ints) \, .
    \end{equation}
    Now we conjugate the Fourier coefficient with $l$ as follows (cf. \eqref{Whittakerorbit})
    \begin{equation}
        \begin{split}
            \MoveEqLeft
            F_{\psi_U}^{(m_1, m_2)}(\chi, g) \coloneqq \intl_{(\rats \bs \ads)^2} E(\chi, 
            \begin{psmallmatrix}
                1 & u_1 & u_2 \\
                0 & 1 & 0 \\
                0 & 0 & 1
            \end{psmallmatrix}
            g) e^{-2\pi i(m_1 u_1 + m_2 u_2)} \, du^2 \\
            &= \intl_{(\rats \bs \ads)^2} E(\chi,
            \begin{psmallmatrix}
                1 & (m_1 u_1 + m_2 u_2)/d & -b u_1 + a u_2 \\
                0 & 1 & 0 \\
                0 & 0 & 1
            \end{psmallmatrix} l g) e^{-2\pi i(m_1 u_2 + m_2 u_2)} \, du^2 \\
            &= \intl_{(\rats \bs \ads)^2} E(\chi,
            \begin{psmallmatrix}
                1 & x_1 & x_2 \\
                0 & 1 & 0 \\
                0 & 0 & 1
            \end{psmallmatrix} l g) e^{-2\pi i d x_1} \, dx^2 = F_{\psi_U}^{(d, 0)}(\chi, lg) 
        \end{split}
    \end{equation}
    where we have made the substitutions $(m_1 u_1 + m_2 u_2)/d \to x_1$ and $-b u_1 + a u_2 \to x_2$. Note that there are other matrices $l \in L(\rats)$ (explicitly given by $m_1$ and $m_2$) that would accomplish similar results, but if $l \in L(\ints)$ the $p$-adic Iwasawa decomposition simplifies when inserting $g = (g_\infty, \id, \ldots, \id)$ as in section \ref{sec:SL3-non-abelian}.

    Expanding further we obtain
    \begin{equation}
        \label{eq:SL3-maximal-parabolic-further-expansion}
        F_{\psi_U}^{(m_1, m_2)}(\chi, g) = \sum_{m_3 \in \rats} \, \intl_{(\rats \bs \ads)^3} E(\chi,
        \begin{psmallmatrix}
            1 & x_1 & x_2 \\
            0 & 1 & x_3 \\
            0 & 0 & 1
        \end{psmallmatrix} l g) e^{-2\pi i( dx_1 + m_3 x_3)} \, d^3 x = \sum_{m_3 \in \rats} W_{\psi_N}^{(d, m_3)}(\chi, lg)
    \end{equation}
    with $d > 0$ and where $W_{\psi_N}^{(d, m_3)}$ is a Whittaker coefficient on $N$ with charges $d$ and $m_3$ for the two simple roots.

    Now restricting to $\chi_\text{min}$, only the maximally degenerate Whittaker coefficients are non-vanishing. This collapses the sum above to $m_3 = 0$ giving
    \begin{equation}
        \label{eq:SL3-maximal-parabolic-in-min-rep}
        F_{\psi_U}^{(m_1, m_2)}(\chi_\text{min}, g) = W_{\psi_N}^{(d, 0)}(\chi_\text{min}, lg) \, .
    \end{equation}

    We note that for $m_2 = 0$ and positive $m_1$ we have that $d = m_1$ and $l = \id$ giving
    \begin{equation} 
        F_{\psi_U}^{(m_1, 0)}(\chi_\text{min}, g) = W_{\psi_N}^{(m_1, 0)}(\chi_\text{min}, g) \, .
    \end{equation}
    The same statement can be made for negative $m_1$ as well and can be derived by directly making a further expansion as in \eqref{eq:SL3-maximal-parabolic-further-expansion} without a conjugation with $l$. We conclude that, in the minimal representation, a maximal parabolic Fourier coefficient charged only on the simple root $\alpha_1$ simplifies to the maximally degenerate Whittaker coefficient charged on the same root.
    
\end{example}

In \cite{Gustafsson:2014iva} it was similarly shown for $SL(3)$ and $SL(4)$ that all non-trivial Fourier coefficients on any maximal parabolic subgroup automorphic forms in the minimal representation simplify to a single translated maximally degenerate Whittaker coefficient on $N$.

This was accomplished by relating Fourier coefficients to the orbit Fourier coefficients defined in definition \ref{def:orbit-coefficient} and which vanish when the orbit does not belong to the wavefront set of the considered automorphic representation. 
Then, the orbit coefficients were expanded as sums of translated Whittaker coefficients which were found to be maximally degenerate for the minimal orbit coefficients and Whittaker coefficients charged under two strongly orthogonal roots for the next-to-minimal orbit coefficients. 
In the minimal representation, the maximal parabolic Fourier coefficients picked up only one of these maximally degenerate Whittaker coefficients in the minimal orbit coefficient. In the same paper the next-to-minimal representation for $SL(4)$ is discussed as well.

Also, it was shown in \cite{MillerSahi} that, for $E_6$ and $E_7$, Fourier coefficients on certain maximal parabolic subgroups of automorphic forms in $\pi_\text{min}$ are determined by maximally degenerate Whittaker coefficients. From their proof, one may also deduce that, concretely, such a Fourier coefficient is exactly a translate of a maximally degenerate Whittaker coefficients similar to the results of \cite{Gustafsson:2014iva} for $SL(3)$ and $SL(4)$.

From this it was conjectured in \cite{Gustafsson:2014iva} that Fourier coefficients on maximal parabolic subgroups for other simply-laced, simple Lie groups simplify in a similar way and that each may be given in terms of a single, translated, maximally degenerate Whittaker coefficient on $N$. In the remaining parts of this section we will explore some applications and verifications of this statement. 

We have seen in sections \ref{sec:unramified-local-Whittaker-coefficients} and \ref{sec:genpsi} that generic Whittaker coefficients factorise
\begin{equation}
    W_{\psi_N}(\chi, g) = \prod_{p \leq \infty} W_{\psi_{N,p}}(\chi_p, g) \, ,
\end{equation}
with $W_{\psi_{N,p}} \in \text{Ind}^{G(\rats_p)}_{N(\rats_p)} \psi_{N,p}$,
but that degenerate Whittaker coefficients, in general, do not and are expressed as sums of factorising terms as seen in \eqref{degW9}. As such, we cannot expect that all Fourier coefficients on any parabolic subgroup should factorise, that is, we cannot a priori expect that
\begin{equation}
    \text{Ind}^{G(\ads)}_{U(\ads)} \psi_U = \bigotimes_{p \leq \infty} \text{Ind}^{G(\rats_p)}_{U(\rats_p)} \psi_{U,p}
\end{equation}

However, in the minimal representation all but the maximally degenerate Whittaker coefficients on $N$ vanish and the remaining simplify, becoming factorisable as seen for $E_6$, $E_7$ and $E_8$ in the tables of example~\ref{ex:degen-whitt-Ed} from appendix A of \cite{FKP2013}. 
This means that if a Fourier coefficient on a maximal parabolic subgroup can be expressed as a single translated maximally degenerate Whittaker coefficient (as in \eqref{eq:SL3-maximal-parabolic-in-min-rep} above), then it does indeed factorise.

This is interesting since, although not much is known in general about $\text{Ind}^{G(\rats_p)}_{U(\rats_p)} \psi_{U,p}$ for non-minimal parabolic subgroups, the image under the embedding
\begin{equation}
    \pi_{\text{min},p} \subset \text{Ind}^{G(\rats_p)}_{P(\rats_p)} \chi_{\text{min},p} \hookrightarrow \text{Ind}^{G(\rats_p)}_{U(\rats_p)} \psi_{U,p} \, ,
\end{equation}
with $\chi_{\text{min},p}$ spherical, has multiplicity one \cite{MR2123125} and the corresponding spherical vectors $f^\circ_{\psi_{U,p}}$ have been computed in several cases using representation theory \cite{DS, Kazhdan:2001nx, KazhdanPolishchuk, SavinWoodbury}. 

Assuming the factorisation of maximal parabolic Fourier coefficients discussed above, it is possible to rederive and extend these results by considering the spherical vectors induced from $\pi_\text{min}$ as coming from local factors of global Fourier coefficients. 

Indeed, in \cite{Gustafsson:2014iva} it was shown that the products of known spherical vectors $f^\circ_{\psi_{U,p}}$ in maximal parabolics for $E_6$, $E_7$ and $E_8$ give exactly the expected translated maximally degenerate Whittaker coefficients in $\pi_\text{min}$ giving strong support for the claim that \eqref{eq:SL3-maximal-parabolic-in-min-rep} can be generalised to all simply-laced simple Lie groups. Let us consider the case $G = E_7$ below.

\begin{example}[$E_7$ spherical vectors]
    Let $G = E_7$ and $P_{7} = LU$ be the maximal parabolic subgroup obtained by removing the simple root $\alpha_7$ using the Bourbaki labelling in figure \ref{fig:E8dynk}. This was one of the parabolic subgroups studied in \cite{MillerSahi}. Then $U$ is abelian and can also be obtained from the 3-grading
    \begin{equation}
        \lie e_7 = \lie g_{-1} \oplus \lie g_0 \oplus \lie g_1 = \mathbf{27} \oplus (\lie e_6 \oplus \mathbf{1}) \oplus \mathbf{27} \, ,
    \end{equation}
    with $\lie u = \lie g_1$.

    The unique spherical vectors in $\text{Ind}_{U(\rats_p)}^{G(\rats_p)} \psi_{U, p}$ in the minimal representation have been computed at the non-archimedean places by \cite{SavinWoodbury} and are here shown evaluated at the identity in $G(\rats_p)$
    \begin{equation}
        \label{eq:E7-p-spherical-vector}
        f^\circ_{\psi_{U,p}} = \frac{1 - p^3 \abs{m}_p^{-3}}{1 - p^3} \, ,
    \end{equation}
    where $m \in \rats^\times$ is the charge of $\psi_{U,p}$ conjugated to the simple root $\alpha_7$ in $U$.
    
    At the archimedean place we have instead, from \cite{DS}, that
    \begin{equation}
        \label{eq:E7-r-spherical-vector}
        f^\circ_{\psi_{U,\infty}} = m^{-3/2}K_{3/2}(m) \, ,
    \end{equation}
    evaluated at the identity in $G(\reals)$. 

    We will now rederive these results by instead viewing the spherical vectors as coming from local factors of a global Fourier coefficient $F_{\psi_{U}}$ of a spherical Eisenstein series in $\pi_\text{min}$. Such an Eisenstein series may be realised from a parabolically induced representation $I_{P_{1}}(\lambda_\text{min})$ with the maximal parabolic subgroup obtained by removing $\alpha_1$ and $\lambda_\text{min} = 2s\Lambda_1 - \rho$ with $s = 3/2$. We consider the Fourier coefficient $F_{\psi_U}$ with a character non-trivial only on the simple root $\alpha_1$. Similar to the examples above, it simplifies to the single maximally degenerate Whittaker coefficient $W_{\psi_N}^{(\alpha_1)}$ charged only on the same root which, in turn, factorises.

    From table A.1 of \cite{FKP2013} we get that
    \begin{equation}
        \begin{split}
            W_{\psi_N}^{(\alpha_1)}(\lambda_\text{min}, \id) &= \frac{2}{\xi(4)} \abs{m}^{-3/2} \sigma_3(m) K_{3/2}(m) = \frac{2}{\xi(4)} \Bigg( \prod_{p < \infty} \frac{1 - p^3 \abs{m}_p^{-3}}{1 - p^3} \Bigg) \Bigg( \abs{m}^{-3/2} K_{3/2}(m) \Bigg)
        \end{split}
    \end{equation}
    where we recognise the first parenthesis as a product of the non-archimedean spherical vectors in \eqref{eq:E7-p-spherical-vector} and the second as the archimedean spherical vector in \eqref{eq:E7-r-spherical-vector}.
 
    In \cite{Gustafsson:2014iva}, the spherical vectors are rederived in a similar way for $E_6, E_7$ and $E_8$ in both the abelian and Heisenberg realisations of the minimal representation with complete agreement.
\end{example}

\chapter[Hecke theory and automorphic \texorpdfstring{$L$}{L}-functions]{Hecke theory and\\ automorphic \texorpdfstring{$L$}{L}-functions}
\label{ch:Hecke}

In this chapter, we outline the theory of Hecke operators and Hecke algebras. In a nutshell, Hecke operators act on the space of automorphic forms  on a group $G$, forming a commutative ring called the Hecke algebra. The representation theory of this algebra carries a wealth of information about automorphic forms and automorphic representations that connect with many of the structures  discussed in the preceding sections. We begin by outlining the Hecke theory in the case of automorphic forms on real arithmetic quotients $G(\mathbb{Z})\backslash G(\mathbb{R})$, providing detailed examples for the case of $SL(2,\mathbb{R})$. After this treatment of the classical Hecke theory, we consider the counterpart in the adelic context. The key object here is the local spherical Hecke algebra $\mathcal{H}_p^{\circ}$ which acts on the space of adelic automorphic forms $\mathcal{A}(G(\mathbb{Q})\backslash G(\mathbb{A}))$. We study the representation theory of the spherical Hecke algebra and show how this relates to automorphic representations via the Satake isomorphism. Our treatment is mainly done in the context of $SL(2,\mathbb{A})$ and $GL(2,\mathbb{A})$, but many results carry over to arbitrary reductive groups. In particular, in section \ref{sec_Lgroup} we give some details on the  generalisation to $GL(n, \mathbb{A})$ and we make contact to the Langlands program. The starting point is the rewriting of the Casselman--Shalika formula that we encountered in section~\ref{sec:CSLD}. Finally, we end this section with a brief discussion of automorphic $L$-functions, which form a cornerstone of the Langlands program.

\section[Classical Hecke operators and Hecke ring: the general idea]{Classical Hecke operators and Hecke ring:\\ the general idea}
\label{sec:Hecke}

Besides the ring of invariant differential operators there is another set of operators that act on the space $\mathcal{A}(G(\mathbb{Z})\bs G(\reals))$ of automorphic functions on the group $G(\reals)$ invariant under the discrete group $G(\mathbb{Z})$. These additional operators are called \emph{Hecke operators} and we sketch their general definition following~\cite{Goldfeld}.  Their power is worked out for $SL(2,\reals)$ in section~\ref{sec:SL2-Hecke} in the classical setting. Hecke operators and algebras can also be introducted in the adelic setting and this will be the topic of sections~\ref{app_adelichecke} and beyond.

Let $g\in G(\reals)$ be a fixed element \emphindex[commensurable group element]{commensurable} with $G(\ints)$, i.e., the intersection $g^{-1}G(\ints)g\cap G(\ints)$ has finite index in both $G(\ints)$ and $g^{-1}G(\ints)g$. We rewrite its double coset with respect to the discrete group $G(\ints)$ as
\begin{align}
\label{eq:dcosets}
G(\ints) g G(\ints) = \bigcup_{i=1}^d G(\ints)g \delta_i.
\end{align}
On the right-hand side we have written the double coset as a finite disjoint union of single cosets with representatives $\delta_i$ for $i=1,\ldots, d$.  The finiteness of this decomposition follows from the commensurability of $g$.

\begin{definition}
The \emphindex{Hecke operator} $T_g$ associated with a $G(\ints)$-commensurable $g\in G(\reals)$ acting on an automorphic function $\varphi$ is defined by
\begin{align}
\label{Heckedef}
(T_g\varphi)(h) = \sum_{i=1}^d \varphi(g \delta_i h) \quad\quad\textrm{with $h\in G(\reals)$,}
\end{align}
where the $\delta_i$ are representatives of the double coset decomposition~\eqref{eq:dcosets}.
\end{definition}

This operator is well-defined as a finite sum. One can check easily that $T_g$ maps $G(\ints)$-invariant functions to $G(\ints)$-invariant functions. 

\begin{remark}
It is often useful to take a slightly larger group than the original $G(\reals)$ if it acts on the same space. For $SL(2,\reals)$ acting on spherical automorphic functions that are defined on the upper half plane $SL(2,\reals)/SO(2,\reals)$ one can also also consider the action of $GL(2,\reals)$ on $\UHP$ and define Hecke operators for elements $g\in GL(2,\reals)$ with respect to $SL(2,\ints)$. This viewpoint will be useful in section~\ref{sec:SL2-Hecke} below.
\end{remark}

\begin{remark}
\label{rmknorm}The normalisation of the Hecke operators in~\eqref{Heckedef} is not uniquely fixed. The one used there yields the Hecke ring  over $\ints$. It can be useful to change the normalisation and then obtain a Hecke algebra over the field $\rats$.
\end{remark}

The Hecke ring is formed by also allowing integer multiples $m T_g$ of Hecke operators for $m\in\ints$ and defining the product of two Hecke operators $T_{g_1}$ and $T_{g_2}$ by representing the combined double coset $G(\ints)g_1 G(\ints)\cdot G(\ints)g_2 G(\ints)$ as the union of double cosets $G(\ints)h G(\ints)$, possibly with multiplicity. It turns out that a finite union suffices and the product of $T_{g_1}$ and $T_{g_2}$ is then the sum over the $T_{h}$ with integer coefficients. This operation turns the set of Hecke operators into a \emphindex{Hecke ring}. 

The Hecke ring is usually defined together with a given choice of \emphindex[semi-group]{semi-group} $\mathcal{S}$  of commensurable elements $g$. A semi-group is a set with an associative product but not all elements in $\mathcal{S}$ need to be invertible. The example to have in mind here is the set of matrices with determinant equal to some positive integer. The semi-group needs to be chosen such that $G(\ints)$ is a (proper) subgroup of $\mathcal{S}$. Importantly, the Hecke ring (for the cases of interest here) turns out to be \emphindex[Hecke ring!commutative]{commutative}. For the precise statement see~\cite[Thm. 3.10.10]{Goldfeld}. 

Furthermore, the Hecke operators also commute with the ring of differential operators. This means that we can seek common automorphic eigenfunctions of the ring of differential operators and the ring of Hecke operators. The action of the operators then puts additional constraints on the Fourier coefficients that appear in the analysis of the automorphic function and in fact captures much of the number-theoretic structure of these coefficients. An example of this is worked out in the following section for the case of $SL(2,\mathbb{R})$.

\section{Hecke operators for \texorpdfstring{$SL(2,\reals)$}{SL(2, R)}}
\label{sec:SL2-Hecke}

In this section, we illustrate some basic features of the Hecke algebra as sketched in the previous section and the way it interacts with the Fourier expansion in the case of $SL(2,\reals)$. The presentation here is based on~\cite{Goldfeld,Apostol2}.

\subsection{Definition of Hecke operators}

Let $f: \UHP \to \reals$ be a \emphindex{Maass wave form}, i.e., an $SL(2,\ints)$ left-invariant function on the upper half plane $\UHP=SL(2,\reals)/SO(2)$ that is also an eigenfunction of the $SL(2,\reals)$ invariant Laplace operator $\Delta$. defined in~\eqref{LapSL2}.  For example $f$ could be the non-holomorphic Eisenstein series $E(s,z)$ as considered in (\ref{PoincareSL2}). As explained in section~\ref{Maass}, Maass wave forms can also be considered spherical automorphic forms on $SL(2,\reals)$.

According to the general discussion of Hecke operators in section~\ref{sec:Hecke}, we have to choose a semi-group $\mathcal{S}$ of $SL(2,\ints)$ commensurable elements. This we do by letting $\mathcal{S}$ be the group of \emph{diagonal} integer $(2\times 2)$-matrices with positive integer determinant $n$. For fixed $n>0$ let 
\begin{align}
g = \begin{pmatrix} m_1m_2&0\\0&m_2\end{pmatrix}
\end{align}
be a parametrisation of such matrices. We will define a Hecke operator $T_n$ not to a single such element but to the union of all diagonal $g$ with determinant equal to $n$. This $T_n$ can be thought of as the sum of \emph{all} the individual $T_g$ defined according to the formula~\eqref{Heckedef}. According to~\eqref{eq:dcosets},  we require the double coset decomposition into right cosets~\cite[Eq. (3.12.2)]{Goldfeld} 
\begin{align}
\label{dcosetRight}
\bigcup_{m_1^2m_2=n}SL(2,\ints) \begin{pmatrix}m_1m_2&0\\0&m_2\end{pmatrix} SL(2,\ints) = \bigcup_{ad=n\atop 0\leq b <d} SL(2,\ints) \begin{pmatrix} a & b \\ 0&d\end{pmatrix}.
\end{align}
in order to define $T_n$. Then to each $n>0$ we can associate a \index{Hecke operator!for $SL(2,\reals)$} Hecke operator $T_n$ acting on a Maass wave form $f(z)$
\begin{align}
\label{heckedefapostol}
(T_n f)(z) &:= \frac1{\sqrt{n}} \sum_{a\geq 1\,;\,ad=n\atop 0\leq b <d} f\left(\left(\begin{array}{cc}a&b\\0&d\end{array}\right)\cdot z\right)\nn\\
&= \frac1{\sqrt{n}}\sum_{d|n} \sum_{b=0}^{d-1} f\left(\frac{nz+bd}{d^2}\right),
\end{align}
which maps $f$ to a new function $T_n f$ on the upper half plane. Here, we have slightly changed the normalisation of the operator compared to the general discussion as anticipated in remark~\ref{rmknorm}. Note that the transformation of the argument is not in $SL(2,\ints)$ but has determinant $n$. Defining the set
\begin{align}
M_2(n) = \left\{ \begin{pmatrix} a&b \\0 & d\end{pmatrix} \st a,b,d\in \ints \textrm{ with } ad=n\right\}
\end{align}
of upper triangular integer $(2\times 2)$-matrices, we can rewrite the Hecke operator also as
\begin{align}
(T_n f)(z) = \frac{1}{\sqrt{n}} \sum_{\gamma_n \in SL(2,\ints) \backslash M_2(n)} f(\gamma_n\cdot z).
\end{align}
The resulting function $T_nf$ is also a Maass wave form since it is $(i)$ invariant and $(ii)$ an eigenfunction of the Laplacian as we will now show.

$(i)$ Invariance requires evaluating
\begin{align}
\label{hecketrm}
(T_n f)(\gamma \cdot z) = \frac1{\sqrt{n}} \sum_{\gamma_n \in SL(2,\ints) \backslash M_2(n)}f\left(\gamma_n\gamma \cdot z\right)
\end{align}
for $\gamma\in SL(2,\ints)$.  Using
\begin{align}
\gamma_n\gamma =  \tilde{\gamma}\tilde{\gamma}_n 
\end{align}
for some other $\tilde{\gamma}\in SL(2,\ints)$ and $\tilde{\gamma_n}\in M_2(n)$ together with invariance of $\varphi$ under $\tilde{\gamma}$ one arrives at
\begin{align}
(T_n f)(\gamma\cdot z) = \frac1{\sqrt{n}}\sum_{\tilde{\gamma}_n \in SL(2,\ints) \backslash M_2(n)} f(\tilde{\gamma}_n\cdot z) = (T_n f)(z)
\end{align}
and the function $T_nf$ is $SL(2,\ints)$ invariant for any positive $n$. (See also chapters 6.8 and 6.9 of~\cite{Apostol2}.)

$(ii)$ Consider the action of the Laplacian (\ref{LapSL2}) on $T_n \varphi$. It is straight-forward to check that
\begin{align}
[\Delta (T_n f)](z) =[ T_n (\Delta f)](z).
\end{align}
Therefore, \emph{the Laplacian commutes with all the Hecke operators} and if $f$ is an eigenfunction of $\Delta$, so is $T_nf$ and with the same eigenvalue.\\

Finally, we study the Fourier expansion of $T_n f$. Suppose that $f$ has a Fourier expansion of the form (cf.~\eqref{EisExpSec4})
\begin{align}
f(z) = f(x+iy) = \sum_{m\in\ints} a_m(y) e^{2\pi i m x }.
\end{align}
Then one finds for $T_n \varphi$
\begin{align}
(T_n f)(z) &= \frac1{\sqrt{n}} \sum_{d|n} \sum_{b=0}^{d-1} f\left(\frac{n}{d^2}x+ \frac{b}{d}+i\frac{n}{d^2}y\right)\nn\\
&= \frac1{\sqrt{n}}\sum_{d|n} \sum_{m\in\ints} a_m\left(\frac{n}{d^2}y\right) e^{2\pi i m n x /d^2} \sum_{b=0}^{d-1} e^{2\pi i m b/d}\nn\\
&= \frac1{\sqrt{n}} \sum_{m\in\ints} \sum_{d|n, d|m} d\, a_m\left(\frac{n}{d^2}y\right) e^{2\pi i m n x/d^2}\nn\\
&= \frac1{\sqrt{n}} \sum_{m\in \ints} \sum_{d|(n,m)} \frac{n}{d} a_{mn/d^2}\left(\frac{d^2}{n}y\right) e^{2\pi i m x}\,,
\end{align}
where we have changed the divisor sum variable from $d$ to $\frac{n}{d}$ in the last step and have relabelled the $m$ sum in between. The Fourier expansion of $T_n f$ is therefore
\begin{align}
\label{fourierhecke}
(T_n f)(z) = \sum_{m\in\ints} \tilde{a}_m(y) e^{2\pi i m x}\quad\text{with}\quad
\tilde{a}_m(y) =  \frac1{\sqrt{n}}  \sum_{d|(n,m)} \frac{n}{d} a_{mn/d^2} \left(\frac{d^2}{n}y\right)\,,
\end{align}
where $d|(n,m)$ means that $d|n$ and $d|m$, i.e. divides the greatest common divisor of $n$ and $m$ which is denoted as usual by $(m,n)=\gcd(m,n)$.

\subsection{Algebra of Hecke operators}

Importantly, the Hecke operators $T_n$ satisfy a simple algebra on the space of Maass wave forms: {\em they all commute}. Moreover, they satisfy the \emphindex[Hecke algebra!for $SL(2,\reals)$]{Hecke algebra}
\begin{align}
\label{heckealgebra}
T_m T_n = T_nT_m = \sum_{d|(m,n)} T_{mn/d^2}\,.
\end{align}
Commutativity is manifest in this expression. To prove (\ref{heckealgebra}), one can first consider the case $(m,n)=1$ and use the explicit definition. In the next step one can consider the case when both $m$ and $n$ are powers of the same prime. A proof can be found in~\cite[ch. 6.10]{Apostol2} where a different normalisation is used.

\subsection{Common eigenfunctions of \texorpdfstring{$T_n$}{Tn} and \texorpdfstring{$\Delta$}{Delta}}
\label{sec:SL2-eigenform}
\index{Hecke operator!common eigenfunctions}

Suppose $f:\UHP\to \reals$ is an eigenfunction of all Hecke operators
\begin{align}
T_n f = c_n f
\end{align}
for some eigenvalues $c_n$ and at the same time a cuspidal Maass wave form with Fourier expansion
\begin{align}
f(z) = \sum_{m\neq 0} a_m(y) e^{2\pi i m x}\,.
\end{align}

Applying $T_n$ to $f$ gives, with (\ref{fourierhecke}),
\begin{align}
c_n f(z) = \frac1{\sqrt{n}} \sum_{m\neq 0} \sum_{d|(n,m)} \frac{n}{d} a_{mn/d^2}\left(\frac{d^2}{n}y\right)e^{2\pi i m x}
\end{align}
Comparing the individual Fourier modes on both sides leads to
\begin{align}
c_n a_m(y) = \frac1{\sqrt{n}} \sum_{d|(n,m)} \frac{n}{d} a_{nm/d^2}\left(\frac{d^2}{n}y\right)
\end{align}
Setting $n=1$ gives $c_1 a_m(y)= a_m(y)$ for all $m\neq 0$, implying $c_1=1$ unless $f$ vanishes. Setting $m=1$ implies 
\begin{align}
\label{HeckeF}
c_n a_1(y) = \sqrt{n} a_n\left(\frac{y}{n}\right).
\end{align} 
If $f$ is not constant, one has $a_1(y)\neq 0$, otherwise all the Fourier coefficients would vanish. 

{}From solving the Laplace condition on $\Delta f=s(s-1)f$ (cf. appendix~\ref{app:SL2Laps}) one knows that the dependence of the Fourier coefficient $a_m(y)$ on $y$ is through the modified Bessel function as
\begin{align}
\label{nonholFCoeffs}
a_m(y) = a_m y^{1/2} K_{s-1/2} (2\pi|m|y)
\end{align}
for some purely numerical coefficient $a_m$ that we will now relate to the Hecke eigenvalues $c_n$.
Rescaling $f$ such that $a_1=1$ (\emphindex{Hecke normalisation}) the relation~\eqref{HeckeF} above implies that the Hecke eigenvalues equal the Fourier coefficients:
\begin{align}
\label{heckefourier}
c_n = a_n.
\end{align}
Obtaining this simple relation was the reason for the choice of normalisation of the Hecke operator $T_n$. \index{Hecke eigenvalues!relation to Fourier coefficients} 
Note that the $c_n$ are only defined for positive $n$ but $a_n$ for any $n$. The reality of $f$ relates $a_n$ to $a_{-n}$. 

By virtue of the Hecke algebra (\ref{heckealgebra}) we have
\begin{align}
T_m T_n f = c_m c_n f = \sum_{d|(m,n)}  c_{mn/d^2} f
\end{align}
so that the Fourier coefficients of a normalised simultaneous eigenfunction satisfy
\begin{align}
\label{heckealgebracoeff}
a_m a_n = \sum_{d|(m,n)}  a_{mn/d^2}\,.
\end{align}
In particular, they must be \index{Fourier coefficient!multiplicative}\emp{multiplicative}, i.e. for coprime $m$ and $n$ one has $a_m a_n = a_{mn}$. This \index{Fourier coefficient!and number theory}number-theoretic property of the Fourier coefficients follows from the action of the Hecke operators and would not have been apparent from $SL(2,\ints)$ invariance alone. Note that the constant term is not captured by these considerations.

The algebra~\eqref{heckealgebracoeff} allows determining \emph{all} Fourier coefficients in terms of the ones for prime numbers $a_p$. We note for later reference that powers of primes can be calculated recursively using the relation
\begin{align}
\label{FCprimepower}
a_{p^{k+1}} = a_{p^k} a_p - a_{p^{k-1}} 
\end{align}
for $k>1$ where again Hecke normalisation $a_1=1$ enters. The Hecke operators $T_p$ determine the full structure of the Hecke algebra and hence are the only relevant ones for the development of the theory. We will see soon that they fit naturally into an adelic framework.

\begin{example}[{Fourier expansion of non-holomorphic Eisenstein series and Hecke algebra}]
\label{sec:HeckeFourier}

We now use the Hecke algebra to rederive the Fourier expansion (\ref{SL2FC2}) of the Eisenstein series $E(s,z)$.
Since the Eisenstein series $E(s,z)$ is defined as a sum over an $SL(2,\ints)$ orbit it is easy to evaluate  the Hecke operators  by multiplying the acting matrices
\begin{align}
\label{poincareheckeevals}
(T_n E)(s,z) &= \frac12 \frac{1}{\sqrt{n}} \sum_{d|n} \sum_{b=0}^{d-1} \sum_{\gcd(p,q)=1} \left[\Im\left(\left(\begin{array}{cc}n/d&b\\0&d\end{array}\right)\left(\begin{array}{cc}*&*\\p&q\end{array}\right)\cdot z\right)\right]^s\nn\\
&= \frac{1}{\sqrt{n}}\sum_{d|n}\sum_{b=0}^{d-1} \left(\frac{n}{d^2}\right)^s E(s,z)
= \sum_{d|n} \left(\frac{n}{d^2}\right)^{s-1/2} E(s,z)\nn\\
&= \underbrace{n^{s-1/2} \sigma_{1-2s}(n)}_{c_n} E(s,z)\,.
\end{align}
We have used that the coset sum $B(\ints)\backslash SL(2,\ints)$ can be parametrised by two co-prime integers $p$ and $q$ and the unspecified top row corresponds to an arbitrary representative of the coset. 
In particular, the Eisenstein series is an eigenfunction of all Hecke operators and the relation (\ref{heckefourier}) between the Fourier coefficients and the Hecke eigenvalues immediately implies the form (\ref{SL2FC2}) for the non-zero Fourier coefficients up to a normalisation factor. The constant term is not fixed by these considerations. However, this method of deriving the Fourier modes did not require any Poisson resummation nor adelic technology.

Let us verify the relation~\eqref{heckealgebracoeff} for the explicit example of the Eisenstein series \eqref{SL2FC2} to check whether it is a simultaneous eigenfunction. There one has for $n>0$
\begin{align}
\label{HeckeNH}
        a_n = c_n= \sum_{d|n} d^{1-2s} n^{s-1/2}
\end{align}
where Hecke normalisation was used. Let $m$ and $n$ be coprime, then
\begin{align}
        a_m a_n =\sum_{d|m}\sum_{\tilde{d}|n} d^{1-2s} m^{s-1/2} \tilde{d}^{1-2s} n^{s-1/2} 
        = \sum_{d|mn} d^{1-2s} (mn)^{s-1/2} = a_{mn}\,.
\end{align}
The more general relation (\ref{heckealgebracoeff}) can also be verified and the Eisenstein series is an eigenfunction of the Hecke operators (with eigenvalues given by the Fourier coefficients). 
\end{example}

\begin{remark}[Hecke operators for holomorphic modular forms]
For holomorphic modular forms $f:\UHP\to\cx $ of weight $k$ one can also define Hecke operators, see for example~\cite[Ch. 6]{Apostol2}. In this case, they act by \index{Hecke operator!for holomorphic modular forms}
\begin{align}
\label{HeckealgebraMod}
(T_nf)(z) = n^{k-1} \sum_{d|n} d^{-k} \sum_{b=0}^{d-1} f\left(\frac{bz+bd}{d^2}\right) \quad\quad\textrm{($f$ holomorphic of weight $k$)}
\end{align}
and map holomorphic modular forms to homomorphic modular forms. Note that the normalisation convention here is slightly different from the non-holomorphic case. The multiplicative law~\eqref{heckealgebra} in this case reads
\begin{align}
T_m T_n =T_n T_m = \sum_{d|(m,n)} d^{k-1} T_{mn/d^2}. \quad\quad
\textrm{(weight $k$ Hecke algebra)}
\end{align}
One can again define Hecke normalised common eigenfunctions. If $f$ is a common eigenfunction with the Fourier expansion $f(z) = \sum_{m\geq 0} a_m q^n$ (with $q=e^{2\pi i z}$ as always), then one has again
\begin{align}
T_n f = a_n f
\end{align}
when $a_1=1$, i.e., the modular form is Hecke normalised. In this case the Fourier coefficients satisfy
\begin{align}
\label{holFC}
a_m a_n = \sum_{d|(m,n)} d^{k-1}  a_{mn/d^2} \quad\quad\textrm{($f$ holomorphic of weight $k$)}
\end{align}
because of~\eqref{HeckealgebraMod}. If $f$ is a cusp form and its Fourier coefficients satisfy the above relation~\eqref{holFC} one can show that it is automatically a common eigenfunction~\cite[Thm. 6.15]{Apostol2}. For non-cuspidal forms this is not guaranteed. We record the following consequence of~\eqref{holFC} for later use
\begin{align}
\label{FCprimepowerHol}
a_{p^\ell} a_p = a_{p^{\ell+1}} + p^{k-1} a_{p^{\ell-1}} \quad\quad \textrm{(Fourier coefficients of weight $k$ modular form)}
\end{align}
for $\ell\geq 0$. This is to be contrasted with~\eqref{FCprimepower} for non-holomorphic forms.
\end{remark}

\section{Hecke operators and Dirichlet series}
\label{app_adelichecke}

Given the powerful applications of Hecke operators demonstrated in the previous subsections, it is natural to wonder about the action of Hecke operators in the adelic setting of 
 automorphic forms on $SL(2,\mathbb{Q})\backslash SL(2,\mathbb{A})$. It turns out that this gives rise to an even richer structure, and provides a link to the theory of automorphic representations. In this section, we take the first steps toward such a theory by studying the Hecke operators $T_p$ for $p$ a prime, based on~\cite{Bump,Goldfeld}. 

\index{Hecke operator!and Dirichlet series}
Hecke's original motivation to study Hecke operators was to find a way to encode the  properties of a holomorphic modular form in terms of its associated \emphindex{Dirichlet series}~\cite{Hecke1,Hecke2}. Given a weight $k$ modular form $f(z)=\sum_{n\geq 0} a_n(f)q^n$ (with $q=e^{2\pi i z}$) one may form the series 
\beq
\label{HeckeL}
L(s, f)=\sum_{n\geq 1} \frac{a_n(f)}{n^s} = \prod_{p<\infty} \sum_{\ell\geq 0} \frac{a_{p^\ell}(f)}{p^{\ell s}},
\eeq
which is called the Dirichlet series attached to $f$.  The rewriting in the second step is the application of prime factorisation under the assumption of absolute convergence of the $L$-series. In the special case when the Fourier coefficients $a_m(f)$ are \emphindex[multiplicative!completely]{completely multipliciative}, i.e. satisfy $a_{m}a_n=a_{mn}$ for any $m,n\in \mathbb{Z}$, then the $L$-function leads to the following Euler product via geometric series
\beq
L(s,f)=\prod_{p<\infty} \frac{1}{1-a_p(f) p^{-s}} , \qquad \quad a_m(f) \, \, \text{completely multiplicative}.
\eeq 
This is called a degree $1$ Euler product since the denominator contains at most the power $p^{-s}$. 
The prime example of a degree 1 Euler product is the Riemann zeta function $\zeta(s)=\prod_{p<\infty}(1-p^{-s})^{-1}$, corresponding to $a_m=1$ for all $m\geq 1$, which, however, is not associated with a holomorphic modular form on $SL(2,\reals)$ but rather with $GL(1,\ads)$ as was explained in section~\ref{RiemannZetaEx}.

Hecke showed that whenever the Fourier coefficients $a_m(f)$ are multiplicative according to~\eqref{FCprimepowerHol}
then the Dirichlet series can be written as an Euler product
\beqa
\label{LEulerProd}
L(s,f)
=\prod_{p<\infty} L_p(s,f) 
=\prod_{p<\infty} \frac{1}{1-a_p(f) p^{-s} + p^{k-1-2s}}.
\eqa
The derivation of this formula is as follows. Let $L_p(s,f)$ be a factor in the Euler product as above. Then~\eqref{FCprimepowerHol} implies
\begin{align}
L_p(s,f) = \sum_{\ell\geq 0} \frac{a_{p^\ell}}{p^{\ell s} }
&= \frac{1}{a_p} \left[ p^s \sum_{\ell\geq 0} \frac{a_{p^{\ell+1}}}{p^{(\ell+1)s}} + p^{k-1-s} \sum_{\ell\geq 0} \frac{a_{p^{\ell-1}}}{p^{(\ell-1)s}}
\right]\nn\\
&= \frac{1}{a_p}\left[  p^s \left(L_p(s,f)-1\right) + p^{k-1-s} L_p(s,f)\right],
\end{align}
which yields~\eqref{LEulerProd} after solving for $L_p(s,f)$. The series $L(s,f)$ is also called the \emph{L-function} of $f$ and $L_p(s, f)$ the \emph{local L-factor}. The $L$-function $L(s, f)$ in (\ref{LEulerProd}) is of degree 2, due to the factor $p^{-2s}$.

\begin{remark}
In the case of non-holomorphic automorphic forms $f$ one can go through the same derivation of an $L$-function. Using the normalisation of Hecke operators defined in~\eqref{heckedefapostol} and the Fourier coefficients $a_m$ defined in~\eqref{nonholFCoeffs} one obtains an $L$-function for a common eigenfunction $f$ in Hecke normalisation of the form
\begin{align}
\label{LNH}
L(s,f) = \prod_{p<\infty} \frac1{1-a_p p^{-s} + p^{-2s}}.
\end{align}
The shifted exponent on the last term in the denominator is due to the normalisation of the Hecke operators. We assume for simplicity that $f$ is even, i.e., $f(-z)=f(z)$. Then one can define a \emphindex[L-function@$L$-function!completed]{completed $L$-function} via
\begin{align}
\label{compL}
L^\compl(s,f) = \pi^{-s} \Gamma\left(\frac{2s+2\nu-1}{4}\right) \Gamma\left(\frac{2s-2\nu+1}{4}\right) L(s,f)
\end{align}
where $\nu$ is the eigenvalue under the Laplacian $\Delta f= \nu(\nu-1) f$. The completed $L$-function satisfies the simple functional relation
\begin{align}
L^\compl(s,f) = L^\compl(1-s,f)
\end{align}
For a proof of this and extensions to odd Maass forms see~\cite[Prop. 3.13.5]{Goldfeld}. One should think of the normalising factors in~\eqref{compL} as arising from the archimedean place $p=\infty$ and the completed $L$-function as a global one.
\end{remark}

The $L$-function (\ref{LEulerProd}) attached to a modular form $f$ therefore characterizes whether or not the Fourier coefficients exhibit a multiplicative behaviour, something which is certainly not 
guaranteed. When does this happen? It turns out that the Fourier coefficients of a modular form $f$ are multiplicative if and only if $f$ is a \emphindex{Hecke eigenform}, i.e., an eigenfunction of the entire ring of  Hecke operators $T_n$~\cite[Thm 6.15]{Apostol2}. As was emphasised above, the ring of Hecke operators is generated by the $T_p$ for $p$ prime and we will now focus on these.

\begin{remark}
Weil~\cite{Weil1,Weil2} has resolved the problem of generalising the $L$-function to automorphic forms for congruence subgroups $\Gamma_0(N)$ of $SL(2,\ints)$. In this case one needs to twist the $L$-function by a \index{Dirichlet character} Dirichlet character.
\end{remark}

\section{The spherical Hecke algebra}

Recall from section \ref{sec_frommodtoaut} that to each modular form $f(z)$ on the upper-half plane $\mathbb{H}$ we have  a corresponding automorphic form $\varphi_f(g)$ on $SL(2,\mathbb{Q})\backslash SL(2,\mathbb{A})$. We now want to find out how  the action of the Hecke operator $T_p$ lifts to the space of automorphic forms $\mathcal{A}(SL(2,\mathbb{Q})\backslash SL(2,\mathbb{A}))$. 

As for the classical case in section~\ref{sec:SL2-Hecke}, the Hecke operators in the adelic context are associated with double cosets of matrices of determinant different from $1$ and hence outside of $SL(2,\rats_p)$. For this reason, we consider the group $GL(2,\rats_p)$. The \emphindex{convolution algebra} on $GL(2,\mathbb{Q}_p)$. This algebra is given by the space of locally constant $\mathbb{C}$-valued functions on $GL(2,\mathbb{Q}_p)$ with the (commutative) product given by convolution:
\beq
\label{convProd}
\left(\Phi_1\star \Phi_2\right)(g) =\int_{GL(2,\mathbb{Q}_p)}\Phi_1(gh) \Phi_2(h^{-1}) dh=\int_{GL(2,\mathbb{Q}_p)} \Phi_1(h)\Phi_2(h^{-1} g) dh, 
\eeq
where $dh$ denotes the bi-invariant Haar measure on the uni-modular group $GL(2,\mathbb{Q}_p)$. Convolution turns the space of such functions into a ring, called the \emphindex[Hecke algebra!local]{(local) Hecke algebra}, commonly denoted by $\mathcal{H}(GL(2,\mathbb{Q}_p))$ or simply $\mathcal{H}_p$ for short. Although it is a ring it has no unit. To see the connection with the classical Hecke algebra generated by the $T_p$'s, we now restrict to \emphindex[function!bi-invariant]{bi-invariant functions} with respect to the maximal compact subgroup $K_p=GL(2,\mathbb{Z}_p)$, i.e. we consider functions in $\mathcal{H}_p$ that satisfy
\beq
\Phi(k g k')=\Phi(g), \qquad k, k'\in K_p, \, g\in GL(2,\mathbb{Q}_p).
\eeq
We then obtain the \emphindex[Hecke algebra!spherical]{spherical Hecke algebra} $\mathcal{H}(GL(2,\mathbb{Q}_p))^{K_p}$ of $K_p$ bi-invariant functions, which we denote by $\mathcal{H}_p^{\circ}$. It is a central result that $\mathcal{H}_p^{\circ}$ forms a commutative ring (see, e.g. \cite{Bump}). If we fix the Haar measure on $GL(2,\mathbb{Q}_p)$ such that $K_p$ has unit volume, then $\mathcal{H}_p^{\circ}$ also has a unit given by the characteristic function on $K_p$:
\begin{equation}
    \text{char}_{K_p}(g) =
    \begin{cases}
        1 & g\in K_p,\\
        0 & \text{otherwise} .\\
    \end{cases}
\end{equation}
To see this we  calculate the convolution product of the characteristic function with any  $\Phi\in \mathcal{H}_p^{\circ}$:
\begin{align}
(\Phi\star \text{char}_{K_p})(g) =\int_{SL(2,\mathbb{Q}_p)} \Phi(gh) \text{char}_{K_p}(h^{-1})dh
= \int_{K_p} \Phi(gh)dh 
= \Phi(g), 
\end{align}
where in the last step we used that $f$ is bi-invariant under $K_p$ and $K_p$ has unit volume. One says that the spherical Hecke algebra is \emphindex{idempotented}.

The spherical Hecke algebra $\mathcal{H}_p^\circ$ acts on the space of $K_p$-spherical functions on $GL(2,\rats_p)$ via right-translation. For any $\Phi\in \mathcal{H}_p^\circ$ and $K_p$-spherical function $\varphi$ on $GL(2,\rats_p)$ we define a new function on $GL(2,\rats_p)$ by
\begin{align}
(\pi(\Phi) \varphi)(g)=\lint_{GL(2,\mathbb{Q}_p)} \Phi(h) \varphi(gh)  dh.
\label{rightregularHeckeaction}
\end{align}
One can check easily that this maps the right-regular action of $GL(2,\rats_p)$ on functions on $GL(2,\rats_p)$ to a representation of the spherical Hecke algebra (with convolution product~\eqref{convProd}) according to 
\begin{align}
\pi\left(\Phi_1\star \Phi_2\right)\varphi = \pi\left(\Phi_1\right)\left(\pi\left(\Phi_1\right)\varphi\right).
\end{align}
The space of $K_p$-spherical functions is therefore a representation of the spherical Hecke algebra $\mathcal{H}^\circ_p$.

\begin{remark}
\label{rmk:globalHecke}
By taking the restricted direct product (with respect to $K_p$, see section \ref{sec_adeles}) over all the local algebras $\mathcal{H}_p^{\circ}$ we obtain the \emphindex[Hecke algebra!global spherical]{global, or adelic, spherical Hecke algebra}
\beq
\mathcal{H}^{\circ}={\bigotimes_{p\leq\infty}}' \mathcal{H}_p^{\circ}. 
\eeq
For $p=\infty$, the spherical Hecke algebra $\mathcal{H}^\circ_\infty$ is given by $K(\reals)$-bi-finite distributions supported on $K(\reals)$~\cite[Lecture 3.1]{Cogdell}. This includes the invariant differential operators on $G(\reals)$ lying in the universal enveloping algebra $U(\mf{g})$.  
The global Hecke algebra $\mathcal{H}^\circ$ acts on $\mathcal{A}(GL(2,\mathbb{Q})\backslash GL(2,\mathbb{A}))$ by the same formula (\ref{rightregularHeckeaction}). Our main interest in the following lies with the spherical Hecke algebra $\mathcal{H}^\circ_p$ at the finite primes $p<\infty$. 
\end{remark}

We now investigate the structure of the (local) spherical Hecke algebra $\mathcal{H}_p^\circ$ in more detail.  More explicitly, we define the elements $\mathbb{T}_p$ and $\mathbb{R}_p\in \mathcal{H}_p^\circ$ by the $K_p$-bi-invariant functions
\begin{align}
\label{pHecke}
\mathbb{T}_p = \textrm{char}_{K_p \left(\begin{smallmatrix} p & \\ & 1\end{smallmatrix}\right)K_p},\quad\quad
\mathbb{R}_p = \textrm{char}_{K_p \left(\begin{smallmatrix} p & \\ & p\end{smallmatrix}\right)K_p}.
\end{align}
It is an important result that $\mathbb{T}_p$, $\mathbb{R}_p$ and $\mathbb{R}_p^{-1}$ together generate the spherical Hecke algebra $\mathcal{H}_p^{\circ}$. A proof of this statement can be found for example in \cite[Prop. 4.6.5]{Bump}.

On functions $\varphi : GL(2,\rats_p)\to \cx$ they act according to~\eqref{rightregularHeckeaction}.  To ease notation we shall simply continue to call them $\mathbb{T}_p$ and $\mathbb{R}_p$ also when acting on spherical functions.
\begin{align}
(\mathbb{T}_p \varphi)(g) &= \int_{K_p \left(\begin{smallmatrix} p & \\ & 1\end{smallmatrix} \right)K_p}\varphi(gh)dh,\\
(\mathbb{R}_p \varphi)(g) &= \int_{K_p \left(\begin{smallmatrix} p & \\ & p\end{smallmatrix} \right)K_p}\varphi(gh)dh.
\end{align}

Even though written in terms of integrals, they act on functions on $GL(2,\rats_p)$ by finite sums after performing a decomposition of the double cosets into a finite union of left cosets, similar to~\eqref{dcosetRight}. This decomposition for the operator $\mathbb{T}_p$ is~\cite[Prop. 4.6.4]{Bump}
\begin{align}
\label{TTp}
K_p \left(\begin{array}{cc} p & \\ & 1\end{array} \right)K_p=\left(\begin{array}{cc} 1 & \\ & p\end{array} \right)K_p \cup \bigcup_{i=0}^{p-1} \left(\begin{array}{cc} p & i\\ & 1\end{array} \right) K_p,
\end{align}
such that for $K_p$-spherical $\varphi$
\begin{align}
(\mathbb{T}_p \varphi)(g) = \varphi\left(g\left(\begin{smallmatrix}1&\\&p\end{smallmatrix}\right)\right) + \sum_{i=0}^{p-1} \varphi\left(g\left(\begin{smallmatrix}p&i\\&1\end{smallmatrix}\right)\right).
\end{align}

The connection with the classical Hecke operators now follows from the fact that if $f:\UHP\to\reals$ is a Maass wave form   with eigenvalue $a_p$ under $T_p$, then the associated adelic automorphic form $\varphi_f\in \mathcal{A}(SL(2,\mathbb{Q})\backslash SL(2,\mathbb{A}))$ defined in section~\ref{Maass} is an eigenform under $\mathbb{T}_p$ with the same eigenvalue, up to a (convention-dependent) factor: 
\beq
T_p f=a_pf\, \,  \longleftrightarrow\, \,  (\mathbb{T}_p\varphi_f)(g)=p^{1/2} a_p\varphi_f.
\label{eigenvaluerelation}
\eeq
This will be verified for Eisenstein series in example~\ref{ex:HeckeConnection} below but it is valid in general.

\begin{remark}
There is also a classical Hecke operator $R_p$ acting on $f$ which lifts to $\mathbb{R}_p$, but we shall not discuss this further here (see \cite{Bump} for more details). 
\end{remark}

\section[Spherical Hecke algebras and automorphic representations]{Spherical Hecke algebras and automorphic\\ representations}

This and the following sections make use of the theory of automorphic representations which was introduced in section~\ref{ch:autoreps}.

Recall that $GL(2,\mathbb{A})$ acts by right-translation on  $\mathcal{A}(GL(2,\mathbb{Q})\backslash GL(2, \mathbb{A}))$, such that at the archimedean places it has the form of a $(\mathfrak{g}_\infty, K_\infty)$-module, while the finite places carry a representation of $GL(2,\mathbb{A}_f)$. The irreducible constituents $(\pi, V)$ in the decomposition of $\mathcal{A}(GL(2,\mathbb{Q})\backslash GL(2, \mathbb{A}))$ are called automorphic representations. But we have also just seen that $\mathcal{A}(SL(2,\mathbb{Q})\backslash SL(2, \mathbb{A}))$ carries an action of the adelic spherical Hecke algebra $\mathcal{H}^{\circ}$. A natural question is then: Is there a relation between these representations? Not surprisingly, the answer is yes, and we shall now sketch how to see this.

Suppose that $(\pi, V)=\otimes_{p\leq \infty} (\pi_p, V_p)$ is an unramified automorphic representation (see definition \ref{def_unramifiedreps}); this implies that $V_p$ contains a spherical vector  ${\sf f}^{\circ}_p$ (unique up to multiplication by a complex scalar, see e.g.~\cite{Bump}), satisfying ${\sf f}^{\circ}_p(k)=1$ for all $k\in K_p$. The spherical vector therefore spans the complex one-dimensional space $V_p^{K_p}$ consisting of $K_p$-invariant vectors in $V_p$. 

We can for example take $\pi_p$ to be the local induced representation with module
\beq
V_p=\text{Ind}_{B(\mathbb{Q}_p)}^{GL(2,\mathbb{Q}_p)} \delta^{1/2} \mu ,
\label{unramifiedrepsSL2}
\eeq
where $\delta$ is the modulus character of the Borel subgroup and the quasi-character $\mu \, :\, B(\mathbb{Z}_p)\backslash B(\mathbb{Q}_p)\to \mathbb{C}^\star$ is defined by 
\beq
\mu(b)=\mu(na)=\mu(a), \qquad n\in N(\mathbb{Q}_p), \quad a\in A(\mathbb{Q}_p).
\eeq
In the notation of~\eqref{eq:sepmod} we have therefore $\chi(g)=\mu(g) \delta^{1/2}(g)$. The explicit separation of the modulus character in \eqref{unramifiedrepsSL2} turns out to be convenient for 
the forthcoming analysis, and also facilitates comparison with the literature. In the notation of that section we would write $\chi(a)=a^{\lambda+\rho}$, so that $\delta^{1/2}(a)=a^{\rho}$ and $\mu(a)=a^{\lambda}$, where $\rho$ is the Weyl vector of the Lie algebra $\mathfrak{gl}(2)$ and $\lambda$ is a (complex) weight. 

The spherical vector ${\sf f}^\circ_p\in V_p$ is the standard section defined by the extension of $\delta^{1/2}\mu$ to all of $GL(2,\mathbb{Q}_p)$ via the Iwasawa decomposition (see also section \ref{standardsectionSL2}):
\beq
{\sf f}^{\circ}_p(g)={\sf f}^{\circ}_p(nak)=\delta^{1/2}(a) \mu(a).
\label{sphericalvectorSL2}
\eeq
The local spherical Hecke algebra $\mathcal{H}_p^{\circ}$ acts on $V_p$ via the action (\ref{rightregularHeckeaction}). By construction this action preserves the one-dimensional space $V_p^{K_p}$ of spherical vectors: indeed for any $\Phi \in \mathcal{H}_p^{\circ}$ we have for all $k\in K_p$
\beq
(\pi(\Phi){\sf f}_p^{\circ})(gk)=\lint_{GL(2,\mathbb{Q}_p)}{\sf f}_p^{\circ}(gkh)\Phi(h)dh
=\lint_{GL(2,\mathbb{Q}_p)}{\sf f}_p^{\circ}(ghk)\Phi(khk^{-1})dh
=(\pi(\Phi){\sf f}_p^{\circ})(g),
\eeq
since ${\sf f}_p^{\circ}$ is spherical and $\Phi$ $K_p$-bi-invariant. This implies that $V_p^{K_p}$ furnishes a representation of $\mathcal{H}_p^{\circ}$. Since the spherical vector ${\sf f}_p^{\circ}$ spans the one-dimensional space $V_p^{K_p}$, we conclude that the action of $\mathcal{H}_p^{\circ}$ must give back ${\sf f}_p^{\circ}$, up to a complex scalar:
\beq
(\pi(\Phi) {\sf f}_p^{\circ})(g)
= \lambda_{\mu}(\Phi) {\sf f}_p^{\circ}(g),
\label{sphericalHeckeeigenvalue}
\eeq
where the eigenvalue $\lambda_{\mu}(\Phi)$ determines a (quasi-)character of the spherical Hecke algebra
\beq
\lambda_{\mu} \, :\, \mathcal{H}_p^{\circ} \longrightarrow \mathbb{C}^{\times}.
\eeq
As we have indicated, this character depends on the choice of $\mu$ in~\eqref{unramifiedrepsSL2}.

To find an explicit description of the characters $\lambda_\mu$ we shall work out the action of the Hecke operator $\mathbb{T}_p$ defined in~\eqref{pHecke}. 
To proceed, we parametrise the Cartan torus $A(\mathbb{Q}_p)\subset GL(2,\mathbb{Q}_p)$ by 
\beq
a=\left(\begin{array}{cc} v_1 & \\ & v_2 \end{array}\right), \qquad v_1, v_2\in \mathbb{Q}_p^{\times}. 
\eeq
We can further describe the unramified character $\mu$ 
explicitly by
\beq
\mu(a) = \mu\left(\left(\begin{array}{cc} v_1 & \\ & v_2 \end{array}\right)\right)
=|v_1|_p^{s_1}|v_2|_p^{s_2}
\label{unramifiedcharacterGL2}
\eeq
where $s_1, s_2\in \mathbb{C}$. Note that the parametrisation in terms of $s_1$ and $s_2$ differs from the one used in chapter~\ref{ch:SL2-fourier}. The reason here is to simplify some of the following expressions. 
The corresponding value of the modulus character on $B\subset GL(2,\mathbb{Q}_p)$ is 
\beq
\delta\left(\begin{array}{cc} v_1 & \\ & v_2 \end{array}\right)=\left|\frac{v_1}{v_2}\right|_p.
\eeq
This implies that the representation $\text{Ind}_{B(\mathbb{Q}_p)}^{SL(2,\mathbb{Q}_p)} \delta^{1/2} \mu$ is in fact completely determined by 
\beq
\label{abvals}
\alpha_p\equiv p^{-s_1}, 
\qquad \beta_p\equiv p^{-s_2}. 
\eeq

Now we wish to compute 
\beq
(\mathbb{T}_p{\sf f}_p^{\circ})(g)=\int_{K_p \left(\begin{smallmatrix} p & \\ & 1\end{smallmatrix} \right)K_p} {\sf f}_p^{\circ}(gh)dh. 
\eeq
Since we know that ${\sf f}_p^{\circ}$ is an eigenfunction and  is normalised so that ${\sf f}_p^{\circ}(1)=1$, it suffices to evaluate the action at the identity $1\in GL(2,\mathbb{Q}_p)$, which then, via (\ref{sphericalHeckeeigenvalue}), directly corresponds to the value of the character $\lambda_\mu$:
\beq
\lambda_{\mu}(\mathbb{T}_p)=\int_{K_p \left(\begin{smallmatrix} p & \\ & 1\end{smallmatrix} \right)K_p} {\sf f}_p^{\circ}(h)dh.
\label{heckeintegral}
\eeq
To evaluate this  we decompose the double coset space as in~\eqref{TTp}.
Plugging the decomposition into the integral (\ref{heckeintegral}) yields
\begin{align}
\lambda_{\mu}(\mathbb{T}_p)&=\sum_{i=0}^{p-1} {\sf f}_p^{\circ}\left(\left(\begin{array}{cc} p & i\\ & 1\end{array} \right)\right)+{\sf f}_p^{\circ} \left(\left(\begin{array}{cc} 1 & \\ & p\end{array} \right)\right)
\nn \\
&= \underbrace{(\delta^{1/2}\mu)\left(\left(\begin{array}{cc} p & \\ & 1\end{array} \right)\right)+\cdots + (\delta^{1/2}\mu)\left(\left(\begin{array}{cc} p & \\ & 1\end{array} \right)\right)}_{p\, \, \text{terms}} + (\delta^{1/2}\mu)\left(\left(\begin{array}{cc} 1 & \\ & p\end{array} \right)\right)
\nn \\
&= p p^{-1/2} p^{-s_1}+p^{1/2} p^{-s_2} 
\nn \\
&= p^{1/2}(\alpha_p+\beta_p).
\label{Heckecalculation}
\end{align}
By a similar analysis one also shows that 
\begin{align}
\label{HeckeRp}
\lambda_\mu(\mathbb{R}_p)&=\int_{K_p \left(\begin{smallmatrix} p & \\ & p\end{smallmatrix} \right)K_p} {\sf f}_p^{\circ}(h)dh=\alpha_p\beta_p,\\
\lambda_\mu(\mathbb{R}_p^{-1})&=\int_{K_p \left(\begin{smallmatrix} p^{-1} & \\ & p^{-1}\end{smallmatrix} \right)K_p} {\sf f}_p^{\circ}(h)dh=\alpha_p^{-1}\beta_p^{-1}.
\end{align}

These results imply that the one-dimensional representation $\lambda_\mu$ of the Hecke algebra $\mathcal{H}_p^{\circ}$ acting on $V_p^{K_p}$ completely determines the unramified character $\mu$, and thereby the automorphic representation $\text{Ind}_{B(\mathbb{Q}_p)}^{GL(2,\mathbb{Q}_p)}\delta^{1/2}\mu$. It is quite remarkable that this infinite-dimensional  automorphic representation  can be encoded in the finite-dimensional representations of $\mathcal{H}_p^{\circ}$. In the next subsection we shall further investigate the consequences of this fact.

\begin{example}[Classical and $p$-adic Hecke operators for Eisenstein series on $SL(2)$]
\label{ex:HeckeConnection}
In this example, we come back to the mentioned relation~\eqref{eigenvaluerelation} between the Hecke eigenvalue of a non-holomorphic Eisenstein series $E(s,z)$ on $\UHP$ under the classical $T_p$ calculated in~\eqref{poincareheckeevals} and the action of $\mathbb{T}_p\in \mathcal{H}_p^\circ$ on the associated adelic Eisenstein series $E(\chi_s,g)\in\mathcal{A}(SL(2,\rats)\bs SL(2,\ads))$ defined in~\eqref{generalSL2AEisenstein}. 

To begin with, we need to relate the parameters $s_1$ and $s_2$ of the $GL(2,\rats_p)$ principal series~\eqref{unramifiedcharacterGL2} to the parameter $s$ occuring in $E(s,z)$ and $E(\chi_s,g)$ via~\eqref{SL2chis}. Elements of the Cartan torus in $SL(2,\rats_p)$ are of the form $a=\textrm{diag}(v,v^{-1})$, so that~\eqref{unramifiedcharacterGL2} yields
\begin{align}
\mu\left(\left(\begin{smallmatrix}v&\\&v^{-1}\end{smallmatrix}\right)\right) = |v|_p^{s_1-s_2}.
\end{align}
This has to be contrasted with (recall the general $\chi = \delta^{1/2}\mu$)
\begin{align}
(\delta^{-1/2}\chi_s)\left(\left(\begin{smallmatrix}v&\\&v^{-1}\end{smallmatrix}\right)\right) = |v|_p^{2s-1}
\end{align}
that follows from~\eqref{SL2chis}. For symmetry reason one therefore deduces
\begin{align}
s_1 = -s_2 = s-\frac12.
\end{align}
Plugging this into $\alpha_p$ and $\beta_p$ in~\eqref{abvals} one therefore finds from~\eqref{Heckecalculation} that
\begin{align}
\lambda_{\mu}(\mathbb{T}_p) = p^{1/2} \left( \alpha_p+ \beta_p\right) = p^{1/2} \left( p^{-s+1/2} + p^{s-1/2}\right) .
\end{align}
This is the eigenvalue of the adelic Eisenstein series under $\mathbb{T}_p$. From~\eqref{poincareheckeevals} one finds that for the classical Hecke operator $T_p$ acting on the classical $E(s,z)$ the eigenvalues is
\begin{align}
(T_p E)(s,z) = p^{s-1/2} \left(1+p^{1-2s}\right) E(s,z) = \left( p^{s-1/2} + p^{-s+1/2}\right) E(s,z).
\end{align}
This confirms the claimed relation~\eqref{eigenvaluerelation} that
\begin{align}
\mathbb{T}_p  \sim p^{1/2} \,T_p,
\end{align}
where we reiterate that the pre-factor is convention dependent.

\end{example}

\section{The Satake isomorphism}

We recall from section \ref{sec:simple-lie-alg}  that the Weyl group $\Weyl=\Weyl(\mathfrak{g})$ acts on the Cartan torus $A$, and consequently it also acts on the characters $\mu$ via 
\beq
w\mu(a)=\mu(w^{-1}aw), \qquad w\in \Weyl.
\eeq
Under this action the unramified principal series remains invariant
\beq
\text{Ind}_{B(\mathbb{Q}_p)}^{GL(2,\mathbb{Q}_p)}\delta^{1/2}w(\mu)\cong \text{Ind}_{B(\mathbb{Q}_p)}^{GL(2,\mathbb{Q}_p)}\delta^{1/2}\mu.
\eeq
This is what the functional relation~\eqref{funrel} for Eisenstein series expresses.

In terms of the parametrisation of $\mu$ by the complex numbers $(\alpha_p, \beta_p)$ the Weyl group $\Weyl=\mathbb{Z}/2\mathbb{Z}$ simply acts by $(\alpha_p, \beta_p)\mapsto (\beta_p, \alpha_p)$. Now notice that the characters $\lambda_{\mu}(\mathbb{T}_p)=p^{1/2}(\alpha_p+\beta_p)$, $\lambda_\mu(\mathbb{R}_p)=\alpha_p\beta_p$ and $\lambda_\mu(\mathbb{R}_p^{-1})=\alpha_p^{-1}\beta_p^{-1}$ are Weyl-invariant. Hence, at the level of the representations of the spherical Hecke algebra we have
\beq
\lambda_{w\mu} = \lambda_\mu, \qquad \forall \, w\in \Weyl.
\eeq
As a consequence, the image of the homomorphism $\mathcal{H}_p^{\circ}\to \mathbb{C}$ lies in the polynomial $\mathbb{C}$-ring of Weyl-invariants
\beq
\mathbb{C}[\alpha_p^{\pm 1}, \beta_p^{\pm 1}]^{\Weyl}\cong \mathbb{C}[\alpha_p+\beta_p, \alpha_p\beta_p, \alpha_p^{-1}\beta_p^{-1}].
\eeq
It is an important result of Satake \cite{Satake} that this homomorphism in fact yields an isomorphism between the spherical Hecke algebra and the ring of Weyl-invariant polynomials in $(\alpha_p, \beta_p)$:
\beq
\mathcal{H}_p^{\circ}\cong \mathbb{C}[\alpha_p+\beta_p, \alpha_p\beta_p, \alpha_p^{-1}\beta_p^{-1}].
\label{SatakeIsomorphism}
\eeq
See \cite{GrossSatake} for a nice survey of the \emphindex{Satake isomorphism} and its applications.

The key step in Satake's analysis was to introduce the \emph{Satake transform}
\beq
\mathcal{S}\, :\, \mathcal{H}_p^{\circ}(GL(2,\mathbb{Q}_p))\longrightarrow \mathcal{H}_p^{\circ}(A(\mathbb{Q}_p))
\eeq
from the spherical Hecke algebra of $GL(2,\mathbb{Q}_p)$ to the spherical Hecke algebra of the Cartan torus $A(\mathbb{Q}_p)$. The Satake transform is defined by
\beq
(\mathcal{S}\Phi)(a)=\delta^{-1/2}(a)\int_{N(\mathbb{Q}_p)} \Phi(na) dn, \qquad \Phi \in \mathcal{H}_p^{\circ}(GL(2,\mathbb{Q}_p)),
\eeq
where $N(\mathbb{Q}_p)$ is the unipotent radical of the Borel subgroup $B(\mathbb{Q}_p)\subset GL(2,\mathbb{Q}_p)$. Satake then proved that the image of $\mathcal{S}$ lies 
in $\mathcal{H}_p^{\circ}(A(\mathbb{Q}_p))^{\Weyl}$, the Weyl invariant elements of the spherical Hecke algebra of $A(\mathbb{Q}_p)$. To see the connection with our previous analysis, we consider again the formula for the eigenvalues $\lambda_\mu$:
\beq
\lambda_\mu(\Phi)=\int_{GL(2,\mathbb{Q}_p)}{\sf f}_p^{\circ}(h)\Phi(h) dh,
\eeq
which is (\ref{sphericalHeckeeigenvalue}) evaluated at the identity $g=1$. We shall now manipulate this expression in order to elucidate the role played by the Satake transform. To the best of our knowledge this calculation was first outlined by Langlands in \cite{LanglandsGodement}, but here we follow the more detailed exposition by Garrett  \cite{Garrett}. We begin by splitting the integral according to the Iwasawa decomposition $GL(2,\mathbb{Q}_p)=B(\mathbb{Q}_p)K_p$:
\beq
\int_{GL(2,\mathbb{Q}_p)}{\sf f}_p^{\circ}(h)\Phi(h) dh=\int_{B(\mathbb{Q}_p)} \int_{K_p} {\sf f}_p^{\circ}(b^{-1}k)\Phi(b^{-1}k) dbdk,
\eeq
where $dk$ and $db$ are  right-invariant Haar measures on $K_p$ and $B(\mathbb{Q}_p)$, respectively. Next, we make the change of variables $b\to b^{-1}$, which brings out a factor of $\delta^{-1}$ from the measure:
\beq
\int_{B(\mathbb{Q}_p)} \int_{K_p} {\sf f}_p^{\circ}(bk)\Phi(bk) \delta(b)^{-1}db dk.
\eeq
Using right $K_p$-invariance of ${\sf f}_p^{\circ}$  and $\Phi$ as well as $\int_{K_p}dk=1$ this further simplifies to
\beq
\int_{B(\mathbb{Q}_p)}{\sf f}_p^{\circ}(b) \Phi(b)\delta(b)^{-1} db=\int_{B(\mathbb{Q}_p)}\Phi(b) (\delta^{-1/2}\mu)(b) db, 
\eeq
where we used \eqref{sphericalvectorSL2}. To proceed we split the integral according to $B(\mathbb{Q}_p)=N(\mathbb{Q}_p)A(\mathbb{Q}_p)$ and use the fact that 
$\delta^{-1/2}\mu$ is trivial on $N(\mathbb{Q}_p)$ acting on the left:
\beq
\int_{A(\mathbb{Q}_p)}\int_{N(\mathbb{Q}_p)} \Phi(na) (\delta^{-1/2}\mu)(na)dn da=\int_{A(\mathbb{Q}_p)}\int_{N(\mathbb{Q}_p)} \Phi(na) (\delta^{-1/2})\mu(a)dn da.
\eeq
After reshuffling the integrand we finally arrive at the result
\beqa
\lambda_\mu(\Phi)&=&\int_{A(\mathbb{Q}_p)}\mu(a)\left[\delta^{-1/2}(a)\int_{N(\mathbb{Q}_p)} \Phi(na) dn\right] da
\nn \\
&=& \int_{A(\mathbb{Q}_p)}\mu(a) (\mathcal{S}\Phi)(a) da.
\eqa
This clearly shows that the Satake transform lies at the heart of the relation between the unramified automorphic representation $\text{Ind}_{B(\mathbb{Q}_p)}^{GL(2,\mathbb{Q}_p)}\delta^{1/2}\mu$ and the one-dimensional representation $\lambda_\mu$ of the spherical Hecke algebra $\mathcal{H}_p^{\circ}$, the essence of which is the Satake isomorphism (\ref{SatakeIsomorphism}).

\section{The \texorpdfstring{$L$}{L}-group and generalisation to \texorpdfstring{$GL(n)$}{GL(n)}}
\label{sec_Lgroup}

It is illuminating to assemble the  parameters $(\alpha_p, \beta_p)$ in a matrix
\beq
A_{\pi_p}=\left(\begin{array}{cc} \alpha_p & \\ & \beta_p \end{array} \right).
\label{SatakeGL2}
\eeq
This matrix belongs to $GL(2,\mathbb{C})$ and since conjugation $A_{\pi_p} \mapsto wA_{\pi_p} w^{-1}$ by an element $w\in \Weyl$ will not alter the result (\ref{SatakeIsomorphism}) 
we find that the representation $\pi_p$ determines a (semi-simple) \emphindex[conjugacy class!semi-simple]{conjugacy class} $[A_{\pi_p}]\subset GL(2,\mathbb{C})$. This conjugacy class is called the \emphindex{Satake parameter} of the local representation $\pi_p$.

The conclusion of the discussion in this and the  previous sections is that \emphindex[automorphic representation!unramified and $L$-group]{unramified automorphic representations $\pi_p$ of $GL(2,\mathbb{Q}_p)$ are in bijection with semi-simple conjugacy classes $[A_{\pi_p}]\subset GL(2,\mathbb{C})$}. The appearance of $GL(2,\mathbb{C})$ in the context of local representations of $GL(2,\mathbb{Q}_p)$ may seem surprising, but is in fact a simple instance of a more general phenomenon envisioned by Langlands \cite{LanglandsProb}. Langlands suggested that to each reductive algebraic group $G$ over a number field $\field$ there exists an associated  complex group ${}^{L}G(\mathbb{C})$, called the \emphindex[L-group@$L$-group]{$L$-group}, or \emphindex{Langlands dual group}. We have already encountered the group ${}^LG$ briefly in section~\ref{sec:CSLD} in our discussion of the Casselman--Shalika formula but we will now put this group into a more general context. 

A precise definition of ${}^LG$ can be found in~\cite{BorelCorvallis}; we here only recall the salient features. For simple groups $G$ the root system of ${}^LG$ is obtained from that of $G$ by interchanging the short and long roots. In other words, the co-weight lattice $\Lambda^{\vee}$ of the Lie algebra $\mathfrak{g}=\text{Lie}\, G$ is identified with the weight lattice ${}^L\Lambda$ of the dual Lie algebra ${}^{L}\mathfrak{g}=\text{Lie}\, {}^LG$. This is captured by the isomorphism
\beq
\text{Hom}({}^LA, U(1))\cong \text{Hom}(U(1), A),
\eeq
between the lattice of characters on ${}^LA$ and the lattice of co-characters on $A$. For example, in the case of $G=GL(n, \mathbb{Q}_p)$ the $L$-group is $GL(n, \mathbb{C})$, and for $G=SL(n, \mathbb{Q}_p)$ we have ${}^LG=PGL(n, \mathbb{C})$. The duality is even more drastic in the case when $G=Sp(n)$ we have ${}^LG=SO(2n+1)$. See also \cite{BumpHecke,CasselmanLGroup} for details. 

\begin{remark}
The group ${}^LG$ we have introduced here is sometimes called the \emphindex[L-group@$L$-group!connected]{connected $L$-group} in order to distinguish it from the $L$-group in the more general context of field extensions. If one considers a finite field extension $\fieldext$ of $\field=\rats_p$ then the $L$-group ${}^LG$ is defined with the inclusion of (finite) Galois group $\textrm{Gal}(\fieldext/\field)$ of the field extension. This more general viewpoint is relevant for the global Langlands conjectures and will be discussed in section~\ref{sec:Langlands}.\label{rmk:fieldext}
\end{remark}

The Satake parameter $A_{\pi_p}$ associated with the automorphic representation $\pi_p$ should thus be viewed as an element of the Cartan torus ${}^LA(\mathbb{C})\subset GL(2,\mathbb{C})$ dual to the original Cartan torus $A(\mathbb{Q}_p)$. In fact, this holds more generally for any (split) reductive algebraic group $G$. From this perspective, one gets the following  reformulation of the Satake isomorphism (adapted from \cite{CasselmanLGroup}): 

\begin{theorem}[reformulated Satake isomorphism]
There is a natural bijection between the Weyl-invariant homomorphism $\mathcal{H}_p^{\circ}(G) \to \mathbb{C}$ and semi-simple ${}^LG(\mathbb{C})$-conjugacy classes in the dual torus ${}^LA(\mathbb{C})$. \index{Satake isomorphism}
\end{theorem}

\begin{remark}
The Satake parameter $A_{\pi_p}\in {}^LA(\cx)$ already appeared in section~\ref{sec:CSLD} where it was denoted by $a_\lambda$ where $\lambda$ parametrises an element of the principal series representation of $G(\rats_p)$ which is here denoted abstractly by $\pi_p$. 
\end{remark}

Let us briefly discuss some details on the generalisation of our analysis to $G=GL(n, \mathbb{Q}_p)$. We take $\pi_p$ to be the unramified principal series with module
\beq
V_p=\text{Ind}_{B(\mathbb{Q}_p)}^{GL(n, \mathbb{Q}_p)} \delta^{1/2}\mu, 
\eeq
where the inducing character is a straightforward generalisation of (\ref{unramifiedcharacterGL2}):
\beq
\mu(a)=\mu\left(\begin{array}{ccc} v_1  &  &  \\ & \ddots & \\  & & v_n 
\end{array}\right)=\prod_{i=1}^{n}|v_i|_p^{s_i}.
\eeq
As before, this representation is determined  by the $n$ complex numbers:
\beq
\alpha_i:= p^{-s_i}, \qquad i=1,\dots, n.
\label{GLnparameters}
\eeq
(Note that the $\alpha_i$ here are for fixed prime $p$ that we do not indicate explicitly unlike in~\eqref{SatakeGL2}.)
Associated with the representation $\pi_p$ we then have the Satake parameter \index{Satake parameter!for $GL(n)$}
\beq
A_{\pi_p}=\left(\begin{array}{ccc} \alpha_1  &  &  \\ & \ddots & \\  & & \alpha_n 
\end{array}\right)\in {}^LA(\mathbb{C})\subset GL(n, \mathbb{C})={}^L GL(n,\rats_p),
\label{SatakeParameterGLn}
\eeq
on which the Weyl group $\Weyl$ acts by permuting the $\alpha_i$'s. The generators of the spherical Hecke algebra $\mathcal{H}_p^{\circ}(G)$ act on elements $\varphi\in V_p$ by~\eqref{rightregularHeckeaction}, \emph{viz}.
\beqa
(\Phi_i \varphi)(g)&=&\int_{GL(n, \mathbb{Q}_p)} \varphi(gh)\text{char}_{K_p \tau_i K_p}(h) dh
\nn \\
&=& \int_{K_p\tau_i K_p} \varphi(gh)dh,
\label{generalHeckeGLn}
\eqa
where we defined~\cite{BumpHecke}
\beq
\tau_i=\left(\begin{array}{cc} p \idm_i & \\ & \idm_{n-i} \end{array}\right),
\eeq
with $\idm_r$ the $r\times r$ identity matrix. We use the convention that for $i=n$ the double coset is $K_p (p \idm_n) K_p$. Thus, in the special case of $n=2$ the definition (\ref{generalHeckeGLn}) reduces to the generators in section \ref{app_adelichecke}, i.e. $\Phi_1=\mathbb{T}_p$ and $\Phi_2=\mathbb{R}_p$. Tamagawa has shown~\cite{Tamagawa} that the operators $\Phi_1,\ldots,\Phi_n$ together with $\Phi_n^{-1}$ (which is the only invertible $\Phi_i$) generate the spherical Hecke algebra $\mathcal{H}_p^\circ$ of $GL(n,\rats_p)$.

As before, the one-dimensional space $V_p^{K_p}=\mathbb{C} \cdot {\sf f}_p^{\circ}$ of $K_p$-invariant vectors in $V_p$ furnishes a representation of the spherical Hecke algebra, such that for any $\Phi\in \mathcal{H}_p^{\circ}$ and any $v^{\circ}\in V_p^{K_p}$ one has
\beq
\pi(\Phi) v^{\circ}=\lambda_\mu(\Phi)v^{\circ},
\eeq
where $\lambda_\mu : \mathcal{H}_p^{\circ}\to \mathbb{C}^\times$ is a (quasi-)character. To evaluate the eigenvalue $\lambda_\mu(\Phi)$ on all the generators $\Phi_i$ we must decompose the double cosets in (\ref{generalHeckeGLn}). The result can be written as follows using the finite Cartan decomposition (see for instance \cite{BumpHecke} for a nice and explicit proof)
\beq
K_p\tau_i K_p=\bigcup_{j} \beta_{i,j} K_p,
\eeq
where 
the matrices $\beta_{,_j}$ are all integral and upper-triangular with diagonal entries are of the form $p^{\eta}$, where $\eta\in \{1, \dots, n\}$ and $j$ ranges over some finite set. These generalise the matrices on the first line of \eqref{Heckecalculation} and similarly to that calculation we must evaluate the spherical vector ${\sf f}_p^{\circ}$ on all $\beta_{i,j}$. Bump shows that this takes the form \cite{BumpHecke}
\beq
{\sf f}_p^{\circ}(\beta_{i,j})=(\delta^{1/2}\mu)(\beta_{i,j})=p^{-\tfrac{i(n+1)}{2}} \prod_{\ell=1}^{i} p^{\eta_\ell} \alpha_{\eta_\ell},
\eeq
where $\eta_\ell\in \{1, \dots, n\}$ are determined by $j$ and ordered such that $\eta_1<\eta_2<\cdots <\eta_i$. The $\alpha_{\eta_\ell}$ are the complex parameters (\ref{GLnparameters}) that determine the representation $\pi_p$. There are furthermore a total number of 
\beq
p^{i(n-i-\tfrac{1}{2})-\sum_{\ell=1}^{i}\eta_\ell}
\eeq
$\beta_{i,j}$ for each $i\in \{1, \dots, n\}$. Combining everything we find that the eigenvalue of the Hecke operator $\Phi_i$ is given by 
\begin{align}
\lambda_\mu(\Phi_i)&= \int_{K_p\tau_i K_p} {\sf f}_p^{\circ}(h)dh
=\sum_{\beta_i\in \Lambda_i} {\sf f}_p^{\circ}(\beta_i)
\nn \\
&= 
\sum_{\eta_1 <\cdots <\eta_i} p^{i(n-i-\tfrac{1}{2})-\sum_{j=1}^{i}\eta_j} p^{-\tfrac{i(n+1)}{2}} \prod_{\ell=1}^{i} p^{\eta_\ell} \alpha_{\eta_\ell}
\nn \\
&= p^{i(n-i)/2} \sum_{\eta_1 <\cdots <\eta_i} \alpha_{\eta_1}\cdots \alpha_{\eta_n}
\nn \\
&= p^{i(n-i)/2} e_i(\alpha_1, \dots, \alpha_n),
\label{CaclculationHeckeGLn}
\end{align}
where $e_i(\alpha_1, \dots, \alpha_n)$ is the $i$th elementary symmetric polynomial in $n$ variables. In fact, the Satake isomorphism can be written in terms of these
elementary symmetric polynomials
\beq
\mathcal{H}_p^{\circ}(GL(n, \mathbb{Q}_p))\cong \mathbb{C}[e_1(\alpha_1, \dots, \alpha_n), \dots , e_n(\alpha_1, \dots, \alpha_n), e_n(\alpha_1, \dots , \alpha_n)^{-1}],
\eeq
corresponding to the values on the generators $\Phi_1,\ldots,\Phi_n$ and $\Phi_n^{-1}$ of the spherical Hecke algebra of $GL(n,\rats_p)$. 
Indeed, for $n=2$ we have 
\beq
e_1(\alpha_1, \alpha_2)=\alpha_1+\alpha_2, \qquad e_2(\alpha_1, \alpha_2)=\alpha_1\alpha_2,
\eeq
thus recovering (\ref{SatakeIsomorphism}).

Let us end this section with a comment on how these results fit into the general theory of automorphic forms. Recall from definition \ref{defauto} that an automorphic form $\varphi$ on the adelic quotient $G(\mathbb{Q})\backslash G(\mathbb{A})$ is required to be $\mathcal{Z}(\mathfrak{g})$-finite, i.e. that $\varphi$ is an eigenfunction of the entire ring of invariant differential operators in the center of $U(\mathfrak{g})$. This can be viewed as a statement about the behavior of $\varphi$ under the action of differential operators in the \emph{real} group $G_\infty=G(\mathbb{R})$. For the case of automorphic forms attached to unramified automorphic representations $\pi = \pi_\infty \otimes \bigotimes_{p<\infty} \pi_p$ the spherical Hecke algebra provides the non-archimedean analogue of this:  for each finite place $p$, $\varphi$ is an eigenfunction of the ring of Hecke operators generated by $\Phi_i\in \mathcal{H}_p^{\circ}$. These statements combine together in the global Hecke algebra as mentioned in remark~\ref{rmk:globalHecke}.

\section{The Casselman--Shalika formula revisited}

\index{Casselman--Shalika formula!and Hecke algebra}
There is  a close relation between the discussion above and the Casselman--Shalika formula for the $p$-adic spherical Whittaker function $W_\psi^{\circ}$. Spherical Whittaker functions were the central objects in chapter~\ref{ch:Whittaker-Eisenstein} and a glimpse of the relation between them and representation theory was already visible in section~\ref{sec:CSLD}. Here, we recall and extend some of the notions in a more general context. For an \index{Whittaker function!unramifed} unramified character $\psi : N(\mathbb{Z}_p)\backslash N(\mathbb{Q}_p) \to U(1)$ we have an embedding 
\beq
W_\psi : \text{Ind}_{B(\mathbb{Q}_p)}^{GL(n, \mathbb{Q}_p)} \delta^{1/2}\mu \longrightarrow \text{Ind}_{N(\mathbb{Q}_p)}^{GL(n, \mathbb{Q}_p)} \psi
\label{WhittakerEmbeddingGLn}
\eeq
of the unramified principal series into the space of functions $W : GL(n, \mathbb{Q}_p)\to \mathbb{C}$ satisfying 
\beq
W_\psi(\delta^{1/2}\mu,ng)=\psi(n)W_\psi(\delta^{1/2}\mu, g), \qquad \forall n\in N(\mathbb{Q}_p),
\eeq
where, as in chapter~\ref{ch:Whittaker-Eisenstein}, the first argument indicates the dependence on the inducing character $\mu$ in the unramified principal series that was written there in terms of  $\chi=\delta^{1/2}\mu$. The image of the space $V_p^{K_p}$ of $K_p$-fixed vectors in $V_p$ is a one-dimensional space of \emphindex[Whittaker function!spherical]{spherical Whittaker functions}. In particular, for the generator ${\sf f}_p^{\circ}\in V_p^{K_p}$ we obtain a canonical spherical Whittaker function via the explicit \emphindex{Jacquet integral} (see chapter~\ref{ch:Whittaker-Eisenstein} for details)
\beq
W_\psi^{\circ}(\delta^{1/2}\mu,g)=\int_{N(\mathbb{Q}_p)} {\sf f}_p^{\circ}(w_0 ng)\overline{\psi(n)} dn,
\eeq
where we used $ {\sf f}_p^{\circ}=\delta^{1/2}\mu$. This satisfies 
\beq
W_\psi^{\circ}(\delta^{1/2}\mu, nak)=\psi(n)W_\psi(\delta^{1/2}\mu, a), 
\label{sphericalWhittakerrelation}
\eeq
and so is completely determined by its restriction to the Cartan torus $A(\mathbb{Q}_p)$. For $GL(n, \mathbb{Q}_p)$ the vanishing properties of $W_\psi^{\circ}$ analysed in section~\ref{VP} can be simplified as follows. Parametrising $a$ according to 
\beq
a=\varpi^J:=\left(\begin{array}{ccc} p^{j_1}  &  &  \\ & \ddots & \\  & & p^{j_n} 
\end{array}\right)\in A(\mathbb{Q}_p)/A(\mathbb{Z}_p),
\label{varpiJ}
\eeq
with $J=(j_1, \dots, j_n)\in \mathbb{Z}^n$, one finds that (see, e.g., \cite{Cogdell})
\beq
W_\psi^{\circ}(\delta^{1/2}\mu, a)=0. \qquad \text{unless} \quad j_1\geq j_2\geq \cdots \geq j_n.
\label{vanishingWhittakerGLn}
\eeq
The map (\ref{WhittakerEmbeddingGLn}) commutes with the Hecke action and therefore the spherical Whittaker function is an eigenfunction of all the Hecke operators with 
the same eigenvalue (\ref{CaclculationHeckeGLn}) as before:
\beq
\Phi_i W_\psi^{\circ}(\delta^{1/2}\mu, a)=\lambda_\mu(\Phi_i)W_\psi^{\circ}(\delta^{1/2}\mu,a). 
\label{HeckeWhittakerEigenvalue}
\eeq
This fact can be used to derive a recursive formula for the value $W_\psi^{\circ}(\delta^{1/2}\mu, a)$ as we will now show. This will give the connection with the Casselman--Shalika formula that we are after. 

The main difference with  the calculation (\ref{CaclculationHeckeGLn}) is of course that \emph{a priori} we do not know the explicit value of $W_\psi^{\circ}$ on $A(\mathbb{Q}_p)$, in contrast to the case of the original spherical vector ${\sf f}_p^{\circ}$ where we had the formula (\ref{sphericalvectorSL2}) at hand. The key is that we should parametrise the decomposition of the cosets $K_p \tau_i K_p$ in such a way that we can make use of the defining relation (\ref{sphericalWhittakerrelation}). Such a parametrisation was given by Shintani \cite{Shintani}; here we follow the treatment by Cogdell \cite{Cogdell}, which reads
\beq
K_p\tau_i K_p= \bigcup_{\epsilon\in I_i} \bigcup_{n\in N_\epsilon} n\varpi^{\epsilon}K_p,
\eeq
where the set $I_i$ is defined as 
\beq
I_i=\{\epsilon=(\epsilon_1, \dots, \epsilon_n)\in \mathbb{Z}^n \, |\, \epsilon_j\in\{0,1\}, \sum_{j=1}^{n} \epsilon_j =i\},
\eeq
and 
\beq
N_\epsilon=N(\mathbb{Z}_p)/(N(\mathbb{Z}_p)\cap \varpi^{\epsilon} K_p \varpi^{-\epsilon}).
\label{Neps}
\eeq
Using this result we can compute the left hand side of (\ref{HeckeWhittakerEigenvalue}) explicitly:
\beqa
\int_{K_p\tau_i K_p} W_\psi^{\circ}(\delta^{1/2}\mu,\varpi^{J} h)dh&=&\sum_{\epsilon\in I_i} \sum_{n\in N_\epsilon} W_\psi^{\circ}(\delta^{1/2}\mu, \varpi^{J} n \varpi^{\epsilon})
\nn \\
&=& \sum_{\epsilon\in I_i} \sum_{n\in N_\epsilon} W_\psi^{\circ}(\delta^{1/2}\mu, \varpi^{J} n\varpi^{-J}\varpi^{J} \varpi^{\epsilon})
\nn \\
&=& \sum_{\epsilon\in I_i} \sum_{n\in N_\epsilon} \psi(\varpi^{J} n\varpi^{-J}) W_\psi^{\circ}(\delta^{1/2}\mu, \varpi^{J} \varpi^{\epsilon}),
\eqa
where we used that $\varpi^{J} n\varpi^{-J}\in N(\mathbb{Q}_p)$ combined with (\ref{sphericalWhittakerrelation}). In fact, because of the constraint (\ref{vanishingWhittakerGLn}), which requires $j_1\geq \cdots \geq j_n$, we have that $\varpi^{J} n\varpi^{-J}\in N(\mathbb{Z}_p)$ and consequently $\psi(\varpi^{J} n\varpi^{-J})=1$. The summand is therefore independent of $n$ and the sum yields only a factor corresponding to the size of the coset space (\ref{Neps}). Cogdell shows that \cite{Cogdell}
\beq
|N_\eps|=p^{i(n-i)/2} \delta^{-1/2}(\varpi^{\epsilon}),
\eeq
so we obtain for all $i$~\cite[Prop.~7.3]{Cogdell}
\beq
\lambda_\mu(\Phi_i) W_\psi^{\circ}(\delta^{1/2}\mu, \varpi^{J})=\sum_{\epsilon\in I_i}p^{i(n-i)/2} \delta^{-1/2}(\varpi^{\epsilon})W_\psi^{\circ}(\delta^{1/2}\mu, \varpi^{J+\epsilon}).
\label{WhittakerRecursion}
\eeq
This is a recursive formula for the spherical Whittaker function $W_\psi^{\circ}(\delta^{1/2}\mu, \varpi^J)$! We recall that all the Hecke eigenvalues $\lambda_\mu(\Phi_i)$ are known from~\eqref{CaclculationHeckeGLn}.

\begin{example}[Unramified Whittaker functions for $GL(2,\rats_p)$]
\label{ex:GL2rec}
Let us determine some unramified spherical Whittaker functions for $GL(2,\rats_p)$ using the recursion relation~\eqref{WhittakerRecursion}, starting from $J=(0,0)$. The recursion relation then reads for the two values $i=1,2$
\begin{align}
p^{1/2} (\alpha_p+\beta_p) W^{(0,0)} &= p^{1/2} \delta^{-1/2}\left(\varpi^{(0,1)}\right) W^{(0,1)} + p^{1/2} \delta^{-1/2}\left(\varpi^{(1,0)}\right) W^{(1,0)},\\
\alpha_p \beta_p W^{(0,0)} &= \delta^{-1/2}\left(\varpi^{(1,1)}\right) W^{(1,1)},
\end{align}
where~\eqref{Heckecalculation} and~\eqref{HeckeRp} were used and we have introduced the short-hand notations
\begin{align}
W^{(j_1,j_2)} \equiv W^\circ_\psi (\delta^{1/2}\mu, \varpi^J)
\quad\textrm{and}\quad
\varpi^{(j_1,j_2)} \equiv \varpi^J.
\end{align}
Since $W^{(0,1)}=0$ according to~\eqref{vanishingWhittakerGLn} we can solve for $W^{(1,0)}$ and $W^{(1,1)}$ in terms of $W^{(0,0)}$ to obtain
\begin{align}
W^{(1,0)} = W^{(0,0)} \delta^{1/2}\left(\varpi^{(1,0)}\right)   (\alpha_p+\beta_p) ,\quad
W^{(1,1)} =  W^{(0,0)} \delta^{1/2}\left(\varpi^{(1,1)}\right) \alpha_p \beta_p.
\end{align}
We note that 
\begin{align}
\alpha_p+\beta_p = \textrm{Tr}\begin{pmatrix}\alpha_p&\\&\beta_p\end{pmatrix}=\textrm{Tr}_{(1,0)}(A_{\pi_p}) 
\quad\textrm{and}\quad
\alpha_p\beta_p =  \textrm{Tr}\begin{pmatrix}\alpha_p\beta_p\end{pmatrix}=
\textrm{Tr}_{(1,1)}(A_{\pi_p}) 
\end{align}
are the characters of the Satake parameter $A_{\pi_p}$ in the two- and one-dimensional representations of $GL(2,\cx)={}^LGL(2,\rats_p)$, respectively, that are labelled here by their Young tableaux indexed by $J=(1,0)$ and $J=(1,1)$. The translation from non-increasing tuples $(j_1,\ldots, j_n)$ to a \emphindex{Young tableau} is such that the $i$th row has $j_i$ boxes. Therefore, we have
\begin{align}
J=(1,0) \longleftrightarrow \yng(1)\quad\quad\textrm{and}\quad\quad
J=(1,1) \longleftrightarrow \yng(1,1)\,\,,\nn
\end{align}
such that $(1,0)$ corresponds to the fundamental two-dimensional representation and $(1,1)$ to the one-dimensional $GL(2,\cx)$ representation of weight $1$ (tensor density). The relation between spherical Whittaker functions and characters is no coincidence as we explain in the text.
\end{example}

\index{highest weight representation!of $GL(n,\cx)$}

The key to solving the recursion relation~\eqref{WhittakerRecursion} for $GL(n,\rats_p)$ is is to note that the set of integers $J=(j_1, \dots, j_n)$,
subject to the condition $j_1\geq j_2\geq \cdots \geq j_n$, is well-known to parametrise the highest weights of irreducible representations $V_J$ of $GL(n, \mathbb{C})$, which, we recall, is the 
$L$-group ${}^LG$ of $GL(n, \mathbb{Q}_p)$. But the analogy goes even further than that. Let $\chi_J=\text{Tr}_{V_J}$ be the character of the representation $V_J$. This is a \emphindex{class function}, meaning that it is invariant under conjugation
\beq 
\text{Tr}_{V_J}(g)=\text{Tr}_{V_J}(hgh^{-1}), \qquad g, h\in GL(n, \mathbb{C}),
\eeq
and so only depends on the conjugacy class of $V_J$. If we take $V_J$ to be the fundamental $n$-dimensional representation of $GL(n, \mathbb{C})$  then we already have a conjugacy class at hand, namely the Satake parameter $A_{\pi_p}\in {}^LA(\mathbb{C})$ (\ref{SatakeParameterGLn}) of $\pi_p$. From this perspective $J$ is a dominant weight in the weight lattice $\Lambda^{\vee}$ of ${}^L\mathfrak{g}$, which is the co-weight lattice of $\mathfrak{g}$. One can then solve the recursion (\ref{WhittakerRecursion}) in terms of the characters $\chi_J$ with the result \cite{Cogdell}
\begin{equation}
    \label{otherwise}
    W_\psi^{\circ}(\delta^{1/2}\mu, \varpi^{J}) =
    \begin{cases}
        \text{const} \times \delta^{1/2}(\varpi^J)\chi_J(A_{\pi_p}) \quad & \text{if } J\in \Lambda^{\vee} \text{ is dominant} \\   
        0 & \text{otherwise.} 
    \end{cases}
\end{equation}
We note that the recursion relation only determines the spherical Whittaker function up to a constant. At first sight this looks very different from the Casselman--Shalika formula (\ref{CSus}) we derived in chapter~\ref{ch:Whittaker-Eisenstein}. To see that they indeed coincide we shall rewrite the formula given there in a way similar to what was done in section~\ref{sec:CSLD}. Setting $a=\varpi^J$ in (\ref{CSus})  and doing some reshuffling we arrive at
\begin{align}
\frac1{\zeta(\lambda)} \sum_{w\in \Weyl} \eps(w\lambda) |a^{w\lambda+\rho}|&= \frac{1}{\zeta(\delta^{1/2}\mu)}a^\rho \sum_{w\in \Weyl}w\left[\frac{a^{\lambda}}{\prod_{\alpha>0} (1-p^{\left<\lambda|\alpha\right>})}\right]
\nn \\
&= \frac{1}{\zeta(\delta^{1/2}\mu)} \delta^{1/2}(\varpi^J)\sum_{w\in \Weyl} w\left[ \frac{\mu(\varpi^J)}{\prod_{\alpha>0} (1-\mu(\varpi^{-\alpha}))}\right],
\label{rewriteCS}
\end{align}
where we rewrote the arguments as follows
\begin{subequations}
\begin{align}
a^{\rho}&= e^{\left<\rho|H(\varpi^J)\right>}=\delta^{1/2}(\varpi^J),\\
a^{\lambda}&=e^{\left<\lambda |H(\varpi^J)\right>}=p^{-\left<\lambda|J\right>}=\mu(\varpi^J), \\
p^{\left<\lambda|\alpha\right>}&=\mu(\varpi^{-\alpha}),
\end{align}
\end{subequations}
To interpret the new form (\ref{rewriteCS}) of the Casselman--Shalika formula we recall that the weight lattice of the $L$-group ${}^LG(\mathbb{C})$ is $\Lambda^{\vee}$, the co-weight lattice of $G$. We now identify this with the character lattice $X^{*}({}^LA)$ according to
\beq
\Lambda^{\vee}\cong X^{*}({}^LA)\cong \text{Hom}({}^LA, U(1))\cong \mathbb{Z}^n.
\eeq
Under this identification a weight $J=(j_1, \dots, j_n)\in \Lambda^\vee$ can be interpreted as a character $J\, :\, {}^LA(\mathbb{C})\to U(1)$. We can in particular evaluate this character on the 
Satake parameter $A_{\pi_p}\in {}^LA(\mathbb{C})$ with the result
\beq
A_{\pi_p}^{J} =J(A_{\pi_p})=J\left(\begin{array}{ccc} \alpha_1  &  &  \\ & \ddots & \\  & & \alpha_n 
\end{array}\right)=\prod_{i=1}^n \alpha_i^{j_i},
\label{Jvalue}
\eeq
which further implies the equality
\beq
A_{\pi_p}^{J} =\mu(\varpi^J).
\eeq
\begin{remark}
The standard notation being used here might be the source for confusions: In general we denote the value of the character $J$ on $a\in {}^LA$ by $a^J$ or $J(a)$ as in (\ref{Jvalue}); however this should \emph{not} be confused with the \emph{matrix} $\varpi^J$, which is defined in (\ref{varpiJ}). We trust that this will not cause any trouble since it should be clear from the context which definition is referred to.
\end{remark}

Next we compare~\eqref{rewriteCS} with the Weyl character formula for a representation $V_J$ of a Lie group $G$ with highest weight $J$. According to~\eqref{charfn}, the character $\chi_J$ evaluated at $z\in A$ is explicitly given by
\beq
\chi_J(z)=\sum_{w\in \Weyl} w\left[\frac{z^J}{\prod_{\alpha >0} (1-z^{-\alpha})}\right], \qquad z\in A.
\eeq
We can therefore rewrite (\ref{rewriteCS}) as 
\beq
\frac{1}{\zeta(\delta^{1/2}\mu)} \delta^{1/2}(\varpi^J)\sum_{w\in \Weyl} w\left[ \frac{\mu(\varpi^J)}{\prod_{\alpha>0} (1-\mu(\varpi^{-\alpha}))}\right]=\frac{1}{\zeta(\delta^{1/2}\mu)} \delta^{1/2}(\varpi^J)\chi_J(A_{\pi_p}).
\eeq
Comparing this with (\ref{otherwise}) we indeed find perfect agreement, provided that we fix the overall constant there to be $\zeta(\mu)^{-1}$. We conclude that the Casselman--Shalika formula for the spherical Whittaker function $W_\psi^{\circ}\in \big(\text{Ind}_{N(\mathbb{Q}_p)}^{G(\mathbb{Q}_p)}\psi\big)^{K_p}$ can be written in terms of the Weyl character formula for an irreducible representation $V_J$ of the Langlands dual group ${}^LG(\mathbb{C})$:
\beq
W_\psi^{\circ}(\delta^{1/2}\mu, A_{\pi_p})=\frac{1}{\zeta(\delta^{1/2}\mu)} \delta^{1/2}(\varpi^J) \chi_J(A_{\pi_p}).
\label{CScharacter}
\eeq

For $GL(n,\cx)$ the characters of $V_J$ is well-known to be given by symmetric polynomials that can be expressed in the basis of Schur polynomials. Examples for $GL(2)$ can be found in~\ref{ex:GL2rec}.

\section{Automorphic \texorpdfstring{$L$}{L}-functions}
\label{sec_autoLfunction}

Equipped with the adelic Hecke technology of the previous sections we shall now revisit the discussion of Dirichlet series of section~\ref{app_adelichecke} in the more general context of $GL(n, \mathbb{A})$. 

Suppose first that $f$ is a Maass form on the upper-half plane $\UHP$ which is an eigenfunction of the classical Hecke operator $T_p$ with eigenvalue $a_p$. For instance, $f$ could be a non-holomorphic Eisenstein series. This lifts to an automorphic form $\varphi_f\in \mathcal{A}(SL(2,\mathbb{Q})\backslash SL(2, \mathbb{A}))$ which is an eigenfunction of $\mathbb{T}_p$ with eigenvalue $\lambda_\mu(\mathbb{T}_p)=p^{1/2} (\alpha_p+\beta_p)$ as we found in (\ref{Heckecalculation}). According to (\ref{eigenvaluerelation}) the relation between the eigenvalues is thus
\beq
a_p=\alpha_p+\beta_p.
\label{alphabeta}
\eeq
This implies that we can rewrite the local factor in the Dirichlet series (\ref{LNH}) as follows 
\beq
(1-a_pp^{-s}+p^{-2s})^{-1} = \left[(1-\alpha_p p^{-s})(1-\beta_p p^{-s})\right]^{-1}=\text{det}\left(\idm-A_{\pi_p}p^{-s}\right)^{-1},
\eeq
where $A_{\pi_p}$ is the semi-simple Satake parameter (\ref{SatakeGL2}) in the fundamental matrix representation. 

The relation (\ref{alphabeta}) has a natural generalisation to higher rank groups. Suppose $\varphi\in \mathcal{A}(G(\mathbb{Q})\backslash G(\mathbb{A}))^{K_\mathbb{A}}$, i.e., $\varphi$ is a spherical automorphic form, is attached to an unramified automorphic representation $\pi$. Suppose also that $\varphi$ is an eigenfunction of the spherical Hecke algebras $\mathcal{H}_p^{\circ}=\mathcal{H}(\mathbb{Q}_p)^{K_p}$. This implies that for $\Phi\in \mathcal{H}_p^{\circ}$ we have $\pi(\Phi)\varphi = \lambda_\pi(\Phi) \varphi$. In this situation there exists a unique Satake class $[A_{\pi_p}]\subset {}^LG(\mathbb{C})$ such that 
\beq
\lambda_\pi(\Phi)=p^{\sharp}\textrm{Tr}_\pi(A_{\pi_p}),
\eeq
where the prefactor is some power of the prime $p$. In particular, for $G=GL(n)$ we see from \eqref{CaclculationHeckeGLn} that 
\beq
\lambda_\pi(\Phi_1)=\lambda_\mu(\Phi_1)=p^{(n-1)/2} (\alpha_1+\cdots +\alpha_n)=p^{(n-1)/2} \textrm{Tr}_\pi(A_{\pi_p}),
\eeq
where the semi-simple conjugacy class $A_{\pi_p}$ is given in (\ref{SatakeParameterGLn}).

We can now generalise the construction of the Dirichlet series to $GL(n, \mathbb{A})$. To this end let $\pi=\bigotimes_{p\leq \infty}\pi_p$ be the unramified principal series $\text{Ind}_{B(\mathbb{A})}^{GL(n, \mathbb{A})}\delta^{1/2}\mu$ and $A_{\pi_p}$ be the corresponding Satake parameter associated with each local factor $\pi_p$. To this data we attach the following \emph{local $L$-factor}:
\beq
L_p(\pi_p, s) =\text{det}\left(\idm-A_{\pi_p} p^{-s}\right)^{-1},
\eeq
and we define the \emph{standard $L$-function} as
\beq
L(\pi, s)=\prod_{p<\infty} L_p(\pi_p, s).
\eeq
Langlands has proven \cite{LanglandsFE} that this can be completed by adding a certain factor for the prime at infinity 
\beq
L^\compl(\pi, s)=L_\infty(\pi_\infty, s)\prod_{p<\infty} L_p(\pi_p, s),
\eeq
which  has an analytic continuation to a meromorphic function in the entire complex $s$-plane, and satisfying a functional equation. This is a vast generalisation of the 
completed Riemann zeta-function $\xi(s)=\xi_\infty(s) \prod_{p<\infty} (1-p^{-s})^{-1}$, where the prime at infinity corresponds to the Gamma-factor $\xi_\infty(s)=\pi^{-s/2}\Gamma(s/2)$. For Maass wave forms on $\UHP$ the factors at infinity were given in~\eqref{compL}.

But Langlands suggested to generalise this even further. Suppose $G$ is a reductive algebraic group over $\mathbb{Q}_p$ and $\pi_p$ is an unramified automorphic representation of $G(\mathbb{Q}_p)$. Let $A_{\pi_p}$ be the associated Satake parameter, giving a semi-simple conjugacy class $[A_{\pi_p}]\subset {}^LG(\mathbb{C})$. Let further 
\beq
\rho \, :\, {}^LG(\mathbb{C})\, \longrightarrow \, GL(n, \mathbb{C})=\textrm{Aut}(\cx^n)
\eeq
be an $n$-dimensional representation of the $L$-group. Note that the representation does \emph{not} depend on the prime $p$. In the case of $G=GL(n, \mathbb{Q}_p)$ and $\rho$ the fundamental representation, $\rho(A_{\pi_p})$ will just be the diagonal matrix (\ref{SatakeParameterGLn}), but in general this need not be the case. 

Moreover, in general one has that for an unramified  global representation $\pi$ of $G(\mathbb{A})$,  only for \emph{all but finitely many} $p$ the local representations $\pi_p$ are spherical, i.e. contain vectors ${\sf f}_p^{\circ}$ fixed under $K_p=G(\mathbb{Z}_p)$. To take care of this complication we let $S$ be a finite set of places such that if $p\notin S$, $\pi_p$ is spherical. The set $S$ always includes the archimedean place $p=\infty$. For this data we now construct the \emphindex[L-function@$L$-function!partial Langlands]{Langlands $L$-function}
\begin{align}
\label{genLL}
L_S(\pi, s, \rho)=\prod_{p\notin S} \frac{1}{\text{det}\left(\idm-\rho(A_{\pi_p})p^{-s}\right)}.
\end{align}
In this situation the analytic continuation is more involved but Langlands has conjectured that $L_S(\pi, s, \rho)$ can be completed at the unramified places $S$ to obtain a 
meromorphic function $L^\compl(\pi, s, \rho)$ of $s$, called the \emphindex[L-function@$L$-function!global Langlands]{global Langlands $L$-function}.

\begin{example}[$L$-function for $G=GL(2,\mathbb{Q}_p)$]
To give a simple example of how such an $L$-function would look like, let us consider $G=GL(2,\mathbb{Q}_p)$ but now take $\rho$ to be the $k$:th symmetric power $Sym^k(\mathbb{C}^2)$ of the fundamental representation $\mathbb{C}^2$ of $GL(2, \mathbb{C})$ (see, e.g., \cite{GelbartMiller} for a nice discussion of this and other examples). The resulting $L$-function  reads 
\beq
\label{eq:symL}
L(\pi, s, Sym^k)=\prod_{p<\infty} \frac{1}{\text{det}(\idm-\rho(A_{\pi_p})p^{-s})}=\prod_{p<\infty}\prod_{j=0}^{k} \frac{1}{1-\alpha_p^{j}\beta_p^{k-j} p^{-s}}.
\eeq
\end{example}
 
Using the formalism outlined above, Langlands thus provided a systematic procedure for attaching $L$-functions to automorphic forms, a task that had previously only been understood in special cases. The relation between automorphic forms on $G$, the $L$-group ${}^LG$ and automorphic $L$-functions provides the cornerstone behind the Langlands program, which are a set of far-reaching conjectures put forward by Langlands, of which only a tiny fraction have been proven. In section \ref{sec:Langlands} we briefly discuss some of the ideas in the Langlands program, and how they relate and extend the theory we have presented in this work.

\section{The Langlands--Shahidi method*}
\label{sec:LSmethod}
\index{Langlands--Shahidi method}

It is important to study the functional properties of $L$-functions such as~\eqref{eq:symL} since these can be used to give estimates on Hecke eigenvalues (or  Fourier coefficients) of cusp forms. This application to number theory is reviewed for example in~\cite{ShahidiKorea,ShahidiBook,GelbartMiller}; we will content ourselves here with explaining the basic construction and its relation to Eisenstein series and Whittaker coefficients.

The starting point is the knowledge of the functional equation~\eqref{funrel} for Eisenstein series on $G$ induced from a representation of the Levi subgroup $L$ of some parabolic subgroup $P=LU$ of $G$. From this functional equation and the knowledge how the $L$-function of interest arises in the Fourier expansion one can then deduce properties of the $L$-function. This method was suggested by Langlands in~\cite{LanglandsFE,LanglandsEP} and then developed in detail by Shahidi~\cite{Shahidi1,Shahidi2,Shahidi3,ShahidiKorea,ShahidiBook}.

To motivate the procedure, we look at the Fourier expansion of the $SL(2,\reals)$ Eisenstein series (cf.~\eqref{SL2FC2})
\begin{align}
\label{FESL2}
E(s,z) = y^s + \frac{\xi(2s-1)}{\xi(2s)} y^{1-s} + \frac{2y^{1/2}}{\xi(2s)} \sum_{n\neq 0} |n|^{s-1/2}\sigma_{1-2s}(n) K_{s-1/2}(2\pi|n| y) e^{2\pi i n x}.
\end{align}
The Eisenstein series satisfies the functional equation (cf.~\eqref{SL2funcrel})
\begin{align}
\label{eq:FrelLS}
E(1-s,z) = \frac{\xi(2s)}{\xi(2s-1)} E(s,z)
\end{align}
and the $L$-function whose properties one is interested in is the completed Riemann zeta function $\xi(k)$. As we have seen in chapters~\ref{ch:SL2-fourier} and~\ref{ch:CTF}, this functional equation can be read off from the constant terms of the Eisenstein series and does not require the knowledge of the completed Riemann zeta function beyond its definition in terms of an Euler product.

Additional properties of $\xi(s)$ can be inferred from the first Fourier coefficient ($n=1$). The functional relation~\eqref{eq:FrelLS} for this Fourier coefficient reads
\begin{align}
\frac{1}{\xi(2(1-s))} K_{1/2-s}(2\pi y) = \frac{\xi(2s)}{\xi(2s-1)} \frac{1}{\xi(2s)} K_{s-1/2}(2\pi y).
\end{align}
Using the property $K_t(x)=K_{-t}(x)$ of the modified Bessel function, one deduces that
\begin{align}
\xi(2s-1) = \xi(2-2s)
\quad \Leftrightarrow \quad
\xi(k) = \xi(1-k).
\end{align}
Thus, the functional equation of the completed Riemann zeta function $\xi(s)$ is a consequence of the functional equation of Eisenstein series. One can also deduce the non-vanishing of $\zeta(s)$ on the line $\Re(s)=1$ from the holomorphy (in $s$) of $E(s,z)$ on the line $\Re(s)=0$ and further properties of $\zeta(s)$ from the study of $E(s,z)$~\cite{GelbartMiller}. (The higher Fourier coefficients $n>1$ provide no additional information.)

The more general realisation of this method relies on Eisenstein series on $G$ induced from a cuspidal automorphic representation $\pi_L$ of the Levi factor $L$ of a maximal parabolic subgroup $P=LU\subset G$. We assume that the representation $\pi_L$ is spherical at almost all places $p$. 

As before, we have that at the spherical finite places $p$ one can characterize the representation by means of its Satake parameter $A_{\pi_p}\in {}^LA$.  Let also $S$ be a set of places that includes all the non-spherical places and the archimedean one. In the everywhere-unramified case one would have $S=\{\infty\}$. Since ${}^LL$ is a complex linear group, it admits standard finite-dimensional complex representations $\rho_L: {}^LL \to GL(n,\cx)$ where $n$ is the dimension of the representation. For any such pair $(\pi_L, \rho_L)$, the \emphindex[L-function@$L$-function!partial Langlands]{partial Langlands $L$-function} is given by
\begin{align}
\label{partialL}
L_S(s,\pi_L, \rho_L) = \prod_{p\notin S} L_p(s,\pi_L,\rho_L) = \prod_{p\notin S} \frac{1}{\det(\idm-\rho_L(A_{\pi_p}) p^{-s})},
\end{align}
where the determinant is taken in the representation associated with $\rho_L$. Formally, this is the same as the definition~\eqref{genLL} above but this time we have emphasised that this is for the Levi part $L$ of a parabolic subgroup $P$ of $G$.
The \emphindex[L-function@$L$-function!global Langlands]{global Langlands $L$-function} requires the definition of factors for the places $S$ that is less uniform and not known in full generality. Important progress for the global $L$-functions for $GL(n)$ and $SO(2n+1)$ can be found in~\cite{Harris1,HarrisTaylor,Henniart,JiangSoudry}. \\

The virtue of these $L$-functions is that they arise in the Fourier expansion of Eisenstein series induced from a \emph{cuspidal} representation $\pi_L$ of $L$. For an automorphic form $\phi\in \pi_L$ we let
\begin{align}
E(s,\phi,g) = \sum_{\gamma \in P(\rats)\backslash G(\rats) } \phi(\gamma g) \delta_P(\gamma g)^s 
\end{align}
be the \emphindex[Eisenstein series!and Langlands--Shahidi method]{Eisenstein series on $G$ induced from $\phi\in\pi_L$}. $\delta_P(g)$ is here (the trivial extension to $G$ of) the modulus character on $P\subset G$ defined by
\begin{align}
d(l u l^{-1}) = \delta_P(l) du
\end{align}
and $\delta_P(ulk)= \delta_P(l)$. It can be given explicitly by $\delta_P(l) = l^{2\rho_P}$, where $\rho_P$ is half the sum of the positive roots contained in $U$.  In the discussion above in section~\ref{nonminEis}, we had taken $\phi=1$ in the non-cuspidal trivial representation.

The Eisenstein series $E(s,\phi,g)$ on $G$ has a Fourier expansion with respect to the unipotent $U$ that is simpler due to the fact that $\phi$ is taken from a \emph{cuspidal} representation of the Levi factor $L$. This arises because in the Bruhat decomposition of $G$ most classes have a vanishing contribution as $\phi$ is cuspidal. This is a collapse mechanism not unsimilar to the one discussed for constant terms in section~\ref{EvalLCF} and for Whittaker coefficients in section~\ref{sec:degpsi}.

Langlands showed~\cite{LanglandsEP} that the constant term of $E(s,\phi,g)$ along $P'$ (the opposite of $P$) is controlled by partial $L$-functions~\eqref{partialL} and Shahidi extended this to non-trivial Fourier coefficients~\cite{Shahidi1,Shahidi2,Shahidi3,ShahidiKorea,ShahidiBook}. Shahidi's work relies also on the Casselman--Shalika formula for (generic) Whittaker coefficients of an Eisenstein series $E(\lambda,g)$ at unramified places~\eqref{CScharacter}.
 
We first explain why Langlands $L$-functions arise from formula~\eqref{CScharacter}. If~\eqref{CScharacter} is evaluated for the special case of $A_{\pi_p}=\id$, corresponding to the trivial representation, one obtains
\begin{align}
W^\circ(\delta^{1/2}\mu,\id) = \frac{1}{\zeta(\delta^{1/2}\mu)} = \prod_{\alpha^\vee>0} (1- p^{-1-\langle\lambda|\alpha\rangle})
= \prod_{\alpha^\vee>0} (1- p^{-1}\mu(\varpi^{\alpha^\vee}))
\end{align}
Now each $\mu(\varpi^{\alpha^\vee})$ corresponds to the adjoint action of the Satake parameter $A_{\pi_p}$ on the root space of $\alpha^\vee$ which is nothing but the representation of the split torus ${}^LA$ on the Lie algebra ${}^L\mf{n}$. Denoting this action by $\rho:{}^LA \to \mathrm{End}({}^L\mf{n})$, we have that
\begin{align}
\prod_{\alpha^\vee>0} (1-p^{-1} \mu(\varpi^{\alpha^\vee}) = \det(\idm- \rho(A_{\pi_p}) p^{-1} ) 
\end{align}
since the representation $\rho$ of ${}^LA$ decomposes into the direct sum of one-dimensional representations labelled by the positive roots (and the determinant therefore factorises). Hence
\begin{align}
\label{WhittL}
W^\circ(\lambda,\id) = \lint_{N(\ints_p)\backslash N(\rats_p)} E(\lambda,n)\overline{\psi(n)} dn =\frac{1}{L_p(1,\lambda, \rho)},
\end{align}
i.e., an $L$-function of the type~\eqref{partialL}. Here, we have labelled the representation $\pi_L$ of the Levi ${}^LA$ of the minimal parabolic (Borel) ${}^LB$ by its quasi-character $\lambda$. 

In the more general case of the Eisenstein series $E(s,\phi,g)$ induced from a cuspidal representation $\pi_L$ of the Levi subgroup of a maximal parabolic $P=UL\subset G$ one has to consider the adjoint action $\rho$ of ${}^LL$ on the Lie algebra ${}^L\mf{u}$ of the unipotent ${}^LU$. Under this action, ${}^L\mf{u}$ decomposes into a finite number of irreducible representations according to 
\begin{align}
\label{udecLS}
\rho = \bigoplus_{j=1}^m \rho_j,
\end{align}
where $m$ is the maximum coefficient (among all roots of $G$) of the simple root defining the maximal parabolic subgroup $P\subset G$. Shahidi showed~\cite{Shahidi1,ShahidiBook} that it is possible to choose $\phi\in\pi_L$ such that the generic Fourier coefficient of $E(s,\phi,g)$ at $g=\id$ is given (for an unramified place $p$) by
\begin{align}
\lint_{U(\ints_p)\backslash U(\rats_p)} E(s,\phi,u) \overline{\psi(u)} du = \prod_{j=1}^m \frac{1}{L_p(1+a j s,\pi_L,\rho_j)},
\end{align}
where $a$ is a fixed number that depends on the choice of parabolic subgroup. The shifts by $s$ in the argument of the $L$-function (compared to~\eqref{WhittL}) is due to the factor $\delta_P(\gamma g)^s$ in the definition of the Eisenstein series. 

In the constant term, the same $L$-functions appear, cf. the intertwining factors $M(w,\lambda)$ in~\eqref{eq:intertwiner}. Langlands showed~\cite{LanglandsEP} that the intertwiner appearing in the constant term for a place $p\notin S$ is
\begin{align}
\prod_{j=1}^m \frac{L_p(ajs,\pi_L,\rho_j)}{L_p(1+ajs,\pi_L,\rho_j)}.
\end{align} 
Due to the cuspidality of $\pi_L$ this is the only non-trivial coefficient appearing in the constant term and it plays the role of the coefficient of $y^{1-s}$ in~\eqref{FESL2} above. It also appears in the functional equation satisfied by $E(s,\phi,g)$ and allows one to deduce a functional equation for the partial $L$-function $L_S$ obtained from all the places $p\in S$:
\begin{align}
 \prod_{j=1}^m L_S(a j s,\pi_L,r_j) =  \prod_{j=1}^m L_S(1-a j s,\pi_L,\tilde{r}_j) 
 \cdot\prod_{p\in S} C(s,\tilde{\pi}_v),
\end{align}
which is called the \emphindex[functional equation!crude]{crude functional equation}~\cite{Shahidi1,Shahidi2}. The factors $C(s,\tilde{\pi}_v)$ appearing in this functional relation are called \emphindex[local factor]{local factors} and they can be determined from the study of the intertwining operator for $p\in S$~\cite{Shahidi1,Shahidi2}. The tildes in the above formula refer to the parabolic subgroup $\tilde{P}$ opposite to $P$. From this identification of the product of $m$ partial $L$-functions as a Fourier coefficient of an Eisenstein series one can also deduce that the product extends to a meromorphic function (in $s$) and does not vanish on the imaginary axis~\cite{ShahidiBook}. Moreover, it is possible to perform induction on $m$ to deduce the same statements for each of the individual factors. This produces a host of non-trivial results for generalised $L$-functions in various representations $r_i$ that arise from all maximal parabolic subgroups~\cite{Shahidi4}. The results described here for split groups can also be extended to so-called \emphindex[group!quasi-split]{quasi-split groups}~\cite{ShahidiBook,Lai}.

\begin{remark}
Besides the Langlands--Shahidi method just outlined, $L$-functions have also been studied using \emphindex[converse theorems]{converse theorems}, most notably those of Cogdell and Piatetski-Shapiro for the general linear group~\cite{CogdellPS1,CogdellPS2}. The virtue of these converse theorems is that they allow to conclude that an $L$-function satisfying certain technical conditions must be a global $L$-function arising from an automorphic form on the general linear group. Such converse theorems can be seen as the extensions of Hecke's results for $L$-functions that were discussed at the end of section~\ref{app_adelichecke}.

Converse theorems make it possible to deduce \emphindex{Langlands functoriality} in some examples. As will be discussed more in section~\ref{sec:Langlands}, Langlands functoriality deals with the question of transferring automorphic forms from a group $G$ to a group $G'$ in which $G$ is a subgroup. Concretely, one starts from an $L$-function that is tentatively associated with the group $G'$ and takes the \emphindex{Rankin--Selberg product} with automorphic $L$-functions on subgroups $G\subset G'$. If certain technical conditions are fulfilled, one can conclude that there must be a (cuspidal) automorphic representation of $G'$ whose $L$-function is the one under study, thereby lifting the representations from $G$ to $G'$. For a nice discussion of this we refer to~\cite{GelbartMiller}.\label{rmk:converse}
\end{remark}

\chapter{Outlook}
\label{ch:outlook}

\vspace{0.5cm}

\epigraph{\emph{It is a deeper subject than I appreciated and, I begin to suspect, deeper than anyone yet appreciates. \\
To see it whole is certainly a daunting, for the moment even impossible, task.}}{--- Robert P. Langlands${}^\S$}

\nonumberfootnote{${}^\S$A review of Haruzo Hida's p-adic automorphic forms on Shimura varieties.}

\vspace{1cm}

\noindent In this concluding chapter, we collect various topics and further directions that we decided not to include in full detail but that  are active fields of research providing interesting context for the study of automorphic forms. The emphasis of the first topics discussed here is mainly on theoretical physics, then we move on to more mathematical areas. We will be much more sketchy in this chapter and refer to the cited literature for additional details.

\section[String scattering amplitudes and automorphic forms]{String scattering amplitudes and\\ automorphic forms}
\label{sec:outlook-strings}

This section is a continuation of the discussion of string theory and automorphic forms in chapter~\ref{ch:intro-strings} which will elaborate on recent research on the topic.
First, we will summarise some of the results from chapter~\ref{ch:intro-strings}. In section \ref{smallreps} we will discuss how automorphic representations are used to specify the coefficients $\mathcal{E}^{\ordd{D}}_{\gra{p}{q}}$ and what this tells us about their Fourier coefficients.
Section \ref{sec:D6R4} treats the $D^6R^4$ term which, as seen in \eqref{eq:intro-LapD6R4}, has an extra source term, and hints at a theory beyond the automorphic forms discussed here.

As discussed in chapter~\ref{ch:intro-strings}, the analytic part of the four-graviton scattering amplitude can be expanded in $\alpha'$ as \eqref{fourgrav} reproduced here for convenience
\begin{align}
    \label{eq:outlook-fourgrav}
    \mathcal{A}^{(D)}(s,t,u,\epsilon_i;g) = 
    \left[\mathcal{E}^{(D)}_{(0,-1)}(g)\frac{1}{\sigma_3}+ 
        \sum_{p\geq 0} \sum_{q \ge 0} \mathcal{E}^{(D)}_{(p,q)}(g) \sigma_2^p \sigma_3^q
    \right] \mathcal{R}^4 \, .
\end{align}
where $\sigma^2 = s^2 + t^2 + u^2$ and $\sigma^3 = s^3 + t^3 + u^3$.

The coefficients $\mathcal{E}^{\ordd{D}}_{\gra{p}{q}}(g)$ are functions on $\mathcal{M} = G(\ints) \bs G(\reals) / K(\reals)$, where $G(\reals)$ is the classical symmetry group, $K(\reals)$ its maximal compact subgroup, and $G(\ints)$ the discrete U-duality subgroup shown in table \ref{tab:CJ}. They satisfy the differential equations \eqref{eq:intro-fourgravLap} 
\begin{subequations}
    \label{eq:outlook-fourgravLap}
    \begin{align}
        \label{eq:outlook-LapR4}
        R^4:& &\left(\Delta_{G/K} - \frac{3(11-D)(D-8)}{D-2} \right) \mathcal{E}^{\scriptstyle{(D)}}_{\scriptstyle{(0,0)}}(g) &= 6\pi\delta_{D,8},&\\
        \label{eq:outlook-LapD4R4}
        D^4R^4:& &\left(\Delta_{G/K} - \frac{5(12-D)(D-7)}{D-2} \right) \mathcal{E}^{\scriptstyle{(D)}}_{\scriptstyle{(1,0)}}(g) &= 40\zeta(2)\delta_{D,7}+7\mathcal{E}_{\scriptstyle{(0,0)}}^{\scriptstyle{(6)}} \delta_{D,6},&\\
        \label{eq:outlook-LapD6R4}
        D^6R^4:&\quad&\left(\Delta_{G/K} - \frac{6(14-D)(D-6)}{D-2} \right) \mathcal{E}^{\scriptstyle{(D)}}_{\scriptstyle{(0,1)}}(g) &= -\left(\mathcal{E}_{\scriptstyle{(0,0)}}^{\scriptstyle{(D)}}\right)^2+ 40\zeta(3) \delta_{D,6}&\\
      &&  &\quad +\frac{55}3\mathcal{E}_{\scriptstyle{(0,0)}}^{\scriptstyle{(5)}} \delta_{D,5} + \frac{85}{2\pi} \mathcal{E}_{\scriptstyle{(1,0)}}^{\scriptstyle{(4)}} \delta_{D,4},&\nn
    \end{align}
\end{subequations}
where $\Delta_{G/K}$ is the Laplace--Beltrami operator on $G/K$, and are well behaved in the limits corresponding to cusps in $G/K$.

As was also covered in chapter~\ref{ch:intro-strings}, there is strong evidence from various consistency checks that they are given by (combinations of) maximal parabolic Eisenstein series defined in section~\ref{nonminEis} and in particular in example~\ref{ex_maximalparabolic}. Specifically, for dimensions $D = 5, 4, 3$ corresponding to tori $T^d$ with $d = 5, 6, 7$ in table~\ref{tab:CJ}, if one considers the maximal parabolic subgroups $P$ of $E_{d+1}$ that have semi-simple Levi parts $SO(d, d)$, then the solutions~\eqref{eq:intro-R4D4R4} 
\begin{subequations}
\label{eq:outlook-R4D4R4}
\begin{align}
\label{eq:outlook-R4fn}
   R^4: & & \mathcal{E}^{(D)}_{(0,0)}(g) &= 2\zeta(3) E(\lambda_{3/2}, P, g),&\\
\label{eq:outlook-D4R4fn}
    D^4R^4: & & \mathcal{E}^{(D)}_{(1,0)}(g) &= \zeta(5)E(\lambda_{5/2}, P, g). &
\end{align}
\end{subequations}
to equations~\eqref{eq:outlook-LapR4} and~\eqref{eq:outlook-LapD4R4} are the conjectured coefficient functions appearing in the four graviton scattering amplitude and have been subjected to numerous consistency checks~\cite{Green:2010kv,Pioline:2010kb,Green:2011vz}.  

The character defining the Eisenstein series is given by the weight \eqref{eq:intro-lambda_s}
\begin{align}
    \label{eq:outlook-lambda_s}
\lambda_s = 2s \Lambda_P -\rho,
\end{align}
where $\Lambda_P$ denotes the fundamental weight orthogonal to the Levi subgroup $L$ of $P=LU$.

\subsection{Small representations and string amplitudes}
\label{smallreps}

As automorphic representations the functions $\mathcal{E}_{\gra{0}{0}}^{\ordd{D}}$ and $\mathcal{E}_{\gra{1}{0}}^{\ordd{D}}$ appearing in~\eqref{eq:outlook-R4D4R4} are attached to so-called \emphindex[representation!small]{small representations}. According to~\eqref{GK-dim2}, the functional dimension of the automorphic representation induced from a parabolic subgroup $P(\mathbb{A})\subset G(\mathbb{A})$ is given by 
\begin{align}
\label{GKformula}
\mathrm{GKdim} I_P(\chi) = \dim G-\dim P = \dim U,
\end{align}
where $\chi$ corresponds to a generic character on the parabolic subgroup $P(\mathbb{A})=L(\ads)U(\ads)$. It turns out that for very special choices of the inducing character $\chi$ there may exist unitarizable submodules of $I_P(\chi)$ with smaller functional dimension. For example, it is well-known that for any real semi-simple Lie group $G(\mathbb{R})$ there exists a \emphindex[representation!minimal]{minimal unitary representation} which has the smallest non-trivial functional dimension among all $G$-representations \cite{MR0342049,MR0404366}. 

The notion of a minimal representation also extends to $p$-adic groups $G(\mathbb{Q}_p)$ \cite{MR2123125} and globally one says that an automorphic representation $\pi=\otimes_{p}\pi_p$ of an adelic group $G(\mathbb{A})$ is minimal if at least one local component $\pi_p$ has smallest non-trivial functional dimension \cite{GRS}. 

Minimal representations of a group $G$ are closely related to minimal nilpotent $G$-orbits. Specifically, via  \emphindex[Kirillov's orbit method]{Kirillov's `orbit method'} one can obtain $\pi_{min}$ through the geometric quantisation of the minimal nilpotent orbit $\mathcal{O}_{min}$ \cite{MR1278630}. This implies that there exists a sequence of small $G$-representations with increasing functional dimension associated with nilpotent orbits of smaller dimension than than the regular orbit. See for example \cite{MR1327538} for an analysis pertaining to exceptional Lie groups of real rank $4$. 

Automorphic forms attached to small representations $\pi$ are interesting both from a mathematical and a physical perspective. It was shown in the seminal paper by Ginzburg--Rallis--Soudry \cite{GRS} that automorphic forms in the minimal representation $\pi_{min}$ have very few non-vanishing Fourier coefficients, a fact that has far-reaching consequences. In particular, it allows to describe the complete Fourier expansion very explicitly, a task which is generally very difficult for Lie groups beyond $SL(2)$. One of the main applications of the theory of small representations has been in the context of the so-called \emphindex[theta correspondence|textbf]{theta correspondence} which is a method of lifting automorphic representations from one group $G$ to another $G'$. Ginzburg has also developed a method which uses small automorphic representations for constructing new automorphic $L$-functions (see \cite{MR3161096} for a survey). 

\subsection{Physical interpretation of non-zero modes}

As shown explicitly in various examples of chapter~\ref{ch:working}, some special Eisenstein series $E(\chi,g)$ have constant terms that are only made up of a very small number of terms, even in the case when the group $G$ on which they are defined is of large rank. This requires an apt choice of the inducing character $\chi$. In a similar way the non-zero Fourier modes of these special Eisenstein series also simplify to a great extent, as for instance shown in example~\ref{ex:degen-whitt-Ed} for $G=E_6,E_7$ and $E_7$. Now we consider the automorphic coefficients~\eqref{eq:outlook-R4D4R4} of the $R^4$ and $D^4R^4$ curvature correction terms which arise in string scattering amplitudes at lowest order in the derivative expansion and which are examples of such special Eisenstein series. We draw together some mathematical results discussed throughout this book which lead to the simplifications in the Fourier expansion and in particular provide a physical interpretation of the special properties of the series following~\cite{Pioline:2010kb,Green:2011vz,Green:2010kv}. 

The automorphic coefficients of the $R^4$ and $D^4R^4$ curvature correction terms satisfy certain physical conditions when approaching special limits of the moduli space $\mathcal M=E_{d+1}/K(E_{d+1})$ which they are functions on. We recall that $\mathcal M$ is the moduli space in the case of compactification on a torus $T^d$. The limits are obtained by considering the Fourier expansion of the coefficients given by Eisenstein series with respect to parabolic subgroups which are not the Borel subgroup. More precisely, special physical meaning is attributed to expansions with respect to the maximal parabolic subgroups $P_{\alpha_1}$, $P_{\alpha_2}$ and $P_{\alpha_{d+1}}$ with a labelling of nodes in the Dynkin diagram as in figure~\ref{fig:CJ} and the associated expansions are referred to as the perturbative, M-theory and decompactification limit respectively. 
We will from now on focus on a discussion of the non-zero Fourier modes of the expansion.
Then, given a maximal parabolic subgroup $P_\alpha=L_\alpha U_\alpha$, the non-zero Fourier modes of the expansion along the unipotent $U_\alpha$ encode non-perturbative instanton effects. The modes are of a form typical for instanton contributions, c.f. the form of the instanton action in equation~\eqref{eq:SUGRA-instanton-action} for the case of $d=0$. The integer instanton charges form conjugate pairs with the variables parameterising the unipotent radical $U_\alpha$ of the maximal parabolic subgroup and thus yield the phases for the non-zero modes. The expansion is then given as a sum over the instanton charges. This can for example be explicitly seen in the case of the $R^4$ coefficient for $d=0$ in equation~\eqref{eq:intro-SL2-fourier-expansion}.

There are, however, physical constraints on the instanton charges contributing to the expansion which stem from supersymmetry. More precisely these physical constraints are the so-called Bogomol'nyi--Prasad--Sommerfield (BPS) conditions which derive from the supersymmetry algebra of the theory~\cite{Obers:1998fb}. 
These BPS conditions provide restrictions on a set of states, in our case instantons, which will transform in a supermultiplet. However the presence of these BPS states also breaks a certain amount of the total supersymmetry of the theory.  
Quantities like the $R^4$ and $D^4R^4$ terms then have such BPS conditions associated with them and one explicitly states the fraction of supersymmetry which is broken. The $R^4$ terms is $\tfrac12$-BPS and the $D^4R^4$ term is $\tfrac14$-BPS. 

The \emphindex[instanton!BPS]{BPS instantons} which are present in the Fourier expansion have different physical meaning depending on the parabolic subgroup with respect to which one is expanding. The instantons present in the case of the decompactification limit are interpreted in terms of black holes carrying BPS charges with their euclidean world-lines wrapped around extended dimensions. In the perturbative limit the instantons encode non-perturbative contributions to string perturbation theory, while in the \emphindex{M-theory} limit they provide non-perturbative contributions to eleven-dimensional supergravity. For a more detailed description of the physical interpretation of these instantons we refer the reader to~\cite{Green:2011vz}.

For each curvature correction term the BPS conditions cause some of the instanton configurations to be trivially realised and thus they will not contribute to the non-zero Fourier modes in the expansion of the respective automorphic coefficient. The instantons which do contribute to the non-zero modes lie in orbits of a representative configuration of instantons under the action of the Levi subgroup $L_\alpha$ of the maximal parabolic. This is in direct connection with the discussion of~\emphindex[character!variety orbit]{character variety orbits} in section~\ref{sec:character-variety-orbits}. There it was shown that the Fourier integral $F_\psi$ of an automorphic function transforms as
\begin{align}
F_\psi(\gamma g) = F_{\psi^\gamma}(g)
\end{align} 
under the discrete Levi factor, where $\gamma\in L(\mathbb Q)$ and $\psi^\gamma(u)=\psi(\gamma u \gamma^{-1})$. The transformation thus acts via adjoint action on the argument of the character $\psi$ and in this way the Fourier modes are related to each other in orbits under the Levi subgroup. It is thus sufficient to calculate a single representative mode for each orbit from which the other modes are then obtained by transformation. As already mentioned these orbits are character variety orbits and they have been classified by Miller and Sahi~\cite{MillerSahi} for all classical and exceptional Lie groups. From the physical point of view the instanton charges which parameterise the phases of the Fourier modes lie in BPS orbits which are carved by the action of the Levi subgroup on a representative configuration of instantons for each orbit.

Let us finally connect this discussion of BPS orbits with the special properties of the $R^4$ and $D^4R^4$ automorphic coefficients. For this we will employ some of the theory of wavefront sets as explained in section~\ref{sec:wavefront-sets}, but for the benefit of the reader we summarise some of the discussion here again. As we have seen in section~\ref{smallreps} the Eisenstein series defining the coefficients of the curvature corrections terms are attached to small automorphic representations. The wavefront sets associated to these automorphic representations only have support on a restricted number of so-called coadjoint nilpotent orbits (see~\ref{sec:wavefront-sets}). For the $R^4$ automorphic coefficient these are the trivial and the~\emphindex[nilpotent orbit!minimal]{minimal nilpotent orbit} and for $D^4R^4$ the trivial, minimal and~\emphindex[nilpotent orbit!next-to-minimal]{next-to-minimal nilpotent orbit}. The character variety orbits discussed above which yield the BPS instanton orbits, lie inside in such coadjoint nilpotent orbits. Finally, theorems by M\oe glin--Waldspurger and Matumoto assert that a given Fourier coefficient $F_{\psi^\gamma}$ can only be non-zero if the character variety orbit which it is associated to intersects one of the coadjoint nilpotent orbits on which the wavefront set has support. This is the reasoning from which the strikingly simple structure of the non-zero Fourier modes of special Eisenstein series follows. The connection between wavefront sets, supersymmetry and curvature corrections was also studied in~\cite{Bossard:2014lra,Bossard:2014aea,Bossard:2015uga} and we will come back to this in section~\ref{sec:WFred}.

\begin{example}[Instantons in four space-time dimensions in the decompactification limit]
In four space-time dimensions one can obtain BPS-instantons by compactifying BPS black holes from five space-time dimensions along their one-dimensional world-line. The relevant symmetry groups are $E_7(\reals)$ in four dimensions and $E_6(\reals)$ in five dimensions that is obtained in the parabolic expansion of $E_7(\reals)$ with respect to node $7$ of the Dynkin diagram. BPS black holes in five dimensions carry generalised electro-magnetic charges $\Gamma$ that transform in the $27$-dimensional representation of $E_6(\reals)$ which is the Levi factor in the parabolic decomposition of $E_7$. The classification of BPS black holes has been carried out in~\cite{Ferrara:1997uz,Ferrara:1997ci,Lu:1997bg}, In five space-time dimensions one finds for example that $\tfrac12$-BPS states are in the coset
\begin{align}
E_6(\reals)/(Spin(5,5)\ltimes \reals^{16}) \quad\textrm{of dimension $17$.}
\end{align}
This coset can be found by considering the stabiliser of a $\tfrac12$-BPS charge in the $27$-dimensional representation. Moreover, it agrees exactly with the intersection of the minimal nilpotent orbit of $E_7(\reals)$ with the character variety orbit of the simple Levi $E_6(\reals)$ acting on the $27$-dimensional unipotent as tabulated in~\cite{MillerSahi}. 

Therefore any automorphic function attached to the minimal representation (like the function multiplying the $R^4$ correction in the string amplitude) can only receive non-perturbative corrections from the $\tfrac12$-BPS states in the decompactification limit. In other limits, other dimensions and for $\tfrac14$-BPS states similar statements, \textit{mutatis mutandis}, hold~\cite{Green:2011vz}.
\end{example}

\subsection{\texorpdfstring{$D^6R^4$}{D6R4}-amplitudes and new automorphic forms} 
\label{sec:D6R4}

The inhomogeneous Laplace equation~\eqref{eq:outlook-LapD6R4} for the $D^6R^4$ coupling does not represent a typical $\mathcal{Z}(\mathfrak{g})$-finiteness condition and therefore the coefficient function $\mathcal{E}^{\ordd{D}}_{\gra{0}{1}}(g)$ is not expected to be an automorphic form in the strict sense of definition~\ref{defauto}. Its solutions have nevertheless been investigated recently in detail by Green, Miller and Vanhove in~\cite{Green:2014yxa} (see~\cite{GreenDualHighDeriv,Bossard:2014lra,Bossard:2014aea,DHoker:2014gfa} for earlier and related work). An $SL(2,\ints)$-invariant solution was found and its Fourier expansion has been studied. 

Green, Miller and Vanhove have also succeeded in expressing the solution as a sum over $G(\ints)$-orbits similar to the standard form of Langlands--Eisenstein series~\cite{Green:2014yxa}
\begin{align}
\label{eq:D6R4sol}
\mathcal{E}^{(10)}_{(0,1)}(g) =\sum_{\gamma\in B(\ints)\backslash G(\ints)} \Phi(\gamma g),
\end{align}
where $\Phi: G \to \reals$ is a right $K=SO(2,\reals)$ invariant function and hence can be interpreted as a function on $B(\reals)$. It is furthermore invariant under $B(\ints)$. However, unlike the case of Eisenstein series, the function $\Phi$ is \emph{not} a character on the Borel subgroup $B(\reals)$ but rather a highly non-trivial function. More precisely, it is given by
\begin{align}
\Phi(z) = \frac23 \zeta(3)^2 y^3 +\frac1{9} \pi^2\zeta(3) y + \sum_{n\neq 0} c_n(y) e^{2\pi i n x},
\end{align}
where the $c_n(y)$ are complicated functions involving Bessel functions and rational functions. The Fourier expansion of such a function is not straight-forward due to complications involving Kloosterman sums and nested integrals over Bessel functions. (See appendix~\ref{app:PSKS} for some details on Fourier expansions of Poincar\'e series and Kloosterman sums.)

A proper framework for $G(\ints)$-invariant functions that satisfy differential equations of the type~\eqref{eq:outlook-LapD6R4} appears to be required in string theory. The class of functions extends the notion of automorphic form discussed elsewhere here. The analysis in~\cite{Green:2014yxa} points in the direction of a relation to \emphindex[automorphic distributions]{automorphic distributions}~\cite{MillerSchmid1,MillerSchmid2,MillerSchmid3,Schmid,Unterberger}. 

The function $\mathcal{E}^{\ordd{10}}_{\gra{0}{1}}(g)$ has the following constant terms~\cite[Eq.~(2.25)]{Green:2014yxa}
\begin{align}
\label{D6R4const}
\lint_{N(\ints)\backslash N(\reals)} \mathcal{E}^{(10)}_{(0,1)}(ng) dn = \frac{2\zeta(3)^2}{3} y^3  + \frac{4\zeta(2)\zeta(3)}{3} y + \frac{4\zeta(4)}{y} + \frac{4\zeta(6)}{27 y^3} + \textrm{non-poly. terms in $y$}.
\end{align}
Here, we have used the usual coordinates from section~\ref{sec:SL2} on $SL(2,\reals)$. 
The non-polynomial terms are of the form $\sum_{n>0} a_n e^{-4\pi n y}/y^2$ and do not have an expansion around weak coupling $y\to\infty$. These terms have an interpretation as instanton/anti-instanton bound states. We see that the structure of constant terms is quite different from that of Eisenstein series where, according to the Langlands constant term formula of theorem~\ref{LCFthm}, one has a sum of polynomial terms in $y$ only and the number is bounded from above by the order of the Weyl group $\Weyl$ which would be $|\Weyl|=2$ here.

In terms of string perturbation theory, the four polynomial terms in~\eqref{D6R4const} correspond to contributions from string world-sheets of genus $h=0,1,2,3$. We recognise the genus $h=0$ contribution from~\eqref{fourgravtree}. The genus $h=2$ contribution predicted here was recently compared to a first principles string theory calculation and found to agree~\cite{DHoker:2013eea,DHoker:2014gfa}, where also remarkably a connection to the so-called \emphindex{Zhang--Kawazumi invariant} on the moduli space of genus $h=2$ Riemann surfaces was found~\cite{Zhang,Kawazumi}. As a consequence of equation~\eqref{eq:outlook-LapD6R4},~\cite{DHoker:2013eea,DHoker:2014gfa} discovered that the Zhang--Kawazumi invariant must satisfy a simple Laplace eigenvalue equation on the moduli space of Riemann surfaces. The genus $h=3$ term in~\eqref{D6R4const} has been verified directly from a string perturbation calculation very recently in~\cite{Gomez:2013sla} in the pure spinor formalism.

In terms of wavefront sets and automorphic representations it seems natural to associate the $D^6R^4$ coupling to the (special) nilpotent orbits of type $3A_1$ and $A_2$~\cite{Bossard:2014lra}. A proper interpretation of these wavefront sets for $SL(2,\reals)$ is missing since the largest nilpotent orbit is the regular $A_1$-type orbit.  $D^6R^4$ correction terms have been analysed recently in various dimensions by different methods~\cite{Basu:2007ck,Green:2008bf,Green:2008uj,Green:2010wi,DHoker:2014gfa,Basu:2014hsa,Pioline:2015yea,Bossard:2015uga,Bossard:2015,Pioline:2015qha,Pioline:2015nfa}.

\subsection[Wavefront sets of curvature corrections and their reduction]{Wavefront sets of curvature corrections and their\\ reduction}
\label{sec:WFred}

\begin{figure}[t]
\centering
\def\deltax{1} 
\def\xmin{1}
\def\ymin{1}
 \begin{tikzpicture}
  \draw[-,draw=black,very thick](\xmin+\deltax,\ymin) -- (\xmin+\deltax,\ymin + 2 );
   \draw[-,draw=black,very thick](\xmin+\deltax - 1.5,\ymin + 1) -- (\xmin+\deltax - 1.5,\ymin + 3.5);
    \draw[-,draw=black,very thick](\xmin+\deltax - 1.5,\ymin + 3.5) -- (\xmin+\deltax,\ymin + 2);
  \draw[-,draw=black,very thick](\xmin+\deltax,\ymin + 2) -- (\xmin+\deltax + 1.5,\ymin + 2.5);
    \draw[-,draw=black,very thick](\xmin+\deltax + 1.5,\ymin + 2.5) -- (\xmin+\deltax,\ymin + 4);
    \draw[-,draw=black,very thick](\xmin+\deltax - 1.5,\ymin + 3.5) -- (\xmin+\deltax,\ymin + 4);
\draw[-,draw=black,very thick](\xmin+\deltax,\ymin + 4) -- (\xmin+\deltax,\ymin + 4.5);
  \draw[-,draw=black,very thick](\xmin+\deltax,\ymin) -- (\xmin+\deltax - 1.5,\ymin + 1);
\draw[-,draw=black,very thick] (\xmin+\deltax,\ymin - 1) -- (\xmin+\deltax,\ymin);
\draw[-,draw=black,very thick] (\xmin+\deltax,\ymin - 2) -- (\xmin+\deltax,\ymin - 1);
 \draw[-,draw=black,very thick](\xmin+\deltax,\ymin+4.5) -- (\xmin+\deltax,\ymin+ 5.5 );
  \draw[-,draw=black,very thick](\xmin+\deltax,\ymin+4.5) -- (\xmin+\deltax+1.5,\ymin+ 5.5 );
  \draw[-,draw=black,very thick](\xmin+\deltax,\ymin+4.5) -- (\xmin+\deltax-1.5,\ymin +5.5 );
  \draw[-,draw=black,very thick](\xmin+\deltax-1.5,\ymin+5.5) -- (\xmin+\deltax-1.5,\ymin +9 );
  \draw[-,draw=black,very thick](\xmin+\deltax,\ymin+5.5) -- (\xmin+\deltax-1.5,\ymin +6.5 );
  \draw[-,draw=black,very thick](\xmin+\deltax+1.5,\ymin+5.5) -- (\xmin+\deltax+1.5,\ymin +9 );
  \draw[-,draw=black,very thick](\xmin+\deltax-1.5,\ymin+9) -- (\xmin+\deltax,\ymin +8 );
  \draw[-,draw=black,very thick](\xmin+\deltax,\ymin+5.5) -- (\xmin+\deltax,\ymin +8 );
   \draw[-,draw=black,very thick](\xmin+\deltax-1.5,\ymin +5.5) -- (\xmin+\deltax-1.2,\ymin + 6);
   \draw[-,draw=black,very thick](\xmin+\deltax-0.975,\ymin +6.375) -- (\xmin+\deltax,\ymin + 8);
   \draw[-,draw=black,very thick](\xmin+\deltax,\ymin +7) -- (\xmin+\deltax+1.5,\ymin + 9);
      \draw[-,draw=black,very thick](\xmin+\deltax,\ymin +8) -- (\xmin+\deltax+1.5,\ymin + 9);
   
\draw[dashed,draw=black,very thick](\xmin+\deltax+1.5,\ymin +9) -- (\xmin+\deltax+1.5,\ymin + 9.5);
\draw[dashed,draw=black,very thick](\xmin+\deltax+1.5,\ymin +9) -- (\xmin+\deltax+1,\ymin + 9.5);
     \draw[dashed,draw=black,very thick](\xmin+\deltax-1.5,\ymin +9) -- (\xmin+\deltax-1.5,\ymin + 9.5);

 \draw (\xmin+\deltax,\ymin) node{\textbullet};
 \draw (\xmin+\deltax,\ymin - 1) node{\textbullet};
  \draw (\xmin+\deltax,\ymin - 2)  node{\textbullet};
  \draw[fill=white,thick] (\xmin+\deltax,\ymin + 2) circle(.5ex);
  \draw (\xmin+\deltax,\ymin + 4)  node{\textbullet};
  \draw (\xmin+\deltax + 1.5,\ymin + 2.5)  node{\textbullet};
   \draw[fill=white,thick] (\xmin+\deltax - 1.5,\ymin + 3.5)  circle(.5ex);
    \draw (\xmin+\deltax - 1.5 ,\ymin + 1) node{\textbullet};
     \draw (\xmin+\deltax,\ymin + 4.5)  node{\textbullet};
        \draw (\xmin+\deltax,\ymin + 5.5)  node{\textbullet};  
             \draw (\xmin+\deltax+1.5,\ymin + 5.5)  node{\textbullet};  
             \draw (\xmin+\deltax-1.5,\ymin + 5.5)  node{\textbullet};  
             \draw (\xmin+\deltax-1.5,\ymin + 6.5)  node{\textbullet};  
             \draw[fill=white,thick] (\xmin+\deltax,\ymin + 7)  circle(.5ex);
             \draw[fill=white,thick] (\xmin+\deltax,\ymin + 8)  circle(.5ex);
             \draw[fill=white,thick] (\xmin+\deltax-1.5,\ymin + 9)  circle(.5ex);
       \draw (\xmin+\deltax+1.5,\ymin + 9)  node{\textbullet};

\draw (\xmin+\deltax+0.4 ,\ymin - 2) node{$0$};
\draw (\xmin+\deltax+0.5,\ymin - 1) node{$A_1$};
\draw (\xmin+\deltax+0.5,\ymin) node{$2A_1$};
\draw (\xmin+\deltax-2.2,\ymin + 1) node{$(3 A_1)^{\prime\prime}$};
\draw (\xmin+\deltax+0.8,\ymin + 1.8) node{$(3 A_1)^{\prime}$};
\draw (\xmin+\deltax+1.9,\ymin + 2.5) node{$A_2$};
\draw (\xmin+\deltax-2.1,\ymin + 3.5) node{$4A_1$};
\draw (\xmin+\deltax+0.9,\ymin + 4) node{$A_2A_1$};
\draw (\xmin+\deltax+1.0,\ymin + 4.5) node{$A_22A_1$};
\draw (\xmin+\deltax+0.5,\ymin + 5.5) node{$2A_2$};
\draw (\xmin+\deltax+2.5,\ymin + 5.5) node{$A_23A_1$};
\draw (\xmin+\deltax-2.0,\ymin + 5.5) node{$A_3$};

\draw (\xmin+\deltax-2.6,\ymin + 6.5) node{$(A_3A_1)^{\prime\prime}$};
\draw (\xmin+\deltax+0.8,\ymin + 7) node{$2A_2A_1$};
\draw (\xmin+\deltax-0.8,\ymin + 7.8) node{$(A_3A_1)^{\prime}$};
\draw (\xmin+\deltax-2.5,\ymin + 9) node{$A_32A_1$};
\draw (\xmin+\deltax+2.4,\ymin + 9) node{$D_4(a_1)$};

\def\deltax{8}

  \draw[-,draw=black,very thick](\xmin+\deltax,\ymin) -- (\xmin+\deltax,\ymin + 2 );
   \draw[-,draw=black,very thick](\xmin+\deltax - 1.5,\ymin + 1) -- (\xmin+\deltax - 1.5,\ymin + 3.5);
    \draw[-,draw=black,very thick](\xmin+\deltax - 1.5,\ymin + 3.5) -- (\xmin+\deltax,\ymin + 2);
  \draw[-,draw=black,very thick](\xmin+\deltax,\ymin + 2) -- (\xmin+\deltax + 1.5,\ymin + 2.5);
    \draw[-,draw=black,very thick](\xmin+\deltax + 1.5,\ymin + 2.5) -- (\xmin+\deltax,\ymin + 4);
    \draw[-,draw=black,very thick](\xmin+\deltax - 1.5,\ymin + 3.5) -- (\xmin+\deltax,\ymin + 4);
\draw[-,draw=black,very thick](\xmin+\deltax,\ymin + 4) -- (\xmin+\deltax,\ymin + 4.5);
  \draw[-,draw=black,very thick](\xmin+\deltax,\ymin) -- (\xmin+\deltax - 1.5,\ymin + 1);
\draw[-,draw=black,very thick] (\xmin+\deltax,\ymin - 1) -- (\xmin+\deltax,\ymin);
\draw[-,draw=black,very thick] (\xmin+\deltax,\ymin - 2) -- (\xmin+\deltax,\ymin - 1);
 \draw[-,draw=black,very thick](\xmin+\deltax,\ymin+4.5) -- (\xmin+\deltax,\ymin+ 5.5 );
  \draw[-,draw=black,very thick](\xmin+\deltax,\ymin+4.5) -- (\xmin+\deltax+1.5,\ymin+ 5.5 );
  \draw[-,draw=black,very thick](\xmin+\deltax,\ymin+4.5) -- (\xmin+\deltax-1.5,\ymin +5.5 );
  \draw[-,draw=black,very thick](\xmin+\deltax-1.5,\ymin+5.5) -- (\xmin+\deltax-1.5,\ymin +9 );
  \draw[-,draw=black,very thick](\xmin+\deltax,\ymin+5.5) -- (\xmin+\deltax-1.5,\ymin +6.5 );
  \draw[-,draw=black,very thick](\xmin+\deltax+1.5,\ymin+5.5) -- (\xmin+\deltax+1.5,\ymin +9 );
  \draw[-,draw=black,very thick](\xmin+\deltax-1.5,\ymin+9) -- (\xmin+\deltax,\ymin +8 );
  \draw[-,draw=black,very thick](\xmin+\deltax,\ymin+5.5) -- (\xmin+\deltax,\ymin +8 );
   \draw[-,draw=black,very thick](\xmin+\deltax-1.5,\ymin +5.5) -- (\xmin+\deltax-1.2,\ymin + 6);
   \draw[-,draw=black,very thick](\xmin+\deltax-0.975,\ymin +6.375) -- (\xmin+\deltax,\ymin + 8);
   \draw[-,draw=black,very thick](\xmin+\deltax,\ymin +7) -- (\xmin+\deltax+1.5,\ymin + 9);
      \draw[-,draw=black,very thick](\xmin+\deltax,\ymin +8) -- (\xmin+\deltax+1.5,\ymin + 9);
   
\draw[dashed,draw=black,very thick](\xmin+\deltax+1.5,\ymin +9) -- (\xmin+\deltax+1.5,\ymin + 9.5);
\draw[dashed,draw=black,very thick](\xmin+\deltax+1.5,\ymin +9) -- (\xmin+\deltax+1,\ymin + 9.5);
     \draw[dashed,draw=black,very thick](\xmin+\deltax-1.5,\ymin +9) -- (\xmin+\deltax-1.5,\ymin + 9.5);

 \draw (\xmin+\deltax,\ymin) node{\textbullet};
 \draw (\xmin+\deltax,\ymin - 1) node{\textbullet};
  \draw (\xmin+\deltax,\ymin - 2)  node{\textbullet};
  \draw[fill=white,thick] (\xmin+\deltax,\ymin + 2) circle(.5ex);
  \draw (\xmin+\deltax,\ymin + 4)  node{\textbullet};
  \draw (\xmin+\deltax + 1.5,\ymin + 2.5)  node{\textbullet};
   \draw[fill=white,thick] (\xmin+\deltax - 1.5,\ymin + 3.5)  circle(.5ex);
    \draw (\xmin+\deltax - 1.5 ,\ymin + 1) node{\textbullet};
     \draw (\xmin+\deltax,\ymin + 4.5)  node{\textbullet};
        \draw (\xmin+\deltax,\ymin + 5.5)  node{\textbullet};  
             \draw (\xmin+\deltax+1.5,\ymin + 5.5)  node{\textbullet};  
             \draw (\xmin+\deltax-1.5,\ymin + 5.5)  node{\textbullet};  
             \draw (\xmin+\deltax-1.5,\ymin + 6.5)  node{\textbullet};  
             \draw[fill=white,thick] (\xmin+\deltax,\ymin + 7)  circle(.5ex);
             \draw[fill=white,thick] (\xmin+\deltax,\ymin + 8)  circle(.5ex);
             \draw[fill=white,thick] (\xmin+\deltax-1.5,\ymin + 9)  circle(.5ex);
       \draw (\xmin+\deltax+1.5,\ymin + 9)  node{\textbullet};

\draw (\xmin+\deltax+0.4 ,\ymin - 2) node{$R$};
\draw (\xmin+\deltax+0.5,\ymin - 1) node{$R^4$};
\draw (\xmin+\deltax+0.7,\ymin) node{$D^4R^4$};
\draw (\xmin+\deltax-2.2,\ymin + 1) node{$D^6R^4$};
\draw (\xmin+\deltax+2.1,\ymin + 2.5) node{$D^6R^4$};

\draw (\xmin - 3 + 0.2 - 1,\ymin - 2) node{$0$};
\draw (\xmin - 3 + 0.3 - 1,\ymin - 1) node{$34$};
\draw (\xmin - 3 + 0.3 - 1,\ymin) node{$52$};
\draw (\xmin - 3 + 0.3 - 1,\ymin + 1) node{$54$};
\draw (\xmin - 3 + 0.3 - 1,\ymin + 2) node{$64$};
\draw (\xmin - 3 + 0.3 - 1,\ymin + 2.5) node{$66$};
\draw (\xmin - 3 + 0.3 - 1,\ymin + 3.5) node{$70$};
\draw (\xmin - 3 + 0.3 - 1,\ymin + 4) node{$76$};
\draw (\xmin - 3 + 0.3 - 1,\ymin + 4.5) node{$82$};
\draw (\xmin - 3 + 0.3 - 1,\ymin + 5.5) node{$84$};
\draw (\xmin - 3 + 0.3 - 1,\ymin + 6.5) node{$86$};
\draw (\xmin - 3 + 0.3 - 1,\ymin + 7) node{$90$};
\draw (\xmin - 3 + 0.3 - 1,\ymin + 8) node{$92$};
\draw (\xmin - 3 + 0.3 - 1,\ymin + 9) node{$94$};

\end{tikzpicture}

\caption{\label{fig:E7closure}\sl The smallest nilpotent orbits of $\mf{e}_7(\cx)$ and their closure ordering. The vertical axis is the dimension of the orbit and on the left they are labelled according to the Bala--Carter classification where we have denoted $2A_2+A_1\equiv 2A_2A_1$ etc. for brevity. The open circles indicate orbits that are not special.  The figure is adapted from~\cite{Spaltenstein,Bossard:2014lra,Bossard:2015}. On the right, the wavefront sets of the various curvature terms appearing in the four-graviton scattering amplitude are shown on the same kind of diagram.
}
\end{figure}
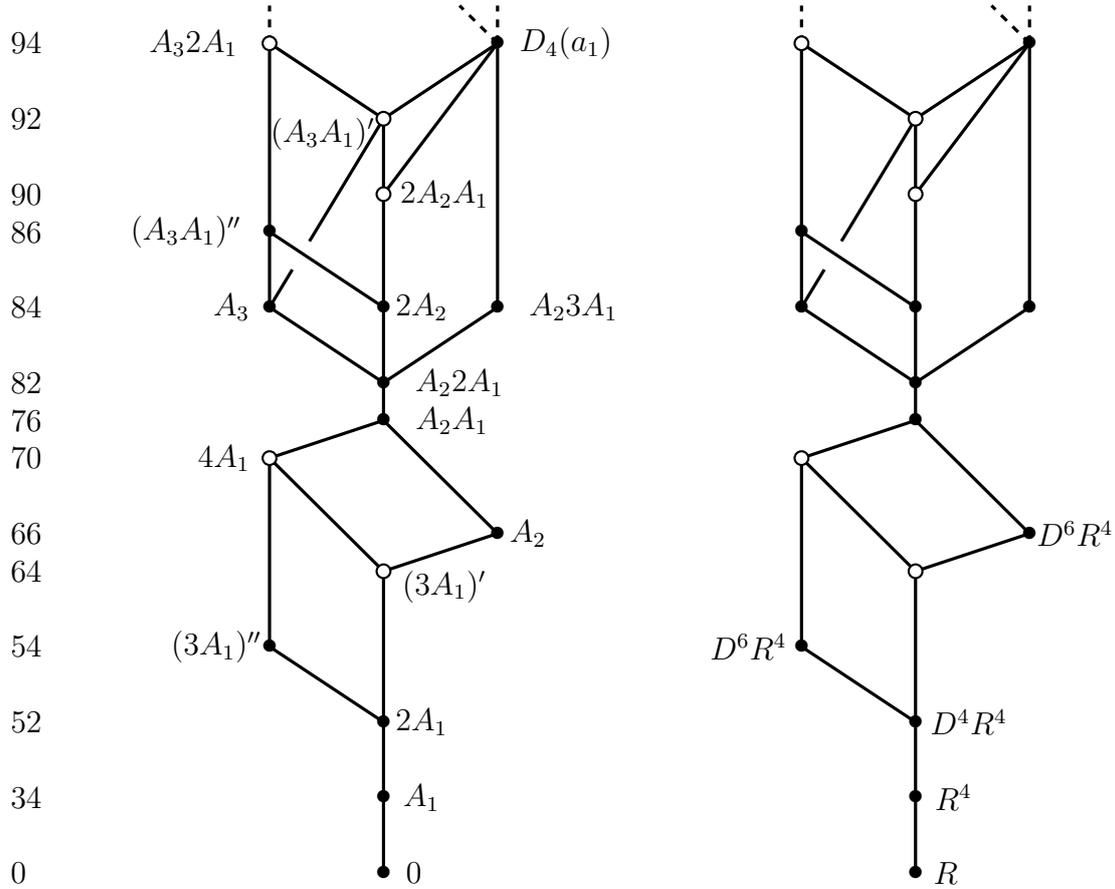

In this section, we would like to collect and systematize some of the remarks on wavefront sets and curvature corrections that have been made in the preceding discussion. We will do this for the case $G=E_{7}(\reals)$ that is relevant for $D=4$ space-time dimensions and maximal supersymmetry. The \emphindex{closure diagram} of nilpotent orbits of $\mf{e}_7(\cx)$ can be found for example in~\cite{Spaltenstein} and that of $\mf{e}_{7}(\reals)$ in~\cite{DjokovicE7}. We display the closure (or \index{Hasse diagram}Hasse) diagram of the smallest nilpotent orbits of $\mf{e}_7(\cx)$ in figure~\ref{fig:E7closure}.

In the figure, we have also shown the wavefront sets of the various types of curvature corrections $D^{2k}R^4$ following~\cite{Pioline:2010kb,Green:2011vz,Bossard:2014lra,Bossard:2015oxa,Bossard:2015}. What is noticeable is that the wavefront sets appear to be associated only with \emp{special} orbits~\cite{Bossard:2015oxa}. Preliminary investigations of higher derivative terms in~\cite{Bossard:2015} suggest that correction terms with more than six derivatives acting on $R^4$ will generically have contributions from the orbit $(A_3+A_1)''$. The expansion in increasing orders of derivatives seems to be related to an expansion in terms of size of the associated wavefront set (with only special orbits as maximal orbits). 

We note also that there can be several maximal nilpotent orbits contributing to a given curvature correction, as in the case of $D^6R^4$. This is related to the fact the U-duality invariant functions $\mathcal{E}_{\gra{p}{q}}^{\ordd{D}}$ that arise are not necessarily automorphic functions of the standard type but more general as discussed in section~\ref{sec:D6R4}. The branching of the diagram is associated in physics with the existence of independent (linearised) supersymmetry invariants~\cite{Bossard:2015uga}.

Let us also relate this discussion back to the analysis of small representations of sections~\ref{sec:CSeval} and~\ref{smallreps}. In the case of $E_7(\reals)$ one has a degenerate principal series representation of functional dimension $33$ that can be realised with Eisenstein series by choosing the weight $\lambda=2s \Lambda_1-\rho$ and we write the associated Eisenstein series as
\begin{align}
E\left( \DEVII{s}{0}{0}{0}{0}{0}{0})\right).
\end{align}
The wavefront set in the case of generic $s$ is of type $A_2$ of dimension $66$. This is twice the dimension of the coset $P_1\backslash E_7$ where $P_1$ is the maximal parabolic subgroup associated with node $1$ in the Dynkin diagram of figure~\ref{fig:CJ}. Reductions occur in this case for the values $s=\tfrac52$ and $s=\tfrac32$ when the wavefront set collapses to the $2A_1$ and $A_1$, respectively. The reduction can be analysed using theorem~\ref{DegWhittThm} as discussed in section~\ref{sec:CSeval}. The two cases where the wavefront set reduces corresponds to the $R^4$ and $D^4R^4$ curvature correction. There is a contribution of this function to the $D^6R^4$ correction for a non-special value of $s$.

\begin{remark}
The most well-studied case of curvature corrections is that of $D=10$ type IIB superstring theory where the symmetry group is $SL(2,\reals)$ with U-duality group $SL(2,\ints)$. The set of nilpotent orbits of $\mf{sl}(2,\cx)$ is very degenerate and consists only of either the trivial or the regular ($A_1$-type) orbit. Nevertheless, the various curvature corrections of type $D^{2k} R^4$ come with very specific orders $s$ of the non-holomorphic Eisenstein series $E(s,z)$. It is an open problem to understand the specific values that appear, notably $s=\tfrac32$ and $s=\tfrac52$, from a mathematical point of view as there seems to be nothing special happening for the automorphic representation for these values.
\end{remark}

\section{Automorphic functions and lattice sums}
\label{sec:latticesums}

As discussed in the introduction, the non-holomorphic Eisenstein series $E(\chi_s,z)$ of $SL(2,\reals)$ (cf.~\eqref{Eisenintro2}) can be equivalently written in terms of a sum over an integral lattice:
\begin{align}
\label{eq:orbitlattice}
E(\chi_s,z) = \sum_{\gamma\in B(\ints)\backslash SL(2,\ints)} \chi_s(\gamma\cdot z) = \frac{1}{\zeta(2s)} \sum_{(c,d)\in \ints^2 \atop (c,d)\neq (0,0)} \frac{y^2}{|cz+d|^{2s}}.
\end{align}
In physics, writing the Eisenstein as a \emphindex[automorphic form!as lattice sum]{lattice sum} can sometimes be interpreted as the sum over the lattice of all possible charges that define the U-duality group (cf.~chapter~\ref{ch:intro-strings}). The sum over the group coset $B(\ints)\backslash G(\ints)$, on the other hand, can be interpreted as the contribution from a single U-duality orbit, if $G(\ints)$ is the U-duality group $G(\ints)$ of table~\ref{tab:CJ} in chapter~\ref{ch:intro-strings}. From the latter point of view, the functions discussed in~\eqref{eq:outlook-R4D4R4} represent simply the U-duality orbit of the perturbative tree level scattering amplitude, whereas the function in~\eqref{eq:D6R4sol} is the U-duality orbit of a finite number of perturbative terms and an infinite number of non-perturbative terms.

Having a representation of an automorphic function as lattice sum can be physically intuitive and it certainly opens up the possibility of employing Poisson resummation for performing the Fourier expansion of the function, as is done in the $SL(2,\reals)$ example in appendix~\ref{app:SL2Fourier}. 

Lattice sums for more general groups $G$ were considered by Obers and Pioline in~\cite{Obers:1999um}. They write the group element $g\in G(\reals)$ in some linear finite-dimensional representation $\mathcal{R}$. In the same representation, a lattice $\Lambda_{\mathcal{R}}$ is embedded that is preserved by the action of $G(\ints)$. This can be constructed for example by starting from the highest weight vector in the representation $\mathcal{R}$. One can form a scalar invariant by considering
\begin{align}
\label{invScalar}
|| g^{-1} \omega||^2,
\end{align}
where $\omega\in \Lambda_{\mathcal{R}}$ and the norm is computed using the $K(\reals)$-invariant inner product on $\mathcal{R}$. In the example~\eqref{eq:orbitlattice} above, this is realised by working in the two-dimensional representation, letting $\omega=\begin{pmatrix}-d\\c\end{pmatrix}\in\ints^2$ and using the Euclidean norm. Then
\begin{align}
\sum_{(c,d)\in \ints^2 \atop (c,d)\neq (0,0)} \frac{y^2}{|cz+d|^{2s}} =\sum_{0\neq \omega\in \Lambda_{\mathcal{R}}} ||g^{-1}\omega||^{-2s}.
\end{align}

The quantity~\eqref{invScalar} is well-defined by construction on $G/K$ for any $\omega\in\Lambda_{\mathcal{R}}$ and one can form a $G(\ints)$-invariant function very generally by letting
\begin{align}
\tilde{\mathcal{E}}_{\mathcal{R},s}(g) =\sum_{0\neq \omega\in \Lambda_{\mathcal{R}}} ||g^{-1}\omega||^{-2s}.
\end{align}
This function is $K$-finite, of moderate growth and $G(\ints)$ invariant. Moreover, it is directly amenable to Poisson resummation on the lattice $\Lambda_{\mathcal{R}}$ and this has been exploited widely to obtain results about the constant terms and also partly the non-constant terms of $\tilde{\mathcal{E}}_{\mathcal{R},s}(g)$~\cite{Lambert:2006ny,Lambert:2010pj,Gubay:2010nd,Bao:2013fga}. 

For some groups $G(\reals)$ and some representations $\mathcal{R}$ it can happen that the function $\tilde{\mathcal{E}}_{\mathcal{R},s}(g)$ is proportional to a (maximal parabolic) Eisenstein series as defined in sections~\ref{sec_Eisenstein} and~\ref{nonminEis}, and as is the case in the $SL(2,\reals)$ example in~\eqref{eq:orbitlattice} above. However, as was already emphasised in~\cite{Obers:1999um}, the function $\tilde{\mathcal{E}}_{\mathcal{R},s}(g)$ will in general \emph{not} be an eigenfunction of the ring of invariant differential operators, i.e., it will not be $\mathcal{Z}(\mf{g})$-finite and hence not a proper automorphic form.

The failure of being an automorphic can be remedied by restricting the lattice sum over $\mathcal{R}$ to an appropriate $G(\ints)$-invariant subset. Such a subset can be found for example by considering the symmetric tensor product $\mathcal{R}\otimes \mathcal{R}$ and then projecting to the largest invariant subspace in there~\cite{Obers:1999um}. The symmetric tensor product arises because~\eqref{invScalar} is computing a symmetric quantity in the $\mathcal{R}$-valued $\omega$. An automorphic form is then given by
\begin{align}
\mathcal{E}_{\mathcal{R},s}(g) =\sum_{0\neq \omega\in \Lambda_{\mathcal{R}}}\delta(\omega\otimes\omega) ||g^{-1}\omega||^{-2s},
\end{align}
where $\delta(\omega\otimes\omega)$ projects on the invariant subspace in $\mathcal{R}\otimes\mathcal{R}$ defined above. In physical applications, this projection has the interpretation of implementing certain conditions that are called \emphindex{BPS conditions} and correspond to considering contributions only from a subset of all instantonic states. The presence of the projection $\delta(\omega\otimes \omega)$ in the sum often makes the direct application of Poisson resummation impossible and renders the Fourier expansion much more difficult. Examples where the full Fourier expansion of a constrained sum was carried out can be found in~\cite{Bao:2009fg,Bao:2010cc}.

Another way of turning $\tilde{\mathcal{E}}_{\mathcal{R},s}(g)$ into an automorphic form is by restricting to a single $G(\ints)$-orbit within $\Lambda_{\mathcal{R}}$ and this leads back to Langlands' definition.

\section{Asymptotics of Fourier coefficients}
\label{sec:asympt}

In applications to physics one is often interested in extracting the \index{Fourier coefficient!asymptotic behavior} asymptotic behavior of Fourier coefficients. Asymptotic here refers to a chosen direction on the symmetric space $G(\reals)/K(\reals)$ on which the Eisenstein series and their Fourier coefficients are defined. The symmetric space often has the interpretation as the \emphindex[moduli!space]{moduli space}, a point of which corresponds to the \emphindex[vacuum expectation value]{vacuum expectation values} of some fields in a physical theory, cf.~the discussion in chapter~\ref{ch:intro-strings}. In an effective description these values are turned into \index{coupling constant} coupling constants and similar parameters of the theory and therefore sending a certain coupling constant to zero corresponds to a limit on the symmetric space $G(\reals)/K(\reals)$. The asymptotic behaviour then reveals the \emphindex[instanton!action]{instanton action} of a particular non-perturbative object in the theory. In the case of $SL(2,\reals)$ this was explained in the introduction in~\eqref{eq:Instaction} and in general one would like to know how the real part of the action behaves asymptotically as a function of the coupling constants. String theory makes definite predictions for these coupling constants that were formalised as conjectures in~\cite[Appendix~A.6]{GarlandMillerPatnaik}. In this section we prove these conjectures by solving the relevant Laplace equations asymptotically.

The interesting coupling constants are associated with directions in the non-compact abelian subgroup $A(\reals)$. Choosing one direction in this space is tantamount to picking a certain weight, or, equivalently, one chooses a maximal parabolic subgroup $P(\reals)\subset G(\reals)$. In order to find the asymptotic behavior of an automorphic function $\varphi(g)$ and its Fourier coefficients in this limit one can analyse the Laplace differential equation that it satisfies. This equation reads
\begin{align}
\Delta_{G/K} \varphi(g) =\mu\varphi(g)
\end{align}
for some eigenvalue $\mu$.

For a maximal parabolic $P(\reals)$ we know from (\ref{maxparLevi}) that there  is a $GL(1,\reals)$ subgroup in the Levi factor $L(\reals)$ and we denote its Cartan subalgebra generator by $d\in\mf{h}$. This element can be used to introduce a grading of $\mf{p}(\reals)=\mf{m}(\reals)\oplus d\reals\oplus \mf{u}(\reals)$ with the properties
\begin{subequations}
\begin{align}
\lb d, m\rb &=0 &&\textrm{for $m\in\mf{m}=\mathrm{Lie}(M)$},\\
\lb d, u_\ell \rb & = \ell u_\ell &&\textrm{for $u_\ell\in\mf{u}_\ell$ with $\mf{u}=\bigoplus_{\ell}\mf{u}_\ell$}.
\end{align}
\end{subequations}
(The decomposition of $\mathfrak{u}$ is the same that arose in~\eqref{udecLS}.) To now consider a `pure instanton' at degree $\ell$ means that we are interested in group elements $g$ of the form $g=e^{\phi d} e^{\chi E_\ell}$ where $E_\ell\in\mf{u}_\ell$ is a chosen fixed generator. Treating the expectation value of $e^\phi$ as the coupling constant, the \emphindex[weak coupling limit]{weak coupling limit} then corresponds to $\phi\to-\infty$. The relevant part of the Laplacian for the directions $\phi$ and $\chi$ is 
\begin{align}
\label{relLap}
\Delta_{G/K} \propto \partial_\phi^2 +\beta \partial_ \phi + e^{-2\ell\phi} \partial_\chi^2 +\ldots,
\end{align}
where $\beta=\sum_\ell \ell \mathrm{dim}(\mf{u}_\ell)$ and we have not fixed the normalisation of the Laplace operator as it can be absorbed into the eigenvalue $\mu$.  For an instanton of charge $m$ we now make the ansatz for the automorphic function that asymptotically for $\phi\to -\infty$
\begin{align}
\varphi_m(\phi,\chi) = e^{-a e^{-b\phi} + 2\pi i m \chi}\left(1 +O(e^\phi)\right).
\end{align}
This correspond to a Fourier coefficient for a character $\psi\left(e^{\chi E_\ell}\right) = e^{2\pi i m \chi}$.

Acting with the relevant part (\ref{relLap}) of the Laplace operator on this ansatz shows that it can only be an eigenfunction (asymptotically) if
\begin{align}
\label{asympt}
a=2\pi |m|\quad\textrm{and}\quad b=\ell.
\end{align}
Note that this reasoning is independent of the eigenvalue $\mu$ and of whether $\varphi$ is an Eisenstein series or any other automorphic function. The important point about (\ref{asympt}) is that it shows that the leading part of the \emphindex[instanton!action]{instanton action} is
\begin{align}
S_{\mathrm{E}}(\phi,\chi) = \log \varphi_m(\phi,\chi) = 2\pi |m| e^{-\ell\phi} + 2\pi i m\chi + \textrm{sub-leading in $e^{\phi}$}.
\end{align}
This is the typical of type $\ell$ instanton where `type' here refers to the degree in $\mf{u}(\reals)$. Making the link to the non-abelian Fourier expansion of section~\ref{sec:FCsec} shows that the more non-abelian a Fourier coefficient is the faster its fall-off in the corresponding weak coupling expansion.

The typical cases encountered in string theory are when $e^\phi= g_{\mathrm{s}}$ is the \emphindex[string coupling]{string coupling}. Instantons with $\ell=1$ are then D-instantons and those with $\ell=2$ are NS-instantons~\cite{Green:1997tv,Green:1997tn,Becker:1995kb,Pioline:2009qt,Bao:2009fg,Bao:2010cc}. In low space-time dimensions one expects also instantons with $\ell>2$~\cite{Obers:1998fb,Englert:2007qb}.

\section[Black hole counting and automorphic representations]{Black hole counting and automorphic\\ representations}
\label{sec:BH-counting}

As explained in chapter~\ref{ch:intro-strings}, string theory compactified on a compact six-dimensional manifold $X$ gives rise to an effective supersymmetric gravitational theory in 4 dimensions. The number of preserved supersymmetries, usually denoted by $\mathcal{N}$, depends on the properties of $X$. Previously we have mainly discussed the case of $X=T^6$, but there are other interesting and relevant manifolds, such as $X=K3\times T^2$ and when $X$ is a Calabi-Yau threefold. The resulting theory has black hole solutions carrying electric and magnetic charges taking values in a lattice $\Gamma$. 

An important observable is the \emph{BPS-index} $\Omega(\gamma)$ which is a function $\Omega : \Gamma \to \mathbb{Z}$ that counts the (signed) degeneracies of a certain class of black holes (called BPS-black holes) with charge vector $\gamma\in \Gamma$. This index then provides a microscopic description of the black hole entropy $S(\gamma)$ via Boltzmann's formula $S(\gamma)=\log \Omega(\gamma)+\cdots$, where the ellipsis denote subleading corrections. 

The BPS-index $\Omega(\gamma)$ holds the key to many interesting connections between string theory and mathematics. The charge lattice $\Gamma$ is nothing but the cohomology lattice $H^{*}(X, \mathbb{Z})$ of the compact manifold $X$ and the index $\Omega(\gamma)$ can roughly be thought of as counting certain submanifolds of $X$ in the cohomology class $[\gamma]\subset H^{*}(X, \mathbb{Z})$. It is therefore naturally related to the enumerative geometry of $X$. Remarkably, the index also provides a link to automorphic forms. To illustrate these statements we shall now consider a few examples. 

\subsection{\texorpdfstring{$\mathcal{N}=8$}{N=8} supersymmetry}

Let us first consider the case when $X=T^6$, the real six-dimensional torus. This leads to $\mathcal{N}=8$ supersymmetry and is the case discussed in section~\ref{sec:outlook-strings}. In four dimensions with electric and magnetic charges taking values in a lattice $\Gamma\cong \mathbb{Z}^{56}$. As reviewed in chapter~\ref{ch:intro-strings} and section~\ref{sec:outlook-strings}, this theory exhibits a classical $E_7(\mathbb{R})$-symmetry which is broken in the quantum theory to the  $E_7(\mathbb{Z})=\{ g\in E_7(\mathbb{R})\, |\, g\Gamma =\Gamma\}$. This implies that the weighted degeneracy $\Omega(\gamma)$ of BPS-black holes of charge $\gamma\in \Gamma$ must be invariant under $E_7(\mathbb{Z})$. 

However, not all black holes have charges supported on the entire lattice $\Gamma$. For example, the $\tfrac12$-BPS black holes preserve half of the supersymmetries of the theory and can only have charges supported on a 28-dimensional (Lagrangian) subspace $\mathcal{C}_{1/2}\subset \Gamma$. Similarly, $\tfrac14$-BPS black holes have support on a 45-dimensional subspace $\mathcal{C}_{1/4}\subset \Gamma$. 

Now denote by $\Omega_{1/A}(\gamma)$ the index counting $\tfrac{1}A$-BPS-black holes ($A=2,4$). Due to the $E_7(\mathbb{Z})$-invariance it is natural to suspect that the index arises as the Fourier coefficient of an automorphic form, constrained so that  $\Omega_{1/A}(\gamma)$ is non-vanishing only when $\gamma\in \mathcal{C}_A$. 

Let us consider the $A=2$ case for illustration. It turns out that all the expected properties are  fulfilled by an automorphic form $\varphi_{min}$ on $E_8(\mathbb{Z})\backslash E_8(\mathbb{R})$ attached to the minimal representation $\pi_{min}$ of $E_8(\mathbb{R})$ \cite{Pioline:2005vi,Gunaydin:2005mx,Pioline:2006ni}. This representation has Gelfand--Kirillov dimension 29 and can thus be realised as the unitary action of $E_8$ on a space of functions of 29 variables \cite{Gunaydin:2001bt}, say $(p,k)\in \mathbb{Z}^{28} \times \mathbb{Z}$. These integers  parametrise characters on the Heisenberg unipotent radical $U_{Heis}\subset E_8$, which has an associated Levi factor $L_{Heis}=E_7\times \mathbb{R}$. The centre  $Z=[U_{Heis}, U_{Heis}]$ is one-dimensional and the integer $k$ parametrises a unitary character $\psi_Z : Z(\mathbb{Z})\backslash \mathbb{Z}(\mathbb{R})\to U(1)$, trivial on the abelianization $Z\backslash U_{Heis}$. On the other hand the electric and magnetic charges $\gamma=(p,q)\in \Gamma\cong \mathbb{Z}^{56}$ parametrise characters $\psi : U_{Heis}(\mathbb{Z})\backslash U_{Heis}(\mathbb{R}) \to U(1)$, trivial on $Z(\mathbb{R})$. Consider the constant term of $\varphi_{min}$ with respect to $Z$:
\beq 
\varphi_{Z, min}=\int_{Z(\mathbb{Z})\backslash Z(\mathbb{R})} \varphi_{min}(zg)dz.
\eeq
This is a function $\varphi_{Z, min} : E_7(\mathbb{R})\to \mathbb{C}$ invariant under $E_7(\mathbb{Z})$. By taking the constant term with respect to $Z$ we have effectively removed the dependence on the variable $k$. The function $\varphi_{Z, min}$ can be expanded further (see section~\ref{ch:fourier})
\beq
\varphi_{Z, min}(g)=\varphi_{U_{Heis}}+ \sum_{\psi\neq 1} F_{\psi}(\varphi_{min}, g),
\eeq
where 
\beq 
F_{\psi}(\varphi_{min}, g)=\int_{U_{Heis}(\mathbb{Z})\backslash U_{Heis}(\mathbb{R})} \varphi(ug)\overline{\psi(u)}du.
\eeq
In general, such a Fourier coefficient might not be Eulerian (i.e. have an Euler product factorisation); however, as we explained in section~\ref{sec-FourWhitt}, for the minimal representation that turns out to be the case:
\beq
F_\psi(\varphi_{min}, g)=F_{\psi_\infty}(\varphi_{min}, g_\infty) \times \prod_{p<\infty} F_{\psi_p}(\varphi_{min}, g_p). 
\eeq

It was shown in \cite{Kazhdan:2001nx,KazhdanPolishchuk} that these Fourier coefficients indeed have support on the Lagrangian subspace $\mathcal{C}_{1/2}$. We can now state the relation to 1/2 BPS black holes: 
\begin{conjecture}[\hspace{1sp}\cite{Pioline:2005vi,Gunaydin:2005mx,Pioline:2006ni}] The index $\Omega_{1/2}(p,q)$ counting charged 1/2 BPS black holes in four-dimensional, $\mathcal{N}=8$ supergravity is given by
\beq
\Omega_{1/2}(p,q)= \prod_{p<\infty} F_{\psi_p}(\varphi_{min}, 1),
\eeq
where $F_{\psi_p}(\varphi_{min}, 1)$ is the $p$-adic spherical vector in the minimal representation $\pi_{min}$ of $E_8$ (obtained in \cite{KazhdanPolishchuk}) and the electric-magnetic charges $(p,q)$ parametrises the character $\psi_p$.
\label{conj:halfBPS}
\end{conjecture}

Similarly, for the 1/4 BPS black holes we have:
\begin{conjecture}[\hspace{1sp}\cite{Pioline:2005vi,Gunaydin:2005mx,Pioline:2006ni}]  The index  $\Omega_{1/4}(p,q)$ counting charged 1/4 BPS black holes is given by
\beq
\Omega_{1/4}(p,q)= \prod_{p<\infty} F_{\psi_p}(\varphi_{ntm}, 1),
\eeq
where $\varphi_{ntm}$ is an automorphic form in the next-to-minimal representation $\pi_{ntm}$ of $E_8$. 
\end{conjecture}

\begin{remark} 
For a certain subclass of 1/2-BPS black holes in $\mathcal{N}=8$ theories there exists an alternative proposal for the degeneracies. More precisely, for black holes which are dual to bound states of D0-D2-D4-D6 branes with charges $(p^\Lambda, q_\Lambda)=(p^0, p^a, q_a, q_0)$, the degeneracies are conjecturally given by \cite{Maldacena:1999bp,Shih:2005qf}
\begin{align}
\Omega_{1/2}(p,q)=\hat{c}\big(I_4(p,q)\big),
\end{align}
where $I_4$ is the quartic $E_7$-invariant, and $\hat{c}$ are the Fourier coefficients of the Jacobi form:
\begin{align}
-\frac{\theta_1^2(\tau, z)}{\eta(\tau)^6}=\sum_{n, \ell} \hat{c}(4n-\ell^2)q^{n}y^\ell, \qquad q:=e^{2\pi i \tau}, \, z:=e^{2\pi i z}.
\end{align}
It would be very interesting to understand the relation between this proposal and Conjecture (\ref{conj:halfBPS}).
\end{remark}

\subsection{\texorpdfstring{$\mathcal{N}=4$}{N=4} supersymmetry}
\label{sec:N=4}

Let us now take $X=K3\times T^2$, where the first factor is a compact $K3$-surface. This yields $\mathcal{N}=4$ supersymmetry in 4 dimensions, which admits $\tfrac12$- and $\tfrac14$-BPS-black holes with electric magnetic charges $\gamma=(p,q)$ taking values in $\Gamma=H^{*}(X, \mathbb{Z})$. The quantum symmetry of this theory is $SL(2,\mathbb{Z}) \times SO(6,22;\mathbb{Z})$~\cite{Schwarz:1993mg,Sen:1994fa,Witten:1995ex} and we are interestested in finding invariant BPS-indices $\Omega_{1/2}(p,q)$ and $\Omega_{1/4}(q,p)$. Mathematically, these indices are counting special Lagrangian submanifolds of $X$ in the class $[\gamma]\subset H^{*}(X, \mathbb{Z})$. As we shall see the counting works quite differently in this case compared to the $\mathcal{N}=8$ theory considered above. 

The 1/2 BPS-states are purely electric $\gamma=(0,q)$ or purely magnetic $\gamma=(p,0)$ and they are known to be exactly counted by \cite{Dabholkar:1989jt,Dabholkar:2004yr}
\beq
\Omega_{1/2}(q,0) = d(q^2/2),
\eeq
where $d(n)$ are the Fourier coefficients of the discriminant function $(\tau\in\UHP$):
\beq
\Delta(\tau)=\frac{1}{\eta(\tau)^{24}} =\sum_{n=-1}^{\infty}d(n) e^{2\pi i n\tau},
\eeq
which is a cusp form of weight 12 for $SL(2,\mathbb{Z})$. Note that the index is automatically invariant under $SO(6, 22;\mathbb{Z})$ since it only depends on the invariant square $q^2=q\cdot q$ of the charge vector $q$. On the other hand the $SL(2,\mathbb{Z})$-part of the quantum symmetry is broken since $(p,q)$ transforms in a doublet. Thus, in order to preserve the full symmetry group we must consider both electric and magnetic charges, as is the case for the $\tfrac14$-BPS-index $\Omega_{1/4}(p,q)$. Moreover, in order to preserve $SO(6,22; \mathbb{Z})$ this can only depend on the invariant combinations $q^2, p^2, p\cdot q$. The answer is that $\Omega_{1/4}(p,q)$ is the Fourier coefficient of the unique weight 10 cusp form for $Sp(4;\mathbb{Z})$, known as the \emphindex{Igusa cusp form} and usually denoted by $\Phi_{10}$. The precise statement is \cite{Dijkgraaf:1996it,Shih:2005uc}:
\beq
\Omega_{1/4}(p,q)=D(q^2/2, p^2/2, p\cdot q),
\eeq
where the numbers $D(m, n, \ell)$ are extracted from the expansion of the inverse of the Igusa cusp form:
\beq 
\frac{1}{\Phi_{10}(\rho, \sigma, \tau)} = \sum_{m,n,\ell}D(m, n, \ell) e^{2\pi i m\sigma} e^{2\pi i n\tau} e^{2\pi i \ell \rho},
\eeq
where $(\rho, \sigma, \tau)$ are complex variables parametrising the \emphindex{Siegel upper half plane}. This can be generalised to orbifolds of $X=K3\times T^2$ by some discrete 
subgroup $\mathbb{Z}_N$, in which case the counting is given by \emphindex{Siegel modular forms} for paramodular groups (see \cite{Sen:2007qy} for a review and further references). 

\subsection{\texorpdfstring{$\mathcal{N}=2$}{N=2} supersymmetry}

Finally, we consider the case when $X$ is a \emphindex[Calabi--Yau space]{Calabi--Yau 3-fold}. This gives rise to $\mathcal{N}=2$ supersymmetry in four space-time dimensions. The lattice $\Gamma$ of electric and magnetic charges is either $H^{\text{even}}(X, \mathbb{Z})$ or $H^3(X, \mathbb{Z})$ depending on whether we consider type IIA or type IIB string theory. According to Kontsevich's \emphindex[homological mirror symmetry]{homological mirror symmetry conjecture} \cite{MR1403918} a BPS black hole with charge $\gamma\in H^{\text{even}}(X, \mathbb{Z})$ can be viewed as a (semi-)stable object in the (bounded) \emphindex{derived category of coherent sheaves} $\mathrm{D^{b}Coh}(X)$, while black holes with charges $\gamma\in H^3(X, \mathbb{X})$ correspond to (semi-)stable objects (special Lagrangians) in the \emphindex[Fukaya category, derived]{derived Fukaya category} $\mathrm{D^{b}Fuk}(X)$. The BPS-index $\Omega : \Gamma \to \mathbb{Z}$ should then be identified with the \emphindex[Donaldson--Thomas invariant]{generalised Donaldson--Thomas invariants} of $X$ \cite{2008arXiv0811.2435K,Gaiotto:2008cd,MR2951762}. String theory predicts that there should be an action of a discrete Lie group $G(\mathbb{Z})$ on the categories $\mathrm{D^{b}Coh}(X)$ and $\mathrm{D^{b}Fuk}(X)$, which is very unexpected from a purely mathematical viewpoint. In general it is not known what the group $G(\ints)$ should be but it must at least contain the ``S-duality'' group $SL(2,\mathbb{Z})$ (see, e.g., \cite{Dabholkar:2005dt,Denef:2007vg}). For certain choices of $X$ there are, however, precise conjectures regarding the nature of $G(\mathbb{Z})$. 

Let $X$ be a rigid Calabi-Yau 3-fold ($h_{2,1}(X)=0$) of CM-type, i.e. admitting complex multiplication by the ring of algebraic integers $\mathcal{O}_d$ in the quadratic number field $\mathbb{Q}(\sqrt{-d})$. In this case the intermediate Jacobian of $X$  is an elliptic curve: 
\beq
H^3(X, \mathbb{R})/H^3(X, \mathbb{Z})\cong \mathbb{C}/\mathcal{O}_d.
\eeq
We then have: 
\begin{conjecture}[\hspace{1sp}\cite{Bao:2009fg,Bao:2010cc}]  For type IIB string theory compactified on a rigid Calabi-Yau 3-fold $X$ of CM-type the ``U-duality group'' $G(\mathbb{Z})$ is the Picard modular group 
\beq
SU(2,1;\mathcal{O}_d):= SU(2,1)\cap GL(3, \mathcal{O}_d).
\eeq
In particular, this group acts on the charge lattice $H^3(X, \mathbb{Z})$ and consequently on $\mathrm{D^{b}Fuk(X)}$.
\end{conjecture}
If correct, this suggests that the BPS-index $\Omega(\gamma)$ should arise as the Fourier coefficient of an automorphic form on $SU(2,1)$ in a similar vein as for $\mathcal{N}=8$ and $\mathcal{N}=4$ supergravity discussed above. Constraints from supersymmetry further imply that there exists a class of 1/2 BPS-states that have support only on charges $\gamma$ such that $\mathcal{Q}_4(\gamma)\geq 0$, where $\mathcal{Q}(\gamma)$ is a quartic polynomial in the charge vector $\gamma$. In other words, the BPS-index is constrained such that 
\begin{equation}
    \Omega(\gamma) =
    \begin{cases}
        n\neq 0 & \mathcal{Q}_4(\gamma) \geq 0 \\ 
        0       & \mathcal{Q}_4(\gamma) < 0.
    \end{cases}
\end{equation}

It turns out that this constraint is precisely satisfied for Fourier coefficients of automorphic forms attached to the quaternionic discrete series of Lie groups $G$ in their quaternionic real form \cite{MR1327538,MR1421947,MR1988198}. This leads to the following:
\begin{conjecture}[\hspace{1sp}\cite{Bao:2009fg,Bao:2010cc,PerssonAuto}]  The generalised Donaldson-Thomas invariants $\Omega(\gamma)$ of a CM-type rigid Calabi-Yau threefold $X$ are captured by the Fourier coefficients of an automorphic form attached to the quaternionic discrete series of $SU(2,1)$.
\end{conjecture}

Another interesting case is when $X$ is a Calabi-Yau threefold with $h_{1,1}=1$. One then expects that the U-duality group is an arithmetic subgroup $G_2(\mathbb{Z})$ of the split real form $G_2(\mathbb{R})$. Automorphic forms on $G_2$ associated with the quaternionic discrete series have been analysed in detail by Gan, Gross, Savin \cite{MR1932327}, and one has:
\begin{conjecture}[\hspace{1sp}\cite{Pioline:2009qt,PerssonAuto}]  There exists Calabi-Yau 3-folds $X$ with $h_{1,1}=1$ whose Donaldson-Thomas invariants $\Omega(\gamma)$ are captured by automorphic forms attached to the quaternionic discrete series of $G_2$, as analysed by Gan-Gross-Savin.
\label{QDS}
\end{conjecture}

\begin{remark}For large values of the charges the index should reproduce the macroscopic entropy of the black hole which is known to be given by $S(\gamma)=\pi \sqrt{\mathcal{Q}_4(\gamma)}$. Translated into mathematics this implies that the Fourier coefficient should have an asymptotic growth given by 
\beq
\Omega(\gamma)\sim e^{\pi\sqrt{\mathcal{Q}_4(\gamma)}} \quad \text{as} \quad \gamma \to \infty.
\label{growth}
\eeq
This gives rise to the following interesting puzzle. In general, Hecke eigenforms always  give rise to Fourier coefficients that grow polynomially, and hence 
the growth in (\ref{growth}) does not seem to be compatible with the fact that the automorphic forms of Gan-Gross-Savin are indeed Hecke eigenforms. One possible resolution to this 
problem is that one should not consider honest automorphic forms in the quaternionic discrete series but rather  some analogue of mock modular forms for $G_2$, a possibility  suggested by Stephen D. Miller. This might also be consistent with the fact that the BPS-index $\Omega(\gamma)$ jumps discontinuously at certain co-dimension one walls in parameter space (known as \emphindex{wall-crossing}) and this phenomenon is closely related to mock modularity (see, e.g., \cite{Manschot:2007ha,Manschot:2009ia,Dabholkar:2012nd}).
\end{remark}

\section{The Langlands program}
\label{sec:Langlands}

Any survey on automorphic forms would be incomplete without at least mentioning some of the key ideas involved in the \emphindex{Langlands program}, 
the collective name given to the visionary conjectures outlined by Langlands in his letter to Weil in 1967 \cite{LanglandsWeilLetter}, and later expanded upon 
in the lecture notes ``Problems in the theory of automorphic forms'' \cite{LanglandsProb}. To give a complete account of these conjectures goes far beyond the 
scope of this survey. 
However, we would like to give a heuristic discussion  of some of the ingredients and 
their implications. This section leans on the discussions in sections~\ref{sec_Lgroup} to~\ref{sec_autoLfunction}. We will also make a few remarks on the geometric version of the Langlands program along with some speculative remarks on the connection with physics.

\subsection{The classical version}

The context of Langlands' letter to Weil was reductive groups $G$ defined over an arbitrary number field $\field$ that can be either local (like $\mathbb{Q}_p$) or global (like $\mathbb{Q}$). 
Let us focus on the global situation. As usual we restrict our treatment to $\field=\mathbb{Q}$, and we let $G$  be a split group over $\mathbb{Q}$; for example $GL(n, \mathbb{Q})$. Recall that being split over $\mathbb{Q}$ means that there exists a maximal torus which is a product of $GL(1, \mathbb{Q})$s.  However, Langlands also
considered groups $G$ that were only \emphindex[group!quasi-split]{quasi-split}, meaning that they contain a Borel subgroup which is defined over $\mathbb{Q}$. Equivalently, a quasi-split group is split over an \emphindex[field extension!unramified]{unramified finite extension} $\fieldext/\field$. We recall that finite extension of a field $\field$ is another field $\fieldext$ that contains $\field$ and which has finite dimension as a vector space over $\field$, so in this case it is a finite-dimensional vector space over $\rats$. The group of automorphisms of the extension $\fieldext$ is called the \emphindex{Galois group} and denoted by $\text{Gal}(\fieldext/\field)$. In this more general context the $L$-group of $G(\mathbb{Q})$ is really defined as the semi-direct product
\beq
{}^LG=\widehat{G}(\mathbb{C})\ltimes \text{Gal}(\fieldext/\field),
\eeq
where the first factor is the complex group that we discussed in section \ref{sec_Lgroup}. In the case when $G$ is split over $\field=\mathbb{Q}$, like for $GL(n, \mathbb{Q})$, the Galois group acts trivially and the $L$-group becomes a direct product ${}^LG=\widehat{G}(\mathbb{C})\times \text{Gal}(E/\mathbb{Q})$. In this situation one can take the representation $\rho\, :\, {}^LG \to GL(n, \mathbb{C})$ that enters in the construction of $L$-functions $L(\pi,s, \rho)$ defined in~\eqref{genLL}, to have a trivial projection on the second Galois factor in ${}^LG$, and we therefore recover the description of $L$-functions in section \ref{sec_autoLfunction} where we had simply assumed ${}^LG=\widehat{G}(\mathbb{C})$,  see also remark~\ref{rmk:fieldext}. 

One of the main parts of Langlands conjectures is the \emphindex{principle of functoriality}. To state it, let  $G$ and $G'$ be reductive groups over $\mathbb{Q}$. The principle of functoriality asserts that whenever we have a group homomorphism between the associated $L$-groups
\beq
\Psi \, :\, {}^LG \, \longrightarrow {}^LG',
\label{Lgrouphom}
\eeq
 there should be a close relation between the associated automorphic forms on $G(\mathbb{Q})\backslash G(\mathbb{A})$ and $G(\mathbb{Q})\backslash G(\mathbb{A})$. What does ``close relationship'' mean? Suppose $ \pi$ is an automorphic representation of $G$ associated with a Satake class $[A_\pi]$ in the dual group ${}^LG$. Functoriality implies that there exists an automorphic representation $\pi'$  of $G'$ with Satake class $[A_{\pi'}]\subset {}^LG'$, such that  
 \beq
 [A_{\pi'}]\cong [\Psi(A_{\pi})].
 \eeq
 It turns out that this has far-reaching consequences even for the case when the first group is taken to be trivial. Suppose for example that $G=\{1\}$ and $G'=GL(n)$. In this situation the dual group of $G$ is simply the Galois group $\text{Gal}(\bar{\mathbb{Q}}/\mathbb{Q})$, where the extension is $\fieldext=\bar{\mathbb{Q}}$, the algebraic closure of $\rats$. (Recall that algebraic closure $\bar{\field}$ of a number field $\field$ is obtaining by adjoining to $\field$ all roots of all polynomials over $\field$. This not a finite extension and so generalises the discussion above.) The $L$-dual group of  $G'$ is the direct product ${}^LG'=GL(n, \mathbb{C})\ltimes \text{Gal}(\bar{\mathbb{Q}}/\mathbb{Q})$. The map $\Psi$ then yields a homomorphism 
 \beq
 \Psi\, :\, \text{Gal}(\bar{\mathbb{Q}}/\mathbb{Q}) \longrightarrow GL(n, \mathbb{C}).
 \label{Galhom}
 \eeq
 This has the remarkable consequence that to each  automorphic representation $\pi$ of $GL(n, \mathbb{A})$ there should exist an associated $n$-dimensional  representation $R$ of $\text{Gal}(\bar{\mathbb{Q}}/\mathbb{Q})$ such that 
 \beq
 L(s, \pi)=L_A(s, R), 
 \label{Lfunctionequality}
 \eeq
 where the object on the left is the standard $L$-function of $\pi$ discussed in section \ref{sec_autoLfunction} (i.e. corresponding to $\rho$ in section \ref{sec_autoLfunction} being the fundamental representation of $GL(n, \mathbb{C})$) and the object on the right is the so-called \emphindex[L-function@$L$-function!Artin]{Artin $L$-function} of the Galois representation $R$. We shall not go into the details of Artin $L$-functions but rather refer to \cite{LanglandsLfunct} for a nice discussion of the two sides of (\ref{Lfunctionequality}) and also to remark~\ref{rmk:converse}. 
 
There are numerous sources which give overviews of  various aspects of the Langlands program; we would like to especially mention \cite{GelbartLanglands,LanglandsLfunct,CasselmanLGroup,KnappLanglandsProgram,ArthurFunctoriality}. See also the two papers by Knapp \cite{KnappPrereq,KnappFirstSteps} which summarises the key references in the field.

\subsection{The Langlands program and physics}
As we have indicated at several occasions in this treatise, automorphic forms occur in abundance in string theory (see \cite{Green:1997tv,Green:1997tn,Green:1998by,Kiritsis:1997em,Obers:1999um,Pioline:2009qt,Bao:2009fg,Bao:2010cc,Green:2010wi,PerssonAuto} for a sample). Despite this fact the physical role of the classical Langlands program remains unclear. We know that automorphic representations play a role in understanding BPS-states and instantons in string theory, but we have no clue as to what is the physical interpretation of the dual side, involving representations of the Galois group. It would be very interesting to find out whether such an interpretation exists. Given that automorphic $L$-functions lie at the heart of the Langlands program, a very natural question, posed by Moore in \cite{MooreFuture}, is the following:

\vspace{.3cm} 

\noindent {\bf Open question:} {\it Is there a natural role for $L$-functions in BPS-state counting problems?} 

\vspace{.1cm} 

For some speculations on this and related issues, see \cite{Moore:1998pn,Miller:1999ag}, and for a conjectured connection between BPS-states in string theory and Galois representations, see \cite{Walcher:2012zk}. Interesting ideas on the connection between the geometry of Calabi--Yau manifolds and modular forms constructed from string world-sheet conformal field theories can be found in~\cite{Schimmrigk:2008mp}.

\subsection{The geometric version}
We should also mention that there exists a version of the Langlands program which does not have its roots in number theory, but rather in the geometry of Riemann surfaces. This is commonly referred to as the \emphindex[Langlands program!geometric]{geometric Langlands program} and was proposed by Beilinson, Deligne, Drinfeld and Laumon (for a nice survey see the lectures notes by Frenkel \cite{FrenkelLanglands}, and for a recent update see this note by Gaitsgory \cite{2016arXiv160609462G}). To each object in the original (or, ``classical'') Langlands program there exists geometric counterparts; for instance, the role of the Galois group is played by the fundamental group of the Riemann surface, while automorphic forms are replaced by certain ``automorphic sheaves'' on the moduli space of principal bundles on the Riemann surface. Remarkably, Kapustin and Witten have shown \cite{KapustinWitten} that the geometric Langlands program can be naturally understood in the context of quantum field theory (more precisely, a twisted version of $\mathcal{N}=4$ supersymmetric quantum field theory in four space-time dimension). In this context the  analogue of the ``Langlands duality'' (\ref{Galhom}) corresponds to a variant of (homological) mirror symmetry. \index{homological mirror symmetry}

\section{Whittaker functions, multiple Dirichlet series and statistical physics}

In this section, we discuss some issues related to a fascinating connection between Whittaker functions and statistical mechanics. Starting from a rewriting of the Casselman--Shalika formula, generalisations of Whittaker functions to metaplectic groups will be given. Their relation to Weyl group multiple Dirichlet series will be discussed and an alternative interpretation in terms of lattice models. This is an active area of research that has received a lot of momentum through the work of Brubaker, Bump, Chinta, Friedberg, Gunnells, Hoffstein and many others. We rely in our exposition mainly on the collection~\cite{MDSbook2012} and on~\cite{BBF} and refer the reader also to~\cite{BBFAMS} for an overview.

\subsection{Generalisations of the Weyl character formula}

The Casselman--Shalika vector for spherical Whittaker functions $W^\circ(\lambda,a)$ on the group $G(\rats_p)$ was discussed in detail in chapter~\ref{ch:Whittaker-Eisenstein} and given an interpretation in terms of characters $\mathrm{ch}_\Lambda$ of the Langlands dual group ${}^LG$ in equation~(\ref{CSWCF}). This formula can actually be inverted to give an alternative formula for highest weight characters of ${}^LG$ through
\begin{align}
\label{WCFCS}
\mathrm{ch}_\Lambda(a_\lambda) = \frac{W^\circ(\lambda,a_\Lambda)\delta^{-1/2}(a_\Lambda)}{\prod_{\alpha>0} (1-p^{-1} a_\lambda^{\alpha})}.
\end{align}
Here, $\lambda$ is a weight of the original group $G$ parametrising the principal series representation and $a_\lambda$ and $a_\Lambda$ distinguished elements of $A$ and ${}^LA$, respectively. These distinguished elements were defined in section~\ref{sec:CSLD}.

Formula (\ref{WCFCS}) resembles the standard Weyl  character formula~(\ref{charfn}), in particular the denominator. Independent of Whittaker functions, Tokuyama~\cite{Tokuyama} considered a one-parameter family of deformations of the Weyl character formula that can be written as\index{Weyl character formula!Tokuyama's deformation}
\begin{align}
\label{WeyldefToku}
\mathrm{ch}_\Lambda(a_\lambda) = \frac{\sum_{v\in\mathcal{B}_{\Lambda+\rho}} G(v,t) a_\lambda^{\wt(v)+\rho}}{\prod_{\alpha>0} (1+t a_\lambda^{\alpha})}
\end{align}
where $t\in \cx$. In this expression, all quantities refer to the Langlands dual group ${}^LG$. The sum here is over all $v$ in the \emphindex[crystal]{crystal} $\mathcal{B}_{\Lambda+\rho}$. The crystal $\mathcal{B}_{\Lambda+\rho}$ is a directed graph with vertices $v$ given by all the weights (with multiplicity) of the irreducible highest weight representation $V_{\Lambda+\rho}$ of ${}^LG$, where the $\rho$ shift is important, and the edges labelled by simple roots. The map $\wt: \mathcal{B}_{\Lambda+\rho}\to \mf{h}$ identifies the vertices with points in the weight lattice of ${}^LG$.  Crystals were introduced by Kashiwara~\cite{Kashiwara} in his study of the \emphindex[universal enveloping algebra!quantum deformation]{quantum deformed universal enveloping algebra} $U_q({}^L\mf{g})$ (closely related to \emphindex[quantum group]{quantum groups}~\cite{LusztigQG}) and possess a \emphindex[canonical basis]{canonical basis} in the sense of Kashiwara and Lusztig~\cite{Kashiwara,LusztigCanonical}. The operators corresponding to the edges are the simple step operators $f_i$ in the limit $q\to 0$. Kashiwara introduced also the crystal $\mathcal{B}_\infty$ that is modeled on the canonical (free) \emphindex{Verma module} of $U({}^L\mf{n}_-)$.

Following~\cite{BBCFH}, we will call the complex function $G(v,t)$ a \emphindex{Tokuyama function} and it is the main object of interest in this expression. In the original paper~\cite{Tokuyama}, the numerator was not written in terms of the crystal $\mathcal{B}_{\Lambda+\rho}$ but in terms of \emphindex{Gelfand--Tsetlin patterns}~\cite{GelfandTsetlin} with top row $\Lambda+\rho$ and the analysis restricted to the special linear group. We will comment later on the status for other groups.

Interesting special cases of the deformed character formula (\ref{WeyldefToku}) are
\begin{itemize}
\item {$t=-1$}. This is the value for the standard character formula~\eqref{charfn}. In this case the sum over the crystal collapses to a sum over the Weyl orbit of the shifted highest weight $\Lambda+\rho$. In other words, $G(v,-1)=0$ unless $\wt(v)=w(\Lambda+\rho)$ for some $w\in\Weyl$ and in that case $G(v,-1)=\eps(w\wlong)$.
\item  {$t=1$}. In this case one obtains a relation to the original formulation of Gelfand--Tsetlin patterns.
\item {$t=0$}. The denominator trivialises and Tokuyama use this case to recover a relation of Stanley's~\cite{Stanley} between Gelfand--Tsetlin patterns and `most singular' values of Hall--Littlewood polynomials. In the crystal formulation, the only contributing terms arise from the embedding of $\mathcal{B}_\Lambda \to \mathcal{B}_{\Lambda+\rho}$~\cite{MDSbook2012} and the sum then becomes the character in the form~\eqref{charactersum}.
\item {$t=-p^{-1}$}. This is the case relevant for the Casselman--Shalika formula and will be discussed in more detail below.
\end{itemize}
The first three cases were originally studied by Tokuyama~\cite{Tokuyama}. The last case in relation to Whittaker functions was first explored by Bump, Brubaker, Bump, Friedberg and Hoffstein~\cite{BBFH}, see also~\cite{HamelKing} for combinatorial aspects.

In the case $t=-p^{-1}$, we can compare~\eqref{WeyldefToku} and~\eqref{WCFCS} to deduce that we have an alternative description of Whittaker functions in terms of a sum over a crystal with a Tokuyama function $G(v,-p^{-1})$:\index{Whittaker function!as sum over crystal}
\begin{align}
\label{WhittToku}
W^\circ(\lambda,a_\Lambda) = \delta^{1/2}(a_\Lambda) \sum_{v\in \mathcal{B}_{\Lambda+\rho}} G(v,-p^{-1}) a_\lambda^{\wt(v)+\rho}.
\end{align}
In order to ease notation, we will suppress the $t$-value in the Tokuyama function in the sequel and will simply write $G(v)$ instead of $G(v,-p^{-1})$.

The identity~\eqref{WhittToku} in some sense defines the Tokuyama function $G(v)$ given the spherical Whittaker function. But it is desirable to have an independent description of the function $G(v)$. This was achieved in crystal form in~\cite{BBF} and can be given in terms of so-called decorated \emphindex[Berenstein--Zelevensky--Littelmann path]{Berenstein--Zelevinsky--Littelmann paths} (BZL paths) in the crystal~\cite{BerensteinZelevinsky,LittelmannC}. A BZL path of a vertex $v\in \mathcal{B}_{\Lambda+\rho}$ is given by first fixing a choice of a reduced expression of the longest Weyl word $\wlong$:
\begin{align}
\label{wlongred}
\wlong = w_{i_1}\cdots w_{i_\ell},
\end{align}
where $\ell=\ell(\wlong)$ is the length of the longest Weyl word and  $w_{i}$ is the $i$-th fundamental reflection. The BZL path $\BZL(v)$ of a crystal vertex $v\in\mathcal{B}_{\Lambda+\rho}$ is then obtained by following the simple lowering operators $f_i$ as far as possible through the crystal in the order given in the reduced expression of $\wlong$. Let $b_{1}$ be the largest integer such that $f_{i_1}^{b_{1}} v \neq 0$, that is, $b_{1}$ is the maximum number of steps in the direction of $f_{i_1}$ one can take in the crystal without leaving it. Starting from the point obtained in this way one then constructs $b_{2}$ as the maximum number of steps in the $f_{i_2}$ direction and so. This yields a sequence of non-negative integers 
\begin{align}
\label{BZLstring}
\BZL(v)=(b_1, b_2,\ldots, b_\ell)
\end{align}
and the endpoint of the crystal always corresponds to the `lowest weight' $v_-\in\mathcal{B}_{\Lambda+\rho}$ with $\wt(v_-)=\wlong(\Lambda+\rho)$. A vertex $v$ is uniquely characterised by its string $\BZL(v)$ (for a fixed choice of reduced expression~\eqref{wlongred}).

\begin{figure}
\centering
\input{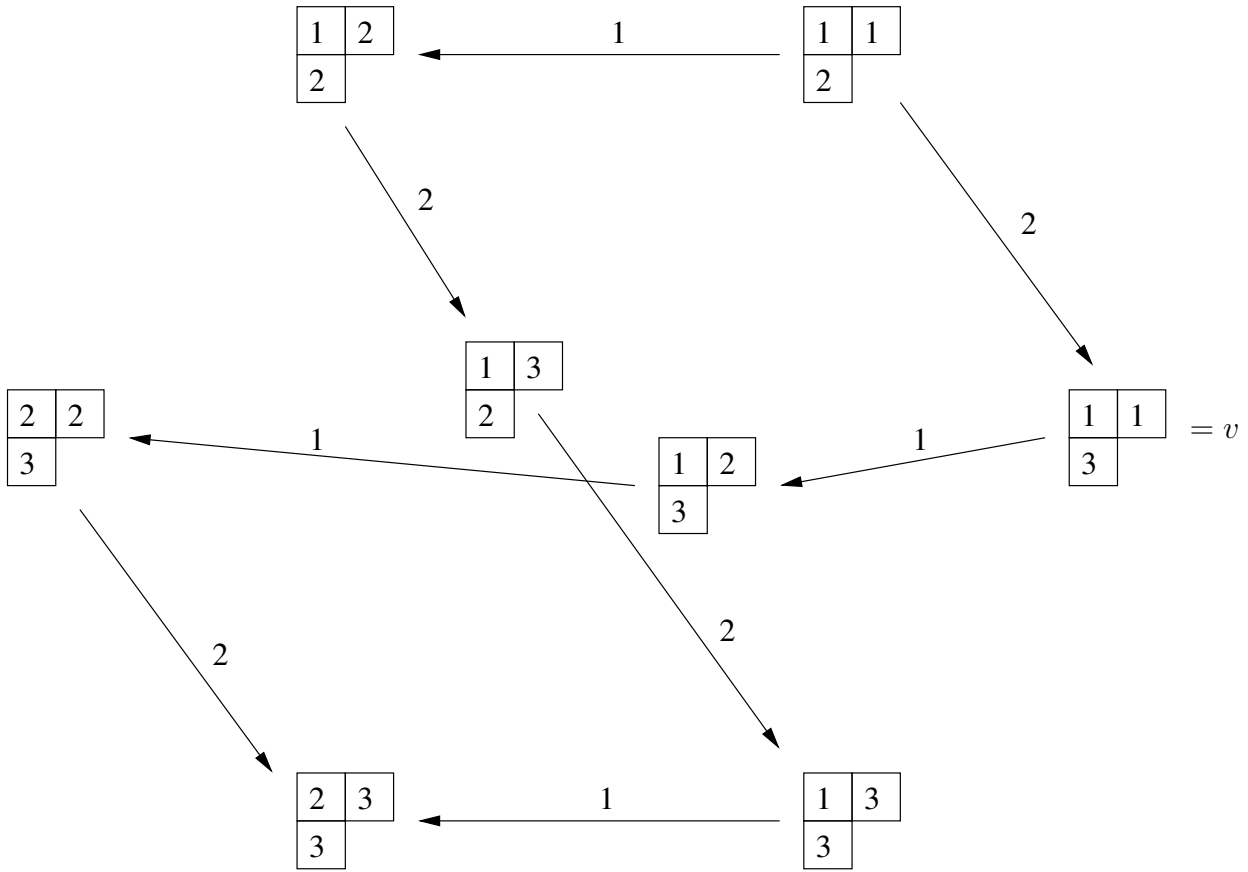}
\caption[Crystal $\mathcal{B}_{\Lambda+\rho}$]{\label{fig:CrystalBZL} \textit{The crystal $\mathcal{B}_{\Lambda+\rho}$ for $SL(3)$ and $\Lambda=0$. The $\rho$-shift turns this into the weight diagram of the adjoint representation and we label the different vertices of the crystal by filled Young tableaux. The arrows with numbers indicate the action of the operators $f_1$ and $f_2$. The two tableaux in the center correspond to the multiplicity two weight space of the adjoint representation associated with the two-dimensional Cartan subalgebra.  For the choice of vertex $v=\scalebox{0.75}{$\young(11,3)$}$, and the reduced expression $\wlong=w_1w_2w_1$, the BZL path is $\BZL(v)=(2,1,0)
=
\scalebox{0.75}{$
\begin{Bsmallmatrix}
    \raisebox{-2pt}{\circled{$0$}} & \young(1)\\
    \young(2) & 
\end{Bsmallmatrix}$}
$, where we have circled and boxed the entries according to the rules described in the text. }}
\end{figure}

For determining the function $G(v)$, the BZL string $(b_1,\ldots, b_\ell)$ needs to be decorated further.\index{Berenstein--Zelevensky--Littelmann path!decorated} In the case of $G=GL(r+1)$ of rank $r$ this is described in~\cite{BBF} for two choices of reduced $\wlong$ words. We will give here the version for 
\begin{align}
\wlong = w_1 w_2 w_1 w_3 w_2 w_1 \cdots w_{r} w_{r-1} \cdots w_2 w_1
\end{align}
and note that $\ell=\frac12r(r+1)$. 

\begin{remark}
There is another common choice of reduced expression $\wlong=w_r w_{r-1} w_rw_{r-1}w_{r-2}\cdots w_r\cdots w_2 w_1$ that is obtained by starting at the other end of the Dynkin diagram~\cite{BBFbook}. We will not use it here---it corresponds to $\Delta$-ice in a statistical mechanics interpretation whereas the choice here corresponds to $\Gamma$-ice.
\end{remark}
 The numbers in the BZL string $\BZL(v)$ of~\eqref{BZLstring} are then arranged in a triangular (Gelfand--Tsetlin-like) pattern according to
\begin{align}
\BZL(v) = \left\{\begin{array}{ccc}
\ldots &\ldots &\ldots\\
b_3 &b_2\\
b_1
\end{array}\right\},
\end{align}
such that the $i$-th column contains all numbers associated with the $w_i$ fundamental reflection. Littelmann proved that the numbers along a fixed row are weakly increasing~\cite{LittelmannC}. Entries in this tableaux now get circle or box decorations according to the following rules: $(i)$ if an entry $b_k$ is equal to its left neighbour (or equal to $0$ if it does not have one) it is circled; $(ii)$ if the crystal point $f_{i_{k-1}}^{b_{k-1}}\cdots f_{i_1}^{b_1} v$ does not have a neighbour in the $e_{i_k}$ direction, i.e. it sits on the boundary, then $b_k$ is boxed. An example of this description is given in figure~\ref{fig:CrystalBZL} for the case $r=3$. The boxing and circling rules can be given a geometrical interpretation in terms of the embedding of $\mathcal{B}_{\Lambda+\rho}$ into $\mathcal{B}_\infty$~\cite{BumpNakasuji}.

The decorated BZL string can then be used to define the Tokuyama function $G(v)$ via~\cite{MDSbook2012}
\begin{align}
\label{crystalrule}
G(v) = \prod_{b_k\in\BZL(v)} \left\{ \begin{array}{cl} 
1 & \textrm{if $b_k$ is circled but not boxed}\\
-p^{-1} & \textrm{if $b_k$ is boxed but not circled}\\
1-p^{-1} & \textrm{if $b_k$ is neither boxed nor circled}\\
0 &\textrm{if $b_k$ is boxed and circled}
\end{array}\right.
\end{align}
More complicated versions of this rule exist for other values of the Tokuyama deformation parameter $t$~\cite{BBF}. An equivalent description of $G(v,t)$ directly in terms of Gelfand--Tsetlin patterns was given in~\cite{BBFGauss,BBFbook,BBFH}.

Tokuyama's formula for the Tokuyama function $G(v)$ of~\eqref{crystalrule} gives a purely combinatorial description of the spherical Whittaker function $W^\circ(\lambda, a_\Lambda)$ for the case $G=GL(r+1,\rats_p)$.  One may wonder whether other choices of Tokuyama function $G(v)$ also correspond to objects related to automorphic forms. An affirmative answer to this was provided by Bump, Brubaker and Friedberg and we will discuss this next in a broader context~\cite{MDSbook2012}. But first we illustrate the crystal picture for the $p$-adic spherical Whittaker function of the $SL(2,\ads)$ Eisenstein series.

\begin{example}[Crystal description of $SL(2,\rats_p)$ Whittaker function]
We consider the case $G=SL(2,\rats_p)$ and verify formulas~\eqref{WhittToku} and~\eqref{crystalrule}. The Whittaker function $W^\circ(\lambda,a_\Lambda)$ is (cf.~\eqref{CS2Lang})
\begin{align}
\label{WSL2Tok}
W^\circ(\lambda,a_\Lambda) = (1-p^{-2s}) \frac{p^{s N-N} - p^{-2s+1-sN}}{1-p^{-2s+1}}.
\end{align}
To work out the crystal sum we fix the longest Weyl word as $\wlong=w_1$ and the highest weight is $\Lambda+\rho=(N+1)\rho$ where we recall that everything refers to  Langlands dual group of $SL(2,\rats_p)$. Then the crystal $\mathcal{B}_{\Lambda+\rho}$ consists of the vertices $v\in\{ (N+1)\rho, (N-1)\rho, \ldots, -(N+1)\rho \}$ that we label $v_k=(N+1-2k)\rho$ for $k=0,\ldots,N+1$. The highest weight representation $V_{\Lambda+\rho}$ is of dimension $N+2$. The BZL path of a vertex $v_k$ is 
\begin{align}
\BZL(v_k) = (N+1-k)
\end{align}
and its single entry is circled for $k=N+1$ and boxed for $k=0$, otherwise it is undecorated. Therefore
\begin{align}
G(v_k) = \left\{\begin{array}{cl}
-p^{-1} & \textrm{for $k=0$}\\
1 & \textrm{for $k=N+1$}\\
1- p^{-1} & \textrm{otherwise}
\end{array}\right.
\end{align}
The right-hand side of equation~\eqref{WhittToku} becomes therefore ($a_\lambda^\rho= p^{-(2s-1)/2}$)
\begin{align}
p^{-N/2} \sum_{k=0}^{N+1} G(v_k) a_\lambda^{(N+2-2k)\rho} 
&= p^{-sN-2s+1} \left(-p^{-1} + (1-p^{-1})\sum_{k=1}^N p^{k(2s-1)} +p^{(N+1)(2s-1)} \right)\nn\\
&= p^{-sN-2s+1} \left(-1 + (1-p^{-1})\frac{1-p^{(N+1)(2s-1)}}{1-p^{2s-1}} +p^{(N+1)(2s-1)} \right)\nn\\
&= p^{-sN-2s} \frac{-1+p^{2s}+p^{(N+1)(2s-1)}-p^{(N+1)(2s-1)+2s}}{1-p^{2s-1}}\nn\\
&= p^{-(N+2)s}(1-p^{2s}) \frac{1-p^{(N+1)(2s-1)}}{p^{2s-1}-1}\nn\\
&= p^{-sN}(1-p^{-2s}) \frac{p^{N(2s-1)}-p^{-2s+1}}{1-p^{-2s+1}},
\end{align}
which agrees with~\eqref{WSL2Tok}.
\end{example}

\subsection{Weyl group multiple Dirichlet series}
\label{sec:WGMDS}

In sections~\ref{app_adelichecke} and~\ref{sec_autoLfunction} we introduced Dirichlet series and automorphic $L$-functions. Both are meromorphic functions of a single complex variable $s$, satisfy functional equations for $s\leftrightarrow 1-s$ and have an Euler product form. They correspond to multiplicative sequences $a_n$ of numbers, in the simplest case of an $SL(2,\ints)$ cuspidal Hecke eigenform $f$ these are just the Fourier coefficients of $f$, cf.~(\ref{HeckeL}), so there is a close connection between Fourier expansions of automorphic forms and Dirichlet series, see sections~\ref{sec_autoLfunction} and~\ref{sec:LSmethod} for more details. 

It is natural to wonder whether these concepts can be generalised to functions of \emph{several} complex variables $s_1,\ldots, s_r$. This is a non-trivial problem and it turns out that multiplicativity of the coefficients cannot be maintained, see~\cite{MDSbook2012} that also discusses the history of the subject. One way of constructing such multiple Dirichlet series is as so-called \emphindex[Dirichlet series!multiple]{Weyl group multiple Dirichlet series}~\cite{BBCFH}. 

To introduce them we again restrict to $G=GL(r+1)$ of rank $r$ and introduce the following additional definitions~\cite{BBF2}. Let $\field$ be a number field that contains the group $\mu_{2n}$ of $2n$-th roots of unity. Let $S$ be a finite set of places of $\field$ such that $S$ includes all archimedean places (e.g. $p=\infty$) and all divisors of $n$. We denote by $\field_p$ the completion of $\field$ at a place $p$ and by $\mf{o}_p$ the corresponding integers for $p$ non-archimedean. The ring of $S$-integers $\mf{o}_S$ in $\field$ are those $x\in\field$ whose component $x_p$ is in $\mf{o}_p$ for all $p\notin S$. We can enlarge, if necessary, $S$ such that the $S$-integers $\mf{o}_S$ are a principal ideal domain. We denote $\field_S = \prod_{p\in S} \field_p$. 

A general form for a multiple Dirichlet series is then given by
\begin{align}
\label{eq:MDS}
Z_\Psi(\mathbf{s},\mathbf{m})
= \sum_{\mathrm{ideals }(C_i)} \Psi(C_1,\ldots,C_r) H(C_1,\ldots,C_r;m_1,\ldots,m_r) |C_1|^{-2s_1} \cdots |C_r|^{-2s_r}
\end{align}
for $\mathbf{s}=(s_1,\ldots, s_r)\in\cx^r$ and $\mathbf{m}=(m_1,\ldots,m_r)\in (\mf{o}_S)^r$. Here, $\Psi : (\field_S^\times)^r \to \cx$ and $H: (\field_S^\times)^r \times (\mf{o}_S)^r\to \cx$ are functions with multiplicativity properties that ensure that the sum over ideals in the principal ideal domain $\mf{o}_S$ is well-defined. These properties rely on the properties of the \emphindex[Hilbert symbol]{$n$-order Hilbert symbol} and on $n$-th order reciprocity. We will make no further use of the precise conditions here and refer the reader to~\cite{BBF2} for the details. We note, however, that the conditions make $\Psi$ a member of a finite-dimensional vector space $\mathcal{M}$ and that this space carries an action of the Weyl group $\Weyl$~\cite{BBF2}. 

The function $H$ satisfies an additional multiplicative property, called \emphindex[multiplicativity!twisted]{twisted multiplicativity}, that ensures that it is completely determined by its values on prime powers $H(p^{k_1},\ldots,p^{k_r}; p^{l_1},\ldots, p^{l_r})$. The parameters $\mathbf{m}$ appearing in~\eqref{eq:MDS} are called the twisting parameters and they enter crucially in the twisted multiplicativity relation. The problem of finding interesting multiple Dirichlet series is then reduced to specifying the $H(p^{k_1},\ldots,p^{k_r}; p^{l_1},\ldots, p^{l_r})$.  One requirement one would naturally impose on them is that, when viewed as a function of one $s_i$ alone, one obtains sums of single Dirichlet functions with standard functional relations. (Since we are working over a field that contains $\mu_{2n}$ these will actually be so-called \emphindex[Dirichlet series!Kubota]{Kubota Dirichlet series}~\cite{Kubota} and we refer again to~\cite{BBF2} for the details.) 

This requirement will ensure that the multiple Dirichlet series $Z_\Psi(\mathbf{s},\mathbf{m})$ will satisfy a functional relation under the fundamental Weyl reflection $w_i$ of the form
\begin{align}
Z_\Psi(w_i\mathbf{s},\mathbf{m}) = Z_{\Psi'}(\mathbf{s},\mathbf{m}),
\end{align}
where $\Psi'$ is some other element of the finite-dimensional space $\mathcal{M}$. By choosing a suitable normalisation
\begin{align}
Z_\Psi^\compl(\mathbf{s},\mathbf{m}) =\left( \prod_{\alpha>0} \zeta_\alpha(\mathbf{s}) G_\alpha(\mathbf{s}) \right) Z_\Psi(\mathbf{s},\mathbf{m})
\end{align}
in terms of factors of the Dedekind zeta function of $\field$ and appropriate $\Gamma$-function factors evaluated at places parametrised by the positive roots $\alpha$ and $\mathbf{s}$ one can bring this functional relation into the nicer form
\begin{align}
Z^\compl_{w\Psi}(w\mathbf{s},\mathbf{m}) = Z^\compl_{\Psi}(\mathbf{s},\mathbf{m})
\end{align}
by using the action of $\Weyl$ on $\mathcal{M}$. A functional equation of this type allows meromorphic continuation of the multiple Dirichlet series from the domain of convergence of (\ref{eq:MDS}) to $\cx^r$ by means of a variant of Bochner's theorem about complex functions in tube domains~\cite{BBCFH}.

Finding $H(p^{k_1},\ldots,p^{k_r}; p^{l_1},\ldots, p^{l_r})$ that satisfy this requirement is a non-trivial combinatorial problem. Essential information on the multiple Dirichlet series is contained in the so-called \emphindex[Dirichlet series!p-part@$p$-part of]{$p$-part} of $Z_\Psi$ that is defined by suppressing $\Psi$:
\begin{align}
\label{eq:ppart}
\sum_{k_i=0}^\infty H(p^{k_1},\ldots,p^{k_r}; p^{l_1},\ldots, p^{l_r}) |p|^{-2k_1s_1-\ldots -2k_r s_r}.
\end{align} 
The $p$-part depends on the $s_i$ and on the twisting parameters that are now given in terms of the integers $l_i$.
We would like to interpret the $p$-part as an expression evaluated on a crystal of ${}^LG$ evaluated at a special point $a_\lambda$ as in~\eqref{WhittToku}. To this end, we consider the case when the non-negative integers $(k_1,\ldots,k_r)$ correspond to a vector $\kappa$ linking a (strongly dominant) highest weight $\Lambda+\rho$ to one of its Weyl images, i.e.,
\begin{align}
\label{eq:stablepoints}
\kappa= \sum_i k_i \alpha_i = \Lambda+\rho - w(\Lambda+\rho)
\end{align}
for some $w\in\Weyl$, such that the second term in~\eqref{eq:ppart} is basically $\delta^{1/2}(a_\Lambda)a_\lambda^{\kappa+\rho}$. The highest weight $\Lambda$ here is determined by the integers $l_i$ through $\Lambda = \sum_i l_i \varpi_i$, where $\varpi_i$ are the fundamental weights of ${}^LG$. Suppose now that the twisting parameters $l_i$ and hence $\Lambda$ are fixed, then one can evaluate for $H(p^{k_1},\ldots,p^{k_r}; p^{l_1},\ldots, p^{l_r})$ for those $k_i$ for which~\eqref{eq:stablepoints} is satisfied as finite product of Gauss sums~\cite{BBFtwist}
\begin{align}
H(p^{k_1},\ldots,p^{k_r}; p^{l_1},\ldots, p^{l_r}) = \prod_{\alpha>0\atop w\alpha<0} g_2(p^{\langle\Lambda+\rho|\alpha\rangle-1},p^{\langle \Lambda+\rho|\alpha\rangle}),
\end{align}
where all elements now refer to the Langlands dual group and $g_2$ is a certain Gauss sum~\cite{BBFtwist}. The points $k_i$ for which one thus has a relatively simple formula for the value of $H$ on the Weyl orbit of the highest weight $\Lambda+\rho$. It is an important observation that the only other values of $k_i$ for which $H(p^{k_1},\ldots,p^{k_r}; p^{l_1},\ldots, p^{l_r})$ can be non-zero are the other points of the crystal $\mathcal{B}_{\Lambda+\rho}$. We are thus in a very similar situation to the discussion of the Weyl character formula and its generalisations above. For the strict Weyl character of the highest weight representation $V_{\Lambda+\rho}$ the crystal sum (\ref{WeyldefToku}) only had support on the Weyl images of $\Lambda+\rho$, but for the spherical Whittaker function one needed to consider also the other points of $\mathcal{B}_{\Lambda+\rho}$, in particular those in its interior. 

Determining $H(p^{k_1},\ldots,p^{k_r}; p^{l_1},\ldots, p^{l_r})$ at the other points of $\mathcal{B}_{\Lambda+\rho}$ is non-trivial and a number of approaches to this problem exist. 
\begin{itemize}
\item One approach is called the \emphindex[averaging method]{averaging method} of Chinta--Gunnells~\cite{ChintaGunnells1,ChintaGunnells2,MDSbook2012} that employs a deformed character constructed using the averaged Weyl group action on the field of rational functions in several variables. This approach works uniformly for any  type of root system. It has also been extended to the affine case in~\cite{LeeZhang}.

\item The approach by Bump, Brubaker and Friedberg~\cite{BBFtwist,BBF} starts from the just mentioned analogy with the crystal sum and finds rules for computing the Tokuyama function $G(v)$ on the crystal. These resemble the rules~\eqref{crystalrule} above but instead one gives different weights to the various parts of the BZL path. These weights are not simple powers of $p$ but instead involve $n$-th order Gauss sums. In this approach one has to treat each type of root system separately. For the various classical and some exceptional types we refer the reader to \cite{BBGGtypeB,BBFtypeC1,BBFtypeC2,FGGtypeG2,GGtypeD,FriedbergZhangBseries,FriedbergZhangTokuyamaB}.

\item Given a Tokuyama function $G(v)$ on the crystal one might wonder, in view of~\eqref{WhittToku}, whether the crystal sum can be interpreted as the Whittaker coefficient of some Eisenstein series. It turns out that in order for this to be true one needs to consider \emphindex[Eisenstein series!metaplectic]{metaplectic Eisenstein series}. These are functions defined over the group $G(\field)$ where $\field$ is the number field that contains the roots of unity $\mu_{2n}$. One can define Eisenstein series over $G(\field)$ and consider their Whittaker coefficients in the same way as we did for $G(\rats_p)$. However, it turns out that for $n>1$ one no longer has a multiplicity one theorem for Whittaker coefficients and that over a global field the Eulerian property is similarly no longer guaranteed. The multiplicative property of standard (non-metaplectic) Whittaker coefficients is replaced by the \emphindex[multiplicativity!twisted]{twisted multiplicativity} that we have already encountered above. The \emphindex[Whittaker function!metaplectic]{metaplectic Whittaker functions} are sources for the coefficients of multiple Dirichlet series as shown by Brubaker, Bump and Friedberg~\cite{BBF,MDSbook2012,McNamara1,McNamara2,BrubakerFriedbergWhittaker}. This approach works for all types of root system.

\item Finally, one can interpret the crystal sum as the partition function of a statistical mechanical mode, as done by Brubaker, Bump and Friedberg in~\cite{BBFYangBaxter,MDSbook2012}. In the case $n=1$, this opens up new tools from the theory of integrable systems, most notably the Yang--Baxter equation. This goes back to work of Kostant on the (quantum) Toda lattice and representation theory~\cite{KostantToda}.  In this approach one has to treat different types of root system differently. 
\end{itemize}
The equivalence of these different approaches has been shown in many cases and we refer to~\cite{MDSbook2012} for an overview.

\section{Quantum Whittaker functions}

\index{Whittaker function!$q$-deformed}

Following Givental's influential work on generalised mirror symmetry for flag varieties \cite{1996alg.geom.12001G}, Gerasimov--Lebedev--Oblezin \cite{GLO1,GLO2,GLO3,2007arXiv0705.2886G} have developed a theory of \emph{$q$-deformed} (or \emph{quantum}) \emph{Whittaker functions}, denoted $\Psi_q$, that have a variety of interesting applications in mathematics and theoretical physics. The quantum Whittaker functions are one-parameter generalisations that interpolate between the archimedean and non-archimedean Whittaker functions. They can be obtained as eigenmodes of the so-called $q$-deformed Toda Hamiltonian associated with $GL(n)$, and reduce in the limit $q\to 1$ to Givental's integral representation  of  the archimedean Whittaker function \cite{1996alg.geom.12001G}, while the limit $q\to 0$  recovers the Casselman--Shalika formula for the $p$-adic Whittaker function.  The quantum Whittaker functions play a role in quantum integrable systems \cite{Gerasimov:2007ys,Gerasimov:2014mfa}, topological string theory and mirror symmetry \cite{Gerasimov:2010xr,2013arXiv1309.5922G}. 

The quantum Whittaker function for $GL(\ell+1)$ can be written in the following explicit form \cite{2006math......8152G}: 
\begin{align}
\Psi_q(p)=\sum_{p_{k, i}\in \mathcal{P}^{(\ell+1)}}\prod_{k=1}^{\ell+1}q^{l_k(\sum_{i=1}^k p_{k,i}-\sum_{i=1}^{k-1} p_{k-1,i})}\frac{\prod_{k=2}^\ell \prod_{i=2}^{k-1}(p_{k, i+1}-p_{k,i})_q!}{\prod_{k=1}^\ell\prod_{i=1}^\ell (p_{k, i}-p_{k+1,i})_q!\, (p_{k+1, i+1}-p_{k,i})_q!}.
\label{qWhit}
\end{align}
Let us briefly explain the notation. The sum runs over a collection of integers $p_{i,j}\in \mathbb{Z}$, $i=1,\dots, \ell+1, \, j=1, \dots, i$, belonging to a set of Gelfand--Tsetlin patterns $\mathcal{P}^{(\ell+1)}$ associated  to an irreducible representation of $GL(\ell+1, \mathbb{C})$. The $q$-deformed factorial function is defined as $(n)_q!=(1-q)(1-q^2)\cdots (1-q^n)$.

In the limit $q\to 1$ the quantum Toda Hamiltonians reduce to the ordinary Toda system, whose solution is given by Givental's integral representation of the archimedean Whittaker function:
\begin{align}
\Psi_\mathbb{R}(x_1, \dots, x_{\ell+1})=\int_{\mathbb{R}^{\ell(\ell+1)/2}}\prod_{k=1}^\ell \prod_{i=1}^k dx_{k,i} e^{\mathcal{F}(x)},
\label{Givental}
\end{align}
where 
\begin{align}
\mathcal{F}(x)=\sum_{k=1}^{\ell+1} l_k\Big(\sum_{i=1}^k x_{k,i}-\sum_{i=1}^{k-1} x_{k-1,i}\Big)-\sum_{k=1}^\ell\sum_{i=1}^k \big(e^{x_{k, i}-x_{k+1,i}} +e^{x_{k+1,i+1}-x_{k,i}}\big).
\end{align}
The new variables are obtained by setting $x_j\equiv q^{p_{\ell+1,j}+j-1}$ while taking the limit $q\to 1$. 

In the other limit $q\to 0$, the quantum Whittaker function can be reduced to
\begin{align}
\Psi_q(p)\big|_{q\to 0}=\sum_{p_{k,i}\in \mathcal{P}^{(\ell+1)}}\prod_{k=1}^{\ell+1} z_k^{(\sum_{i=1}^kp_{k,i}-\sum_{i=1}^{k-1}p_{k-1},i)},
\label{padiclimit}
\end{align}
where we defined $z_i=q^{l_i}$, $i=1,\dots, \ell+1$. This has precisely the form of the characters of irreducible representations of $GL(\ell+1)$ in the Gelfand--Tsetlin bases. 

To get a better handle on the somewhat unwieldy expressions above we now consider a simple example.

\begin{example}
For $\ell=1$ we have simply $GL(2)$, for which the lattice is $(p_{2,1}, p_{2,2})\in \mathbb{Z}^2$. Then the general quantum Whittaker function becomes:
\begin{align}
\psi_q(p_{2,1}, p_{2,2})=\sum_{p_{2,1}\leq p_{1,1}\leq p_{2,2}}\frac{q^{l_1 p_{1,1}}q^{l_2(p_{2,1}+p_{2,2}-p_{1,1})}}{(p_{1,1}-p_{2,1})_q!(p_{2,2}-p_{1,1})_q!}
\end{align}
In the $q\to 1$ limit the integral representation \eqref{Givental} reduces to a Bessel integral, and up to normalisation one obtains the real spherical vector $\tilde{f}_\infty^\circ$ in \eqref{FourierSpherical}. Similarly, for $q\to 0$ the expression \eqref{padiclimit} can be readily identified with a standard character of $GL(2,\mathbb{C})$ and after appropriate identification of variables one obtains the $p$-adic spherical vector $\tilde{f}_p^\circ$ in \eqref{FourierSpherical}.
\end{example}

\begin{remark} 
The quantum Whittaker function $\Psi_q$  interpolates between archimedean and non-archimedean Whittaker functions which appear as local factors in the Fourier coefficient of an Eisenstein series. This fact suggests the following:

\vspace{.3cm}

\noindent {\bf Open question:} {\it Is there some kind of ``quantum automorphic form'' for which $\Psi_q$ appears naturally as a Fourier coefficient, or some generalisation thereof?}
\end{remark}

\section{Theta correspondences}
\label{sec:ThetaCorr}
Roughly, a theta correspondence (or theta lift) is a method for lifting (or transferring) an automorphic form $\varphi$ on a group $H$ to an automorphic form $\varphi'$ on another group $H'$, where $H\times H'$ is a subgroup of bigger group $G$. In the prototypical example, $H$ is an orthogonal group, $H'$  is a symplectic group and $G$ is a larger symplectic group. The correspondence is realised through an integral formula with kernel given by a classical theta function on $G$, restricted to the product $H\times H'$, thereby explaining the appearance of the word ``theta'' in the name of the correspondence. Theta correspondences have played an important role in establishing concrete examples of Langlands functoriality principle, but have also appeared ubiquitously in string theory. In this section we shall review some basic features of the theta correspondence and make various remarks on interesting connections with other parts of this treatise, both physical and mathematical. 

\subsection{The correspondence in a nutshell}
Let $G$  be an algebraic group defined over a number field $\field$. 

\begin{definition}[Dual pair]
A \emph{dual  pair} $(H, H^\prime)$  is a pair of subgroups  of $G$, such that $H\times H'\subset G$ and where $H$ (resp. $H'$) is the centraliser of $H'$ (resp. $H$) inside $G$ \cite{MR546602,MR1159101}. The dual pair is called reductive if, under the action of $H$ (and of $H'$), the defining representation of $G$ splits into the direct sum of two complementary invariant subspaces.
\end{definition}

Let $\pi$ be an automorphic representation of $G(\mathbb{A}_\field)$ of small Gelfand--Kirillov dimension. Let $\theta\in \pi$ be an automorphic form on $G(F)\backslash G(\mathbb{A}_\field)$. Restricting the representation $\pi$ to the product $H\times H'\subset G$ then yields an $H(\field)\times H'(\field)$-invariant function $\theta(h, h')$ on $H(\mathbb{A}_\field)\times H'(\mathbb{A}_\field)$  known as the \emph{theta kernel}. Starting from the theta kernel we can construct two different theta lifts:
\begin{definition}[Theta lift]
For $\varphi$ an automorphic form on $H(\field)\backslash H(\mathbb{A}_\field)$ we define the \emph{theta lift} from $H$ to $H'$  by the integral
\begin{align}
\Theta_\varphi(h')=\int_{H(\field)\backslash H(\mathbb{A}_\field)} \varphi(h) \theta(h,h') dh,
\end{align}
where $dh$ is the invariant measure on $H(\field)\backslash H(\mathbb{A}_\field)$. Similarly, for $\varphi'$ an automorphic form on $H'(\field)\backslash H'(\mathbb{A}_\field)$ we can define the theta lift in the other direction, i.e. from $H'$ to $H$, by the integral
\begin{align} 
\Theta_{\varphi'}(h)=\int_{H'(\field)\backslash H'(\mathbb{A}_\field)} \varphi'(h') \theta(h,h') dh'.
\end{align}
\end{definition}

\noindent Clearly, the function $\Theta_\varphi$ is defined on $H'(\ads_\field)$ and invariant under $H'(\field)$. The prototypical questions that one wishes to study in this context are the following:
\begin{itemize}
\item Does $\Theta_\varphi$ converge?
\item Is $\Theta_\varphi$ non-zero?
\item Is $\Theta_\varphi$ cuspidal?
\item If $\Theta_\varphi\neq0$, what kind of object is it?
\end{itemize}
If $\Theta_\varphi$ can be identified with an automorphic form on $H'(\mathbb{A}_\field)$ then the theta lift can be viewed as a map
\begin{align}
\Theta \, :\, \mathcal{A}(H(\field)\backslash H(\mathbb{A}_\field))\longrightarrow \mathcal{A}(H'(\field)\backslash H'(\mathbb{A}_\field)),
\end{align}
and therefore may in principle provide a concrete realisation of Langlands functoriality (see \cite{MR1145805,MR1109355} for some examples). 

Currently these problems cannot be addressed in general, but below we shall look at some examples where much is known.

\subsection{The Siegel--Weil formula}

We shall now consider a specific example of the theta correspondence. Let $V$ be an orthogonal vector space, i.e. a finite-dimensional vector space over $\field$ with non-degenerate inner product $(\, , \, )$. Let $V'$ be a symplectic vector space. Take $H=O(V)$, the orthogonal group of $V$ and $H'=Sp(V')$, the symplectic group over  $\field$. Then $W=V\otimes V'$  is naturally a symplectic vector space, and $(O(V), Sp(V'))$ forms a dual reductive pair inside $G=Sp(W)$.  The (metaplectic cover of) the adelic group $Sp(W, \mathbb{A}_\field)$ acts on the space of Schwartz--Bruhat functions $\mathcal{S}(V(\mathbb{A}_\field))$ via the Weil representation $\rho$:
\begin{align}
\rho : Sp(W, \mathbb{A}_\field)\longrightarrow \mathcal{S}(V(\mathbb{A}_\field)).
\end{align}
 Recall that this representation can be realized as follows
\begin{align}
\rho\left(\begin{pmatrix} A & 0 \\ 0 & (A^{-1})^T\end{pmatrix}\right) f(X) &=\sqrt{|\det(A)|}f(A^T X)
\nn \\
\rho\left(\begin{pmatrix} I & B \\ 0 & I\end{pmatrix}\right) f(X) &=\psi\left(\frac{X^TBX}{2}\right)f(X)
\nn \\
\rho\left(\begin{pmatrix} 0 & I \\ -I & 0\end{pmatrix}\right) f(X) &=\gamma \widehat{f(X)},
\end{align}
where $\psi : \field\backslash \mathbb{A}_\field \rightarrow \mathbb{C}^{\times}$ is a non-trivial character.  This action commutes with the natural action of $O(V)$ on $V$. 

Now, for any $\phi\in \mathcal{S}(V(\mathbb{A}_\field))$ define a distribution 
\begin{align}
\theta \quad :\quad  \mathcal{S}(V(\mathbb{A}_\field)) & \longrightarrow \mathbb{C}
\nn \\
\quad  \quad \phi & \longmapsto \sum_{x\in V(\field)} \phi(x).
\end{align}
From this distribution we define a function, the \emph{theta kernel}, by restriction from $Sp(W)$ to $Sp(V')\times O(V)$ as discussed above:
\begin{align}
\theta_\phi \quad :\quad Sp(V', \field)\backslash Sp(V', \mathbb{A}_\field)\times O(V)(\field)\backslash O(V)(\mathbb{A}_\field)\longrightarrow \mathbb{C}
\end{align}
by the formula 
\begin{align}
\theta_\phi(g,h)=\sum_{x\in V(\field)} \rho(g)\cdot \phi(h^{-1}x), \quad g\in Sp(V', \mathbb{A}_\field), \, h\in O(V)(\mathbb{A}_\field).
\end{align}
We want to study the associated theta lift 
\begin{align}
\Theta_\varphi(g, \phi)=\int_{O(V)(\field)\backslash O(V)(\mathbb{A}_\field)} \varphi(h)\theta_\phi(g,h) dh,
\label{thetaSiegel}
\end{align}
where $\varphi\in \mathcal{A}(O(V)(\field)\backslash O(V)(\mathbb{A}_\field))$. The Siegel--Weil formula is a special case of this theta lift where one takes $\varphi$ to be the identity function \cite{MR946349,MR961164}. To state it we need some additional data on Eisenstein series on symplectic groups.  

From now on we fix $Sp(V', \mathbb{A}_\field)$ to be $Sp(n, \mathbb{A}_F)$. Let $P=LU$ be the Siegel maximal parabolic subgroup of $Sp(n)$ with Levi factor 
\begin{align} 
L=\left\{l(a)= \begin{pmatrix} a & \\ & (a^T)^{-1}\end{pmatrix}\, \Big|\, a\in GL(n)\right\},
\end{align}
and unipotent radical 
\begin{align} 
U=\left\{u(b)= \begin{pmatrix} 1_n & b\\ & 1_n \end{pmatrix}\, \Big|\, b=b^T\right\}.
\end{align}
We have the (non-unique) decomposition $Sp(n, \mathbb{A}_\field)=U(\mathbb{A}_\field)L(\mathbb{A}_\field)K_{\mathbb{A}_F}$ and we write
\begin{align}
g=u\, l(a) \, k\in Sp(n),
\end{align}
with $a\in GL(n, \mathbb{A}_\field)$. We further introduce the quantity 
\begin{align}
|a(g)|=|\det(a)|.
\end{align}
Now define a standard section of the adelic principal series according to
\begin{align}
\Phi(g, s, \phi)=\rho(g)\cdot \phi(0) |a(g)|^{s-s_0(m,n)},
\end{align}
where 
\begin{align}
s_0(m,n)=\frac{m}{2}-\frac{n+1}{2}.
\end{align}
We now have all the data required to define the \emphindex{Siegel--Eisenstein series}
\begin{align}
E(g,s,\phi)=\sum_{\gamma\in P(\field)\backslash Sp(n, \field)} \Phi(\gamma g, s, \phi).
\end{align}
Provided that the function $\phi$ is $K_{\mathbb{A}_\field}$-finite this converges absolutely for $\text{Re}(s)>(n+1)/2$ and admits a meromorphic continuation and functional equation \cite{MR946349}.

Finally we can state the Siegel--Weil formula:

\begin{theorem}[Siegel--Weil formula \cite{MR0223373,MR946349,MR961164}]
Let $\Theta_\varphi(g, \phi)$ be the theta lift defined in (\ref{thetaSiegel}) and let $E(g,s,\phi)$  be the Siegel--Eisenstein series defined above. The we have:
\begin{align}
 \Theta_1(g, \phi)=\kappa E(g, s, \phi), 
 \end{align}
 where 
 \begin{align}
 \kappa=\left\{\begin{array}{cc} 1 & m>n+1 \\
 2 & m\leq n+1. \end{array}\right.
 \end{align}
 \end{theorem}
 
 The Siegel--Weil formula thus states that the theta lift of the trivial  function $\varphi=1$ on $O(V)(\field)\backslash O(V)(\mathbb{A}_\field)$ is a Siegel--Eisenstein on $Sp(n, \field)\backslash Sp(n, \mathbb{A}_\field)$. 
 
 \begin{example}
 Consider now the special case of the Siegel--Weil formula when $n=1$, corresponding to the dual reductive pair $SL(2)\times O(V)$. In this case the Siegel--Weil formula provides a transfer of automorphic forms on the orthogonal group $O(V)$ to $SL(2)$. The Siegel--Eisenstein series reduces  to the standard principal Eisenstein series on $SL(2,\mathbb{A}_\field)$ discussed at length in section~\ref{sec:nonspherical}:
 \begin{align}
 E({\sf f}_\lambda, g)=\sum_{\gamma\in B(\field)\backslash SL(2, \field)} {\sf f}_\lambda(\gamma g),
 \end{align} 
 where ${\sf f}_\lambda=\Pi_\nu {\sf f}_{\lambda, \nu}$ is a standard section of the adelic principal series which is non-spherical at the archimedean place:
 \begin{align}
 {\sf f}_{\lambda, \infty}\left(g_\infty \left(\begin{array}{cc}
\phantom{-}\cos\theta & \phantom{-}\sin \theta \\
-\sin \theta & \phantom{-}\cos \theta \\
\end{array} \right)\right)=e^{i w\theta} {\sf f}_{\lambda, \infty}(g_\infty).
\end{align}
The Siegel--Weil formula now reads explicitly
\begin{align}
\sum_{\gamma\in B(\field)\backslash SL(2, \field)} {\sf f}_\lambda(\gamma g)=\int_{O(V)(\field)\backslash O(V)(\mathbb{A}_\field)} \theta_\phi(g,h) dh.
\end{align}
 \end{example}
 
 \begin{remark}
 The Siegel-Weil formula was proven by Weil \cite{MR0223373} assuming certain convergence properties of the theta integral. In a series of works, Kudla and Rallis \cite{MR946349,MR961164,MR1289491} developed the theory of the \emphindex{regularised Siegel--Weil formula}, which holds when these convergence properties are not satisfied. Recently, the results of Kudla--Rallis were extended by Gan--Qiu--Takeda \cite{MR3279536}, who proved the regularized Siegel--Weil formula in full generality.
 \end{remark}
 
\subsection{Exceptional theta correspondences}
There exist versions of theta correspondences which do not correspond to orthogonal-symplectic dual pairs. In particular, one can consider the case when the ambient group $G$ is an exceptional Lie group \cite{MR1455531,MR1457344}. At the level of complex Lie groups we have the following examples of dual pairs
\begin{align}
PGL(3)\times G_2 & \subset E_6
\nn \\
 PGSp(6) \times G_2&\subset E_7
\nn \\
F_4 \times G_2 & \subset E_8.
\end{align}
These examples fit into the ``tower of theta correspondences'' constructed in \cite{MR1455531}. The tower consists of a sequence of dual pairs $(H, H')$, where the first member is always $H=G_2$ and the second varies. 
In this setting the Weil representation is replaced by the  minimal automorphic representation $\pi_{min}$ of $G$. This representation is unique, factorizable $\pi_{min}=\otimes_\nu \pi_{min, \nu}$, and has functional dimension
\begin{align}
E_6\quad & :\quad \text{GKdim}(\pi_{min})=11
\nn \\
E_7 \quad & :\quad\text{GKdim}(\pi_{min})=17
\nn \\
E_8 \quad & :\quad\text{GKdim}(\pi_{min})=29.
\end{align}

Let  $\theta_G(g), \, g\in G(\mathbb{A}_\field),$ be an automorphic form in the space of $\pi_{min}$. This may be viewed as the exceptional analogue of the Siegel theta series of $Sp(2n)$ considered in the previous section. By restriction to the dual pair $H\times H'\subset G$, this yields a function $\theta_G(h,h')$ on $H(\field)\backslash H(\mathbb{A}_\field)\times H'(\field)\backslash H'(\mathbb{A}_\field)$, i.e. the analogue of the theta kernel. Let $\pi$  be an automorphic representation of $H'(\mathbb{A}_\field)$. For any $\varphi\in \pi$ we then have the exceptional theta lift 
\begin{align} 
\Theta_\varphi(h)=\int_{H(\field)\backslash H(\mathbb{A}_\field)} \varphi(h)\theta_G(h,h')dh.
\end{align}
For $H'=G_2$, this yields a tower of maps
\begin{align}
\Theta_H^{G_2} \, :\, \mathcal{A}(H(\field)\backslash H(\mathbb{A}_\field))\longrightarrow \mathcal{A}(G_2(\field)\backslash G_2(\mathbb{A}_\field)),
\end{align}
and we denote the space of such theta lifts by $\Theta_G^{G_2}(\pi)$ as $f$ varies in $\pi$ and $\theta_G$ varies in $\pi_{min}$. The space $\Theta_G^{G_2}(\pi)$ is thus a 
subspace of the space of automorphic forms on $G_2(\mathbb{A}_\field)$.

Consider now the case when $G$ is the quaternionic real form of either $E_7$ or $E_8$. These  real forms are characterized by having real rank 4, and the restricted root system is of type $F_4$ (in the literature they are sometimes denoted by $E_{7(-5)}$ and $E_{8(-24)}$).  The maximal compact subgroup is of the form $SU(2)\times M$, where $M$ is the compact form of $D_6$ for $G=E_7$ and the compact form of $E_7$  for  $G=E_8$. The first member $H$ in the dual pair   corresponds to the automorphism group of the exceptional Jordan algebra of $3\times 3$ hermitian matrices (see \cite{MR1387237} for details) with coefficients in the quaternionic division algebra $\mathbb{H}$ (for $E_7$), or the octonionic division algebra $\mathbb{O}$ (for $E_8$). For both $E_7$ and $E_8$ the second member $H'$  in the dual pair is the split real form of $G_2$. 

Gan \cite{MR1767400} has constructed an automorphic realization of the minimal representation $\pi_{min}$ of $G(\mathbb{A}_\field)$, where $G(\field_\nu)$  is split for all finite places $\nu$, while the archimedean factor $G(\field_\infty)$ is the quaternionic real form of $E_7$ or $E_8$. This was then used  to study the theta lift from automorphic forms on $H(\field)\backslash H(\mathbb{A}_\field)$ to $G_2(\field)\backslash G_2(\mathbb{A}_\field)$. Let $P=LU$  be the Heisenberg parabolic subgroup of $G$ (corresponding to node 1 on the Dynkin diagram) with modulus character $\delta_P$. The Levi factor $L$  has derived group of type  $D_6$  for  $G=E_7$  and type  $E_7$  for  $G=E_8$. For $s\in \mathbb{C}$,  let ${\sf f}_s=\prod_{\nu}{\sf f}_{s,\nu}$ be a standard section of the induced representation $\text{Ind}_{P(\mathbb{A}_\field)}^{G(\mathbb{A}_\field)}\delta_P^{s}$ and construct the corresponding Eisenstein series
\begin{align}
E(g, {\sf f}_s)=\sum_{\gamma\in P(\field)\backslash G(\field)} {\sf f}_s(\gamma g).
\end{align}

\noindent Gan then proved the following: 

\begin{theorem}[Automorphic realisation of $\pi_{min}$ for quaternionic real forms\cite{MR1767400,MR1801660}]
Let $G(\ads_\field)$ be the adelic group associated with the quaternionic real form of $E_7$ or $E_8$ as above. For any standard section ${\sf f}_s\in \text{Ind}_{P(\mathbb{A}_\field)}^{\, G(\mathbb{A}_\field)}\delta_P^{s}$ the Eisenstein series $E(g, {\sf f}_s)$ has at most a simple pole at $s=s_0$, with 
\begin{align}
s_0=\left\{ \begin{array}{cc} 3/17 & \text{for}\quad G=E_7\\ 
5/29 & \text{for}\quad G=E_8. \end{array} \right.
\end{align}
 At this value of $s$ the global minimal representations $\pi_{min}$ of the quaternionic real forms of the exceptional groups arise as a square-integrable subspace
\begin{align}
\pi_{min} \subset \text{Ind}_{P(\mathbb{A}_\field)}^{\, G(\mathbb{A}_\field)}\delta_P^{s}\Big|_{s=s_0}, 
\end{align}
and its automorphic realization is spanned by the residues of the Eisenstein series $E(g, {\sf f}_s)$ at $s=s_0$. 
\end{theorem}

\begin{remark}
The above theorem may be viewed as the quaternionic version of an earlier result by Ginzburg--Rallis--Soudry \cite{GRS} who constructed an automorphic realisation of the minimal representation of $E_8(\mathbb{A}_\field)$ when the archimedean component is the split real form of $E_8$.
\end{remark}

\begin{remark} In contrast to the split case considered in \cite{GRS}, the minimal representation $\pi_{min}=\otimes_\nu \pi_{min, \nu}$ of the  $G(\mathbb{A}_\field)$ constructed by Gan may be non-spherical for some finite places $\nu$.
\end{remark}

We shall now proceed to discuss theta correspondences in the context of quaternionic exceptional groups. Let us denote by $\theta_{G}(g)$ the (residue of the) Eisenstein series $E(g, {\sf f}_s)$ at $s=s_0$. By restriction of $g\in G$ to the dual pair $H\times G_2$, through $g=hh'$ we obtain the theta kernel 
\begin{align}
\theta_{G}(h, h')=E(hh', {\sf f}_{s_0}).
\end{align}
For $\pi$ an automorphic representation of $H(\mathbb{A}_\field)$ and $\varphi\in \pi$ we then write the theta lift as
\begin{align}
\Theta_\varphi(h', {\sf f}_{s_0})=\int_{H(\field)\backslash H(\mathbb{A}_\field)} \varphi(h) E(hh', {\sf f}_{s_0})dh.
\end{align}
Gan proved that $\Theta_\varphi$ is non-zero and cuspidal under certain conditions on $\pi$. In particular:
\begin{itemize}
\item $\Theta_\varphi$ is cuspidal if and only if 
\begin{align}
\int_{C(\field)\backslash C(\mathbb{A}_\field)} \varphi(c) dc=0,
\end{align}
where $C$ is the stabilizer of a particular character $\psi$ of the unipotent radical of a maximal parabolic of $G$ (see  \cite{MR1767400} for the precise definition).
\item $\Theta_\varphi$ is zero if and only if
\begin{align}
\int_{Z(\field)\backslash Z(\mathbb{A}_\field)} \varphi(zh)dz=0,
\end{align}
where $Z=[U,U]$ is the centre of the Heisenberg unipotent radical $U$ of $P$.
\end{itemize}

For the case of $G=E_8$, Gan has also established a Siegel--Weil formula for the dual pair $H\times G_2\subset E_8$ \cite{MR1801660}. In this case, $H$ is a real form of $F_4$. The Siegel-Weil formula arises the theta lift of the trivial function $\varphi=1$ on $F_4(\mathbb{A}_\field)$  to $G_2(\mathbb{A}_\field)$:
\begin{align}
\Theta_1(h', {\sf f}_{s_0})=\int_{F_4(\field)\backslash F_4(\mathbb{A}_\field)} E(hh', {\sf f}_{s_0})dh.
\label{trivialthetaG2}
\end{align}
In this situation the theta lift does not give rise to cusp forms, 
but analogously to the classical case of orthogonal-symplectic pairs, it yields certain Eisenstein series. To state it we first need some information on Eisenstein series on $G_2$. Let $P_2=L_2U_2$ be the maximal parabolic subgroup of $G_2$ with unipotent radical $U_2$ a 5-dimensional Heisenberg group, and Levi subgroup $L_2=SL(2)\times GL(1)$. Denote the modulus character of $P_2$ by $\delta_{P_2}$ and for $s\in \mathbb{C}$  consider the global degenerate principal series
\begin{align}
I(s)=\text{Ind}_{P_2(\mathbb{A}_\field)}^{G_2(\mathbb{A}_\field)}\delta_{P_2}^{s}=\bigotimes_{\nu}I_\nu(s).
\end{align}
For any standard section ${\sf f}_s\in I(s)$ we have the Eisenstein series
\begin{align}
E(g, {\sf f}_s)=\sum_{\gamma\in P_2(\field)\backslash G_2(\field)}{\sf f}_s(\gamma g).
\end{align}
We will not be needing the Eisenstein series for any such section, but rather a special case thereof, obtained by restriction from the minimal representation $\pi_{min}$ of $E_8$. Since $\pi_{min}$ is a submodule of the degenerate principal series $\text{Ind}_{P(\mathbb{A}_\field)}^{ E_8(\mathbb{A}_\field)}\delta_P^{5/29}$ one can restrict it to a subspace of $\text{Ind}_{P_2(\mathbb{A}_\field)}^{G_2(\mathbb{A}_\field)}\delta_{P_2}^{5/2}$. More precisely, one has a restriction map \cite{MR1801660}:
\begin{align}
\text{Res} \, \,  :\,  \pi_{min} & \longrightarrow \text{Ind}_{P_2(\mathbb{A}_\field)}^{G_2(\mathbb{A}_\field)}\delta_{P_2}^{5/2}
\nn \\
\,  f& \longmapsto \text{Res}(f).
\end{align}
A standard section obtained by restriction in this way is called a \emphindex{Siegel--Weil section} \cite{MR946349}. These are the relevant ones for the Siegel--Weil formula. To understand the structure of these sections we need the following result:

\begin{proposition}[Restriction of the minimal representation of $E_8$ \cite{MR1387237}]
Let $\pi_{min}=\pi_{min, \infty} \otimes \hat{\bigotimes}_\nu \pi_{min, \nu}$ be the minimal representation of $E_8(\mathbb{A}_\field)$, with archimedean component being the quaternionic real form of $E_8$. Then the restriction of $\pi_{min}$ to its $F_4$-invariant subspace  is a quaternionic discrete series representation $\pi$ of the split real form $G_2(\mathbb{R})$, in the sense of Gross--Wallach \cite{MR1327538,MR1421947}. 
\end{proposition}

\noindent Given the above result one has the following consequence:

\begin{corollary}
A standard section ${\sf f}={\sf f}_\infty \times \hat{\prod}_{\nu}{\sf f}_\nu\in \text{Ind}_{P_2(\mathbb{A}_\field)}^{\, G_2(\mathbb{A}_\field)}\delta_{P_2}^{5/2}$ is a Siegel--Weil section if and only if ${\sf f}_\infty \in \pi$.
\end{corollary}

A priori the theta lift $\Theta_1(h)$ in \eqref{trivialthetaG2} can have a cuspidal component (in the discrete spectrum) and an Eisenstein component (in the continuous spectrum):
\begin{align}
 \Theta_1(h, {\sf f})=(\Theta_1(h, {\sf f}))_{\text{Eisenstein}} + (\Theta_1(h, {\sf f}))_{\text{cusp}}.
 \end{align}
 The main result of \cite{MR1801660} is then:
 
 \begin{theorem}[Exceptional Siegel--Weil formula for $G_2$ \cite{MR1801660}]
 For a section ${\sf f} \in \pi_{min}$ of $E_8(\mathbb{A}_F)$, one has 
 \begin{align}
 (\Theta_1(h, {\sf f}))_{\text{cusp}} & = 0
 \nn \\
(\Theta_1(h, {\sf f}))_{\text{Eisenstein}} & = E\big(h, \text{Res}({\sf f})\big).
 \end{align}
 \end{theorem}

\noindent That is, the theorem implies the following explicit form of the Siegel--Weil formula
\begin{align}
\int_{F_4(\field)\backslash F_4(\mathbb{A}_\field)} E(hh', {\sf f})dh=\sum_{\gamma\in P_2(\field)\backslash G_2(\field)} \text{Res}({\sf f})(\gamma h').
\end{align}

\begin{remark} The above result is expected to have applications in string theory. In fact, the Eisenstein series $E\big(h, \text{Res}({\sf f})\big)$ is a natural candidate for the type of automorphic form on $G_2$ that is relevant for conjecture~\ref{QDS} above.
\end{remark}

\section{Moonshine}
\label{sec:Moonshine}

In mathematics and theoretical physics  the term `moonshine' refers to webs of connections between different areas, such as number theory, group theory, and string theory. The most famous example is known as `monstrous moonshine', and its impact can hardly be overstated; it led to numerous developments in mathematics which crucially relied on results in string theory. In many respects, moonshine is similar to the Langlands program in the sense that it connects many a priori distinct areas under a single framework. Moreover, automorphic forms play a pivotal role in moonshine and therefore clearly deserves a place in this treatise. However, as we will see, the automorphic forms that appear in moonshine are of a different kind and the relation to automorphic representations is not clear. We shall comment on this in a little more detail at the end of this section. 

\begin{remark}
In order to stay closer to the original literature and to avoid confusion when talking about Jacobi forms, we change some of our standard notation in this section. The typical element of the upper half plane $\UHP$ will be called $\tau$ (not $z$!) and typically the dependence on it will be holomorphic. This implies that $SL(2,\ints)$-invariant objects will admit a \emphindex[q-expansion@$q$-expansion]{$q$-expansion} with $q=e^{2\pi i \tau}$.
\end{remark}

\begin{figure}
\begin{center}
    \begin{tikzpicture}[thick]
        \fill[lightgray!30] (0.38, -3.29) coordinate (A) -- (1.53, -3.62) coordinate  (B) -- (1.88, -0.25) coordinate (D);
        \draw (A) -- (B) -- (3.81, -3.27) coordinate (C) -- (D) -- cycle;
        \draw (B) -- (D);
        \node[below, left=2.5mm] at (A) {modular forms};
        \node[below=2.5mm, align=center] at (B) {infinite-dimensional \\[-0.1em] Lie algebras};
        \node[below, right=2.5mm,align=center] at (C) {string theory \\[-0.1em] VOAs};
        \node[above=2.5mm] at (D) {finite groups};
    \end{tikzpicture}
\end{center}
\caption{Pictorial overview of moonshine.}
\end{figure}
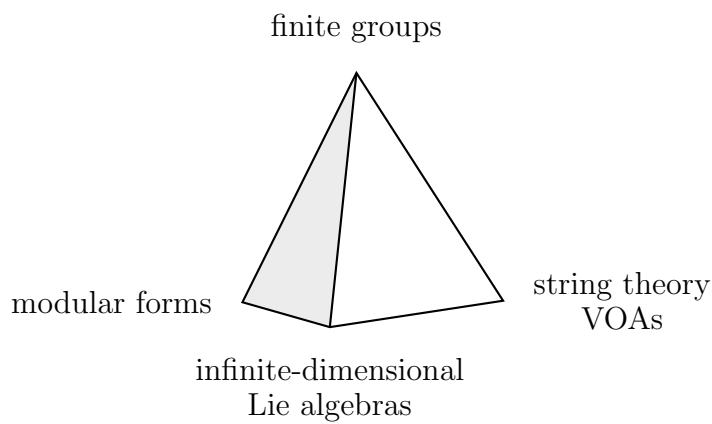

\subsection{Monstrous moonshine}
Monstrous moonshine pertains to a relation between the largest of the 26 finite sporadic simple groups known as the \emphindex[monster group]{monster} $\mathbb{M}$, and a certain  class of modular functions. The starting point of monstrous moonshine was McKay and Thompson's observation that the first few Fourier coefficients in the $q$-expansion of the unique (up to the constant term) holomorphic $SL(2,\mathbb{Z})$-invariant $J(\tau)$-function 
\beq
J(\tau)=\sum_{n=-1}^{\infty} c(n) q^n=q^{-1} + 196884q + 21493760 q^2 + \cdots
\label{Jfunction}
\eeq 
correspond to sums of dimensions of irreducible representations of $\mathbb{M}$:
\beq
196884=1 + 196883, \qquad 21493760=1 + 196883+12196876.
\eeq
 This suggests the existence of an infinite-dimensional graded representation $V=\bigoplus_{n=-1}^{\infty}V_n$ of $\mathbb{M}$, such that each graded piece $V_n$ decomposes into direct sums of irreducible representations of $\mathbb{M}$, and the Fourier coefficients can be reinterpreted as $c(n)=\text{dim}\, V_n$; in other words, the $J$-function is the \emph{graded dimension}Ê of $V$.

One can then consider a more general class of objects where the $c(n)$ are replaced by characters $\text{Tr}_{V_n}(g)$; this yields the \emphindex{McKay--Thompson series} 
\beq
T_g(\tau)=\sum_{n=-1}^{\infty} \text{Tr}_{V_n}(g) q^n.
\eeq
That is, to each element $g\in \mathbb{M}$ one associates a holomorphic function
\beq
T_g : \mathbb{H} \longrightarrow \mathbb{C}.
\eeq
By virtue of the cyclic properties of the trace these are really \emphindex[class function]{class functions}, meaning that they only depend on the conjugacy class of $g$: 
\beq
T_{g}(\tau)=T_{xgx^{-1}}(\tau), \qquad \forall x\in \mathbb{M}.
\eeq

The $J$-function corresponds to the identity element $e\in \mathbb{M}$ for which $\text{Tr}_{V_n}(e)=\text{dim}\, V_n$. Conway and Norton \cite{MR554399} proposed that all McKay--Thompson series are special types of modular forms known as \emph{Hauptmoduln} for  special \emph{genus zero} subgroups $\Gamma_g$ of $SL(2,\mathbb{R})$. A subgroup $\Gamma\subset SL(2,\mathbb{R})$ is genus zero if the closure of the quotient $\Gamma\backslash \mathbb{H}$ has the topology of a sphere when viewed as a Riemann surface. A hauptmodul for $\Gamma$ then generates the field of rational functions on $\Gamma\backslash \mathbb{H}$.  The conjecture by Conway and Norton is indeed true for the $J$-function which is the hauptmodul for $\Gamma_e=SL(2,\mathbb{Z})$. Frenkel--Lepowsky--Meurman \cite{MR996026} constructed a  \emphindex[vertex operator algebra]{(conformal) vertex  algebra}$V^{\natural}$, with $\text{Aut}(V^{\natural})=\mathbb{M}$, central charge $c=24$, and with graded dimension ($L_0=$ Virasoro generator)
\beq
J(\tau)=\text{Tr}_{V^{\natural}}(q^{L_0-1}),
\eeq
 providing an explanation for monstrous moonshine. 

\begin{remark}
We recall that a \emphindex{vertex algebra} is a vector space with some additional structure, such as a product obeying certain axioms; a \emph{conformal} vertex algebra also carries an action of the \emphindex{Virasoro algebra}, whose zeroth generator $L_0$ gives rise to the grading. Conformal vertex algebras are also known as \emph{chiral algebras} or \emphindex[vertex operator algebra (VOA)!holomorphic]{holomorphic vertex operator algebras} (VOAs).
\end{remark}

\begin{remark}
The monster vertex algebra $V^{\natural}$  has a natural interpretation in physics: it corresponds to a holomorphic conformal field theory describing bosonic string theory on the orbifold $(\mathbb{R}^{24}/\Lambda_{\text{Leech}})/\mathbb{Z}_2$, where $\Lambda_{\text{Leech}}$ is the Leech lattice. From this point of view the modularity of the McKay--Thompson series becomes completely natural: it originates from the symmetry of the path integral over maps from the world-sheet torus, with modified boundary conditions along the spatial direction, into space-time. One missing piece is the so-called \emph{genus zero property}, namely the fact that the invariance groups $\Gamma_g$  of the McKay-Thompson series are all genus zero subgroups of $SL(2,\mathbb{R})$, meaning that $\overline{\Gamma_g\backslash \mathbb{H}}$ has the topology of a sphere. This can not be explained by conformal field theory methods alone since $\Gamma_g$ contains elements, so-called \emph{Fricke elements}, which act as $\tau\to-1/(N\tau)$, $N=\mathrm{order}(g),$ and these lie outside of $SL(2,\mathbb{Z})$. Recently, a new type of string theory compactification was proposed \cite{Paquette:2016xoo}, in which these Fricke symmetries appear naturally as space-time dualities, thus providing a complete physical explanation of monstrous moonshine. This construction also relates monstrous moonshine to BPS-state counting, since in the proposal of \cite{Paquette:2016xoo}, the McKay--Thompson series are reinterpreted as supersymmetric indices counting BPS-states. 
\end{remark}

The complete Conway--Norton conjecture was subsequently proven by Borcherds \cite{MR1172696} who made use of $V^{\natural}$ combined with an earlier  `no-ghost'-theorem in string theory \cite{Goddard:1972iy}. In the course of his proof, Borcherds also introduced new ingredients, notably the \emph{monster Lie algebra} $\mathfrak{m}$,  carrying an action of $\mathbb{M}$ by Lie algebra automorphisms. The construction of $\mathfrak{m}$ was  inspired by physics: it arises as the (degree 1) BRST cohomology of a certain extension of $V^{\natural}$. 
The monster Lie algebra is an example of a \emphindex{Borcherds--Kac--Moody algebra}, which is an infinite-dimensional Lie algebra constructed from a Cartan matrix by generators and relations. The main difference from ordinary Kac--Moody algebras is that Borcherds--Kac--Moody-algebras allows for \emph{imaginary simple roots} (i.e. with non-positive norm). A few years after the Conway--Norton conjectures were put forward, Norton proposed a vast generalisation that he called \emphindex{generalised moonshine}. This pertains to a wider class of modular functions $T_{g,h}(\tau)$ associated to \emph{commuting pairs} of elements $g,h\in \mathbb{M}$. He gave a precise list of properties that these functions should satisfy; in particular for the identity element $g=e$ they  reproduce the McKay--Thompson series. 

Norton's generalised moonshine conjecture has recently been proven in a series of impressive papers by  Carnahan  \cite{MR2728485,MR2904095,2012arXiv1208.6254C}, in particular by constructing a  class of Borcherds--Kac--Moody-algebras $\mathfrak{m}_g$, labelled by  classes $[g]\subset \mathbb{M}$, admitting an action of the centralizer $C_{\mathbb{M}}(g)$ of $g$ in $\mathbb{M}$. Physically, the generalised moonshine functions can conjecturally be understood as follows \cite{Dixon:1988qd}
\beq
T_{g,h}(\tau)=\text{Tr}_{V_g^{\natural}}(hq^{L_0-1}),
\eeq
where $V_g^{\natural}$ denotes the so-called $g$-twisted sector in the orbifold of $V^\natural$ by the element $g\in \mathbb{M}$. For this formula to make sense, $h$ must be a symmetry of $V_g^\natural$, which has automorphism group the centraliser $C_{\mathbb{M}}(g)$, and so indeed $h$ and $g$ must commute.  
 
\subsection{The singular theta correspondence}

We now want to focus on a particular aspect of Borcherds' proof of monstrous moonshine. The starting point was the following famous product formula for the modular invariant $J$-function:
\beq
J(p)-J(q)=p^{-1}\prod_{m>0, n\in \mathbb{Z}}(1-p^m q^n)^{c(mn)},
\label{denominatorJ}
\eeq
where $q=e^{2\pi i \tau}$, $p=e^{2\pi i \sigma}$, $z, \sigma\in \mathbb{H}$. The exponent on the right hand side involves the Fourier coefficients $c(n)$ of the $J$-function itself \eqref{Jfunction}. Borcherds key insight was to interpret this formula as the denominator formula (see \eqref{denominatorformula}) of an infinite-dimensional Lie algebra, called the monster Lie algebra $\mathfrak{m}$. The monster Lie algebra is an example of a class of Lie algebras that are now known as {\it Borcherds--Kac--Moody algebras}, or BKM-algebras for short. The main difference from ordinary Kac-Moody algebras is that BKM-algebras are allowed to have (a possibly infinite number of) imaginary simple roots, beyond the ordinary real simple roots. The denominator formula of a BKM-algebra has the same structure as the one of an ordinary Lie algebra \eqref{denominatorformula}, but with a correction factor due to the imaginary simple roots:
\begin{align}
\sum_{w\in\Weyl}\eps(w) w(e^{\rho} \Sigma)=e^\rho \prod_{\alpha>0} (1-e^{-\alpha})
\label{denominatorformulaBKM}
\end{align}
The factor $\Sigma$ contains a certain combination of the imaginary simple roots, including information about their multiplicities. The monster Lie algebra has only a single real simple root and therefore the Weyl group (which is defined with respect to reflections in the real simple roots only) is $\Weyl =\{\pm 1\}$. The left hand side of the denominator formula \eqref{denominatorformulaBKM} for $\mathfrak{m}$ then has only two terms and one can show that this gives precisely the difference of $J$-functions on the left hand side of \eqref{denominatorJ}. Similarly, the product over positive roots in the denominator formla is identified with the product over $m>0, n\in \mathbb{Z}$ in \eqref{denominatorJ}.

The infinite product then corresponds to a product over all the positive roots of $\mathfrak{m}$. The formula may also be viewed as a lift from a modular form on $\mathbb{H}$, namely the $J$-function, to an automorphic form on $SO(2,2)/(SO(2)\times SO(2))\cong \mathbb{H}\times \mathbb{H}$. It is an example of a more general process known as a \emph{Borcherds lift}, or \emph{multiplicative lift}. This, in turn, can be understood in the context of the theta correspondence. 

Consider the following integral 
\beq
\Theta_J(\sigma, \rho) = \int_{SL(2,\mathbb{Z})\backslash \mathbb{H}} \theta_{\Pi^{2,2}}(\tau;\sigma, \rho) J(\tau)\frac{dx dy}{y^2},
\eeq
where $\theta_{\Pi^{2,2}}$ is the weight 0 Siegel theta series for the four-dimensional indefinite lattice $\Pi^{2,2}$ of signature $(2,2)$, and $(\sigma, \rho)\in SO(2,2; \mathbb{R})/(SO(2)\times SO(2))\cong \mathbb{H}\times \mathbb{H}$. The measure $dxdy/y^2$  is invariant so the integrand is indeed $SL(2,\mathbb{Z})$-invariant for $\tau=x+iy$. This integral has exactly the form of the theta integrals discussed in section \ref{sec:ThetaCorr} above. However, there is a crucial difference. The seed function is now the modular invariant $J$-function $J(\tau)=q^{-1}+\cdots $ which has a pole at the cusp $\tau\to i \infty$. Hence,  the integral is divergent and does not fit the requirements for the theta correspondence. However, it can be evaluated after proper regularisation as was first done by Harvey and Moore in \cite{Harvey:1995fq}. The result is 
\beq
\int_{SL(2,\mathbb{Z})\backslash \mathbb{H}} \theta_{\Pi^{2,2}}(\tau;g) J(\tau)\frac{dx dy}{y^2}=-2\log\Big|e^{-2\pi i \sigma} \prod_{m>0\atop n\in \mathbb{Z}}\big(1-e^{2\pi i (m\sigma+n\rho)}\big)^{c(mn)}\Big|.
\eeq
By comparing with \eqref{denominatorJ} we see that the result of the integral involves the product side of the denominator formula for the $J$-function. By replacing the seed function by the McKay--Thompson series $T_g(\tau)$ one can obtain similar product formulas for elements $g\in \mathbb{M}$. Borcherds later generalised the Harvey-Moore construction into the so-called \emphindex[theta correspondence!singular]{singular theta correspondence} \cite{MR1625724} which deals with arbitrary integrals of the form
\beq
\Theta_\Lambda(F):=\int_{SL(2,\mathbb{Z})\backslash \mathbb{H}} \frac{dx dy}{y^2}\big( \overline{\theta}_\Lambda(\tau), \, F(\tau)\big).
\label{thetalift}
\eeq
Here $\Lambda$ is a lattice, $\tau=x+iy$, $F$ is a  vector-valued modular form of weight $k$ (valued in the group ring $\mathbb{C}[\Lambda^{\vee}/\Lambda]$), $\overline{\theta}_\Lambda$ is a  vector-valued 
Siegel theta series of weight $-k$, and  $(\, , \, )$ denotes the scalar product in the vector space $\Lambda^{\vee}/\Lambda$. Let $e_\gamma$, $\gamma\in \Lambda^\vee/\Lambda$, be a basis for the group ring $\mathbb{C}[\Lambda^\vee/\Lambda]$ so that 
\beq
e_{\gamma}e_{\gamma'}=e_{\gamma+\gamma'},
\eeq
with inner product 
\beq
(e_\gamma, e_{\gamma'})=\delta_{\gamma+\gamma',0}.
\eeq
In what follows we shall restrict to the case of weight $k=0$ for simplicity. A vector-valued modular function $F$ for a congruence subgroup $\Gamma\subset SL(2,\mathbb{Z})$ is then defined as 
\beq
F(\tau)=\sum_{\gamma\in \Lambda^\vee/\Lambda} F_\gamma(\tau)e_\gamma,
\eeq
where the components $F_\gamma(\tau)$ are modular functions for $\Gamma$, transforming in the metaplectic representation of (the double cover of) $SL(2,\mathbb{Z})$ on $\mathbb{C}[\Lambda^\vee/\Lambda]$:
\beqa
{} F_\gamma(\tau+1)&=& e^{2\pi i Q(\gamma)} F_\gamma(\tau),
\nonumber \\ 
{} F_\gamma (-1/\tau)&=& \frac{1}{\sqrt{|\Lambda^\vee/\Lambda|}} \sum_{\gamma'\in \Lambda^\vee/\Lambda} e^{-2\pi i (\gamma, \gamma')} F_{\gamma'}(\tau),
\eqa
where $Q:\Lambda\to \mathbb{Z}$ is the quadratic form on $\Lambda^\vee/\Lambda$ and we also defined the even bilinear form
\beq
(\gamma, \gamma'):=Q(\gamma+\gamma')-Q(\gamma)-Q(\gamma').
\eeq
The vector-valued theta series is defined similarly
\beq
\Theta_\Lambda(\tau)=\sum_{\gamma\in \Lambda^\vee/\Lambda} \theta_{\Lambda+\gamma} e_\gamma, 
\eeq
where $\theta_{\Lambda+\gamma}$ is the ordinary Siegel theta series for the shifted lattice $\Lambda+\gamma$.

In the course of his proof of generalised monstrous moonshine, Carnahan has computed the integral $\Theta_\Lambda$ for the following choice of data. Let $\Lambda=\Pi^{1,1}(N)\times \Pi^{1,1}$, where $\Pi^{1,1}$ is the standard two-dimensional Lorentzian lattice, and $\Pi^{1,1}(N)\subset \Pi^{1,1}$ is the sublattice where the quadratic form is rescaled by $N>0$. Thus, as a free abelian group of rank 4 we can view the lattice as 
\beq
\Lambda\cong \mathbb{Z}\times N\mathbb{Z}\times \mathbb{Z}\times \mathbb{Z}
\eeq
 and for any vector $u=(a,b,c,d)\in \Lambda$ the quadratic form is 
 \beq
 Q(u)=Q(a,b,c,d)=ab+cd.
 \eeq
  The quotient $\Lambda^\vee/\Lambda$ is then identified with $\mathbb{Z}/N\mathbb{Z}\times \mathbb{Z}/N\mathbb{Z}$ and for any vector $v=(a,b)\in \Lambda^\vee/\Lambda$ the quadratic form descends to $Q(v)=Q(a,b)=ab/N$. One can construct a vector-valued modular form $F_g$ by starting from the generalised Moonshine functions $T_{g^{i}, g^{j}}$ in the following way. Let $\gamma=(i,k)\in  \Lambda^\vee/\Lambda$ and construct the components $F_\gamma=F_{i,k}$ by discrete Fourier transform:
\beq
F_{\gamma}(\tau)=F_{i,k}(\tau)=\frac{1}{N}\sum_{j\in \mathbb{Z}/N\mathbb{Z}} e(-jk/N)T_{g^{i}, g^{j}}(\tau).
\label{defDF}
\eeq
One can show under modular transformations this satisfies:
\beqa
F_\gamma(\tau+1)=F_{i,k}(\tau+1)&=& e(ik/N) F_{i,k}(z)=e(Q(\gamma))F_\gamma(\tau)
\nonumber \\
\nonumber \\
F_\gamma(-1/\tau)=F_{i,k}(-1/\tau)&=& \frac{1}{N} \sum_{j,l \in \mathbb{Z}/N\mathbb{Z}}e\left((-jk-il)/N\right) F_{j,l}(\tau)
\nn\\ 
&=&\frac{1}{N}\sum_{\delta\in M^\vee/M} e(-(\gamma, \delta))F_\delta(\tau),
\nonumber \\
\eqa
where $\delta=(j,l)$ so that $(\gamma, \delta)=Q(i+j, k+l)-Q(i,k)-Q(j,l)=(il+jk)/N$. Thus
\beq
F_g(\tau)=\sum_{i,k\in \mathbb{Z}/N\mathbb{Z}} F_{i,k}(\tau)e_{i,k} 
\eeq
is a vector-valued modular form of weight $0$. 
Given this choice of data we have the following theorem due to Carnahan:

\begin{theorem}[Singular theta lift of McKay--Thompson series \cite{MR2904095}] {\it The singular theta lift (\ref{thetalift}) of the vector-valued modular function $F_g(z)$ defined in (\ref{defDF}) is given by}
\beq
\Theta_\Lambda(F_g)=-4 \log \big|T_{1,g}(\rho)-T_{g,1}(\sigma)\big|,
\eeq
where $T_{1,g}$ and $T_{g,1}$  are specializations of the generalised moonshine functions $T_{g,h}$.
\end{theorem}

\begin{remark}[Product formula for McKay-Thompson series]
After evaluating the integral $\Theta_\Lambda(F_g)$ one also obtains an infinite product formula generalising the one for the $J$-function:
\beq
T_{1,g}(\rho)-T_{g,1}(\sigma)=e^{-2\pi i \rho}\prod_{m>0, \,  n\in \tfrac{1}{N}\mathbb{Z}}\big(1-e^{2\pi i (m\rho+n\sigma)}\big)^{c_g(m,n)},
\eeq
where $c_g(m,n)$ is the $q^{mn}$:th coefficient of the vector-valued function $F_{m,Nn}$. For each element $g\in \mathbb{M}$ this is the denominator formula 
for a Borcherds-Kac-Moody algebra $\mathfrak{m}_g$ \cite{MR2728485,MR2904095}.
\end{remark}

\subsection{Mathieu moonshine}

In 2010, a completely new moonshine phenomenon was conjectured by physicists Eguchi,  Ooguri and  Tachikawa \cite{Eguchi:2010ej}, which involves another of the finite sporadic simple groups, namely the \emphindex{Mathieu group} $M_{24}$. Here the role of the $J$-function is played by an object known as the \emphindex{elliptic genus of K3-surfaces}, denoted $\phi_{K3}$. This is a  topological invariant of a family of 2-dimensional complex surfaces, and it is also an example of a modular form, or, more precisely, a \emphindex{Jacobi form}, which is a 2-variable generalisation of a modular form. The key to revealing the presence of $M_{24}$ is to notice that $\phi_{K3}$ affords the following decomposition \cite{Eguchi:2009cq}:
\beq
\phi_{K3}(\tau, z)=\frac{\theta_1(\tau, z)^2}{\eta(\tau)^3}\Big(24\mu(\tau, z)+H(\tau)\Big), \quad \tau\in \mathbb{H}, \, z\in \mathbb{C},
\label{phiK3}
\eeq
where $\theta_1$ is the Jacobi theta function, $\eta$ is the Dedekind eta function and $\mu$ is an Appell-Lerch sum. The key object for us is the function $H(\tau)$ which turns out to be a so-called mock modular form of weight 1/2. See \cite{Dabholkar:2012nd} for  details on mock modular forms and the decomposition (\ref{phiK3}) of $\phi_{K3}$.  The function $H(\tau)$ can be represented explicitly in terms of a Fourier expansion:
\beq
H(\tau)=q^{-1/8}\big(-2+\sum_{n=1}^\infty A_n q^n\big), \qquad q=e^{2\pi i \tau}.
\eeq
One notices now that the first few coefficients $A_n$ are given by:
\beq
A_1=2\times 45, \qquad A_2=2\times 231, \qquad A_3=2\times 770, \qquad A_4=2\times 2277, \dots
\eeq
where $45, 231, 770,\dots $ are all dimensions of representations of $M_{24}$. This observation thus bears many similarities with monstrous moonshine, although the modular object is of a different type.  The `Mathieu moonshine' conjecture has by now been extended and to a certain extent established \cite{Cheng:2010pq,Gaberdiel:2010ch,Gaberdiel:2010ca,Eguchi:2010fg,Eguchi:2011aj,Gaberdiel:2011fg,Gannon:2012ck,Gaberdiel:2012gf,Persson:2013xpa}. However, despite this amazing progress, the main question remains unanswered: \emph{Why does the elliptic genus carry information about the Mathieu group?} This is analogous to McKay's original question that  initiated monstrous moonshine. 

There are several indications that the resolution to this problem can be found in the context of BPS black holes in string theory. Consider string theory compactified on a manifold $X$ which is a product of a K3-surface $S$ and a torus $T^2$. Recall from section \ref{sec:N=4} that this yields a theory on four dimensions with $\mathcal{N}=4$ supersymmetry. In (type IIA) string theory on $X=S\times T^2$ the charge lattice is $\Gamma=\Gamma^{6,22}\oplus \Gamma^{6,22}$ and generic BPS-states have  charge vector $\gamma=(Q,P)$, with electric-magnetic charges $Q,P\in \Gamma^{6,22}$. A generic BPS-state is a \emph{dyon} and preserves 1/4 of the full $\mathcal{N}=4$ supersymmetry. As mentioned in section \ref{sec:N=4} the number of such states is counted (with signs) by the \emph{sixth helicity supertrace}, which is an index defined as follows 
\beq
\Omega_{1/4}(Q,P) = \text{Tr}_{\mathcal{H}_{\mathrm{BPS}}(Q,P)}((-1)^{2J} (2J)^6), 
\eeq
where $J$ is the helicity. This index is invariant under $SL(2,\mathbb{Z})\times SO(6,22;\mathbb{Z})$, where the electric and magnetic charge vectors $Q$ and $P$ transform as a doublet under $SL(2,\mathbb{Z})$, and can only enter through the $SO(6, 22;\mathbb{Z})$-invariant combinations $Q^2, P^2$ and $Q\cdot P$. Remarkably, the partition function that counts these states is known, and is given by the reciprocal of a Siegel modular form, known as the \emphindex{Igusa cusp form} $\Phi_{10}$ \cite{Dijkgraaf:1996it,Shih:2005uc}:
\beq
\frac{1}{\Phi_{10}(\sigma, \tau, z)} = \sum_{m, n, \ell} D(m, n, \ell) p^m q^n y^{\ell},
\eeq
with the identification 
\beq
\Omega_{1/4}(Q,P)  = D\left(\frac{Q^2}{2}, \frac{P^2}{2}, Q\cdot P\right).
\eeq
This gives an exact formula for the degeneracies of $\tfrac14$-BPS black holes in $\mathcal{N}=4$ string theory, and the macroscopic \emphindex{Bekenstein--Hawking--Wald entropy} is reproduced from 
the asymptotic behaviour of the Fourier coefficients for large charges:
\beq
S_{\text{BHW}}=\mathrm{log} \,   \Omega_{1/4}(Q,P)\sim \pi \sqrt{Q^2 P^2- (Q\cdot P)^2} +\cdots, 
\eeq
where the ellipsis indicate quantum corrections to the classical entropy. 

The Igusa cusp form furthermore arises from a singular theta correspondence. One has the following result due to Kawai:
\begin{theorem}[Singular theta  representation of the Igusa cusp form \cite{Kawai:1995hy}]
\beq
\int_{SL(2,\mathbb{Z})\backslash \mathbb{H}} \textrm{Tr}_{\mathcal{H}_{\mathrm{BPS}}}((-1)^{2J}(2J)^6q^{L_0}\bar{q}^{\bar{L}_0}) \frac{dxdy}{y^2}=\log (Y^{10}\big|\Phi_{10}(\sigma, \tau, z)\big|^2),
\eeq
where $Y=\det\, \begin{pmatrix} \sigma & z \\ z & \tau \end{pmatrix}$ is the determinant of the genus 2 period matrix.
\end{theorem}
Evaluating the integral directly yields the product representation 
\beq
\Phi_{10}(\sigma, \tau, z)=e^{2\pi i (\sigma+\tau-z)} \prod_{n,m\geq 0, \ell\mathbb{Z}\atop \ell>0 \text{ when } m=n=0}\big(1-e^{2\pi i (n\sigma+m\tau+\ell z)}\big)^{c(4mn-\ell^2)},
\label{productPHI10}
\eeq
where the exponents $c(4mn-\ell^2)$ are the Fourier coefficients of the K3 elliptic genus
\beq
\phi_{K3}(\tau, z)=\sum_{m\geq 0, \ell\in \mathbb{Z}} c(m,\ell)q^m y^\ell.
\eeq
The infinite product \eqref{productPHI10} corresponds to the right hand side of the denominator formula \eqref{denominatorformulaBKM} of a Borcherds--Kac--Moody-algebra with Cartan matrix (associated with the real simple roots) \cite{1995alg.geom..4006G}:
\beq
\begin{pmatrix} \phantom{-}2 & -2 & -2 \\ -2 & \phantom{-}2 & -2 \\ -2 & -2 & \phantom{-}2 
\end{pmatrix}
\label{CartanmatrixBKM}
\eeq

The product form of $\Phi_{10}$ corresponds to an automorphic (Borcherds) lift of the K3 elliptic genus, analogously to the Borcherds lift of the $J$-function which yields the product in~\eqref{denominatorJ}. 
The full denominator formula is given by the following expression:
\begin{align}
\sum_{m>0} e^{2\pi i m\tau}T_m\phi_{10,1}(\sigma, z)=e^{2\pi i (\sigma+\tau-z)} \prod_{n,m\geq 0, \ell\mathbb{Z}\atop \ell>0 \text{ when } m=n=0}\big(1-e^{2\pi i (n\sigma+m\tau+\ell z)}\big)^{c(4mn-\ell^2)},
\label{denominatorPhi10}
\end{align}
where $\phi_{10,1}$ is the weak Jacobi form of weight 10, index 1, and $T_m$ is the $m$:th Hecke operator (to be precise, the denominator formula for the BKM with Cartan matrix~\eqref{CartanmatrixBKM} is actually the \emph{square} of~\eqref{denominatorPhi10}). See \cite{EichlerZagier} for the precise definition of how the Hecke operator $T_m$ acts on Jacobi forms. 
Just as for the monster Lie algebra, the left hand side encodes the structure of the imaginary simple roots, whose multiplicities are in this case captured by the Fourier coefficients of the Jacobi form $\phi_{10,1}$. The sum side is known as the \emph{additive automorphic lift} (or Saito-Kurokawa-Maass lift from its inventors), in contrast to the right hand side which is a multiplicative lift. 

\begin{remark}
We have seen that the multiplicative lift can be understood within the singular theta correspondence. It is an interesting open question whether there is some alternative version of the (singular) theta lift that encodes also the additive automorphic lift.
\end{remark}

It has been speculated (starting with \cite{Harvey:1995fq,Harvey:1996gc}) that the space of BPS-states in $\mathcal{N}=4$ string theory carries an action of the Mathieu group $M_{24}$, thereby providing a possible explanation for the Mathieu moonshine phenomenon. Further support for this was given in \cite{Cheng:2010pq,Persson:2013xpa,Persson:2015jka}, where a class of functions 
\beq
\Phi_{g,h}\, : \mathbb{H}^{(2)} \longrightarrow \mathbb{C}
\eeq
on the Siegel upper half plane $\mathbb{H}^{(2)}$ was defined for each commuting pair of elements $g,h\in M_{24}$ (see \cite{Persson:2013xpa} for the precise definition of $\Phi_{g,h}$). It was proven that these functions are Siegel modular forms with respect to subgroups $\Gamma_{g,h}^{(2)}\subset Sp(4;\mathbb{R})$, and they play a key role in generalised Mathieu moonshine \cite{Gaberdiel:2012gf}. When restricted to $(g,h)=(1,1)$, where $1$ is the identity element of $M_{24}$, this reproduces the Igusa cusp form $\Phi_{10}=\Phi_{1,1}$. A key observation was that $\Phi_{g,h}$ satisfies the surprising property:
\beq
\Phi_{g,h}(\sigma, \tau, z)=\Phi_{g,h}(\tfrac{\tau}{N}, \sigma, z),
\label{emduality}
\eeq
where $N$ is the order of the element $g$. What could be the physical interpretation of the functions $\Phi_{g,h}$? In the follow-up paper \cite{Persson:2015jka} it was demonstrated that the symmetry (\ref{emduality}) corresponds precisely to a novel type of electric-magnetic duality in certain orbifold compactifications of heterotic string theory on $T^6$. In these theories the most general index counting $h$-twisted black hole states in the $g$-orbifold theory \cite{Sen:2010ts}:
\beq
\Omega_{1/4}^{g,h}(Q,P)=\text{Tr}_{\mathcal{H}_{\mathrm{BPS},g}(Q,P)}((-1)^{2J} (2J)^6 h),
\eeq
where the trace is taken over the $g$-orbifold sector $\mathcal{H}_{\mathrm{BPS},g}(Q,P)$. This suggests a very natural interpretation of all the functions $\Phi_{g,h}$. Fourier expanding its reciprocal one obtains:
\beq
\frac{1}{\Phi_{g,h}(\sigma, \tau, z)}=\sum_{m, n, \ell} D_{g,h}(m, n, \ell) p^m q^n y^{\ell},
\eeq
This leads to the following conjecture:

\begin{conjecture}[Counting of twisted dyons \cite{Persson:2013xpa,Persson:2015jka}]
For each $g\in M_{24}$ and $h\in C_{M_{24}}(g)$ the twisted black hole states in (generalised) $g$-orbifold CHL-models are counted by the Fourier coefficients of $1/\Phi_{g,h}$, i.e.
\beq
\Omega_{1/4}^{g,h}(Q,P)=D_{g,h}\left(\frac{Q^2}{2}, \frac{P^2}{2}, Q\cdot P\right).
\label{CHLconjecture}
\eeq
\end{conjecture}

These results certainly indicate that BPS-black holes will ultimately play a role in explaining Mathieu moonshine. 

\begin{remark}
A possible connection between BPS-states and Mathieu moonshine was also observed in recent work of Katz, Klemm and Pandharipande \cite{Katz:2014uaa} (see also \cite{Cheng:2015kha}). They studied certain stable pair invariants (Gopakumar--Vafa invariants) of type II string theory on a K3-surface, and found a numerical coincidence with dimensions of $M_{24}$-representations. This also prompted recent ongoing work by Harvey and Moore, aimed at proving (or disproving) the assertion that $M_{24}$ is the automorphism group of the algebra of space-time BPS-states of string theory on K3 \cite{MooreTalk}.
\end{remark}

\subsection{Representation theoretic interpretation?}
It is natural to ask whether the automorphic forms which appear in monstrous moonshine can be understood in the context of automorphic representations. Put differently, can parts of monstrous moonshine be understood in the context of the Langlands program? These questions have in particular been raised by Borcherds on several occasions. The first thing to note is that the automorphic infinite products that arise in the singular theta correspondence has poles at cusps and are therefore not Hecke eigenforms in general. This means that the Hecke algebra acts freely and Flath's theorem (see theorem~\ref{thm:Flath}) does not apply. Hence, if there were an underlying automorphic representation it is most likely not factorizable. Another way to see this is to try to construct $L$-functions associated with automorphic infinite products. If one naively defines these in the standard way via Dirichlet series formed from their Fourier coefficients one ends up with divergent series. Borcherds has proposed that instead one might define $L$-functions via some regularised version of the Mellin transform. However, to the best of our knowledge, this has not yet been realised. We are therefore forced to admit that a satisfactory interpretation of the singular theta correspondence in general, and  moonshine in particular, is still an open problem.

\section{Extension to Kac--Moody groups}
\label{sec:KM}

In this article we have concentrated on the study of automorphic forms and in particular Eisenstein series defined on finite-dimensional Lie groups, as categorised in the Cartan classification. In this section we give a short summary of what happens as one makes the extension to infinite-dimensional\index{group!Kac--Moody} \emphindex{Kac--Moody groups}, which are generated by an infinite number of raising and lowering operators. A complete classification of Kac--Moody groups is at present not known and we will restrict our attention here mainly to Eisenstein series defined on\index{Kac--Moody groups!hyperbolic}  \emphindex[Kac--Moody groups!affine]{affine, hyperbolic and Lorentzian Kac--Moody groups}. Full accounts of Kac--Moody algebras can be found in the books~\cite{Kac,MoodyPianzola,WakimotoBook,KumarBook,GoddardOlive}. The motivation for studying Eisenstein series defined on infinite-dimensional Kac--Moody groups is twofold.

\subsection{String theory motivation: infinite-dimensional U-duality}

In string theory, Kac--Moody groups appear, for example,  in the list of  U-duality groups encoding discrete symmetries of type II string theory compactified on a $(10-D)$-dimensional torus from ten down to $D$ space-time dimensions. The list of these groups was given in table~\ref{tab:CJ} in chapter~\ref{ch:intro-strings} and consists of the groups in the exceptional series of the Cartan classification, where in $D\geq3 $ dimensions the respective U-duality group is given by the finite-dimensional and discrete group $E_{11-D}(\mathbb{Z})$. In $D=2,1$ and $0$ dimensions, however, the corresponding U-duality groups are infinite-dimensional and are conjectured to be given by the affine, hyperbolic and Lorentzian Kac--Moody groups $E_9(\mathbb{Z})$, $E_{10}(\mathbb{Z})$ and $E_{11}(\mathbb{Z})$, respectively~\cite{Hull:1994ys}. In particular,  the groups $E_{10}$ and $E_{11}$ are of special relevance~\cite{Damour:2002cu,West:2001as}, since they have been conjectured as fundamental symmetries of \emphindex{M-theory}~\cite{Witten:1995ex,Polchinski}, a theory whose low-energy limit is eleven-dimensional supergravity and from which the five different known types of string theories can be derived as particular limits.

\subsection{Gravity motivation: quantum cosmological billiards}

There is also an amusing connection of the theory of Kac--Moody automorphic forms to quantum cosmology. Classical cosmology, as the study of the evolution of the Universe in the form of a dynamical space-time subject to Einstein's field equations, often assumes the existence of a space-like singularity like the \emphindex{big bang}. Similar space-like singularities arise in the interior of Schwarzschild and Kerr black holes. In seminal work devoted to the study of Einstein's partial differential equations, Belinskii, Khalatnikov and Lifshitz (BKL)\index{BKL analysis}\index{Belinskii--Khalatnikov--Lifshitz analysis} discovered that $(i)$ there is most probably a regime in the vicinity of the singularity where one effectively use ordinary differential equations and $(ii)$ the behaviour of the four-dimensional empty (matterless) Universe in this set-up becomes chaotic~\cite{Belinsky:1970ew,Belinsky:1982pk}, see also~\cite{Misner:1969hg}. 

This was later generalised to the context of maximal supergravity by Damour and Henneaux~\cite{Damour:2000wm,Damour:2000hv} where it was found that similar phenomena persist and moreover a Weyl group structure related to the hyperbolic Kac--Moody group $E_{10}$ was found, see also~\cite{Damour:2002et,Henneaux:2007ej} for reviews. The Weyl group determines the shape of a certain effective billiard system, giving rise to the name \emphindex{cosmological billiard} for this field of theoretical cosmology. Moreover, in~\cite{Damour:2002cu} it was shown that there is a full $E_{10}(\reals)$ structure for part of the ordinary differential equations arising in this analysis, lending strong support to the conjecture that Kac--Moody symmetries are fundamental for string theory and M-theory.

At the quantum level, it is interesting to contemplate how much of the Kac--Moody structure remains relevant. In a first step, one can try to quantise the cosmological billiard system associated with the Weyl group acting on some generalised upper half plane. This was undertaken in~\cite{Kleinschmidt:2009cv,Kleinschmidt:2009hv,Kleinschmidt:2010bk} where it was shown that the spectrum of the Laplace operator in this context is such that one obtains a quantum-mechanical resolution of the space-like singularity, see also~\cite{Graham:1990jd,Ivashchuk:1994fg,Kirillov:1995wy,Pioline:2002qz,Forte:2008jr}. This \emphindex{quantum cosmological billiards} analysis, however, should properly be embedded into a full quantum gravitational analysis that goes beyond the crude approximations used in these first investigations. At any rate, it is to be expected that automorphic forms will play a central role in the final formulation of this model.

\subsection{Mathematical motivation: new automorphic \texorpdfstring{$L$}{L}-functions}

From a mathematician's perspective one motivation to study Kac--Moody Eisenstein series is to consider them as a potential source for deriving new $L$-functions in a Fourier through a Fourier expansion of the series. It is however not precisely clear if this extension of the theory of Eisenstein series will yield necessarily to new types of $L$-functions and the focus of the discussion has so far been on series defined on affine groups. In fact in~\cite{ShahidiInfinite} an argument was provided that no new functions will be found, while in~\cite{GarlandLS} a new method for obtaining such functions was devised. This new method relies on an expansion of the series with respect to ``lower triangular parabolics'', instead of only ``upper triangular parabolics''. 

In recent years there has been some work, developing the theory of Eisenstein series for Kac--Moody groups. The most well developed part is that of Eisenstein series defined on affine Kac--Moody groups which was started by Garland~\cite{Garland1,GarlConvergence}. While for the finite-dimensional groups convergence with respect to the (complex) defining weight $\lambda$ was proven over almost all of the complex plane, c.f.~\eqref{absoluteconvergence}, for the infinite-dimensional Kac--Moody groups convergence is restricted and the defining weight has lie inside the Tits cone~\cite{Kac}. Furthermore, a restriction on the group element forming the argument of the series has to be imposed~\cite{GarlConvergence}. First steps towards a definition of Eisenstein series on hyperbolic Kac--Moody groups have been made in~\cite{CarboneHyp}, where the case of rank $2$ hyperbolics was considered. Furtheromore, in~\cite{FK2012,FKP2013} Eisenstein series defined on the hyperbolic $E_{10}$ group (along with $E_9$ and $E_{11}$), were discussed. A general proof of convergence of Eisenstein series on general hyperbolics remains to be developed, however.

\subsection{Fourier coefficients and small representations}
Despite the absence of a mathematically rigorous definition of Kac--Moody Eisenstein series, quite a bit can be said about the Fourier expansion of these series. The foundation for this work was laid in~\cite{Garland1}, where the analogue of Langlands' formula~\eqref{LCF} for the constant term, was developed for the case of affine Kac--Moody Eisenstein series. While for Eisenstein series on finite-dimensional groups we have explained how to evaluate Langlands' formula in section~\ref{EvalLCF}, it is not clear how to apply this formula in the case of affine Kac--Moody Eisenstein series. The reason for this is that the sum over Weyl words appearing in the formula is an infinite sum due to the infinite-dimensional nature of affine groups and their associated Weyl groups. The same problem of course also appears when considering the extension of Langlands' formula for the cases of hyperbolic and Lorentzian Kac--Moody groups. This question was taken up in~\cite{FK2012}, see also~\cite{FK2012Summary} for a summary of this work, where it was shown that for special types of Kac--Moody Eisenstein series, the naively infinite sum `collapses' to finite sum and can be explicitly~\textit{computed}. On a more technical level, to evaluate Langlands' formula, one proceeds just as in the case of a finite-dimensional group and one successively constructs Weyl words in the set $\mathcal{C}(\lambda)$ by the orbit method, c.f. section~\ref{EvalLCF}. It can then be shown that for particular types of Eisenstein series, which we will discuss in a moment, only the coefficients $M(w,\lambda)$ associated with the first few Weyl words $w$ in the carefully constructed orbit, are non-zero. All other coefficients associated with the infinite number of Weyl words that follow in the orbit are however zero and therefore do not contribute to the constant term.  

The Eisenstein series for which this collapse of the constant term happens are of the general form~\eqref{ESeriesMaxParabDefn} defined on the groups $E_9$, $E_{10}$ and $E_{11}$. More specifically, the defining weight $\lambda$ is of the form $\lambda=2s\Lambda_1-\rho$, such that the series is defined with respect to the maximal parabolic subgroup $P_1$ associated with the first node of the $E_{11-D}$ Dynkin diagram in Bourbaki labelling. In order to observe collapse, the generically complex parameter $s$ has to take real values of $s=\tfrac32$, $\tfrac52$, \ldots, see~\cite{FK2012}  for an extended list. 

\begin{remark}
Let us mention as an aside that these particular Eisenstein series appear as the automorphic couplings of the $R^4$ and $D^4R^4$ curvature correction terms in the low-energy effective action of type IIB string theory in $D=2,1$ and $0$ dimensions, in line with the discussion in section~\ref{sec:outlook-strings}. Since the $R^4$ and $D^4R^4$ term are $\frac12$- and $\frac14$-BPS protected terms, it is a reassuring confirmation that the corresponding constant terms only contain a finite number of perturbative contributions. Furthermore, these Eisenstein series are associated with small representations which we discuss in some detail in section~\ref{smallreps}.
\end{remark}

Developing an understanding of the structure of Fourier coefficients of Kac--Moody Eisenstein series is an open problem and part of ongoing work. However, it is possible to apply formula~\eqref{degW1} for the degenerate Whittaker coefficients also to Kac--Moody Eisenstein series (with slight modifications in the affine case) and use this to make some statements about the Fourier coefficients. This was done in~\cite{FKP2013}, where the formula for the degenerate Whittaker coefficients~\eqref{degW1} was derived, applied to the cases of the particular maximal parabolic Eisenstein series just mentioned above and explicit expressions for the Whittaker coefficients were computed. In particular it was found that in the case of $s=\tfrac32$ the (abelian) Fourier coefficients are completely determined by maximally degenerate Whittaker coefficients. The collapse property discussed for the case of the constant term above, also plays a central role in computing these Whittaker coefficients.   
Related work, with a focus on the rank $2$ affine case, can be found in~\cite{LoopGroupsLiu}. We formalise the observations of~\cite{FK2012,FKP2013} as:
\begin{conjecture}[Small representations for Kac--Moody groups]
Kac--Moody groups possess a minimal unitary representation that can be realised automorphically. In the case of $E_n(\reals)$ (for $n\geq 9$)  this can be achieved by inducing from the maximal parabolic subgroup of with semi-simple Levi group of type $D_{n-1}$ that is obtained by deleting the first node of the $E_n$ Dynkin diagram. The canonically associated Eisenstein series for $s=\tfrac32$ (obtained by analytic continuation) is the spherical vector in the minimal automorphic representation. The wave-front set is of Bala--Carter type $A_1$.

A similar next-to-minimal representation is obtained for $s=\tfrac52$ and its wave-front set is of Bala--Carter type $2A_1$.
\end{conjecture}

\subsection{Langlands program for Kac--Moody groups?}
Braverman and Kazhdan have also started to develop the local theory for affine Kac--Moody groups (see \cite{2012arXiv1205.0870B} for a survey). In particular, they have constructed the local spherical Hecke algebra \cite{MR2846488}, as have Gaussent and Rousseau~\cite{GaussentRousseau}. With Patnaik they have proven an affine version of the Satake isomorphism \cite{2014arXiv1403.0602B}. These results were recently used in \cite{2012arXiv1212.6473B} to prove a Gindikin--Karpelevich formula for affine Kac--Moody groups and in \cite{2011arXiv1101.4912K,2014arXiv1407.8072P} the Casselman--Shalika formula has been generalised to the affine setting.

As formulated in \cite{2012arXiv1205.0870B}, \emph{the dream} is to have a fully developed representation theory and an associated Langlands correspondence for any (symmetrizable) Kac--Moody group. Although at present this remains a dream, the recent developments reviewed above certainly provides hope that such a theory is within reach.

\appendix
\cftaddtitleline{toc}{chapter}{Appendices}{}

\chapter{Fourier expansion of \texorpdfstring{$SL(2,\reals)$}{SL(2, R)} series by Poisson resummation}
\label{app:SL2Fourier}

In this appendix we perform the Fourier expansion of the series (\ref{Eisenintro})
\begin{align}
\label{Eisenapp}
f_s(z) = \sum_{(c,d)\in\ints^2\atop(c,d)\neq(0,0)} \frac{z_2^s}{|cz+d|^{2s}}
\end{align}
that is related to the standard $SL(2,\reals)$ Eisenstein series through $f_s(z)=2\zeta(2s) E(s,z)$, cf. (\ref{SL2lattice}). Here, $z=x+i y$ lies on the upper half plane $\UHP$ as defined in section~\ref{sec:SL2}.

The invariance of $f_s(z)$ under shifts $z\to z+1$ implies that it should have a Fourier expansion
\begin{align}
f_s(z)
=  C(y) + \sum_{m\neq 0} a_m(y) e^{2\pi i m x} \,.
\end{align}
The `constant term(s)' $C(y)$ and the non-zero Fourier coefficients $a_m(t)$ are determined in the following. We suppress the label $s$ on the constant terms and Fourier coefficients for ease of notation.

The technique to be used rests on Poisson resummation whose fundamental equation here is (cf. \cite[Eqn. (8.2.210]{Polchinski})
\begin{align}
\label{poisson1}
\sum_{m\in\ints} \exp(-\pi a m^2 + 2\pi i b m) = a^{-1/2} \sum_{\tilde{m}\in\ints} \exp\left(-\frac{\pi(\tilde{m}-b)^2}{a}\right)\,.
\end{align}
Another useful form of this same formula is
\begin{align}
\label{poisson2}
\sum_{m\in\ints} \exp\left(-\frac{\pi}{t} (m+n x)^2\right) 
  = t^{1/2} \sum_{\tilde{m}\in\ints} \exp\left(-\pi t \tilde{m}^2 - 2\pi i \tilde{m}nx\right)\,.
\end{align}
Note that the sums are over all integers and not constrained to a single $SL(2,\ints)$-orbit.

We will also use the following representation of powers for $\Re(s)>0$ and $\Re(M)>0$
\begin{align}
\label{intpower}
M^{-s} = \frac{\pi^s}{\Gamma(s)} \int_0^\infty \frac{dt}{t^{s+1}} e^{-\frac{\pi}{t} M}\,.
\end{align}
Finally, we require the following integral representation of the modified Bessel function for real $a,b\neq 0$
\begin{align}
\label{BesselK}
\int_0^\infty  \frac{dt}{t^{s+1}} e^{-\pi t a^2 -\frac{\pi}{t} b^2} = 2 \left|\frac{a}{b}\right|^{s} K_s(2\pi|a b|)\,.
\end{align}

\section{Constant term(s)}

First extract the term $c=0$ from (\ref{Eisenapp}). Then $d\neq 0$ and 
\begin{align}
\label{const1}
f_s(z) &= y^s \sum_{d\neq 0} |d|^{-2s} + \underbrace{y^s \sum_{c\neq 0} \sum_{d\in\ints} |cz+d|^{-2s}}_{f_s^{(1)}(z)}
= 2\zeta(2s) y^s+ f_s^{(1)}(z)\,.
\end{align}

The power $|cz+d|^{-2s}$ appearing in the second term can be rewritten as an integral using (\ref{intpower}). Then one can Poisson resum over $d\in\ints$ using (\ref{poisson2}):
\begin{align}
\label{afterpoisson}
f_s^{(1)}(z) &= \frac{\pi^s}{\Gamma(s)}y^s\sum_{c\neq 0}\sum_{d\in\ints} \int_0^\infty \frac{dt}{t^{s+1}} \exp\left(-\frac{\pi}{t} |cz+d|^2\right)\\
 &= \frac{\pi^s}{\Gamma(s)}y^s\sum_{c\neq 0}\sum_{d\in\ints} \int_0^\infty \frac{dt}{t^{s+1}} \exp\left(-\frac{\pi}{t} ((cx+d)^2 +(cy)^2)\right)\nn\\
&= \frac{\pi^s}{\Gamma(s)}y^s\sum_{c\neq 0} \sum_{\tilde{d}\in\ints}\int_0^\infty \frac{dt}{t^{s+1}} t^{1/2} \exp\left(-\pi t \tilde{d}^2 -2\pi i \tilde{d}cx -\frac{\pi}{t} (cy)^2\right)\,.\nn
\end{align}

In the final line of (\ref{afterpoisson}) one can separate out the term with $\tilde{d}=0$ by
\begin{align}
\label{const2}
f_s^{(1)} (z)=  \frac{\pi^s}{\Gamma(s)}y^s\sum_{c\neq 0} \int_0^\infty\frac{dt}{t^{s+1/2}} \exp\left(-\frac{\pi}{t} (cy)^2\right) + f_s^{(2)}(z)
\end{align}
since it does not have any $x$ dependence and where $f_s^{(2)}$ are the terms with $\tilde{d}\neq 0$:
\begin{align}
\label{nonzeromodes}
f_s^{(2)}(z) = \frac{\pi^s}{\Gamma(s)}y^s\sum_{c\neq 0} \sum_{\tilde{d}\neq 0}\int_0^\infty \frac{dt}{t^{s+1/2}}  \exp\left(-\pi t \tilde{d}^2 -2\pi i \tilde{d}cx -\frac{\pi}{t} (cy)^2\right)\,.
\end{align}
The integral in the term with $\tilde{d}=0$ can be undone using (\ref{intpower}) and the sum over $c\neq 0$ can be carried out afterwards. Hence the first term in (\ref{const2}) becomes
\begin{align}
 &\quad\frac{\pi^s}{\Gamma(s)}y^s\sum_{c\neq 0} \int_0^\infty\frac{dt}{t^{s+1/2}} \exp\left(-\frac{\pi}{t} (cy)^2\right)
 &&=  \frac{\pi^s}{\Gamma(s)}\frac{\Gamma(s-1/2)}{\pi^{s-1/2}} y^{s-2(s-1/2)} \sum_{c\neq 0} c^{-2(s-1/2)}\nn\\
 &=2\zeta(2s) \frac{\pi^{-(s-1/2)}\Gamma(s-1/2)\zeta(2s-1)}{\pi^{-s}\Gamma(s)\zeta(2s)} y^{1-s}
 &&= 2\zeta(2s) \frac{\xi(2s-1)}{\xi(2s)}y^{1-s}\,,
\end{align}
where we have pulled out the same  overall factor as in (\ref{const1}) and regrouped the $\pi$-factors to use the definition of the completed Riemann zeta function $\xi(k)=\pi^{-k/2} \Gamma(k/2) \zeta(k)$.

\section{Non-zero Fourier modes}

The current status of the Fourier expansion is then
\begin{align}
f_s(z) = 2\zeta(2s)\left( y^s + \frac{\xi(2s-1)}{\xi(2s)}y^{1-s}\right) + f_s^{(2)}(z)\,,
\end{align}
with the non-zero mode part $f_s^{(2)}$ given by (\ref{nonzeromodes}). The $t$-integral appearing in that expression is a Bessel integral and can be evaluated using (\ref{BesselK}) as
\begin{align}
f_s^{(2)}(z) &= \frac{2\pi^s}{\Gamma(s)}y^s \sum_{c\neq 0}\sum_{\tilde{d}\neq 0} 
\left|\frac{\tilde{d}}{ny}\right|^{s-1/2}K_{s-1/2}(2\pi|\tilde{d}c|y) e^{-2\pi i \tilde{d}cx}\nn\\
&= \frac{2\pi^s}{\Gamma(s)} y^{1/2} \sum_{c\neq 0}\sum_{\tilde{d}\neq 0} 
\left|\frac{\tilde{d}}{c}\right|^{s-1/2}K_{s-1/2}(2\pi|\tilde{d}c|y) e^{-2\pi i \tilde{d}c x}\,.
\end{align}
To find the Fourier coefficient $a_m(y)$ of a mode $e^{2\pi i m x}$ we transform the double summation to one over $m\neq 0$ and the (positive) divisors  $d|m$. Then
\begin{align}
f_s^{(2)}(z) &= \frac{4\pi^s}{\Gamma(s)} y^{1/2} \sum_{m\neq 0}\sum_{d|m}
d^{1-2s} |m|^{s-1/2} K_{s-1/2}(2\pi|m|y) e^{2\pi i m x}\nn\\
&=2\zeta(2s) \frac{2y^{1/2}}{\xi(2s)} \sum_{m\neq 0} |m|^{1/2-s} \sigma_{2s-1}(m) K_{s-1/2}(2\pi|m|y) e^{2\pi i m x}\,,
\end{align}
again pulling out the same overall factor $2\zeta(2s)$ and using the divisor sum\index{divisor sum}
\begin{align}
\sigma_s(m) = \sum_{d|m} d^s
\end{align}
where only positive divisors are included. 

The full Fourier expansion is therefore given by
\begin{align}
\label{fullfourier}
f_s(z) = 2\zeta(2s)\left[ y^s + \frac{\xi(2s-1)}{\xi(2s)}y^{1-s}  + \frac{2y^{1/2}}{\xi(2s)} \sum_{\neq 0}|m|^{1/2-s} \sigma_{2s-1}(m) K_{s-1/2}(2\pi |m|y) e^{2\pi i m x}\right].
\end{align}
The term in the square brackets is the full expansion of the Eisenstein series $E(s,\tau)$ for $SL(2,\reals)$. This agrees with the adelic derivation of theorem~\ref{SL2AESexp}.

\chapter{Laplace operators on \texorpdfstring{$G/K$}{G/K} and automorphic forms}
\label{App:Laplace}

In this appendix, we briefly review the connection between the \emphindex[Laplace operator]{scalar Laplace operator} on the symmetric space $G(\reals)/K(\reals)$ and the quadratic Casimir. We do this first for a general simple, simply-laced split group $G(\reals)$ and then give a very explicit analysis for $G=SL(2,\reals)$.

\section{Scalar Laplace operator and quadratic Casimir}

For a simple, simply-laced split $G(\reals)$ we denote by $\mf{h}$ a fixed Cartan subalgebra of the Lie algebra $\mf{g}(\reals)$ of $G(\reals)$. With respect to $\mf{h}$ and a choice of simple roots $\alpha_i$ ($i=1,\ldots,r$ with $r=\dim_\reals(\mf{h})$) the remaining generators arrange into positive and negative step operators, cf.~\eqref{SL2triplets}. We denote by $E_\alpha$ the step operator of a given root $\alpha$. In \emphindex{Iwasawa gauge} we choose to write an arbitrary element $g\in G(\reals)/K(\reals)$ as
\begin{align}
\label{Iwagauge}
g = na= \exp\left(\sum_{\alpha>0} u_\alpha E_\alpha\right)\prod_{i=1}^r v_i^{h_i},
\end{align}
where $h_i$ are the Cartan generators associated with the choice of simple roots~cf.~\eqref{simpleChev}.   The variables $v_i$ (for $=1,\ldots,r$) and $u_\alpha$ (for $\alpha\in\Delta_+$) are coordinates on the symmetric space $G(\reals)/K(\reals)$. 

The $G(\reals)$-invariant metric on the symmetric space can be constructed from
\begin{align}
\label{cosmet}
ds_{G/K}^2 = 2\langle \mathcal{P} | \mathcal{P} \rangle,
\end{align}
where we chose a convenient normalisation and
\begin{align}
\mathcal{P} = \frac12\left(g^{-1} dg - \theta(g^{-1} dg)\right)
\end{align}
is the coset projection of the \emphindex{Maurer--Cartan form} $g^{-1}dg$ associated with the vector space decomposition $\mf{g}=\mf{p}\oplus \mf{k}$. Here, $\mf{k}$ is the Lie algebra of $K$. The \emphindex[Cartan involution]{(Cartan) involution} $\theta$ leaving $\mf{k}$ fixed can be defined by
\begin{align}
\theta (E_{\alpha}) = -E_{-\alpha},\quad
\theta (h_i) =-h_i.
\end{align}
With this convention, $\mf{k}$ and $\mf{p}$ have the bases
\begin{align}
\mf{k} &=  \langle E_{\alpha}-E_{-\alpha}\,|\, \alpha>0 \rangle,\nn\\
\mf{p} &=  \langle E_{\alpha}+E_{-\alpha}\,|\, \alpha>0 \rangle \oplus \langle h_i\,|\, i=1,\ldots, r\rangle.
\end{align}
We further choose the normalisation ($A_{ij}$ is the Cartan matrix (\ref{CartanMatrix}) of the simply-laced $\mf{g}(\reals)$)
\begin{align}
\label{normD}
\langle E_{\alpha} | E_{-\beta} \rangle = \delta_{\alpha,\beta},\quad
\langle h_i | h_j \rangle = A_{ij}.
\end{align}
Working out the Maurer--Cartan form for the element (\ref{Iwagauge}) one finds
\begin{align}
\label{MCF}
g^{-1} dg &= \sum_{i=1}^r v_i^{-1}dv_i h_i + a^{-1} \left(\sum_{\alpha>0} Du_\alpha E_{\alpha}\right) a\nn\\
&= \sum_{i=1}^r v_i^{-1} dv_i h_i + \sum_{\alpha>0}  a^{-\alpha} Du_\alpha E_\alpha
\end{align}
where  $D u_\alpha = du_{\alpha} + \ldots$
and the dots represent finitely many terms coming from commutator terms when expanding out the \emphindex{Baker--Campbell--Hausdorff identity} 
\begin{align}
e^{-X} d(e^{X}) = dX -\frac1{2!} \lb X, dX \rb +\frac1{3!} \lb X, \lb X, dX\rb\rb + \ldots
\end{align}
for the nilpotent $E_{\alpha}$. The expression (\ref{MCF}) together with (\ref{normD}) leads to a block-diagonal metric of the form
\begin{align}
ds_{G/K}^2 = g_{\mu\nu} dx^\mu dx^\nu =2 \sum_{i,j=1}^r v_i^{-1} v_j^{-1} dv_i dv_j A_{ij} +  \sum_{\alpha>0} a^{-2\alpha} (Du_\alpha)^2.
\end{align}

The scalar Laplacian associated with this metric is ($\partial_\mu \equiv\frac{\partial}{\partial x^\mu}$ and $\sqrt{g}=\sqrt{\det(g_{\mu\nu})}$)
\begin{align}
\label{GKlapl}
\Delta_{G/K} &= \frac1{\sqrt{g}} \partial_{\mu} \left(\sqrt{g}g^{\mu\nu} \partial_{\nu}\right)\nn\\
&= \frac12\sum_{i,j=1}^r (A^{-1})^{ij} a^{2\rho}v_i \partial_i \left(a^{-2\rho}v_j\partial_j \right) + \sum_{\alpha>0} a^{2\alpha} \partial_\alpha^2 +\ldots,
\end{align}
where the dots come from inverting the metric in the $du_\alpha du_\beta$ sector and $(A^{-1})^{ij}$ is the inverse of the Cartan matrix. We have used the relation $\sum_{\alpha>0} \alpha= 2\rho$ for the Weyl vector, cf.~\eqref{Weylvec}. The Laplace operator~\eqref{GKlapl} is $G(\reals)$-invariant since the Maurer--Cartan form trivially is: The transformation of $g\in G(\reals)/K(\reals)$ is given by $g\to g_0 g k^{-1}$ with constant $g_0\in G(\reals)$ and $k\in K(\reals)$ such that $g^{-1}dg\to k (g^{-1}dg) k - dk k^{-1}$ is independent of $g_0$.

We can evaluate the eigenvalue of the Laplacian~\eqref{GKlapl} when acting on an Eisenstein series $E(\lambda,g)$ as defined in ~\eqref{generalEisenstein}. Due to the invariance of the Laplacian, it suffices to evaluate it on the summand $\chi(g)=\chi(a)=a^{\lambda+\rho}$, corresponding to $\gamma=\id$. For this term, the derivatives $\partial_\alpha$ with respect to the coordinates $u^\alpha$ vanish and one finds
\begin{align}
\Delta_{G/K} a^{\lambda+\rho} &= \Delta_{G/K} \prod_{i=1}^r  v_i^{2s_i}
= \frac12\sum_{i,j=1}^r (A^{-1})^{ij}  2s_i (2s_j-2) a^{\lambda+\rho}\nn\\ 
&= \frac12\left(\langle\lambda|\lambda\rangle - \langle\rho|\rho\rangle\right) a^{\lambda+\rho},
\end{align}
where we stress that we assumed $\mf{g}$ to be simply-laced. As already mentioned, $G(\reals)$-invariance implies that this is also the eigenvalue for the full Eisenstein series:
\begin{align}
\Delta_{G/K} E(\lambda,g) &= \frac12 \left(\langle\lambda|\lambda\rangle - \langle\rho|\rho\rangle\right) E(\lambda,g).
\end{align}
This agrees up to a factor with the standard quadratic Casimir evaluated on a lowest weight representation with lowest weight $\Lambda=\lambda+\rho$~\cite{Kac}.

\section{Automorphic forms on \texorpdfstring{$SL(2,\reals)$}{SL(2, R)} as Laplace eigenfunctions}
\label{app:SL2Laps}

For the case of $SL(2,\reals)$ we can give fully explicit expressions. Using
\begin{align}
g = na = \exp( u e) v^h
\end{align}
one finds from~\eqref{cosmet}
\begin{align}
ds^2_{G/K} = 4 v^{-2}dv^2 + e^{-4v} du^2
\quad \Rightarrow \quad
\Delta_{G/K} = \frac14 e^{2v} v\partial_v \left(e^{-2 v} v\partial_v\right) + e^{4v} \partial_u^2.
\end{align}
This can be brought into a more familiar form by using $v= y^{1/2}$ and $u=x$, cf.~\eqref{SL2IWA:app}. This leads to
\begin{align}
\Delta_{G/K} =y^2 \left( \partial_x^2+ \partial_y^2\right),
\end{align}
which agrees with the Laplacian on the upper half plane $\UHP$ given in~\eqref{eq:SL2Lapapp}. 

Consider now a real eigenfunction $\varphi(z)$ of the Laplace operator with eigenvalue $s(s-1)$:
\begin{align}
\label{laplapp}
\Delta_{G/K} \varphi(z) = s(s-1) \varphi(z).
\end{align}
If the function is furthermore invariant under $SL(2,\ints)$, this implies $\varphi(z)=\varphi(z+1)$ and one has a Fourier expansion of the form
\begin{align}
\label{Fapp}
\varphi(z) =\sum_{m\in\ints } a_m(y) e^{2\pi i m x} 
\end{align}
where $m\in\ints$ denotes the `instanton charge' of the character in the terminology of section~\ref{sec:FCsec} and $a_0(y)$ is the constant term. This Fourier expansion is due to the translations $x\to x+1$ contained in the action of $SL(2,\ints)$ acting on $SL(2,\reals)$. Reality of $\varphi(z)$ implies that $a_m(y) = a_{-m}(y)$ for all $m>0$. We therefore restrict to $m\geq 0$.

Plugging the Fourier expansion~\eqref{Fapp} into the Laplace equation~\eqref{laplapp} and analysing each mode individually leads to the following equations
\begin{subequations}
\begin{align}
\label{sl2zero}
m=0 &:& y^2 \partial_y^2 a_0(y) &= s(s-1) a_0(y),\\
\label{sl2nonzero}
m\neq 0 &:& y^2 \left(\partial_y^2 -4\pi^2m^2\right)a_m(y) &= s(s-1) a_m(y).
\end{align}
\end{subequations}
The equation (\ref{sl2zero}) for the constant term has two linearly independent solutions 
\begin{subequations}
\begin{align}
s\neq \frac12&:& a_0(y) &= y^s &\textrm{or}&& a_0(y) &= y^{1-s},\\
\label{shalfsols}
s=\frac12&:& a_0(y) &= y^{1/2} &\textrm{or}&& a_0(y) &= y^{1/2}\log y.
\end{align}
\end{subequations}
All these solutions are at most power laws when $y$ approaches any cusp, e.g. $y \to\infty$.

Equation (\ref{sl2nonzero}) for the non-zero modes becomes more familiar when one uses $a_m(y) = y^{1/2} b_m(y)$ which leads to
\begin{align}
y^{2} \partial_y^2 b_m(y) + y \partial_y b_m(y) -\left(4\pi m^2y^2+\left(s-\tfrac12\right)^2 \right) b_m(y) =0.
\end{align}
After a rescaling of the $y$ coordinate this becomes the modified Bessel equation with the two modified Bessel functions $K_{s-1/2}$ and $I_{s-1/2}$ as linearly independent solutions. Translated back to $a_m(y)$ these are
\begin{align}
a_m(y) = y^{1/2} K_{s-1/2}(2\pi |m| y)\quad\textrm{or}\quad
a_m(y) = y^{1/2} I_{s-1/2}(2\pi |m| y).
\end{align}
If one insists on at most power law growth near the cusp $y\to\infty$ the solution involving the function $I_{s-1/2}$ is disallowed. This is an instance of the `multiplicity one theorem' mentioned in chapter~\ref{ch:fourier}.

Putting everything together, we see that any real function $\varphi(g)$ on $SL(2,\reals)$ that is right-invariant under $SO(2,\reals)$ and satisfies the three conditions stated for automorphic forms in the introduction can be expanded as
\begin{align}
\varphi(z) = a_0^{(s)} y^s + a_0^{(1-s)} y^{1-s} + y^{1/2}\sum_{m\neq 0} a_m K_{s-1/2}(2\pi|m| y) e^{2\pi i m x}
\end{align}
with $a_m= a_{-m}$ and these are purely numerical coefficients. For cusp forms one has that the numerical coefficients $a_0^{(s)}$ and $a_0^{(1-s)}$ vanish identically. The above expansion is valid for $s\neq \tfrac12$; for $s=\tfrac12$ one has to replace the constant terms by the solutions of~\eqref{shalfsols}.

As shown in section~\ref{sec:SL2-Hecke}, the coefficients $a_m$ can also be determined for cusp forms if one demands in addition to the Laplace condition that $\varphi(z)$ is also an eigenfunction of all the Hecke operators. These can be thought of as the analogues of the Laplace operator for finite $p<\infty$ and therefore an automorphic function that obeys simple equations for all $p\leq \infty$ is uniquely fixed (up to an overall normalisation), cf. remark~\ref{rmk:globalHecke} and example~\ref{ex:HeckeConnection}.

\chapter{Local-to-global principle}\label{appendixLocGlob}
\label{sec_Hasse}
In this appendix we provide some background on the local-to-global principle, also known as Hasse's principle. This principle forms the basis of a powerful approach to problems in arithmetic which we are going to illustrate in the following.

The local-to-global principle is nicely motivated by the study of algebraic equations. For instance, consider the polynomial
\begin{align}
f(x)=7x^3-2x+2
\end{align}
which has integer coefficients.  Then the question of whether the equation $f(x)=0$ has any integer solutions can be answered by considering a reduction of the polynomial's coefficients modulo $3$ leaving us with
\begin{align}
\tilde f (x)=x^3+x+2\,.
\end{align}
Since it is only necessary to test three integers, one quickly verifies that $\tilde{f}(x)=0$ possesses no solution in $\mathbb Z/3\mathbb Z$. Now, since $\mathbb Z\rightarrow \mathbb Z/n\mathbb Z$ is a ring homomorphism, any solution of $f(x)=0$ in $\mathbb Z$ is mapped to a solution of $\tilde f(x)=0$ in $\mathbb Z/3\mathbb Z$ we have also shown that $f(x)=0$ has no integer solutions.

Considering an algebraic equation over $\mathbb Z$, such as the one above, modulo a prime number $p$ is referred to as seeking solutions~\textit{locally} in other words in $\mathbb Z/p\mathbb Z$.  If no solution is found locally it is then possible to deduce also that no solution exists~\textit{globally} namely in $\mathbb Z$. The approach presented above provides a simple example for what is known as the local-to-global principle. Although the method has worked nicely in this simple example it actually has limited applicability. In order to develop a more powerful method that realises the local-to-global principle it is useful to introduce $p$-adic numbers.

To this end consider for instance the quartic equation
\begin{align}
7y^4=2^2x^6+2\cdot 7^3 x^2+7^4
\end{align}
and we ask the question whether this equation possesses any non trivial solutions $(x,y)\in\mathbb Q^2$. In a first attempt we may try to proceed in an analogous way to above and reduce the equation modulo $2$ or modulo $7$. In the first case this leaves us with
\begin{align}
y^4=1
\end{align}
which implies that $(x,\pm1)$ is a solution for all $x\in\mathbb Q$.
From the second reduction we obtain
\begin{align}
4x^6=0
\end{align}
implying as solution $(0,y)$ for all $y\in\mathbb Q$. Using the simple method of reducing the equation with respect to prime numbers thus does not help us in answering the question posed. Instead we will now use $p$-adic numbers to show that no solution to the equation exists in $\mathbb Q_p$ and since $\mathbb Q\hookrightarrow \mathbb Q_p$ also no solution can exist in $\mathbb Q$.

Specifically we will work with $p=7$ and define the $7$-adic valuations of the variables $x$ and $y$ as $\nu_7(x)=n$ and $\nu(y)_7=m$ with $m,n\in\mathbb Z$. Taking the $7$-adic valuation of the left-and right-hand-sides of our algebraic equation we thus find
\begin{align}
\nu_7(\text{l.h.s.})\equiv\nu_7(7y^4)=1+4m
\end{align}
and from the right-hand-side
\begin{align}
\nu_7(\text{r.h.s.})\equiv\nu_7(2^2x^6+2\cdot 7^3 x^2+7^4)\geq\min(4,3+2n,6n)\,,
\end{align}
where we have used property~\eqref{valprop} of the $p$-adic valuation. In fact, since none of the arguments of the  minimum function are equal, the inequality sharpens to an equality, $\nu_7(\text{r.h.s.})=\min(4,3+2n,6n)$. The value of the $7$-adic valuation of the left-hand-side $\nu_7(\text{l.h.s.})$ is odd and we would thus require that $\nu_7(\text{r.h.s.})=3+2n$.  
However note that for $n\geq1$ we have that $\nu(\text{r.h.s.})=4$ and for $n<1$ we have $\nu_7(\text{r.h.s.})=6n$. Either case is in contradiction with the value of the $7$-adic valuation of the left-hand-side and we thus conclude that no solution in $\mathbb Q_7$ exists. As a consequence also no solution in $\mathbb Q$ exists, providing another realisation of the local-to-global principle.

Even though the introduction of $p$-adic numbers improves our ability to analyse algebraic equations, there is still a major  limitation to our analysis. In particular neither of the above methods is able to prove the existence of a global solution, meaning a solution in $\mathbb Q$. In other words, if we find a solution of an algebraic equation locally in $\mathbb Q_p$ for some prime $p$, this does not imply that a solution in $\mathbb Q$ exists. In fact one can give examples of algebraic equations where a local solution exists for every prime $p$, but no global solution exists. See for instance~\cite{Selmer}.

Nevertheless, in some cases one can go further and prove theorems which provide information about the existence of global solutions from the local analysis. An example of such a case is the Hasse--Minkowski theorem for quadratic forms. Although the theorem holds for general number fields $K$ we will content ourselves to stating the theorem for quadratic forms over the rationals which we also define for completeness. 
\begin{definition}[Quadratic form]
A quadratic form $f$ over $\mathbb Q$ is a polynomial of degree two in the variables $x_i\in K$ with $i=1,...,n$, where $n$ is called the rank of the form. Given some $y\in \mathbb Q$, we say that the quadratic form $f$ represents $y$ if there exists a solution $(X_1,...,X_n)\in \mathbb Q^n$ with $(X_1,...,X_n)\neq(0,...,0)$, such that
\begin{align}
y=f(X_1,...,X_n)\,.
\end{align}
\end{definition}
With this definition we can then state the Hasse-Minkowski theorem which applies to non-degenerate quadratic forms. 
\begin{theorem}[Hasse--Minkowski]
Let $f$ be a quadratic form over $\mathbb Q$ and for a prime number $p$ let $f_p$ be the form over $\mathbb Q_p$. Then $f$ represents zero if and only if $f_p$ represents zero for all prime numbers $p$ including the prime at infinity. 
\end{theorem}
Put differently, the theorem states that for $f$ to have a global zero it is necessary and sufficient for $f$ to have a local zero at all places. We refer the reader to~\cite{MR0344216} for more details and a proof of this theorem.

\chapter{Poincar\'e series and Kloosterman sums}
\label{app:PSKS}

In this appendix, we briefly review some classic material related to the Fourier expansion of Poincar\'e series on the upper half plane that are invariant under $G(\ints) = SL(2,\ints)$. References for this exposition include~\cite{Iwaniec}.

\section{Poincar\'e series}

The primary ingredient for a Poincar\'e series is the \emphindex{seed function} $\sigma\,:\, \UHP \to \cx$ that is invariant under discrete shifts:
\begin{align}
\sigma(z+1) = \sigma(z),
\end{align}
or more group-theoretically
\begin{align}
\sigma(\gamma z) = \sigma(z) \quad
\textrm{for all }\quad \gamma\in B(\ints)=\left\{\begin{pmatrix}\pm 1 & n \\&\pm 1\end{pmatrix} \st n\in\ints\right\}\subset SL(2,\ints).
\end{align}
Note that the seed function can depend on both the real and imaginary part of $z$ and it is therefore more general than the character that enters in the construction of non-holomorphic Eisenstein series in~\eqref{PoincareSL2}. The periodicity assumption, however, entails that the seed $\sigma(z)$ has a Fourier expansion
\begin{align}
\label{eq:seedFE}
\sigma(z) = \sigma_0(y) + \sum_{n\neq 0} \sigma_n(y) e^{2\pi i n x}.
\end{align}

The general form of a \emphindex[Poincar\'e series for $SL(2,\ints)$]{Poincar\'e series} is 
\begin{align}
f(z)\equiv P(\sigma,z) = \sum_{\gamma\in B(\ints)\backslash G(\ints)} \sigma(\gamma z).
\end{align}
This sum makes sense for suitably fast decaying seeds. The Poincar\'e series is invariant under $\Gamma=SL(2,\ints)$ by construction and we would like to say as much as possible about its Fourier expansion
\begin{align}
f(z) = \sum_{n\in \ints} f_n(y) e^{2\pi i n x}.
\end{align}

\section{Fourier expansion}

We will tackle the Fourier expansion of $f(z)$ by starting from the Poincar\'e series representation and performing yet another right quotient to study
\begin{align}
B(\ints) \backslash SL(2,\ints) /B(\ints).
\end{align}
Choosing a unique representative then includes a sum over $k\in\ints$ for the right quotient that will be very useful.

The left quotient can be parametrised in standard way in terms of co-prime $c$ and $d$ such that
\begin{align}
\begin{pmatrix}a&b\\c&d\end{pmatrix}
\end{align}
is any fixed matrix representing the single coset class (i.e. one fixes arbitrarily any $a$ and $b$ that satisfy $ad-bc=1$ over $\ints$). Performing the right quotient then allows to bring $d$ in the range between $0$ and $c-1$ since it includes translates by $c$ in the lower right component. This is true for $c\neq 0$; for $c=0$ one has to have $d=1$ and the coset is represented by the identity. Let  $(\ints/c\ints)^\times$ denote the invertible integers mod $c$; these are all the numbers $d$ co-prime with $c$. 

Therefore
\begin{align}
f(z) &= \sum_{\gamma\in B(\ints)\bs SL(2,\ints)} \sigma(\gamma z) \nn\\
&= \sigma(z) + \sum_{c>0} \sum_{d\in (\ints/c\ints)^\times} \sum_{k\in \ints} \sigma\left(\frac{a(z+k)+b}{c(z+k)+d}\right)\nn\\
&=\sigma(z) + \sum_{c>0} \sum_{d\in (\ints/c\ints)^\times} \sum_{k\in \ints} \sigma\left(\frac{a}{c} -\frac{1}{c(c(z+k)+d)}\right)
\end{align}
where the first term comes from the $c=0$ coset. The $k$-sum represents the right quotient. $a$ here is a chosen upper left entry of the matrix in the left quotient. 

Next we use the Poisson summation formula
\begin{align}
\sum_{k\in\ints} \varphi(z+k) = \sum_{n\in\ints} \int_{\reals} \varphi(z+\omega) e^{-2\pi i \omega n} d\omega.
\end{align}
Applied to $f(z)$ we get
\begin{align}
f(z) =\sigma(z) + \sum_{c>0} \sum_{d\in (\ints/c\ints)^\times} \sum_{n\in \ints} \int_\reals e^{-2\pi i \omega n}\sigma\left(\frac{a}{c} -\frac{1}{c^2((z+\omega)+d/c))}\right)d\omega .
\end{align}
Now letting $z=x+iy$ and introducing $\tilde\omega= \omega+x+d/c$ one gets
\begin{align}
f(z) &=\sigma(z) + \sum_{c>0} \sum_{d\in (\ints/c\ints)^\times} \sum_{n\in \ints} e^{2\pi i n (x+d/c)} \lint_\reals e^{-2\pi i \tilde\omega n} \sigma\left(\frac{a}{c} - \frac1{c^2(\tilde\omega+iy)}\right) d\tilde\omega.
\end{align}
This rearrangement has completely brought out the Fourier mode $e^{2\pi i n x}$ in the second term and so one finds the Fourier modes of the Poincar\'e sum to be
\begin{align}
\label{eq:fFE2}
f_n(y) = \sigma_n(y) + \sum_{c>0} \sum_{d\in(\ints/c\ints)^\times} e^{2\pi i n \frac{d}{c}} \int_\reals e^{-2\pi i \tilde\omega n} \sigma\left(\frac{a}{c} - \frac1{c^2(\tilde\omega+iy)}\right) d\tilde\omega.
\end{align}

In the next step we use the Fourier expansion of the seed $\sigma(z)$ given in~\eqref{eq:seedFE} that here leads to
\begin{align}
\sigma\left(\frac{a}{c} - \frac1{c^2(\tilde\omega+iy)}\right) = \sum_{m\in\ints} \sigma_m\left(\frac{y}{c^2(\tilde\omega^2+y^2)}\right) e^{2\pi i m \frac{a}{c}}e^{-2\pi i m \frac{\tilde\omega}{c^2(\tilde\omega^2+y^2)}}.
\end{align}
Substituting this back into~\eqref{eq:fFE2} yields
\begin{align}
f_n(y) &= \sigma_n(y) + \sum_{c>0}\sum_{m\in\ints} \sum_{d\in(\ints/c\ints)^\times} \!\!\!\!\!e^{2\pi i n \frac{d}{c}+2\pi i m\frac{a}{c}} \!\lint_\reals e^{-2\pi i \tilde\omega n-2\pi i m \frac{\tilde\omega}{c^2(\tilde\omega^2+y^2)}}\sigma_m\left( \frac{y}{c^2(\tilde\omega^2+y^2)}\right) d\tilde\omega\nn\\
&= \sigma_n(y) + \sum_{c>0}\sum_{m\in\ints} S(m,n;c)\lint_\reals e^{-2\pi i \tilde\omega n-2\pi i m \frac{\tilde\omega}{c^2(\tilde\omega^2+y^2)}}\sigma_m\left( \frac{y}{c^2(\tilde\omega^2+y^2)}\right) d\tilde\omega,
\end{align}
where we have carried out the sum over $d\in(\ints/c\ints)^*$ in terms of the \emphindex{Kloosterman sum}
\begin{align}
S(m,n;c) = \sum_{d\in(\ints/c\ints)^\times} e^{2\pi i (mx +nx^{-1}) /c}.
\end{align}
Here, we have used that $ad-bc =1$ which means that $ad\equiv 1\mod c$ and therefore $a\equiv d^{-1}\mod c$. (This is the meaning of $x^{-1}$ in the definition above.) Note that $S(m,n;c)=S(n,m;c)$.

For the zero mode $f_{n=0}(y)$ of $f(z)$ we therefore obtain
\begin{align}
\label{eq:zeromode}
f_0(y) = \sigma_0(y) + \sum_{c>0}\sum_{m\in\ints} S(m,0;c)\int_\reals e^{-2\pi i m \frac{\tilde\omega}{c^2(\tilde\omega^2+y^2)}}\sigma_m\left( \frac{y}{c^2(\tilde\omega^2+y^2)}\right) d\tilde\omega.
\end{align}
For the zero mode, the Kloosterman sum simplifies to the \emphindex{Ramanujan sum}. However,  this integral is not possible to evaluate in general.

\section{The case of Eisenstein series}

Let evaluate the Fourier modes for the case of a standard non-holomorphic Eisenstein series. In this case, the seed is a character on the Borel and has the Fourier expansion
\begin{align}
\sigma(z) = y^s \quad\quad
\Rightarrow \quad\quad
\sigma_m(y) =\left\{ \begin{array}{cl} y^s & \textrm{if $m=0$}\\ 
0 &\textrm{if $m\neq0$} \end{array}\right.
\end{align}

The constant term~\eqref{eq:zeromode} then collapses to the $m=0$ term and reads
\begin{align}
f_0(y) &= y^s + \sum_{c>0} S(0,0;c) c^{-2s} \int_\reals \left(\frac{y}{\tilde\omega^2+y^2}\right)^s d\tilde\omega\nn\\
&= y^s +  \sqrt{\pi}\frac{\Gamma(s-1/2)}{\Gamma(s)} y^{1-s} \sum_{c>0} \varphi(c) c^{-2s} 
\end{align}
The integral converges for $\textrm{Re}(s) >\frac12$. The Euler totient function $\varphi(c)=S(0,0;c)$ gives the cardinality of $(\ints/c\ints)^\times$, i.e., the number of non-zero integers less than $c$ that are also co-prime with $c$. It is a standard result for the Dirichlet series that~\cite[Eq.~(2.28)]{Iwaniec} 
\begin{align}
 \sum_{c>0} \varphi(c) c^{-2s} = \frac{\zeta(2s-1)}{\zeta(2s)}.
\end{align}
In this way one recovers the standard constant terms. The sum converges for $\textrm{Re}(s)>1$. For $s=1$ it diverges. Combining the constant terms leads to
\begin{align}
f_0(y) = y^s + \frac{\xi(2s-1)}{\xi(2s)} y^{1-s}
\end{align}
with the completed Riemann zeta function $\xi(s) = \pi^{-s/2} \Gamma(s/2) \zeta(s)$.

For seeing the $y$-dependence of the constant term, the following change of variables is sufficient ($t=\tilde\omega/y$)
\begin{align}
\int_\reals \left(\frac{y}{\tilde\omega^2+y^2}\right)^s d\tilde\omega = y^{1-s} \int_\reals (1+t^2)^{-s} dt.
\end{align}
(Note that the integrand coincides with the spherical vector of the principal series representation that the Eisenstein series belongs to.)

For the non-zero mode one also recognises quickly the Bessel integral. The number-theoretic sum in this case is more complicated but yields of course the standard divisor sum~\cite[Eq.~(2.27)]{Iwaniec}. In some more detail
\begin{align}
f_n(y) &= \sum_{c>0}  S(0,n;c)c^{-2s} \int_\reals e^{-2\pi i n \tilde\omega} \left(\frac{y}{\tilde\omega^2+y^2}\right)^sd\tilde\omega\nn\\
&= \sum_{c>0}  S(0,n;c) c^{-2s} y^{1-s}\int_\reals e^{-2\pi i n y t} \left(1+t^2\right)^{-s} dt\nn\\
&= \sum_{c>0}  S(0,n;c) c^{-2s} y^{1-s} \frac{2\pi^s}{\Gamma(s)} |ny|^{s-1/2} K_{s-1/2}(2\pi|n| y)\nn\\
&= \frac{2\pi^s}{\Gamma(s)}  |n|^{s-1/2} y^{1/2} K_{s-1/2}(2\pi|n| y) \sum_{c>0}  S(0,n;c) c^{-2s} \nn\\
&= \frac{2\pi^s}{\Gamma(s)}  |n|^{s-1/2} y^{1/2} K_{s-1/2}(2\pi|n| y) \frac{\sigma_{1-2s}(n)}{\zeta(2s)}\nn\\
&= \frac{2}{\xi(2s)} \sigma_{1-2s}(n)|n|^{s-1/2} y^{1/2} K_{s-1/2}(2\pi|n| y),
\end{align}
the familiar result for $SL(2)$ Eisenstein series.

%BIBTEX

\cleardoublepage
\phantomsection
\addcontentsline{toc}{chapter}{References}
{\small
\bibliography{reviewbib}
\bibliographystyle{utphys-sort}
}

\newpage
\fancyhead[CE]{\textit{Index}}
\printindex

\end{document}